\documentclass[11pt]{livreMetMat}
\ifx\pdfpageheight\undefined\PassOptionsToPackage{dvips}{graphicx}\else%
\PassOptionsToPackage{pdftex}{graphicx}
\PassOptionsToPackage{pdftex}{color}
\fi
\ifx\pdfpageheight\undefined\newcommand{\texorpdfstring}[2]{#1}
\else
\usepackage[bookmarksopen=false,pdftex=true,breaklinks=true,backref=page,%
    pagebackref=true,plainpages=false,%
    hyperindex=true,pdfstartview=FitH,colorlinks=true,pdfpagelabels=true,%
    colorlinks=true,linkcolor=blue,citecolor=red]{hyperref}
\fi
\usepackage{graphicx}
\usepackage{amsfonts}
\usepackage[french]{minitoc}
\usepackage{array}
\usepackage{makeidx}
\usepackage[french]{varioref}
\usepackage[francais]{babel}
\usepackage{float}  
\usepackage{A--definitions}

\date{\bf Version mise \`a jour en, avril 2012}  
\author{ABDELJAOUED Jouna\"{\i}di \\[1mm]
Ma\^{\i}tre Assistant \\[1mm]
Universit\'e  de Tunis\\[5mm]
LOMBARDI Henri \\[1mm]
Ma\^{\i}tre de Conf\'erences \\[1mm]
Universit\'e  de Franche-Comt\'e, Besan\c{c}on\\[10mm]}  

\title{\bf\Huge M\'ethodes Matricielles \\[2mm]  
Introduction \`a la Complexit\'e Alg\'ebrique\\[5mm]} 
   
\makeindex  

\begin{document}  
\renewcommand\figurename{Figure}
\renewcommand\listfigurename{Liste des figures}
\renewcommand\listtablename{Liste des tableaux}
\renewcommand\indexname{Index des termes}
\renewcommand\tablename{Tableau}

\thispagestyle{empty}
\maketitle
\pagenumbering{roman}
\setcounter{tocdepth}{2}
\setcounter{page}{0}
\thispagestyle{empty}
\cleardoublepage 
\pagestyle{headings}



\chapter*{Avant-Propos}
\addcontentsline{toc}{chapter}{Avant-Propos}
\markboth{Avant-Propos}{Avant-Propos}

L'\agr \lin nous a sembl\'e \^etre le
terrain id\'eal pour une introduction simple et
p\'edagogique aux outils modernes de la \com
\agq d\'evelopp\'es durant les trois derni\`eres d\'ecennies.

\ms Le tournant en mati\`ere de \coag  en \agr \lin fut la d\'ecouverte par Strassen
\cite{Str}, en 1969, d'un fait d'une simplicit\'e d\'econ\-cer\-tante, mais d'une
port\'ee consid\'erable, \`a savoir que la \mul  de deux \macas d'ordre deux pouvait
se faire avec seulement sept \muls (non commutatives) au lieu de huit dans l'anneau
de base. Ce qui ramenait la \com asymptotique de la \mul de deux \macas d'ordre
\,$n$\, \`a \,$\O(n^{\log_{2}7})$\, au lieu de \,$\O(n^3)$\, et faisait descendre
pour la premi\`ere fois l'exposant de \,$n$\, au-dessous de $3$, alors que les
recherches ant\'erieures n'avaient r\'eussi qu'\`a r\'eduire le \coe de \,$n^3$\,
dans le nombre d'\oparis \ncrs pour calculer le produit de deux \macas d'ordre
\,$n$\, (\cf \cite{Cop}).

\ms Depuis, de nombreux outils ont \'et\'e d\'evelopp\'es. Des notions nouvelles
sont apparues comme celles de \com \bil et de rang tensoriel utilis\'ees de
mani\`ere intensive notamment par Bini, Pan, Sch\"onhage, Strassen, Winograd et
d'autres  (\cf \cite{Bin,BCLR,Win,Pan,Scho2,Stra2}) pour r\'eduire l'exposant
\,$\alpha$: \`a l'heure actuelle, on sait que  \,$\alpha < 2,376$. Il est cependant
conjectur\'e que la borne inf\'erieure des exposants \,$\alpha$\, acceptables serait
$2$, \cad que pour tout \,$\vep >0$\, le produit de deux \macas d'ordre \,$n$\,
pourrait \^etre calcul\'e  par un \cari de taille \,$\O(n^{2+\vep})$\, et de \prof
\,$\O(\log n)$. Cependant ces  \metsz, d'un int\'er\^et th\'eorique certain, sont
\`a l'heure actuelle inapplicables \`a cause notamment de la constante d\'emesur\'ee
que le \gui{$grand~\O$}  cache (\cf \cite{Knu} \S~4.6.4). Par contre la \met de
Strassen a pu trouver une impl\'ementation concr\`ete \cite{Chan}, et elle commence
\`a battre la \mul \usle (dite \cvlez) \`a partir de \,$n=70$.

\ms Le {calcul \paralz} est une technique, en plein d\'eveloppement, qui distribue
un calcul \`a faire sur un grand nombre de processeurs travaillant au m\^eme moment,
en
\paralz. Pour la \mul rapide de \macasz, si le nombre de processeurs disponibles est
suffisamment grand (de l'ordre de \,$\O(n^{\alpha})$), le temps de calcul est alors
extr\^{e}mement faible (de l'ordre de  \,$\O(\log n)$\, pour des matrices sur un
corps fini).

\ms La \mul rapide des \macas a de nombreuses applications en \agr \lin sur les
corps, par exemple l'inversion d'une \maca peut se faire en  \,$\O(n^{\alpha})$\,
avec le m\^eme exposant. Cependant, contrairement \`a la \mul rapide des matrices,
ces \algos ne sont pas bien  adapt\'es au \emph{calcul \paralz}. Ainsi l'agorithme
d'inversion d'une \maca auquel on vient de faire allusion, et que nous \'etudierons
dans la section \ref{sec detinv}, ne voit jamais son temps de calcul descendre en
dessous d'un \,$\O(n\,\log n)$.

\ms C'est sur la base de r\'esultats parfois anciens qu'on a pu exhiber, en  \agr
\linz, des \algos bien adapt\'es au calcul
\paralz, s'appuyant sur la \mul rapide des matrices. Ces \algos
sont en outre des \algos sans divisions (ou presque) et s'appliquent donc \`a des
\acomsz. 

\ss C'est le cas en particulier de la \met d\'evelopp\'ee en 1847 par
l'astronome fran\c{c}ais Le Verrier am\'elior\'ee, un si\`ecle plus tard, par
Souriau, Frame et Faddeev qui l'utilisent pour le calcul des d\'eterminants, du
\polcarv pour l'inversion des matrices, et pour la r\'eso\-lu\-tion des \slisz.
Cette \met s'est av\'er\'ee porteuse d'un \algo tr\`es bien adapt\'e au calcul
\paralz, d\^u \`a Csanky, qui en 1976 a construit, dans le cas d'un \acom contenant
le corps des rationnels, une \fam de \carisz,  pour calculer en \,$\O(\log^2n)$\,
\'etapes \parals les \coes du \polcarz.

\ss Une autre \metz, dite \emph{de partitionnement} (\cite{Gas} pp.~291--298) et
attribu\'ee \`a Sa\-muel\-son \cite{Sam} (1942), a eu un regain d'int\'er\^et avec
l'\agb \cite{Ber}, qui fournit un calcul rapide, \paral et sans division, du
\polcarz. Cet \algo a permis de \gnr aux \acomas le r\'esultat de Csanky concernant
la \com \paralz, par une voie tout \`a fait diff\'erente. Nous en pr\'esenterons une
version \paral am\'elior\'ee (section \ref{sec Berkopar}).

La version \sqle la plus simple de l'\agb n'utilise pas de produits de matrices
mais seulement des produits d'une matrice par un vecteur.\\ Elle s'est  av\'er\'ee
tout \`a fait efficace sur les ordinateurs usuels, et particu\-li\`erement bien
adapt\'ee au cas des matrices creuses.

\ms Nous pr\'esentons dans cet ouvrage les principaux \algos en \agr \linz, et
donnons plus particuli\`erement un aper\c{c}u d\'etaill\'e des diff\'e\-ren\-tes
\mets utilis\'ees pour le calcul du \polcarz, avec des r\'esultats r\'ecents.

\ss L'int\'er\^et port\'e au \polcar d'une matrice est justifi\'e par le fait que la
d\'etermination de ses coefficients suffit \`a conna\^{\i}tre le d\'eterminant de
cette matrice et \`a calculer son adjointe. Dans le cas des corps cela permet de
calculer son inverse et de r\'esoudre les \syses \linsz. Il r\'eside \egmt dans les
renseignements que cela donne sur une forme quadratique, comme par exemple sa
signature dans le cas du corps des r\'eels.

\bni\textbf{Plan de l'ouvrage}

\ms Nous faisons quelques rappels d'\agr \lin dans le chapitre
\ref{Chap RapAlgLin}.

\ss Le chapitre \ref{chap BasicAlgoAlin} contient quelques \mets
classiques cou\-ramment utilis\'ees pour le calcul du \polcarz:
l'\ajbz, la \mhbz, la \miLz, l'\algo de Le Verrier et son
am\'e\-lio\-ration par Souriau-Faddeev-Frame, la \met de Samuelson
modifi\'ee \`a la \Berz, en \gnl la plus efficace, la \met de
Chistov qui a des performances voisines, et enfin des \mets
reli\'ees aux \srlsz, les plus efficaces sur les corps finis.

\ss  Le chapitre \ref{chap circuits} d\'eveloppe
le formalisme des \caris (ou \prevsz) pour une description
formelle des   calculs \agqsz. Nous y expliquons la
technique importante d'\elidz, elle aussi invent\'ee par
Strassen.

\ss Dans le chapitre \ref{ChapNoCo}
nous  donnons un aper\c{c}u des principales notions
de \com les plus couramment utilis\'ees.
Ces notions constituent une tentative de th\'eoriser les
calculs
sur ordinateur, leur temps d'ex\'ecution et
l'espace m\'emoire qu'ils occupent.

\ss
Dans le chapitre \ref{chap divpar} nous expliquons la
strat\'egie
\gnle \gui{diviser pour gagner},
bien adapt\'ee  au calcul \paralz.
Nous donnons quelques
exemples de base.

\ss Le chapitre \ref{chap multipoly} est consacr\'e \`a la \mul
rapide des \polsz, avec la \met de Karatsuba et la Transform\'ee
de Fourier Discr\`ete.

\ss Le chapitre \ref{chap multimat} est consacr\'e
\`a  la \mul rapide des  matrices.
Nous y abordons notamment les notions fondamentales
de \com \bilz, de rang tensoriel et
de \cabasz.

\ss Le chapitre \ref{chap AlinCorps} est consacr\'e \`a des \algos dans lesquels
intervient la \mul rapide des matrices, mais sans que l'ensemble de l'\algo soit
bien adapt\'e au calcul \paralz. 

On obtient ainsi en \gnl les \pcds les plus
rapides connues en ce qui concerne le \emph{temps \sql asymptotique}, pour la
plupart des \pbs classiques li\'es \`a l'\agr \linz. Ces performances sont en \gnl
obtenues uniquement sur les corps. Seule la derni\`ere section du chapitre,
consacr\'ee \`a l'\algo de Kaltofen-Wiedemann concerne le calcul sur un \acomaz.

\ss
Le chapitre \ref{chap par Lever}  pr\'esente
les \parans
de la \met de Le Verrier, qui s'appliquent dans tout \acom
o\`u les entiers sont non diviseurs de z\'ero
et o\`u la division par un entier, quand elle est possible,
 est explicite.

\ss Le chapitre \ref{chap PolCarAnn} est consacr\'e
 aux \mets \parals de Chistov
et de Berkowitz qui s'appliquent en toute \gnlz it\'e.

\ss Le chapitre \ref{chap experim} pr\'esente tout d'abord
quelques tableaux r\'ecapitulatifs des \coms des diff\'erents
\algos \'etudi\'es, \sqls ou \paralsz, pour le calcul du \deter et
celui du \polcarz. Nous donnons ensuite les r\'esultats des tests
exp\'erimentaux concernant quelques \mets \sqles du calcul du
\polcarz. Ces r\'esultats montrent des performances sup\'erieures
pour les \algos de Chistov et de Berkowitz avec un l\'eger
avantage pour ce dernier.

\ss Les deux derniers chapitres sont consacr\'es aux travaux
de Valiant
sur un analogue \agq de la conjecture $\P\neq \NP$,
dans lesquels le d\'eterminant et le permanent occupent une
place centrale.
Bien qu'on ait tr\`es peu d'id\'ees sur la mani\`ere
de r\'esoudre la conjecture de Valiant  $\VP\neq \VNP$,
celle-ci semble quand m\^eme moins hors de port\'ee que
la conjecture \algq dont elle s'inspire.

\ss L'annexe contient
 les codes {\sc Maple} des \algos exp\'eriment\'es.
  Nous avons choisi le logiciel de Calcul Formel
 {\sc Maple} essentiellement pour des raisons de commodit\'e.
Le langage de programmation qui
 lui est rattach\'e est proche de celui de nombreux autres langages
 classiques, permettant de d\'efinir et de pr\'esenter
 de mani\`ere lisible et efficace les \algos
 consid\'er\'es.
Les autres langages de calcul formel \gnlz istes auraient
pu aussi bien faire l'affaire.
Il n'y aura d'ailleurs aucun mal \`a impl\'ementer dans un de
ces langages les \algos
pr\'esent\'es dans ce livre. Une liste r\'ecapitulative en est donn\'ee
dans la table page \pageref{listeAg}.

\bni\textbf{L'esprit dans lequel est \'ecrit cet ouvrage}

\ms Nous avons en \gnl donn\'e des preuves compl\`etes de nos r\'esultats, en
accordant une grande place aux exemples. Mais il nous est aussi arriv\'e de ne
donner qu'une id\'ee de la preuve, ou de ne la donner compl\`etement que sur un
exemple, ou de renvoyer \`a une r\'ef\'erence. Nous assumons tr\`es consciemment ce
que nous avons sacrifi\'e de la rigueur formelle au profit de la compr\'ehension de
\gui{ce qui se passe}. Nous avons essay\'e de donner  dessins et figures pour
illustrer notre texte, tout en ayant conscience d'en avoir fait bien trop peu.

\ss Nous avons aussi essay\'e de rapprocher cet expos\'e de la pratique concr\`ete
des \algosz, en d\'eveloppant chaque fois que nous l'avons pu des calculs de \com
dans lesquels nous explicitons les constantes \gui{cach\'ees dans le grand $\O$},
sans la connaissance desquelles les r\'esultats th\'eoriques 
n'ont pas de r\'eelle port\'ee pratique, et peuvent \^{e}tre trompeurs.

\ss Le niveau requis pour lire ce livre est seulement une bonne familiarit\'e avec
l'\agr \linz. Le mieux serait \'evidemment d'avoir lu auparavant cette perle rare
qu'est le livre de Gantmacher \cite{Gan}. On peut recommander aussi le grand
classique (toujours disponible) \cite{LM} de Lancaster \& Tismenetsky. Il est
naturellement pr\'ef\'erable, mais pas indispensable, d'avoir une id\'ee des
concepts de base de la \cobi pour lesquels nous recommandons les ouvrages \cite{Bal}
et \cite{Ste}. 

\ss Enfin, sur les \algos en \gnlz, si vous n'avez pas lu le livre de Knuth
\cite{Knu} parce que vous comprenez mal l'anglais ou que vous \^etes plut\^ot
habitu\'es \`a la langue de Voltaire, avant m\^eme de commencer la lecture de notre
ouvrage, \'ecrivez une lettre \`a tous les \'editeurs scientifiques en leur
demandant par quelle aberration la traduction en fran\c{c}ais n'a pas encore \'et\'e
faite.

\ss Pour aller au del\`a en Calcul Formel nous recommandons les 
livres de von zur Gathen \& Gerhard
\cite{VonZurbook}, Bini  \& Pan \cite{BP},
B\"{u}rgisser, Clausen \& Shokrollahi \cite{Bur},
B\"{u}rgisser  \cite{Bur2} et le Handbook of Computer Algebra~\cite{GKW}.

\ss Nous esp\'erons que notre
livre contribuera \`a mieux faire saisir l'impor\-tan\-ce
de la \com \agq \`a un moment o\`u les math\'ematiques constructives et les
solutions \algqs se d\'eveloppent de mani\`ere rapide et commencent \`a occuper de
plus en plus une place essentielle dans l'enseignement des Math\'ematiques, de
l'Informatique et des Sciences de l'ing\'enieur.

\bni\textbf{Remerciements}
Nous remercions Marie-Fran\c{c}oise Roy et Gilles Villard pour leur 
relecture attentive et leurs suggestions pertinentes, ainsi que Peter B\"urgisser
pour son aide concernant les deux derniers chapitres. 
Et enfin Fran\c{c}ois P\'etiard qui nous a fait b\'en\'eficier avec une 
patience infinie de son expertise en LaTeX.

\newpage \thispagestyle{empty}

\newpage
\thispagestyle{empty}\cleardoublepage 

\rdb
\addcontentsline{toc}{chapter}{Table des mati\`eres:  vous y \^etes!
\hspace{5.2cm}~}

\tableofcontents  

\newpage  \thispagestyle{empty}
\cleardoublepage 

\pagenumbering{arabic}



\chapter{Rappels d'alg\`ebre lin\'eaire}
\label{Chap RappelsAlgLin}
\label{Chap RapAlgLin}
\minitoc
\subsection*{Introduction }
Ce chapitre est consacr\'e \`a des rappels d'\agr \lin insistant
sur quelques \idas li\'ees aux \deters et \polcarsz. Notre but est
double. D'une part fixer les notations et donner les formules qui
justifieront les \algos de calcul de ces objets. D'autre part,
donner une id\'ee de l'\'etendue des applications qui pourront
\^{e}tre tir\'ees de ces calculs.

La section \ref{sec AlgLin} fixe les notations et rappelle la
formule de Binet-Cauchy ainsi que les \idts de Cramer et de
Sylvester.
La section \ref{subsec polycar} est consacr\'ee au
\polcar et \`a la formule de Samuelson. Dans la section
\ref{subsec polmin} nous \'etudions le \polmin et les sous-espaces
de Krylov. La section \ref{subsec suites rec} est consacr\'ee aux
\srlsz. Nous rappelons les \idts li\'ees aux sommes de Newton dans
la section \ref{subsec Newton}. La section \ref{secHadamard}
aborde les \mets du calcul modulaire. Enfin la section
\ref{secUnifSli} est consacr\'ee \`a l'\iMP et \`a ses
g\'en\'eralisations.

\section {Quelques propri\'et\'es g\'en\'erales}
\label{sec AlgLin}

\subsection{Notations}
\label{subsecNotachap1}

Dans cet ouvrage \,$\A$\, est un \acom et unitaire\footnote{~Si
\,$\A$\, n'est pas unitaire, on peut toujours le plonger dans un
anneau avec unit\'e (\cf \cite{Jac}).} et  \,$\K$\, un corps
commutatif. \indexnota{Aa@$\A$}  \indexnota{K@$\K$} Pour deux
entiers positifs quelconques \,$m$\, et \,$n$, \,$\A^{m\times
n}$\, \indexnota{Amn@$\A^{m\times n}$} d\'esigne l'ensemble des
matrices de \,$m$\, lignes et \,$n$\, colonnes \`a \coes dans
\,$\A$.

Soit \,$A=(a_{ij})\in \A^{n\times n}$\, une \maca
  d'ordre \,$n\geq 2$ \`a \coes dans \,$\A$
et soit un entier $r\;(1 \leq r \leq n)$.
On adopte les notations suivantes:
\begin{itemize}
\item $\Ir$ est la matrice unit\'e d'ordre $r$.
\item $\det(A)$\,  d\'esigne le
\emph{\deterz} de \,$A$: par \dfnz, c'est le \pol en les $a_{ij}$
d\'efini par la m\^eme formule que dans le cas d'une matrice \`a
coefficients dans un corps, autrement dit $\,\det{(A)}$\, est
donn\'e par:
$$ \det{(A)}=\sum_{\sigma} \,{{\varepsilon(\sigma)} \prod_{i=1}^n
{a_{i,\sigma(i)}}} $$
o\`u $\sigma$ parcourt l'ensemble des permutations de l'ensemble
d'indices $\{1,\alb\ldots,\alb n\}$ et o\`u $\varepsilon(\sigma)$
est la signature de la permutation $\sigma$. Lorsque cela ne
pr\^ete pas \`a confusion, nous noterons parfois \,$|A|$\, au lieu
de $\,\det{(A)}\,$ le \deter de \,$A$\,.
\item $\Tr{(A)}$\, est la \emph{trace} \index{matrice!trace d'une}
de \,$A$, \cad la somme de ses
\elts diagonaux.
\item La \emph{comatrice} de \,$A$\, \index{comatrice}
est la matrice $(d_{ij})_{1\leq i,j\leq n}$ o\`u chaque
\,$d_{ij}$\, est le cofacteur de l'\elt en position $(i,j)$ dans
\,$A$, \cad
$$ d_{ij}=(-1)^{i+j}\det{(B_{ij})} $$
o\`u \,$B_{ij}$\, est la matrice obtenue \`a
partir de \,$A$\, en supprimant la
\,$i\,$\eme ligne et la
\,$j\,$\eme colonne.
\item On a alors les formules de d\'eveloppement de \,$\det{(A)}$\,
(suivant la \,$i\,$\eme ligne ou suivant la \,$j\,$\eme colonne)
valables sur un \acomaz:
$$\det{(A)}=\sum_{k=1}^na_{ik}\,d_{ik}=\sum_{k=1}^na_{kj}\,d_{kj}
~~~(1\leq i,j\leq n)\,.$$
\item $\Adj(A)$\, d\'esigne la matrice adjointe
\index{matrice!adjointe} de \,$A$.
C'est la transpos\'ee de la comatrice de \,$A$.
Rappelons qu'elle v\'erifie
la double \egtz:
\begin{equation} \label{EqAdj}
A\,\Adj(A) =  \Adj(A)\,A=\det{(A)}\;\I_n\,.
\end{equation}
\end{itemize}

\sni Maintenant \,$A=(a_{ij})\in \A^{m\times n}$\, d\'esigne une
matrice quelconque  \`a \coes dans \,$\A$\, et \,$r$\, un entier
$\in\{1,\ldots, \min(m,n)\}$.
\begin{itemize}
\item Un \emph{mineur} d'ordre \,$r\,$\index{mineur!d'une matrice}
de  \,$A$\, est le \deter d'une \maca  extraite de
  \,$A$\,  en supprimant \,$m-r$\, lignes et \,$n-r$\,
  colonnes.
\item \refstepcounter{bidon}
\label{Minijhk} $A_{i..j,h..k}$\, d\'esigne la matrice
extraite de \,$A$\, sur les lignes de \,$i$\, \`a \,$j\,$ et sur
les colonnes de \,$h$\, \`a \,$k$\, ($i\leq j\leq m,\,h\leq k\leq
n$). \indexnota{Aijhk@$A_{i..j,h..k}$} Si \,$j-i=k-h$\, on note
\,$a^{h..k}_{i..j}=\det{(A_{i..j,h..k})}\,$ le mineur
correspondant. \indexnota{aijhk@$a^{h..k}_{i..j}$}
\item Plus g\'en\'eralement, si
\,$\alpha=\{\alpha_1,\ldots,\alpha_r\}$\, avec \,$1\leq\alpha_1<
\cdots< \alpha_r\leq m $\, et
\,$\beta=\{\beta_1,\ldots,\beta_s\}$\, avec \,$1\leq\beta_1<
\cdots< \beta_s\leq n $\, on notera \,$A_{\alpha,\beta}$\, la
matrice extraite de \,$A$\, sur les lignes (resp. les colonnes)
dont les indices sont en ordre croissant dans \,$\alpha $\, (resp.
dans \,$\beta$). \indexnota{Aalpha@$A_{\alpha,\beta}$}
\item Les \emph{sous-matrices principales}
\index{sous-matrice!principale} sont les sous-matrices dont la
diagonale principale est extraite de celle de \,$A$. On appellera
\emph{\mips de}  \,$A\,$\index{mineur!principal} les \deters des
sous-matrices principales de \,$A$. Comme cas particulier, $A_r$
d\'esigne la sous-matrice $\,A_{1..r,1..r}~$: nous dirons que
c'est une \emph{\smpdz} \index{sous-matrice!principale dominante}
de \,$A$.
  Son \deter est appel\'e un \emph{\mip dominant.}
\item $A_{i,h..k}=(a_{ih},\dots,a_{ik})=A_{i..i,h..k}$\, est la
matrice-ligne extraite de la  \,$i\,$\eme ligne de \,$A$\, sur les
colonnes de \,$h$\, \`a \,$k$. On pose de m\^{e}me
$A_{i..j,h}=\tra{(a_{ih},\dots,a_{jh})}=A_{i..j,h..h}.$
  \indexnota{Aihk@$A_{i,h..k}$}  \indexnota{Aihk@$A_{i..j,h}$}
\item \refstepcounter{bidon} \label{minrij}
On d\'efinit
\,$a_{\,ij}^{(r)}=\left|\begin{array}{cc} A_r &
\,A_{1..r,j} \\
A_{i,1..r} & a_{ij}\\
\end{array}\right|$\,
pour tous entiers \,$r,\,i,\,j$\, tels que
\,$1\leq r\leq \min(m,n)-1$\, et \,$r< i\leq m,\,r< j\leq n$,
le mineur d'ordre \,$r+1$\, de \,$A$\, obtenu
en bordant la \smpd
\,$A_r$\, par les \coes correspondants de
la \,$i\,$\eme ligne et de la
\,$j\,$\eme colonne de \,$A$,
ce qui fait par exemple que
\,$|A_r|=a_{\,rr}^{(r-1)}$.
Par convention on pose \,$a_{\,ij}^{(0)}= a_{ij}$\,
et \,$a_{\,00}^{(-1)}=1$.
%
\end{itemize}

\ms Parmi les propri\'et\'es du \deter qui restent valables dans
un \acom arbitraire, on doit citer en premier la lin\'earit\'e par
rapport \`a chaque ligne et par rapport \`a chaque colonne. La
deuxi\`eme propri\'et\'e la plus importante est son caract\`ere
\emph{altern\'e}, \cad qu'il s'annule si deux lignes ou si deux
colonnes sont \'egales. On d\'eduit de ces propri\'et\'es que le
\deter d'une \maca ne change pas si on ajoute \`a une ligne (resp.
\`a une colonne) une \coli des autres lignes (resp. des autres
colonnes).

Plus \gnltz, on peut citer toutes les propri\'et\'es qui
rel\`event d'\idasz. Par exemple l'\egt
\,$\det{(AB)}=\det{(A)}\det{(B)}$, ou encore le \tho de
Cayley-Hamilton (qui peut \^etre vu, pour une \maca d'ordre \,$n$,
comme une famille de \,$n^2$\, \idasz). Ces \idts sont
v\'erifi\'ees lorsque les \coes sont r\'eels, elles sont donc
vraies dans tous les \acoms (un \pol en \,$k$\, variables \`a
\coes entiers est identiquement nul \ssi il s'annule sur $\ZZ^k$).

\subsection{Formule de Binet-Cauchy}
\label{subsecBinCau}

C'est une formule qui \gnlz ise l'\egt
$\,\det{(AB)}=\det{(A)}\det{(B)}\,$ au cas d'un produit \,$AB$\,
avec \,$A\in\A^{m{\times}n}$\, et \,$B\in\A^{n{\times}m}$\, ($\,m
< n\,$): pour chaque \,$m$-uple \,$\beta =\beta_1,\ldots
,\beta_m$\, extrait en ordre croissant de \,$\{1,\ldots ,n\}$\,,
on consid\`ere \,$A_{1..m,\beta }$\, la matrice extraite de
\,$A$\, sur les colonnes \,$\beta_1,\ldots ,\beta_m$\, et
\,$B_{\beta,1..m}$\, la matrice
extraite\,\footnote{~$A_{1..m,\beta }$\, et \,$B_{\beta,1..m}$\,
sont des \macas d'ordre $\,m$.}  de \,$B$\, sur les lignes
\,$\beta_1,\ldots ,\beta_m$, alors on a la \emph{formule de
Binet-Cauchy} (\cf \cite{Gan} p. 9):
\begin{equation} \label{eqBC}
\det{(AB)}=\sum\nolimits_{\beta}\;\det{(A_{1..m,\beta})}
\,\det{(B_{\beta,1..m})} \index{formule!de Binet-Cauchy}
\end{equation}
(la somme comporte \,$\rC^m_n$\, termes).

\ss Pour le v\'erifier, on pose \,$A=(a_{ij})\in
\A^{m{\times}n}$\,, \,$B=(b_{ij})\in \A^{n{\times}m}$\, et on
utilise les propri\'et\'es \'el\'ementaires des  \deters pour
obtenir la suite d'\'egalit\'es imm\'ediates suivantes (avec des
notations \'evidentes): $$\det{(AB)} = \left|\begin{array}{ccc}
\sum_{i_1=1}^n\,a_{1,\,i_1}b_{i_1,1} & \ldots &
\sum_{i_m=1}^n\,a_{1,\,i_m}b_{i_m,m}
\\ \vdots & \vdots  & \vdots \\ \sum_{i_1=1}^n\,a_{m,\,i_1}b_{i_1,1}
& \ldots  & \sum_{i_m=1}^n\,a_{m,\,i_m}b_{i_m,m}
\end{array}\right|~~~~~~~~~~~~~~~~~~~~~~~~~~~~~~$$

$$~~~~~ = ~~\sum_{1\leq i_1\,,\,i_2\,,\, \ldots\,,i_m\leq n}~
\left|\begin{array}{ccc}a_{1,\,i_1}b_{i_1,1} & \ldots &
a_{1,\,i_m}b_{i_m,m}
\\ \vdots & \vdots  & \vdots \\ a_{m,\,i_1}b_{i_1,1}
& \ldots  & a_{m,\,i_m}b_{i_m,m}
\end{array}\right|~~~~~~~~~~~~~~~~~~~~~~~~~~~~~~~~~$$

$$~~~~~ = ~~\sum_{1\leq i_1\,,~i_2\,,\,\ldots\,,~i_m\leq n}~
b_{i_1,1}\times b_{i_2,2}\times \cdots \times b_{i_m,m}\times
\left|\begin{array}{ccc}a_{1,\,i_1} & \ldots & a_{1,\,i_m}\\
\vdots & \vdots  & \vdots \\ a_{m,\,i_1} & \ldots & a_{m,\,i_m}
\end{array}\right|.~~~~~~~~~~~$$

\sni Parmi les $\,n^m\,$ termes de cette somme, il n'y a que
$\,m!\,\rC^m_n\,$ termes qui risquent de ne pas \^etre
nuls\footnote{~\`a cause du fait qu'un d\'eterminant ayant deux
colonnes identiques est nul~; c'est d'ailleurs la raison simple
pour laquelle $\,\det{(AB)}=0\,$ lorsque $\,m>n$.}. Ce sont,
pour chacun des $\,\rC^m_n\,$ multi-indices
$\,\beta=(\,\beta_1,\,\beta_2,\,\ldots,\,\beta_m\,)\,$ avec
$\,1\leq \beta_1<\beta_2<\ldots<\beta_m\leq n\,,$ les $\,m!\,$
termes correspondant aux multi-indices
$\,(i_1,\,i_2,\,\ldots,\,i_m)\,$ tels que:
$$\{\,i_1,\,i_2,\,\ldots,\,i_m\, \}=
\{\,\beta_1,\,\beta_2,\,\ldots,\,\beta_m\, \}\,.$$ 
On regroupe ces termes en $\,\rC^m_n$\, sommes partielles. Chaque
somme partielle correspond \`a une valeur $\,\beta\,$ du
multi-indice $\,(\beta_1,\,\beta_2,\,\ldots,\,\beta_m)\,$ et
comporte $m!$ termes dans lesquels on peut mettre en facteur:
$$\left|\begin{array}{ccc}a_{1,\,\beta_1} & \ldots &
a_{1,\,\beta_m}\\ \vdots & \vdots  & \vdots \\ a_{m,\,\beta_1} &
\ldots & a_{m,\,\beta_m}\end{array}\right|~=~
\det{(A_{1..m,\beta})}\,.$$
Il suffit alors d'utiliser la d\'efinition du \deter de la matrice
$\,B_{\beta,1..m}\,$ pour voir que la somme partielle
correspondant au multi-indice $\,\beta\,$ n'est autre que:
$\,\left|\begin{array}{ccc} b_{\beta_1\,,\,1} & \ldots &
b_{\beta_1\,,\,m}\\ \vdots & \vdots & \vdots \\ b_{\beta_m\,,\,1}
& \ldots & b_{\beta_m\,,\,m}\end{array}\right|\,
\det{(A_{1..m,\beta})}$.  Ce~qui~donne:
$$\det{(AB)} =
\sum_{1\leq \beta_1< \beta_2<\ldots< \beta_m\leq
n}~\det{(A_{1..m,\beta})}\,\times\, \det{(B_{\beta,1..m})}\,.
~~~~~~~~~~~~~\Box$$

\subsection{Rang,  \deter et \idts de Cramer}
\label{subsecRaDe}

Une matrice carr\'ee est dite
\emph{r\'eguli\`ere}\index{matrice!r\'eguli\`ere} si
elle est inversible, \cad si son \deter
est inversible dans \,$\A$, et \emph{singuli\`ere}
\index{matrice!singuli\`ere} si son \deter est
nul. Une matrice (non \ncrt carr\'ee) est dite
\emph{\fregz}
\index{matrice!fortement r\'eguli\`ere} si toutes ses
\smpds sont \regs
(\ie tous les \mips dominants sont
inversibles dans l'anneau de base consid\'er\'e).

Lorsque \,$\A$\,
est suppos\'e \emph{int\`egre}, on d\'esignera par
\,$\mathrm{rg}(M)$\, le rang d'une
matrice\index{matrice!rang d'une} \,$M$\, quelconque
\`a \coes dans \,$\A$,
\cad l'ordre\indexnota{rgM@$\mathrm{rg}(M)$, rang d'une matrice}
maximum des mineurs non nuls de \,$M.$

\medskip Utilisant les notations ci-dessus,
nous rappelons maintenant quelques r\'esultats \elrs
d'\agr \lin dont certains seront
accom\-pa\-gn\'es de br\`eves d\'emons\-trations.

Comme nous travaillerons souvent avec un \acoma  \,$\A$,
nous aurons besoin de la notion de \emph{\Amoz},
\index{module!sur un anneau} \index{Amo@\Amoz}
qui est la \gnn aux anneaux de la notion
d'\evc sur un corps. Un \Amo \,$M$\, est par \dfn
un groupe ab\'elien (la loi de groupe est not\'ee $+$) muni d'une loi
externe
\,$\A\times M\rightarrow M,\;(a,x)\mapsto a.x$\, qui v\'erifie les
axiomes
usuels (pour tous \,$a,b\in\A$\, et \,$x,y\in M$):
$$\begin{array}{rcl}
1.x&=&x  \\
a.(b.x)&=&(ab).x  \\
(a+b).x&=&a.x+b.x  \\
a.(x+y)&=&a.x+a.y  \\
\end{array}$$
Nous ne consid\'ererons dans cet ouvrage que des
\emph{modules libres de dimension finie}, \index{module!libre}
\cad isomorphes \`a \,$\A^n$, ou parfois le module \,$\A^{\N}$.

Dans le cas o\`u il est int\`egre, l'anneau \,$\A$\, peut
\^{e}tre plong\'e dans son corps des
fractions qui est not\'e \,$\FA$, et tout \Amo
libre de dimension finie (isomorphe \`a un module \,$\A^n$)
peut \^etre consid\'er\'e comme inclus dans un
\FAev (isomorphe \`a $({\FA})^n$).
Le rang d'une matrice est alors \'egal \`a son rang usuel si on
consid\`ere que ses \coes
sont dans le corps \,$\FA$.

Dans la suite, chaque fois que l'hypoth\`ese d'int\'egrit\'e
sur \,$\A$\, doit intervenir,
elle sera  clairement soulign\'ee.

\begin{propri}\label{p1}
Soit \,$\A$\, un anneau int\`egre.
Pour toutes \macas \,$M$\, et \,$N$\,
d'ordre \,$n\geq 2$, \`a \coes dans
\,$\A$, on a:

$\mathrm{(i)}\, \,  MN=0\, \implique \,
\rg(M)+\rg(N)\leq n;$

$\mathrm{(ii)}\, \, \Adj(M)=0
\,\equiva\, \rg(M)\leq n-2.$
\end{propri}
\preuve (i) provient du fait que si \,$u$\, et
\,$v$\, sont deux \endos d'un \FAev de dimension finie \,$n$, alors:
$\mathrm{dim}(\Im\, u) + \mathrm{dim}(\Ker\, u)
= n$ ~et~ ($u\circ v=0\impliq \Im\, v
\subseteq \Ker\, u$). \\
(ii) d\'ecoule du fait que si $\Adj(M)=0$,
tous les mineurs d'ordre \,$n-1$\, de \,$M$\, sont
nuls, et r\'eciproquement. \qed

\begin{propri}\label{p2} ~
$\A$\, \'etant un anneau int\`egre, le rang de la matrice adjointe
de toute matrice carr\'ee \emph{singuli\`ere} \,$M\in \A^{n\times
n}~(n\geq 2)$\, est au plus \'egal \`a $1$, et on a les
\'equivalences: $$\Adj(M)\neq 0 \; \equiva\;  \rg(M)=n-1 \equiva
\rg(\Adj(M))=1.$$
\end{propri}
\preuve C'est une cons\'equence de la propri\'et\'e \ref{p1}
sachant que, par hypoth\`ese, \,$M\,\Adj(M)=\det{(M)}\,\I_n=0$.
\qed

\begin{propri}\label{p3}
~Soit \,$\A$\, un \acom arbitraire. Pour tous entiers \,$n,\,p
\geq 2$\, et toutes matrices \,$P\in\A^{p\times n},\
M\in\A^{n\times n} ~et~Q\in \A^{n\times p}$, on a l'implication:
$$\det{(M)}=0 \implique \det{(P\,\Adj(M)\,Q)}=0$$
\end{propri}
\preuve~\\
$\bullet~$ Supposons tout d'abord l'anneau
\,$\A$\, int\`egre. \\
La propri\'et\'e (\ref{p2})
nous permet alors d'\'ecrire: \\
$\det{(M)}=0 \impliq \rg(\Adj(M)) \leq 1 \impliq \rg(P\,\Adj(M)\,
Q) \leq 1 < p$ \\
et de conclure.

\sni $\bullet~$ Voyons maintenant le cas o\`u \,$\A$\, n'est pas
int\`egre, et commen\c{c}ons par le cas g\'en\'erique.
Consid\'erons pour cela l'anneau
\,$\ZZ[a_1,a_2,\ldots,a_t]=\ZZ[\s{a}]$\, o\`u \,$t=n^2+2np$\, et
o\`u  les $a_i$ sont des variables repr\'esentant les entr\'ees
des matrices \,$M,\,P$\, et \,$Q$. Il est bien connu que
\,$\det{(M)}\in \ZZ[\s{a}]$\, est un \poly  irr\'eductible dans
cet anneau et que donc \,$\aqo{\ZZ[\s{a}]}{\det{(M)}}$\, est un
anneau int\`egre.
Comme dans cet anneau \,$\det{(M)}=0\,$, on est ramen\'e au cas
pr\'ec\'edent et l'on a: $\det{(P\,\Adj(M)\,Q)}=0$\, dans
\,$\aqo{\ZZ[\s{a}]}{\det{(M)}}$.
Ceci prouve l'existence d'un \poly
\,$f(\s{a})\in\ZZ[\s{a}]$\, v\'e\-ri\-fiant l'\idtz:
$$\det{(P\,\Adj(M)\,Q)} =  f(\s{a})\;\det{(M)}$$
(le \poly \,$f(\s{a})$\,  peut \^etre calcul\'e
par un \algo de division exacte dans l'anneau
\,$\ZZ[\s{a}]$).
Cette \ida est v\'erifi\'ee dans tout
\acomz, ce qui permet de conclure dans
le cas \gnlz.
\qed

\begin{proposition}\label{propICRA}
\emph{(Identit\'es de Cramer)}\index{Cramer}
Soit \,$\A$\, un \acom arbitraire et \,$A\in\A^{m{\times}n}$.
Notons \,$C_j$\, la \,$j$\eme colonne de \,$A$\, ($C_j=A_{1..m,j}$)
et \,$B_j$\, la matrice extraite de \,$A$\, en supprimant la colonne
\,$C_j$.
\begin{enumerate}
\item Si \,$n=m+1$\, et si \,$\mu_j=\det{(B_j)}$\, on a:
\begin{equation} \label{eqIC1}
\sum\nolimits_{j=1}^{n}\,(-1)^j\,\mu_j\,C_j\,=\,0
\index{Cramer!identif\'es de}
\end{equation}
\item Supposons  \,$n\leq m$\, et que tous les mineurs d'ordre \,$n$\,
de \,$A$\, sont nuls. Soit
\,$\alpha=\alpha_1,\ldots,\alpha_{n-1}$\, extrait en ordre
croissant de \,$\{1,\ldots ,\alb m\}$. Soit
\,$D_j=(B_j)_{\alpha,1..n-1}\in\A^{(n-1){\times}(n-1)}$\, la
matrice extraite de \,$B_j$\, en gardant les lignes de \,$\alpha$,
et soit \,$\nu_j=\det{(D_j)}$. Alors on a:
\begin{equation} \label{eqIC2}
  \sum\nolimits_{j=1}^{n}\,(-1)^j\,\nu_j\,C_j\,=\,0
\end{equation}
\end{enumerate}
\end{proposition}
\preuve~\\
Pour le premier point: la coordonn\'ee \num$k$\, de
\,$ \sum_{j=1}^{n}\,(-1)^j\,\mu_j\,C_j$\, est \'egale au \deter de la
matrice obtenue en collant au dessous de \,$A$\, la ligne \num$k$\, de
\,$A$\,
(ceci se voit en d\'eveloppant ce \deter selon la derni\`ere ligne).
Cette coordonn\'ee est donc nulle.

\noindent Le deuxi\`eme point se prouve de mani\`ere analogue.  \qed

\medskip Ces deux \'egalit\'es peuvent \^{e}tre relues sous la forme
\gui{solution d'un \sliz}. Pour la premi\`ere on consid\`ere
\,$A\in\A^{n{\times}n}$\, et \,$V\in\A^{n{\times}1}$\, on note
\,$E_j$\, la matrice obtenue \`a partir de \,$A$\, en
rempla\c{c}ant la colonne \,$C_j$\, par \,$V$\,, et
\,$\mu'_j=\det{(E_j)}$. Alors (\ref{eqIC1}) se relit sous la
forme plus classique:
\begin{equation} \label{eqIC3}
\det{(A)}\,V\;=\; \sum\nolimits_{j=1}^{n}\,\mu'_j\,C_j \;=\;
A\,\cmatrix{\mu'_1\cr\vdots\cr\mu'_n}\;=\;A\cdot\Adj(A)\cdot V
\end{equation}
On peut \egmt relire (\ref{eqIC2}) comme suit. On consid\`ere une
matrice   \,$A\in\A^{m{\times}n}$\, et \,$V\in\A^{m{\times}1}$\,
avec \,$m>n$. \emph{On suppose que tous les mineurs d'ordre
\,$n+1$\, de \,$(A|V)$\, sont nuls.} On note toujours
\,$C_j=A_{1..m,j}$\, la \,$j\,$\eme colonne de \,$A$. On choisit
\,$\alpha=\alpha_1,\ldots,\alpha_{n}$\, extrait en ordre croissant
de \,$\{1,\ldots , m\}$. On consid\`ere la matrice \,$F_{j}$\,
obtenue \`a partir de \,$A_{\alpha,1..n}$\, en rempla\c{c}ant la
\,$j\,$\eme colonne par \,$V$\, et on pose
\,$\nu_{\alpha,j}=\det{(F_{j})}$. Si on applique (\ref{eqIC2})
avec la matrice \,$\left(A|V\right)\in \A^{m{\times}(n+1)}$\, on
obtient:
\begin{equation} \label{eqIC4}
\det{(A_{\alpha,1..n})}\,V\;=\;
\sum\nolimits_{j=1}^{n}\,\nu_{\alpha,j}\,C_{j}
\end{equation}

\begin{proposition}\label{p5}
~Soit \,$\A$\, un \acom arbitraire non tri\-vial, \cad dans lequel \,$1_\A\neq
0_\A$, et \,$A\in\A^{m{\times}n}$. \Propeq
\begin{enumerate}
\item Pour tout \,$V\in \A^{m{\times}1}$\, il existe \,$X\in
\A^{n{\times}1}$\, tel que \,$A\,X=V$. Autrement dit, l'\ali
\,$\varphi:\A^n\rightarrow \A^m$\, d\'efinie par \,$A$\, est surjective.
\item Il existe \,$B\in\A^{n{\times}m}$\, tel que \,$A\,B=\I_m$.
\item On a \,$n\geq m$\, et il existe une \coli des mineurs d'ordre
\,$m$\, de \,$A$\, qui est \'egale \`a \,$1$.
\end{enumerate}
\end{proposition}
\preuve~\\
(1)\,$\Rightarrow $\,(2)\, Soit \,$e_j$\, le \,$j\,$\eme vecteur
de la base canonique de \,$\A^m$\, et soit \,$X_j$\, un vecteur de
  \,$\A^n$\, tel que \,$A\,X_j=e_j$. On prend pour matrice \,$B$\, la
matrice dont les colonnes sont les \,$X_j$.

\noindent (2)\,$\Rightarrow $\,(1)\, On prend \,$X=B\,V$.

\noindent (2)\,$\Rightarrow $\,(3)\, Montrons que \,$n<m$\, est
impossible. Si tel est le cas on rajoute \,$m-n$\, colonnes nulles
\`a droite de \,$A$\, et \,$m-n$\, lignes nulles en dessous de
\,$B$, on obtient deux \macas \,$A'$\, et \,$B'$\, pour lesquelles
on a \,$\det{(A')}=\det{(B')}=0\,$ et $\,A'\cdot B'=A\cdot B =
\I_m\,$, ce qui donne \,$0_\A=1_\A$\,.

\noindent Pour la \coliz, on applique la formule de Binet-Cauchy
(\ref{eqBC}) avec \,$A\,B=\I_m$.

\noindent (3)\,$\Rightarrow $\,(2)\, Supposons \,$\sum_{\beta}
c_\beta \det{(A_{1..m,\beta})}=1$. La somme est \'etendue \`a tous
les \,$\beta=\{\beta_1,\ldots ,\beta_m\}$\, o\`u \,$1\leq\beta_1<
\cdots < \beta_m\leq n $. On a
\,$A_{1..m,\beta}=A\cdot(\I_n)_{1..n,\beta}$. Posons
\,$B_\beta=(\I_n)_{1..n,\beta}\cdot\Adj(A_{1..m,\beta})$. Alors
\,$A\cdot B_\beta=\det{(A_{1..m,\beta})}\,\I_m$. Il suffit donc de
prendre \,$B=\sum_\beta \,c_\beta\,B_\beta $.
  \qed

\subsection{Identit\'es de Sylvester}
\label{subsecIdSylv}

\begin{propri}\label{ber}
\emph{(Une \idt de Sylvester)} \\ Soit \,$\A$\, un \acom
quelconque. Pour tout entier \,$n\geq 2$\, et toute matrice
\,$A=(a_{ij})\in\A^{n\times n}$, on a: $$\
\det{(A)}=a_{nn}\,\det{(A_{n-1})} - A_{n,1..n-1}\,[\Adj(A_{n-1})]
\,A_{1..n-1,n}\,.$$
\end{propri}
\preuve
Pour obtenir cette formule, il suffit
de d\'evelopper
$$\det{(A)}=\left|\begin{array}{lr} A_{n-1} & A_{1..n-1,n} \\
A_{n,1..n-1} & a_{nn}\\
\end{array}\right|$$
suivant la derni\`ere ligne puis chacun des
cofacteurs des \elts de \,$A_{n,1..n-1}$\,
intervenant dans ce \dev
suivant la derni\`ere colonne qui n'est autre
que \,$A_{1..n-1,n}$. \qed

\ss
Nous allons voir maintenant que la propri\'et\'e
pr\'ec\'edente peut \^etre
\gnlz is\'ee \`a d'autres partitions de
\,$A$. \'Etant donn\'es en effet deux entiers
\,$r$\, et \,$n$\, avec \,$n\geq 2$\, et
\,$1\leq r<n$, on associe \`a toute matrice
\,$A\in \A^{n\times n}$\, la partition suivante
de la matrice \,$A$\, en blocs:
$$
A= \cmatrix{
A_{r} & A_{12} \cr
A_{21} & A_{22} }
$$
o\`u $\,A_{12}=A_{1..r,r+1..n}\in\A^{r\times (n-r)}\,,$
$\,A_{21}=A_{r+1..n,1..r}\in \A^{(n-r)\times r}\,$ et
$\,A_{22}=A_{r+1..n,r+1..n}\in\A^{(n-r) \times (n-r)}$.

\ss On a alors le r\'esultat suivant, valable pour
tout \acom et unitaire \,$\A$.

%
\begin{prop}[Identit\'es de Sylvester]\label{syl}
\index{Sylvester!identit\'es de}
Avec les notations ci-dessus, et pour tous entiers
\,$n$\, et \,$r$\, tels que \,$1\leq r\leq n-1$,
on a les \idts suivantes, dans lesquelles on a pos\'e
$B_r=\Adj A_{r}$:
\begin{equation}\label{eq1}
|A_{r}|\, A = \cmatrix{ A_{r} & 0 \cr
A_{21} &\I_{n-r} }
\cmatrix{ |A_{r}|\,\I_r &
B_r\, A_{12}\cr
0 &|A_{r}|\,A_{22}-A_{21}\,B_r\,A_{12}
}
\end{equation}
\begin{equation}\label{eq2}
|A_{r}|^{n-r-1}|A|= \det{(|A_{r}|\,A_{22}-A_{21}\,B_r\,A_{12})}
\end{equation}
\begin{equation}\label{eq3}
  \left(a_{\,rr}^{(r-1)}\right)^{n-r-1}\,a_{\,nn}^{(n-1)}
=\left|\begin{array}{ccc}
a_{r+1,r+1}^{\,(r)} & \ldots  & a_{r+1,n}^{\,(r)} \\
\vdots &   & \vdots \\
a_{n,r+1}^{\,(r)} & \ldots  & a_{n,n}^{(r)}
\end{array}\right|
\end{equation}
\begin{equation} \label{eq3bis}
  a_{n-2,n-2}^{\,(n-3)}\,a_{\,nn}^{(n-1)}
=\left|\begin{array}{cc}
a_{n-1,n-1}^{\,(n-2)} & a_{n-1,n}^{(n-2)} \\
a_{n,n-1}^{(n-2)} & a_{\,n,n}^{(n-2)}
\end{array}\right|
\end{equation}
\end{prop}
\preuve

\noindent -- L'\egt matricielle (\ref{eq1})
r\'esulte de l'\idt
$\;A_{r}\;B_{r}=|A_{r}|\;\I_r$.

\noindent --
Pour d\'emontrer (\ref{eq2}) qui est une \idt
\agqz, on peut se restreindre au cas o\`u
les \coes \,$a_{ij}$\, de \,$A$\, sont
des \idtrsz. Les mineurs de \,$A$\,
peuvent alors \^{e}tre vus comme des \pols non
nuls dans l'anneau int\`egre
\,$\ZZ[(a_{ij})_{1\leq i,j \leq n}]$.
L'\egt (\ref{eq2}) est alors obtenue en
prenant les \deters des deux membres de
l'\egt  (\ref{eq1}) et en simplifiant par
le \poly $|A_{r}|^{r+1}$.

\noindent -- Le premier membre de l'\idt (\ref{eq3}) est le m\^eme
que celui de (\ref{eq2}). L'\egt des seconds membres provient du
fait que l'\elt de la \,$(i-r)\,$\eme ligne et \,$(j-r)\,$\eme
colonne ($r+1\leq i,j\leq n$) de la matrice:
$$\det{(A_{r})}\,A_{22}-A_{21}\,
[\Adj(A_{r})]\,A_{12}\in\A^{(n-r)\times (n-r)} $$
est \'egal \`a:
$$~a_{ij}\,\det{(A_{r})}- A_{i,1..r}\,[\Adj(A_{r})]\,
A_{1..r,j}\,$$
qui n'est autre que $~a_{\,ij}^{(r)}~$ d'apr\`es la
propri\'et\'e (\ref{ber}) appliqu\'ee \`a la matrice
\,$\left[\begin{array}{lr} A_r & A_{1..r,j} \\
A_{i,1..r} & a_{ij}\end{array}\right]$.

\noindent -- L'\egt (\ref{eq3bis}) est un cas
particulier de (\ref{eq3}) pour \,$r=n-2$. \qed

\mni
\rem \refstepcounter{bidon}
\label{r1} Si l'on fait \,$r=n-1$\, et par cons\'equent
$$\,A=\left[\begin{array}{lr} A_{n-1} & A_{1..n-1,n} \\
A_{n,1..n-1} & a_{nn}\end{array}\right],$$ l'\egt (\ref{eq2})
donne exactement la formule de la propri\'et\'e (\ref{ber}), ce
qui permet d'affirmer que celle-ci est une \idt de Sylvester
particuli\`ere.

\ms Les \idts de Sylvester seront utilis\'ees dans la section
\ref{sec JorBar} pour le calcul des \deters  par la \mjbz.

\section {Polyn\^ome caract\'eristique}
\label{subsec polycar}

\sni On appelle \emph{matrice \caraz}
\index{matrice!caract\'eristique} d'une matrice \,$A\in
\A^{n\times n}$\, la matrice $A-X\I_n\in (\A[X])^{n\times n}$\,
(\,$X$\, d\'esigne une \idtr sur \,$\A$\,).
\\ Le \emph{\polcarz} \index{polyn\^ome!caract\'eristique} de
\,$A$\, est, par d\'efinition, le \deter de sa matrice
\caraz~({\footnote{~Il serait en fait plus pratique de d\'efinir
comme le fait Bourbaki le \polcar de \,$A$\, comme le \deter de
\,$X\I_n-A$, mais nous nous en tenons \`a l'usage le plus
r\'epandu.}}). On le notera \,$\rPA$:
$$ \rPA(X)= \det{(A-X\I_n)}=p_0X^n + p_1X^{n-1} + \cdots +
p_{n-1}X + p_n\,. $$
Notons que \,$p_0=(-1)^n\,,~p_n=\det{(A)}$\, et que pour \,$1\leq
k\leq n-1$, le \coe \,$p_k$\, est le produit par \,$(-1)^{n-k}$\,
de la somme de tous les mineurs diagonaux d'ordre \,$k$\, de
\,$A$\,(\footnote{~Pour s'en convaincre, on peut examiner la
formule qui donne par \dfn \,$\rPA(-X)= \det{(A+X\I_n)}$\, et voir
quels sont les produits qui contiennent  \,$X^{n-k}$. On peut
aussi faire une preuve par \recu sur \,$n$\, en d\'eveloppant
\,$\det{(A+X\I_n)}$\, suivant la premi\`ere colonne.}). En
particulier:
$$(-1)^{n-1}\,p_1\,=\,\sum_{i=1}^n \,a_{ii}=\Tr(A).$$

\subsection{Matrice caract\'eristique adjointe}
\label{subsecMacad} 

On appelle \emph{matrice \cara adjointe} de \,$A\in
\A^{n\times n}$\, \index{matrice!caract\'eristique adjointe} la
matrice:
$$Q(X)=\Adj(X\I_n-A)=(-1)^{(n-1)}\,
        \Adj{(A-X\I_n)}\,.$$
C'est, \`a un signe pr\`es, l'adjointe de la matrice \cara de \,$A$\,.
Elle peut \^etre vue comme un \pol matriciel de degr\'e \,$n-1$\,
en $\,X$, \`a \coes dans \,$\A^{n\times n}$: en effet, la
matrice des cofacteurs de \,$X\I_n-A$\, est une matrice \,$n\times
n$\, dont les \elts diagonaux sont des \pols de degr\'e \,$n-1$\,
en \,$X$\, et les autres des \pols de degr\'e \,$n-2$\, en
\,$X$\,. \\
Ce qui fait que la \emph{matrice \cara adjointe} de
\,$A$\, s'\'ecrit:
\begin{equation}
\label{eqMacad0} Q(X)=
B_0X^{n-1}+B_1X^{n-2}+\cdots+B_{n-2}X+B_{n-1} \in \A^{n\times
n}[X]\,.~~
\end{equation}
D'apr\`es l'\'equation (\ref{EqAdj}) page \pageref{EqAdj}, on a
l'\egtz:
$$ (X\In - A)\,Q(X)=P(X)\,\In ~~\mathrm{o\grave{u}}~~P(X)=
(-1)^n\,\rPA{(X)}. $$
$\mathrm{On~posera~:}~~~P(X)= X^n-[\,c_1X^{n-1} + \cdots +
c_{n-1}X + c_n\,]\,.$

\sni Ainsi le \poly \,$P(X)\,\In$\, est divisible par le \poly
\,$(X\In - A)$\, au sens de la division euclidienne dans l'anneau
de \pols \,$\A^{n\times n}[X].$ Pour obtenir les \coes
(matriciels) du quotient \,$Q(X)$\, dans cette division, on
applique la \pcd de 
Horner{\footnote{~La \pcd (ou encore sch\'ema)
de Horner
\label{schemaHorner} 
n'est rien d'autre qu'une mise en
forme \algq de la division euclidienne d'un \pol
\,$P(X)=\sum_{i}^nc_iX^i$\, par un \pol \,$X-a$. Le reste, \'egal
\`a \,$P(a)$, est alors obtenu sous la forme
\,$c_0+a\,(c_1+a\,(c_2+\cdots +a\,c_n)\cdots )$. Ceci consitue une
\eva efficace de \,$P(a)$, utilisant un minimum de \mulsz. Cette
\met est en fait identique \`a celle de Ch'in Chiu-Shao employ\'ee
en Chine m\'edi\'evale. Elle a \'et\'e red\'ecouverte par Ruffini
(1802) et Horner (1819). Voir l'Encyclopedia of Mathematics, chez
Kluwer (1996).}} au \poly matriciel \,$P(X)\,\In\in\A^{n\times
n}[X]$\, avec la constante \,$A\in\A^{n\times
n}\,$\index{Horner!sch\'ema de}. Le reste $\,B_{n}=P(A)\,$ est
nul~; ce qui donne l'identit\'e \,$\rPA{(A)}=0\,$ et fournit en
passant une d\'emonstration (\'el\'egante) du \tho de
Cayley-Hamilton. \index{Cayley-Hamilton!th\'eor\`eme de}

\ms Le proc\'ed\'e de Horner peut \^etre repr\'esent\'e par le
sch\'ema suivant:

\ms
\begin{tabular}{|c|c|c|c|c|c|c|c|} \hline $P(X)\,\In$ & $\In$
& $-c_1\In$ & $-c_2\In$ & \dots & $-c_{n-1}\In$ & $-c_n\In$ & \\
\hline $A$ & 0 & $B_0$ & $B_1$ & \dots & $B_{n-2}$ & $B_{n-1}$ &
$B_{n}=0$ \\ \hline
\end{tabular}

\begin{equation}
\label{eqMacad1} {\mathrm {avec}}~:~~ B_0=\In~~~{\mathrm {et}}~~~
B_k=AB_{k-1}-c_k\In~~(1\leq k\leq n)\,.
\end{equation}

Cela donne une m\'ethode rapide et efficace pour calculer, \`a
partir des \coes du \polcarv la matrice \cara adjointe, et fournit
les relations d\'etaill\'ees suivantes utilis\'ees pour \'etablir
la formule de Samuelson (\S~\vref{subsecSamuelson}):
\begin{equation}
\label{eqMacad2} \left\{
\begin{array}{lcl} B_1 & = & A-c_1\In \\ B_2 &
= & A^2-c_1A-c_2\In \\ \;\vdots & \vdots & \quad \quad \vdots \\
B_k & = & A^k-c_1A^{k-1}-\dots-c_{k-1}A-c_k\In \\ \;\vdots &
\vdots &\quad  \quad \vdots \\ B_n & = &
A^n-c_1A^{n-1}-\dots-c_{n-1}A-c_n\In = 0
\end{array}\right.
\end{equation}

Notons qu'\`a la fin de la \pcd de Horner, on obtient \,$B_n=0$\,
\cad \,$A\,B_{n-1}-c_n\,\In=0$\, ou encore:
$$ A\,B_{n-1}=c_n\,\In=(-1)^{n-1}\,\det{(A)}\,\In\,.$$
Si \,$\det{(A)}$\, est inversible dans \,$\A$, alors \,$A$\,
poss\`ede un inverse qui peut \^{e}tre calcul\'e par la formule
\,$A^{-1}=(c_n)^{-1}\,B_{n-1}$. \\ Notons que
$\,B_{n-1}=Q(0)=(-1)^{n-1}\,\Adj(A)$. Donc
$$
\Adj(A)=(-1)^{(n-1)}\,[\,A^{n-1}-c_1\,A^{n-2}-\cdots-c_{n-2}\,A
-c_{n-1}\,\In\,]\,.
$$

Ainsi le calcul de la matrice \cara adjointe nous donne l'adjointe
de $\,A$. Il nous permet aussi d'obtenir l'inverse de $\,A\,$
(s'il existe) notamment dans les cas o\`u la \mpg s'av\'ererait
impraticable, ce qui se produit lorsque l'anneau contient des
diviseurs de z\'ero. La matrice \cara adjointe de $\,A\,$ sert
\'egalement, comme nous le verrons plus loin (voir \S \vref{sec
faddeev} et suivantes), \`a calculer le \polcar de $\,A\,$ par la
\met de Faddeev-Souriau-Frame et, dans certains cas, des vecteurs
propres non nuls.

Nous allons \`a pr\'esent l'utiliser pour \'etablir un r\'esultat
important pour la suite et faisant l'objet du paragraphe suivant.

\subsection{Formule de Samuelson}
\label{subsecSamuelson} \index{Samuelson}

\sni Utilisant l'\egt \ref{eqMacad0} et les relations
\ref{eqMacad2} dans lesquelles on remplace les \coes $\,c_i\,$ par
les \coes $\,p_i\,$ du \polcar de \,$A$\,, sachant que
$\,p_i=(-1)^{n+1}\,c_i$, on obtient l'\idaz:
\begin{equation} \label{eqAdj2}
\Adj(A-X\,\I_n)\,=\,-\;\sum_{k=0}^{n-1}\;
\left(\sum\nolimits_{j=0}^{k}\; p_j\,A^{k-j} \right)\,X^{n-1-k}\,.
\end{equation}
Cette \egt nous sert \`a d\'emontrer la \emph{formule de
Samuelson} \cite{Sam} (voir \cite{Gas}, \met de partitionnement
pp.~291--298, \cite{Fad})\,. \index{formule!de Samuelson}
\begin{prop}[Formule de Samuelson]\label{Samu}~\\
Soit \,$\A$\, un \acomaz, \,$n$\, un entier $\geq 2$\,,
\,$A=(a_{ij})\in\A^{n\times n}$\, et \,$r=n-1$\,. Notons
\,$P_{r}(X)=\sum_{i=0}^{r} q_{r- i}X^i=\det{(A_{r}-X\I_r)}$\, le
\polcar de la \smpd \,$A_{r}$. Posons \,$R_r:=A_{n,1..r}\,$ et
\,$S_r:=A_{1..r,n}$\, de sorte que la matrice \,$A$\, est
partitionn\'ee comme suit:
$$ A:=\cmatrix{
A_r
&S_r\cr R_r&a_{nn}}\,.
$$
Alors on a:
$$ \rPA(X)\,=\,(a_{n,n} - X)\, P_{r}(X)\, +\, \sum_{k=0}^{r-1}
\;\left(\sum\nolimits_{j=0}^{k} q_j\,(R_r\,A_r^{k-j} S_r)  \right)
X^{r-1-k}\,. $$
C'est-\`a-dire encore:
\begin{equation}\label{EqSamu}
\rPA(X)= \left\{
\begin{array}{lll}
(a_{n,n} - X)\, P_{n-1}(X)\, +
\\[1mm]
\sum_{k=0}^{n-2} \,\left[q_0\,(R_r\,A_r^{k}\,S_r) \,+ \cdots
+\,q_{k}\, (R_r\, S_r) \right] X^{n-2-k}
\end{array}
\right.
\end{equation}
\end{prop}
\prv Tout d'abord on applique l'\idt de Sylvester donn\'ee en
\ref{ber} \`a la matrice \,$(A-X\I_{n})$. On obtient: $$
\rPA(X)=(a_{n,n} - X)\,\det{(A_r-X \I_r)}-
R_r\,\Adj(A_r-X\I_r)\,S_r\,. $$ Ensuite on applique (\ref{eqAdj2})
en rempla\c{c}ant \,$A$\, par \,$A_r$\, et \,$n$\, par \,$r$:
$$\Adj(A_r-X \I_r)\,= \,-\;\sum_{k=0}^{r-1}\,
\left(\sum\nolimits_{j=0}^{k} q_j\,A_r^{k-j} \right)\,
X^{r-1-k}\,. ~~~~~~~~~ \Box$$

\sni La formule de Samuelson sera utilis\'ee dans l'\agbz.

\subsection{Valeurs propres de \texorpdfstring{\,$f(A)$}{f(A)}}

Supposons l'anneau \,$\A$\, int\`egre et soit \,$\FA$\, son
corps des fractions et \,$\L$\, une extension de
\,$\FA$\, dans laquelle \,$\rPA$\, se d\'ecompose
en produit de facteurs du premier degr\'e.
Une telle extension est ce qu'on appelle un
\emph{corps de \decoz} de \,$\rPA$.
Tout z\'ero de \,$\rPA$\, dans \,$\L$\, est
appel\'e \emph{valeur propre} de \,$A$,\index{valeur
propre} et sa \mut est, par \dfnz, la
\emph{\mut \agqz}
\index{multiplicit\'e!alg\'ebrique}
\index{valeur propre!multiplicit\'e alg\'ebrique d'une}
de la valeur propre de \,$A$.
Plus \gnltz, m\^eme si \,$\A$\, n'est pas
int\`egre, il est parfois utile de l'envoyer dans un corps
\,$\K$\, par un \homo d'anneaux unitaires
\,$\varphi: \A\rightarrow \K$.
Les valeurs propres de la matrice \,$A$\, dans
\,$\K$\, seront alors par \dfn les z\'eros
du \poly
$$\varphi(\rPA) \,\eqdef\,
\sum\nolimits_{i=1}^n\varphi(p_i)X^{n-i}\in\K[X]\,.
$$

Soit un \pol
\,$f=a_0X^m+a_1X^{m-1}+\cdots+a_{m-1}X+a_m\in \L[X].$
Consid\'erons la matrice
$~f(A)=a_0A^m+a_1A^{m-1}+\cdots+a_{m-1}A+a_m\I_n\in
\L^{n\times n}~$ et  les \polcars
\,$\rPA$\, et \,$\rP_{\!f(A)}$\, de \,$A$\,
et de \,$f(A)$.
Le lemme suivant exprime alors en particulier le lien entre les
valeurs propres de la matrice \,$A\in\A^{n\times n}$\,
et celles de la matrice \,$f(A)$, dans le
cas o\`u l'anneau de base \,$\A$\, est int\`egre.

\begin{lem}\label{trac}
Soit \,$\A$\, un anneau
int\`egre, \,$\L$\, une extension du corps des
fractions de \,$\A$, et \,$f$\, un \pol de
\,$\L[X]$. Si le \polcar de \,$A$\, s'\'ecrit:
$$\rPA(X)=(-1)^n\,(X-\l_1)\,(X-\l_2) \cdots (X-\l_n)$$
avec les \,$\l_i\in \L$, alors le \polcar de
\,$f(A)$\, s'\'ecrit:
$$\rP_{\!f(A)}=(-1)^n\,(X-f(\l_1))\,(X-f(\l_2)) \cdots
(X-f(\l_n)).~$$
En particulier, ${\Tr}(f(A))=\sum_{i=1}^nf(\l_i)$\, et pour tout
\,$k\in \NN$,   ${\Tr}(A^k)=\sum_{i=1}^n\l_i^k$.
\end{lem}
\preuve
Il suffit de montrer le premier point.
Dans le corps \,$\L$, qui contient toutes les valeurs
propres de la matrice \,$A$, celle-ci peut \^etre
ramen\'ee \`a une forme \trgz. C'est-\`a-dire
qu'il existe une \matg
\,$A'\in \L^{n\times n}$\, avec
les \,$\l_i$\, sur la diagonale et une matrice
\,$M\in \L^{n\times n}$\, inversible telles que
\,$A'=M^{-1}AM$. Comme \,$f$\, est un \polz,
et que \,$A'$\, est \trg de la forme
$$A'= \cmatrix{
\l_{1} & 0 &  \cdots & 0 \cr
\times & \ddots & \ddots& \vdots  \cr
\vdots & \ddots& \ddots & 0 \cr
\times & \cdots &\times & \l_{n}
}$$
la matrice \,$f(A')$\, sera aussi \trgz, de la forme
$$f(A')= \cmatrix{
f(\l_{1}) & 0 &  \cdots & 0 \cr
\times & \ddots & \ddots& \vdots  \cr
\vdots & \ddots& \ddots & 0 \cr
\times & \cdots &\times & f(\l_{n}) \cr
}\,.$$ \\
De plus, puisque \,$\Theta_M:A\vers M^{-1}A\:M$\, est un
automorphisme{\footnote{~Pr\'ecis\'ement,
\,$\B=\L^{n\times n}$\, est une \,$\L$\,-\,\agrz, \cad un anneau
  muni d'une loi externe
\,$(x,A)\mapsto x.A$\,
(produit d'une matrice par un scalaire) qui v\'erifie les \idts
\,$(x+y).A=x.A+y.A$, \,$x.(A+B)=x.A+x.B$, \,$x.(y.A)=xy.A$,
\,$x.(AB)=(x.A)B$\, et \,$1.A=A$. Et \,$\Theta_M$\, est
un automorphisme de cette structure: un \homo bijectif d'anneau qui
v\'erifie en plus \,$\Theta_M(x.A)=x.\Theta_M(A)$.
On en d\'eduit que pour tout
\pol \,$f\in\A[X]$\, on a: \,$\Theta_M(f(A))= f(\Theta_M(A))$.}}
de l'anneau
\,$\L^{n\times n}$, on a
\,$f(A')=M^{-1}f(A)\:M$, et par suite
\,$\rP_{\!f(A)}=\rP_{\!f(A')}=
(-1)^n\prod_{i=1}^n(X-f(\l_i)).$ \qed

\ms En fait, il n'est pas \ncr que les deux
matrices \,$A$\, et \,$A'$\, soient semblables pour
conclure que \,$\rP_{\!f(A)}=\rP_{\!f(A')}$. Il
suffit pour cela que \,$\rPA=\rP_{\!A'}$,
comme l'indique le r\'esultat suivant, valable
dans un \acoma \,$\A$.
\begin{propri}
Soient \,$A$\, et \,$A'$\, deux matrices
carr\'ees \`a \coes  dans \,$\A$\, ayant m\^eme
\polcar \,$\rPA=\rP_{\!A'}$. Alors, pour tout
\,$f\in\A[X]$\, on a: \,$\rP_{\!f(A)}=\rP_{\!f(A')}.$
\end{propri}
\dem
Soit \,$B$\, la \emph{\maco du \pol unitaire} \refstepcounter{bidon}
\index{matrice!compagnon d'un polyn\^ome unitaire}
\index{compagnon!matrice}\label{matrcomp}
\,$(-1)^n\rPA(X)=  \alb
X^n-\alb(c_1X^{n-1}\alb +\cdots\alb +c_{n-1}X+c_n)$, \cad
la matrice
$$B = \cmatrix{
0      & \cdots & \cdots & 0      & c_n \cr
1      &   0    & \cdots & 0      & c_{n-1} \cr
0      & \ddots & \ddots & \vdots & \vdots \cr
\vdots & \ddots & \ddots & 0      & \vdots \cr
0      & \cdots & 0      & 1      & c_{1}
}\in \A^{n\times n}\,.$$
On v\'erifie sans difficult\'e que \,$\rPA(X)=\rP_{\!B}(X)$.
Il suffit de montrer que \,$\rP_{\!f(A)}=
\rP_{\!f(B)}$\, pour conclure (puisqu'alors
on a aussi \,$\rP_{\!f(A')}=
\rP_{\!f(B)}$).
On peut \'ecrire
$$\,\rP_{\!f(B)}=
Q_n(c_1,\ldots,c_n,X) = (-1)^nX^n +
\sum\nolimits_{i=0}^{n-1}q_{n,i}(c_1,\ldots,
c_n)\,X^i\,.$$
Il s'agit alors de montrer que
$$\,\rP_{\!f(A)}=(-1)^nX^n+
\somm_{i=0}^{n-1}q_{n,i}(c_1,c_2,\ldots,c_n)X^i\,.$$
  Si l'on consid\`ere les \coes
de \,$f$\, et les en\-tr\'ees de \,$A$\, comme
\idtrs
\,$x_1,\ldots,x_\ell$\, sur \,$\ZZ$\,
(on a \,$\ell=n^2+1+{\deg}\,f$), \'etablir
l'\egt pr\'e\-c\'edente  revient \`a
d\'emontrer \,$n$\,  \idas dans
\,$\ZZ[x_1,\ldots,x_\ell].$
Ces \idas sont ensuite valables dans tout \acom
\,$\A,$ en rempla\c{c}ant les variables
formelles $x_{i}$ par des \elts $\xi_i$ de
\,$\A$\, et elles donnent le r\'esultat souhait\'e. \\
Or, pour d\'emontrer ces \idasz, il
suffit de les v\'erifier sur un ouvert \,$U$\, de
\,$\CC^\ell$, c'est-\`a-dire lorsqu'on substitue \`a
\,$(x_1,\ldots,x_\ell)$\, un \elt
\,$(\xi_1,\ldots,\xi_\ell)$\, arbitraire de \,$U$.\\
Pour cela, on consid\`ere par exemple l'ouvert
correspondant \`a des matrices
\gui{suffisamment proches} de la matrice diagonale suivante:
$$\cmatrix{
1      & 0      & \cdots & 0 \cr
0      & 2      & \ddots & \vdots  \cr
\vdots & \ddots & \ddots & 0 \cr
0      & \cdots & 0      & n
}\,.$$
Ces matrices sont diagonalisables puisque leurs valeurs
propres restent distinctes dans l'ouvert consid\'er\'e.
Dans ce cas-l\`a, \,$A$\, et \,$B$\, sont
diago\-nali\-sables avec les m\^emes valeurs propres,
et le r\'esultat est trivial. \qed

\section {Polyn\^ome minimal}
\label{subsec polmin}

Soit \,$\K$\, un corps, \,$E$\, un \Kev de dimension \,$n>0$\, et
\,$\varphi:E\dans E$
un op\'erateur \linz, repr\'esent\'e en \gnl
par une matrice \,$A$\, dans une base donn\'ee de \,$E$.

\subsection[Sous-espaces de Krylov]{Polyn\^{o}me minimal et sous-espaces
de Krylov}
\label{subsecKrylov}

Le \emph{\polmin} de \,$\varphi\,$\index{polyn\^ome!minimal}
(ou de \,$A$) est par \dfn le \poly unitaire
\,$\rP^\varphi\in\K[X]$\, engendrant l'id\'eal des \pols
\,$f\in\K[X]$\, tels que \,$f(\varphi)=0$.
A priori, ce \pol peut \^{e}tre calcul\'e par les \mets \usles
d'\agr \lin en cherchant la premi\`ere relation de \dlin entre les
\gui{vecteurs} successifs: \,$\Id_E,\,\varphi,\,\varphi^2,
\,\varphi^3\ldots, $\, de l'\evc \,$\End(E)\simeq \K^{n{\times}n}$\, des
\endos de \,$E$.

Puisque \,$\rP_{\!\varphi}(\varphi)=0$ (\tho de Cayley-Hamilton), on 
obtient que
\,$\rP^\varphi$\, divise \,$\rP_{\!\varphi}$.

\ss Par ailleurs, \'etant donn\'e un vecteur \,$v\in E$, le sous-espace
\emph{$\,\varphi$-engendr\'e} par \,$v$\, (on dit aussi \emph{sous-espace de
Krylov\index{Krylov!sous-espace de} pour le couple \,$(\varphi,v)\,$}) est, par
\dfnz, le sous-espace de \,$E$\, engendr\'e par le \sys de vecteurs
\,$(\varphi^i(v))_{\,i\in \NN}$. On le note \,$\Kr(\varphi,v)$.

  Sa dimension
n'est autre que le degr\'e du \poly unitaire \,$\rP^{\varphi,v}\in\K[X]$\, qui
engendre l'id\'eal \,$\{f\in\K[X]\,|\,f(\varphi)(v)=0\}$\, de 
\,$\K[X]$. Ce \poly
s'appelle \emph{le \poly $\varphi$-minimal de} \,$v$. Les deux \pols
\,$\rP^\varphi$\, et \,$\rP_{\!\varphi}$\, appartiennent \'evidemment \`a cet
id\'eal et on a:
$$
\Kr(\varphi,v)=E \;\equiva \; \dg{\rP^{\varphi,v}}=n
\;\equiva \; \rP^{\varphi,v}=\rP^\varphi= (-1)^n\rP_{\!\varphi}\,.
$$
Ainsi, l'existence d'un vecteur \,$v$\, qui \,$\varphi$-engendre
\,$E$\, suffit pour que le \polcar de l'\endom
\,$\varphi$\, de \,$E$\, soit \'egal (\`a un signe pr\`es)
\`a son \polminz:
$$
\exists v\in E,\,\dg{\rP^{\varphi,v}}=n\; \implique\;
\rP^\varphi=(-1)^n\rP_{\!\varphi}\,.
$$

Mais la r\'eciproque est aussi vraie, comme nous allons le voir
bient\^{o}t.

Remarquons que, comme le \pol  \,$\rP^{\varphi}$, le \pol
\,$\rP^{\varphi,v}$\,
peut \^{e}tre calcul\'e par les \mets \usles
d'\agr \linz.

Rappelons maintenant
la propri\'et\'e classique suivante relative \`a la d\'ecomposition
de \,$E$\, en sous espaces \,$\varphi$-stables:

\begin{propri}\label{sdpo}
Soit \,$\varphi$\, un \endom de \,$E$\,
et \,$f\in \K[X]$\, tel que \,$f=f_1f_2\cdots f_r$\, avec
\,$\pgcd(f_i,f_j)\,=1$\, pour \,$i\neq j$\, et \,$f(\varphi)=0$.
Posons \,$g_i=f/f_i$, \,$\theta_i=f_i(\varphi)$,
\,$\psi_i=g_i(\varphi)$, \,$E_i=\Ker\,{\theta_i}$. Chaque   \,$E_i$\, est
un sous-espace \,$\varphi\,$-stable, on note \,$\varphi_i:E_i\rightarrow E_i$\,
la restriction de \,$\varphi$.
Alors:
\begin{itemize}
\item [a)] $E =  E_1\oplus \cdots \oplus E_r$.
\item [b)] $E_i= \Im\,{\psi_i}$,  et la restriction  de
\,$\psi_i$\, \`a \,$E_i$\,  induit  un automorphisme de
l'\evc~$\,E_i$.
\item [c)] On a les \egts
\,$\rP_{\!\varphi}=\prod_{i=1}^r\,\rP_{\!\varphi_i}$\,
et  \,$\rP^{\varphi}=\prod_{i=1}^r\,\rP^{\varphi_i}$.
\item [d)] Si  \,$(v_1,\ldots ,v_r)\in E_1\times\cdots {\times}E_r$, et
\,$v=v_1+\cdots +v_r,$ on a l'\egt
\,$\rP^{\varphi,v}=\prod_{i=1}^r\,\rP^{\varphi_i,v_i}$.
\end{itemize}
\end{propri}
\prv
Tout d'abord on remarque que deux \endos de la forme \,$a(\varphi)$\, et
  \,$b(\varphi)$\, commutent toujours puisque \,$a(\varphi)\circ
b(\varphi)=(ab)(\varphi)$. Ensuite il est clair que tout sous-espace
du type \,$\Ker\,a(\varphi)$\, ou \,$\Im\,a(\varphi)$\, est
\,$\varphi$-stable.

La preuve des points a) et b) est alors bas\'ee sur les
\egts \,$\varphi_i\circ\psi_i=\psi_i\circ\varphi_i=f(\varphi)=0_E$\, et
sur l'\idt de B\'ezout \,$u_1\,g_1+\cdots +u_r\,g_r=1$, qui implique
\,$u_1(\varphi)\circ g_1(\varphi)+\cdots +u_r(\varphi)\circ
g_r(\varphi)=\Id_E$, ce qui se lit: \,$\psi_1\circ\alpha_1+\cdots +
\psi_r\circ\alpha_r=\Id_E$, o\`u \,$\alpha_i = u_i\,(\varphi)$.

Si on choisit dans chaque \,$E_i$\, une base \,$B_i$\, leur r\'eunion
est une base \,$B$\, de \,$E$\, et la matrice de \,$\varphi$\, sur
\,$B$\, est diagonale par blocs, chaque bloc \'etant la matrice de
  \,$\varphi_i$\, sur \,$B_i$.
Ceci implique \,$\rP_{\!\varphi}=\prod_{i=1}^r\,\rP_{\!\varphi_i}$.

Soit \,$h_i=P^{\varphi_i}$.
  Il est clair que \,$(h_1\cdots h_r)(\varphi)=0_E$
(car nul sur chaque \,$E_i$), et
que les \,$h_i$\, sont premiers entre eux 2 \`a 2
(car \,$h_i$\, divise  \,$f_i$). Si \,$g(\varphi)=0$, a fortiori
\,$g(\varphi_i)=0$, donc \,$g$\, est multiple des \,$h_i$, et par suite
multiple de \,$h=h_1\cdots h_r$. Ceci termine la preuve du point c).

Et la preuve du point d) est analogue, en rempla\c{c}ant
\,$h_i=P^{\varphi_i}$\, par \,$P^{\varphi_i,v_i}$\,. \qed

\begin{propri}\label{polmin}
Soit \,$\rP^\varphi = P_1^{m_1}\,P_2^{m_1}\,\ldots\,P_r^{m_r}$\,
la \deco du \polmin de \,$\varphi$\, en facteurs
premiers distincts \,$P_1$, \,$P_2$, \,$\ldots$, \,$P_r.$
On reprend les notations de la propri\'et\'e \ref{sdpo} avec 
\,$f_i=P_i^{m_i}$. On
pose en outre \,$Q_i=\rP^\varphi/P_i,$  \,$G_i=\Ker(Q_i(\varphi)),$
\,$F_i=\Ker(P_i^{m_i-1}(\varphi)),$ pour
\,$1\leq i\leq r$.
Alors \,$E=E_1\oplus\cdots\oplus E_r,$  chaque \,$F_i$\,  est
strictement inclus dans \,$E_i$\, et pour tout
\,$v=v_1+\cdots+v_r$ ($ v_i\in E_i$) on a les \'equivalences suivantes:
$$ v\notin\bigcup_{i=1}^r  G_i \equiva
P^{\varphi,v}=P^{\varphi} \equiva \bigwedge_{i=1}^r\;
P^{\varphi_i,v_i}=P^{\varphi_i}
  \equiva \bigwedge_{i=1}^r\; v_i\notin F_i
  .$$
\end{propri}
\prv Le fait que  \,$E=E_1\oplus\cdots\oplus E_r$\, r\'esulte de la
propri\'et\'e \ref{sdpo}~a).
La deuxi\`eme \'equivalence r\'esulte des
propri\'et\'es \ref{sdpo}~c) et d).
La premi\`ere \'equivalence est claire: la
premi\`ere condition signifie exactement
que le  \pol  \,$\varphi$-minimal de \,$v$\, ne divise pas strictement
le \polmin  de \,$\varphi$.
M\^{e}me chose pour la derni\`ere \'equivalence (tout diviseur strict de
\,$P_i^{m_i}$\, divise  \,$P_i^{m_i-1}).$   Cette remarque montre aussi
que l'inclusion \,$F_i\subset E_i$\, est stricte.
\qed

\begin{corollary} \label{corpolmin}
  Il existe toujours un vecteur \,$v$\, tel que
$\rP^{\varphi,v}=\rP^{\varphi}$.
En particulier, si le \polmin de \,$\varphi$\,  a pour degr\'e la
dimension de \,$E$\, (autrement dit si, au signe pr\`es, il est \'egal
au
\polcarz)
il existe des vecteurs qui \,$\varphi$-engendrent \,$E$.
\end{corollary}
\prv
Il suffit de prendre \,$v=v_1+\cdots+v_r$ avec
chaque \,$v_{i}$\, dans \,$E_i\setminus F_i.$\qed

\ss Ainsi, sauf exception, \cad si l'on
choisit \,$v$\, en dehors de la r\'eunion d'un petit nombre
de sous-\evcs stricts de \,$E$,
les \,$\Ker{(Q_i(\varphi))},$  (cette r\'eunion ne remplit jamais
l'espace, m\^eme si le corps de base est fini), le vecteur
\,$v$\, convient.

Notons que la preuve pr\'ec\'edente est peu satisfaisante parce
qu'en \gnl on ne sait pas
calculer la \deco en facteurs premiers d'un \polz. Il
s'ensuit que notre preuve de
l'existence d'un vecteur qui \,$\varphi$-engendre \,$E$\, reste
plus th\'eorique que pratique.

Voici une mani\`ere de contourner cet obstacle.

Premi\`erement on \'etablit le lemme suivant.

\begin{lemma} \label{lemFactor}
Si on conna\^{\i}t un  \pol \,$Q$\, qui est un facteur strict
d'un \pol \,$P\,$
de \,$\K[X]$\, on peut, ou bien d\'ecomposer \,$P$\, en un produit
\,$P_1\,P_2$\, de deux \pols \'etrangers, ou bien \'ecrire
\,$P$\, et \,$Q$\, sous la forme \,$P_1^k$\, et  \,$P_1^\ell$\,
($1\leq \ell<k$).
\end{lemma}

La \pcd (que nous ne donnons pas en d\'etail)
consiste \`a partir de la factorisation intitiale
\,$P=QQ_1$\, et \`a la raffiner au maximum en utilisant les pgcd's.

Ensuite on rappelle que les \pols $\rP^{\varphi,v}$ et $\rP^{\varphi}$
peuvent facilement \^etre calcul\'es par les \mets classiques
d'\agr \linz.

On d\'emarre alors avec un \,$v$\, non nul arbitraire.
Si $\rP^{\varphi,v}$ est \'egal \`a $\rP^{\varphi}$ on a termin\'e.
Sinon, on applique le lemme pr\'ec\'edent
avec \,$P=\rP^{\varphi}$ et \,$Q=\rP^{\varphi,v}$.
Dans le premier cas, on applique la propri\'et\'e \ref{sdpo}~c) et d)
avec les \pols \,$P_1$\, et \,$P_2$.
On est ramen\'e \`a r\'esoudre le m\^eme \pb s\'epar\'ement
dans deux espaces de dimensions plus petites.
Dans le deuxi\`eme cas, on choisit un nouveau \,$v$\, dans
le compl\'ementaire de $\Ker(P_1^{k-1}(\varphi))$, ce qui fait
que le degr\'e de son \polmin augmente strictement.


\ms\rem 
Toute \deco
\,$E=\bigoplus_{i=1}^kE_i$\,  en sous-espaces
\,$E_i$\, \,$\varphi\,$-stables
donne une forme r\'eduite de \,$\varphi$,
\cad la repr\'esentation de \,$\varphi$\, dans une base convenable
par une matrice diagonale par blocs de la forme:
$$A = \cmatrix{
A_1 & 0 & \cdots  & 0 \cr
0 & A_2 & \ddots  & \vdots  \cr
\vdots & \ddots  & \ddots & 0 \cr
0 & \cdots  & 0 & A_k
}\,.$$
o\`u \,$A_1,\,A_2,\,\ldots,\,A_k$\, sont les matrices des
\endos \,$\varphi_i$\, induits par \,$\varphi$\, dans les
sous-espaces \,$E_i$.
Certaines de ces formes r\'eduites sont dites canoniques, comme
la r\'eduction de Jordan dans le cas o\`u le \polcar se factorise en
facteurs lin\'eaires sur \,$\K$. Il existe d'autres formes r\'eduites
canoniques enti\`erement rationnelles, \cad qui n'utilisent pour le
changement de base que des expressions rationnelles en les \coes de la
matrice de d\'epart.
Sur le sujet des formes normales r\'eduites (et sur bien d'autres) nous
recommandons
le livre de Gantmacher
(\cite{Gan}) dont on attend toujours la r\'e\'edition \`a un prix
abordable.

\subsection{Cas de matrices \`a \coes dans \texorpdfstring{$\ZZ$}{Z}.} 
\label{pageAIC}
Dans certains \algos que nous aurons \`a d\'evelopper par la suite, nous
partirons d'une matrice \,$C$\, \`a \coes dans un anneau int\`egre
\,$\A$\, et bien souvent il sera avantageux qu'aucun des calculs \itmds
ne produise
des \elts qui seraient dans le corps des fractions de \,$\A$\, sans
\^{e}tre dans \,$\A$.

Se pose alors naturellement la question suivante:  les \pols
\,$\rP^{C}$\, et \,$\rP^{C,v}$\, que nous pouvons \^{e}tre amen\'es \`a
envisager comme r\'esultats de calculs \itmds en vue de trouver le
\polcar de \,$C$,  sont-ils toujours des \pols \`a \coes dans \,$\A$?

Dans le cas de \macas \`a \coes dans  $\ZZ$, ou dans un anneau de \pols
$\ZZ[x_1,\ldots ,x_n]$\, ou $\QQ[x_1,\ldots ,x_n]$,  la r\'eponse est
positive.

Ce r\'esultat n'est pas \'evident. Il est bas\'e sur les \dfns et les
propri\'et\'es qui suivent, pour lesquelles on peut consulter
les livres classiques d'alg\`ebre (par exemple \cite{Gob} ou
\cite{MRR}).

\begin{definition}
\label{defAIC} Un anneau int\`egre \,$\A$\, est dit
\emph{int\'egralement clos} si tout diviseur unitaire
dans \,$\FA[X]$\, d'un \pol unitaire de \,$\A[X]$\, est
dans  \,$\A[X]$.
\end{definition}

Avec de tels anneaux les \pols \,$\rP^{C}$\, et \,$\rP^{C,v}$\,
sont donc automatiquement \`a \coes dans \,$\A$.

\begin{definition}
\label{defAFa} Un anneau int\`egre \,$\A$\, est dit \emph{anneau \`a 
pgcds} si tout
couple d'\elts \,$(a,b)$\, admet un pgcd, \cad un \elt \,$g\in \A$\, 
tel que: $$
\forall x\in \A\quad  ((x\;\mathrm{divise}\;a\;\mathrm{et}\;b)
\;\Longleftrightarrow\;x\;\mathrm{divise}\;g )\,. $$
\end{definition}

\begin{propri}
\label{propriAIC1}
Pour qu'un anneau int\`egre soit int\'egralement clos, il suffit
que la propri\'et\'e qui le d\'efinit soit v\'erifi\'ee pour les
diviseurs de degr\'e \,$1$. Autrement dit, tout z\'ero dans \,$\FA$\,
d'un \pol unitaire de \,$\A[X]$\, est dans \,$\A$.
\end{propri}

\begin{propri}
\label{propriAIC2}
Tout anneau \`a pgcds est int\'egralement clos.
\end{propri}

\begin{propri}
\label{propriAIC3}
Si \,$\A$\, est un anneau \`a pgcds il en va de m\^{e}me pour
$\A[x_1,\ldots ,x_n]$.
\end{propri}


\section {Suites r\'ecurrentes lin\'eaires}
\label{subsec suites rec}

\subsection{Polyn\^{o}me g\'en\'erateur, op\'erateur de d\'ecalage}

Soit \,$E$\, un \Kev (resp.  un \Amoz).
On consid\`ere une suite \,$(a_n)_{n\in \NN}$\, d'\elts de \,$E\,$
et un entier \,$p\in \NN.$
Une \emph{\rrl d'ordre \,$p\,$}
pour cette suite est d\'efinie par la donn\'ee de \,$p+1$\,
\elts \,$c_0,$ \,$c_1,\ldots,$ \,$c_p$\, de \,$\K$\,
(resp. de \,$\A$) v\'erifiant:
%
\begin{equation} \label{EqRecLin}
\forall n\in \NN ~~~~ c_0a_n + c_1a_{n+1} + \cdots + c_pa_{n+p}
= 0\quad \quad
\end{equation}
%
Le \poly \,$h(X)=\sum_{i=0}^pc_iX^i$\,  dans \,$\K[X]$\,
(resp. dans \,$\A[X]$) est appel\'e
\emph{un \polgz} de la suite
\,$(a_n)$.\index{polyn\^ome!g\'en\'erateur d'une suite}
Lorsque le \coe \,$c_p$\, est inversible, la suite est alors
d\'etermin\'ee par la donn\'ee de \,$a_0,$ \,$a_1,\ldots,$
\,$a_{p-1}$\, car elle peut ensuite \^{e}tre construite
par \recuz: elle est en quelque sorte \gui{engendr\'ee}
par le \pol \,$h$, ce qui justifie la terminologie adopt\'ee.
Une \emph{\srlz}
\index{suite r\'ecurrente lin\'eaire}
dans  \,$E$\, est une suite
\,$(a_n)_{n\in \NN}$\, d'\elts de \,$E$\,
qui poss\`ede un \polg dont le \coe
dominant est inversible.

On interpr\`ete  cette situation de la mani\`ere
suivante en \agr \linz.
On appelle \,$\cS$\, le \Kev (resp. le \Amoz) \,$E^\N$\, form\'e de
toutes les suites \,$(u_n)_{n\in \NN}$\, \`a valeurs dans
  \,$E$. On note \,$\Phi:\cS\rightarrow \cS$\, \emph{l'op\'erateur de
d\'ecalage} qui donne pour image de  \,$(u_n)_{n\in \NN}$\,
la suite d\'ecal\'ee d'un cran  \,$(u_{n+1})_{n\in \NN}.$
Il est clair que l'op\'erateur de d\'ecalage est un op\'erateur
\linz. Dire que la suite  \,$\s{a}=(a_n)_{n\in \NN}$\, v\'erifie
la \rrl (\ref{EqRecLin}) se traduit
exactement, en langage un peu plus abstrait, par:
$$ \s{a}\;\in\;\Ker(h(\Phi))\,.
$$
Cela montre que les \pols g\'en\'erateurs
d'une \srl donn\'ee  forment un
id\'eal de \,$\K[X]$\, (resp. de \,$\A[X]$).
Dans le cas d'un corps et d'une \srlz, comme \,$\K[X]$\,
est un anneau principal, cet id\'eal (non nul) est
engendr\'e par un \poly unitaire unique
qu'on appelle \emph{le \polgminz} ou simplement
\emph{le \polminz} de la suite.
\index{polyn\^ome!minimal d'une suite r\'ecurrente lin\'eaire}
\index{polyn\^ome!g\'en\'erateur minimal d'une suite
r\'ecurrente lin\'eaire}
Nous le noterons \,$\rP^{\s{a}}$.
Nous allons voir plus loin que ce \pol peut effectivement \^{e}tre
calcul\'e, d\`es qu'on conna\^{\i}t un \polg de la \srlz.

Consid\'erons maintenant un \pol unitaire fix\'e
$$ \,f(X)=X^p-\somm_{i=0}^{p-1}b_iX^i\,
$$
dans \,$\K[X]$\,
(resp. dans \,$\A[X]$). Le \Kev (resp. le \Amoz) \,$\Ker(f(\Phi))$\,
form\'e des \srls dans \,$E$\, pour lesquelles
\,$f$\, est un \polg sera not\'e
\,$\cS_{\!f}.$
Il est isomorphe \`a \,$\K^p$\,
(resp. \,$\A^p$) et une base canonique
est fournie par les \,$p$\, suites\footnote{~Chacune de ces suites est
\'evidemment d\'efinie par ses  \,$p$\, premiers termes.}
\,$e^{(i)}\;(i=0,\ldots ,p-1)$ telles
que \,$e^{(i)}(j)=\delta_{ij}$\, pour \,$j=0,\ldots ,p-1$,
o\`u $\delta_{ij}$ est le \refstepcounter{bidon}\label{Kron} 
\emph{symbole de Kronecker}
($\delta_{ij}=1$
\indexnota{delta.ij@$\delta_{ij}$ (symbole de Kronecker)} si
$i=j$ et $\delta_{ij}=0$ si $i\neq j$).
Pour une
\srl arbitraire \,$\s{a}=(a_n)_{n\in \NN}$\,
dans \,$\cS_{\!f}$\,
on a alors: \,$\s{a}=\sum_{j=0}^{p-1}a_je^{(j)}$.

Il est clair que \,$\cS_{\!f}$\, est stable par \,$\Phi$. Notons
\,$\Phi_{\!f}$\, la restriction de \,$\Phi$\, \`a  \,$\cS_{\!f}$.
On constate \immt que la matrice de \,$\Phi_{\!f}$\, sur
la base canonique $(e^{(p-1)},\ldots,e^{(1)},e^{(0)})$
est la \maco du \polz~$f$:
$$\rC_f =
\cmatrix{
\matrix{0 & \cdots & 0}
           & b_0                                \cr
  \I_{p-1} & \matrix{b_{1} \cr \vdots \cr b_{p-1}}
}
=\cmatrix{
0      & \cdots & \cdots & 0      & b_0 \cr
1      & 0      & \cdots & 0      & b_{1} \cr
0      & \ddots & \ddots & \vdots & \vdots \cr
\vdots & \ddots & \ddots & 0      & \vdots \cr
0      & \cdots & 0      & 1      & b_{p-1}
}\,.$$
En particulier
\,$\rP^{\Phi_{\!f},e^{(1)}}=\rP^{\Phi_{\!f}}=\rP_{\Phi_{\!f}}=f$.
En outre, par simple application des \dfns on obtient
\,$\rP^{\Phi_{\!f},\s{a}}=\rP^{\s{a}}$.

Comme exemples importants de \srlsz, on peut citer:

-- la \srl form\'ee des
puissances \,$(A^n)_{n\in\NN}$\, d'une matrice
\,$A\in\K^{m\times m}$\, dont le \poly
g\'en\'erateur minimal n'est autre que le \poly
minimal \,$\rP^A$\, de la matrice;

-- pour des vecteurs donn\'es
\,$u,v\in\K^{m\times 1}$, les \srls
  \,$(A^n\,v)_{n\in \NN}$\, et
\,$(\tra{u}\,A^n\,v)_{n\in \NN}$\, dont les \pols
minimaux, not\'es respectivement \,$\rP^{A,v}$\,
et \,$\rP_{u}^{A,v}$, sont alors tels que:
$~\rP_{u}^{A,v}$\, divise \,$\rP^{A,v}~$
et $~\rP^{A,v}$\, divise \,$\rP^A$.

\subsection[Matrices de Hankel]{\Srls et matrices de Hankel}

Comme le montre
la discussion qui suit,
l'\'etude des \srls est \'etroitement li\'ee
\`a celle des
\emph{matrices de Hankel}.
\index{Hankel!matrice de}\index{matrice!de Hankel}

Une matrice de Hankel est une matrice (pas \ncrt carr\'ee)
\,$H=(v_{ij})$\, dont les
\coes sont constants sur les diagonales montantes:
\,$v_{ij}=v_{hk}$\, si \,$i+j=h+k$.

Les matrices de Hankel fournissent un exemple
de \emph{matrices structur\'ees\index{matrice!structur\'ee}}.
L'autre exemple le plus important est celui des
\index{Toeplitz!matrice de}\index{matrice!de Toeplitz}
\emph{matrices de Toeplitz,}
celles dont les
\coes sont constants sur les diagonales descendantes:
\,$v_{ij}=v_{hk}$\, si \,$i-j=h-k$.

Remarquons qu'une matrice de Hankel carr\'ee d'ordre \,$n$\, est une
matrice sym\'etrique et que les produits \,$H\,\Jn$\, et \,$\Jn\,H$\,
d'une
matrice de Hankel \,$H$\, carr\'ee d'ordre \,$n$\, par la matrice de
Hankel
particuli\`ere \,$\Jn$:
$$\Jn = \left[
\begin{array}{ccccc}
0 & 0 & \cdots  & 0 & 1 \\
0 & 0 & \cdots & 1 & 0 \\
\vdots & \vdots  & ~ & \vdots & \vdots \\
0 & 1 & \cdots  & 0 & 0 \\
1 & 0 & \cdots  & 0 & 0
\end{array} \right] $$
sont des \mtosz. Cette matrice de permutation d'ordre \,$n$\,
  permet de renverser l'ordre des
\,$n$\, colonnes (resp. des \,$n$\, lignes) d'une matrice lorsque
celle-ci
est multipli\'ee \`a droite (resp. \`a gauche) par la matrice \,$\Jn$:
c'est pourquoi on l'appelle \emph{matrice de renversement} ou encore
\emph{matrice d'arabisation} du fait qu'elle permet d'\'ecrire de
droite \`a gauche les colonnes que l'on lit de gauche \`a droite
et inversement.

Inversement, les produits \,$\Jn\, T$\, et \,$T\,\Jn$\,
d'une \mto \,$T$\, carr\'ee d'ordre \,$n$\, par la matrice
\,$\Jn$\, sont des matrices de Hankel.

Une matrice structur\'ee est d\'etermin\'ee par la donn\'ee de
beaucoup moins de \coes qu'une matrice ordinaire de m\^{e}me
taille. Par exemple une matrice de Hankel (resp. de Toeplitz)
de type $(n,p)$ est d\'eter\-mi\-n\'ee par la donn\'ee de \,$n+p-1$\,
\coesz: ceux des premi\`ere ligne et derni\`ere (resp. premi\`ere)
colonne.
Cela rend ces matrices particuli\`erement importantes
pour les \gui{grands calculs}  d'\agr \linz.

Si \,$\s{a}=(a_n)_{n\in \NN}$\, est une suite arbitraire
et si \,$i,$ $r,$ $p\in\N$\,  nous noterons
\,$\rH^{\s{a}}_{i,r,p}$\, la matrice de Hankel suivante,
qui poss\`ede \,$r$\, lignes et \,$p$\, colonnes:
$$\rH^{\s{a}}_{i,r,p}=
\left[ \begin{array}{ccccc}
a_i & a_{i+1} & a_{i+2}   &\ldots & a_{i+p-1}  \\
a_{i+1} & a_{i+2} &       &       & a_{i+p}  \\
a_{i+2} &     &       &       &     \\
\vdots &  &       &       & \vdots   \\
a_{i+r-1}&a_{i+r}&\ldots&\ldots & a_{i+r+p-2}
\end{array}\right]\,.
$$

Le fait suivant est une simple constatation.
\begin{fact} \label{factSRL}
  Reprenant les notations ci-dessus,
une suite  \,$\s{a}$\, est une \srl avec
\,$f$\, comme \pol g\'en\'erateur
\ssi sont v\'erifi\'ees pour tous  \,$i,$ $r\in\N$\, les \'equations
matricielles
\begin{equation} \label{EqSRL1}
\rH^{\s{a}}_{i,r,p}\,\rC_f=\rH^{\s{a}}_{i+1,r,p}
\end{equation}
Ou ce qui revient au m\^{e}me, en transposant,
$ \,\tra{\rC_f}\,\rH^{\s{a}}_{i,p,r}=\rH^{\s{a}}_{i+1,p,r}$.
Naturellement, il suffit que ces \'equations soient v\'erif\'ees
lorsque \,$r=1$.
\end{fact}
On en d\'eduit
\begin{equation} \label{EqSRL2}
\rH^{\s{a}}_{i,r,p}\,(\rC_f)^k=\rH^{\s{a}}_{i+k,r,p}
\end{equation}
et donc on a aussi:
\begin{fact} \label{factSRL2}
  Sous les m\^{e}mes hypoth\`eses,
dans toute matrice \,$\rH^{\s{a}}_{i,r,p+k}$\, les
\,$k$\, derni\`eres colonnes sont \colis des \,$p$\, premi\`eres.
Et, par transposition, dans toute matrice \,$\rH^{\s{a}}_{i,p+k,s}$\,
les \,$k$\, derni\`eres lignes sont \colis des \,$p$\, premi\`eres.
\end{fact}

On en d\'eduit la proposition suivante.
\begin{proposition} \label{propSRL}
Avec les notations ci-dessus, et dans le cas d'un corps
\,$\K,$
si \,$\s{a}$\, est une \srl qui admet \,$f$\, pour
\pol g\'en\'erateur,
  le degr\'e \,$d$\, de son \pol g\'en\'erateur minimal
\,$\rP^{\s{a}}$\, est \'egal au rang de la matrice de Hankel
\,$\rH^{\s{a}}_{0,p,p}.$
Les \coes de
\,$\rP^{\s{a}}(X)=X^d-\sum_{i=0}^{d-1}g_iX^i\in\K[X]\,$
sont l'unique solution de l'\'equation
\begin{equation} \label{EqSRL3}
\rH^{\s{a}}_{0,d,d}\,\rC_{\rP^{\s{a}}}=\rH^{\s{a}}_{d,d,1}
\end{equation}
\cad encore l'unique solution du \sli
$$
\left[ \begin{array}{ccccc}
a_0 & a_{1} & a_{2}   &\cdots & a_{d-1}  \\
a_{1} & a_{2} &       & \adots  & a_{d}  \\
a_{2} &     & \adots & \adots & \vdots    \\
\vdots &\adots   & \adots  &       & \vdots   \\
a_{d-1}&a_{d}&\cdots&\cdots & a_{2d-2}
\end{array}\right] \;
\left[ \begin{array}{c}
g_0 \\ g_1 \\ g_2 \\  \vdots   \\  g_{d-1}
\end{array}\right]
=
\left[ \begin{array}{c}
a_d \\  a_{d+1} \\ a_{d+2} \\  \vdots \\ a_{2d-1}
\end{array}\right] \,.
$$
\end{proposition}
\prv
Consid\'erons la premi\`ere \rdl entre les colonnes de
  \,$\rH^{\s{a}}_{0,p,p},$ et appelons \,$g$\, le \pol
unitaire correspondant, qui v\'erifie
  $$\,\tra{\rC_g}\,\rH^{\s{a}}_{0,d,p}=\rH^{\s{a}}_{1,d,p}\,.$$
Cela donne l'\'equation
\,$\tra{\rC_g}\,\rH^{\s{a}}_{0,d,p}=\rH^{\s{a}}_{0,d,p}\,\rC_f$\, d'o\`u
\immt par \recu sur \,$k\;$:
$$\,(\tra{\rC_g})^k\,\rH^{\s{a}}_{0,d,p}=
\rH^{\s{a}}_{0,d,p}\,(\rC_f)^k=\rH^{\s{a}}_{k,d,p}\,,$$
et donc
$$\,\tra{\rC_g}\,\rH^{\s{a}}_{k,d,p}=\rH^{\s{a}}_{k+1,d,p}\,.$$
A fortiori, pour tout \,$k$\,  on a:
$$\,\tra{\rC_g}\,\rH^{\s{a}}_{k,d,1}=\rH^{\s{a}}_{k+1,d,1}\,,$$
et cela signifie que \,$g$\, est \egmt un \pol
g\'en\'erateur pour la suite \,$\s{a}$. \\
On laisse le soin \`a la lectrice et au lecteur de finir la preuve.
\qed

\section {Polyn\^omes sym\'etriques et relations de Newton}
\label{subsec Newton}

Soit \,$\A$\,  un \acom unitaire
et \,$\A\,[x_1, \dots , x_n]$\, l'\agr sur
\,$\A$\, des \pols \`a \,$n$\,
\idtrs $x_1,\dots , x_n$.

Tout \pol \,$f\in\A\,[x_1, \dots , x_n]$\,
s'\'ecrit de mani\`ere unique comme une somme finie
de mo\-n\^o\-mes distincts \,$a_J\,x^J=a_J\, x_1^{j_1}\cdots
x_n^{j_n}$\, o\`u $J=(j_1,\dots,j_n)\in \NN^n$\, et
\,$a_J\in\A$:
$$\,f=\sum\nolimits_{J} a_J\; x_1^{j_1}\cdots
x_n^{j_n}\,
$$
o\`u la somme porte sur une
partie finie de $\NN^n$.
On a souvent int\'er\^et \`a donner un bon ordre
sur les termes \,$x^J=x_1^{j_1}\cdots  x_n^{j_n}\,$
pour faire des preuves par induction.
Il suffit par exemple d'ordonner les \idtrs
($x_1<x_2<\cdots <x_n$) pour d\'efinir un bon ordre
sur les termes \,$x^J$,
par exemple l'ordre lexicographique, ou l'ordre lexicographique
subordonn\'e au degr\'e total.

On \'ecrit souvent aussi le \pol comme somme de ses composantes
homog\`enes
$$f=\sum\nolimits_{h=0}^n\, f_h~,
$$
o\`u \,$n$\, est le degr\'e
total de \,$f$\, et \,$f_h$\, la composante homog\`ene
de degr\'e \,$h$\, de \,$f$.
Une mani\`ere simple de voir les composantes
homog\`enes d'un \poly \,$f\in\A\,[x_1,\dots,x_n]$\,
est de consid\'erer une nouvelle \idtr
\,$z$\, et le \poly en \,$z$:
$g(z)=f(x_1z,\alb \dots,\alb x_nz)\in \A\,[x_1,\alb
\dots,\alb  x_n]\,[z].$
La composante homog\`ene
  \,$f_h$\, n'est autre que le
\coe de \,$z^h$\, dans \,$g(z)$.

\ss D\'esignant par \,${\cal S}_n$\, le groupe des
permutations de \,$\{1,2,\dots,n\}$,
un \pol \,$f\in\A\,[x_1,\dots, x_n]$\,
est dit \emph{\symz} si son stabilisateur par
l'action
$$f= f(x_1, \dots , x_n) \stackrel{\tau}{\vers}
\tau f = f(x_{\tau(1)}, \dots , x_{\tau(n)})$$
de \,${\cal S}_n$\, sur \,$\A\,[x_1, \dots , x_n]$\, est
le groupe \,${\cal S}_n$\, tout entier. \Cad encore si
les mon\^omes d'une m\^eme orbite de \,${\cal S}_n$\,
figurent dans l'expression de \,$f$\, avec le m\^eme
\coez.

\ss Des \pols \syms importants sont les
  \emph{sommes de Newton \`a \,$n$\,
\idtrsz}: \index{Newton!somme de}\index{somme de Newton}
\begin{equation} \label{DefSN}
S_k(x_1, \dots , x_n) = \sum\nolimits_{i=1}^n
x_i^k~~\in \A\,[x_1, \dots ,x_n]~~~(k\in \NN)
\end{equation}

On notera \,$\A\,[x_1,\dots,x_n]_{sym}$\,
l'ensemble des \pols \syms en
\,$x_1, \dots , x_n$\, sur \,$\A$. C'est une
sous-\agr propre de \,$\A\,[x_1,\dots, x_n]$.
Il est bien connu (la d\'emonstration peut se faire par
\recu en utilisant l'ordre lexicographique) que
tout \pol \sym s'exprime de mani\`ere unique
comme \pol en \,$\sigma_1,\sigma_2,\dots,\sigma_n$\,
o\`u les \,$\sigma_p~(1\leq p\leq n)$\, sont les
\emph{\pols \syms \elrsz\/} en
\,$x_1,x_2,\dots , x_n$:
$$\sigma_p=\sum_{1\leq i_1<i_2<\cdots<i_p\leq n}
x_{i_1}x_{i_2}\cdots x_{i_p}~\cdot
$$
Cela signifie que l'\homo  de \,$\A$\,--\agrs
$$~\varphi:\A\,[y_1,\dots,y_n]\longrightarrow
\A\,[x_1, \dots ,x_n]_{sym}\,
$$
  d\'efini par
$$\varphi(f(y_1,\dots,y_n)):=
f(\sigma_1,\sigma_2,\dots,\sigma_n)
$$
est un  \isoz.

  En outre, lorsque \,$\A$\, est
int\`egre, \,$\FA$\, d\'esignant le corps des fractions
de \,$\A$, cet \iso \,$\varphi$\, se prolonge de
mani\`ere unique en un \iso de
\,$\FA\,$--\,\agrsz, de \,$\FA(y_1,\dots,y_n)$\,
vers \,$\FA (x_1,\dots , x_n)_{sym}$, qui est par
\dfn la sous-\agr de \,$\FA (x_1,\dots ,
x_n)$\, form\'ee des fractions rationnelles invariantes
par permutation des variables. Autrement dit, toute
fraction rationnelle sur \,$\FA$, \sym en
\,$x_1,\dots , x_n$, s'\'ecrit de mani\`ere unique
comme une fraction rationnelle
en $\sigma_1,\dots ,\sigma_n$ sur \,$\FA$.
On exprime ce fait en disant que
\,$(\sigma_1,\dots , \sigma_n)$\, est un
\emph{\sys fondamental} de \pols \syms
en \,$x_1,\dots,x_n$\, sur le corps \,$\FA$.
Plus \gnltz:

\begin{defi} ~ \label{defsfps}
\'Etant donn\'e un \acom
unitaire \,$\A$\, (resp. un corps \,$\K$), et
\,$n$\, \idtrs $x_1,\dots,x_n$ sur
cet anneau (resp. ce corps), on appelle
\emph{\sys fondamental de \pols \syms en
  \,$x_1,\dots,x_n$\, sur l'anneau \,$\A$\,}
(resp. \emph{fractions rationnelles \syms en
  \,$x_1,\dots,x_n$\,  sur le corps $\K$}) tout \sys
$(f_1,\dots, f_n)$ de \,$n$\, \elts de
$\A\,[x_1,\alb\dots,x_n]_{sym}$\,
(resp. de $\K(x_1,\alb  \dots,\alb x_n)_{sym}$)
v\'erifiant \,$\A\,[x_1,\alb \dots,\alb x_n]_{sym}\alb =\alb
\A\,[f_1,\alb \dots,\alb f_n]$\,  (resp.
$~\K(x_1,\alb \dots,\alb x_n)_{sym}=\K(f_1,\alb\dots,f_n)$).
\index{syst\`eme fondamental!de \pols sym\'etriques}
\index{syst\`eme fondamental!de fractions rationnelles
sy\-m\'e\-tri\-ques}
\end{defi}

Attention \`a l'ambiguit\'e de langage: un \sys
fondamental \emph{sur le corps} \,$\K$\, n'est pas
\ncrt un \sys fondamental \emph{sur l'anneau}
\,$\K$, m\^eme s'il est form\'e de \polsz.
Par contre, un \sys fondamental sur l'anneau
int\`egre \,$\A$\, est toujours un \sys fondamental
sur le corps \,$\FA$.

La \dfn d'un \sys fondamental  sur un corps \,$\K$\, implique
l'in\-d\'e\-pen\-dan\-ce \agq
du \sys \,$(f_1,\alb f_2,\alb \dots,\alb f_n)$\, et garantit
l'uni\-cit\'e de l'expression rationnelle, dans ce \sys
fondamental, de toute fraction rationnelle \sym sur
\,$\K$.

\ss Les relations dites
de Newton permettent  d'exprimer les sommes
de Newton dans le \sys fondamental
des \pols \syms \elrsz.

\begin{prop}\label{Newt}
  \emph{(Relations de Newton)}\,
  Les \pols de Newton \`a \,$n$\,
\idtrs \,$(S_k)_{k\in \NN}$\, sont
re\-li\'es aux \,$n$\, \pols \syms
\elrs \,$(\sigma_k)_{1\leq k \leq n}$\,
par les relations suivantes:
\index{Newton!relation de}\index{relation de Newton}
\begin{itemize}
\item[\emph{(i)}] $S_0=n$;
\item[\emph{(ii)}] pour $1\leq k\leq n~:~~~S_k +
\sum_{i=1}^{k-1}(-1)^{i} \,\sigma_i\,S_{k-i} +
(-1)^k\,k\,\sigma_k=0$;
\item[\emph{(iii)}] pour $k>n~:~~~S_k +
\sum_{i=1}^{n}(-1)^{i}\, \sigma_{i}\,S_{k-i}=0$.
\end{itemize}
\end{prop}
\prv On pose \,$\alpha_k =(-1)^{k+1}\,\sigma_k$\,
($k=1\ldots ,n$). Les \pols \syms
\elrs apparaissent dans le \dev
du \pol
$$
Q(X)\,=\,\prod_{i=1}^n\,(1-x_iX) \,=\,
1-\,\sum_{k=1}^{n}\alpha_k\,X^{k}\,.
$$
Par d\'erivation \logq formelle, on obtient
les \egts suivantes dans l'\agr de
s\'eries formelles \,$\FA(x_1,\dots,x_n)\,[[X,{X}^{-1}]]\,$:
$$
\frac{Q'(X)}{Q(X)}=\sum_{i=1}^n\,\frac{-x_i}{1-x_iX}=
\frac{1}{X}\,\sum_{i=1}^n
\left[\,1-\frac{1}{1-x_iX}\,\right]=
-\frac{1}{X}\,\sum_{k=1}^\infty S_k\,X^k
$$
ou encore dans
\,$\FA(x_1,\dots,x_n)\,[[X]]\,$:
\begin{equation}\label{Newt1}
- X\,Q'(X)\,=\,
Q(X)\,\,\sum_{k=1}^\infty S_k\,X^k
\end{equation}
avec \,$Q(X)= 1 -\left[\alpha_1\,X +
\dots  + \alpha_n\,X^n\right] $\, et
\,$-XQ'(X)\,=\,\alpha_1\,X+2\,\alpha_2\,X^2+\dots +
n\,\alpha_{n}\,X^{n}$.
Identifiant dans l'\'equation (\ref{Newt1}) les termes
en \,$X^{k}$\, pour \,$1\leq k\leq n$, on obtient les
formules (ii). Les formules (iii) sont obtenues par
identification des termes de degr\'e sup\'erieur \`a \,$n$.

\ms \rem En notant \,$\alpha_i =(-1)^{i+1}\,\sigma_i$\,
les relations de Newton s'\'ecrivent sous la forme matricielle suivante:
\begin{equation}\label{Newt3}
\left[\begin{array}{ccccc}
1      & 0      & \dots  & \dots & 0      \\
S_1    & 2      & \ddots &       & \vdots \\
\vdots & \ddots & \ddots & \ddots& \vdots \\
S_{n-2}&        & \ddots & n-1   & 0      \\ [1mm]
S_{n-1}& S_{n-2}&    & S_{1} & n  \\[1mm]
S_{n}& S_{n-1}&    &   & S_{1}  \\
\vdots &   &   &   & \vdots \\
S_{n+k}&  &    &   & S_{k+1}\\
\vdots &  \vdots & \vdots   & \vdots & \vdots
\end{array}\right]
\begin{array}{c}
     \left[\begin{array}{c} \alpha_1 \\[1,5mm]
     \alpha_2 \\ [1mm]
     \vdots \\ 
     \vdots \\ [1mm]
     \alpha_n
     \end{array}\right]\\
\begin{array}{c} ~ \\ ~ \\ ~ \\ ~ \\[1mm] ~
\end{array}
\end{array}
=\left[\begin{array}{c}
    S_1 \\[1,5mm]
    S_2 \\[1mm]
    \vdots \\ \vdots \\[1,5mm]
    S_n \\[1mm]
    S_{n+1} \\ \vdots \\ S_{n+k+1}\\ \vdots
\end{array}\right]
\end{equation}
Les relations (ii) (qui correspondent aux \,$n$\,
premi\`eres lignes
dans la matrice infinie ci-dessus) et (iii)
(qui correspondent aux lignes \num$n$\, et suivantes)
donnent la m\^eme
formule si l'on fait \,$k=n$.
D'autre part, les relations (iii) peuvent \^etre obtenues directement:
en multipliant par \,$X^{k-n}~(k>n)~$ le \pol
$$
P(X)=\prod_{i=1}^n\,(X-x_i)=X^n\,-\,\sum_{i=1}^{n}\,\alpha_i\,X^{n-i}~,
$$
on obtient
$$
X^{k-n}\,P(X)=X^{k}\,-\,\sum\nolimits_{i=1}^{n}\,\alpha_i\,X^{k-i}\,.
$$
Lorsque \,$k>n$\, les \idts
\,$x_j^{k-n}\,P(x_j)=0$\, donnent alors par sommation:
$$
\sum_{j=1}^{n}\,x_j^{k}\,-\,\sum_{i=1}^{n}\,
\alpha_i\,\left[\sum\nolimits_{j=1}^{n}x_j^{k-i}\right]
=S_k -  \sum_{i=1}^{n}\,\alpha_i\,S_{k-i}=0\,.
$$

\begin{cor}
Si \,$\A$\, est un \acom
o\`u les entiers \,$1,$ $2,$ $\ldots,$ $n$\, sont inversibles,
alors les sommes de Newton
\,$(S_k)_{1\leq k\leq n}$\, en les \,$n$\,
\idtrs \,$x_1,\dots,x_n$\, forment
un \sys fondamental de \pols sym\'etri\-ques
en \,$x_1,\dots,x_n$\, sur l'anneau \,$\A$.
\end{cor}
\prv  Le \sys \trg form\'e par les \,$n$\, premi\`eres
\'equations dans (\ref{Newt3}) admet clairement une solution unique
en les \,$\alpha_i$.
\qed

\begin{cor}\label{corSNT}
Soit \,$P=X^n-\left[a_1X^{n-1}+\dots+a_{n-1}X+a_n\right]$\,
un \pol unitaire \`a une \idtr sur un
anneau int\`egre \,$\A$, et soient
\,$\lambda_1,\alb \lambda_2,\alb\dots,\alb\lambda_n$\, les \,$n$\,
racines de ce \pol (distinctes ou non) dans un
corps de \deco de $P$\, (une extension de
$\FA$). Si l'on pose \,$s_k=\sum_{i=1}^n\lambda_i^k$\,
pour \,$1\leq k\leq n$, on a les relations:
\begin{equation}\label{Newt2}
\left\{\begin{array}{rcl} s_1&=&a_1 \\
s_2&=&s_1\,a_1+2\,a_2 \\ \vdots&\vdots&~~~~~~~\vdots \\
s_n&=&s_{n-1} a_{1}+\dots+ s_1 a_{n-1}+n\,a_n
\end{array}\right.
\end{equation}
\end{cor}

Ainsi les  \coes du \pol \,$P$\,
sont d\'etermin\'es de mani\`ere unique dans
\,$\A$\, par la donn\'ee de ses sommes de Newton
si  \,$n\,!$\, est non diviseur de
z\'ero dans \,$\A$.
Et si la division exacte par
chacun des entiers \,$2,\ldots,n$\, est explicite
lorsqu'elle est possible, le calcul des \,$s_i$\, \`a partir
des \,$a_i$\, est lui aussi explicite.

\ss Dans la \dfn suivante, on \gnlz ise les
sommes de Newton pour les z\'eros du \pol \,$P$\,
au cas d'un \acom arbitraire.
\begin{defi}  \label{SNe}
Soit
\,$P=X^n -\left[a_1X^{n-1}+\dots+a_{n-1}X+a_n\right]$\, un \pol
unitaire \`a une \idtr sur un \acom \,$\A$. Alors
les \elts \,$s_i$\, donn\'es par les \'equations
(\ref{Newt2})  s'appellent les sommes de Newton de \,$P$\,
(ou de tout \pol \,$u\,P$\, avec \,$u$\, non diviseur de z\'ero
  dans $\A$).
\end{defi}

Un fait important est le suivant.

\begin{lem}\label{SNTr}
Soit \,$\A$\, un \acomaz,
\,$A$\, une \maca d'ordre \,$n$\, sur \,$\A$\,
et \,$\rPA$\, son \polcarz.
Les sommes de Newton de \,$\rPA$\, donn\'ees dans
la \dfn pr\'ec\'edente sont les traces
des puissances de la matrice \,$A$:
$s_k=\Tr(A^k)$ pour $1\leq k\leq n$.
\end{lem}

\prv On remarque que les \egts
\`a d\'emontrer sont des \idas en
les entr\'ees de la matrice \,$A$\, (consid\'er\'ees
comme des \idtrsz). Il suffit donc de traiter
le cas d'un anneau int\`egre.
Dans ce cas, le r\'esultat est donn\'e par
le corollaire \ref{corSNT} et le lemme \ref{trac}.
\qed


\section {In\'egalit\'e de Hadamard et calcul modulaire}
\label{secHadamard}

Nous d\'ecrivons d'abord quelques majorations
utiles en \agr \linz.

\subsection{Normes matricielles}

Si \,$A=(a_{ij})\in\A^{m\times n}$\, est une matrice \`a \coes r\'eels
ou complexes,
on d\'efinit classiquement les \emph{normes} suivantes
(\cf \cite{Cia}, \cite{Gol})
$$
\norme{A}_1\,=\max_{1\leq j\leq
n}\left[\sum\nolimits_{i=1}^m|a_{ij}|\right]\,, \quad
\norme{ A}_\infty\,=\max_{1\leq i\leq m}
\left[\sum\nolimits_{j=1}^n|a_{ij}|\right]
$$
$$
\norme{A}_F\,=
\sqrt{\,
\sum\nolimits_{i=1}^m\sum\nolimits_{j=1}^n|a_{ij}|^2\,}
$$
Chacune de ces normes v\'erifie les relations classiques
$$
\norme{c \, A}\,=|c|\norme{A}\,,
\quad \norme{ A+B } \,\,\leq\,\, \norme{A} +\norme{B}
$$
(si \,$A$\, et \,$B$\, ont m\^{e}mes dimensions) et
$$
\norme{A\,B}\,\,\leq\,\, \norme{A} \,\norme{B}
$$
(si le produit \,$A\,B$\, est d\'efini).

Consid\'erons maintenant des matrices \`a \coes entiers.
La \emph{taille d'un entier} \,$x$\, est l'espace qu'il occupe lorsqu'on
l'implante sur machine. Si le codage des entiers est standard,
cela veut dire que la taille de \,$x$\, est correctement appr\'eci\'ee
par \,$1+\esup{\log_{2}(1+|x|)}$.

\begin{notation}
\label{notaLog}
{\rm  Dans tout cet ouvrage \,$\log{x}$\,
\indexnota{log@$\log$} d\'enote
\,$\max(1,\log_{2}{|x|})$.
}
\end{notation}

Lorsque \,$x$\, est entier, cela repr\'esente
donc la taille de \,$x$\, \`a une constante pr\`es.
Si \,$\lambda(A) = \log(\norme{ A })$\, avec l'une des
normes pr\'ec\'edentes, la taille de chaque \coe de \,$A$\, est
clairement major\'ee par \,$\lambda(A)$\,
(\`a une constante additive pr\`es)
et en outre les relations pr\'ec\'edentes impliquent
\immt que
$$\lambda(A\,B)\leq \lambda(A)+\lambda(B),\quad
\lambda(A+B)\leq \max(\lambda(A),\lambda(B))+1.$$

Ces relations sont souvent utiles pour calculer des majorations
de la taille des entiers qui interviennent comme r\'esultats de calculs
matriciels.

\subsection*{L'in\'egalit\'e de Hadamard}

L'in\'egalit\'e de
Hadamard
\index{Hadamard!in\'egalit\'e de}
s'applique aux matrices \`a \coes r\'eels.
La valeur absolue du \deter repr\'esente le
volume ($n\,$-\,dimensionnel) du \parap construit
sur les vecteurs colonnes de la matrice, \gui{et donc}
elle est major\'ee par le produit des longueurs de ces vecteurs:
\begin{equation} \label{borneHadamard1}
\left\vert\,\det{((a_{ij})_{1\leq i,j\leq n})}\right\vert\leq
\prod_{j=1}^n{\sqrt{\sum_{i=1}^n{a_{ij}^2}}}
\end{equation}

Il y a \'evidemment des preuves rigoureuses de
ce fait intuitif. Par exemple le processus d'orthogonalisation
de Gram-Schmidt remplace la matrice par une matrice de m\^eme
\deter dont les vecteurs-colonnes sont deux \`a deux
orthogonaux et ne sont pas plus longs que ceux de la
matrice initiale. La signification g\'eom\'etrique de cette
preuve est la suivante: le processus d'orthogonalisation
de Gram-Schmidt remplace le \parap (construit
sur les vecteurs colonnes) par un
\parap droit
de m\^{e}me volume dont les cot\'es sont devenus plus courts.
Ce m\^eme raisonnement donne l'\ine dans le cas d'une
matrice \`a \coes complexes en rempla\c{c}ant \,$a_{ij}^2$\, par
\,$\left|a_{ij}\right|^2$\, (mais l'interpr\'etation g\'eom\'etrique
directe dispara\^{\i}t).

Avec les normes \,$\norme{~~~}_1$\, et \,$\norme{~~~}_\infty$\, on
obtient pour une \macaz:
\begin{equation} \label{borneHadamard2}
     \left\vert\,\det{(A)}\right\vert\leq
\left(\norme{ A }_1\right)^n\,,\quad
\left\vert\,\det{(A)}\right\vert\leq \left(\norme{ A
}_\infty\right)^n
\end{equation}

Avec la norme de Frobenius $\norme{ A }_F$
on obtient la majoration suivante (un produit de \,$n$\,
r\'eels positifs dont la somme est constante est maximum
lorsqu'ils sont tous \'egaux):
\begin{equation} \label{borneHadamard3}
\left\vert\,\det{(A)}\right\vert\leq \left(\frac{\norme{ A
}_F}{\sqrt{n}}\right)^n
\end{equation}

\subsection[Th\'eor\`eme chinois et applications]{Le th\'eor\`eme
chinois et son application aux calculs modulaires}
\index{Th\'eor\`eme chinois} \label{Calcmod}
Soient \,$p_{1},p_{2},\ldots,p_{r}$\, des entiers positifs deux \`a
deux premiers entre eux. On pose \,$M=p_{1}p_{2}\cdots p_{r}$.
Pour toute suite \,$x_{1},x_{2},\ldots,x_{r}$\, de \,$r$\, entiers
relatifs, il existe un entier \,$x$\, (unique modulo \,$M$)
v\'erifiant: $x\equiv x_{i} ~[{\rm mod}~p_{i}]~(i=1,\ldots, r)$.
On peut calculer cet entier (modulo $M$) en remarquant que, pour
tout \,$i$\, compris entre 1 et \,$r$, les nombres \,$p_{i}$\, et
\,$q_{i}=M/p_{i}$\, sont premiers entre eux et que par cons\'equent
il existe des entiers \,$u_{i}$\, et \,$v_{i}$\, (relation de B\'ezout)
tels que
\,$p_{i}u_{i}+q_{i}v_{i}=1$ ($i=1,\ldots ,r$). Le nombre \,$x$\,
recherch\'e
n'est autre que \,$\sum_{i=1}^rx_{i}q_{i}v_{i}$ (modulo \,$M$). Il
est facile de v\'erifier qu'il r\'epond bien \`a la question.

\ss L'une des cons\'equences importantes du Th\'eor\`eme chinois en
calcul formel est son utilisation pour le calcul de
\coes entiers \,$x\in \ZZ$\, dont on sait majorer la valeur
absolue par un entier \,$B$\, strictement positif. Il
arrive souvent que les calculs \itmdsz, lorsqu'ils sont
effectu\'es (avec ou sans division) dans \,$\ZZ$, donnent des
\coes dont la taille explose rapidement, ce qui risque de
rendre ces calculs impraticables ou trop co\^{u}teux, alors
que la taille du r\'esultat final est bien plus petite.
Supposons que l'on ait \`a calculer un \,$x\in \ZZ$\, tel que
\,$-B\leq x\leq B$\, par un \algo sans divisions.

On commence par choisir
des entiers positifs $p_{1},\,p_{2},\ldots,p_{r}$ deux \`a deux
premiers entre eux dont le produit d\'epasse strictement $2B$.
Au lieu de calculer directement \,$x$,
on effectue tous les calculs modulo
\,$p_{i}$\, s\'epar\'ement pour chaque $i$ ($i=1, \ldots, r$). Les
r\'esultats \,$x_{1},\,x_{2},\ldots,\, x_{r}$\, ainsi obtenus sont tels
que $x$ est dans la classe de $x_{i}$ modulo \,$p_{i}$\, (pour $i=1
\ldots r$). Utilisant les m{\^e}mes notations que ci-dessus pour
les coefficients de B\'ezout relatifs aux couples
\,$(p_{i},q_{i})$, on r\'ecup\`ere ensuite le r\'esultat principal
\,$x$\, \`a partir des r\'esultats partiels \,$x_{1},\,x_{2},\ldots,\,
x_{r}$\, en remarquant que \,$x$\, est l'entier relatif de plus petite
valeur absolue congru \`a \,$\sum_{i=1}^rx_{i}q_{i}v_{i}$\, modulo
\,$p_{1}p_{2}\cdots p_{r}$\, (puisque $-B\leq x\leq B$). Dans le cas
d'un \algo avec divisions, les facteurs $p_{i}$ doivent
\^etre choisis de mani\`ere \`a ce qu'ils soient premiers avec
les diviseurs intervenant dans les calculs.

Pour le calcul des d\'eterminants de matrices
\`a \coes entiers, par exemple, on peut utiliser
l'\inH
(\ref{borneHadamard1}) pour faire fonctionner la \met
modulaire. On prendra pour borne \,$B=M^nn^{n/2}$\, o\`u
\,$M=\max_{1\leq i,j\leq n}{|a_{ij}|}$. Si cela s'av\`ere
pr\'ef\'erable, on peut choisir une des bornes donn\'ees dans les
\'equations (\ref{borneHadamard2}) et (\ref{borneHadamard3}).

Il en est de m{\^e}me pour le calcul du \polcar $\rPA$
d'une matrice \,$A\in \ZZ^{n\times n}$\, car chacun des \coes de
$\rPA$ est une somme de mineurs diagonaux de la matrice
\,$A$. On peut donc l\`a aussi utiliser l'\inH pour majorer les valeurs
absolues des \coes en vue du
traitement modulaire.

Plus pr\'ecis\'ement, si l'on prend
\,$M=\max_{1\leq i,j\leq n}{|a_{ij}|}$\, comme ci-dessus et si l'on
d\'esigne par \,$m_{k,j}~(1\leq j\leq C_n^k)$\, les mineurs
de \,$A$\, diagonaux d'ordre \,$k$\, (pour $k$ donn\'e entre $1$ et
\,$n$),
alors le \coe \,$\mu_{k}$\, du terme de degr\'e \,$n-k$\, de $\rPA$
est major\'e en valeur absolue comme suit :
\begin{equation} \label{borneHadamard4}
|\mu_{k}|=\left\vert\;\sum_{j\in\{1,\ldots,C_n^k\}}\!\!\!
m_{k,j}\;\right\vert
\leq C_n^k M^k k^{k\over 2} \leq {(2M)}^nn^{n\over 2}
\end{equation}
puisque \,$C_n^k\leq 2^n$.

\subsection*{Quelques consid\'erations pratiques}

L'id\'ee principale dans l'utilisation du calcul modulaire
est de remplacer un \algo dans \,$\ZZ$\, permettant de
r\'esoudre un \pb donn\'e par plusieurs \algos modulo
des nombres premiers.

Pour \^{e}tre vraiment efficace,
cette \met doit \^{e}tre appliqu\'ee
avec des listes  de nombres premiers
\,$p_{1},\,p_{2},\ldots,\,p_{r}$\,
d\'ej\`a r\'epertori\'ees et pour lesquelles on a
d\'ej\`a calcul\'e les \coes \,$q_{i}v_{i}$\, correspondants
qui permettent de r\'ecup\'erer \,$x$\, \`a partir des \,$x_{i}$.

Ces nombres premiers peuvent \^etre choisis par rapport
\`a la taille des mots trait\'es par les processeurs.
Par exemple, pour des processeurs
qui traitent des mots \`a 64 bits, on prend des nombres premiers
compris entre $2^{63}$ et $2^{64}-1$: il y en a suffisamment
(bien plus que $10^{17}$ nombres!) pour r\'esoudre dans la
pratique tous les \pbs de taille humainement raisonnable
et r\'ealiste~\cite{VonZurbook}.
En outre on poss\`ede des tests rapides pour savoir si un
nombre est premier, et cela a permis d'\'etablir des listes
\,$p_{1},p_{2},\ldots,p_{r}$\, avec la liste des \coes \,$q_{i}v_{i}$\,
correspondants, qui r\'epondent \`a tous les cas qui se posent
en pratique.

Chaque \opari \elr
modulo un tel nombre premier se fait alors en
temps constant, ce qui r\'eduit consid\'erablement le temps de
calcul.
En outre la \deco du \pb en \algos modulaires
offre la possibilit\'e d'utiliser plusieurs processeurs en \paralz.

\section{R\'esolution uniforme des \slis}
\label{secUnifSli}

Nous expliquons ici comment le \polcar permet de d\'eterminer le rang
d'une matrice et de r\'esoudre
\emph{uniform\'ement} (avec une seule formule, du type Cramer) les \slis
ayant un format donn\'e et un rang fix\'e. Et ceci sur un corps
arbitraire.

Cette solution uniforme (\cf \cite{DiGL}) constitue une extension d'un
r\'esultat de Mulmuley
\cite{Mul} qui ne traite que la question du rang.

Naturellement, le rang d'une matrice peut \^{e}tre calcul\'e par la
\mpgz. Mais la \met n'est pas uniforme et, a priori, ne se laisse pas
bien \pararz.

Les applications des formules et \algos que nous allons d\'ecrire ici
seront de deux ordres: d'une part en calcul \paralz, d'autre part
lorsqu'on doit traiter des \slis d\'ependant de \paratsz.

Dans ce deuxi\`eme cas de figure, la \mpg produit un arbre de calcul qui
risque de comporter un tr\`es grand nombre de branches, correspondant
\`a un grand nombre de formules distinctes, lorsque les \parats prennent
toutes les valeurs possibles.
Le cas extr\^{e}me est celui o\`u toutes les entr\'ees d'une matrice
sont des \parats
ind\'e\-pen\-dants. Par exemple avec une matrice de rang maximum de
format \,$n\times 2n$\, la solution du \sli correspondant par la \mpg
d\'epend du mineur maximal non nul qu'on extrait, et ce dernier peut
\^{e}tre n'importe lequel des \,${C_{2n}^n}>2^n$\, mineurs d'ordre
\,$n$\, de la matrice.

En analyse num\'erique matricielle, avec des matrices \`a \coes r\'eels
ou complexes, une formule uniforme compacte en rang fix\'e est obtenue
par l'utilisation des \coes de Gram de la matrice correspondant au \sli
homog\`ene: dans le cas r\'eel, le Gram d'ordre \,$k$\, d'une matrice
\,$A$\, est \'egal \`a la somme des carr\'es de tous les mineurs d'ordre
\,$k$\, de \,$A$, son annulation signifie que le rang de la matrice
n'exc\`ede pas \,$k-1$.

Les \idts que nous allons obtenir sont des \gnns directes des formules
usuelles qui expriment l'\iMP en fonction des \coes de Gram de la
matrice. L'\'etonnant est que, m\^{e}me sur un corps fini, un petit
nombre de sommes de carr\'es de mineurs suffit \`a contr\^{o}ler le rang
d'une matrice, et que des formules semblables aux formules usuelles
fonctionnent encore.

Il y a cependant un prix \`a payer, non n\'egligeable, qui est
d'introduire un \parat \supt dans les calculs.

\subsection[L'\iMP]{Les \coes de Gram et l'\iMP dans le cas r\'eel ou
complexe}
\label{subsecIMP}
\subsubsection*{Th\'eorie g\'en\'erale}
Dans toute la section \ref{subsecIMP} \,$A$\, est une matrice dans
\,$\KK^{m \times n}$, avec  \,$\KK=\CC\;\mathrm{ou}\;\RR$,
repr\'esentant sur des bases orthonorm\'ees une \ali \,$\varphi
:E\rightarrow F$\, entre espaces vectoriels hermitiens ou euclidiens de
dimension finie. Nous noterons
\,$\gen{x,y}$\, le produit scalaire des vecteurs \,$x$\, et \,$y$.
Nous notons \,$A\sta$\, la transpos\'ee de la conjugu\'ee de \,$A$\,
(dans le cas r\'eel on a \,$A\sta=\tra{\!A}$).
  La matrice \,$A\sta$\, repr\'esente sur les m\^{e}mes bases l'\ali
adjointe{\footnote{~A ne pas confondre avec la matrice adjointe
\,$\Adj(A)$. Cette ambig\"uit\'e dans la terminologie, en fran\c{c}ais,
est ennuyeuse.}}  \,$\varphi\sta$, caract\'eris\'ee par:
\begin{equation} \label{eqApAdj}
\forall x\in E\quad \forall y\in F\quad \quad \gen{\varphi(x),y
}_F=\gen{x,\varphi\sta(y)}_E
\end{equation}
Les matrices \,$AA\sta$\, et  \,$A\sta A$\,  sont des \macas
hermitiennes positives (sym\'etriques r\'eelles positives dans le cas
r\'eel), en \gnl non r\'eguli\`eres.
Si \,$H$\, est un sous-\evc de \,$E$\, nous noterons
\,$\pi_H$\, la \pro \orte de \,$E$\, sur \,$H$, vue comme
\ali de \,$E$\, dans \,$E$.

D'un point de vue de pure \agr \lin tous les r\'esultats de la
\gui{th\'eorie \gnlez} qui suit sont bas\'es sur la \deco des espaces
\,$E$\, et  \,$F$\, en sommes directes de noyaux et d'images de
\,$\varphi$\, et  \,$\varphi\sta$.
\begin{lemma}
\label{lemSDO}
Nous avons deux sommes directes:
\begin{equation} \label{eqSDO}
   \Im\,\varphi \oplus \Ker\,\varphi\sta= F, \; \; \Ker\,\varphi\oplus
\Im\,\varphi\sta= E
\end{equation}
\end{lemma}

Cela r\'esulte du fait que
\,$\Ker\,\varphi\sta$\, (resp. \,$\Ker\,\varphi$) est le
sous-espace \ort de
\,$\Im\,\varphi$\,  (resp. \,$\Im\,\varphi\sta$), ce qui est une
cons\'equence directe de l'\egtz~(\ref{eqApAdj}).

Nous en d\'eduisons les faits suivants.

\begin{fact}
\label{factSDO}~
\begin{enumerate}
\item \label{factSDO.2}  L'\ali \,$\varphi$\, se restreint en un \iso
\,$\varphi_0$\, de \,$\Im\,\varphi\sta$\, sur \,$\Im\,\varphi$\, et
\,$\varphi\sta$\, se restreint en un \iso \,$\varphi_0\sta$\, de
\,$\Im\,\varphi$\, sur \,$\Im\,\varphi\sta$.
\item \label{factSDO.3} En outre:
\begin{equation} \label{SDO4}
\begin{array}{rclrrcl}
\Im\,\varphi  & =  &\Im\,\varphi\varphi\sta\,,   & \quad&
\Ker\,\varphi\sta  & = & \Ker\,\varphi\varphi\sta\,, \\[1mm]
\Ker\,\varphi   & =  & \Ker\,\varphi\sta\varphi\,,  &
&\Im\,\varphi\sta & = & \Im\,\varphi\sta\varphi\,.
\end{array}
\end{equation}
\item \label{factSDO.4} Soit \,$\varphi_1:\Im\,\varphi\rightarrow
\Im\,\varphi$\,
l'\auto \lin d\'efini par
\,$\varphi_1=\varphi_0\varphi_0\sta$. C'est la restriction de
\,$\varphi\,\varphi\sta$\, \`a
\,$\Im\,\varphi$.
Nous avons:
$$ \formule{\det{(\Id_{\Im\,\varphi} + Z \,\varphi_1)} & = &
\det{(\Id_F + Z \,\varphi \varphi\sta)} \\ & = & 1+ a_1 Z+ \cdots
+ a_r Z^r\,.} $$
o\`u \,$r=\rg(\varphi)=\rg(\varphi_1)$\, et  \,$a_r\neq 0$. De la
m\^eme fa\c{c}on, nous avons l'\auto
\,$\varphi_1\sta=\varphi_0\sta\,\varphi_0$\, de
\,$\Im\,\varphi\sta$\, et $$ \formule{\det{(\Id_{\Im\,\varphi\sta}
+ Z \,\varphi_1\sta)} & = & \det{(\Id_E + Z \,\varphi\sta
\varphi)}\\ & = & 1+ a_1 Z+ \cdots+ a_r Z^r\,.} $$
\end{enumerate}
\end{fact}

Ce sont des cons\'equences directes du lemme \ref{lemSDO}.

Peut-\^etre cela sera plus clair si nous
repr\'esentons \,$\varphi\sta$\,  et
\,$\varphi$\, dans les sommes \ortes (\ref{eqSDO}):
\begin{equation} \label{SDO2}
\varphi=
\pmatrix{
\varphi_0    &  0_{K,I}     \cr
  0_{I\sta,K\sta}   &
0_{K,K\sta}   },
\quad
\varphi\sta=
\pmatrix{
\varphi_0\sta    &  0_{K\sta,I\sta}     \cr
  0_{I,K}   &  0_{K\sta,K}
}
\end{equation}
o\`u \,$K=\Ker\,\varphi$, \,$K\sta=\Ker\,\varphi\sta$,
\,$I=\Im\,\varphi$, \,$I\sta=\Im\,\varphi\sta$.
\begin{definition}\label{defGk}
  Les {\em \coes de Gram
de \,$A$\, (ou de \,$\varphi$)} \index{coefficient de Gram!d'une
matrice}
\index{Gram!coefficient de}  sont les  \,$\cG_k(A)=\cG_k(\varphi)=a_k$\,
\indexnota{Gk(A)@$\cG_k(A)$} donn\'es par la formule
\begin{equation} \label{eqDefGram}
\det{(\I_m \, + \, Z \, AA\sta)} \; = \; 1 \; + \; a_1 \; Z \; +
\; \cdots \; + \; a_m \; Z^m.
\end{equation}
Nous d\'efinissons aussi
\,$\cG_0(A)=1{\quad\rm et\quad}\cG_\ell(A)=0
{\quad\rm  pour\quad} \ell>m.$
\end{definition}

Notez que le \polcar de \,$B=AA\sta$\, est \'egal \`a
\,$(-1)^mZ^mQ(1/Z)$\, o\`u \,$Q(Z)=\det{(\I_m+Z\,B)}$. Les \coes
de Gram de \,$\varphi$\, sont donc, au signe pr\`es, les \coes du
\polcar de \,$\varphi\varphi\sta$.
\begin{lem}
\label{lemGk} \emph{(Conditions de Gram pour le rang)}
\begin{enumerate}
\item L'\ali \,$\varphi$\, est de rang \,$\leq r$\, \ssi
\,$\cG_{k}(\varphi)=0$\, pour  \,$r<k\leq n$. Elle est de rang \,$r$\,
si en outre \,$\cG_{r}(\varphi)\neq 0.$
\item Le \coe de Gram \,$\cG_k(A)$\, est un nombre r\'eel positif ou
nul, \'egal \`a la somme des carr\'es des modules des mineurs d'ordre
\,$k$\, de la matrice \,$A$. En cons\'equence,
\,$\cG_{r+1}(\varphi)=0$\, suffit pour certifier que le rang est \,$\leq
r$.
\end{enumerate}
\end{lem}
\preuve
Le premier point est une cons\'equence directe du fait
\ref{factSDO}-\ref{factSDO.4}.

\noi Il pourrait aussi \^etre vu comme une cons\'equence du second
point, que nous d\'emontrons maintenant.

\noi Le \coe \,$a_k$\, est la somme des \mips d'ordre \,$k$\, de
\,$AA\sta$. Chaque \mip d'ordre \,$k$\, est obtenu comme \deter de la
matrice correspondante, qui est \'egale \`a \,$A_\alpha
(A_\alpha)\sta$\, o\`u \,$\alpha$\, d\'esigne  un \,$k$-uple
\,$\alpha_1< \cdots < \alpha_k$\, extrait de \,$\{1,\ldots ,m\}$\, et
\,$A_\alpha$\, est la matrice
extraite de \,$A$\, en gardant seulement les \,$k$\, lignes
correspondant \`a \,$\alpha$.
La formule de Binet-Cauchy (\ref{eqBC}) nous indique alors que ce
d\'eterminant est la somme des carr\'es des modules des mineurs d'ordre
\,$k$\, extraits de~$\,A_\alpha$.
\qed

\medskip Nous supposons d\'esormais
\,$r=\rg(\varphi)$\, (donc \,$a_r=\cG_r(\varphi)\neq 0$). Puisque
$$ \det{(\Id_{\Im\,\varphi} + Z \,\varphi_1)} \; = \; \det{(\Id_F
+ Z \,\varphi \varphi\sta)} \; = \; 1+ a_1 Z+ \cdots+ a_r Z^r $$
le th\'eor\`eme de Cayley-Hamilton nous donne
\begin{equation} \label{eqCHfi1}
  \varphi_1^r - a_1 \varphi_1^{r-1}+ \cdots+ (-1)^ra_r
\Id _{\Im\,\varphi } \,= \,0\,.
\end{equation}
Par suite, on obtient en rempla\c{c}ant
\,$\varphi_1$\, par \,$\varphi\varphi\sta$\, dans la formule
pr\'ec\'edente
$$  (\varphi\varphi\sta)^r - a_1 (\varphi\varphi\sta)^{r-1}+ \cdots+ (-
1)^ra_r
\pi_{\Im\,\varphi } \,= \,0\,.
$$
Ainsi:

\begin{lem}
\label{lem proj orth} {\rm  (\pros \ortes sur l'image et sur le noyau)}
\begin{enumerate}
\item La  \pro \orte \,$\pi_I$\, sur le
sous-espace \,$I=\Im\,\varphi \subseteq F$\,
est \'egale \`a:
\begin{equation} \label{eqProIm}
a_r^{-1}\left(a_{r-1}\varphi\, \varphi \sta-a_{r-2}(\varphi\,
\varphi \sta)^2+ \cdots
+(-1)^{r-1}(\varphi\, \varphi \sta)^r\right).
\end{equation}

\item La  \pro \orte \,$\pi_{I\sta}$\, sur le
sous-espace
\,$I\sta=\Im\,\varphi\sta \subseteq E$\,
est \'egale \`a:
\begin{equation} \label{eqProKer}
a_r^{-1}\left(a_{r-1}\varphi \sta \varphi -a_{r-2}(\varphi
\sta \varphi )^2+ \cdots
+(-1)^{r-1}(\varphi \sta \varphi )^r\right).
\end{equation}
Et la \pro \orte sur le noyau de \,$\varphi $\, est
$\Id_E-\pi_{I\sta }$.
\end{enumerate}
\end{lem}

En outre l'\'equation (\ref{eqCHfi1}) implique que l'inverse de
\,$\varphi_1$\, est donn\'e par
$$
\varphi_1^{-1} \; = \; a_r^{-1}\left(a_{r-1}{\rm  Id}_{\Im\,\varphi}
-a_{r-2} \varphi_1+ \cdots +(-1)^{r-1} \varphi_1^{r-1}\right) \,.
$$
De m\^{e}me  on a
$$
{\varphi_1\sta}^{-1} \; = \; a_r^{-1}\left(a_{r-1}{\rm
Id}_{\Im\,\varphi\sta}-a_{r-2} {\varphi_1\sta}+ \cdots + (-1)^{r-1}
{\varphi_1\sta}^{r-1}\right) \,.
$$
et puisque
  \,$\varphi_1\sta=\varphi_0\sta\,\varphi_0,$
cela donne
\begin{equation} \label{eqfi1sta}
\formule{\forall y \in \Im\,\varphi\sta \quad{\varphi_1\sta}^{-1} (y)\;
= \; \\[1mm]
\;\;\;a_r^{-1}\left( a_{r-1}{\rm Id}_{E}-a_{r-2} ({\varphi\sta\varphi})+
\cdots +
(-1)^{r-1} ({\varphi\sta\varphi})^{r-1}\right)(y) }
\end{equation}

\begin{definition}
\label{defIMP}
Supposons que \,$\varphi$\, est de rang \,$r$.  \emph{L'\iMP de
\,$\varphi$\, (en rang \,$r$)}
\index{inverse de Moore-Penrose!en rang \,$r$} est  l'\ali
\index{Moore-Penrose!inverse de}
\,$\varphi^{[-1]_r}:F\rightarrow E$\,  d\'efinie par:
$$\forall y \in F\qquad
\varphi^{[-1]_r}(y)= \varphi_0^{-1}( \pi_{\Im\,\varphi}(y)).$$
\end{definition}

\rem Nous n'avons pas \'ecrit \,$\varphi^{[-1]_r} = \varphi_0^{-1}\circ
\pi_{\Im\,\varphi}$\,
parce que le deuxi\`eme membre est a priori \gui{mal d\'efini}:
\,$\pi_{\Im\,\varphi}$\, est une application de
\,$F$\, dans \,$F$,  \,$\varphi_0^{-1}$\,
est une application de \,$\Im\,\varphi$\, dans \,$\Im\,\varphi\sta$\, et
\,$\varphi^{[-1]_r}$\,
est une application de \,$F$\, dans \,$E$.

\medskip D'apr\`es  (\ref{eqSDO})   et  (\ref{SDO2})  on voit que
$$\begin{array}{ccrcl}
\forall y \in F& \quad  & \pi_{\Im\,\varphi}(y)  & =  &
{\varphi_0\sta}^{-1}(\varphi\sta (y)) \\[1mm]
\forall y \in F& \quad   & \varphi^{[-1]_r}(y) &  = &
\varphi_0^{-1}({\varphi\sta_0}^{-1}(\varphi\sta(y))).
\end{array}$$
et puisque \,$\varphi_0^{-1}\circ {\varphi_0\sta}^{-
1}={\varphi_1\sta}^{-1}$
nous obtenons
\begin{equation} \label{SDO5}
\forall y \in F\quad \varphi^{[-1]_r}(y) =
{\varphi_1\sta}^{-1}(\varphi\sta(y))
\end{equation}

  En appliquant (\ref{eqfi1sta}) on obtient alors une formule uniforme en
rang \,$r$\, qui donne une solution  des \slis en analyse num\'erique
matricielle:

\begin{prop}
\label{propMoPe}
        \emph{(Inverse de Moore-Penrose)}
  Soit \,$v\in F$. Soit
\,$\varphi\oplus v$\, l'\ali \,$E\oplus \KK\rightarrow F$
      d\'efinie par \,$(\varphi\oplus v) (x,\lambda) =
     \varphi(x)+\lambda v$.

\sni  1) L'\iMP
\,$\varphi ^{[-1]_r}\in\L(F,E)$
est donn\'e par:
\begin{equation} \label{eqMoPe}
\formule{
a_r^{-1}\left(a_{r-1}\Id_E-a_{r-2}(\varphi \sta \varphi)  + \cdots+
(-1)^{r-1}(\varphi \sta \varphi )^{r-1}\right)  \varphi \sta \\
\quad =
\\
a_r^{-1}\, \varphi \sta \left(a_{r-1}\Id_F-a_{r-2}
(\varphi  \varphi\sta)  + \cdots+
(-1)^{r-1}(\varphi  \varphi\sta )^{r-1}\right)}
\end{equation}
o\`u \,$a_k=\cG_k(\varphi )$.

\sni 2) Nous avons \,$v\in \Im(\varphi)$
  \ssi \,$\cG_{r+1}(\varphi  \oplus v)=0$\,
  \ssi
\begin{equation} \label{eqMoPe2}
    v \; = \; \varphi  \; \varphi ^{[-1]_r} (v)
\end{equation}
Dans ce cas \,$x=\varphi ^{[-1]_r} (v)$\, est
l'unique solution dans  \,$\Im(\varphi \sta)$.
\end{prop}

\begin{remark}
\label{remMatrices}
Voici la formulation matricielle du lemme
  \ref{lem proj orth} et de la proposition
\ref{propMoPe}.
Soient
  \,$m,\, n,\, r>0$\, dans \,$\NN$\,
avec \,$r \le \min(m, n)$\, et
  \,$A\in\KK^{m \times n}$\, une matrice de rang \,$r$.
Posons  $a_k=\cG_k(A)$. Soit \,$V\in\KK^{m \times 1}$.

\sni 1) La matrice de la \pro \orte sur le
sous-espace
\,$\Im\,A \subseteq \KK^m$\,
est \'egale \`a
$$
P=a_r^{-1}\left(a_{r-1}AA\sta-a_{r-2}(AA\sta)^2+ \cdots
+(-1)^{r-1}(AA\sta)^r\right).
$$

\sni 2) La matrice de la \pro \orte  sur le sous-espace
\,$\Im\,A\sta \subseteq \KK^n$\,
est \'egale \`a:
$$
P\sta=a_r^{-1}\left(a_{r-1}A\sta A-a_{r-2}(A\sta A)^2+ \cdots
+(-1)^{r-1}(A\sta A)^r\right).
$$
Et celle de la \pro \orte sur le noyau de \,$A$\, est
\,$\I_n-P\sta$.

\sni 3) La matrice
  \,$A^{[-1]_r}\in\KK^{n \times m}$\, (\iMP de
\,$A$\, en rang \,$r$) est \'egale \`a:
\begin{equation} \label{eqMoPeM}
a_r^{-1}\left(a_{r-1}\I_n-a_{r-2}A\sta A + \cdots+
(-1)^{r-1}(A\sta A)^{r-1}\right) \; A\sta \,
\end{equation}

\sni 4) Le \sli \,$AX=V$\,
admet une solution \ssi \,$\cG_{r+1}(A \vert V)=0$\,
($(A\vert V)$\, est la matrice obtenue en juxtaposant la colonne \,$V$\,
\`a droite de la matrice\,$A$)  \ssi on a l'\'egalit\'e:
\begin{equation} \label{eqMoPe2M}
    V \; = \; A \; A^{[-1]_r} \; V
\end{equation}
Dans ce cas \,$X=A^{[-1]_r} \; V$\, est l'unique solution dans
l'espace \,$\Im\,A\sta$.
\end{remark}

Notez que la matrice \,$A^{[-1]_r}\in\KK^{n \times m}$\, est bien
d\'efinie par la formule (\ref{eqMoPeM}) d\`es que \,$A$\, est de rang
\,$\geq r$. Cela est utile en analyse num\'erique et de mani\`ere plus
\gnle chaque fois que les
\coes de \,$A$\, sont des r\'eels connus avec seulement une pr\'ecision
finie (ce qui peut introduire  une incertitude sur le rang de la
matrice).

\subsubsection*{Cas des matrices hermitiennes}
Lorsque \,$E=F$\,  et  \,$\varphi=\varphi\sta$, l'\endom \,$\varphi$\,
est dit hermitien.
Alors on a une \deco \orte
\,$E=\Ker\,\varphi\oplus\Im\,\varphi$\, et la restriction
\,$\varphi_0$\,
de \,$\varphi$\,  \`a
  \,$\Im\,\varphi$\, est un \auto \lin de  \,$\Im\,\varphi$.
Nous posons
\begin{equation} \label{eqDefGH}
\det{(\Id_E \, + \, Z \, \varphi )} \, = \det{(\I_m \, + \, Z \,
A)} \, = \, 1 \, + \, b_1 \, Z \, + \, \cdots \, + \, b_n \, Z^n.
\end{equation}
Au signe pr\`es, les \,$b_i$\, sont donc les \coes du
\polcar de \,$\varphi$.
Si le rang de  \,$\varphi$\, est \'egal \`a \,$r$\, alors  \,$b_r\neq
0$,
  \,$b_{r+1}=\ldots= b_n= 0$\, et
$$ \det{(\Id_E \, + \, Z \, \varphi )}\,=\,
\det{(\Id_{\Im\,\varphi} \, + \, Z \, \varphi_0 )}\, = \, 1 \, +
\, b_1 \, Z \, + \, \cdots \, + \, b_r \, Z^r. $$ Ainsi par
Cayley-Hamilton $$
\varphi_0^r-b_1\,\varphi_0^{r-1}+b_2\,\varphi_0^{r-2}+ \cdots+
(-1)^{r-1}b_{r-2}\,\varphi_0+(-1)^{r-1}b_r\,\Id_{\Im\,\varphi}=0
$$ Et en rempla\c{c}ant \,$\varphi_0$\, par  \,$\varphi$\, nous
obtenons: $$ \varphi^r-b_1\,\varphi^{r-1}+b_2\,\varphi^{r-2}+
\cdots+
(-1)^{r-1}b_{r-2}\,\varphi+(-1)^{r-1}b_r\,\pi_{\Im\,\varphi}=0 $$

Ceci donne, pour le cas des matrices hermitiennes, une version
simplifi\'ee des r\'esultats pr\'ec\'edents plus g\'en\'eraux.
Elle se trouve dans l'ouvrage \cite{BP} de Bini et Pan.
\begin{prop}\label{propMoPeH}
          \emph{(\iMPz,  cas hermitien)}

\sni 1) La \pro \orte
\,$\pi_\Im\,\varphi$\, sur le sous-espace \,$\Im\,\varphi$\, est \'egale
\`a:
\begin{equation} \label{eqProImH}
b_r^{-1}
\left(b_{r-1}\,\varphi-b_{r-2}\,\varphi^2+ \cdots+
(-1)^{r-1}b_1\,\varphi^{r-1}+ (-1)^{r} \varphi^r\right)
.
\end{equation}


\sni 2)  L'\iMP  $\varphi ^{[-1]_r}\in\L(E,E)$\, est \'egal \`a:
\begin{equation} \label{eqMoPeH}
b_r^{-1}\left(b_{r-1}\,\pi_{\Im\,\varphi} -b_{r-2}\,\varphi   + \cdots+
(-1)^{r-1}b_1\,\varphi^{r-2}+ (-1)^{r} \varphi^{r-1}\right)
\end{equation}
\end{prop}

Remarquez que l'\'equation (\ref{eqProIm}) peut \^etre d\'eduite de
(\ref{SDO4})
et  (\ref{eqProImH}).

\subsubsection*{Interpr\'etation g\'eom\'etrique}
Si \,$\lambda_1\geq \lambda_2\geq \cdots \geq \lambda_r$\, sont les
\emph{valeurs singuli\`eres} non nulles de \,$\varphi$, \cad les racines
carr\'ees des valeurs propres \,$>0$\, de \,$\varphi \varphi\sta$\, il
existe des bases orthonorm\'ees de \,$E$\, et \,$F$\, par rapport
auxquelles la matrice de \,$\varphi $\, est \'egale \`a \,$L$:
$$ L\;=\;\cmatrix{
  \lambda_1&   0         & \cdots &\cdots&0  &\cdots& 0  \cr
     0    &  \lambda_2   &   0    &     &        && \vdots  \cr
  \vdots  &  \ddots  &\ddots&\ddots&     &    &  \vdots \cr
  \vdots   &             &\ddots& \lambda_r&0 &   &   \vdots \cr
  0   &  \cdots    & \cdots & 0 &    0  &       &   0    \cr
  \vdots   &             &      &  &      &    &  \vdots \cr
  0   &  \cdots          &  & \cdots &      &     \cdots   &   0
}
$$
(\cf \cite{Cia,Gol,LM}).

Matriciellement on obtient \,$A=ULV$\, o\`u \,$U$\, et \,$V$\, sont des
matrices unitaires (\ortes dans la cas r\'eel) convenables. Ceci
s'appelle la d\'ecomposition de \,$A$\, en valeurs singuli\`eres (la SVD
en anglais).

On voit que \,$\varphi $\, transforme la sph\`ere unit\'e de \,$E$\, en
un ellipso\"{\i}de dans \,$\Im\,\varphi$\, avec pour longueurs des axes
principaux \,$2\lambda_1,\ldots ,2\lambda _r$. Dans ces conditions la
matrice de  \,$\varphi\sta $\, est \'egale \`a \,$L\sta=\tra{\!L}$\, et
celle de \,$\varphi^{[-1]_r}$\, est \'egale \`a
$$ L^{[-1]_r}\;=\;\carray{cccccccccc}{
  (\lambda_1)^{-1}&   0    & \cdots &   \cdots  &0  &\cdots& 0 \\
     0    &  (\lambda_2)^{-1}   &   0    &     &   && \vdots  \\[1mm]
  \vdots   &  \ddots\;\;&\ddots\;\;&\ddots&     &   &   \vdots \\[1mm]
  \vdots   &             &\ddots\;& (\lambda_r)^{-1}&0 &&\vdots \\[2mm]
  0   &  \cdots    & \cdots & 0 &    0  &        &   0 \cr
  \vdots   &             &      &  &      &        &   \vdots \cr
  0   &  \cdots          &  & \cdots &      &     \cdots   &   0
}
$$
de m\^{e}me format que \,$\tra{\!L}$.

Bien que les matrices \,$L$\, et \,$A^{[-1]_r}$\, soient attach\'ees de
mani\`ere
unique \`a \,$A$\, et d\'ependent continument de \,$A$\, (sous
l'hypoth\`ese que le rang est fix\'e), il n'en va pas de m\^{e}me pour
les matrices \,$U$\, et \,$V$\, de la d\'ecomposition en valeurs
singuli\`eres, qui sont fondamentalement instables.

Que le vecteur \,$v$\, appartienne ou non \`a \,$\Im\,\varphi$\, on a
toujours \,$\varphi^{[-1]_r}(v)\in \Im\,\varphi\sta  $, qui est le
sous-espace \ort \`a \,$\Ker\,\varphi $\, et \,$\varphi (\varphi^{[-
1]_r}(v))$\, est la \pro \orte de \,$v$\, sur
\,$\Im\,\varphi $. Ainsi lorsque \,$v$\, n'est pas dans l'image, l'\iMP
fournit une \emph{solution approch\'ee} \,$x=\varphi^{[-1]_r}(v)$\, qui
donne pour \,$\varphi(x)$\, la meilleure approximation (au sens des
moindres carr\'es) de \,$v$. En outre \,$x$\, est la plus petite en
norme (parmi les solutions qui r\'ealisent cette meilleure
approximation).

Ce qui est remarquable est qu'on arrive \`a calculer
(essentiellement \`a l'aide du \polcar de \,$\varphi \varphi\sta$) les
\pros \ortes et l'\iMP par une formule uniforme (plus exactement, par
une formule qui  ne d\'epend que du rang, lequel se lit sur le \polcar
en question) sans qu'on ait besoin de calculer les bases orthonorm\'ees
dans lesquelles se r\'ev\`ele la g\'eom\'etrie de l'\ali \,$\varphi$.

\subsection{G\'en\'eralisation sur un corps arbitraire}
\label{subsecIMPG}
\subsubsection*{Th\'eorie g\'en\'erale}
Dans le  cas r\'eel ou complexe, on a vu qu'en termes d'\agr \lin tout
le paragraphe \gui{th\'eorie \gnlez} est gouvern\'e par les sommes
directes (\ref{eqSDO}) entre les noyaux et images de \,$\varphi $\, et
\,$\varphi\sta$\,
(lemme \ref{lemSDO}):
$$
   \Im\,\varphi\oplus \Ker\,\varphi\sta= F,
\; \; \Ker\,\varphi\oplus \Im\,\varphi\sta= E\,.
$$
Il suffit en effet, lorsqu'on parle de \pro \ortez, de remplacer par
exemple l'expression \emph{\pro \orte sur \,$\Im\,\varphi$\,} par
\emph{\pro sur \,$\Im\,\varphi$\, \paralm \`a
\,$\Ker(\varphi\sta)$}.

Nous allons voir maintenant que ces relations  (\ref{eqSDO}) peuvent
\^{e}tre r\'ealis\'ees de mani\`ere automatique sur un corps arbitraire
\,$\K$\, \`a condition d'introduire, \`a la place de  \,$A\sta$, une
matrice \,$A^{\circ}$\, \`a \coes dans le corps \,$\K(t)$\, o\`u
\,$t$\, est une \idtrz.

Pour cela nous nous limitons au point de vue purement matriciel, (c'est
le point de vue o\`u des bases ont \'et\'e fix\'ees dans \,$E$\, et
\,$F$).
Nous consid\'erons une forme quadratique \,$\Phi_{t,n}$\, sur
\,$E'=\K(t)^n$\, et une forme quadratique \,$\Phi_{t,m}$\, sur
\,$F'=\K(t)^m$:
$$
\begin{array}{rcl}
\Phi_{t,n}(\xi_1,\ldots ,\xi_n)& =  &{\xi_1}^2+t\, {\xi_2}^2 + \cdots +
t^{n-1}\,{\xi_n}^2    \\[1mm]
\Phi_{t,m}(\zeta_1,\ldots ,\zeta_m)& =  & {\zeta_1}^2 + t\,{\zeta_2}^2 +
\cdots + t^{m-1}\,{\zeta_m}^2
\end{array}
$$

Nous notons les \gui{produits scalaires}  correspondants par
\,$\gen{\cdot,\cdot}_{E'}^{t}$\, et \,$\gen{\cdot,\cdot}_{F'}^t$. Nous
notons \,$Q_n$\,
et \,$Q_m$\, les matrices (diagonales) de ces formes sur les bases
canoniques.

Toute \ali \,$\varphi:E\rightarrow F $\, donne lieu \`a une \ali
\,$E'\rightarrow F'$\, que nous notons encore \,$\varphi$\, et qui est
d\'efinie par la m\^{e}me matrice sur les bases canoniques. Il existe
alors une unique
\ali  \,$\varphi^{\circ}:F'\rightarrow E' $\, qui r\'ealise les \egts
(\ref{eqApAdj}) dans ce nouveau contexte:
\begin{equation} \label{eqApAdjbis}
\forall x\in E'\quad \forall y\in F'\quad \quad \gen{\varphi(x),y
}_{F'}^t=\gen{x,\varphi^{\circ}(y)}_{E'}^{t}
\end{equation}
La matrice \,$A^{\circ}$\, de \,$\varphi^{\circ}$\,
sur les bases canoniques est alors
$$
A^{\circ}\,=\, {Q_n}^{-1}\,\tra{\!A}\;Q_m\,,
$$
puisqu'on doit avoir pour tous
\,$X\in\K(t)^{n{\times}1},\;Y\in\K(t)^{m{\times}1}$\,
$$
\tra{(A\,X)}\,Q_m\,Y\;=\tra{\!X}\,{Q_n}\,(A^{\circ}\,Y)\,.
$$
En pratique si \,${A}=(a_{i,j})$\, on obtient
\,$A^{\circ}=(t^{j-i}\,a_{j,i})$, par exemple:
$$ A\;=\;\cmatrix{
a_{11}&a_{12}&a_{13}&a_{14}&a_{15}\cr
a_{21}&a_{22}&a_{23}&a_{24}&a_{25}\cr
a_{31}&a_{32}&a_{33}&a_{34}&a_{35}
}
,\quad A^{\circ}\;=\;\cmatrix{
a_{11}&t\,a_{21}&t^2\,a_{31}\cr
t^{-1}\,a_{12}&a_{22}&t\,a_{32}\cr
t^{-2}\,a_{13}&t^{-1}\,a_{23}&a_{33}\cr
t^{-3}\,a_{14}&t^{-2}\,a_{24}&t^{-1}\,a_{34}\cr
t^{-4}a_{15}&t^{-3}a_{25}&t^{-2}a_{35}
}
$$

Comme nous l'avons d\'ej\`a indiqu\'e, pour pouvoir reproduire (avec les
l\'eg\`eres variations \ncrsz) le fait \ref{factSDO},
les \dfns \ref{defGk} et \ref{defIMP}, les lemmes \ref{lemGk},
\ref{lem proj orth} et la proposition
\ref{propMoPe} il nous suffit de d\'emontrer l'analogue du lemme
\ref{lemSDO}.
\begin{lemma}
\label{lemMul}
Avec les notations ci-dessus on a pour toute matrice
\,$M\in\K^{m{\times}n}$, des sommes directes
\ortes dans les espaces \,$F'=\K(t)^{m}$\, et  \,$E'=\K(t)^{n}$\,
\begin{equation} \label{eqSDObis}
   \Im\,\varphi \oplus \Ker\,\varphi^{\circ}= F'\,, \; \;
\Ker\,\varphi\oplus \Im\,\varphi^{\circ}= E'\,
\end{equation}
\end{lemma}
\preuve
Les dimensions conviennent et il suffit de montrer que l'intersection
est r\'eduite \`a 0. Prenons par exemple la premi\`ere.
La relation (\ref{eqApAdjbis}) implique que l'\ort de
\,$ \Im\,\varphi$\, au sens de la forme bi\lin
\,$\gen{\cdot,\cdot}_{F'}$\, est \'egal \`a \,$\Ker\,\varphi^{\circ}$.
Il nous suffit donc de montrer que si \,$H$\, est
un sous-\evc de \,$F'=\K(t)^m$\, \emph{d\'efini sur \,$\K$\,}, son \ort
\,$H^{\perp}$\, dans \,$\K(t)^{m}$\, au sens du produit scalaire
\,$\gen{\cdot,\cdot}_{F'}^t$\, ne le coupe qu'en 0.
Soit donc \,$(p_1(t),\ldots ,p_m(t))\in H\cap H^{\perp}$.  Il existe
\,$v_1,\ldots ,v_r\in H$\, et \,$a_1(t),\ldots,a_r(t)\in\K(t)$\, tels
que \,$(p_1(t),\ldots ,p_m(t))=\sum_i a_i(t)\,v_i$.
Quitte \`a multiplier par le produit des \denos
on peut supposer que les \,$a_i$\,
sont des \pols et donc aussi les \,$p_i$. On peut introduire une
nouvelle \idtr \,$u$\, et travailler dans \,$\K[t,u]$. Alors puisque
\,$\sum_i a_i(t)\,v_i$\, est \ort \`a tous les \,$v_i$\, il est \ort
\`a \,$\sum_i a_i(u)\,v_i=(p_1(u),\ldots ,p_m(u))$\, et cela donne
$$ P(t,u)=\sum\nolimits_{i=1}^m\,p_i(t)\,p_i(u)\,t^{i-1}=0\,.
$$
Il nous reste \`a voir que cette relation implique que les \,$p_i$\,
sont tous nuls. Supposons l'un des \,$p_i$\, non nul. Soit \,$d\geq 0$\,
le plus grand des degr\'es des \,$p_{i}$. Soit \,$k$\, le plus grand
indice pour lequel \,$\deg p_{k}=d$\, et \,$a_k$\, le \coe dominant de
\,$p_k$. Alors on v\'erifie facilement que
le \coe de \,$u^d\,t^{d+k-1}$\, dans \,$P(t,u)$\,
est \'egal \`a \,$a_k^2$, et donc \,$P$\, est non nul.
\qed

\medskip Nous nous contentons maintenant de reproduire les \dfns et
r\'esultats dans notre nouveau cadre.

\begin{fact}
\label{factSDObis}~
\begin{enumerate}
\item \label{factSDO.2bis}  L'\ali \,$\varphi:E'\rightarrow F'$\, se
restreint en un \iso \,$\varphi_0$\, de \,$\Im\,\varphi^\circ$\, sur
\,$\Im\,\varphi$\, et \,$\varphi^\circ$\, se restreint en un \iso
\,$\varphi_0^\circ$\, de \,$\Im\,\varphi$\, sur \,$\Im\,\varphi^\circ$.
\item \label{factSDO.3bis} En outre:
\begin{equation} \label{SDO4bis}
\begin{array}{rclrrcl}
\Im\,\varphi  & =  &\Im\,\varphi\varphi^\circ\,,   & \quad&
\Ker\,\varphi^\circ  & = & \Ker\,\varphi\varphi^\circ\,, \\[1mm]
\Ker\,\varphi   & =  & \Ker\,\varphi^\circ\varphi\,,  &
&\Im\,\varphi^\circ & = & \Im\,\varphi^\circ\varphi\,.
\end{array}
\end{equation}
\item \label{factSDO.4bis} Soit \,$\varphi_1:\Im\,\varphi\rightarrow
\Im\,\varphi$\,
l'\auto \lin d\'efini par
\,$\varphi_1=\varphi_0\varphi_0^\circ$. C'est la restriction de
\,$\varphi\,\varphi^\circ$\, \`a
\,$\Im\,\varphi$.
Nous avons:
$$ \formule{\det{(\Id_{\Im\,\varphi} + Z \,\varphi_1)} & = &
\det{(\Id_F + Z \,\varphi \varphi^\circ)} \\ & = & 1+ a_1 Z+
\cdots + a_r Z^r\,.} $$
o\`u \,$r=\rg(\varphi)=\rg(\varphi_1)$\, et  \,$a_r\neq 0$. De la
m\^eme fa\c{c}on, nous avons l'\auto
\,$\varphi_1^\circ=\varphi_0^\circ\,\varphi_0$\, de
\,$\Im\,\varphi^\circ$\, et $$
\formule{\det{(\Id_{\Im\,\varphi^\circ} + Z \,\varphi_1^\circ)} &
= & \det{(\Id_E + Z \,\varphi^\circ \varphi)} \\ & = & 1+ a_1 Z+
\cdots+ a_r Z^r\,.} $$
\end{enumerate}
\end{fact}

Les \coes de la matrice \,$A\,A^\circ$\, sont des
\emph{\pols de Laurent}, autrement dit des \elts
de \,$\K[t,1/t]$.

\begin{definition}\label{defGkbis}
  Les {\em \pols de Gram (g\'en\'eralis\'es)
de \,$A$} \index{polyn\^{o}me de Gram!(g\'en\'eralis\'e)}
  \index{Gram!polyn\^{o}me de} sont les
\pols de Laurent \,$\cG'_{k}(A)(t)=a_{k}(t)\in\K[t,1/t]$,
\indexnota{Gpk@$\cG'_{k}(A)\in\K[t,1/t]$} et  les {\em \coes de Gram
\gnes de \,$A$}
\index{Gram!coefficient g\'en\'eralis\'e de} \index{coefficient de
Gram!g\'en\'eralis\'e} sont les \coes
\,$\cG'_{k,\ell}(A)=a_{k,\ell}$\, \indexnota{Gpkl@$\cG'_{k,\ell}(A)$}
donn\'es par la formule
\begin{equation} \label{eqDefGrambis}
\left\{
\begin{array}{rcll}
\det{(\I_m  +  Z \, AA^{\circ})} &  = & 1 \, + \, a_1(t) \, Z \, +
\, \cdots \, + \, a_m(t) \, Z^m     \\[2mm] a_{k}(t)&  = &
t^{-k(n-k)} \,\left(\sum_{\ell=0}^{k(m+n-2k)} \,a_{k,\ell}\,t^\ell
\right)
\end{array}
\right.
\end{equation}
Nous d\'efinissons aussi
\,$\cG'_0(A)=1{\quad\rm et\quad}\cG'_\ell(A)=0
{\quad\rm  pour\quad} \ell>m.$
\end{definition}

On v\'erifie ais\'ement que \,$\cG'_{k}(A)(t)=\cG'_{k}(\tra{\!A})(1/t)$.

Dans la suite de cette section, nous dirons,
pour abr\'eger, \gui{\polz} \`a la place de
\gui{\pol de Laurent} en laissant au lecteur le soin de d\'eterminer
selon le contexte si des puissances n\'egatives de la variable sont
pr\'esentes ou non.

Notons que les \coes de Gram usuels sont donn\'es par
\begin{equation} \label{eqGG}
\cG_{k}(A)=\cG'_{k}(A)(1)=\sum\nolimits_{\ell} \,a_{k,\ell}
\end{equation}

Les
\coes de Gram \gnlz is\'es sont des sommes de carr\'es de mineurs et ils
permettent de contr\^{o}ler le rang de la matrice en vertu du lemme
suivant, qui est l'analogue du lemme \ref{lemGk}.
\begin{lem}
\label{lemGkbis} \emph{(Conditions de Gram \gnees pour le rang)}
Soit
  \,$A\in\K^{m \times n}\subseteq \K(t)^{m \times n}$.

\sni $(1)$\, La matrice \,$A$\, est de rang \,$\leq r$\, \ssi les \pols
\,$\cG'_{k}(A)(t)$\, pour \,$k>r$\, sont identiquement nuls.
Elle est de rang \,$r$\, si en outre  $\cG'_{r}(A)\neq 0.$

\sni $(2)$\, Le \coe de Gram  \,$a_{k,\ell}=\cG'_{k,\ell}(A)$\, est
\'egal \`a la somme des carr\'es des  mineurs \,$\mu_{\alpha,\beta}$\,
d'ordre \,$k$\, de la matrice \,$A$\, extraits sur les lignes et les
colonnes
correspondant aux multi-indices
\,$\alpha=(\alpha_1,\ldots ,\alpha_k)$\,  et
\,$\beta=(\beta_1,\ldots ,\beta _k)$\, pour toutes les paires
\,$(\alpha,\beta)$\,  qui v\'erifient l'\egt \,$\sum_{i=1}^k\alpha_i-
\sum_{j=1}^k\beta_j=\ell-k(n-k)$.
\end{lem}

\medskip En particulier \,$\cG'_k(A)=0$\, si \,$k>p=\inf(m,n)$.

En posant \,$p=\inf(m,n)$\, et \,$p'=\sup(m,n)$\,  le nombre
total des
\coes de Gram \gnlz is\'es est  \'egal \`a:
$$
\sum_{k=1}^p\,\left(k\, (m+n-2k)+ 1 \right)=p+\frac{1}{6}\,p\,(p+1)\,
(3p'-p-2)\leq \frac{1}{2}p(p+1)p'\,.
$$

Nous avons les analogues du lemme \ref{lem proj orth} et de la
proposition \ref{propMoPe}.
\begin{lem}
\label{lem proj orth bis} {\rm  (\pros sur l'image et sur le noyau)}

\noi
Soient
  \,$m,\, n,\, r>0$\, dans \,$\NN$\, avec \,$r \le \min(m, n)$,
  \,$A\in\K^{m \times n}$\, une matrice de rang \,$r$.
Posons  $a_k(t)=\cG'_k(A)$.

\sni $(1)$\, La matrice de la \pro  sur le sous-espace
\,$\Im\,A \subseteq \K(t)^m$\, \paralm \`a \,$\Ker\,A^{\circ}$\,
est \'egale \`a
\begin{equation} \label{eqProImbis}
P=a_r^{-1}\left(a_{r-1}AA^{\circ}-a_{r-2}(AA^{\circ})^2+
\cdots +(-1)^{r-1}(AA^{\circ})^r\right)
\end{equation}

\noindent  $(2)$\, La matrice de la \pro  sur le
sous-espace
\,$\Im\,A^{\circ} \subseteq \K(t)^n$\,
  \paralm \`a \,$\Ker\,A$\,
est \'egale \`a
\begin{equation} \label{eqProKerbis}
P^{\bullet}=a_r^{-1}\left(a_{r-1}A^{\circ}A-a_{r-2}(A^{\circ}A)^2+
\cdots +(-1)^{r-1}(A^{\circ}A)^r\right)
\end{equation}
Et la matrice de \pro sur le noyau de \,$A$\,
  \paralm \`a \,$\Im\,A^{\circ}$\, est \,$\In-P^{\bullet}$.
\end{lem}

Notez qu'il s'agit de \pros \ortes par rapport aux
formes \bils \,$\gen{.,.}_{E'}^{t}$\, et \,$\gen{.,.}_{F'}^{t}$.

\ss\rem
En fait chaque formule peut \^{e}tre sp\'ecialis\'ee en rempla\c{c}ant
\,$t$\, par n'importe quelle valeur \,$\tau\in \K\setminus\{0\}$\, qui
n'annule pas le \deno \,$a_r(t)$\, (ce qui est toujours possible si le
corps poss\`ede au moins \,$r(m+n-2r)+1$\, \eltsz).

\begin{definition}
\label{defIMPbis}
Supposons que \,$\varphi$\, est de rang \,$r$.  \emph{L'\iMP \gne de
\,$\varphi$\, (en rang \,$r$)}
\index{inverse Moore-Penrose!g\'en\'eralis\'e en rang \,$r$} est  l'\ali
\,$\varphi^{[-1]_{r,t}}:F'\rightarrow E'$\,  d\'efinie par:
$$\forall y \in F'\qquad
\varphi^{[-1]_{r,t}}(y)= \varphi_0^{-1}( \pi_{\Im\,\varphi}(y)).$$
\end{definition}

\begin{prop}
\label{propMoPebis}
        \emph{(\iMP g\'en\'eralis\'e)}

\noi Soient
  \,$m, n, r>0$\, dans \,$\NN$\, avec \,$r \le \min(m, n)$,
  \,$A\in\K^{m \times n}$ de rang \,$r$, \,$V\in\K^{m\times 1}$.

\sni $(1)$\, L'\iMP \gne de \,$A$\, en rang \,$r$\,
est la matrice
  \,$A^{[-1]_{r,t}}\in\K(t)^{n \times m}$\,  \'egale \`a
\begin{equation} \label{eqMoPebis}
a_r^{-1}\left(a_{r-1}\I_m-a_{r-2}A^{\circ}A + \cdots+
(-1)^{r-1}(A^{\circ}A)^{r-1}\right) \; A^{\circ}\,
\end{equation}
o\`u \,$a_k=\cG'_k(A)$.

\sni $(2)$\, Le \sli \,$AX=V$\, admet une solution \ssi  le \pol
\,$\cG'_{r+1}(A \vert V)$\, est identiquement nul \ssi
\begin{equation} \label{eqMoPe2bis}
  V \; = \; A \; A^{[-1]_{r,t}} \; V
\end{equation}
Dans ce cas \,$X=A^{[-1]_{r,t}} V$\, est
l'unique solution dans  \,$\Im(A^\circ)$.
\end{prop}

\rem
Si \,$\varphi$\, est injective et \,$v\in\Im\,\varphi$\, est
repr\'esent\'e par un vecteur colonne \,$V\in\K^{m{\times}1}$\,
l'\elt \,$ A^{[-1]_{r,t}} \, V $\, est l'unique solution du \sli
correspondant. En cons\'equence, il ne d\'epend pas de \,$t$\,
et les fractions rationnelles donn\'ees par le calcul
des coordonn\'ees de
\,$ A^{[-1]_{r,t}} \, V $\, se simplifient en des constantes.

\subsubsection*{Cas des matrices sym\'etriques}

Dans ce paragraphe \,$E'=F'=\K(t)^{n\times n}$,
\,$\varphi$\, est d\'efini par une matrice
sym\'etrique \,$A=\tra{\!A}\in\K^{n\times n}$\,
et \,$\rg(\varphi)=r$. Soit \,$\lambda$\,
l'\auto \lin de \,$E'$\, defini par \,$Q_n$\,
par rapport \`a la base canonique.

D\'efinissons
\,$\whf= \lambda^{-1}\circ\varphi$,  \,$\wtf=\varphi\circ\lambda$.

La
matrice de \,$\whf$\, est \,$\wha=Q_n^{-1}A$, celle de
\,$\wtf$\, est \,$\wta=A\,Q_n$.

Puisque
\,$A$\, est sym\'etrique, on a:
$$\wha(t)=\tra{(\wta)}(1/t)\;\;\mathrm{et}\;\;\varphi^\circ=
\lambda^{-1}\circ\varphi\circ\lambda=
\whf\circ\lambda= \lambda^{-1}\circ\wtf\,.$$

On en d\'eduit:
$$\begin{array}{lll}
\Im\,\whf=\Im\,\varphi^\circ,& \Ker\,\whf=\Ker\,\varphi,  &
\Im\,\wtf=\Im\,\varphi,     \\[1mm]
\Ker\,\wtf=\Ker\,\varphi^\circ,&
\Im\,\varphi=\lambda(\Im\,\varphi^\circ), &
\Ker\,\varphi=\lambda(\Ker\,\varphi^\circ)      \\
\end{array}$$

Donc l'\'equation (\ref{eqSDObis})
peut \^{e}tre r\'e\'ecrite comme deux \decos
\ortes (par rapport \`a
la forme \bil \,$\gen{.,.}_{E'}^{t}$):
$$\begin{array}{c}
\Im\,\whf  \oplus \Ker\,\whf\;=\;\Im\,\varphi^\circ
\oplus \Ker\,\varphi\;=\; I^\circ
\oplus K\;=\; E',\\[1mm]
\Im\,\wtf  \oplus \Ker\,\wtf\;=\;\Im\,\varphi
\oplus \Ker\,\varphi^\circ\;=\; I
\oplus K^\circ\;=\; E'.
\end{array}$$

Nous notons \,${\wtfo}$\,  l'\auto
de \,$I$\, obtenu par restriction de \,$\wtf$.

Les \alis \,$\varphi$,   \,$\wtf$\, et \,$\wtfo$\,
ont m\^{e}me rang \,$r$, et la somme directe
\,$\Im\,\wtf\oplus\Ker\,\wtf=I
\oplus K^\circ=E'$\, entra\^{\i}ne que:
\begin{equation} \label{eqDefGHbis}
\det{(\I_n + Z \, \wtf )}  = \det{(\Id_{I}
  +  Z \, {\wtfo})} =
  1 +  b_1(t) \, Z  +
\cdots  +  b_r(t)  Z^r
\end{equation}
avec \,$b_r\neq 0$\, et \,$b_{r+1}=\ldots= b_n= 0$. Les \,$b_i(t)$\,
sont au signe pr\`es les \coes du
\polcar de \,$\wtf$.

On vient de d\'emontrer la
version simplifi\'ee du lemme \ref{lemGkbis}.
Ceci cons\-titue le r\'esultat cl\'e de Mulmuley
\cite{Mul}.
\begin{lem}
\label{lemGkbisSym} \emph{(Conditions de Mulmuley pour le rang
d'une matrice sym\'etrique)} Soit \,$\wta=A\,Q_n$\, avec
\,$A\in\K^{n{\times}n}$\, sym\'etrique. Soit \,$c_k=c_k(t)$\, le
\coe de \,$Z^{n-k}$\, dans le \polcar \,$\rP_{\!\wta}(Z)$\, de
\,$\wta$. \\
Alors, la matrice \,$A$\, est de rang \,$\leq r$\, \ssi les \pols
\,$c_k(t)$\, pour \,$k>r$\, sont identiquement nuls. \\
Elle est de
rang \,$r$\, si en outre  $c_{r}(t)\neq 0.$
\end{lem}

Puisqu'on a la somme \orte \,$\Im\,\wtf\oplus\Ker\,\wtf=E'$\,
on peut reproduire les calculs donn\'es
dans le cas des matrices sym\'etriques r\'eelles
et on obtient le r\'esultat suivant, qui simplifie
ceux obtenus
dans le cas d'une matrice arbitraire,
de fa\c{c}on similaire \`a la proposition \ref{propMoPeH}.
\begin{prop}
          \emph{(Un \ing d'une matrice sym\'etrique)} \\
Soit \,$A\in\K^{n{\times}n}$\, sym\'etrique de rang \,$r$,
\,$E'=\K(t)^n$\, et \,$\varphi:E'\rightarrow E'$\,
l'\ali d\'efinie par \,$A$.
On consid\`ere la matrice  \,$\wta=A\,Q_n$.
Les \coes \,$b_i$\, sont d\'efinis par l'\egt $(\ref{eqDefGHbis})$.
Dans la suite l'or\-thogonalit\'e s'entend par rapport \`a
la forme \bil \,$\gen{.,.}_{E'}^{t}$.

\sni $(1)$ La \pro \orte
\,$\pi_{\Im\,\varphi}$\,
sur le sous-espace \,$\Im\,\varphi$\, de
\,$E'$\, a pour matrice:
$$
P=b_r^{-1}
\left(b_{r-1}\,\wta-b_{r-2}\,\wta\,^2+ \cdots+
(-1)^{r-1}b_1\,\wta\,^{r-1}+
  (-1)^{r}\wta\,^r\right).
$$

\noi $(2)$ L'\iMP \gne
\,$\wtf\,^{[-1]_{r,t}}\in\L(E',E')$\, de
\,$\wtf$\, a pour matrice:
$$
\wta\,^{[-1]_{r,t}}={b_r}^{-1}  \,\left(b_{r-1}\,P -
b_{r-2}\,\wta   + \cdots+
(-1)^{r-1}b_1\,\wta\,^{r-2}+ (-1)^{r} \wta\,^{r-1}\right)
$$

\noi $(3)$ L'\endom
\,$\psi =\lambda\circ \wtf\,^{[-1]_{r,t}}$\, de
\,$E'$, dont la matrice est
\,$B=Q_n\,\wta\,^{[-1]_{r,t}}$, est un \emph{\ingz} de \,$\varphi $\,
au sens suivant:
$$ \varphi \circ \psi \circ \varphi =\varphi \quad \mathrm{et}
\quad
\psi \circ \varphi \circ \psi = \psi\,.
$$
L'\ali \,$ \varphi \circ \psi$\, est la \pro sur \,$\Im\,\varphi$\,
\paralm \`a \,$\Ker\,\psi$\, et
\,$\psi \circ \varphi $\, est la \pro sur \,$\Im\,\psi$\,
\paralm \`a \,$\Ker\,\varphi$.\\
Pour tout vecteur colonne \,$V$\, le \sli \,$AX=V$\, admet une solution
\ssi \,$ABV=V$\, et en cas de r\'eponse positive, \,$BV$\, est 
l'unique solution
dans l'espace \,$\Im\,\psi$.
\end{prop}
\prv
Il reste \`a prouver le point 3. Posons \,$\theta=\wtf\,^{[-1]_{r,t}}$.
  On sait que
$$\wtf\circ\theta\circ\wtf=\wtf\quad \mathrm{et}
\quad \theta\circ\wtf\circ\theta=\theta\,.$$
Puisque \,$\wtf=\varphi \circ \lambda$\, et \,$\psi =\lambda\circ\theta$,
cela donne tout de suite les deux \egts demand\'ees pour \,$\psi$\,
et \,$\varphi$. Tout le reste suit sans difficult\'e.
\qed

\ss Pour la th\'eorie des \ings nous recommandons les livres
\cite{Bha} et \cite{LM}.

\subsubsection*{Interpr\'etation par les \idts de Cramer}

Supposons la matrice \,$A$\, de rang \,$r$\, et \,$V$\, dans
l'espace engendr\'e par les colonnes de \,$A$. Appelons \,$C_j$\,
la \,$j$\,\eme colonne de \,$A$. Soit
$\mu_{\alpha,\beta}=\det{(A_{\alpha,\beta})}$\, le mineur d'ordre
\,$r$\, de la matrice \,$A$\, extrait sur les lignes
\,$\alpha=\{\alpha_1,\ldots ,\alpha_r\}$\, et les colonnes
\,$\beta=\{\beta_1,\ldots ,\beta _r\}$. Pour \,$j=1,\ldots ,r$\,
soit \,$\nu_{\alpha,\beta,j}$\, le \deter de la m\^{e}me matrice
extraite, \`a ceci pr\`es que la colonne \,$j$\, a \'et\'e
remplac\'ee par la colonne extraite de \,$V$\,  sur les lignes
\,$\alpha$. Alors on obtient pour chaque couple
\,$(\alpha,\beta)$\, de multi-indices une \idt de Cramer:
\begin{equation} \label{eqMPC1}
\qquad \mu_{\alpha,\beta}\;V=\sum\nolimits_{j=1}^r
\nu_{\alpha,\beta,j}\,C_{\beta_j}\qquad
\end{equation}
due au fait que le rang de la matrice \,$(A_{1..m,\beta}|V)$\, est
inf\'erieur ou \'egal \`a \,$r$\, (\cf l'\egt \vref{eqIC4}). Ceci peut
se relire comme suit:
\begin{eqnarray} \label{eqMPC1rc}
\qquad \mu_{\alpha,\beta}\;V&=& \left[
\begin{array}{ccc}
\nonumber
     C_{\beta_1} & \ldots  & C_{\beta_r}
\end{array}
\right] \cdot  \left[
\begin{array}{c}
       \nu_{\alpha,\beta,1} \\
     \vdots  \\
      \nu_{\alpha,\beta,r}
\end{array}
\right]=\\
\nonumber
  &=&\left[
\begin{array}{ccc}
     C_{\beta_1} & \ldots  & C_{\beta_r}
\end{array}
\right] \cdot
\Adj(A_{\alpha,\beta})
  \cdot
  \left[
  \begin{array}{c}
      v_{\alpha_1}  \\
      \vdots   \\
      v_{\alpha_r}
  \end{array} \right]=\\[2mm]
  &=&A \cdot (\I_n)_{1.. n,\beta} \cdot
\Adj(A_{\alpha,\beta}) \cdot (\I_m)_{\alpha,1.. m}\cdot V
\end{eqnarray}

Notons
\,$|\alpha|=\sum_{i=1}^r\alpha_i$,
\,$|\beta|=\sum_{i=1}^r\beta_i.$
Rappelons que
$$\cG'_{r}(A)(t)\;=\;\sum_{\alpha,\beta}\; t^{|\alpha|-
|\beta|}\;\;\mu_{\alpha,\beta}^2.$$
Si nous multiplions chaque \egt (\ref{eqMPC1rc}) par
\,$\mu_{\alpha,\beta}\,t^{|\alpha|-|\beta|}$\, et si nous additionnons
toutes ces \egts nous obtenons une expression de la forme:
$$\cG'_r(A)\cdot V =
A \cdot   \left( \sum\nolimits_{\alpha,\beta}\;
\mu_{\alpha,\beta}\,t^{|\alpha|-|\beta|} \cdot (\I_n)_{1.. n,\beta}
\cdot
  \Adj(A_{\alpha,\beta}) \cdot
(\I_m)_{\alpha,1.. m} \right) \cdot V
$$
Cette formule ressemble beaucoup trop \`a (\ref{eqMoPe2bis})
donn\'ee dans la proposition \ref{propMoPebis}:
$$V \; = \; A \; A^{[-1]_{r,t}} \; V
$$
pour ne pas \^{e}tre due \`a une \egt
\begin{equation} \label{eqMPC3rcbis}
\cG'_r(A) \,A^{[-1]_{r,t}}\,=\,\sum_{\beta,\alpha}\;
\mu_{\alpha,\beta}\,\,t^{|\alpha|-|\beta|}  \, (\I_n)_{1.. n,\beta}
\cdot\Adj(A_{\alpha,\beta}) \cdot
(\I_m)_{\alpha,1.. m}.
\end{equation}
Ainsi l'\iMP \gne peut \^{e}tre interpr\'et\'e
comme une somme pond\'er\'ee d'\idts de Cramer.

Nous ne prouverons cependant pas cette derni\`ere \egtz.
On peut la trouver, d\'emontr\'ee dans un cadre diff\'erent (plus
g\'en\'eral) et formul\'ee diff\'eremment, comme l'\egt 2.13
dans \cite{PB} ou, avec la m\^{e}me formulation qu'ici, dans
\cite{DiGLQ}.

\newpage \thispagestyle{empty}
 

\chapter{Algorithmes de base en alg\`ebre lin\'eaire}
\label{chap BasicAlgoAlin}
\minitoc

\subsection*{Introduction}

Il s'agit dans ce chapitre de d\'ecrire et d'analyser
certaines \mets \sqlesz, plus ou
moins classiques, pour le calcul du \deter et
du \polcar \`a \coes dans un anneau commutatif.

L'objectif recherch\'e est de comparer ces \algos
 \sqls et de d\'egager le meilleur
 possible, \cad le plus rapide
th\'eoriquement et pratiquement, occupant le moins
d'espace  m\'emoire possible (en \'evitant notamment
l'explo\-sion de la taille des  r\'esultats \itmdsz),
le plus facilement impl\'ementable sur machine
\sqle et le plus \gnlz, \cad
applicable dans un \acom arbitraire.

Nous introduirons plus loin (chapitre \ref{ChapNoCo})
des notions
pr\'ecises de \comz. Dans ce chapitre nous nous
contenterons de la
notion informelle de \com \arit donn\'ee par le compte
du nombre
d'\oparis dans l'anneau de base lors de l'ex\'ecution
de l'\algo
consid\'er\'e. Nous ferons \egmt quelques
commentaires, souvent
informels, sur le bon contr\^ole (ou non) de la taille
des r\'esultats
\itmdsz.

Nous commen\c{c}ons par l'\apg pour le calcul du
\deterz. C'est l'\algo
d'\agr \lin le plus classique. Il fonctionne sur un
corps et
poss\`ede de nombreuses applications (solutions de
\slis, calcul de
l'inverse, $\LU$-\deco \ldots). \emph{La \met du pivot
pour la r\'esolution des \slis est en fait due aux savants
chinois}:
on pourra consulter
\`a ce sujet la notice historique du chapitre 3 dans
l'ouvrage de Schrijver \cite{Schri} ainsi que l'\'etude
plus r\'ecente de Karine Chemla{\footnote{~D\`es
le troisi\`eme si\`ecle de notre \'ere, on trouve dans 
les commentaires de Liu Hui sur le texte classique 
\emph{Les neuf Chapitres} ce qu'il semble l\'egitime
d'appeler une preuve de correction de l'\apg pr\'esent\'e
dans ce texte ancien.}} \cite{Che1,Che2}.

Nous continuons avec un \algo qui a
pour m\'erite sa tr\`es grande simplicit\'e: l'\ajb
qui peut \^{e}tre vu comme une adaptation de la \mpgz,
avec un meilleur comportement des \coes \itmdsz.
Cet \algo fonctionne sur un \acom int\`egre,
\`a condition que les divisions exactes ne soit pas
trop co\^uteuses.
Dans le cas du calcul du \polcarv il devient un
\algo sans division et s'applique sur un \acomaz.
C'est ce que nous appelons la \mjb modifi\'ee. Une
variante de l'\ajb due \`a Dodgson (alias Lewis
Caroll) offre des perspectives int\'eressantes dans 
le cas des matrices structur\'ees.

Nous \'etudions ensuite l'\aghbz,
couramment uti\-li\-s\'e en analyse num\'erique.
Il utilise des divisions par des \elts non nuls
arbitraires
(on suppose donc qu'on travaille sur un corps). Mais
en calcul
formel, o\`u l'on veut des r\'esultats exacts, se pose
s\'erieusement le \pb de la croissance de la taille
 des  r\'esultats \itmdsz.

Nous signalons la \miL dans laquelle le calcul du
\polcar d\'epend
du calcul de plusieurs \detersz.

Nous examinons ensuite des \mets qui utilisent des
divisions
uniquement par des nombres entiers (de petite taille).
Il
s'agit de la \met de Le Verrier et de son
am\'elioration \`a la
Souriau-Faddeev-Frame.

Nous continuons avec les \mets sans division
 de Samuelson-\Ber et de Chistov.
Plus sophistiqu\'ees et nettement plus efficaces que
la \mjb modifi\'ee
elles fonctionnent \egmt sur un \acomaz.
L'\algo de Chistov pr\'esente les m\^emes \caras que
celui de Samuelson-\Ber mais r\'ev\`ele un l\'eger
handicap par rapport \`a ce dernier dans les tests
exp\'erimentaux.

Nous terminons avec les \mets qui utilisent les
\srlsz.
Celle que nous appelons \met de Frobenius a de bonnes
\caras tant du point de vue du nombre d'\oparis que
de la taille des r\'esultats \itmdsz. Elle ne
s'applique cependant
pas en toute \gnlz it\'e, et, hormis le cas des corps
finis,
elle est en pratique surpass\'ee par l'\agbz, sans
doute parce
que ce dernier n'utilise pas de division, et a besoin
de
moins d'espace m\'emoire
(meilleur contr\^{o}le des r\'esultats \itmdsz).
Nous exposons \egmt une variante due \`a Wiedemann.

\ss Dans ce chapitre, nous nous int\'eressons
seulement \`a des
versions assez simples des \algosz.

Nous dirons qu'une version d'un \algo est \emph{\elrz}
si les \muls de
matrices, de \pols ou de nombres entiers qui
interviennent
en son sein sont ex\'ecut\'ees selon la \met classique
\usle
(dite parfois \gui{na\"{\i}ve}). Pour les matrices et
les \pols
la \mul \usle  consiste \`a appliquer
simplement la formule d\'efinissant le produit.
Pour la \mul des entiers, il s'agit de l'\algo qu'on
apprend \`a
l'\'ecole primaire.

D'autre part, nous parlons de \emph{versions \sqlesz}
dans la mesure
o\`u  les \mets qui cherchent \`a acc\'el\'erer
l'ex\'ecution lorsque de nombreux processeurs
sont utilis\'es en \paral  ne sont pas envisag\'ees.

Dans ce chapitre, nous ne d\'eveloppons que des
versions \sqles \elrsz.

\ss Rappelons enfin la convention importante suivante:
dans tout cet ouvrage la notation \,$\log{n}$\, 
signifie
\,$\max(1,\log_{2}{n})$.


\section{M\'ethode du pivot de Gauss} 
\label{subsec Gauss}

\index{Gauss!m\'ethode du pivot}
C'est la \met la plus r\'epandue et la plus
courante aussi bien pour le calcul exact que pour
le calcul approch\'e des \deters lorsque
les \coes appartiennent \`a un corps
\,$\K$\, dans lequel les \ops de base
\,$(\,+\,,\,\,-\,,\,\times\,,\,/\,)$\, ainsi que le
test d'\egt
\`a 0 s'effectuent par des \algosz.

Son int\'er\^et r\'eside non seulement dans le fait
qu'elle poss\`ede plusieurs variantes
(symbo\-li\-ques ou num\'eriques) jouant un r\^ole
important dans la r\'eduction et l'inversion des
matrices et dans la r\'esolution des
\slisz, mais aussi dans le fait que la
technique du pivot est utilis\'ee dans d'autres
\mets de r\'eduction comme celle de
\JBz, ou pour le calcul du \polcar
 comme nous le verrons plus
loin avec, par exemple, les \mets de
\gui{\JB modifi\'ee} ou de Hessenberg.

\subsection{Transformations
\'el\'ementaires}\label{subsec.trel}

Une matrice est dite \emph{\tgu} (resp. \tgiz)
\index{matrice!triangulaire}
\index{triangulaire!matrice}
si  les \elts situ\'es au-dessous
de  (resp. au dessus de) la diagonale principale sont
nuls.
On dit \matg lorsque le contexte rend clair de quelle
variante
il s'agit. Une \matg est dite \emph{\utgz}
\index{matrice!unitriangulaire}
\index{unitriangulaire!matrice} si les \coes sur la
diagonale
principale sont tous \'egaux \`a $1$.

Bas\'ee sur l'id\'ee des \'eliminations
successives des inconnues dans la
r\'esolution d'un \sliz, la \mpg consiste \`a
r\'eduire
une matrice $A \in \K^{m\times n}$ \`a une matrice
\tgu par une
succession de transformations \elrs sur
les lignes (et \'eventuellement sur les colonnes)
de \,$A$.

Les \emph{transformations \elrsz} sur les lignes
\index{transformation \'el\'ementaire}
\index{elementaire@\'el\'ementaire!transformation}
d'une matrice sont de trois types:
\begin{itemize}
\item[(i)] multiplier une ligne
par un \elt non nul de \,$\K$;
\item[(ii)] \'echanger deux lignes;
\item[(iii)] ajouter \`a une ligne le produit
d'une autre ligne par un \elt de~~$\K$.
\end{itemize}

On d\'efinit de mani\`ere analogue les transformations
\elrs sur les
colonnes.

On associe \`a toute transformation \elr (sur les
lignes ou sur les colonnes) d'une matrice $A \in
\K^{m\times n}$
la matrice (dite \elrz)
\index{matrice!elementaire@\elr}
\index{elementaire@\'el\'ementaire!matrice}
obtenue en effectuant cette m\^eme
transformation \elr de la matrice unit\'e (matrice
unit\'e
d'ordre $m$ ou $n$ selon
le cas).
Toute transformation \elr sur les lignes (resp.
colonnes) de \,$A$\, revient alors \`a multiplier \`a
gauche
(resp. \`a droite) la  matrice \,$A$\, par
la matrice \elr correspondante.
Ceci est d\^u simplement au fait que si~ $L_1,L_2 \in
\:\K^{1\times n}$, on a pour tout $\lambda\in \K$:
$$
\left[\begin{array}{cc} \lambda & 0 \\ 0 & 1
\end{array}\right]
\left[\begin{array}{c} L_1 \\ L_2 \end{array}\right] =
\left[\begin{array}{c}
\lambda L_1
\\ L_2 \end{array}\right]~;~~~~
\left[\begin{array}{cc} 0 & 1 \\
1 & 0
\end{array}\right]
\left[\begin{array}{c} L_1 \\ L_2 \end{array}\right] =
\left[\begin{array}{c} L_2
\\ L_1
\end{array}\right]~;
$$
$$
 \mathrm{et} \: \: \left[\begin{array}{cc}
1 & 0 \\
\lambda & 1
\end{array}\right] \left[\begin{array}{c} L_1 \\ L_2
\end{array}\right] =
\left[\begin{array}{c} L_1 \\ L_2 + \lambda L_1
\end{array}\right]\cdot~(\footnote{~\'{E}galit\'{e}s
analogues pour les colonnes.})
$$

Il est clair que l'inverse d'une transformation \elr
sur les  lignes (resp. colonnes) est une
transformation
\elr de m\^eme  type sur les lignes (resp.
colonnes). Pr\'ecis\'ement:
$$
\left[\begin{array}{cc} \lambda & 0 \\ 0 & 1
\end{array}\right]^{-1} = \left[\begin{array}{cc}
\lambda^{-1} & 0 \\ 0 & 1 \end{array}\right]~~~,
~~~~\left[\begin{array}{cc} 0 & 1 \\ 1 & 0
\end{array}\right]^{-1} =
\left[\begin{array}{cc} 0 & 1 \\ 1 & 0
\end{array}\right]
$$
$$
~~\mathrm{et} ~~\left[\begin{array}{cc} 1 & 0 \\
\lambda & 1 \end{array}\right]^{-1} =
\left[\begin{array}{cc} 1 & 0 \\
-\lambda & 1 \end{array}\right]
~~\mathrm{pour~tout}
~\lambda \in \K.
$$

\ss Une matrice (et l'\ali correspondante)
est dite \refstepcounter{bidon}
\emph{unimodulaire} si elle est de
d\'eterminant 1.
\index{matrice!unimodulaire}\label{unimod}
Lorsqu'on veut se limiter aux transformations \elrs
qui correspondent au produit par une matrice
unimodulaire,
on a droit seulement \`a celles du troisi\`eme type.
N\'eanmoins, il est facile de voir qu'une succession
de
trois telles transformations permet d'obtenir un
\emph{\'echange sign\'e de lignes (ou de colonnes)}
du type \,$(L_{i},L_{j})\leftarrow (L_{j},- L_{i})$,
qui est consid\'er\'e comme la variante unimodulaire
des transformations \elrs du deuxi\`eme type.
Les \'echanges sign\'es et les transformations \elrs
du trois\`eme type sont appel\'ees
\emph{transformations \elrs unimodulaires}.
\index{transformation \'el\'ementaire!unimodulaire}
\index{unimodulaire!transformation \elrz}

\ss L'\'elimination de Gauss proprement dite que nous
consid\'erons
ici est essentiellement une succession de
transformations
\elrs du troisi\`eme type sur les lignes: des
\'echanges de lignes ou de colonnes
n'interviennent que s'il y a lieu de chercher un pivot
non nul pour le ramener au bon endroit.  Chaque
\'etape de
l'\algo de Gauss consiste donc
\`a traiter le pivot (non nul) issu de l'\'etape
pr\'ec\'edente,
en faisant appara\^{\i}tre des z\'eros au-dessous de
ce pivot,
et \`a d\'eterminer ensuite le pivot de l'\'etape
suivante pour
le placer sur la diagonale cons\'ecu\-ti\-vement au
pivot
pr\'ec\'edent.
Si on rempla\c{c}ait les \'echanges (de lignes ou de
colonnes)
sign\'es, on obtiendrait donc une r\'eduction
n'utilisant
que des transformations \elrs unimodulaires.

\ss En fait, il est bien connu, et c'est une
cons\'equence de la
\mpgz, que toute \maca
inversible est \'egale \`a un produit de matrices
\elrsz.
Et qu'en cons\'equence toute matrice (de n'importe
quel format)
peut \^{e}tre ramen\'ee par manipulations \elrs de
lignes et de
colonnes \`a une forme canonique du type suivant
$$
\left[\begin{array}{ccc} \Ir&|&0\\ -&-&-\\0&|&0
\end{array}\right]
$$
avec la possibilit\'e de lignes ou de colonnes vides.

Cette r\'eduction est d'une importance th\'eorique
capitale. Citons
par exemple  Gabriel \& Roiter \cite{GaRo} page 5, qui
donnent
d'ailleurs dans leur
chapitre 1 des extensions tr\`es int\'eressantes de la
\metz: [\ldots]
\emph{
en d\'epit de son \'evidence et de sa simplicit\'e, ou
peut-\^{e}tre gr\^{a}ce \`a elles,
cette r\'eduction est tr\`es utile, et son usage
r\'ep\'et\'e conduit \`a des r\'esultats profonds.}

Si on se limite aux transformations \elrs
unimodulaires, alors
la forme r\'eduite est la m\^{e}me que ci-dessus dans
le cas d'une
matrice
rectangulaire ou carr\'ee non inversible, et pour une
\maca inversible
il faut
modifier la forme r\'eduite en prenant son dernier
\coe diagonal non
\ncrt
\'egal \`a 1.

\subsection{La {\LU\,-\,\decoz}}
\label{subsec.ludec}
Lorsque le processus de triangulation d'une matrice
\,$A \in
\K^{m\times n}\,$
aboutit sans qu'aucune permutation de lignes ou de
colonnes
n'intervienne
--- ce qui a lieu si les \,$r$\, premi\`eres \smpds
de \,$A$, \,$r$\, \'etant le rang de \,$A$, sont
\regs --- et si l'on
garde en m\'emoire les matrices \elrs
associ\'ees aux transformations effectu\'ees, la \mpg
permet d'obtenir,
en
m\^eme temps que la triangulation de
\,$A$, ce qu'il est convenu d'appeler une
\LU-\emph{\decoz}\index{LU-decom@\LU-d\'ecomposition},
\cad une 
fa\c{c}on d'\'ecrire \,$A$\, sous la forme:
\,$A=LU$, o\`u \,$U \in \K^{m\times n}$\, est une
\matgu
(c'est la forme triangulaire recherch\'ee de \,$A$),
et \,$L
\in \K^{m\times m}\,$
une \mutgiz:
\,$L$\, n'est autre que l'inverse du produit des
matrices
\elrs
correspondant aux transformations successives
effectu\'ees sur
les lignes de \,$A$.

Pour une \maca \regz, l'existence
d'une telle \deco \'equivaut au fait que le processus
de
triangulation arrive \`a son terme sans aucun
\'echange de lignes
ni de colonnes. Elle
\'equivaut aussi \`a la compl\`ete r\'egularit\'e de
la matrice
puisque les mineurs principaux de la matrice
consid\'er\'ee ne
sont autres que les produits successifs des pivots
rencontr\'es
au cours du processus. Enfin, toujours dans le cas
d'une \maca
\regz, l'existence de la \deco implique son unicit\'e.
Cela ne serait plus le cas pour une matrice
singuli\`ere comme
on peut le voir ici
$$
\left[
{\begin{array}{rrrrrr} 1  &  0   &   0   \cr
  2 &  1   &    0  \cr
   3 &  4   & 1\end{array}}
 \right]\,\,
\left[
{\begin{array}{rrrrrr} 1  &  3   &   1   \cr
  0 &  0   &    1  \cr
   0 &  0   & 1\end{array}}
 \right]
=\left[
{\begin{array}{rrrrrr} 1  &  3   &   1   \cr
  2 &  6   &    3  \cr
   3 &  9   & 8\end{array}}
 \right]
=\left[
{\begin{array}{rrrrrr} 1  &  0   &   0   \cr
  2 &  1   &    0  \cr
   3 &  0   & 1\end{array}}
 \right]
\left[
{\begin{array}{rrrrrr} 1  &  3   &   1   \cr
  0 &  0   &    1  \cr
   0 &  0   & 5\end{array}}
 \right].
$$

\ss Nous donnons maintenant \`a voir le r\'esultat de
la \mpg avec des matrices \`a coefficients entiers.

\begin{examples}
\label{exaGauss1} {\footnotesize
\emph{Nous montrons deux exemples caract\'erisitiques,
o\`u tous les pivots
qui se pr\'esentent sur la diagonale sont non nuls.
Nous donnons les matrices \,$L$\, et \,$U$.
Le premier est celui d'une matrice
dont les coefficients entiers ne prennent pas plus que
$2$ chiffres. Sur la premi\`ere ligne les matrices
\,$M_1$\, et \,$L_1$,
ensuite la ma\-tri\-ce~$\,U_1$.
\[
\left[
{\begin{array}{rrrrrr}
9 & 7 & 8 & 11 & 13 & 4 \\
19 & 4 & 56 & 84 & 73 & 10 \\
35 & 62 & -13 & 17 & 23 & 11 \\
20 & 3 & 6 & 7 & 5 & 9 \\
49 & 23 & 50 & 42 & 2 & 17
\end{array}}
 \right] ,\,
 \left[
{\begin{array}{rrrrr}
1 & 0 & 0 & 0 & 0 \\
{\ds \frac {19}{9}}  & 1 & 0 & 0 & 0 \\ [2ex]
{\ds \frac {35}{9}}  & {\ds \frac {-313}{97}
}  & 1 & 0 & 0 \\ [2ex]
{\ds \frac {20}{9}}  & {\ds \frac {113}{97}}
 & {\ds \frac {-5562}{7963}}  & 1 & 0 \\ [2ex]
{\ds \frac {49}{9}}  & {\ds \frac {136}{97}}
 & {\ds \frac {-4694}{7963}}  & {\ds \frac {
-21433}{244718}}  & 1
\end{array}}
 \right],
\]
\[U_1=
 \left[
{\begin{array}{rccccc}
9 & 7 & 8 & 11 & 13 & 4 \\
0 & {\ds \frac {-97}{9}}  & {\ds \frac {352}{
9}}  & {\ds \frac {547}{9}}  & {\ds \frac {
410}{9}}  & {\ds \frac {14}{9}}  \\ [2ex]
0 & 0 & {\ds \frac {7963}{97}}  & {\ds
\frac {16523}{97}}  & {\ds \frac {11586}{97}}  &
{\ds \frac {45}{97}}  \\ [2ex]
0 & 0 & 0 & {\ds \frac {244718}{7963}}  &
{\ds \frac {51521}{7963}}  & {\ds \frac {
-10965}{7963}}  \\ [2ex]
0 & 0 & 0 & 0 & {\ds \frac {-15092695}{244718}}  &
{\ds \frac {-1665525}{244718}}
\end{array}}
 \right].
\]
Le deuxi\`eme exemple est celui d'une matrice \`a
coefficients dans
$\QQ$.
Le \nume et le \deno n'ont qu'un chiffre, mais la
croissance
de la taille des coefficients est spectaculaire.
Sur la premi\`ere ligne \,$M_{2}$\, et \,$U_2$, sur la
seconde \,$L_2$.
$$
 \left[
{\begin{array}{ccccc}
{\ds \frac {1}{6}}  & {\ds \frac {3}{2}}  &
{\ds \frac {-9}{5}}  & {\ds \frac {7}{6}}  &
{\ds \frac {-7}{6}}  \\ [2ex]
{\ds \frac {3}{2}}  & {\ds \frac {-9}{8}}  &
1 & {\ds \frac {5}{7}}  & {\ds \frac {2}{9}}
 \\ [2ex]
{\ds \frac {-1}{9}}  & {\ds \frac {2}{3}}  &
{\ds \frac {7}{6}}  & {\ds \frac {1}{8}}  &
{\ds \frac {9}{5}}  \\ [2ex]
{\ds \frac {7}{8}}  & {\ds \frac {3}{4}}  &
{\ds \frac {1}{4}}  & {\ds \frac {-7}{9}}  &
{\ds \frac {-4}{3}}  \\ [2ex]
{\ds \frac {-1}{3}}  & -1 & {\ds \frac {-7}{6
}}  & {\ds \frac {-4}{9}}  & {\ds \frac {6}{7
}}  \\ [2ex]
{\ds \frac {9}{8}}  & {\ds \frac {-1}{2}}  &
{\ds \frac {2}{5}}  & {\ds \frac {-5}{9}}  &
{\ds \frac {9}{8}}  \\ [2ex]
{\ds \frac {2}{3}}  & {\ds \frac {-8}{7}}  &
-1 & {\ds \frac {4}{9}}  & {\ds \frac {-3}{7}
}
\end{array}}
 \right] ,
\left[
{\begin{array}{ccccc}
{\ds \frac {1}{6}}  & {\ds \frac {3}{2}}  &
{\ds \frac {-9}{5}}  & {\ds \frac {7}{6}}  &
{\ds \frac {-7}{6}}  \\ [2ex]
0 & {\ds \frac {-117}{8}}  & {\ds \frac {86}{
5}}  & {\ds \frac {-137}{14}}  & {\ds \frac {
193}{18}}  \\ [2ex]
0 & 0 & {\ds \frac {6763}{3510}}  & {\ds
\frac {-4175}{19656}}  & {\ds \frac {35446}{15795}} 
\\
 [2ex]
0 & 0 & 0 & {\ds \frac {-3391183}{1704276}}  &
{\ds \frac {-959257}{486936}}  \\ [2ex]
0 & 0 & 0 & 0 & {\ds \frac {25849022797}{10254937392}}
 \\ [2ex]
0 & 0 & 0 & 0 & 0 \\
0 & 0 & 0 & 0 & 0
\end{array}}
 \right]
$$
$$
L_{2}= \left[
{\begin{array}{cccccrr}
1 & 0 & 0 & 0 & 0 & 0 & 0 \\
9 & 1 & 0 & 0 & 0 & 0 & 0 \\
{\ds \frac {-2}{3}}  & {\ds \frac {-40}{351}
}  & 1 & 0 & 0 & 0 & 0 \\ [2ex]
{\ds \frac {21}{4}}  & {\ds \frac {19}{39}}
 & {\ds \frac {4635}{6763}}  & 1 & 0 & 0 & 0 \\ [2ex]
-2 & {\ds \frac {-16}{117}}  & {\ds \frac {
-8475}{6763}}  & {\ds \frac {-969733}{6782366}}  & 1 &
0 & 0 \\ [2ex]
{\ds \frac {27}{4}}  & {\ds \frac {85}{117}}
 & {\ds \frac {381}{13526}}  & {\ds \frac {
8966489}{13564732}}  & {\ds \frac {251177120859}{
258490227970}}  & 1 & 0 \\ [2ex]
4 & {\ds \frac {400}{819}}  & {\ds \frac {
-54066}{47341}}  & {\ds \frac {-3752551}{23738281}}  &
{\ds \frac {64239864618}{129245113985}}  & 0 & 1
\end{array}}
 \right]  \,,
$$
}
}
\end{examples}

Nous allons comprendre ces comportements typiques en
exprimant
pr\'ecis\'ement les \coes calcul\'es dans la \mpg en
fonction de \deters
extraits de la matrice initiale. Nous avons pour cela
besoin de
pr\'eciser les notations.

\begin{nota}
\label{nota Gauss}
Soit \,$A \in \K^{m\times n}$ une matrice de rang $r$.
On suppose que la
triangulation de Gauss aboutit
\`a son terme sans \'echange de ligne ni de colonne.
Dans ces conditions on pose \,$A^{[0]}=A$, on note
\,$A^{[p]}$\, la matrice
transform\'ee de \,$A$\, \`a l'issue de l'\'etape
\,$p$\,
($p\leq r$)
et on note \,$L^{[p]}$\, le produit
des matrices \elrs correspondant aux transformations
effectu\'ees
au cours de l'\'etape \,$p$, de sorte que
\,$A^{[p]}=L^{[p]}A^{[p-1]}$\, et
\,$L^{[p]}$\, est une matrice qui ne
diff\`ere de la matrice unit\'e \,$\I_m$\, que des
\elts
de la $p\,$\eme colonne situ\'es au-dessous
de la diagonale principale.
On note \,$a_{\,ij}^{[p]}$\,  l'\elt en  position
\,$(i,j)$\, de la
matrice \,$A^{[p]}$ et \,$l_{\,ij}^{[p]}$\, celui de
la matrice
\,$L^{[p]}$.
\end{nota}

Le symbole de Kronecker est d\'efini
page \pageref{Kron} et la notation
\,$a_{\,ij}^{(p)}\,$
 page \pageref{minrij}. On a alors:

\begin{propri}
\label{propri Gauss}
Avec les notations pr\'ec\'edentes, les \'el\'ements
\,$l_{\,ij}^{[p]}$, \,$a_{\,ij}^{[p]}$\, et
$a_{\,ij}^{(p)}$ sont li\'es par les relations
suivantes
(dans $(\ref{EqGauss1})$ on a $1\leq p\leq r$,  $p<
j\leq n$ et $p<
i\leq
m$):
\begin{equation} \label{EqGauss1}
a_{\,ij}^{[p]} = a_{~ij}^{[p-1]}-
\frac{a_{~ip}^{[p-1]}}{a_{~pp}^{[p-1]}}\,a_{~pj}^{[p-1]}=
\frac{ a_{~ij}^{[p-1]} \, a_{~pp}^{[p-1]} -
a_{~ip}^{[p-1]} \, a_{~pj}^{[p-1]} } { a_{~pp}^{[p-1]}
}.
\end{equation}
\begin{equation} \label{EqGauss2}
l_{\,ip}^{[p]}=-\frac{a_{~ip}^{[p-1]}}{a_{~pp}^{[p-1]}}
~~\mathrm{ si} ~i > p,~~~
l_{\,ij}^{[p]}=\delta_{ij}~~ \mathrm{sinon}.
\end{equation}
\begin{equation} \label{EqGauss3}
a_{\,11}^{[0]}  a_{\,22}^{[1]} \cdots
a_{~pp}^{[p-1]}= a_{~pp}^{(p-1)}.
\end{equation}
\begin{equation} \label{EqGauss4}
a_{\,ij}^{[p]}=
\frac{a_{\,ij}^{(p)}}{a_{~pp}^{(p-1)}}.
\end{equation}
\begin{equation} \label{EqGauss5}
l_{\,ip}^{[p]}=-\frac{a_{~ip}^{(p-1)}}{a_{~pp}^{(p-1)}}
~~\mathrm{ si} ~i > p,~~~
l_{\,ij}^{[p]}=\delta_{ij}~~ \mathrm{sinon}.
\end{equation}
\end{propri}
\prv
Les deux premi\`eres \'equations correspondent
exactement aux affectations de  l'\algo de
Gauss. Les deux suivantes correspondent au fait que
les \deters des sous matrices
correspondantes de \,$A$\, sont inchang\'es par les
transformations \elrs de
lignes utilis\'ees dans l'\algoz.
La derni\`ere r\'esulte de la deuxi\`eme et la
quatri\`eme.
\qed

\ss Il est clair que la matrice \,$U=A^{[r-1]}$\,
(o\`u \,$r=\rg(A)$) obtenue \`a l'issue
de la derni\`ere \'etape de l'\algo de Gauss dans ce
cas,
est bien la forme
\tgu recherch\'ee de la matrice \,$A$, et
que \,$A=LU$\, o\`u
\begin{equation} \label{EqGaussL1}
L=\left[l_{ij}\right]=\left[L^{[1]}\right]^{-1}
\left[L^{[2]}\right]^{-1}\cdots
\left[L^{[r-1]}\right]^{-1}
\end{equation}
est une \matgi avec en outre
\begin{equation} \label{EqGaussL2}
l_{\,ij}=-l_{\,ij}^{[j]}=\frac{a_{\,ij}^{(j-1)}}{a_{\,jj}^{(j-1)}}
~~\mathrm{ si} ~ m\geq i > j\geq 1,~~~
l_{\,ij}=\delta_{ij}~~ \mathrm{sinon}.
\end{equation}
En effet la matrice
\,$\left[L^{[p]}\right]^{-1}$\, ne diff\`ere de
\,$L^{[p]}$\, que des
\elts
\,$l_{\,ip}^{[p]}$\, pour \,$1\leq p<i\leq m$\, qui
doivent \^etre remplac\'es par leurs oppos\'es, et la
\mul \`a gauche par
\,$\left[L^{[p-1]}\right]^{-1}$\, du produit
\,$\left[L^{[p]}\right]^{-1}\,\alb\left[L^{[p+1]}\right]^{-1}
\,\alb\cdots\,\alb \left[L^{[r-1]}\right]^{-1}$\,
n'affecte que la
$(p-1)\,$\eme colonne de ce dernier (identique \`a la
$(p-1)\,$\eme colonne de $\I_m$) et revient tout
simplement \`a la remplacer par la
$(p-1)\,$\eme colonne de
$\left[L^{[p-1]}\right]^{-1}$.

Remarquons aussi que la relation (\ref{EqGauss3})
montre comment
l'\apg permet de calculer les \mips de la matrice
\,$\A$\, (et donc
son \deter lorsqu'elle est carr\'ee).

\ss Nous comprenons maintenant dans le cas d'une
matrice initiale \`a
\coes entiers le comportement typique de la taille des
\coes calcul\'es
dans
la \mpg (\cf la matrice $M_1$ de l'exemple
pr\'ec\'edent).
On voit sur les relations (\ref{EqGauss4}),
(\ref{EqGauss2}) et
(\ref{EqGaussL2}) que tous ces \coes peuvent \^{e}tre
\'ecrits comme des
fractions
dont le \nume et le \deno sont des mineurs de la
matrice initiale.
En outre les mineurs sont major\'es (en valeur
absolue, donc aussi en
taille si ce sont des entiers) en utilisant l'\inHz.
Grosso modo, en partant d'une matrice \`a $k$ lignes
avec des \coes de
taille $\tau$, on obtient dans l'\apg des \coes de
taille $k\tau$.
Pour ce qui concerne une matrice \`a \coes dans $\QQ$
(comme  $M_2$), pour obtenir
une majoration de la taille des \coes calcul\'es, nous
devons remplacer
$M_2$ par une matrice \`a \coes entiers $M'_2=c\,M_2$ 
(o\`u $c$ est le
ppcm
des \denosz). Grosso modo, en partant d'une matrice
\`a $k$ lignes
avec des \coes dont le \deno et le \nume sont de
taille $\tau$, on
obtient
maintenant dans l'\apg des \coes de taille $k\tau^2$.

\subsubsection*{Algorithme  du pivot de Gauss
simplifi\'e}

L'\algo simplifi\'e pour la \mpg s'applique pour les
matrices \fregsz.
Dans ce cas, il n'y a pas de
recherche de pivot
et la matrice est de rang maximum
$\inf(m,n)$. Cet \algo  remplace la matrice \,$A$\,
par une matrice de m\^emes
dimensions dont la partie \supee (diagonale
principale comprise) est celle de la
 matrice \,$U$\, et la partie  \infee (sans la
diagonale){\footnote{~Et sans les
\elts nuls en position \,$(i,j)$\, avec \,$i>j>n$\,
lorsque \,$m\geq n+2$.}}
celle de la matrice \,$L$\, de la \LU-\deco de
\,$A$. On obtient l'\algo \vref{apgs}.

\begin{algorH}[Algorithme  du pivot de Gauss
simplifi\'e
   (sans recherche de pivot)  et
\LU-\decoz.]\label{apgs}
\acl{apgs}{Algorithme  du pivot de Gauss simplifi\'e
et \LU-\decoz}
\Entree Une matrice $A=(a_{ij})\in\K^{m\times n}$
\fregz.
\Sortie La matrice $A$ transform\'ee ainsi que les
matrices $L$
et $U$ comme expliqu\'e ci-dessus.
\Debut
\Varloc $i$, $j$, $p\in\N\;$; $\;piv\in\A\;$;
\hsu \pour{p}{1}{\inf(m,n)}
\hsd $piv:=a_{pp}\;;$
\hsd \pour{i}{p+1}{m}
\hst $a_{ip}:=a_{ip} / piv\;;$
\hst \pour{j}{p+1}{n} $a_{ij}:=a_{ij}- a_{ip}*a_{pj}$
\hst \finpour
\hsd \finpour
\hsu \finpour
\fin
\end{algorH}

En fait la derni\`ere \'etape ($\,p=\inf(m,n)$) de la 
boucle
principale ne s'ex\'ecute que si \,$m>n$\, et elle ne
modifie
alors que les valeurs des \,$a_{in}$\, pour \,$i>n$.
On aurait donc pu \'ecrire
\textsf{\pour{p}{1}{\inf(m,n)-1}\ldots }
mais il aurait fallu rajouter \`a la fin:

\begin{agH}{
\hspace*{4mm}\textsf{\sialors{m>n}
\hsd $piv:=a_{nn}\;;$
\hsd \pour{i}{n+1}{m}
\hst $a_{in}:=a_{in} / piv\;;$
\hsd \finpour
\hsu \finsi}}
\end{agH}

Un calcul \elr donne le r\'esultat suivant.

\begin{proposition}
\label{propComArG} Le nombre d'\oparis dans 
\,$\K$\, lorsqu'on ex\'ecute l'\apg simplifi\'e, est
major\'e par:

\smallskip \centerline{$
n(m-1)(2n-2m+1) +{1\over 6}m(m-1)(4m-5)
$} 

\sni ce qui donne pour \,$m=n$\, la majoration
\,${2\over 3} n^3 - {1\over 2}n^2 - {1\over 6} n$.
\end{proposition}

Si la matrice \,$A$\, est de rang \,$r$\, et si les
\,$r$\, premiers \mips dominants sont
non nuls, l'\algo pr\'ec\'edent, modifi\'e pour
s'arr\^{e}ter lorsque le pivot \,$piv$\, est nul,
fournit encore la \LU-\deco de \,$A$.
Cela donne l'\algo \vref{apgs2}.

\begin{example}
\label{exaGauss2} {\footnotesize
\emph{Voici une matrice  \,$M_{3}\in\ZZ^{6\times5}$\,
de rang $4$, suivie des matrices \,$L_{3}$\, et
\,$U_3$\,
obtenues \`a partir de l'\algo \ref{apgs2}.
$$
\left[
{\begin{array}{rrrrr}
-73 & -53 & -30 & 45 & -58 \\
21 & -54 & -11 & 0 & -1 \\
72 & -59 & 52 & -23 & 77 \\
33 & 55 & 66 & -15 & 62 \\
-41 & -95 & -25 & 51 & -54 \\
14 & 55 & 35 & -5 & 25
\end{array}}
 \right], \,
 \left[
{\begin{array}{ccccrr}
1 & 0 & 0 & 0 & 0 & 0 \\
{\ds \frac {-21}{73}}  & 1 & 0 & 0 & 0 & 0 \\ [2ex]
{\ds \frac {-72}{73}}  & {\ds \frac {8123}{
5055}}  & 1 & 0 & 0 & 0 \\ [2ex]
{\ds \frac {-33}{73}}  & {\ds \frac {-2266}{
5055}}  & {\ds \frac {220594}{272743}}  & 1 & 0 & 0 \\
 [2ex]
{\ds \frac {41}{73}}  & {\ds \frac {4762}{
5055}}  & {\ds \frac {52277}{272743}}  &
{\ds \frac {1193}{949}}  & 1 & 0 \\ [2ex]
{\ds \frac {-14}{73}}  & {\ds \frac {-1091}{
1685}}  & {\ds \frac {83592}{272743}}  &
{\ds \frac {1052}{949}}  & 0 & 1
\end{array}}
 \right]
$$
$$
\left[
{\begin{array}{rcccc}
-73 & -53 & -30 & 45 & -58 \\
0 & {\ds \frac {-5055}{73}}  & {\ds \frac {
-1433}{73}}  & {\ds \frac {945}{73}}  & {\ds
\frac {-1291}{73}}  \\ [2ex]
0 & 0 & {\ds \frac {272743}{5055}}  & {\ds
\frac {196}{337}}  & {\ds \frac {243716}{5055}}  \\
 [2ex]
0 & 0 & 0 & {\ds \frac {2911532}{272743}}  &
{\ds \frac {-3038698}{272743}}  \\ [2ex]
0 & 0 & 0 & 0 & 0 \\
0 & 0 & 0 & 0 & 0
\end{array}}
 \right]
$$
}
}
\end{example}

\begin{algor}[Deuxi\`eme \apg simplifi\'e (sans
recherche de pivot)
    et \LU-\decoz.] \label{apgs2}
\acl{apgs2}{Deuxi\`eme \apg simplifi\'e}
\Entree Une matrice
$A=(a_{ij})\in\K^{m\times n}$.
\Sortie La matrice $A$
transform\'ee ainsi que le rang $r$ de $A$ lorsque
celui-ci est
\'egal \`a l'ordre du dernier \mip dominant non nul.
On obtient
\egmtz, dans ce cas, la \LU-\deco de la matrice $A$
comme dans
l'\algo \ref{apgs}. \Debut \Varloc $i$, $j$, $p$,
$r\in\N$ ;
               $\;piv\in\K$ ;
 \hsu  $p:=1$ ;  $r:=\inf{(m,n)}$ ;
 \hsu  \tantque{p\leq  \inf{(m,n)}}
 \hsd    $piv:=a_{pp}$ ;
 \hsd    \sialors{piv= 0} $r:=p-1$~;
$p:=\inf{(m,n)}\;$ \sinon
 \hst      \pour{i}{p+1}{m}
 \hsq         $a_{ip}:=a_{ip} / piv$~;
 \hsq         \pour{j}{p+1}{n}
 \hsc             $a_{ij}:= a_{ij} - a_{ip}* a_{pj}$
 \hsq         \finpour
 \hst      \finpour
 \hsd    \finsi ;
 \hsd    $p:= p+1$
 \hsu  \fintantque
\fin
\end{algor}

\subsection[Recherche de pivot non nul]{Algorithmes
avec recherche de
pivot non nul}\label{subsec.rpnn}

Si on rencontre un pivot nul sur la diagonale
principale au
cours du processus
de triangulation on doit proc\'eder \`a des \'echanges
de lignes
et/ou de colonnes pour ramener un pivot en position
convenable
(s'il reste un \elt non nul dans le coin sud-est).
Alors ce n'est pas une \LU-\deco de
\,$A$\, que l'on obtient avec la \mpgz, mais une
\textit{\~{P}LUP}-\deco (voir par exemple
\cite{Aho,BP}),
\cad  une \LU-\deco du produit
\`a droite et \`a gauche
de la matrice \,$A$\, par  des matrices de
permutation.

De mani\`ere plus pr\'ecise, si \`a l'issue de
l'\'etape
\,$p-1$\, du processus de
triangulation, on obtient un pivot nul ($a_{~pp}^{[p-
1]}=0$), alors de deux choses
l'une: ou bien \,$a_{~ij}^{[p-1]}=0$\, pour tous
\,$i,j\geq p$\,
auquel cas le  rang de \,$A$\, est \'egal \`a
\,$p-1$, et le processus est termin\'e, ou bien on
peut
trouver des entiers \,$i\,$
et \,$j$\,  $\geq p$\, et tels que
\,$a_{~ij}^{[p-1]}\neq 0$.
Dans ce cas, une permutation de lignes et/ou de
colonnes doit
intervenir pour
remplacer le pivot nul par l'\elt $a_{~ij}^{[p-1]}$:
ce
qui revient \`a
remplacer la matrice \,$A^{[p-1]}$\, par la matrice
\,$E_{ip}(m)\:A^{[p-1]}\:E_{jp}(n)$\, o\`u
\,$E_{kl}(h)$\,
d\'esigne la matrice
\elr obtenue \`a partir de \,$\I_h$\, par \'echange
des lignes \,$k$\, et
\,$l$\, (ou, ce qui revient au m\^eme, par \'echange
des
colonnes \,$k$\, et \,$l$).
Cette \opz, qui pr\'epare \,$A^{[p-1]}$\, \`a subir
avec
succ\`es l'\'etape
\,$p$, n'alt\`ere pas les \,$p-1$\, premi\`eres lignes
et les
\,$p-1$\, premi\`eres
colonnes de cette matrice.
Plus pr\'ecis\'ement, elle commute avec les \ops de
type \gui{traitement
d'un pivot} d\'ej\`a effectu\'ees (qui correspondent
au produit \`a
gauche par une \matgiz).
Par exemple, si sur une matrice $6\times 6$ on doit
faire des \'echanges
de lignes et de colonnes avant de traiter les pivots
\num$3$ et $5$ on
obtiendra la \deco suivante

\smallskip  \centerline{$ 
L^{[5]}\,Q_5\,L^{[4]}\,L^{[3]}\,Q_3\,L^{[2]}\,L^{[1]} 
= 
L^{[5]} \, \tilde{L}^{[4]}\, \tilde{L}^{[3]}\, 
\tilde{L}^{[2]}\, \tilde{L}^{[1]}\,Q_5\,Q_3 = L_5 \,Q_5\,Q_3 
$} 

\sni  o\`u 

\centerline{$\,\tilde{L}^{[4]}= Q_5\,L^{[4]}\, Q_5\,$, 
$\,\tilde{L}^{[3]}= Q_5\,L^{[3]}\,Q_5 \,$,} 

\smallskip 
\centerline{$\,\tilde{L}^{[2]}= Q_5\, Q_3\,L^{[2]}\,Q_3 \,Q_5 \,$, 
$\,\tilde{L}^{[1]}= Q_5\, Q_3\,L^{[1]}\, Q_3 \,Q_5\,$}

\sni   et 

\centerline{$\,L_5 = L^{[5]} \, \tilde{L}^{[4]}\, \tilde{L}^{[3]}\, 
\tilde{L}^{[2]}\, \tilde{L}^{[1]}$,}




\sni et donc

\smallskip \centerline{$
A=  \,Q_3\,Q_5(L_5)^{-1}\, U \,P_5\, P_3=
\tilde{P}\,L\,U\,P\,.
$}

\smallskip  Ainsi le processus de triangulation de Gauss,
lorsqu'une recherche
de pivots intervient, se ram\`ene \`a un processus
sans
recherche de pivot
sur le produit \`a droite et \`a gauche de la matrice
\,$A$\,
par des matrices
de permutation. Cela montre aussi que l'\apgz,
appliqu\'e
\`a la matrice \,$A$, donne, en m\^eme temps que sa
\textit{\~{P}LUP}-\decoz, \index{PLUP@\~{P}LUP-\deco}
le rang de la matrice \,$A$.

Notons aussi que la \met avec recherche du pivot
permet
de calculer dans tous les cas le \deter de la matrice
\,$A$\,
si elle est carr\'ee. Il suffit de garder en m\'emoire
et de mettre \`a jour \`a chaque
\'etape la parit\'e des
permutations de lignes et de colonnes d\'ej\`a
effectu\'ees.

\subsubsection*{\LUP-\deco d'une matrice surjective}
\refstepcounter{bidon}
\label{LUPdeco}\index{LUP@\LUP-\deco}
Un cas particulier est donn\'e par les matrices
surjectives.
Un pivot non nul existe toujours
sur la ligne voulue.
Cela donne l'\algo \vref{apgmatsur}.

\begin{algor}[\LUP-\deco d'une matrice
surjective.]\label{apgmatsur}
\acl{apgmatsur}{\LUP-\deco d'une matrice surjective}
\Entree Une matrice $A=(a_{ij})\in\K^{m\times n}$
surjective.
\Sortie La matrice $A$ transform\'ee (elle
donne les matrices $L$ et $U$ comme dans l'\algo
\ref{apgs}),
la matrice de permutation $P$ et sa signature
$e\in\{-1,1\}$.
\Debut
\Varloc $i$, $j$, $p\in\N$ ; $piv\in\K$ ;
\hsu $P:= \In$ ; $e:=1$ ;
\hsu \pour{p}{1}{m}
\hsd     $piv:=a_{pp}$ ; $j:=p$ ;
\hsd     \sialors{piv=0}
\hst         \tantque{piv=0}
\hsq              $j:=j+1$ ; $piv:= a_{pj}$ ;
\hst         \fintantque  ;
\hst         EchCol$(A,p,j)$ ; EchCol$(P,p,j)$ ;
$e:=-e$ ;
\hsq         \# EchCol$(A,p,j)$ est une \pcd qui
             \'echange
\hsq         \# les colonnes $p$ et $j$ de la matrice
$A$.
\hsd     \finsi ;
\hsd \pour{i}{p+1}{m}
\hst $a_{ip}:=a_{ip} / piv\;;$
\hst \pour{j}{p+1}{n} $a_{ij}:=a_{ij}- a_{ip}*a_{pj}$
\hst \finpour
\hsd \finpour
\hsu \finpour
\fin
\end{algor}

Ainsi lorsqu'une matrice \,$A\in {\cal K}^{m\times
n}$\,
est surjective (\cad si son  rang est \'egal
au nombre de ses lignes), on peut
la d\'ecomposer  en un produit de trois matrices
\,$L,\, U,\, P$\, o\`u \,$L\in {\cal K}^{n\times n}$\,
est une
\matgi avec des $1$ sur la
diagonale, \,$U\in {\cal K}^{m\times n}$\, une \matgu
et
\,$P\in {\cal K}^{n\times n}$\, une matrice de
permutation.

La \LUP-\deco permet de r\'esoudre des
probl\`emes comme le calcul du d\'eterminant ou la
r\'esolution d'un \syse \linsz.
En effet, pour r\'esoudre le \sys
\,$A\,x=b$\, avec \,$A=L\,U\,P$, on commence par
r\'esoudre le \sys \,$L\,z=b$\, puis le \sys
\,$U\,y=z$\, et enfin le \sys \,$P\,x=y$.
Les deux premiers \syss sont des \syss
\trgs que l'on peut r\'esoudre par substitutions
successives des inconnues (en \,$\O(n^2)$\, \oparis
donc)
et le dernier \sys est
une simple permutation des inconnues.
Enfin \,$\det\,{A}=\pm\,\det\,U$\, (selon
la parit\'e de la permutation repr\'esent\'ee par la
matrice \,$P$).

Il faut remarquer qu'une matrice non surjective
n'admet pas toujours de \LUP-\deco comme
par exemple la matrice \,$\left[\begin{array}{cc}
0 & 0 \\ 1 & 1 \end{array}\right]$.\\
Par ailleurs la   \LUP-\deco d'une
matrice surjective n'est pas unique
comme on peut le voir sur la matrice
\,$A=\left[\begin{array}{ccc} 2 & 1 & 1\\ 0 & 3 & 4
\end{array} \right]$\, qui admet les deux
\LUP-\decos \,$A=L\,U\,P$\, avec
\,$L=\mathrm{I}_2$, \,$P=\mathrm{I}_3$\, et \,$U=A$,
ou encore \,$A=L\,U\,P$\, avec
$$
\,L=\left[\begin{array}{cc}
1 & 0 \\ 3 & 1 \end{array} \right],\;\;
U=\left[\begin{array}{rrr}
1 & 2 & 1 \\ 0 & -6 & 1 \end{array} \right]
\;\; \mathrm{et} \;\;
P=\left[\begin{array}{ccc} 0 & 1 & 0 \\
1 & 0 & 0 \\ 0 & 0 & 1 \end{array}\right]\,.
$$
Notons enfin que la matrice \,$U$\, obtenue dans la
\deco \,$A=L\,U\,P$\, est une matrice
surjective et \fregz.

\subsubsection*{R\'esolution de \slis et calcul de
l'inverse}

La \mpg permet de r\'esoudre un ou plusieurs
\slis associ\'es \`a la m\^eme matrice,
en triangulant
la matrice \'elargie aux seconds membres.

Dans sa variante \gui{Gauss-Jordan}, qui consiste \`a
poursuivre le
processus d'\'elimination de Gauss \gui{de bas en
haut} et \gui{de
droite \`a
gauche} sur les lignes de la matrice \,$U$\, de
fa\c{c}on
\`a annuler les \elts au dessus de la diagonale
principale,  la \mpg sert \egmt
\`a calculer l'inverse d'une \maca  inversible
lorqu'on l'applique
\`a cette matrice
\'elargie (\`a droite) avec la matrice unit\'e de
m\^eme ordre,
moyennant un
co\^ut l\'eg\`erement \supe qui fait passer la
constante
dans \,$\O(n^3)$\, de
\,${2\over 3}$\, \`a \,${4\over 3}$.

\section{M\'ethode de \JB}
\label{secDJB}\label{sec bareiss}
\label{sec JorBar}
\index{Jordan-Bareiss!m\'ethode de}

La \mpg
est une \met de traitement automatique des
\syses \lins dont les
\coes et les inconnues sont dans un corps
donn\'e $\K$.
Cette \met fonctionne bien dans le cas de
matrices \`a \coes dans un corps fini
(et dans une moindre mesure, dans le
cas du corps \,$\QQ$). Mais hormis le cas des
corps finis, elle poss\`ede l'inconv\'enient majeur
de n\'ecessiter une simplification syst\'ematique
des fractions calcul\'ees si on ne veut pas voir
la taille des \coes exploser, ce qui
entra\^{\i}ne souvent un temps de calcul prohibitif,
par
exemple lorsqu'on travaille avec un corps de fractions
rationnelles \`a plusieurs variables. En outre cette
\met utilise des divisions et ne s'applique donc
pas si la matrice a ses \coes dans un anneau
arbitraire.

Nous allons voir dans cette section que la \met
connue aujourd'hui sous le nom de \gui{\met de
Bareiss}, qui
peut \^etre consid\'er\'ee comme une adaptation de
la \mpg classique,
permet dans une certaine mesure de pallier aux
inconv\'enients pr\'esent\'es par cette derni\`ere.

La \met de Bareiss (\cf \cite{Bar}, 1968) \'etait
connue
de Jordan (\cf \cite{Dur}), et elle semble avoir
\'et\'e
d\'ecouverte par Dodgson
(plus connu sous le nom de Lewis Caroll)
qui en a donn\'e une variante dans \cite{Dodgson}.
Nous la d\'esignerons d\'esormais sous
le nom  de \emph{\mjbz}.

Nous r\'eservons le nom
de \emph{\mddz} \`a la variante de Lewis Caroll que
nous
exposons \`a la fin de la section.

La \mjb est valable dans le cas d'un anneau
int\`egre \,${\cal A}$\, o\`u l'\egt peut
\^etre test\'ee par un \algoz, et l'addition,
la \mul et la division \gui{exacte} (quand
il y a un quotient exact) peuvent \^etre
effectu\'ees par des \algosz. Cela signifie,
pour la division exacte, qu'il y a un \algo
prenant en entr\'ee un couple
\,$(a,b)\in\A^2,~b\neq 0$, et donnant en sortie
l'unique \elt \,$x\in\A$\, v\'erifiant
\,$ax=b$, dans le cas o\`u il existe.

\subsection[Formule de \DJB]{Formule de \DJB et
variantes}
\label{subsec.fdjb}

Soit \,$A$\, une matrice dans \,$\A^{m\times n}$.
Reprenant les relations donn\'ees \`a
la propri\'et\'e \ref{propri Gauss}, et puisque tous
les
\coes $a_{\,ij}^{[p]}$ s'\'ecrivent
$a_{\,ij}^{(p)}/a_{~pp}^{(p-1)}$
(relation (\ref{EqGauss4})) avec le m\^{e}me \deno
pour
un $p$ fix\'e, l'id\'ee est de calculer
directement les \numes de mani\`ere r\'ecursive.
L'\'equation
(\ref{EqGauss1}) se relit
alors sous la forme
\begin{equation}\label{eqJB}
a_{\,ij}^{(p)} =
\frac{a_{~ij}^{(p-1)}\,a_{~pp}^{(p-1)}-
a_{~ip}^{(p-1)}\,a_{~pj}^{(p-1)}}{a_{p-1,p-1}^{(p-2)}}~
\end{equation}
C'est ce que nous appellerons la \emph{formule de
\DJBz}.
On peut obtenir ce m\^{e}me r\'esultat en
appliquant  l'\idt de Sylvester (\ref{eq3bis}) de la
proposition
\vref{syl} \`a la matrice
$$\left[\begin{array}{lr} A_p & A_{1..p,j} \\
A_{i,1..p} & a_{ij}\\
\end{array}\right]$$
avec \,$p\in [1..\min(m,n)-1]\,,~
i\in [p+1..m]\,,~j\in [p+1..n]\,. $
Cela donne:
%

\begin{proposition}
\label{formule J-B}\emph{(Formule de \DJBz)}
\index{Dodgson-Jordan-Bareiss!formule de}\\
Soit \,$\A$\, un \acom arbitraire.
Pour toute matrice \,$A=(a_{ij})\in\A^{m\times n}$,
on a la relation:
%
\begin{equation}\label{eq5}
a_{ij}^{(p)}\times a_{p-1,p-1}^{(p-2)} =
\left|\begin{array}{cc} a_{pp}^{(p-1)} &
a_{pj}^{(p-1)}
\\ a_{ip}^{(p-1)} & a_{ij}^{(p-1)}
\end{array}\right|
\end{equation}
avec les conventions \usles
$~a_{~00}^{(-1)}=1~$ et $~a_{ij}^{(0)}=a_{ij}$.
\end{proposition}

On a \egmt la variante suivante.
Si l'on applique la formule (\ref{eq3}) de
la proposition \ref{syl}  \`a la
matrice $~\left[\begin{array}{lr} A_p & A_{1..p,j} \\
A_{i,1..p} & a_{ij}\\ \end{array}\right]$, on obtient:

\begin{proposition}
\label{propB2} \emph{(Formule de Bareiss \`a plusieurs
\'etages)}\\
Soit \,$\A$\, un
\acom arbitraire. Pour toute matrice
\,$A\in\A^{m\times n}$\, et tout entier \,$p\geq 2$,
on a lorsque \,$1\leq r\leq p-1,$ \,$p+1\leq i\leq m,$
et \,$p+1\leq j\leq n$:
\begin{equation}\label{eq4}
\left(a_{rr}^{(r-1)}\right)^{p-r}
    \,a_{ij}^{(p)}
\: = \left|\begin{array}{cccc}
a_{r+1,r+1}^{(r)} & \ldots  & a_{r+1,n}^{(r)} &
a_{r+1,j}^{(r)} \\
\vdots & \ddots & \vdots & \vdots \\
a_{p,r+1}^{(r)} & \ldots  & a_{p,p}^{(r)} &
a_{p,j}^{(r)}
\\ a_{i,r+1}^{(r)} & \ldots  & a_{i,p}^{(r)} &
a_{i,j}^{(r)}
\end{array}\right|.
\end{equation}
\end{proposition}

Dans son article, Bareiss a remarqu\'e qu'on pouvait
utiliser cette identit\'e avec \,$p-r=2$\, pour
calculer les \,$a_{ij}^{(p)}$\, de proche en proche,
lorsque l'anneau est int\`egre et poss\`ede un \algo
de division exacte.

En fait la \gui{\met de Bareiss}
couramment utilis\'ee aujourd'hui est plut\^ot bas\'ee
sur
la premi\`ere formule (celle de \DJBz).
L'\'equation (\ref{eqJB}) permet en effet de calculer
les \,$a_{ij}^{(p)}$\, de proche en proche.

La \mjb est donc une
adaptation de la \mpg qui garantit,
tout au long du processus de triangulation  de la
matrice
trait\'ee, l'appartenance des \coes \`a l'anneau
de base. L'efficacit\'e de cet \algo tient \`a ce que
les
\coes calcul\'es sont tous des \deters
extraits de la matrice initiale, et donc restent de
taille
raisonnable pour la plupart des anneaux \uslsz.

\subsubsection*{Algorithme de \JB}

En utilisant la relation (\ref{eqJB}) on obtient
l'\ajb
\vref{ajb},
dans sa version de \emph{l'\'elimination \`a un seul
\'etage}
d\'echarg\'ee de la recherche du pivot.

Rappelons les conventions
\,$a_{\,00}^{(-1)}=1$\, et \,$a_{~pp}^{(p-1)}=0$ pour
\,$p>\inf(m,n)$.

\begin{algor}[Algorithme de Jordan-Bareiss]\label{ajb}
\acl{ajb}{Algorithme de \JBz}
\Entree Une matrice $\,A=(a_{ij})\in\A^{m\times n}$.
L'anneau $\A$ est suppos\'e int\`egre avec un \algo
de division exacte.
\Sortie La matrice $\,A\,$ transform\'ee.
Si les $r$ premiers \mips dominants sont non nuls,
et si le $(r+1)$-\`eme est nul, elle contient en
position
\,$(i,j)\,$ le mineur $\,a^{(p)}_{i,j}\,$ avec
\,$p=\inf(r,i-1,j-1)$. L'entier $r$ est aussi
calcul\'e.
Si en outre $\,r=\rg(A)\,$ on retrouve facilement la
\LU-\deco de $\,A\,$ \`a partir de la sortie, comme
expliqu\'e avant l'exemple \ref{exaJB1}.
\Debut
\Varloc $i$, $j$, $p\in\N\;$; $piv$, $den$,
$coe\in\A\;$;
\hsu $p:=1$ ; $den:= 1$ ;  $r:=\inf(m,n)$ ;
\hsu \tantque{p < \inf(m,n)}
\hsd    $piv:=a_{pp}$ ;
\hsd    \sialors{piv= 0}   $p:=\inf(m,n)$ ; $r:=p-1\;$
\sinon
\hst      \pour{i}{p+1}{m}
\hsq         $coe:=a_{ip}$;
\hsq         \pour{j}{p+1}{n}
\hsc            $a_{ij}:=(piv*a_{ij}-coe*a_{pj})\,/ \,
den$
\hsq         \finpour
\hst      \finpour
\hsd    \finsi ;
\hsd    $p:= p+1$ ;
\hsd    $den:=piv$
\hsu \fintantque
\fin
\end{algor}

On retrouve facilement la \LU-\deco de \,$A$\,
\`a partir de la matrice retourn\'ee par l'\algo
pr\'ec\'edent
en utilisant les formules (\ref{EqGauss4})
(propri\'et\'e \ref{propri Gauss})
et (\ref{EqGaussL2}) page \pageref{EqGaussL2}: notons
\,$c_{ij}$\, les
\coes de cette
matrice\,; alors
pour la matrice \,$L$\, on a
\,$l_{ij}=c_{ij}/c_{jj}$\, si \,$1\leq
j<i\leq m$\, ($l_{ij}=\delta_{ij}$\, sinon) et pour la
matrice  \,$U$\,
on a \,$u_{ij}=c_{ij}/c_{i-1,i-1}$\, si \,$1\leq
i\leq j$\, ($u_{ij}=0$\, sinon).
On peut le voir sur l'exemple suivant.

\begin{example}
\label{exaJB1} {\footnotesize
\emph{Dans cet exemple on reprend la matrice \,$M_3$\,
de l'exemple
\ref{exaGauss2} et on donne ses transform\'ees par les
\algos de Jordan-Bareiss  et de Gauss.
$$
M3\,=\,\left[
{\begin{array}{rrrrr}
-73 & -53 & -30 & 45 & -58 \\
21 & -54 & -11 & 0 & -1 \\
72 & -59 & 52 & -23 & 77 \\
33 & 55 & 66 & -15 & 62 \\
-41 & -95 & -25 & 51 & -54 \\
14 & 55 & 35 & -5 & 25
\end{array}}
 \right]
$$
$$
 \left[
{\begin{array}{rrrrr}
-73 & -53 & -30 & 45 & -58 \\
21 & 5055 & 1433 & -945 & 1291 \\
72 & 8123 & 272743 & 2940 & 243716 \\
33 & -2266 & 220594 & 2911532 & -3038698 \\
-41 & 4762 & 52277 & 3660124 & 0 \\
14 & -3273 & 83592 & 3227536 & 0
\end{array}}
 \right]
$$
$$
 \left[
{\begin{array}{ccccc}
-73 & -53 & -30 & 45 & -58 \\
{\ds \frac {-21}{73}}  & {\ds \frac {-5055}{
73}}  & {\ds \frac {-1433}{73}}  & {\ds
\frac {945}{73}}  & {\ds \frac {-1291}{73}}  \\ [2ex]
{\ds \frac {-72}{73}}  & {\ds \frac {8123}{
5055}}  & {\ds \frac {272743}{5055}}  & {\ds
\frac {196}{337}}  & {\ds \frac {243716}{5055}}  \\
 [2ex]
{\ds \frac {-33}{73}}  & {\ds \frac {-2266}{
5055}}  & {\ds \frac {220594}{272743}}  &
{\ds \frac {2911532}{272743}}  & {\ds \frac {
-3038698}{272743}}  \\ [2ex]
{\ds \frac {41}{73}}  & {\ds \frac {4762}{
5055}}  & {\ds \frac {52277}{272743}}  &
{\ds \frac {1193}{949}}  & 0 \\ [2ex]
{\ds \frac {-14}{73}}  & {\ds \frac {-1091}{
1685}}  & {\ds \frac {83592}{272743}}  &
{\ds \frac {1052}{949}}  & 0
\end{array}}
 \right]
$$
}
}
\end{example}

Comparons l'\ajb \`a l'\apg
dans le cas de l'anneau \,$\A=\ZZ[X,Y]$.

Lorsqu'on utilise l'\apg dans le corps des fractions
\,$\FA=\QQ(X,Y)$\, sans r\'eduire les fractions au
fur et \`a mesure qu'elles sont calcul\'ees (ce qui
est
tr\`es co\^uteux), il n'est pas difficile de voir
que les degr\'es des \numes et \denos ont
un comportement exponentiel. Avec l'\ajbz, par contre,
les degr\'es ont seulement une croissance \linz.

\bni
\rem
Dans le \emph{cas non int\`egre\/}, le
fonctionnement de l'\ajb
sans recherche du pivot reste possible si tous les
mineurs principaux rencontr\'es \,$a_{pp}^{(p-1)}$\,
au
cours du processus de triangulation sont non diviseurs
de 0, et si \emph{les divisions exactes:}
$$
\frac{a_{pp}^{(p-1)}a_{ij}^{(p-1)}-
a_{ip}^{(p-1)}a_{pj}^{(p-1)}}{a_{p-1,p-1}^{(p-2)}}
~~~~\mathit{peuvent~se~faire~algorithmiquement}\,.$$
Cette condition est satisfaite lorsqu'on remplace la
matrice
carr\'ee \,$A$\, par sa matrice \cara
\,$A-X\I_n \in\A\,[X]^{n\times n}$\, o\`u tous les
pivots
rencontr\'es sont des \pols unitaires (au signe
pr\`es).
Donc l'\ajb appliqu\'e \`a
\,$A-X\I_n$\, ne fait intervenir que  la structure
d'anneau
de \,$\A$\, et ne n\'ecessite aucune division dans
\,$\A$.
C'est l'objet du paragraphe suivant.

\subsection[M\'ethode de \JB modifi\'ee]{Cas d'un
\acoma: \mjb
modifi\'ee} \label{subsec.mjbmo}
\index{Jordan-Bareiss!m\'ethode modifi\'ee }

C'est la \mjb appliqu\'ee
\`a la matrice \cara \,$A-X\I_n$\, d'une
\maca \,$A\in \A^{n\times n}$.
Les \coes de \,$A-X\I_n$\, sont dans l'anneau
\,$\A\,[X]$. M\^eme si \,$\A$\, n'est pas int\`egre,
les divisions exactes requises sont ici des divisions
par
des \pols unitaires qui ne n\'ecessitent par
cons\'equent aucune division dans \,$\A$, mais
uniquement
des additions, soustractions, et \mulsz.
En parti\-culier, aucune permutation de lignes ou de
colonnes n'intervient au cours du processus de
triangulation.

La \mjb
modifi\'ee permet donc de calculer  le \polcar
 de la matrice \,$A$, et par
cons\'equent son \deterz, son adjointe, et, au
cas o\`u elle est inversible, son inverse.

Cette \met a \'et\'e propos\'ee en 1982 par
Sasaki \& Murao \cite{Sas}.
Les auteurs remarquent \egmt que dans un
calcul de base de l'\algo (du type \gui{produit en
croix divis\'e par le pivot pr\'ec\'edent}):
$$
f(X):=\frac{a(X)\,c(X)-b(X)\,d(X)}{e(X)}\,,
$$
les degr\'es en \,$X$\, sont \'egaux \`a \,$k$\, ou
\,$k+1$\, pour \,$f,$
\`a \,$k$\, ou
\,$k-1$\, pour \,$a,b,c,d$\, et \`a \,$k-1$\, pour
\,$e$.
On peut donc se passer de calculer les \coes
des mon\^omes de degr\'e \,$<k-1$\, dans \,$ac-bd$\,
et le calcul du quotient ne doit pas non plus
s'encombrer
des termes de degr\'es \,$<k-1$\, dans les restes
successifs
(pour l'\algo \usl de division des \polsz).
Ceci conduit pr\'ecis\'ement aux r\'esultats suivants.
\begin{itemize}
\item Les \coes des mon\^{o}mes de degr\'e \,$k-1$\,
\`a \,$2k$\, dans le
produit de deux \pols de degr\'e \,$k$\, se calculent
(en utilisant la
\met \uslez) en \,$k^2+2\,k-2$\, \oparisz.
\item  La division exacte d'un \pol de degr\'e
\,$2k$\, par un \pol
unitaire de degr\'e \,$k-1$\, se calcule (en utilisant
la \met \uslez)
en \,$k^2+3\,k-1$\, \oparisz.
\end{itemize}
On en d\'eduit qu'une affectation
\,$f:=\frac{a\,d-b\,c}{e}$\, dans
l'\ajb modifi\'e, lorsque \,$e=e(X)$\, est le pivot
unitaire de degr\'e
\,$k-1$, consomme \,$3\,k^2+\O(k)$\, \oparis
dans l'anneau de base  (et en tout cas au plus
\,$3\,k^2+8\,k-4$).

Pour l'ensemble de l'\algo on obtient un nombre total
d'\oparis
inf\'erieur \`a
$$\sum_{k=2}^{n} \,(3\,k^2+8\,k-4)\,(n-k+1)^2\,\leq\,
{\frac {1}{10}} \,n^{5} + {\frac {7}{6}} \,n^{4} +
{\frac {7}{3}}
\,n^{3}\,.$$

\begin{proposition}
\label{propJBC} Soit \,$A\in\A^{n{\times}n}$\, une
\maca sur un \acomaz.
L'\ajb appliqu\'e \`a la matrice \cara \,$A-X\,\In$\,
s'ex\'ecute
(en utilisant la \met \uslez) en \,${\frac {1}{10}}
\,n^{5} + \O(n^4)$\,
\oparis dans l'anneau \,$\A$.
\end{proposition}

\subsection{La \mdd}\label{subsec.mdd}

La \mdd est une variante \'el\'egante et \sym de la
\mjbz.
Cependant son but n'est pas le calcul de la \LU-\deco
d'une matrice,
mais seulement celui de ses \emph{mineurs connexes},
\index{mineur!connexe}\cad les mineurs
\,$a_{i..j}^{h..k}$\, (avec
\,$j-i=k-h$).
En particulier elle peut \^{e}tre utilis\'ee pour
le calcul du \deter d'une \macaz.

Une variante de la formule \ref{eq5} (apr\`es un
\'echange de lignes
et un \'echange de colonnes)
est la formule suivante concernant les mineurs
connexes
\begin{equation} \label{EqLewisCaroll}
a_{i+1..j-1}^{h+1..k-1}\cdot a_{i..j}^{h..k} =
\left\vert\matrix{
 a_{i..j-1}^{h..k-1}   &     a_{i..j-1}^{h+1..k}     
\cr
 a_{i+1..j}^{h..k-1}   &     a_{i+1..j}^{h+1..k}
}\right\vert
\end{equation}
Cela donne les affectations correspondantes dans
l'\agdz.
Mais ce dernier fonctionne uniquement si tous les
mineurs connexes
appel\'es \`a servir de \deno sont non nuls:
contrairement \`a la \mpg
et \`a la \mjbz, la \mdd ne poss\`ede pas de variante
connue
efficace dans le cas o\`u une
affectation  \,$x:=0/0$\, est produite par
l'\algoz({\footnote{~Lewis
Caroll propose dans sa communication d'op\'erer des
permutations circulaires sur les lignes et les
colonnes de
la matrice. Voici un contre-exemple montrant que la
\mdd ne
s'applique pas toujours. La matrice
$
A=\left[ \begin{array}{cccc} 1 & 0 & 1 & 1 \\
1 & 0 & 1 & 0 \\ 0 & 1 & 0 & 0 \\ 1 & 0 & 0 & 1
\end{array} \right]
$
est une matrice inversible
de d\'eterminant \,$1$, lequel ne peut pas
se calculer par la \met
de Lewis Carrol, m\^eme lorsqu'on effectue des
permutations
circulaires de lignes et de colonnes.}}).

Pour voir plus clairement ce que signifie
l'\'equation de Lewis Caroll (\ref{EqLewisCaroll})
appelons \,$B$\,  la
matrice extraite \,$A_{i+1..j-1,h+1..k-1,}$\, et 
notons
\,$p$, \,$q$, \,$u$, \,$v$ les indices
\,$i+1$, \,$j-1$, \,$h+1$, \,$k-1$. L'\'equation se
r\'e\'ecrit alors:

{\footnotesize
$$
\left\vert\matrix{
 B
}\right\vert
\cdot
\left\vert\matrix{
a_{i,h}    &  A_{i,u..v}   &  a_{i,k}    \cr
A_{p..q,h}    &  B         &   A_{p..q,k}  \cr
a_{j,h}    &  A_{j,u..v}   &  a_{j,k}
}\right\vert
=\quad \quad \quad \quad \quad \quad
\quad \quad \quad \quad \quad $$
$$
\left\vert\matrix{
a_{i,h}    &  A_{i,u..v}  \cr
A_{p..q,h}    &  B
}\right\vert
\cdot
\left\vert\matrix{
 B         &   A_{p..q,k}  \cr
 A_{j,u..v}   &  a_{j,k}
}\right\vert
-
\left\vert\matrix{
A_{p..q,h}    &  B          \cr
a_{j,h}    &  A_{j,u..v}
}\right\vert
\cdot
\left\vert\matrix{
  A_{i,u..v}   &  a_{i,k}          \cr
  B         &   A_{p..q,k}
}\right\vert
$$
Un exemple:
$$
\left\vert\matrix{
  c_2  &   c_3        \cr
  d_2  &   d_3
}\right\vert
\cdot
\left\vert\matrix{
b_1    &  b_2   &  b_3  &  b_4 \cr
c_1    &   c_2  &  c_3  &  c_4 \cr
d_1    &  d_2   &  d_3  &  d_4 \cr
e_1    &  e_2   &  e_3  &  e_4
}\right\vert
=\quad \quad \quad \quad \quad \quad
\quad \quad \quad \quad \quad $$
$$
\left\vert\matrix{
 b_1   &   b_2  &   b_3   \cr
 c_1   &   c_2  &   c_3   \cr
 d_1   &   d_2  &   d_3
}\right\vert
\cdot
\left\vert\matrix{
  c_2  &  c_3   &  c_4    \cr
  d_2   &   d_3 &   d_4   \cr
 e_2   &  e_3 &   e_4
}\right\vert
-
\left\vert\matrix{
c_1    &   c_2  &  c_3 \cr
d_1    &  d_2   &  d_3 \cr
e_1    &  e_2   &  e_3
}\right\vert
\cdot
\left\vert\matrix{
b_2   &  b_3  &  b_4 \cr
 c_2  &  c_3  &  c_4 \cr
d_2   &  d_3  &  d_4
}\right\vert
$$
}

\ss Dans la \mjb sans recherche de pivot on calcule
\`a l'\'etape \num$p$ tous les mineurs
\,$a^{(p)}_{ij}$\, $((i,j>p)$ d'une
matrice \,$A$. Dans la \mpg on calcule les quotients
\,$a^{[p]}_{ij}=a^{(p)}_{ij}/a^{(p-1)}_{p,p}.$
Si la matrice a une \gui{structure interne}
comme dans le cas des matrices de Hankel
ou de Toeplitz la structure est perdue
d\`es la premi\`ere \'etape.

Dans la \mddz, on calcule
\`a l'\'etape \num$p$ tous les mineurs connexes
d'ordre \,$p+1$\,
de la matrice \,$A$.
Il s'ensuit que dans le cas d'une matrice
structur\'ee,
les matrices \itmds calcul\'ees
par la \mdd sont \egmt structur\'ees. Ceci
diminue tr\`es s\'erieusement le nombre d'\oparis \`a
effectuer et le fait passer de $\O(n^3)$ \`a
$\O(n^2)$.

Dans le cas d'un anneau int\`egre o\`u les divisions
exactes
sont faisables par un \algoz, on obtient les m\^{e}mes
avantages que dans l'\ajb concernant la taille des
\coes \itmdsz.

\subsubsection*{Algorithme de Dodgson pour une matrice
de Hankel}

Nous donnons ici une version pr\'ecise de l'\algo de
Dodgson pour
les matrices de Hankel dont tous les mineurs connexes
sont non
nuls. C'est l'\algo \vref{agdh}. \\ L'entr\'ee est une
liste
\,$L=(a_i)$\, contenant les $\,m+n-1\,$ \coes de la
matrice de
Hankel initiale \,$H\in\A^{m\times n}$\, 
($h_{i,j}=a_{i+j-1}$).
La sortie est un tableau  \,$T=(t_{r,j})$\,
($r=0,\ldots,\inf(m,n)$, $j=r,\alb\ldots ,\alb m+n-r$)
qui
contient tous les mineurs connexes de la matrice 
\,$H$,
calcul\'es en suivant l'\algo de Dodgson. Pour
l'initialisation,
sur la ligne $0$ il y a des $1$ (les \gui{mineurs
connexes d'ordre
0}) et sur le ligne $1$ les \coes de \,$H$\, (les
\gui{mineurs
connexes d'ordre 1}). Sur la ligne \,$r\geq 2$\, il y
a les \coes
de la matrice de Hankel form\'ee par les mineurs
connexes d'ordre
\,$r$\, de  \,$H$. Dans la colonne \,$j$\, il y a les
\deters des
sous-\macas de  \,$H$\, qui ont le \coe \,$a_j$\, sur
leur
diagonale ascendante.

\begin{algor}[Algorithme de Dodgson pour une matrice
    de Hankel]\label{agdh}
\acl{agdh}{Algorithme de Dodgson pour une matrice de
Hankel}
\Entree Deux entiers \,$m,n\in\N$ et une liste
\,$L=(a_{i})\in\A^{m+n-
1}$.
Cette liste contient les \coes d'une matrice de Hankel
$\,H\in\A^{m\times n}$.
L'anneau $\A$ est suppos\'e int\`egre avec un \algo
de division exacte.
\Sortie Un tableau $T=(t_{r,j})$ rempli
d'\elts de $\A$ pour
$r\in\left\{0,\ldots,\inf(m,n)\right\}$,
$j\in\left\{r,\alb\ldots ,\alb m+n-r\right\}$.
Il contient sur la ligne \,$r$\, les mineurs connexes
d'ordre \,$r\,$ de la matrice  $\,H$, suppos\'es
tous non nuls.
\Debut
\Varloc $r$, $j$, $q\in\N$ ;
\hsu $q:=\inf(m,n)$ ;
\hsu  $T:=$TableauVide$(0..q,1..m+n-1)$ ;
\hsq     \# on a cr\'e\'e $\,T\,$ tableau vide de
taille voulue
\hsu \pour{j}{1}{m+n-1}
\hsd    $t_{0,j}:=1$ ; $t_{1,j}:=a_{j}$  ;
\hsu \finpour;
\hsq     \# fin de l'initialisation
\hsu \pour{r}{1}{q-1}
\hsd   \pour{j}{r+1}{m+n-r-1}
\hst      $
t_{r+1,j}:=(t_{r,j-1}\,t_{r,j+1}-t_{r,j}^2)/t_{r-1,j}$
\hsd   \finpour
\hsu \finpour
\fin
\end{algor}

L'\algo est pratiquement le m\^{e}me dans le cas d'une
\mto
\,$Z$\, (il suffit de changer le signe dans
l'affectation de
 \,$t_{r+1,j}$\,) et il peut s'appliquer pour le
calcul du
\polcarz.

\begin{examples}~\\
\label{exaJBD}{\footnotesize\emph{%
Dans le premier exemple, on consid\`ere la matrice de
Hilbert
d'ordre $5$, qui est un exemple classique de matrice
de Hankel
mal conditionn\'ee (le \deter de la matrice,
$1/266716800000$, est
l'inverse d'un entier tr\`es grand).\\
\[
\mathit{A} =  \left[
{\begin{array}{cccccc}
     1 &   1/2  &   1/3  &   1/4  &   1/5   \\
  1/2  &   1/3  &   1/4  &   1/5  &   1/6       \\
  1/3  &   1/4  &   1/5  &   1/6  &   1/7       \\
  1/4  &   1/5  &   1/6  &   1/7  &   1/8    \\
  1/5  &   1/6  &   1/7  &   1/8  &   1/9
\end{array}}
 \right]
\]
Voici alors la sortie de l'\agd \ref{agdh}
(on a supprim\'e la ligne des \,$1$):
\[
{\begin{array}{ccccccccc}
1 & {\displaystyle \frac {1}{2}}  & {\displaystyle
\frac {1}{3}}
 & {\displaystyle \frac {1}{4}}  & {\displaystyle
\frac {1}{5}}
 & {\displaystyle \frac {1}{6}}  & {\displaystyle
\frac {1}{7}}
 & {\displaystyle \frac {1}{8}}  & {\displaystyle
\frac {1}{9}}
 \\ [2ex]
\  & {\displaystyle \frac {1}{12}}  & {\displaystyle
\frac {1}{72
}}  & {\displaystyle \frac {1}{240}}  & {\displaystyle
\frac {1}{
600}}  & {\displaystyle \frac {1}{1260}}  &
{\displaystyle
\frac {1}{2352}}  & {\displaystyle \frac {1}{4032}}  &
\  \\
 [2ex]
\  & \  & {\displaystyle \frac {1}{2160}}  &
{\displaystyle
\frac {1}{43200}}  & {\displaystyle \frac {1}{378000}}
 &
{\displaystyle \frac {1}{2116800}}  & {\displaystyle
\frac {1}{
8890560}}  & \  & \  \\ [2ex]
\  & \  & \  & {\displaystyle \frac {1}{6048000}}  &
{\displaystyle \frac {1}{423360000}}  & {\displaystyle
\frac {1}{
10668672000}}  & \  & \  & \  \\ [2ex]
\  & \  & \  & \  & {\displaystyle \frac
{1}{266716800000}}  & \
 & \  & \  & \
\end{array}}
\]
Voici ensuite un exemple de la sortie
de l'\algo avec une matrice de Hankel car\-r\'ee
d'ordre $7$
\`a \coes entiers (lisibles sur la premi\`ere ligne):
\[
{\begin{array}{ccccccrcccccc}
1 & 7 & 7 & 1 & 2 & 2 & 4 & 3 & 5 & 3 & 7 & 2 & 4 \\
  & -42 & -42 & 13 & -2 & 4 & -10 & 11 & -16 & 26 &
-43 & 24 &
  \\
  &   & -330 & -85 & 24 & 2 & -14 & 13 & 6 & 4 & -175
&   &
 \\
  &   &   & -1165 & 373 & -85 & 17 & -23 & -1 & -41 & 
 &
 &   \\
  &   &   &   & -1671 & -442 & -119 & -42 & 157 &   & 
 &
  &   \\
  &   &   &   &   & -41 & 259 & 889 &   &   &   &   &
 \\
  &   &   &   &   &   & 870 &   &   &   &   &   &
\end{array}}
\]
}}
\end{examples}

\section{M\'ethode de Hessenberg}\label{sec hessen}

Toutes les matrices consid\'er\'ees ici sont
\`a \coes dans un corps commutatif \,$\K$.

\subsubsection*{Matrices quasi-triangulaires}

\begin{defi}  Une \maca
\,$H=(h_{ij})\in \K^{n\times n}~(n\in
\NN^*)$\, est dite \qtgu (resp. \qtgiz) si
$~h_{ij}= 0 \, \mbox{ d\`es que }\,i-j \geq
2$\, (resp. d\`es que \,$j-i \geq 2$). On dit
encore que $H$ est une matrice de Hessenberg.
\end{defi}

Une matrice \qtgu
$H$ est donc une matrice de la forme:
\begin{equation} \label{hess1}  H =
\left[\begin{array}{ccccc} h_{11} & h_{12} &
\ldots &
 h_{1,n-1} & h_{1n} \\ h_{21} & h_{22} &
\ldots & h_{2,n-1} & h_{2n} \\ 0 & h_{32} &
\ddots & & h_{3n} \\
\vdots & \ddots & \ddots & \ddots &  \vdots\\
0  & \ldots& 0 & h_{n,n-1} & h_{nn} \\
\end{array}\right]
\end{equation}

On d\'emontre par une  \recu
\imme sur \,$k\ (1\leq k\leq n)$\, la
propri\'et\'e suivante des matrices de
Hessenberg:

\begin{prop}
Soit \,$H=(h_{ij})$\, une matrice
de Hessenberg (\supee ou \infeez).
On d\'esigne par \,$H_k\  (1\leq k\leq n)$\,
la \smpd d'ordre
\,$k$\, de \,$H$\, et par \,$D_k$\, le
\deter de \,$H_k$. On pose \,$D_0=1$.
La suite \,$(D_k)_{1\leq k\leq n}$\,
(des \mips dominants de \,$H$)
v\'erifie alors la relation de \recuz:
$$
D_k=h_{kk}D_{k-1} +\sum_{i=1}^{k-1}
(-1)^{k-i} h_{k,k-1}h_{k-1,k-2}\ldots
h_{i+1,i}\, h_{ik}\,D_{i-1}\,\cdot
$$
\end{prop}

Pour le voir, il suffit de d\'evelopper
\,$D_k$\, suivant la derni\`ere ligne (resp.
la derni\`ere colonne) de \,$H_k$\, si
celle-ci est une matrice quasi-\tgu (resp. \infeez).

\ms Appliquant ce r\'esultat \`a la
matrice \,$H-X\I_n$, elle-m\^eme
\qtgz, dont les \mips dominants
sont les \polcars \,$P_k(X)$\, des \smpds \,$H_k$\, de
\,$H~(1\leq k\leq n)$, on obtient les relations
de \recu suivantes dites \emph{relations
de Hessenberg} permettant de calculer de
proche en proche les \polcars
\,$P_k(X)$\, de
\,$H_k$\, pour \,$2\leq k\leq n$\,
sachant que \,$P_0(X)=1$,
\,$P_1(X)=h_{11}-X$; et
\begin{equation} \label{hess2}
P_k(X)=\left\{ \,
\begin{array}{l}
(h_{kk}-X)\,P_{k-1}(X) \;+\\[1mm]
\sum_{i=1}^{k-1}
{\left(\left[\prod_{j=i+1}^k(-h_{j,j-1})\right]
\,h_{ik}\, P_{i-1}(X)
\right)}
\end{array}
\right.
\end{equation}

\subsubsection*{La \met de Hessenberg}

Elle consiste \`a calculer le \polcar  d'une \maca
\,$A$\, d'ordre \,$n \geq 2 $, dont les
\elts
\,$a_{ij}$\, appartiennent \`a un corps
\,$\K$, en la r\'eduisant \`a la forme
\,(\ref{hess1}) c'est-\`a-dire \`a une
matrice de Hessenberg \,$H$\, semblable \`a
\,$A$\, dont les \elts \,$h_{ij}$\,
appartiennent \egmt \`a \,$\K$.

\begin{algorH}[\Algo de Hessenberg ($\K$ est un
corps)]
\label{algohessen}
\acl{algohessen}{\Algo de Hessenberg}
\index{Hessenberg!algorithme de}
\Entree Un entier $n\geq 2$ et une matrice
$A=(a_{ij})\in\K^{n\times n}$.
\Sortie Le \polcar  de $A$: $\rPA(X)$.
\Varloc  $jpiv$, $ipiv$, $iciv$, $i$, $m\in\N$ ; $piv,
c, \in \K$ ; \\
$H:=(h_{{ij}})\in \K^{n\times n}:$
les matrices transform\'ees successives de $A$; \\
$P=(P_{i}):$ liste des \polcars successifs dans
$\K[X]$ ;
\Debut
\hsu $P_{0} :=1$ ; $H:=A$ ; \hspace{1cm} \#
Initialisations
\hsc \# R\'eduction de $H$ \`a la forme de Hessenberg
\hsu \pour{jpiv}{1}{n-2}
\hsd $ipiv:=jpiv+1$ ; $iciv:=ipiv$ ;
$piv:=h_{iciv,jpiv}$ ;
\hsd \tantque{piv=0\; \ET \; iciv<n}
\hst $iciv:=iciv+1$ ; $piv:=h_{iciv,jpiv}$
\hsd \fintantque ;
\hsd \sialors{piv \neq 0}
\hst \sialors{iciv>ipiv}
\hsq EchLin$(H,ipiv,iciv)$ ; \# Echange de lignes
\hsq EchCol$(H,ipiv,iciv)$ \,  \# Echange de colonnes
\hst \finsi ;
\hst \pour{i}{iciv+1}{n}
\hsq $c:=h_{i,jpiv}/piv$ ;
\hsq AjLin$(H,ipiv,i,-c)$ ; \# Manipulation de lignes
\hsq AjCol$(H,i,ipiv,c)$ \,\,\,\,\, \# Manipulation de
colonnes
\hsd \finsi
\hsu \finpour ;
\hsc \# Calcul du \polcar
\hsu \pour{m}{1}{n}
\hsd $P_{m}:=(h_{mm}-X)\cdot P_{m-1}$ ; $c:=1$ ;
\hsd \pour{i}{1}{m-1}
\hst  $c:=-c\cdot h_{m-i+1,m-i}$ ;
      $P_{m}:=P_{m}+c\cdot h_{m-i,m}\cdot P_{m-i-1}$
\hsd \finpour
\hsu \finpour ;
\hsu $\rPA(X):=P_{n}(X)$ \hspace{1cm} \# le \polcar de
$A$.
\fin
\end{algorH}

\ss Pour ce faire, on applique la \mpg aux lignes de
la matrice
donn\'ee \,$A$\, en prenant comme pivots les
\elts sous-diagonaux de la matrice
trait\'ee, et en prenant bien soin
d'effectuer les transformations \gui{inverses} sur
les colonnes de \,$A$\, pour que la matrice
et sa transform\'ee soient semblables.

\ss
Plus pr\'ecis\'ement, l'\'etape \,$p~(1\leq
p\leq n-2)$\, consiste tout d'abord
\`a voir si l'\elt en position
\,$(p+1,p)$\, est nul, auquel cas il faut
chercher un \elt non  nul au-dessous
de lui (sur la colonne $p$): si un tel
\elt n'existe pas, on passe \`a
l'\'etape suivante \,$p+1$. Sinon par une
permutation de lignes, on  ram\`ene le pivot
non nul au bon endroit, \cad \`a la position
\,$(p+1,p)$, ce qui revient \`a multiplier
\`a gauche la matrice trait\'ee par la
matrice de permutation
\,$E_{i,p+1}~(i> p+1)$\, obtenue en
permutant les lignes $i$ et $p+1$ (ou les
colonnes $i$ et $p+1$, ce qui revient au
m\^eme) de la matrice unit\'e \,$I_n$.
On multiplie
\`a droite par la m\^eme matrice de
permutation (qui est ici \'egale
\`a son inverse) afin que la matrice obtenue
\`a l'issue de chaque \'etape reste semblable
\`a la matrice de d\'epart.

\ss  Le pivot non nul \'etant alors au bon
endroit, on ach\`eve l'\'etape \,$p$\, en
utilisant ce pivot pour faire appara\^{\i}tre des
z\'eros au-dessous de lui dans sa colonne, ce
qui revient \`a multiplier \`a gauche la
matrice trait\'ee par une matrice du type:

\begin{samepage}
$$\begin{array}{ccccccc}
\;~~ & \; & ~p         & \,~p+1 &   &   &    \\
 & & ~\downarrow & \,~\downarrow &~~~\; & ~~~\; &
~~~~~~~~~~~
\end{array} $$
$$\,L = \left[\begin{array}{ccccccc}  1 &
\ldots & 0 & 0 & 0 & \ldots & 0 \\
\vdots & \ddots & \vdots & \vdots & \vdots &   &
\vdots \\
0 & \ldots & 1 & 0 & 0 &
\ldots & 0 \\  0 & \ldots & 0 & 1 & 0 &
\ldots & 0 \\ 0 & \ldots & 0 & l_{p+2,p+1} &
1 & \ldots & 0 \\
\vdots & \ddots & \vdots & \vdots & \vdots &
\ddots & \vdots \\ 0 & \ldots & 0 & l_{n,p+1}
& 0 & \ldots & 1
\end{array}\right]\,
\begin{array}{l}
\; \\ \; \\ \aff~p \\ \aff~p+1 \\ \aff~p+2
\\ \; \\ \;
\end{array} $$
\end{samepage}

\sni et multiplier ensuite  \`a droite
la matrice trait\'ee par la matrice \,$L^{-1}$\,
(obtenue \`a partir de \,$L$\, en changeant le signe
des \elts
sous la diagonale).

Il est clair que ces op\'erations, qui
d\'efinissent l'\'etape \,$p$\, et qui  sont
effectu\'ees sur la matrice provenant de
l'\'etape pr\'ec\'edente (appelons-la
\,$A^{(p-1)}$), n'affectent pas les
\,$p-1$\, premi\`eres colonnes de
\,$A^{(p-1)}$\, et donnent une matrice
\,$A^{(p)}$\, semblable \`a \,$A^{(p-1)}$.

Ceci donne l'\aghbz:
Primo, calculer, \`a l'aide
de la \pcd d\'ecrite ci-dessus,  une
matrice de Hessenberg \,$H$\, semblable \`a la
matrice donn\'ee \,$A\in \K^{n\times n}$.
Secundo, calculer le
\polcar de \,$H$\, (qui est aussi celui de \,$A$)
en utilisant les relations de Hessenberg
(\ref{hess2}).

On obtient ainsi l'\algo
\vref{algohessen}. Dans cet \algo la \pcd
\textsf{AjLin}$(H,i,j,c)$ op\`ere une manipulation de
lignes sur la matrice \,$H$: on ajoute \`a la ligne
\,$j$\, la ligne \,$i$\,  multipli\'ee par \,$c$.

\begin{example}
\label{exaHessen} {\footnotesize
\emph{Dans cet exemple on montre une \gui{petite}
matrice \`a
\coes entiers et sa r\'eduction \`a la forme de
Hessenberg.
$$
A :=  \left[
{\begin{array}{rrrrrrrr}
-3 & 3 & 0 & -2 & -3 & 1 & -2 & 2 \\
1 & 2 & 1 & 1 & 2 & 3 & 2 & 0 \\
2 & 2 & 3 & -3 & 3 & 0 & -2 & -3 \\
-2 & 0 & 1 & -2 & 0 & 0 & -1 & 2 \\
3 & -3 & 3 & 3 & 2 & 3 & 0 & 3 \\
3 & 2 & 3 & -3 & 1 & 2 & -1 & -2 \\
-3 & 3 & 3 & 2 & 3 & 1 & -2 & 0 \\
0 & -1 & -3 & -1 & -1 & 1 & -1 & -3
\end{array}}
 \right]
$$
Voici la liste des lignes de la forme r\'eduite de
Hessenberg.
Nous n'avons pas indiqu\'e les $0$ en position $(i,j)$
lorsque
\,$i>j+1$.
$$\begin{array}{l}
\left[\,   -3\,, \,7\,, \,{\ds \frac {-122}{7}}
\,, \,{\ds \frac {-26680}{6639}} \,, \,{\ds
\frac {-4080544}{1522773}} \,, \,{\ds \frac
{4747626797
}{1757764263}} \,, \,{\ds \frac {109259596132466}{
234026268743849}} \,, \,2 \,  \right]
\\[4mm]
%
\left[ \,  1\,, \,11\,, \,{\ds \frac {87}{7}} \,
, \,{\ds \frac {24037}{13278}} \,, \,{\ds
\frac {17521799}{7613865}} \,, \,{\ds \frac
{1473144559
}{1757764263}} \,, \,2\,, \,0 \,  \right]
\\[4mm]
%
\left[\,7\,, \,{\ds \frac {-415}{7}} \,
, \,{\ds \frac {-54333}{4426}} \,, \,{\ds
\frac {-1294739}{2537955}} \,, \,{\ds \frac {689762552
}{585921421}} \,, \,{\ds \frac {-2270125812893340}{
234026268743849}} \,, \,-3 \,  \right]
\\[4mm]
%
\left[  \,{\ds \frac {13278}{49}}
\,, \,{\ds \frac {1670911}{30982}} \,, \,
{\ds \frac {13199211}{1973965}} \,, \,{\ds
\frac {-11965124859}{4101449947}} \,, \,{\ds \frac {
79329636778655517}{1638183881206943}}\,,
\,{\ds \frac {107}{7}}\,\right]
\\[4mm]
%
 \left[ \,{\ds \frac {
-17765685}{19589476}} \,, \,{\ds \frac {-2532182353}{
2246597766}} \,, \,{\ds \frac {2798215923779}{
2593288209346}} \,,
{\ds \frac {6108776229950083011}{
1035800265460275674}} \,, \,{\ds \frac {25553}{4426}}
\,
 \right]
\\[4mm]
%
\left[  \,{\ds
\frac {-2593288209346}{1288243116405}} \,, \,{\ds
\frac {954443884297868}{1487042200034055}} \,,
 {\ds \frac {17689012510838333947}{
28283244709037870895}} \,, \,{\ds \frac {600431}{
2537955}}  \,
\right]
\\[4mm]
%
 \left[  \,
{\ds \frac
{13198847530884339751}{5149558673799888615}} \,,
 {\ds \frac {-10729114442300396518997896}{
2056815059005858366341435}} \,, \,{\ds \frac {
-64207585234}{26366463945}}  \, \right]
\\[4mm]
%
 \left[ \,{\ds \frac {306462654496531416683963262645}{
54768294462168235404375334801}} \,,
 {\ds \frac {481086736521535}{234026268743849}}
\, \right]
\end{array}
$$
On voit appara\^{\i}tre des fractions de grande
taille:
les \coes de la matrice initiale sont major\'es par
$3$ en valeur
absolue, et le \nume le plus grand dans
la matrice transform\'ee est environ \'egal \`a
$3^{61,8}$. Une
\'etude exp\'erimentale dans les m\^{e}mes conditions
avec des
\macas d'ordre \,$n$\, variant entre 8 et 32 donne une
taille  des \coes
\itmds de type quadratique: le \nume ou \deno de
taille maximum est
de l'ordre de  $3^{2(n-2)^2}$.
Il s'agit donc ici d'un cas typique d'une \met
qui ne s'applique efficacement de mani\`ere directe,
en calcul formel,
que dans le cas d'un corps  fini.}
}
\end{example}

\begin{remark}
\emph{Une \matg \,$A\in \K^{n\times n}$\,
 est une matrice de Hessenberg
particuli\`ere qui a ses valeurs propres dans \,$\K$.
Mais une
matrice de Hessenberg qui a ses valeurs propres dans
\,$\K$\,
n'est pas \ncrt triangulaire ni semblable \`a une
matrice
triangulaire, comme on le voit avec la matrice~}
$\left[\begin{array}{ccc}
1&0&1\\1&1&1\\0&0&1 \end{array}\right].$

\end{remark}

\subsubsection*{Nombre d'\oparis}

\noi $\bullet~$ La phase 1 de r\'eduction \`a la
forme de Hessenberg est compos\'ee de
\,$n-2$\, \'etapes. Chacune des \'etapes
\,$p\ (1\leq p \leq n-2)$\, comporte un
travail sur les lignes avec \,$(n-p-1)$\,
divisions, \,$(n-p-1)\,(n-p)$\,
\muls et autant d'additions.
L'op\'eration inverse sur les colonnes
comporte
\,$(n-p-1)\,n$\, \muls et autant
d'additions. \\
Ce qui donne \,${1\over 6}\,(n-1)\,
(n-2)\,(5n+3)\,\asymp {5\over 6}\,n^3$\,
\mulsz/divisions et
\,${5\over 6}\,n\,(n-1)\,(n-2)$\,
\adsosz,
\cad un nombre total d'\oparis dans \,$\K$\, qui est
asymptotiquement de l'ordre de \,${5\over 3}\,n^3$.

\mni
$\bullet~$ La phase 2 qui consiste \`a calculer
les \polcars \,$P_k(X)~(2\leq k\leq n)$\, des
\smpds de la
r\'e\-duite de Hessenberg s'effectue par
\recu sur \,$k$\, \`a partir de
\,$P_0(X)=1$\, et \,$P_1(X)=h_{11}-X$. Si
l'on d\'esigne par \,$S(k)$\, le nombre de
\mulsz/divisions (resp.
\adsosz) permettant de
calculer le \polcar \,$P_k(X)$\, de
\,$H_k$, l'utilisation des relations de
Hessenberg conduit aux relations de
\recu suivantes, vraies pour
\,$2\leq k\leq n\,$:
$$
\left\{ \begin{array}{ll}
S(k)= & S(k-1)+(k-1)+\somm_{i=1}^{k-1}i \\
\; & \mbox{{\small pour les
\mulsz/divisions}} \\
S(k)= & S(k-1)+2\,(k-1)+\somm_{i=1}^{k-2}i \\
\; & \mbox{{\small pour les additions/soustractions}}
\end{array}\right.$$
c'est-\`a-dire, dans les
deux cas (que \,$S(k)$\, d\'esigne le
nombre de \mulsz/divi\-sions
ou celui des \adsosz):
$$
 S(k) = S(k-1) + {1\over
2}\,(k-1)\,(k+2) \,.$$

Comme \,$S(1) = 0$, cela donne par
sommation: \,${1\over 6}\,n\,(n-1)\,(n+4)\,
\asymp {1\over 6}\,n^3$\, \mulsz/divisions et autant
d'\adsos dans le corps \,$\K$\, (la phase de
quasi-triangularisation  est donc la plus
co\^uteuse, asymptotiquement cinq fois plus
ch\`ere en nombre d'\oparis que la phase
de calcul du \polcar de la matrice  \qtgz).

D'o\`u le r\'esultat:

\begin{prop}
L'\aghb calcule les \polcars de toutes
les \smpds  d'une matrice
\,$n\times n$\, sur un corps \,$\K$\,
avec moins de
$~(n+1)\,(n-1)^2~$ \mulsz/divisions
et $~n\,(n-1)^2~$ \adsosz, soit en tout
\,$2\,n^3-3\,n^2+1$\, \oparisz.
\end{prop}

\begin{remark}
\emph{
Ce que la \mhb gagne en
\com \arith par rapport aux pr\'ec\'edentes
\mets de calcul du \polcarz,
elle le perd sur un
aspect essentiel au plan pratique.
Celui de l'absence de
contr\^ole raisonnable de la taille des
\coes \itmdsz. La formule permettant
d'exprimer dans la \mpg
chaque \coe \itmd comme quotient
de deux \deters extraits de la  matrice de
d\'epart (voir propri\'et\'e \vref{propri Gauss}),
ne s'applique plus
dans le processus de quasi-triangularisation de
Hessenberg. En effet les transformations
subies par les lignes sont ici suivies par des
transformations inverses sur les colonnes.
Et on ne dispose pas actuellement pour la
\mhbz, pourtant la plus
rapide en temps \sql si on ne prend
en compte que le nombre d'\oparisz,
d'une formule analogue qui
permette de conclure sur la question
de la taille des \coes \itmdsz.
Cela est confirm\'e par les
r\'esultats exp\'erimentaux que nous avons pu avoir
(voir l'exemple \ref{exaHessen} et le chapitre
\ref{chap experim}).
Dans le cas de matrices \`a \coes entiers,
il y a la possibilit\'e de rem\'edier \`a ce \pb
en utilisant le calcul modulaire (\cf section
\ref{secHadamard}
page~\pageref{Calcmod}).
}
\end{remark}
\begin{remark}
\emph{
Signalons l'existence d'une version modifi\'ee
r\'e\-cen\-te de l'\aghb
sur un anneau int\`egre, d\'evelopp\'ee
dans \cite{JoMa}, qui permet de garder les \coes
\itmdsz, tout au long des calculs, dans l'anneau de
base, suppos\'e int\`egre. Elle semble
bien adapt\'ee au
calcul modulaire sur les anneaux de
\pols \`a \coes entiers.
}
\end{remark}

\section{M\'ethode d'interpolation de Lagrange}
\label{sec lagrange}

Elle ram\`ene le calcul du \polcar d'une \maca
 \,$A\in \A^{n\times n}$\, au calcul de
$n+1$ \detersz.

On est donc suppos\'e \^etre dans une situation
o\`u le calcul des \deters ne pose pas
\pbz: cela peut \^etre le cas par exemple
lorsque la \mpg ne se heurte
pas \`a des \pbs graves de simplification de
fractions, ou lorsque l'on dispose d'un \algo
efficace et sans division pour le calcul des
\deters (comme celui du \dev
suivant une ligne ou une colonne si la matrice
donn\'ee est creuse).
La \met consiste \`a appliquer la formule
d'\iL au \polcar
\,$\rPA(X)=\det{(A-X\In)}$,
c'est-\`a-dire la formule bien connue:
$$
\rPA(X)=\sum_{i=0}^n  \;
{ \left(\ P(x_i)\!
\prod_{i\in\left\{0,\ldots,n\right\}}^{i\neq k}
\!\frac{\,X-x_i\,}{x_k-x_i}\; \right)}
$$
o\`u \,$x_0, x_1, \ldots, x_n \,$ sont
\,$n+1$\, \elts distincts de \,$\A$,
avec la restriction suivante:
les \,$(x_i-x_j)$\, (pour $i\neq j$) doivent
\^etre non diviseurs  de z\'ero dans \,$\A$\,
et on doit disposer d'un \algo de division
exacte par les \,$(x_i-x_j)$\, dans \,$\A$.

C'est par exemple le cas avec \,$x_i=i\times
1_\A$\, lorsque \,$\A$\, est de ca\-rac\-t\'eristique
nulle, ou finie \'etrang\`ere \`a \,$n!$, ou plus
\gnlt lorsque la division exacte par les
\emph{les entiers de \,$\A$\, \infes
ou \'egaux \`a \,$n$}\, (\cad les \elts
$1_\A$, $1_\A+1_\A$, \ldots, $n\,1_\A$), si elle est
possible, est
unique et
r\'ealisable par un \algoz.

En effet, si l'on choisit \,$x_k=k$\, pour
\,$0\leq k\leq n$, la formule d'interpolation
s'\'ecrit:
$$
\det{(A-X\In)}=\sum_{k=0}^n \
\left(\ (-1)^{k}\; \frac{\det{(A-k\,\In)}}{k\,!\,
(n-k)\,!}\!\prod_{i\in\left\{0,\ldots,n\right\}}^{i\neq
k}\!(X-i)\
\right)
$$
 ce qui exige la
possibilit\'e d'effectuer des divisions (exactes) par
les entiers de \,$\A$\, \infes ou \'egaux \`a \,$n$.

En fait, \,$\rPA(X)=(-1)^nX^n+Q(X)$\, avec
\,$\deg(Q)\leq n-1$\, et il
suffit d'appliquer la \miL \`a \,$Q(X)$, ce qui
revient \`a calculer la
valeur de \,$\rPA$\, en \,$n$\, points au lieu de
\,$n+1$.

\mni  {\bf Le nombre  d'\oparisz}
lors de l'ex\'ecution
de cet \algo est
\`a peu pr\`es \,$n$\, fois celle du calcul
d'un \deter d'ordre \,$n$.
Si, pour le calcul des \deters
$P(x_i)=\det(A-x_i\,\I_n)$, on choisit d'utiliser
l'\apg (ou l'\ajb si on est dans une situation
o\`u il s'av\`ere \^etre pr\'ef\'erable \`a
l'\algo de Gauss\footnote{~Dans les deux
cas, nous avons vu que le nombre d'\oparis  est
$\O(n^3)$.}), on obtient donc pour la \miL un
\,$\O(n^4)$.
En fait les meilleurs \algos
sans division dont on dispose actuellement pour
calculer les
\deters passent par le calcul du \polcarz,
ce qui rend caduque la \miLz.
Celle-ci, avec le calcul du
\deter dans l'anneau de base  abandonn\'e
\`a la sagacit\'e de \textsc{Maple},
sera compar\'ee \`a ces autres \algos sur
quelques exemples test\'es sur machine
(voir chapitre \ref{chap experim}).

\section{M\'ethode de Le Verrier et variantes}
\label{sec Leverrier}

Cette \metz, d\'ecouverte en 1848 par
l'astronome fran\c{c}ais Le Verrier \cite{Lev},
repose sur les relations  de Newton entre les
sommes de Newton et les
\pols \syms \elrs
dans l'\agr des \pols \`a \,$n$\,
\idtrs \,$x_1,\dots , x_n$\,
sur un \acom \,$\A$.

De mani\`ere g\'en\'erale, la \met de
Le Verrier appliqu\'ee \`a une matrice
$n\times n$ r\'eclame qu'on soit dans un
anneau o\`u les entiers $1,2,\ldots ,n$
sont non diviseurs de z\'ero. Les seules
divisions requises sont des divisions
exactes par l'un de ces entiers.

\subsection{Le principe g\'en\'eral}

La \met de  Le Verrier consiste
pr\'ecis\'ement \`a d\'eduire le calcul
des \coes du \polcar du calcul de ses sommes de
Newton. Celles-ci sont en effet \'egales aux
traces des puissances de \,$A$\, comme le montre le
lemme \vref{SNTr}.

Rappelons que l'anneau de base n'a pas besoin
d'\^{e}tre
int\`egre, puisque les sommes de Newton peuvent
\^{e}tre
d\'efinies sans recours aux valeurs propres, en
utilisant  les \'equations
(\ref{Newt2}) page \pageref{Newt2}
(\cf \dfn \ref{SNe}).

\ms Ceci donne l'\algo \vref{ALV}.

\begin{algor}[\Algo de Le Verrier]\label{ALV}
\acl{ALV}{\Algo de Le Verrier}
\index{Le Verrier!algorithme de}
\Debut
\Etap{1} {}{Calculer les puissances
\,$A^2,\dots,A^{n-1}$\, de la matrice \,$A$\, ainsi
que les}
\hsu \elts diagonaux  de la matrice \,$A^n$;
\Etap{2} {}{Calculer les traces des matrices
\,$A^1,\,A^2,\dots,\,A^{n}$;}
\Etap{3} {}{Calculer les \coes
\,$p_k$ $(1\leq k\leq n)$  en utilisant les
\'equations (\ref{Newt2}).}
\fin
\end{algor}

\subsubsection*{Nombre d'\oparis}
Pour un anneau \,$\A$\, fix\'e par le contexte,
\emph{nous noterons \,$\mu_M(n)$\,
\indexnota{muM@$\mu_M(n)$} le nombre
d'\oparis \ncrs pour la \mul de deux \macas d'ordre
\,$n$} (on trouvera
une d\'efinition plus pr\'ecise dans la notation
\vref{IDNConst}).
Lorsqu'on utilise la \met \usle de \mul des \macas on
a
\,$\mu_M(n)=n^2\,(2n-1)$.

\ss Pour l'\algo de Le Verrier, le compte est le
suivant:

\ss
\begin{itemize}
\item  l'\'etape 1 utilise \,$(n-2)\,\mu_M(n) +
n\,(2n-1)$\,
\oparis
\item  les \'etapes 2 et 3 utilisent
\,$n^2+2\somm_{k=1}^n(k-1) = 2\,n^2 - n$\,
op\'erations.
\end{itemize}

\begin{proposition}
\label{propComaLV} ~Le nombre total d'\oparis lors de
l'ex\'ecution de l'\algo de Le Verrier, si on utilise
la \mul
usuelle des matrices, est  $2\,n^4+\O(n^3)=\O(n^4)$
(pr\'ecis\'ement \'egal \`a
\,$2\,n\,(n-1/2)\,(n^2-2n+2)$).
\end{proposition}

\ss
 Des \algos d\'eriv\'es de l'\algo de Le
Verrier  ont \'et\'e propos\'es par de nombreux
auteurs, avec des
am\'eliorations concernant la \com (\cf
\cite{Csan,GalPan,PrSa,Sou} et \cite{Fad}). Nous les
\'etudions
dans la suite de cette section et dans le chapitre 
\ref{chap par
Lever}.

\subsection{M\'ethode de Souriau-Faddeev-Frame}
\label{sec faddeev}

Cette \metz, d\'ecouverte s\'epar\'ement par Faddeev
\& Sominskii
\cite{Fadso}, Souriau \cite{Sou} et Frame \cite{Fram},
est une
am\'elioration astucieuse de l'\algo de Le Verrier.

Comme dans le cas de l'\algo de
Le Verrier, l'anneau \,$\A$\, est suppos\'e tel que
\emph{la division par un entier, quand elle est
possible, est unique} (autrement dit, les
entiers ne sont pas des diviseurs de z\'ero)
\emph{et r\'ealisable par un \algoz}.
Cette \met permet  de calculer:
\begin{itemize}
\item le \polcar \,$\rPA$\, d'une \maca
\,$A\in \A^{n\times n}$;
\item l'adjointe de la matrice $A$
et son inverse (s'il existe);
\item un vecteur propre non nul relatif \`a une valeur
propre
donn\'ee de \,$A$\, si l'on suppose de plus que
l'anneau $\A$ est
int\`egre.
\end{itemize}

\sni Posant
\,$P(X)=(-1)^n\rPA(X)=X^n-\left[c_1X^{n-1}+ \cdots
+c_{n-1}X+c_n\right]\,,$ la \met consiste \`a calculer
les \coes
$\,c_k\,$ pour en d\'eduire le \polcar de \,$A$\,. On
utilise pour
cela le calcul de la matrice \cara adjointe de \,$A$\,
tel que
d\'evelopp\'e dans \ref{subsecMacad}. Rappelons la
d\'efinition de
la matrice \cara adjointe de $\,A\,$: c'est la matrice
\,$Q(X) =
\Adj{(X\In -A)}=\sum_{k=0}^{n-1}B_kX^{n-1-k}\,$
(formule
\ref{eqMacad0} page \pageref{eqMacad0}) dans laquelle
les matrices
$\,B_k\,$ sont donn\'ees par les relations
\ref{eqMacad1} (page
\pageref{eqMacad1}): $$B_k=AB_{k-1}-c_k\In~~(1\leq
k\leq
n)~~~\mathrm{avec}~~~B_0=\In\,.$$ On d\'emontre, en
utilisant les
relations de Newton (\ref{Newt2} page
\pageref{Newt2}), que:
$$c_k={1\over k}\,\Tr(A\,B_{k-1})\quad
\mathrm{pour}~1\leq k\leq
n\,.$$
En effet, partant des \'equations suivantes (voir
\ref{eqMacad2}
page \pageref{eqMacad2}) qui d\'ecoulent des relations
\ref{eqMacad1} rappel\'ees ci-dessus:
$$B_k=A^k-c_1A^{k-1}-\dots-c_{k-1}A-c_k\In~~\mathrm{pour~tout~entier}~~
k\in \left\{1,\ldots,n\right\},$$
on consid\`ere les traces des deux membres dans
chacune de ces
\,$n$\, \egts matricielles pour obtenir:
$$
\Tr(B_k)=s_k-c_1s_{k-1}-\dots -c_{k-1}s_1-
nc_k\quad (1\leq k\leq n)\,.
$$
Mais \,$s_k=c_1s_{k-1} + \dots + c_{k-1}s_1 + kc_k$\,
(ce sont les
relations de Newton pour le \poly \,$P(X)$). Comme
\,$\Tr(B_k)=\Tr(AB_{k-1})-nc_k$\, (\`a cause de l'\egt
\,$B_k=AB_{k-1}-c_k\In$), on obtient
\,$\Tr(AB_{k-1})=kc_k\,$.
\qed

\ms Notons par ailleurs (comme nous l'avons fait au \S
\vref{subsecMacad}) que
\,$B_{n}=A\,B_{n-1}-c_n\,\In=0$, \cad
$ A\,B_{n-1}=c_n\,\In=(-1)^{n-1}\det(A)\,\In.$ 
Ce qui montre que si \,$\det(A)$\, est inversible dans
\,$\A$,
alors \,$A$\, poss\`ede un inverse donn\'e par
\,$A^{-1}=(c_n)^{-1}\,B_{n-1}$.

\ss Rappelons \'egalement que
$\,B_{n-1}=(-1)^{n-1}\,\Adj(A)\,$ ce qui
donne, sans autre calcul, l'adjointe de la matrice
\,$A$\,.

\ms Nous traduisons la \met de Souriau-Faddeev-Frame
qui vient
d'\^etre d\'evelopp\'ee par l'\algo \vref{afad}, dans
lequel
$\,B\,$ d\'esigne successivement les matrices
\,$B_0=\In$,
\,$B_1$,\ldots, \,$B_n$, et \,$C$\, les matrices
\,$A$,
\,$A\,B_1$,\ldots, \,$A\,B_n\,$.

\begin{algor}[Algorithme de
Souriau-Faddeev-Frame]\label{afad}
\acl{afad}{Algorithme de Souriau-Faddeev-Frame}
\index{Souriau-Faddeev-Frame!algorithme de}
\Entree Un entier \,$n$\, et une matrice
\,$A\in\A^{n\times n}$.
L'anneau $\A$ est suppos\'e avoir un \algo de division
exacte
par les entiers \,$\leq n$.
\Sortie Le \polcar $\rPA$ de $\,A$.
\Debut
\Varloc $k\in\N$; $c\in\A$; $C,B\in\A^{n\times n}$;
$P\in\A[X]$.
\hsu $Id:=\In$; $B:= Id$; $P:=X^n$;
\hsu \pour{k}{1}{n-1}
\hsd  $C:=B\cdot A$;  $c:=\Tr(C)/k$;
\hsd  $P:=P-c\cdot X^{n-k}$; $B:=C-c\cdot Id$
\hsu \finpour;
\hsu $c:=\Tr(B\cdot A)/n$; $P:=P-c$;
\hsu  $\rPA:=(-1)^n\,P$
\fin
\end{algor}

\subsubsection*{Calcul de vecteurs propres}
Dans le cas o\`u \,$\A$\, est int\`egre, si $\l$ est
une valeur
propre simple de \,$A$\, (\cad une racine simple de
\,$P(X)$), le
m\^eme calcul (donnant entre autres la matrice \cara
adjointe
\,$Q(X) = \Adj{(X\In - A)}$) nous permet d'obtenir un
vecteur
propre non nul associ\'e \`a $\l$. \\[2mm] En effet
\,$Q(X)=\somm_{i=0}^{n-1}B_iX^{n-1-i}$\, donc
\,$\Tr(Q(X))=\somm_{i=0}^{n-1}\Tr(B_i)\,X^{n-1-i}$. \\
Mais
\,$\Tr(B_i)=(i-n)\,c_i$, donc
$$
\Tr(Q(X))=-\somm_{i=0}^{n-1}(n-i)\,c_i\,X^{n-i-1}
= P'(X)
$$
 o\`u \,$P'(X)$\, d\'esigne le
\poly d\'eriv\'e de \,$P(X)$.
Ainsi
\,$\Tr(Q(\l))=P'(\l)\neq 0$\, puisque $\l$
est une racine simple de \,$P$. Par cons\'equent
la matrice \,$Q(\l)$\, n'est pas nulle.

Mais l'\egt \,$(X\In - A)\,Q(X)=P(X)\,\In$\, donne
\,$(A-\l\In)\,Q(\l)=-P(\l)\,\In=0$. Ce qui prouve que
n'importe
quelle colonne non nulle \,$v$\, de \,$Q(\l)$\,
v\'erifie \,$Av=\l
v$\, et c'est donc un vecteur propre non nul de
\,$A$\, relatif
\`a la valeur propre $\l$.

Si l'on d\'esigne par \,$\ell$\, le num\'ero de la
colonne
pr\'esum\'ee non nulle de la matrice
\,$Q(\l)=B_0\l^{n-1} +
B_1\l^{n-2}+\cdots+B_{n-2}\l+B_{n-1}$\, et par
\,$b_k^{\,\ell}$\,
la colonne \num$\ell$\, de \,$B_k$, le calcul de ce
vecteur propre
peut se faire de la mani\`ere suivante:

$\bullet$~~ Poser $\ v_0 = e_\ell\ \ $
(colonne num\'{e}ro \,$\ell$\, de la matrice $~\In$);

$\bullet$~~ Faire $~v_k = \l\,v_{k-1} + b_k^{\,\ell}~$
(pour \,$k$\, allant de 1 \`{a} \,$n-1$).

\sni Le vecteur propre recherch\'e n'est autre que
\,$v=v_{n-1}$.

\ss Plus \gnltz, si la \mut \geoqz{\footnote{~La \mut
\geoq d'une
valeur propre est par \dfn la dimension du sous-espace
propre correspondant.}} de $\l$ est \'egale \`a 1,
\index{multiplicit\'e!g\'eom\'e\-tri\-que}
\index{valeur propre!multiplicit\'e g\'eom\'e\-tri\-que d'une} la
matrice
\,$Q(\l)$\, n'est pas nulle (elle est de rang $1$) et
n'importe
quelle colonne non nulle de \,$Q(\l)$\, repr\'esente
un vecteur
propre non nul de \,$A$\, pour la valeur propre $\l$.

Si par contre la \mut \geoq (et par cons\'equent la
\mut \agqz\footnote{~C'est la \mut de \,$\l$\, en tant que
z\'ero du \polcarz.}) de la valeur propre \,$\l$\, 
est \supee ou \'egale \`a
2, non seulement la trace, mais la matrice \,$Q(\l)$\,
elle-m\^eme est nulle d'apr\`es la propri\'et\'e \vref{p2},
puisque le rang de la matrice singuli\`ere \,$A-\l\In$\, 
est, dans ce cas, au plus \'egal \`a \,$n-2$. 
La matrice \,$Q(\l)$, dans ce cas, ne donne
donc aucun vecteur propre non nul de $A$.

On montre alors que ce sont les matrices d\'eriv\'ees
successives
(par rapport \`a \,$X$) de la matrice \,$Q(X)$\, qui
permettent de
calculer des vecteurs propres non nuls relatifs \`a
\,$\l$.

Consid\'erons en effet pour
\,$k\in\N$\, l'op\'erateur \,$\Delta^{[k]}\,:\,\A[X]
\dans \A[X],$\,
\,$P\mapsto\Delta^{[k]}P$,  o\`u
\,$\Delta^{[k]}P\,(X)=P^{[k]}(X)$\, est d\'efini par
l'identit\'e
$$
P(X+Y) =
\somm\nolimits_{k\geq 0}P^{[k]}(X)\,Y^k\,.
$$

Remarquons que \,$P^{[0]}(X)=P(X)$,
\,$P^{[1]}(X)=P'(X)$\, et qu'en \cara nulle
\,$\Delta^{[k]} = {1\over {k!}}\,D^k$\, o\`u \,$D^k$\,
est l'op\'erateur de d\'erivation d'ordre \,$k$\, dans
\,$\A\,[X]$\, ($P^{[k]}(X)={1\over{k!}}P^{(k)}(X)$).
Il est facile de voir qu'en \cara quelconque,
\,$\Delta^{[k]}$\, est un \,$\A$\,-\,\endom de
l'\agr \,$\A\,[X]$\, qui v\'erifie une formule
analogue \`a la formule de Leibnitz, mais plus simple:
$$
\Delta^{[k]}(P_1\,P_2)=\sum_{i=0}^k
\Delta^{[i]}(P_1)\;\Delta^{[k-i]}(P_2)
$$

En outre un \poly \,$P$\, admet $\l$ comme racine
d'ordre
\,$k\geq 1$\, si et seulement si
\,$\,P^{[0]}(\l)=P^{[1]}(\l)= \dots=P^{[k-1]}(\l)=0$\,
et \,$P^{[k]}(\l)\neq 0$.
Appliquant successivement les op\'erateurs
\,$\Delta^{[k]}$\, pour \,$k$\, allant de 1 \`a
\,$m$\, (o\`u $m$ est la \mut \agq
de la valeur propre $\l$) \`a l'\egt
matricielle:
$$
\,(X\In - A)\,Q(X)=P(X)\,\In\,,
$$
on obtient la suite d'\egtsz:
$$
Q^{[k-1]}(X)+(X\In - A)\,Q^{[k]}(X) =
P^{[k]}(X)\,\In\quad (1\leq k\leq m).$$

Rempla\c{c}ant dans ces \egts $X$ par la valeur propre
$\l$, et
tenant compte du fait que \,$Q(\l)=P(\l)=P^{[1]}(\l)=
\cdots=P^{[m-1]}(\l)=0$, on obtient le \sysz:
$$
\left\{ \begin{array}{l} Q(\l)=0 \\
(\l\In - A)\,Q^{[1]}(\l) = 0 \\
Q^{[1]}(\l) + (\l\In - A)\,Q^{[2]}(\l) = 0 \\
~~~\vdots~~~\vdots~~~\vdots ~~~ \vdots ~~~ \vdots \\
 Q^{[m-2]}(\l) + (\l\In - A)\,Q^{[m-1]}(\l) = 0 \\
 Q^{[m-1]}(\l) + (\l\In - A)\,Q^{[m]}(\l) =
 P^{[m]}(\l)\,\In\,.
\end{array}\right.$$

Soit \,$r\in \NN$\, le plus petit entier tel que
\,$Q^{[r]}(\l)
\neq 0$. \\ Alors \,$r<m$\, car sinon, on aurait
\,$Q(\l)=Q^{[1]}(\l)=\cdots=Q^{[m-1]}(\l)=0$\, et
\,$(\l\In -
A)\,Q^{[m]}(\l) = P^{[m]}(\l)\,\In$\, avec
\,$P^{[m]}(\l) \neq 0$,
ce qui contredit (nous sommes dans un anneau
int\`egre) le fait
que la matrice \,$\l\In - A$\, est singuli\`ere.

Donc \,$(\l\In - A)\,Q^{[r]}(\l)=0$\, et toute colonne
non nulle
de \,$Q^{[r]}(\l)$\, est un vecteur propre non nul de
\,$A$\, pour
la valeur propre multiple $\l$.

\subsubsection*{Nombre d'\oparis }
 L'\algo de Souriau-Faddeev-Frame
 consiste \`a calculer, pour \,$k$\, allant de
1 \`a \,$n$, le produit matriciel \,$A_k=AB_{k-1}$,
le \coe \,$c_k=
{1\over k}\Tr(A_k)$\,
et enfin la matrice \,$B_k=A_k-c_k\In$.

Rappelons qu'on d\'esigne par \,$\mu_M(n)$\, le nombre
d'\oparis
dans l'anneau de base pour la \mul de deux matrices
carr\'ees
d'ordre \,$n$. Le co\^ut de l'\algo de
Souriau-Faddeev-Frame
s'\'el\`eve \`a $$(n-2)\,\mu_M(n) + n\,(2n-1) +
2\,n\,(n-1)\,.$$
C'est \`a tr\`es peu pr\`es le m\^eme co\^ut que pour 
l'algo de
Le Verrier (on gagne \,$n$\, \oparisz).

Outre la plus grande simplicit\'e, l'avantage est que
l'on a aussi
calcul\'e la matrice adjointe. En particulier le
calcul de la
matrice inverse, si elle existe ne co\^ute que
\,$n^2$\, divisions
\supts dans \,$\A$. Enfin le calcul d'un vecteur
propre non nul
relatif \`a une valeur propre donn\'ee de \mut \geoq
\'egale \`a 1
se fait moyennant \,$2\,n\,(n-1)$\, \oparis  \suptsz.

\begin{proposition}
\label{propComAFSF}
Avec la \met de
Faddeev-Souriau-Frame, le calcul du \polcar de la
matrice \,$A$,
de son \deterz, de sa matrice adjointe,
de son inverse quand elle existe, ainsi que des
sous-espaces propres de dimension $1$ (quand on
conna\^{\i}t la valeur propre correspondante)
se fait en \,$2\,n^{4}+\O(n^3)=O(n^4)$\, \oparisz.
Pour le calcul du seul \polcar on en effectue
pr\'ecis\'ement
\,$2\,n\,(n-1)\,(n^{2}-3n/2+1/2)$ \opsz.
\end{proposition}

\subsection{M\'ethode de Preparata \& Sarwate}
\label{secPrepSar}

La \met de \PrSa est une acc\'el\'eration astucieuse 
de la \met de Le Verrier, bas\'ee sur la remarque
simple suivante. Pour calculer la trace d'un produit
\,$AB$\, de deux \macas d'ordre \,$n$, 
il suffit d'\'ex\'ecuter \,$2n^2$\, \oparis
puisque \,$\Tr\,AB=\som_{k,\ell}a_{k,\ell}\,b_{\ell,k}$.
Or le calcul le plus co\^{u}teux dans la \met de Le Verrier
est celui des traces des puissances successives de
la matrice \,$A$\, dont on veut calculer le \polcarz.

Posons donc \,$r=\esup{\sqrt{n}\,}$, \,$B_0=C_0=\In$,
\,$B_1=A$, et calculons
les \,$B_{i}=A^i$\, pour \,$i=2,\ldots,r$,
puis les \,$C_{j}={B_r}^j$\, pour \,$j=1,\ldots,r-1$.
Ce calcul consomme \,$(2r-3)(n^3-n^2)\simeq 2\,n^{3,5}$\,
\oparis dans \,$\A$.

On a alors 
\,$\Tr\,A^{rj+i}=\Tr\,B_iC_j
=\som_{k,\ell}b_{i,k,\ell}\,c_{j,\ell,k}$,
et les valeurs \,$rj+i$\, pour \,$0\leq i,j\leq r-1$\,
parcourent l'intervalle \,$[0,r^2-1]$.
Si \,$r^2=n$\, on doit calculer en outre \,$S_n=\Tr\,C_1C_{r-1}$.
On obtient donc toutes les sommes de Newton \,$S_m=\Tr\,A^m$,
\,$1\leq m\leq n$\, pour un peu moins que 
\,$2n^3$\, \oparis \supts dans \,$\A$.

\ss Ceci donne l'\algo \vref{APSa}. Comme la r\'ecup\'eration
des \coes du \polcar \`a partir des sommes de Newton
r\'eclame \,$\O(n^2)$\, \oparis on obtient la proposition
suivante.

\begin{algor}[\Algo de \PrSaz, version \sqle simple.]
\label{APSa}
\acl{APSa}{\Algo de \PrSaz, version \sqlez}
\index{Preparata\&Sarwate!algorithme de}
\Entree Une \maca  \,$A\in\A^{n\times n}$\, o\`u \,$\A$\, 
est un anneau v\'erifiant les hypoth\`eses de l'\aglvz.
\Sortie Le \polcar 
\,$\rPA=(-1)^n\,(X^n+\som_{k=1}^{n}p_{k}X^{n-k})$.
\Debut
\Varloc $i$, $j$, $r\in\N$; $B_i$, $C_{j}\in\A^{n\times n}$ 
$(i,j=1..r-1)$,
$S_i\in\A$ $(i=1..n)$;
\Etap{1} {Calcul des puissances \,$A^i$\, 
pour \,$i< r=\esup{\sqrt{n}\,}$.} 
{$r:=\esup{\sqrt{n}\,}$; \,$B_1:=A$; 
\,$S_0:=n$; \,$S_1:=\Tr\,A$;}
\hsu \pour{i}{1}{r-2} 
\hsd $B_{i+1}:=A\,B_{i}$; \,$S_{i+1}:=\Tr\,B_{i+1}$
\hsu \finpour;
\Etap{2} {Calcul des puissances \,$A^{rj}$\, 
pour \,$j< r$.}
{$C_{1}:=A\,B_{r-1}$; \,$S_{r}:=\Tr\,\,C_1$;}
\hsu \pour{j}{1}{r-2} 
\hsd $C_{j+1}:=C_{1}\,C_{j}$; \,$S_{(j+1)r}:=\Tr\,\,C_{j+1}$
\hsu \finpour;
\Etap{3}{Calcul des sommes de Newton.}
{\pour{i}{1}{r-1}}
\hsd {\pour{j}{1}{r-1}}
\hst \,$S_{jr+i}:=\Tr\,B_iC_{j}$
\hsd \finpour
\hsu \finpour
\hsu \sialors{n=r^2}{$S_{n}:=\Tr\,\,C_{1}C_{r-1}$} \finsi;
\Etap{4} {Calcul des \coes de $\rPA$.}
{Calculer les \coes \,$p_k$ $(1\leq k\leq n)$  
en utilisant les \'equations (\ref{Newt2}).}
\fin
\end{algor}

\begin{proposition}
\label{propCoAPSa} ~Supposons que l'\acom \,$\A$\,
satisfasse les hypoth\`eses de l'\aglvz:
\emph{la division par un entier, quand elle est
possible, est unique, et r\'ealisable par un \algoz}.
\index{Le Verrier!hypoth\`eses de l'\algo de}
Le nombre total d'\oparis lors de
l'ex\'ecution de l'\algo de \PrSaz, si on utilise
la \mul \usle des matrices, est  
\'egal \`a $2\,n^{3,5}+\O(n^3)=\O(n^{3,5})$.
\end{proposition}

\section{M\'ethode de Samuelson-Berkowitz}
\label{sec berkoseq}

Elle est bas\'ee sur la \met de partitionnement
(\cite{Gas}  pp.~291--298, \cite{Fad}), attribu\'ee
\`a Samuelson
\cite{Sam}, et elle a l'avantage de s'appliquer
\`a un \acomaz.

Berkowitz \cite{Ber} en donne une version \paral
de laquelle nous
extrayons une \met \sqle
particuli\`erement simple et
efficace.
Elle montre l'int\'er\^et pratique de cet \algo
pour les machines \sqlesz,
en le pla\c{c}ant parmi les
\algos les plus performants actuellement pour
le calcul sans division du \polcar (\cf les test
exp\'erimentaux
pr\'e\-sen\-t\'es au chapitre \ref{chap experim}).

\subsection[Principe g\'en\'eral de l'\algo]{Principe
g\'en\'eral
de l'\algo}\label{subsec.samu}

Soit \,$A=(a_{ij})\in\A^{n\times n}$\, une \maca
d'ordre \,$n \geq 2$\,
sur un
\acoma  \,$\A$.
Conform\'ement aux notations introduites dans la
section
\ref{sec AlgLin},  pour tout entier
\,$r~(1 \leq r \leq n-1)$, on d\'esigne par \,$A_r$\,
la \smpd d'ordre \,$r$\,de \,$A$.  On partitionne
comme suit la
matrice \,$A_{r+1}$:
$$ A_{r+1}\,=\,
\cmatrix{
A_r
&A_{1..r,r+1}\cr A_{r+1,1..r}&a_{r+1,r+1}}\,=\,
\cmatrix{
A_r
&S_r\cr R_r&a_{r+1,r+1}}\,.
$$

Le \polcar \,$P_{r+1}(X)$\, de \,$A_{r+1}$\, est 
reli\'e au
\polcar \,$P_r(X)=\sum_{i=0}^r p_{r-i}X^i$\, de
\,$A_r$\, par la
formule de Samuelson (\ref{EqSamu}) (proposition
\vref{Samu}) que
l'on peut r\'e\'ecrire sous la forme suivante:
\begin{equation}\label{EqSam}
P_{r+1}= \left\{
\begin{array}{lll}
(a_{r+1,r+1} - X)\, P_r(X)\, +
\\[1mm]
\sum_{k=2}^{r+1} \;\left[(R_rA_r^{k-2}S_r) p_0 +
\cdots
+ (R_r S_r) p_{k-2}\right]\, X^{r+1-k}
\end{array}
\right.
\end{equation}

\ss Notons \,$Q_{r+1}$\, le \polz: $$
-X^{r+1}+a_{r+1,r+1}\,X^r+R_rS_r\,X^{r-1}+
R_rA_rS_r\,X^{r-2}+\cdots+R_rA_r^{r-1}S_r\,. $$ Pour
calculer
\,$P_{r+1}(X)$\, selon la formule de Samuelson on
peut:
\begin{itemize}
\item  effectuer le produit \,$P_{r}\,Q_{r+1}$~;
\item  supprimer les termes de degr\'e $<r$~;
\item  et enfin diviser par \,$X^r\,$.
\end{itemize}
\noindent On peut aussi d\'ecrire ce calcul sous la
forme:
\begin{equation}\label{sam}
\overrightarrow{P_{r+1}}=\Toep(Q_{r+1})\times
\overrightarrow{P_r}
\end{equation}
o\`u \,$\overrightarrow{P_r}$\, est le vecteur colonne
\,$\tra{(p_0,p_{1},\ldots ,p_r)}$\, des \coes du \poly
\,$P_r$\, et
\,$\Toep(Q_{r+1})\in\A^{(r+2)\times(r+1)}$\,
est la \mto suivante d\'efinie \`a partir
du \pol \,$Q_{r+1}$:
$$
\Toep(Q_{r+1}) =
\left[\begin{array}{ccccc}
-1          & 0 & \cdots & \cdots & 0 \\[1mm]
a_{r+1,r+1} & -1 & \ddots &   & \vdots \\[1mm]
R_rS_r & \ddots & \ddots & \ddots & \vdots \\[2mm]
\vdots & \ddots & \ddots & \ddots & 0 \\[2mm]
R_rA_r^{r-2}S_r &  &  \ddots
 & \ddots & -1 \\[2mm]
R_rA_r^{r-1}S_r &  R_rA_r^{r-2}S_r & \dots & R_rS_r &
a_{r+1,r+1}
\end{array}\right]
$$
\begin{algor}[\Algo de \Berz, principe g\'en\'eral.]
\label{algoBerG}
    \acl{algoBerG}{\Algo de \Berz, principe \gnlz}
\index{Berkowitz!algorithme de}
\Entree Une matrice $A\in\A^{n\times n}$.
\Sortie Le \polcar $\rPA(X)$ de $A$.
\Debut
\Etap{1}{}{Pour
\,$k<r$\, dans \,$\{1,\ldots,n\}$\,  calculer
les produits \,$R_r\,(A_r)^k\,S_r$,}
\hsu ce qui donne les \pols \,$Q_{r+1}$\, et les
matrices
\,$\Toep(Q_{r+1})$,
\Etap{2}{}{Calculer le produit
\,$\Toep(Q_{n})\,\Toep(Q_{n-1})\,\cdots\,\Toep(Q_{2})\,
\overrightarrow{P_1}$:}
\hsu on obtient \,$\overrightarrow{\rPA}$.
\fin
\end{algor}
\sni Avec $\overrightarrow{P_1}=
\left[\begin{array}{c}-1\\
a_{1,1}\end{array}\right]$, on obtient  l'\agb
informel
\vref{algoBerG}.

\subsection{Version \sqle}  \label{subsec.berkoseq}
Dans la version \sqle la plus simple de l'\agbz, le
calcul des
\coes de la matrice \,$\Toep(Q_{r+1})$ se fait
naturellement par
l'utilisation exclusive de produits scalaires ou de
produits de
matrices par des vecteurs. De m\^{e}me dans l'\'etape
2, le
produit s'effectue de droite \`a gauche, donc
n'utilise que des
produits matrice par vecteur. Cela donne l'\algo
\vref{algoberko}.

\begin{algorH}[\Algo de \Berz, version \sqle simple.]
    \label{algoberko}
    \acl{algoberko}{\Algo de \Berz, version \sqlez}
\Entree Un entier $n\geq 2$ et une matrice
$A=(a_{ij})\in\A^{n\times n}$.
\Sortie Le \polcar de $A$:
$\rPA(X)$.
\Debut
\Varloc $i$, $j$, $k$, $r\in\N$ ;
$v=(v_{i})$, $c=(c_{i})$, $s=(s_{i})$, $P=(p_{i})$:
listes de longueur variable
$r$\, $(1\leq r\leq n+1)$ dans $\A$.
\hsu \# initialisation de $c$ et de $v$
\hsu $c_{1}:=-1$ ; $v:=(-1,a_{11})$ ;
\hsu \# Calcul des \polcars des matrices principales
\hsu \# dominantes d'ordre $\geq 2$ (les listes
successives
        dans $P=(p_{i})$)
\hsu \pour{r}{2}{n}
\hsd $(s_{i})_{i=1..r-1}:=(a_{ir})_{i=1..r-1}$;
\hsd $c_{2}:=a_{rr}$ ;
\hsd \pour{i}{1}{r-2}
\hst $c_{i+2}:=\sum_{j=1}^{r-1}a_{rj}\,s_{j}$ ;
\hst \pour{j}{1}{r-1}
$p_{j}:=\sum_{k=1}^{r-1}a_{jk}\,s_{k}\;$
     \finpour;
\hst $(s_{j})_{j=1..r-1}:=(p_{j})_{j=1..r-1}$
\hsd \finpour;
\hsd $c_{r+1}:=\sum_{j=1}^{r-1}a_{rj}\,s_{j}$ ;
\hsd \pour{i}{1}{r+1}
\hst 
$p_{i}:=\sum_{j=1}^{\min{(r,i)}}c_{i+1-j}\,v_{j}$
\hsd \finpour;
\hsd $(v_{i})_{i=1..r+1}:=(p_{i})_{i=1..r+1}$
\hsu \finpour;
\hsu $\rPA(X):=\sum_{i=0}^nv_{i+1}X^{n-i}$
\fin
\end{algorH}

Ainsi, sans calculer des puissances de matrices
\,$A_r^{k-1}$\,
($3\leq k\leq r$) on commence par calculer
\,$R_rS_r$\, puis
successivement, pour \,$k$\, allant de 2 \`a \,$r$, le
produit
(matrice par vecteur) $\:A_r^{k-1}S_r\:$  suivi du
produit
scalaire $\:R_rA_r^{k-1}S_r\:$, ce qui se traduit par
\,$2r^3+r^2-3r+1$\, \oparis pour chaque \,$r~(1\leq r
\leq n-1)$.
On en d\'eduit que le nombre d'\oparis (dans l'anneau
de base
\,$\A$) intervenant dans ce calcul est \'egal \`a:
$$\,\somm_{r=1}^{n-1}(2r^3+r^2-3r+1)={1\over 2}n^4 -
{2\over 3}n^3
 -{3\over 2}n^2
+ {8\over 3}n -1 \,.$$

Il en est de m\^eme pour la \mul des \mtos
\,$\Toep(Q_r)$. On
commence par multiplier la premi\`ere matrice
\,$\Toep(Q_1)=\overrightarrow{P_1}$\, \`a gauche par
la matrice
\,$\Toep(Q_2)$\, pour obtenir le vecteur
\,$\overrightarrow{P_2}$\, qui est un vecteur
\,$3\times 1$\, et
ainsi de suite jusqu'\`a
\,$\overrightarrow{P_n}=\Toep(Q_n)\times
\overrightarrow{P_{n-1}}$. \\
Comme chaque \mul \,$\overrightarrow{P_{r}}=
\Toep(Q_{r})\times
\overrightarrow{P_{r-1}}$\, (d'une matrice
sous-triangulaire
\,$(r+1)\times r$\, avec des $-1$ sur la diagonale par
un vecteur
\,$r\times 1$) co\^ute  \,$r\,(r-1)$\, \oparis dans
\,$\A$, le
calcul de \,$\overrightarrow{P_n}$\, se fait en
\,$\somm_{r=2}^n
(r^2 -r) = {1\over 3}(n^3-n)$\, \oparis de
base\footnote{~Ce
calcul peut \^{e}tre acc\'el\'er\'e en utilisant une
\mul rapide
des \pols (\cf chapitre \ref{chap multipoly}), mais
cela ne change
pas substantiellement le r\'esultat global qui reste
de \,$\O(n^4)$\, \oparis avec la m\^{e}me constante
asymptotique.
Nous n'avons pas impl\'ement\'e cette am\'elioration
lors de nos tests exp\'erimentaux.}.

\begin{proposition}
\label{propberkoseq} Le co\^ut total de l'\algo \sql
simple de
Berkowitz s'\'el\`eve \`a
$$
 {1\over 2}n^4 - {1\over 3}n^3
 -{3\over 2}n^2
+ {7\over 3}n -1 \; \; \leq \; \; {1\over 2}n^4 -
{1\over 3}n^3 $$
\oparis dans l'anneau de base.
\end{proposition}

\section{M\'ethode de Chistov}
\label{secChistov}

\subsection{Le principe
g\'en\'eral}\label{secChistovGene}
La \met de Chistov \cite{Chis} consiste \`a calculer le \polcar
\,$\rPA(X)$\,
d'une \maca \,$A\in\A^{n\times n}$\, ($n \geq 2$) en
le
ramenant \`a l'inversion du \pol formel
\,$Q(X)=\det{(\In-XA)}$\, dans l'anneau des
s\'eries formelles
\,$\A[[X]]$.

Ce \poly est, \`a un signe pr\`es, le \poly
r\'eciproque du \polcar puisque
$$
(-1)^nX^nQ\left({1\over X}\right)=
(-X)^n\det{\left(\In-{1\over X}A\right)}=
\det{(A-X\In)}=\rPA(X).
$$

Comme les \pols \,$\rPA$\, et \,$Q$, r\'esultat
final du calcul, sont de degr\'e \,$n$, tous
les calculs peuvent se faire modulo \,$X^{n+1}$\,
dans l'anneau des s\'eries formelles \,$\A[[X]]$,
\cade peuvent se faire dans \emph{l'anneau des \dlis
\`a l'ordre \,$n$\, sur \,$\A$}:
\,$\aqo{\A[X]}{X^{n+1}}$.
\index{anneau des  \dlis \`a l'ordre \,$n$}

Dans la suite, nous noterons souvent cet anneau
\,$\A_n$.
\indexnota{An@$\A_n=\aqo{\A[X]}{X^{n+1}}$ !anneau des
\dlisz}

L'\algo de Chistov\index{Chistov!algorithme de}
utilise le fait que, pour toute \maca
\,$B$\, d'ordre \,$n$\, \regz,  le $n\,$\eme \elt
de la diagonale  de la matrice inverse $B^{-1}$,
not\'e \,$(B^{-1})_{n,n}$, est \'egal \`a:
$$
(B^{-1})_{n,n}=\frac{\det\,{B_{n-1}}}{\det\,{B}}
\quad \mathrm{ou\,\, encore}\quad
(B_{n}^{-1})_{n,n}=\frac{\det\,{B_{n-1}}}{\det\,{B_n}}\,,
$$
o\`u \,$B_{r}$\, d\'esigne la \smpd d'ordre \,$r$\, de
\,$B$
(avec la convention \,$\det\,{B_0}=1$).
Ceci permet d'\'ecrire lorsque \,$B$\, est \fregz:
$$
\left(B_n^{-1} \right)_{n,n}\times \left(B_{n-1}^{-1}
\right)_{n-1,n-1}
\times \cdots \times \left(B_{1}^{-1} \right)_{1,1} =
\frac{1}{\det{B}}\,.
$$

Appliquant ce fait \`a la matrice
\,$B=\In-XA_n\in\A[X]^{n\times n}$\,
qui est \freg puisque
tous ses \mips dominants sont des \pols
de terme constant \'egal \`a 1 et sont donc
inversibles dans
l'anneau \,$\A[[X]]$, on obtient:
\begin{equation} \label{EqChis1}
Q(X)^{-1}=\left[\,\det{(\In-XA_n)}\,\right]^{-1}=
\prod_{r=1}^n\left(B_{\,r}^{-1} \right)_{\\ r,r}.
\end{equation}

Mais on a un \iso canonique \,$\A[[X]]^{r\times
r}\simeq
\A^{r\times r}[[X]]$\,  et la matrice
\,$B_r=\Ir-XA_r$\,
est aussi inversible dans
l'\agr des s\'eries formelles sur l'anneau de
matrices \,$\A^{r\times r}$, et son inverse est la
matrice:
%
\begin{equation} \label{EqChis2}
\,B_{r}^{-1}=\Ir+\somm_{k=1}^{\infty}(A_r)^k\,X^k
\in \A^{r\times r}[[X]]\,.
\end{equation}

Donc en notant \,$\Er$\, la \,$r\,$\eme colonne de
$\Ir$:
\begin{equation} \label{EqChis3}
\left(B_{r}^{-1} \right)_{r,r} \ \mathrm{mod} \  
X^{n+1}=
1+\sum_{k=1}^n\left(\tra{\,\Er}\,(A_r)^k\,\Er\right)\,X^k\,.
\end{equation}

Par cons\'equent en notant \,$
\widetilde{Q}(X) =  Q(X)^{-1}\ \mathrm{mod} \ 
X^{n+1}$\,
on obtient:

$$\widetilde{Q}(X) = \prod_{r=1}^n\,
\left[\,1+\sum_{k=1}^n
\left(\tra{\,\Er}\,(A_r)^k\,\Er\right)\, X^k \,\right]
\ \mathrm{mod} \  X^{n+1}
$$
et donc,
$$
Q(X) =\left\{\prod_{r=1}^n\,\left[\,1+\sum_{k=1}^n
\left(\tra{\,\Er}\,(A_r)^k\,\Er\right)\, X^k\,
\right]  \right\}^{-1}\;\,\mathrm{dans}\;
\aqo{\A[X]}{X^{n+1}}
$$

Ainsi \,$Q(X)$\, est l'inverse modulo \,$X^{n+1}$\, du
produit modulo
\,$X^{n+1}$\, de \,$n$\, \pols de terme constant
\'egal
\`a 1 et de degr\'e \infe ou \'egal \`a \,$n$.

\ss Rappelons que le \polcar \`a calculer
\,$\rPA(X)$\,
est le produit par \,$(-1)^{n}$\, du \poly
r\'eciproque
\index{polyn\^ome!r\'eciproque} \`a l'ordre \,$n$\, de
\,$Q(X)$. On obtient alors l'\algo de Chistov
\vref{algoChisG}.

\begin{algor}[\Algo de Chistov, principe g\'en\'eral.]
\label{algoChisG}
    \acl{algoChisG}{\Algo de Chistov, principe \gnlz}
\Entree la matrice $A\in\A^{n\times n}$.
\Sortie le \polcar $\rPA(X)$ de $A$.
\Debut
\Etap{1}{}{Calculer  pour
\,$r,k\in \{1,\ldots,n\}$\,
les produits \,$\tra{\,\Er}\,(A_r)^k\,\Er$,}
\hsu ce qui donne les
\pols \,$\left(B_{\,r}^{-1} \right)_{\\ r,r}$
(formule $(\ref{EqChis3})$).
\Etap{2}{}{Calculer le produit des $n$ \pols
pr\'ec\'edents modulo \,$X^{n+1}$,}
\hsu ce qui donne
\,$Q(X)^{-1}\ \mathrm{mod} \  X^{n+1}$
(formule $(\ref{EqChis1})$).
\Etap{3}{}{Inverser modulo \,$X^{n+1}$\,
le \poly pr\'ec\'edent: on obtient \,$Q(X)$.}
\Etap{4}{}{Prendre le \poly r\'eciproque
\`a l'ordre \,$n$\, du \poly \,$Q(X)$.}
\hsu On obtient \,$\rPA(X)$\, en multipliant par
$(-1)^{n}$.
\fin
\end{algor}

Nous d\'etaillons
maintenant la version \sqle
\elr de cet \algoz.

\subsection{La version \sqle}\label{secChisSeq}

Dans la  version \sqle la plus simple on obtient
l'\algo
\vref{algochistov}.
\begin{algor}[Algorithme de Chistov, version
    \sqle simple.]\label{algochistov}
    \acl{algochistov}{Algorithme de Chistov, version
    \sqle simple}
\Entree Un entier $n\geq 2$ et une matrice
$A=(a_{i,j})\in\A^{n\times n}$.
\Sortie Le \polcar de $A$, $\rPA\in\A[X]$.
\Debut
\Varloc $i$, $j$, $k$, $r\in\N$;
$q=(q_{i}),\,b=(b_{i}),\,c=(c_{i})\,\in\A^{n+1}$
o\`u  $0\leq i\leq n$;
$v=(v_{i})\,,w=(w_{i})\,\in\A^{r}$ o\`u
$1\leq i\leq r$,  $Q\in \A[X]$.
\hsu $q:=(1)_{i=0..n}$; $c:=q$; 
\,\,\,\,(initialisation)
\hsu \pour{i}{1}{n} $q_{i}:=q_{i-1}\,a_{1,1}\;$
\finpour;
\hsu \pour{r}{2}{n}
\hsd $v:=(a_{i,r})_{i=1..r}$;   $c_{1}:=v_{r}$;
\hsd \pour{i}{2}{n-1}
\hst \pour{j}{1}{r}
      $w_{j}:=\sum_{k=1}^{r}a_{j,k}\,v_{k}\;$ 
\finpour;
\hst $v:=w$;   $c_{i} := v_{r}$
\hsd \finpour;
\hsd $c_{n}:=\sum_{k=1}^{r}a_{r,k}\,v_{k}$;
\hsd \pour{j}{0}{n}
     $b_{j}:=\sum_{k=0}^{j}c_{j-k}\,q_{k}\;$ \finpour;
\hsd $q := b$
\hsu \finpour
\hsu $Q:=1/\left(\sum_{k=0}^{n}q_{k}X^k\right)\;{\rm
mod}\;X^{n+1}$;
\hsu $\rPA:=(-1)^nX^nQ(1/X)$
\fin
\end{algor}


On d\'emontre maintenant que le co\^ut en
nombre d'\oparis  dans l'anneau de
base pour cette version \elr
 est asymptotiquement de l'ordre de
\,$(2/3)\,n^4$.

Reprenons en effet les quatre \'etapes
dans l'\algo de Chistov \gnl \vref{algoChisG}.

\ss
$\bullet$ {\bf L'\'etape 1}, la plus co\^uteuse,
se ram\`ene en fait \`a calculer successivement
les produits \,$(A_r)^k\Er$\, pour \,$1\leq r\leq n$\,
et pour \,$1\leq k\leq n$, puisque
\,$\tra{\,\Er}\,(A_r)^k\,\Er$\, n'est autre que la
\,$r\,$\eme composante du vecteur
\,$(A_r)^k\,\Er$.

Pour chaque valeur de \,$r~(1\leq r\leq n)$,
on commence par calculer \,$A_r\Er$\, puis,
pour $k$ allant de 2 \`a $n$, le produit de
la matrice \,$A_r$\, par le vecteur
\,$(A_r)^{k-1}\Er$, ce qui se traduit par
\,$n\,(2r^2 - r)$\, \oparis
(\adsos et \mulsz) pour
chaque $r$ compris entre 1 et $n$. On en
d\'eduit que le nombre d'\ops dans ce
calcul est \'egal \`a:
$$n\,\somm_{r=1}^{n}\,(2r^2 - r) \,=\,
{1 \over 6}\,n^2\,(n+1)(4n-1)
\,=\,
(2/3)\,n^4+\O(n^3)\,.$$

\ss
$\bullet$ {\bf L'\'etape 2} revient \`a calculer
le produit tronqu\'e de \,$n$\, \pols de degr\'e au
plus \'egal
\`a \,$n$, ce qui se fait en
\,$\O(n^3)$\, \ops de base.

\ss
$\bullet$ {\bf L'\'etape 3} consiste \`a
calculer le produit tronqu\'e \`a l'ordre
\,$n$\, des \,$\esup{\log{n}}$\, \pols
\,$1+R,\,1+R^2,\,\ldots,\,1+
R^{2^{\esup{\log{n}}}}$\, eux-m\^emes
obtenus \`a l'issue de \,$\esup{\log{n}}$\,
\'el\'evations successives au carr\'e, tronqu\'ees
\`a l'ordre $n$, de \pols de degr\'e \,$n$.
Cela fait un nombre d'\oparis de l'ordre de
\,$n^2\esup{\log{n}}$.

\ss
$\bullet$ {\bf L'\'etape 4} a un co\^ut 
n\'egligeable.

\begin{proposition}
\label{propchistoseq}
Le co\^ut total de l'\algo \sql \elr
de Chistov s'\'el\`eve \`a
 \,$(2/3)\,n^4+\O(n^3)$\,  \oparis dans l'anneau de
base.
\end{proposition}

\ss Nous pr\'esentons ci-dessous un r\'esum\'e de la
discussion
sur le nombre d'\oparisz:

\begin{center}
\label{seqchis}
\begin{tabular}{|l|c|c|}
\hline Etape~ & Co\^ut \\
\hline Etape~1 & $(2/3)\,n^4+\O(n^3)$ \\
\hline Etape~2 & $\O(n^{3})$ \\
\hline Etape~3 & $\O(n^{2}\log{n})$ \\
\hline Etape~4 & n\'egligeable \\
\hline
\end{tabular}

\ms {\bf {\small Tableau  \ref{seqchis}}} \\
{\small Complexit\'e de la version
\sqle  de l'\algo de Chistov}
\end{center}

\section{M\'ethodes reli\'ees aux \srls}
\label{sec algo-srl}

Dans la section \ref{subsec.Melr} nous donnons un
\algo de
calcul du
\polcar d'une matrice \,$A$\, bas\'e sur la
consid\'eration
des
transform\'es successifs de vecteurs de la base
canonique
par \,$A$.

Dans la section \ref{subsec.BerMa} nous pr\'esentons
un
\algo
d\^{u} \`a Berlekamp qui permet de calculer le
\polgmin
d'une \srl dans un corps
lorsqu'on sait qu'elle v\'erifie une \rrl d'ordre
\,$n$\,
et qu'on conna\^{\i}t les \,$2n$\, premiers termes de
la
suite.

Dans la section \ref{subsec.wied} on d\'ecrit l'\algo
de
Wiedemann
qui utilise celui de Berlekamp pour trouver avec une
bonne
probabilit\'e le \polcar d'une matrice sur un corps
fini.

\subsection{L'\algo de Frobenius}  \label{subsec.Melr}

Nous donnons ici un \algo
qui est bas\'e sur une  description de nature \geoq
pour un \endom d'un \Kevz.

Comme cons\'equence, on calcule
le \polcar de l'\endom avec essentiellement le
m\^{e}me
nombre d'\oparis que dans la \mpg (qui ne calcule que
le
\deterz), sans les inconv\'enients
que pr\'esentait la \mhb (hormis le cas des
corps finis) concernant la taille des
\coes \itmdsz.

\subsection*{Le cas usuel}

Nous aurons besoin de la proc\'edure \vref{ajbsol}
(d\'eriv\'ee de
l'\ajbz) \`a laquelle nous donnons le nom de
\textbf{JorBarSol}.
Elle calcule, \`a la \JBz, la \rdl exprimant la
derni\`ere colonne
en fonction des premi\`eres dans une matrice \freg
ayant une
colonne de plus que de lignes. La fin du calcul,
apr\`es la
triangulation, reste dans l'anneau \,$\A$\, si la \rdl
est \`a
\coes dans \,$\A$.

\ss Consid\'erons une \maca \,$A\in\ZZ^{5\times 5}$\,
d'ordre 5
prise au hasard,
donn\'ee par exemple par {\sc Maple}.
Elle d\'efinit un \endom \,$h_{\!A}$\, de \,$\QQ^5$.
On note \,$(f_i)_{1\leq i\leq 6}$\, le premier vecteur
de
la base canonique de \,$\QQ^5$\, et ses 5
transform\'es
successifs par
\,$A$.  Ceci fournit une matrice \,$B\in\ZZ^{5\times
6}$.

Voici un exemple typique
{\footnotesize
$$
{\begin{array}{cc}
A\,=\,\left[
{\begin{array}{rrrrrrrr}
57 & -82 & -48 & -11 & 38 \\
-7 & 58 & -94 & -68 & 14 \\
-35 & -14 & -9 & -51 & -73 \\
-73 & -91 & 1 & 5 & -86 \\
43 & -4 & -50 & 50 & 67
\end{array}}
 \right]\,,
 \\~
 \\
B\,=\,\left[
\begin {array}{rrrrrrrr}
1 & 57 & 7940 & 55624 & -46831857 & -22451480858 \\
\noalign{\smallskip}
0 & -7 & 8051 & 1071926 & 199923276 & 14745797441 \\
\noalign{\smallskip}
0 & -35 & -998 & -245490 & 54032957 & 9123769947 \\
\noalign{\smallskip}
0 & -73 & -7622 & -1648929 & -128141849 & -10372211183
\\
\noalign{\smallskip}
0 & 43 & 3460 & 209836 & -58008810 & -15808793525
\end {array}\right]\,.
\end{array}}
$$}

\normalsize
En \gnl les vecteurs  \,$(f_i)_{1\leq i\leq 5}$\,
sont ind\'ependants, et m\^{e}me, la matrice \,$B$\,
est
\fregz. C'est le cas ici.

Concernant la taille des \coesz, ceux de la matrice
initiale
sont major\'es par \,$100$\, en valeur absolue, et
ceux de
la matrice \,$B$\, dans la \,$k$\,-\`eme colonne,
sont major\'es par \,$M^{k-1}$\, si \,$M$\, est une
des
normes de \,$A$\, d\'ecrite en section
\ref{secHadamard}, par exemple
\,$M<270$\, pour la norme de Frobenius.

\begin{algor}[Algorithme JorBarSol]\label{ajbsol}
\acl{ajbsol}{Algorithme JorBarSol}
\Entree Une matrice $A=(a_{ij})\in\A^{n\times (n+1)}$
\fregz. L'anneau $\A$ est suppos\'e int\`egre avec un
\algo
de division exacte.
\Sortie $L=(\ell_j)\in\A^n$: on a, en notant $C_j$
la $j$-\`eme colonne de $A$, $C_{n+1}=\sum_{j=1}^n
\ell_j\,
C_j$. La fin du calcul reste dans $\A$ si les $\ell_i$
sont
dans $\A$.
\Debut
\Varloc $i$, $j$, $p$, $q\in\N$; $piv$, $den$,
$coe\in\A$;
\hsc        \# \LU-\deco \`a la \JB
\hsu $den:= 1$ ;  $m:=n+1$;
\hsu \pour{p}{1}{n-1}
\hsd    $piv:=a_{pp}$ ;
\freg
\hsd      \pour{i}{p+1}{n}
\hst         $coe:=a_{ip}$;
\hst         \pour{j}{p+1}{m}
\hsq            $a_{ij}:=(piv*a_{ij}-coe*a_{pj})\,/ \,
den$
\hst         \finpour
\hsd      \finpour
\hsd    $den:=piv$
\hsu \finpour
\hsc        \# calcul des \coes $\ell_i$
\hsu \pour{q}{1}{n-1}
\hsd     $p:=n-q$;   $\ell_p:=a_{pm}/a_{pp}$;
\hsd     \pour{i}{1}{p-1}
\hst       $a_{im}:=a_{im}-\ell_p\,a_{ip}$
\hsd     \finpour
\hsu \finpour
\fin
\end{algor}

On peut calculer la \rdl
\,$f_6=\sum_{i=1}^5\alpha_i\,f_i$, qui se
relit
\,${h_{\!A}}^5(f_1)=\sum_{i=0}^4\alpha_i\,{h_{\!A}}^i(f_1)$.
La matrice de \,$h_{\!A}$\, sur la base
\,$(f_1,\ldots, f_5)$\,
est alors clairement la \maco du \pol
\,$P(X)=X^5-\sum_{i=0}^4\alpha_i\,X^i$\,
 (\cf page \pageref{matrcomp}) et
on obtient le \polcar de \,$A$\, par la formule
$\rPA(X)=(-1)^5P(X)$.

Cet \algo qui calcule le \polcar de \,$A$\, fonctionne
lorsque le
\polgmin \,$P$\, de la \srl
\,$(A^n\,f_{1})_{n\in\NN}$\, dans
\,$\QQ^5$\, est de degr\'e \,$\geq 5$. Dans ce cas le
\pol \,$P$\,
est en effet \'egal au \polmin et au \polcar de
\,$A$\, (au signe
pr\`es).

Pour calculer la \rdl 
\,$f_6=\sum_{i=1}^5\alpha_i\,f_i$\, on
applique la \pcd  \textbf{JorBarSol}. Elle commence
par la
triangulation \`a la \JB de la matrice \,$B$. Cette
triangulation
(en fait, une \LU-\decoz) ne change pas les deux
premi\`eres
colonnes et donne les 4 derni\`eres suivantes.
{\footnotesize
$$
\begin {array}{rrrrrr}
7940 & 55624 & -46831857 & -22451480858 \\
\noalign{\smallskip}
8051 & 1071926 & 199923276 & 14745797441 \\
\noalign{\smallskip}
288771 & 39235840 & 6619083961 & 452236520806 \\
\noalign{\smallskip}
641077 & -110921281313 & -32874613863452 &
-5984786805270056
\\
\noalign{\smallskip}
-370413 & -114147742050 & -8244227015780803785 &
-1467472408808983073730
\end {array}
$$}

\normalsize
Si on traite une \maca \,$A$\, d'ordre \,$n$\,
dont une norme est major\'ee par \,$M$,
le \coe en position \,$(i,j)$\, dans la matrice ainsi
obtenue est \'egal au mineur \,$b^{(k-1)}_{i,j}$\, de
la
matrice \,$B=[f_1|Af_1|\cdots|A^nf_1]$\, avec
\,$k=\min(i,j)$.
A priori (comme dans l'exemple ci-dessus, d'ailleurs)
le
plus grand \coe serait
en position \,$(n,n+1)$, major\'e
par \,$M^{1+\cdots+(n-2)+n}$, \cad
\,$M^{(n^2-n+2)/2}$, ce qui reste
raisonnable, en tout cas bien meilleur que dans
l'\aghbz.

L'\algo termine en donnant la \coli recherch\'ee par
le
calcul successif des \coes
\,$\alpha_5,\,\alpha_4,\ldots,\alpha_1$.

\ss Une \emph{matrice de Frobenius d'ordre \,$n$\,}
\index{Frobenius!matrice
de} \index{matrice!de Frobenius}
 est un autre nom donn\'e \`a
une \maco
d'un \poly \,$P(X)$\, de degr\'e \,$n$.
On peut l'interpr\`eter comme
la matrice de l'application \lin \gui{\mul par \,$x$}
(la classe de \,$X$) dans l'\agr quotient
\,$\aqo{\K[X]}{P(X)}$\,
sur la base canonique \,$1,\,x,\,\ldots,\,x^{n-1}$.

L'\algo que nous venons de d\'ecrire n'ayant pas de nom
officiel, nous l'appellerons \emph{\algo de 
Frobenius},\index{Frobenius!algorithme de} c'est
l'\algo \vref{afrob}.

\begin{algor}[Algorithme de Frobenius (le cas
simple)]\label{afrob}
\acl{afrob}{Algorithme de Frobenius.}
\Entree Une matrice $\,A=(a_{ij})\in\A^{n\times n}$.
L'anneau $\A$ est suppos\'e int\`egre avec un \algo
de division exacte.
\Sortie Le \polcar \,$\rPA(X)$.
L'\algo ne fonctionne que si le premier vecteur de
base et
ses \,$n-1$\, transform\'es successifs par \,$A$\,
sont lin\'eairement ind\'ependants.
\Debut
\Varloc $i$, $k$, $m\in\N$;
$B=(b_{i,j})\in\A^{n\times (n+1)}$;
\hsz
$V=(v_i)\in\A^{n\times 1}$, $L=(\ell_j)\in\A^n$;
\hsu  $m:= n+1$; $V:=$ premi\`ere colonne de $A$;
\hsu  $B:=0$ dans $\A^{n\times (n+1)}$; $b_{1,1}:=1$;
\hsu  $2$-\`eme colonne de $B\,:=\,V$;
\hsu  \pour{k}{3}{m}
\hsd      $V:=A\,V$;
\hsd  $k$-\`eme colonne de $B\,:=\,V$
\hsu  \finpour;
\hsu  $L:=\textsf{JorBarSol}(B)$;
\hsu  $\rPA:=(-1)^n(X^n-\sum_{k=1}^n \ell_k\,X^{k-1})$
\fin
\end{algor}

\subsubsection*{Nombre d'\oparis }

Avec une \maca d'ordre \,$n$\, l'\algo de Frobenius
donne un calcul en
\,$\O(n^3)$\, \oparisz,
ce qui est du m\^{e}me ordre de grandeur que pour la
\mpgz.

Plus pr\'ecis\'ement l'\algo \ref{afrob} se
d\'ecompose en deux grandes
\'etapes. La premi\`ere \'etape cr\'ee la matrice
\,$B$\, et ex\'ecute
\,$(n-1)\,n^2$\, \muls et \,$(n-1)^3$\, additions.
La deuxi\`eme \'etape applique l'\algo
\textbf{JorBarSol} \`a la matrice
\,$B$. Vue la premi\`ere colonne de celle-ci, cet
\algo
utilise
$$ \sum_{p=2}^{n-1} (n-p)\,(n-p+1)+\sum_{p=1}^{n-1}
(p-1)=
{\displaystyle \frac {1}{3}} \,n^{3} -
{\displaystyle \frac {1}{2}} \,n^{2} - {\displaystyle
\frac {5}{6
}} \,n + 1
$$
\adsos et
$$ 3\,\sum_{p=2}^{n-1} (n-p)\,(n-p+1)+\sum_{p=1}^{n-1}
(p-1)+n-2=
n^{3} - {\displaystyle \frac {5}{2}} \,n^{2} +
{\displaystyle \frac {3}{2}} \,n - 1
$$
\mulsz/divisions (les divisions sont toutes exactes).
Ceci donne le r\'esultat suivant.
\begin{proposition}
\label{propFrobCompt} L'\algo de Frobenius (dans le
cas \usl
simple) appliqu\'e \`a une \maca d'ordre \,$n$\, sur
un anneau
int\`egre dans lequel les divisions exactes sont
explicites
demande en tout $${\displaystyle \frac {10}{3}}
\,n^{3} - 7\,n^{2}
+ {\displaystyle \frac {11}{3}} \,n - 1$$
 \oparis dans l'anneau.
Plus pr\'ecis\'ement cet \algo ex\'ecute
\,${\displaystyle \frac {4}{3}} \,n^{3} -
{\displaystyle \frac {7}{2}} \,n^{2} + {\displaystyle
\frac {13}{
6}} \,n$\, \adsos et
\,$2\,n^{3} - {\displaystyle \frac {7}{2}} \,n^{2} +
{\displaystyle \frac {3}{2}} \,n - 1$\,
\mulsz/divisions.
\end{proposition}

En pratique, sur un corps fini, les \algos de
Hessenberg et
de Frobenius s'av\`erent meilleurs que tous les
autres, ce qui correspond au fait qu'ils fonctionnent en
ex\'ecutant seulement \,$\O(n^3)$\, \oparisz. Mais d\`es qu'on
passe \`a des matrices \`a \coes dans \,$\ZZ$, l'\agb devient
plus performant, car ses \,$\O(n^4)$\, \oparis sont
ex\'ecut\'ees sur des entiers de taille mieux contr\^{o}l\'ee. 
Si on passe \`a des anneaux tels que \,$\ZZ[t,u]$, l'\agb
b\'en\'eficie en plus du fait qu'il n'utilise pas de divisions.
Enfin sur des anneaux non int\`egres, les \algos de Hessenberg
et de Frobenius, m\^{e}me dans leurs variantes avec
recherche de pivot non nuls, ne fonctionnent plus en
toute \gnlz it\'e.
\subsection*{Le cas difficile: triangularisation par
blocs}
\label{subsec.TriFro}\label{sec.TriFro}

La \met que nous d\'ecrivons maintenant est
l'adaptation de la pr\'ec\'edente pour le cas 
le plus difficile, qui se pr\'esente
cependant rarement. 
Cette \metz, comme celle d\'ecrite
pour le cas \usl (celui o\`u le \polcar de  \,$A$\, est \'egal \`a
son \polmin et o\`u le premier vecteur de base 
\,$A$\,-engendre l'espace
\,$\K^n$) fait partie de l'usage, et nous ne savons
pas \`a qui l'attribuer.

Soit \,$A$\, une \maca dans \,${\K}^{n\times n}$.
Notons \,$a=(e_{1},\ldots,e_{n})$\, la base canonique
de \,${\K}^{n}$\, (on identifiera \,${\K}^{n\times 1}$
avec \,${\K}^{n}\,$) et \,$h_{A}$\, l'\endom de
\,${\K}^{n}$\, ayant pour matrice \,$A$\, dans cette
base.

Nous allons construire une nouvelle base
\,$b=(f_{1},\ldots,f_{n})$ dans laquelle l'\endom
\,$h_{\!A}$\,
aura une matrice suffisamment sympathique, dont le
\polcar sera
facile \`a calculer. Il s'agit pr\'ecis\'e\-ment de
r\'eduire la
matrice de \,$h_{\!A}$\, \`a une forme triangulaire
par blocs avec
des blocs diago\-naux ayant la forme de Frobenius.

Voyons ce qui se passe sur un exemple.
\subsubsection*{\bf Un exemple dans \,$\K^7=\QQ^7$:}
Soit $(e_{1},e_{2},e_{3},e_{4},\alb e_{5},\alb
e_{6},\alb
e_{7})$\, la base canonique de \,${\K}^{7}$\,, \,$v$\,
un vecteur
de \,${\K}^{7}\,$ et \,$A$\, la \maca d'ordre 7
ci-dessous.
Rappelons que le sous-espace de Krylov\footnote{~On
dit parfois aussi \emph{sous espace cyclique}.} 
associ\'e au couple
\,$(A,v)$\, (ou au couple \,$(h_{\!A},v)\,$), not\'e
\,$\Kr_{A,\,v}$, est le sous-espace de \,${\K}^{7}$\,
engendr\'e par la \srl 
\,$(A^n\,v)_{n\in\NN}$.
{\footnotesize $$
 A=\left
[\begin {array}{rrrrrrr}
1&5&43&683&794&206&268 \\
\noalign{\smallskip}
-1&-2&-26&-458&-554&-148&-186 \\
\noalign{\smallskip}
0&0&1&14&18&5&6 \\
\noalign{\smallskip}
1&3&24&387&469&125&157 \\
\noalign{\smallskip}
1&4&21&300&357&92&119 \\
\noalign{\smallskip}
2&5&53&888&1082&292&363 \\
\noalign{\smallskip}
-7&-23&-163&-2547&-3074&-813&-1028
\end {array}\right]\,, 
$$}

\normalsize Pour obtenir une base du sous-espace
\,$\Kr_{A,\,e_1}$\, (de dimension \,$k_{1}$) nous
calculons
successivement les vecteurs \,$e_{1}$, $Ae_{1}$,
$A^2e_{1}$,
\ldots\, en nous arr\^etant au dernier vecteur qui ne
soit pas \coli de ceux qui le pr\'ec\`edent ou, ce
qui revient
au m\^eme, nous construisons successivement les
matrices
$[e_{1}]$, $[e_{1}\,|\,Ae_{1}]$,
$[e_{1}\,|\,Ae_{1}\,|\,A^2e_{2}]\,,\dots$\, en nous
arr\^etant \`a
la derni\`ere matrice dont le rang est \'egal au
nombre de
colonnes. Dans notre cas, cela donne la suite de
matrices:
{\footnotesize $$
 \left
[\begin {array}{c}
1\\
0\\
0\\
0\\
0\\
0\\
0\\
\end {array}\right ]\,,\,\,
\left [\begin {array}{cr} 1&1\\
0&-1\\
0&0\\
0&1\\
0&1\\
0&2\\
0&-7
\end{array}\right]\,,\,\,
\left [\begin {array}{crr} 1&1&9\\
0&-1&-5\\
0&0&0\\
0&1&5\\
0&1&5\\
0&2&10\\
0&-7&-35
\end {array}\right]\,,~~\ldots
$$
}

\noi\normalsize Il faut s'arr\^eter \`a la deuxi\`eme
matrice car
la matrice \,$[\,e_{1}\,|\,Ae_{1}\,|\,A^2e_{1}\,]$\,
est de rang
$2$. On remarque en effet que \,$A^2e_{2}=4\,e_{1} +
5\,Ae_{1}$.
Une base du sous-espace \,$V_{1}=\Kr_{A,\,e_1}$\, est
donc
form\'ee du couple \,$(e_{1},Ae_{1})$\, et
\,$\dim{V_{1}}=k_{1}=2$. La matrice correspondante est
not\'ee
\,$U_1=[\,e_{1}\,|\,Ae_{1}\,]$.

\ss Passons au second vecteur de la base cano\-nique.
On remarque que \,$e_{2}$\, n'est pas dans le
sous-espace
\,$V_{1}$\, et l'on poursuit la construction de la
base
recherch\'ee avec les matrices
$[\,U_{1}\,|\,e_{2}\,]$ {puis}
$[\,U_{1}\,|\,e_{2}\,|\,Ae_2\,]$ {puis}
$[\,U_{1}\,|\,e_{2}\,|\,Ae_2\,|\,A^2e_2\,]$ \ldots\,\,
 jusqu\`a obtenir une matrice dont la derni\`ere
colonne est \coli des autres.

Ici c'est le vecteur \,$A^3e_{2}$\, qui est
combinaison
\lin de ceux qui le pr\'ec\`edent (\cad qu'il
appartient
au sous-espace \,$V_{2}=V_{1}+\gen{e_2,Ae_2,A^2e_2}$)
et on obtient la  suite de  matrices:
{\footnotesize
$$
\left [\begin {array}{ccc} 1&1&0\\
0&-1&1
\\
0&0&0\\
0&1&0\\
0&1&0\\
0&2&0\\
0&-7&0\end {array}\right ]~,~~
\left [\begin {array}{cccc} 1&1&0&5\\
0&-1&1&-2
\\
0&0&0&0\\
0&1&0&3
\\
0&1&0&4\\
0&2&0&5
\\
0&-7&0&-23\end {array}\right ]~,~~
U_2=\left [\begin {array}{ccccc}
1&1&0&5&86\\
0&-1&1&-2&-53\\
0&0&0&0&1\\
0&1&0&3&50\\
0&1&0&4&48\\
0&2&0&5&103\\

0&-7&0&-23&-347\end {array}\right]\,.
$$}

\noi \normalsize Et la matrice
$[\,e_{1}\,|\,Ae_{1}\,|\,e_{2}\,|\,Ae_2\,|\,A^2e_2\,|
\,A^3e_2\,]= [\,U_{2}\,|\,A^3e_2\,]$:
{\footnotesize $$
~~~~\left [\begin {array}{cccccc}
1&1&0&5&86&348\\
\noalign{\smallskip}
0&-1&1&-2&-53&-200\\
\noalign{\smallskip}
0&0&0&0&1&-2\\
\noalign{\smallskip}
0&1&0&3&50&209\\
\noalign{\smallskip}
0&1&0&4&48&214\\
\noalign{\smallskip}
0&2&0&5&103&411\\
\noalign{\smallskip}
0&-7&0&-23&-347&-1471
\end {array}\right]
$$}

\noi\normalsize est de rang 5 puisque sa derni\`ere
colonne
\,$A^3e_{2}$\, est combinaison des autres. On peut le
voir par
exemple par la \mpgz, qui fournit la \rdl
\,$A^3e_2=209e_1+306Ae_1+2e_2+Ae_2-2A^2e_2$. On passe
ensuite au
troisi\`eme vecteur de la base canonique. On remarque
que
\,$e_{3}$\, n'est pas dans le sous-espace \,$V_{2}$.
On construit
alors une base de \,$V_{3}=V_{2}+\Kr_{A,\,e_3}$. On
poursuit donc
la construction de la base recherch\'ee avec les
nouvelles
matrices $[\,U_{2}\,|\,e_3\,]$ et
$[\,U_{2}\,|\,e_3\,|\,Ae_3\,]$:
{\footnotesize
$$
\left [\begin {array}{cccccc}
1&1&0&5&86&0\\
\noalign{\smallskip}
0&-1&1&-2&-53&0\\
\noalign{\smallskip}
 0&0&0&0&1&1\\
\noalign{\smallskip}
0&1&0&3&50&0\\
\noalign{\smallskip}
0&1&0&4&48&0\\
\noalign{\smallskip}
0&2&0&5&103&0\\
\noalign{\smallskip}
0&-7&0&-23&-347&0
\end {array}\right]\,,~
\left [\begin {array}{ccccccc}
1&1&0&5&86&0&43\\
\noalign{\smallskip}
0&-1&1&-2&-53&0&-26\\
\noalign{\smallskip}
0&0&0&0&1&1&1\\
\noalign{\smallskip}
0&1&0&3&50&0&24\\
\noalign{\smallskip}
0&1&0&4&48&0&21\\
\noalign{\smallskip}
0&2&0&5&103&0&53\\
\noalign{\smallskip}
0&-7&0&-23&-347&0&-163
\end {array}\right]
$$}

\noi\normalsize
La derni\`ere matrice, que nous notons \,$U_3=U$, est
de
rang 7. C'est la matrice de passage de
la base canonique \`a la base que nous venons de
construire
\,$b = (e_{1}, \alb Ae_{1},\alb e_{2},Ae_{2},\alb
A^2e_{2},e_{3},Ae_{3})$.

Dans cette
nouvelle base, il est clair que la matrice de l'\endom
\,$h_{\!A}$\, est une matrice \tgu par blocs, les
blocs
diagonaux \'etant form\'es de matrices de Frobenius.

On peut d'ailleurs le v\'erifier, en calculant le
produit
matriciel \,$U^{-1}A\,U$\, pour obtenir:
{\footnotesize
$$
U^{-1}\,A\,U\,=\,
\left [\begin {array}{ccccccc} 0&4&0&0&209&0&179\\
\noalign{\smallskip}
1&5&0&0&306&0&291\\
\noalign{\smallskip}
0&0&0&0&2&0&6\\
\noalign{\smallskip}
0&0&1&0&1&0&-17\\
\noalign{\smallskip}
0&0&0&1&-2&0&-4\\
\noalign{\smallskip}
0&0&0&0&0&0&1\\
\noalign{\smallskip}
0&0&0&0&0&1&5
\end {array}\right]
$$}

\noi\normalsize
dont le \polcar (celui aussi de
\,$A\,$) est \'egal au produit des \polcars
des blocs diagonaux de
Frobenius, \cad :
\,$\left({X}^{2}-5\,X-4\right)
\left({X}^{3}+2\,{X}^{2}-X-2\right)
\left({X}^{2}-5\,X-1\right).$

En fait la matrice \,$U^{-1}A\,U$,
et par suite les \polcars des blocs diagonaux de
Frobenius,
peuvent \^etre retrouv\'es \`a partir des \rdls
d\'ej\`a
calcul\'ees
et de la relation qui exprime le vecteur
\,$A^2e_{3}$\,
comme
\coli des vecteurs de la base \,$b$, qui peut \^{e}tre
obtenue en
appliquant la \mpg \`a la matrice
\,$[\,U_3\,|\,A^2e_{3}\,]$.

\subsubsection*{Description g\'en\'erale de
l'algorithme}

On prend \,$f_{1}=e_{1}$\, puis \,$f_{2}=Ae_{1}$,
sauf si \,$Ae_{1}$\, est co\lin avec \,$e_{1}$, auquel
cas on prend \,$f_{2}=e_{2}$.

Pr\'ecis\'ement, on d\'efinit l'entier \,$k_{1}\in
\{1,\ldots,n\}$
comme suit:
les vecteurs \,$e_{1},$ $Ae_{1},\ldots,$ $A^{k_{1}-
1}e_{1}$\, sont
ind\'e\-pen\-dants, mais \,$A^{k_{1}}e_{1}$\, d\'epend
\lint des pr\'ec\'edents. Ceci d\'efinit le d\'ebut
$$
(f_{1},\ldots,f_{k_{1}})=(e_{1},Ae_{1},\ldots,A^{k_{1}-
1}e_{1})
$$
de notre base.

Les tests de d\'ependance \lin dont nous avons eu
besoin
peuvent \^etre obtenus en appliquant le pivot de
Gauss, avec
\'eventuels \'echanges de lignes mais sans \'echange
de
colonnes,
sur les matrices successives  \,$(e_{1},Ae_{1})$,
\,$(e_{1},Ae_{1},A^2e_{1})$\, etc.
Cette \met fournit aussi la relation de d\'ependance
\lin lorsque l'entier \,$k_{1}$\, est atteint.
Notez aussi que nous n'avons pas besoin de calculer
les
puissances successives de la matrice \,$A$\, mais
seulement
les transform\'es successifs du vecteur  \,$e_{1}$\,
par
\,$A$.

\smallskip Si \,$k_{1}=n$\, notre base \,$b$\, est
trouv\'ee,
et en exprimant \,$Af_{k_{1}}$\, sur la base \,$b$\,
nous
obtenons en m\^eme temps la matrice de $h_{\!A}$ sur
\,$b$\,
sous forme d'une matrice de Frobenius.
$$
\left [\begin {array}{ccccccc}
  0       &  \cdots   &  0   &   a_{0} \cr
  1       &  \ddots   &  \vdots    &   a_{1}      \cr
  \vdots  &  \ddots   &  0    &    \vdots     \cr
  0       &  \cdots   &  1   &   a_{n-1}
\end {array}\right]
$$
dont le \polcar est, au signe pr\`es,
$$
P(X)=X^n-(a_{n-1}X^{n-1}+\cdots+a_{1}X + a_{0})
$$

\smallskip Si \,$k_{1}<n$, nous cherchons le premier
vecteur \,$e_{i}$\, ($i>1$) \lint ind\'ependant
de \,$f_{1},\ldots,f_{k_{1}}$. Ceci nous fournit le
vecteur \,$f_{k_{1}+1}$.
Le calcul de l'indice $i$ et donc du vecteur
\,$f_{k_{1}+1}$\, peut de nouveau \^etre obtenu par la
m\'ethode du pivot de Gauss (sans \'echange de
colonnes)
appliqu\'ee aux matrices
\,$(f_{1},\ldots,f_{k_{1}},e_{i})$.
On d\'efinit ensuite l'entier $k_{2}\in \{1,\ldots,n-
k_{1}\}$
comme suit: les vecteurs
$$
f_{1}, \ldots, f_{k_{1}+1},Af_{k_{1}+1}, \ldots,
A^{k_{2}-1}f_{k_{1}+1}\,
$$ sont ind\'ependants, mais
\,$A^{k_{2}}f_{k_{1}+1}$\, d\'epend \lint
des pr\'ec\'edents. Ceci d\'efinit le nouveau d\'ebut
de notre base,
$$
(f_{1},\ldots,f_{k_{1}+k_{2}})=
(f_{1},\ldots,f_{k_{1}+1},Af_{k_{1}+1},\ldots,
A^{k_{2}-1}f_{k_{1}+1}).
$$

\smallskip Si \,$k_{1}+k_{2}=n$\, notre base \,$b$\,
est
trouv\'ee,
et en exprimant \,$Af_{k_{1}+k_{2}}$\,
sur la base \,$b$\, nous obtenons en m\^eme temps la
matrice
de \,$h_{\!A}$\, sur \,$b$\, sous forme d'une \matg
par blocs, ayant pour
blocs diagonaux deux matrices de Frobenius.

$$
\left[
  \begin {array}{cccccccccc}
            0       &  \cdots  &  0   &   a_{0} & 0   
  &\cdots &\cdots
&  0   &   c_{0} \cr
            1       &  \ddots  &\vdots&   a_{1}
&\vdots  &       &
&\vdots&   c_{1} \cr
            \vdots  &  \ddots  &  0   &  \vdots
&\vdots  &       &
&\vdots&  \vdots \cr
            0       &  \cdots  &  1   &  a_{n-1}& 0   
  & \cdots&\cdots
&  0   &   c_{n-1}\\[1mm]
            0       &  \cdots  &\cdots&   0     & 0   
  &\cdots &\cdots
&  0   &   b_{0} \cr
            \vdots  &          &      &   \vdots&1    
  &\ddots &
&\vdots&   b_{1}\cr
            \vdots  &          &      &  \vdots &0    
  &\ddots &\ddots
&\vdots&  \vdots \cr
            \vdots  &          &      &  \vdots
&\vdots  &\ddots &\ddots
&  0   &  \vdots \cr
            0       &  \cdots  &\cdots&   0     & 0   
  &\cdots &     0
&  1   &   b_{m-1}
\end {array} \right]
$$
dont le \polcar est, au signe pr\`es,
$$
P(X)=\left(X^n-\sum\nolimits_{i=0}^{n-1}{a_{i}X^{i}}\right)
\cdot\left(X^m-\sum\nolimits_{j=0}^{m-1}{b_{j}X^{j}}\right)\,.
$$

\smallskip Si \,$k_{1}+k_{2}<n$\, nous cherchons le
vecteur
\,$f_{k_{1}+k_{2}+1}$\, parmi les vecteurs restants de
la
base canonique, et nous continuons le processus.
En fin de compte, en ayant calcul\'e un nombre
relativement
restreint (certainement $<2n$) de produits du type
matrice
fois vecteur \,$Ag$, et en ayant appliqu\'e le pivot
de
Gauss un nombre relativement restreint de fois nous
avons
obtenu une nouvelle base \,$b$\, ainsi que la matrice
de
\,$h_{\!A}$\, sur cette base sous la forme d'une \matg
par
blocs, ayant sur la diagonale des blocs form\'es de
matrices
de Frobenius.
Le \polcar de la matrice est donc \'egal au produit
des
\polcars \,$P_{i}(X)$\, des blocs diagonaux, qui sont
donn\'es par simple lecture de la derni\`ere colonne
de la
matrice de Frobenius.

\bigskip Notons pour terminer qu'il est facile de
v\'erifier
sur une telle forme r\'eduite que chacun des vecteurs
\,$f_{j}$\,
est annul\'e par l'\endom
\,$\prod_{i}{P_{i}(h_{\!A})}$,
ce qui fournit une preuve \geoq \elr
du \tho de Cayley-Hamilton.
Pour la preuve de ce \tho il suffit d'ailleurs
de constater le fait pour le vecteur \,$f_{1}$, car
celui-ci
est simplement le premier vecteur d'une base, et donc
n'importe quel vecteur non nul a priori.

\subsubsection*{Domaine de validit\'e et nombre
d'\oparis}
Dans cet \algo le nombre d'\oparis est encore un
\,$\O(n^3)$.

Son domaine de validit\'e est celui des corps, et plus
\gnlt celui des
anneaux int\`egres et
int\'egra\-le\-ment clos, que nous avons envisag\'es
\`a l'occasion de
l'\'etude du \polmin (voir section \vref{pageAIC}).

En effet avec un tel anneau, si
\,$C\in\A^{n{\times}n}$\, les
\pols \,$\rP^{C}$\, et \,$\rP^{C,v}$\, sont
automatiquement \`a
\coes dans \,$\A$. Il s'ensuit que les proc\'edures du
type
\textbf{JorBarSol} que nous utilisons au cours de
l'\algo ne
calculent que des \elts dans \,$\A$.

Dans le cas d'une matrice \`a \coes dans
\,$\ZZ$\, on a les m\^{e}mes majorations des \coes
\itmds
que celles que nous avons esquiss\'ees dans le cas
facile le
plus \uslz.

\subsection{Algorithme de
Berlekamp/Massey}\label{subsec.BerMa}

On donne dans un corps \,$\K$\, les  \,$2n$\,
premiers \elts d'une \srl \,$(a_k)_{k\in\NN}$\,
pour laquelle on sait qu'il existe un \polg de
degr\'e \,$n$. Le probl\`eme est de calculer
le \polgmin \,$g$\,  de la suite.

Une telle solution est donn\'ee par l'\algo de
Berlekamp/Massey \cite{Dor} qui donne en sortie le
degr\'e
\,$d$\, ainsi que les \coes d'un \poly \,$f=c_d\,g$\,
associ\'e  au \pol \,$g$.
Ce \pol \,$g$\, est alors obtenu en divisant par
\,$c_d$.
\index{Berlekamp/Massey!algorithme de}

L'\algo de Berlekamp/Massey utilise les propri\'et\'es
de la suite des triplets \,$(R_i,U_i,V_i)$\, form\'ee
des restes et des multiplicateurs de B\'ezout
successifs
dans l'\algo d'Euclide \'etendu pour le couple de
\pols
\,$(R_{-1},R_0)$\, o\`u \,$R_{-1}=X^{2n}$\, et \,$R_0
=
\sum_{i=0}^{2n-1} a_iX^i$.

Posant \,$V_{-1}=U_0=0$\, et \,$U_{-1}=V_0=1$,  ces
triplets
v\'erifient, pour tout
\,$i\geq 0$, les relations:
$$\begin{array}{rrclclr}
&R_{i-1}&=&R_i\,Q_i+R_{i+1}&\mbox{o\`u }\dg R_{i+1}
< \dg R_{i} \\
&         U_{i+1}       &=& U_{i-1} - Q_i\,U_i,\\
&         V_{i+1}       &=& V_{i-1} - Q_i\,V_i,&\mbox{
d'o\`u:}  \\
&         R_i           &=& U_i\,R_{-1} + V_i\,R_0, \\
\mbox{de plus:}&U_i\,V_{i-1} - V_i\,U_{i-1}&=& (-1)^i
\qquad \mbox{  et}&  \dg R_{i} < 2n - \dg V_{i}\,.
\end{array}$$

Les deux derni\`eres relations se v\'erifient
facilement par \recu sur~$\,i$. 

On arr\^ete le processus au premier reste,
disons \,$R_m$, de degr\'e plus bas que
\,$n$, pour obtenir:
$$\,U_{m}\,X^{2n} +
V_m\,R_0=R_m~~\mbox{ avec }~~\dg R_m<n.$$

Posons \,$d=\sup(\dg V_m, 1+\dg R_m)$\, et $P=X^dV_{m}(1/X)$.
Alors on peut montrer que \,$P$\, divis\'e par son \coe dominant 
est le \polgmin de la suite \,$(a_k)$
(cf. \cite{VonZurbook} et \cite{Dor}).
Par exemple dans le cas o\`u \,$\dg V_m=n$\, et \,$V_{m}(0)\neq 0$, 
en \'ecrivant que les termes de degr\'e
compris entre \,$n$\, et \,$2n-1$\, du \poly
\,$V_m(X)\,R_0(X)$\, sont nuls, on constate que \,$P(X)$\, est bien un
\polg de la suite \,$(a_k)$.

Ceci donne pr\'ecis\'ement l'\algo \vref{algoBMA} (dans lequel
\,${\rm cd}(P)$\, d\'esigne le \coe dominant de \,$P$).

\begin{algor}[Algorithme de Berlekamp-Massey]\label{algoBMA}
    \acl{algoBMA}{\Algo de Berlekamp-Massey}
\Entree Un entier $n\geq 1$. Une liste non nulle d'\'el\'ements du
corps $\K$, $[a_0,a_1,\ldots,a_{2n-1}]$: les $2n$ premiers termes
d'une \srlz, sous l'hypoth\`ese qu'elle admet un \polg de degr\'e
$\leq n$. 
\Sortie Le \polgmin $P$ de la \srlz. 
\Debut
\Varloc $R,R_0,R_1,V,V_0,V_1,Q$ : \pols en $X$  
\hst \# initialisation 
\hsu   $R_0:=X^{2n}$ ;
$R_1:=\sum_{i=0}^{2n-1}a_iX^i$ ;  $V_0:=0$  ;
$V_1:=1$ ; 
\hst \# boucle 
\hsu    \tantque{n \leq  \deg(R_1)} 
\hsd
$(Q,R) :=$ quotient et reste de la division de $R_0$ par $R_1$ ;
\hsd         $V := V_0-Q*V_1$ ; 
\hsd 
$V_0:=V_1$ ; $V_1:=V$ ; $R_0:=R_1$ ; $R_1:=R$ ; 
\hsu
\fintantque \hst \# sortie \hsu $d:=\sup(\deg(V_1),1+\deg(R_1)) $
; $P:=X^dV_1(1/X)$ ; 
\hsu Retourner $P:=P/{\rm cd}(P)$. \fin
\end{algor}

Cet \algo est d\^{u} \`a Berlekamp, mais sous une
forme o\`u la relation avec l'\algo d'Euclide \'etendu
\'etait invisible. C'est Massey qui a fait le rapprochement.
Pour plus de d\'etails sur la
relation entre cet \algo et l'\algo d'Euclide
\'etendu, on pourra consulter  \cite{Dor}.

\subsection{M\'ethode de Wiedemann}\label{subsec.wied}

L'\algo de Wiedemann \cite{Wie} pour la r\'esolution
des \slis sur un corps \,$\K$\, est un
\algo probabiliste, avec divisions, qui est bas\'e sur
la th\'eorie des \srlsz. Il est particul\`erement
efficace dans le cas des matrices creuses sur les corps finis.
\index{Wiedemann!algorithme de}

Il utilise le fait que si le \polmin \,$\rP^A$\,
d'une matrice \,$A\in \K^{n\times n}$\, est de degr\'e
\,$n$\, (donc \'egal, \`a un signe pr\`es, au \polcar
\,$\rP_A$\, de \,$A\,$)\,, alors il existe toujours un
vecteur \,$v\in \K^{n\times 1}$\, pour lequel le
\polgmin de
la \srl \,$(A^k\,v)_{\,k\in \NN}$\, est \'egal \`a
\,$\rP^A$. Il suffit en effet, comme nous l'avons
vu dans la section \ref{subsec polmin} (corollaire
\ref{corpolmin}),
de prendre un vecteur de \,$\K^{n}$\, en dehors d'une
r\'eunion finie de sous-espaces de \,$\K^{n}$.

L'\algo de Wiedemann choisit au hasard une forme \lin
\,$\pi:\K^{n}\rightarrow \K$\,
et un vecteur \,$v\in\K^{n}$\,  puis il
calcule les \,$2n$\, premiers termes de la \srl
\,$(\pi(A^k\,v))_{\,k\in \NN}$\, dans \,$\K$. Enfin le
 \polgmin de cette suite est obtenu par l'\algo de
Berlekamp-Massey.

Dans le cas \'etudi\'e par Wiedemann, le corps
\,$\K$\, est fini de cardinal
\,$q$, et on a une mesure de probabilit\'e naturelle,
en postulant une \'equiprobabilit\'e des \'el\'ements
du corps.
Si le \polmin de la matrice \,$A$\, est \'egal \`a son
\polcarz,
la probabilit\'e pour trouver un vecteur \,$v$\,
convenable
apr\`es \,$k$\, essais successifs est \supee \`a
$~1-\log{\frac{q^{k-1}}{q^{k-1}-1}}\,\geq
1-\frac{1}{q^{k-1}-1}~$ (\cf \cite{Wie}).

\ss Si on compare avec l'\algo de Frobenius, on voit
que l'on doit calculer \,$2n$\, vecteurs \,$A^k\,v$\, au
lieu de \,$n$. Par contre le calcul du \polgmin est ensuite
beaucoup plus rapide.  En outre, dans le cas des matrices
creuses, m\^{e}me le calcul des  \,$2n$\, vecteurs \,$A^k\,v$\,
est tr\`es rapide.

\ss Notons enfin que les \algos de Frobenius et de
Wiedemann
peuvent \^{e}tre acc\'el\'er\'es tr\`es
significativement au moyen
de la \mul rapide des \pols et de la \mul rapide des
matrices (\cf
sections \ref{sec kellerseq} et \ref{sec kalto}).

\newpage \thispagestyle{empty}
 
\chapter{Circuits arithm\'etiques}\label{chap circuits} 
\minitoc
\acvide
 
\subsubsection*{Introduction}   

Dans ce chapitre nous introduisons la notion fondamentale de \cari
qui est le cadre \gnl dans lequel se situe l'analyse des 
\algos d\'evelopp\'es dans cet ouvrage.

La \coag peut \^etre vue comme une th\'eorie qui 
cherche \`a analyser les \algos qui acceptent de se mettre sous 
forme de \fams de \carisz.

Dans un \cari les instructions de branchement ne sont pas 
autoris\'ees, ce qui semble une limitation assez s\'ev\`ere. 
Les \algos usuels d'\agr \lin sont en effet ordinairement 
\'ecrits en utilisant des tests d'\egt \`a 0.
N\'eanmoins, il s'av\`ere que dans beaucoup de cas, cette limitation 
apparente n'en est pas une, notamment en raison de la \pcd 
d'\elid de Strassen que nous exposons dans la section 
\ref{sec elidiv}.
Par contre le cadre un peu strict fourni par les \caris s'av\`ere 
tr\`es f\'econd. C'est gr\^ace \`a lui que l'on peut mettre en place 
la strat\'egie \gnle \gui{diviser pour gagner}.

Lorsqu'on envisage les \algos li\'es \`a la g\'eom\'etrie 
\agq r\'eelle, la n\'ecessit\'e des tests de signe, et donc 
des instructions de branchement, devient souvent
imp\'erieuse, et une autre 
branche de la  \coag est \ncrz, avec la 
th\'eorie des \emph{r\'eseaux \arithsz} que nous ne d\'evelopperons 
pas ici. 

\ss Dans la section \ref{sec circprev} nous donnons les 
d\'efinitions pr\'ecises
des \caris et de leur variante \gui{programm\'ee}, les \prevs 
(straight-line programs en anglais). C'est l'occasion d'introduire
quelques mesures de \com pour ces \algosz.

\ss Dans la section \ref{sec elidiv} nous introduisons l'\elid 
selon la \met de Strassen et nous 
\'etablissons quelques uns des r\'esultats les plus importants qui la 
concernent.

\ss Dans la section \ref{secpartialder} nous donnons une \met qui 
transforme un \cari \,$\Gamma$\, qui calcule une fraction rationnelle 
\,$f$\, en un \cari \,$\Gamma'$\,  de taille comparable (la taille est 
multipli\'ee par au plus $5$), qui calcule \`a la fois la fonction 
\,$f$\, et toutes ses d\'eriv\'ees partielles.

\section{Circuits arithm\'etiques et \prevs}
\label{sec circprev}  
  
Un \cari constitue une  
fa\c{c}on naturelle et simple de  
repr\'e\-senter les calculs \agqs  
dans un anneau arbitraire, dans le cas o\`u un \algo
n'utilise pas d'instructions de branchements,
et uniquement des boucles du type 

\centerline{\textsf{\pour{i}{m}{n}}\, \ldots\,\,\, \textsf{\finpour}.}

\noindent Si la taille de l'entr\'ee est fix\'ee, ces boucles peuvent
\^{e}tre \gui{mises \`a plat} et on obtient un \prg
dont les seules instructions sont des affectations.

La plupart des \algos pr\'esent\'es dans le chapitre 
\ref{chap BasicAlgoAlin}
sont de ce type et donnent donc lieu, lorsque les dimensions
des matrices sont fix\'ees, \`a des \prevsz. 

\subsection{Quelques d\'efinitions}

Par exemple l'id\'ee d'un \cari est donn\'ee par le calcul du \deter 
d'une \maca par l'\apg  
simplifi\'e, dans le cas des matrices  
\fregsz, pour des matrices de taille fix\'ee.
  Le calcul est alors toujours exactement  
le m\^eme et peut \^etre repr\'e\-sen\-t\'e comme une suite 
d'affectations qu'on peut disposer \sqlt ou 
dessiner au moyen d'un graphe plan. 

Par exemple pour une matrice  
\,$4\times 4$, en donnant un nom diff\'erent \`a  
chaque r\'esultat de calcul \elr  
(addition, soustraction, \mul ou division),  
et en reprenant la notation $a^{[k]}_{ij}$ introduite  
\`a la section \ref{subsec Gauss}, on obtient la mise  
\`a plat sous la forme du \prev \vref{prevLU4} 
dans lequel 
toutes les affectations situ\'ees \`a une m\^eme \prof peuvent en 
principe \^etre ex\'ecut\'ees simultan\'ement.

Pour une \prof donn\'ee, les calculs sont faits avec des  
variables d\'efinies aux \'etages pr\'ec\'edents.  
Ce calcul comprend $37$ \ops arithm\'etiques, et sa  
\emph{\profz} est \'egale \`a $10$, sa \emph{largeur}
est \'egale \`a $9$.

\begin{proevaH}[Calcul du \deter et de la \LU-\deco d'une 
 \maca d'ordre $4$ par la \mpg (sans recherche de 
 pivot).]\label{prevLU4}
\acl{prevLU4}{\Prodeva\, du \deter et de la \LU-\deco d'une 
 \maca d'ordre $4$}
\Entree Une matrice $A=(a_{ij})\in\K^{4\times 4}$ \`a \coes
dans un corps $\K$.
\Sortie Les \coes $l_{ij}$ en dessous de la diagonale de
la matrice \,$L$, les \coes $u_{ij}=a^{[i-1]}_{ij}$ ($j\geq i$)
de la matrice
\,$U$, le \deter $d_4$ de \,$A$.
\Debut 
\Pro{1}{Traitement du premier pivot}
{l_{21}:=a_{21}/a_{11} \;; \;
l_{31}:=a_{31}/a_{11}\;; \;
l_{41}:=a_{41}/a_{11}}
\Pro{2}{largeur 9} 
{b^{[1]}_{22}:=l_{21}\,a_{12}\;;\;	
b^{[1]}_{23}:=l_{21}\,a_{13} \;;\; 
b^{[1]}_{24}:=l_{21}\,a_{14} \;;}
\hsud
{b^{[1]}_{32}:= l_{31}\,a_{12} \;;\; 
b^{[1]}_{33}:=l_{31}\,a_{13} \;;\;
b^{[1]}_{34}:=l_{31}\,a_{14} \;;}
\hsud 
{b^{[1]}_{42}:= l_{41}\,a_{12} \;;\; 
b^{[1]}_{43}:=l_{41}\,a_{13} \;;\;
b^{[1]}_{44}:=l_{31}\,a_{44}  
}
\Pro{3}{largeur 9} 
{a^{[1]}_{22}:=a_{22}-b^{[1]}_{22}\;;\;
a^{[1]}_{23}:=a_{23}-b^{[1]}_{23} \;;\;
a^{[1]}_{24}:=a_{24}-b^{[1]}_{24} \;;  }
\hsud	
{a^{[1]}_{32}:=a_{32}-b^{[1]}_{32} \;;\;
a^{[1]}_{33}:=a_{33}-b^{[1]}_{33}\;;\;
a^{[1]}_{34}:=a_{33}-b^{[1]}_{34}\;; }
\hsud	
{a^{[1]}_{42}:=a_{42}-b^{[1]}_{42} \;;\;
a^{[1]}_{43}:=a_{43}-b^{[1]}_{43}\;;\;
a^{[1]}_{44}:=a_{43}-b^{[1]}_{44}}
\Pro{4}{Traitement du deuxi\`eme pivot} 	
{l_{32}:=a^{[1]}_{32}/a^{[1]}_{22}\;;\;
l_{42}:=a^{[1]}_{42}/a^{[1]}_{22}\;;\;  
d_2:=a_{11}a^{[1]}_{22}}
\Pro{5}{largeur 4} 
{b^{[2]}_{33}:=l_{32}\,a^{[1]}_{23}\;;\;
b^{[2]}_{34}:=l_{32}\,a^{[1]}_{24}\;;  }
\hsud	
{
b^{[2]}_{43}:=l_{42}\,a^{[1]}_{23}\;;\;
b^{[2]}_{44}:=l_{42}\,a^{[1]}_{24}}
\Pro{6}{largeur 4} 
{a^{[2]}_{33}:=a^{[1]}_{33}-b^{[2]}_{33}\;;\;
a^{[2]}_{34}:=a^{[1]}_{34}-b^{[2]}_{34}\;;\;  }
\hsud	
{
a^{[2]}_{43}:=a^{[1]}_{43}-b^{[2]}_{43}\;;\;
a^{[2]}_{44}:=a^{[1]}_{44}-b^{[2]}_{44}
}
\Pro{7}{Traitement du troisi\`eme pivot} 
{l_{43}:=a^{[2]}_{43}/a^{[2]}_{33}\;;\;
l_{44}:=a^{[2]}_{44}/a^{[2]}_{33}\;;\;
d_3:=d_2\,a^{[2]}_{33}}
\Pro{8}{}
{b^{[3]}_{44}:=l_{43}\,a^{[2]}_{44} }
\Pro{9}{}
{a^{[3]}_{44}:=a^{[2]}_{44}-b^{[3]}_{44}}
\Pro{10}{}
{d_4:=d_3\,a^{[3]}_{44}}
\fin
\end{proevaH}

Plus \gnltz.  
  
\begin{defi} \label{defi prevalrit}  
Un \emph{\prev \arith}\index{programme d'\'evaluation!arithm\'etique}  
\,$P$\, sans  
division (\rsp avec division) avec constantes dans  
\,$C$, o\`u \,$C$\, est une partie (cod\'ee) d'un  
anneau \,$\A$\, ou d'un corps \,$\K$, est  
la donn\'ee:\\  
--- d'un ensemble de variables $x_{p,u}$ o\`u $p$  
est un entier $\ge 0$ donnant la \prof de la variable  
et \,$(p,u)$\, est l'identificateur de la variable,\\  
--- d'une suite d'instructions d'affectations de  
l'un des types suivants:  
  
--$~~x_{p,u}:= a \circ b$\, o\`u \,$a$\, et  
\,$b$\, d\'esignent ou bien une variable \,$x_{q,u}$\,  
avec \,$q<p$\, ou bien une constante \,$c\in C$, et  
\,$\circ\in\{+,-,\times \}$\,  
(\rsp $\,\circ\in\{+,-,\times,/ \}$).  
  
--$~~x_{p,u}:= a$\, (avec les m\^emes conventions  
pour \,$a$). \\  
Mises \`a part les variables \,$x_{0,u}$\, qui sont les  
entr\'ees du \prevv toutes les variables  
\,$x_{p,u}$\, sont affect\'ees exactement une fois dans le  
\prgz. Ce sont \emph{les variables d'affectation}
\index{programme d'\'evaluation!variable d'affectation d'un}   
du \prgz.\\  
Les constantes sont consid\'er\'ees comme de \prof nulle.  
En con\-s\'e\-quen\-ce on note \,$\pr(a)=0$\, si \,$a\in C$\,  
et \,$\pr (x_{p,u})=p$.  
\end {defi}  

Quelques commentaires sur cette d\'efinition.

Dans le cas d'un \prev avec divisions, l'\eva
peut \'echouer pour certaines valeurs des variables 
d'entr\'ee dans le corps $\K$. 
Souvent l'ensemble \,$C$\, est vide ou r\'eduit \`a  
\,$\{0,1\}$. Le \prg peut alors \^etre \'evalu\'e  
sur un corps arbitraire (sur un anneau arbitraire s'il  
est sans division).

Naturellement, tous les identificateurs doivent  
\^etre distincts.  
Les affectations du type \,$x_{p,u}:= a\,$
sont pr\'evues uniquement pour le cas o\`u on d\'esire 
respecter certaines contraintes dans une
gestion pr\'ecise des \'etapes \paralsz.

Ordinairement, on demande que dans une affectation  
\,$x_{p,u}:= a \circ b$\, on ait 
\,$\pr(x_{p,u})=1+\max(\pr(a),\pr(b))$\,  
et dans une affectation \,$x_{p,u}:= a$\, on ait  
\,$\pr(x_{p,u})=1+\pr(a)$.  
On peut aussi demander  
que dans une affectation \,$x_{p,u}:= a \circ b$\,  
on ait \\ \,$a\in C$\, ou \,$p=1+\pr(a)$, et  
\,$b\in C$\, ou \,$p=1+\pr(b)$. 

Le texte du \prev doit normalement pr\'eciser 
quelles sont les variables repr\'esentant les sorties. Mais on  
peut demander que les sorties soient exclusivement les  
variables de \prof maximum. On peut demander aussi que  
toute variable de \prof non maximum soit utilis\'ee. 
   
\begin{remark} 
\label{rem defi prevalrit2}  
\emph{De mani\`ere plus \gnlez, un \prev peut \^etre  
d\'efini pour n'importe quel type de structure  
\agqz, une fois qu'ont \'et\'e pr\'ecis\'es les  
op\'erateurs de base dans la structure, qui peuvent  
\^etre de n'importe quelle arit\'e. Par exemple un  
\emph{\prev \bolz}\index{programme d'\'evaluation!bool\'een} 
correspond \`a la structure d'\agr  
de Boole avec les op\'erateurs \bols usuels.  
Autre exemple, dans les anneaux commutatifs, on peut  
d\'efinir une notion de \emph{\prev avec \detersz}  
si on introduit en tant qu'op\'erateurs de base les  
\,$\det_n$\, comme \ops d'arit\'e \,$n^2$\, qui donnent le  
\deter d'une matrice \,$n\times n$\, en fonction  
de ses entr\'ees.  
}  
\end{remark}  

\begin{defi} 
\label{defLongueurs} Nous utiliserons la terminologie suivante  
concernant les \prevsz: 
\begin{itemize}
\item  
Le nombre des entr\'ees dans un \prev est en   
\gnl contr\^ol\'e par un ou plusieurs 
\parats  qu'on appelle 
\emph{les \parats d'entr\'ee}\index{param\`etre 
d'entr\'ee!d'un programme d'\'evaluation}
\index{programme d'\'evaluation!param\`etre d'entr\'ee d'un}
du \prgz. Par exemple, dans un \prev qui calcule le produit de 
deux matrices \,$n\times n$, on prend  
l'entier \,$n$\, comme \parat d'entr\'ee et dans un  
\prev associ\'e \`a la r\'esolution d'un \sys de  
\,$m$\, \'equations \polles \`a \,$n$\,  
\idtrsz, de degr\'e maximum \,$d$, \'ecrites en  
\rpn dense{\footnote{~Un \label{labdense} \pol est \emph{cod\'e en  
\rpn dense} 
\index{repr\'esentation!dense} lorsque le codage donne la liste de tous 
les \coes des mon\^{o}mes en dessous d'un degr\'e donn\'e, dans un ordre 
convenu.
Il est \emph{cod\'e en \rpn creuse} 
\index{repr\'esentation!creuse} lorsque le codage donne la liste des 
paires \,$(a_m,m)$\, o\`u \,$m$\, code un mon\^{o}me (par exemple 
\,$x^2y^3z^5$\, peut \^{e}tre cod\'e par \,$(2,3,5)$\, et \,$a_m$\, son 
\coe (non nul) dans le \polz.}}, les \parats d'entr\'ee
sont \,$m,n,d$.
\item Dans une affectation du type \,$x_{p,u}:=  
a \circ b$, \,$a$\, et \,$b$\, sont les  
\emph{ant\'e\-c\'e\-dents} de \,$x_{p,u}$.
\item  La \emph{profondeur} du \prev est  
\index{programme d'\'evaluation!profondeur d'un}
\index{profondeur!d'un programme d'\'evaluation}
la \prof maximum de ses variables d'affectation;  
not\'ee \,$\pr(P)$, elle corres\-pond au nombre  
d'\'etapes \parals du \prev \,$P$.
\item  La \emph{taille\/}\index{programme d'\'evaluation!taille d'un} 
\index{taille!d'un programme d'\'evaluation}ou 
\emph{longueur\/}\index{programme d'\'evaluation!longueur d'un} 
\index{longueur!d'un programme d'\'evaluation} du \prev  
d\'esignera le nombre total de toutes les \oparisz,
\cad les affectations  
du type \,$x_{p,u}:= a \circ b$.
\item Pour chaque \'etape \,$p~(1\leq p\leq \pr(P))$, on  
consid\`ere le nombre \,$\tau_p$\, d'\ops effectu\'ees  
durant cette \'etape. On appelle \emph{largeur}\index{programme 
d'\'evaluation!largeur d'un} du  
\prev le plus grand de ces nombres , \cad  
\,$\max{\{ \tau_p~|~1\leq p\leq \pr(P)\}}$.
\item Lors de l'\eva d'un  \prev \arith sur un
anneau ou sur un corps dont les \elts sont cod\'es, 
les entr\'ees $x_i$ ont a priori n'importe quelle 
taille tandis que les constantes
du \prg ont une taille fix\'ee une fois pour toutes. 
Du point de vue du calcul concret sur des objets cod\'es,
on est donc souvent en droit d'estimer 
que seules importent vraiment
les \emph{affectations sans scalaires},
\index{affectation sans scalaires} 
\cad celles du type \,$x_{p,u}:= a \pm b$\, et 
\,$x_{p,u}:= a \times  b$\, ou aucun des deux ant\'ec\'edents 
n'est une constante, ainsi que  
\,$x_{p,u}:= a /  b$\, o\`u $b$ n'est pas une constante. Ceci 
donne lieu aux notions de \emph{longueur stricte} 
\index{longueur!stricte} et de
\emph{\prof stricte}, \index{profondeur!stricte} dans lesquelles seules 
sont prises en 
compte les affectations sans scalaires. Une \mul ou division 
sans scalaire dans un \prev \arith
est encore dite \emph{essentielle}.\index{essentielle!multiplication}
\item  \label{profmul} Variation sur le th\`eme pr\'ec\'edent. 
Dans la mesure o\`u on consi\`ere que les additions ainsi que 
les 
\muls ou divisions par des constantes sont 
relativement peu co\^uteuses (ou \'eventuellement 
pour des raisons plus profondes d'ordre th\'eorique)
on est int\'eress\'e par la 
\emph{longueur \muvz}\index{longueur!multiplicative} 
d'un \prev et par sa 
\emph{\prof \muvz}\index{profondeur!multiplicative}
qui sont d\'efinies comme la longueur et la \prof 
mais en ne tenant compte que des 
\muls et divisions essentielles. 
\end{itemize}
\end{defi}

Par exemple le \prev \ref{prevLU4} a une \prof \muv \'egale \`a $6$
et une largeur \muv \'egale \`a $9$. 

\subsection{Circuit \arith vu comme un graphe}   

On peut \egmt repr\'esenter un \prev sous forme d'un dessin plan.  
Par exemple pour le calcul du \deter par l'\apg avec une matrice  
\freg \,$3\times 3$, on peut le repr\'esenter par le dessin  
du circuit \vref{F6}.  
  
\ms  
\setcounter{agh}{\arabic{algorithm}}
\stepcounter{algorithm}
\begin{agh}  
\begin{center}
\includegraphics*[width=11cm]{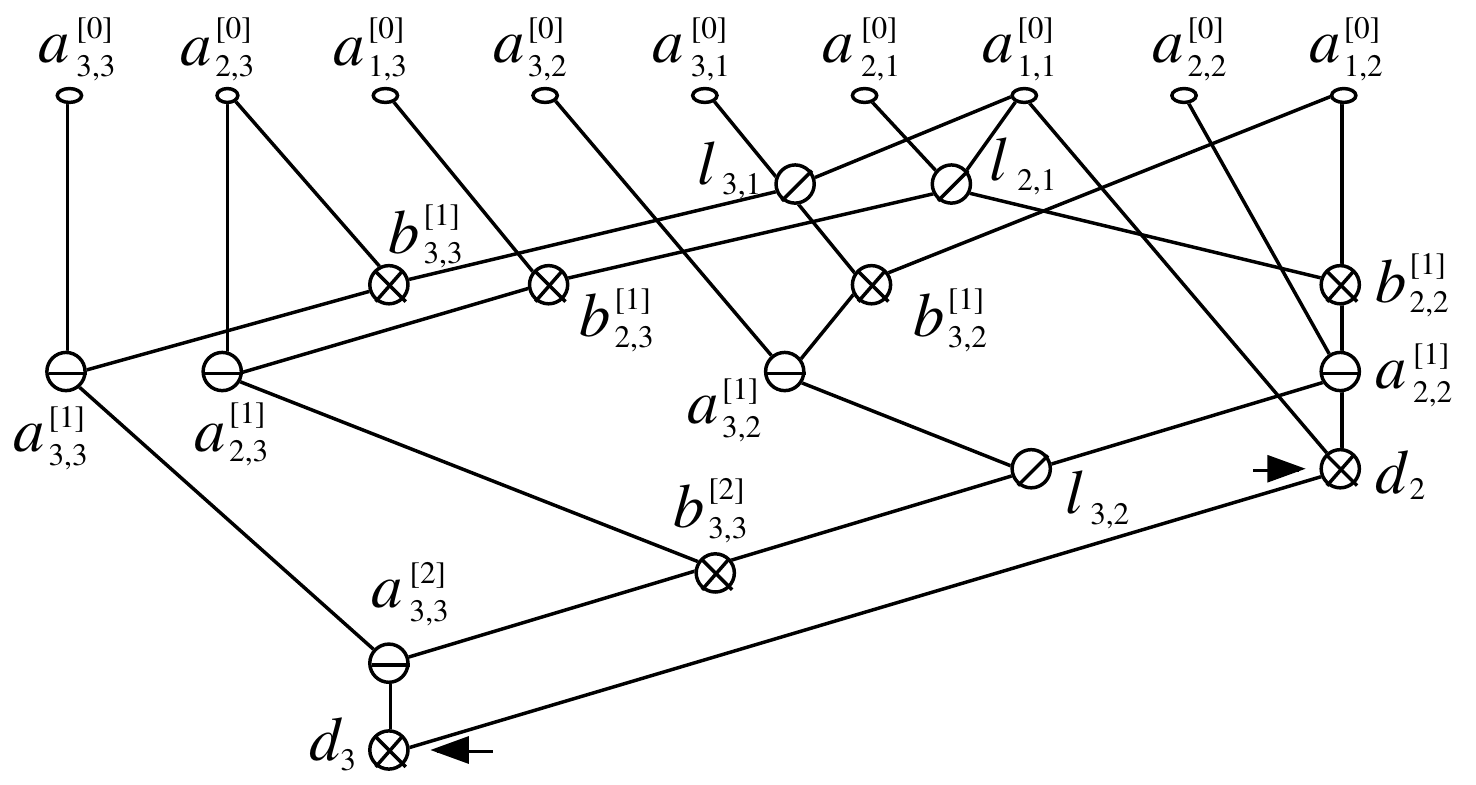}
\end{center}
\caption[Circuit de l'\apg simplifi\'e]{\label{F6} 
Pivot de Gauss  
simplifi\'e pour une matrice \,$3\times 3\,$}  

\centerline{($\oslash$: division; \ \ $\otimes$:  
\mulz; \ \ $\ominus$: soustraction)}  
  
\centerline{{\footnotesize Pour la division et la  
soustraction}}
 
\centerline{{\footnotesize le brin entrant gauche 
repr\'esente le premier terme de l'op\'eration}}  
\end{agh} 
  
\noindent
Pour une matrice \,$n\times n$, on obtiendra un  
\cir de \prof $3n-2$ avec un nombre de  
portes, en tenant compte des \,$n-1$\, affectations  
\,$d_p:=d_{p-1}\,a^{[p-1]}_{~pp}~~(2\leq p\leq n)~$  
qui donnent les \mips dominants de la  
matrice, \'egal \`a:  
$$
(n-1) + \som_{k=1}^{n-1}k\,(2k+1) =  
{1\over 6}\,(n-1)\,(4n^2+n+6)\,.$$

Si on veut formaliser ce genre de dessin qui visualise  
bien la situation, on peut adopter la d\'efinition  
 suivante.  

\begin {defi}  
Un \emph{\cari avec divisions}\index{circuit!arithm\'etique}
 \index{circuit!arithm\'etique!avec divisions}
\index{circuit!arithm\'etique!sans division} (\rsp sans 
division) \gui{avec  
constantes dans \,$C$}, o\`u \,$C$\, est une partie  
(cod\'ee) d'un anneau \,$\A$\, ou d'un corps \,$\K$,  
est un graphe acyclique orient\'e et \'etiquet\'e, 
chaque n{\oe}ud qui n'est  
pas une porte d'entr\'ee ayant exactement deux  
ant\'ec\'edents. Le \cir est \'etiquet\'e  
de la mani\`ere suivante:\\  
--- chaque porte d'entr\'ee est \'etiquet\'ee par un  
triplet \,$(0,n,c)$\, o\`u \,$0$\, est la \profz,  
\,$(0,n)$\, est le nom qui identifie le n{\oe}ud et  
\,$c\in \{x\}\cup C$\, (avec \,$x\notin C$).  
Ici un triplet \,$(0,n,x)$\, repr\'esente la variable  
\,$x_n$\, et un triplet \,$(0,n,c)$\, avec \,$c\in C$\,  
repr\'esente l'\elt (cod\'e par) \,$c$\, de  
\,$\A$.\\  
--- chaque n{\oe}ud interne et  chaque porte de sortie est  
\'etiquet\'e par un triplet \,$(m,n,\circ)$\, o\`u  
\,$m$\, est sa \profz, \,$(m,n)$\, est son  
identificateur, et \,$\circ\in\{+,-,\times,/ \}$\,  
(\rsp \,$\circ\in\{+,-,\times\}$)  
d\'esigne une \opariz.\\  
--- enfin, dans le cas des op\'erateurs $/$ et $-$ (et  
dans le cas de l'op\'erateur $\times$ si le \cir est  
destin\'e \`a \^etre \'evalu\'e dans un anneau non  
commutatif) il faut \'etiqueter de mani\`ere \`a les  
distinguer (gauche, droite) les deux arcs qui aboutissent  
\`a un n{\oe}ud correspondant \`a cet op\'erateur. \\  
--- les portes de sortie correspondant aux r\'esultats du  
calcul sont sp\'ecifi\'ees par une marque distinctive dans  
leur identificateur.  
\end{defi}  
  
En fait, dans toute la suite, nous utiliserons  
indiff\'e\-rem\-ment \gui{\cariz} et  
\gui{\prev \arithz}, tout en sous-entendant que,  
pour ce qui est du codage, nous choisissons toujours  
un codage correspondant \`a la d\'efinition d'un  
\prev \arithz.  De la m\^eme mani\`ere, nous consid\'ererons comme 
synonymes \prev \bol et \ciboz\index{circuit!bool\'een}.  
  
Dans un \cari on peut interpr\'eter chaque n{\oe}ud com\-me  
repr\'esentant un \poly de \,$\A\,[(x_i)_{i\in I}]$\, ou  
$\K\,[(x_i)_{i\in I}]$ (dans le cas sans division) ou une  
\frat de $\K((x_i)_{i\in I})$ (dans le cas avec division).  

\subsection{Circuits \ariths homog\`enes} 

\begin{definition} 
\label{defCariH} 
On appelle \emph{\carihz} \index{circuit!arithm\'etique!homog\`ene} un 
\cari sans division dont tous les noeuds
repr\'esentent des \pols \hogs et qui a la structure suivante.
\begin{itemize}
\item Les \pols de degr\'e \,$d$\, sont calcul\'es apr\`es ceux de 
degr\'es strictement inf\'erieurs.
\item Le calcul des \pols de degr\'e \,$d$\, se fait en deux
phases. Dans la premi\`ere phase, en une seule \'etape \paralz,
on effectue des produits de \pols
pr\'ec\'edemment calcul\'es (de degr\'es \,$d'<d$\, et \,$d-d'$). Dans 
la deuxi\`eme phase on calcule des
\colis des pr\'ec\'edents.
\end{itemize}
\end{definition}

\begin{prop}  
\label{prop Eli Div4}  
Tout \cari sans division qui calcule une  \fam de \polys de 
degr\'e $\leq d$  peut \^etre r\'eorganis\'e en un \carih
qui calcule toutes les composantes \hogs des \pols en sortie.
Le \cari \hog obtenu est de \prof \muv \,$d-1$. 
Par rapport au \cir initial, la \prof a \'et\'e  multipli\'ee
par \,$\O(\log\,d)$, 
la longueur \muv  a \'et\'e au plus multipli\'ee par 
\,$d(d-1)/2$, et la longueur totale a \'et\'e au plus multipli\'ee
par \,$(d+1)^2$.   
\end{prop}  
\preuve  
Chaque noeud \,$y_j$\, du \cir initial repr\'esente un \pol en les 
entr\'ees
\,$x_i$\, qu'on d\'ecompose en somme de composantes \hogsz. 
$$ y_j:=y_j^{[0]}+y_j^{[1]}+\cdots+y_j^{[d]}+ \mathrm{des\; 
composantes\; 
sans\; importance}
$$
On analyse alors le calcul qui est fait sur les composantes
 \hogs de degr\'e $\leq d$.\\ 
Lorsqu'on a dans le \cari original une affectation correspondant 
\`a une addition $y_\ell:=y_h+y_k$  on obtient sur les
composantes \hogs au plus \,$d+1$\, additions qui peuvent \^{e}tre
ex\'ecut\'ees en \,$\esup{\log(d+1)}$\, \'etapes \paralsz.\\ 
Lorsqu'on a dans le \cari original une affectation correspondant 
\`a une \mul essentielle $y_\ell:=y_h\,y_k$  on obtient 
$$\begin{array}{rclll} 
y_\ell^{[0]}& = & y_h^{[0]}\, y_k^{[0]}    \\ 
y_\ell^{[1]}& = & y_h^{[0]}\, y_k^{[1]}+ y_h^{[1]}\, y_k^{[0]}   \\ 
y_\ell^{[2]}& = & y_h^{[0]}\, y_k^{[2]}+ y_h^{[1]}\, y_k^{[1]} + 
y_h^{[2]}\, y_k^{[0]}    \\ 
\vdots~~& \vdots & ~~~~\vdots    \\ 
y_\ell^{[d]}& = & y_h^{[0]}\, y_k^{[d]} + y_h^{[d-1]}\, y_k^{[1]} 
+ \cdots + y_h^{[d]}\, y_k^{[0]}    
\end{array}$$
ce qui correspond \`a (au plus)  \,$d(d-1)/2$\, \muls essentielles 
entre les composantes \hogsz,  \,$2d+1$\, \muls scalaires, et 
\,$d(d+1)/2$\, additions, soit \,$(d+1)^2$\, \oparis en tout.\\
 Par ailleurs on peut r\'eorganiser 
l'ensemble du calcul de mani\`ere que tous les \pols \hogs de degr\'e 
$k$ soient
calcul\'es apr\`es ceux de degr\'e $ < k$.
\qed

\subsection[Le \pb des divisions]{Le \pb des divisions dans les \caris} 

Certains \caris avec division comportent une division par  
une \frat identiquement nulle, et ne repr\'esentent plus  
aucun calcul raisonnable. Implicitement on suppose  
toujours qu'on n'est pas dans ce cas.  
  
\ss  
Le cas de l'\ajb sans recherche de pivot est un peu plus subtil. Il  
correspond \`a un \cari \gui{avec divisions exactes},  
\cad que, lorsqu'on le regarde comme produisant \`a  
chaque porte un \elt du corps des fractions  
rationnelles, on reste en fait toujours dans l'anneau  
des \polsz: les divisions ont toujours pour  
r\'esultat un \poly et non une fraction rationnelle.  
  
\ss  
Si les portes de sortie d'un \cari sont des \pols  
(en les entr\'ees) il est a priori pr\'ef\'erable que  
le \cir soit sans division. Il pourra en effet \^etre  
\'evalu\'e dans n'importe quel anneau.  
  
Dans le cas d'un \cir avec divisions \'evalu\'e dans  
un corps, il se peut que certaines divisions soient  
impossibles, non parce qu'on doit diviser par une \frat  
identiquement nulle, mais parce que les valeurs des \,$x_i$\,  
annulent la \frat du \denoz.  
C'est encore une raison qui fait qu'on pr\'ef\`ere les  
\cirs sans divisions.  
  
Une autre raison est que si le corps \,$\K$\, est de  
\cara nulle ou s'il contient des \elts  
transcendants, l'addition de deux fractions est une  
affaire bien encombrante. 
L'addition dans  
\,$\QQ$\, par exemple, r\'eclame, dans \,$\ZZ$\, 
tout d'abord 3 \muls et une addition, suivies d'une simplification
de fraction, qui r\'eclame un calcul de pgcd, donc les divisions 
successives de l'\algo d'Euclide. 
Ainsi, lorsque les entr\'ees sont dans \,$\ZZ$\, par  
exemple, on pr\'ef\`ere que tout le calcul reste dans \,$\ZZ$.  
  
\ss Si on essaie d'\'evaluer un \cir avec divisions dans  
un anneau arbitraire \,$\A$, la situation est encore un peu  
compliqu\'ee. Toute division par un diviseur de z\'ero est  
impossible. Et si on divise par un non diviseur de z\'ero,  
on se retrouve naturellement dans l'anneau total des  
fractions de \,$\A$, d\'efini de la m\^{e}me mani\`ere que
le corps des fractions  
d'un anneau int\`egre, mais en autorisant comme  
\denos uniquement des non diviseurs de z\'ero dans  
\,$\A$. Naturellement, les calculs dans ce nouvel anneau  
\,$\A'$\, sont nettement plus compliqu\'es que ceux dans \,$\A$\,  
(cf. la discussion \`a propos de la \met du pivot  
dans $\QQ$.)  
  
\ss  
La \prof d'un \cir est un \parat pertinent \`a  
plus d'un titre. Tout d'abord, la \prof repr\'esente en  
quelque sorte le \gui{temps de calcul \paralz} si on  
donne une unit\'e de temps pour chaque \op \arith  
et si on dispose de suffisamment de \gui{processeurs} entre  
lesquels on r\'epartit les calculs \`a faire.  
Ensuite, la \prof permet un contr\^ole de la taille des  
objets \itmds lorsque le calcul est effectu\'e  
comme une \eva par exemple dans \,$\ZZ$, \,$\QQ$, ou  
\,$\QQ(x,y)$. Grosso modo, la taille double au maximum  
lorsque la \prof augmente de 1. Dans le cas des \cirs  
sans divisions \'evalu\'es par exemple dans \,$\ZZ$\, ou  
\,$\ZZ[x,y]$\, la \prof multiplicative est de loin la plus  
importante pour le contr\^ole de la taille des objets  
\itmdsz.  
  
\ss  
Tout ceci a conduit \`a attacher une importance toute  
particuli\`ere aux \cirs sans division et de faible  
\profz.

\section{\'Elimination des divisions \`a la Strassen}  
\label{subsec elidivStrassen} \label{sec elidiv} 

Lorqu'on dispose d'une \pcd utilisant les divisions  
dans le corps des fractions rationnelles pour calculer un \poly  de  
degr\'e d\'etermin\'e \`a \coes dans un anneau  
int\`egre, une technique de Strassen (\cite{Stra}) 
bas\'ee sur une id\'ee  
tr\`es simple permet d'\'eliminer toutes les divisions dans  
cette \pcdz.  

\subsection{Le principe g\'en\'eral}    
\label{subsecelidivPG}
L'id\'ee de base est que la division par un \poly de la  
forme \,$(1-u)$\, o\`u \,$u = u((x_i))$\, peut \^etre  
remplac\'ee par le produit par la s\'erie formelle  
$$1+u+\cdots+u^m+\cdots$$  
\`a condition d'\^etre dans une situation o\`u on sait  
qu'on peut ne consid\'erer qu'une partie finie bien  
contr\^ol\'ee des s\'eries formelles en jeu.  
  
Nous allons expliquer cette id\'ee fondamentale sur  
l'exemple du calcul du \deter d'une matrice  
carr\'ee \,$A$\, par l'\apg simplifi\'e (\cad sans recherche  
de pivot) mis sous forme d'un \cari pour les matrices  
\,$n\times n$\, pour une valeur fix\'ee de \,$n$.  
On consid\`ere ce \cir comme un \prev dans un anneau  
arbitraire $\A$ (on peut se limiter au sous-anneau  
engendr\'e par les \coes de la matrice, ou plut\^ot  
\`a l'anneau total des fractions de ce sous-anneau).  
Naturellement un obstacle appara\^{\i}t \'eventuellement  
lors d'une affectation \,$v_h:=v_k/v_{\ell}$\, si  
$v_{\ell}$ est diviseur de z\'ero.  
Il y a cependant des cas o\`u cet obstacle  
n'appara\^{\i}t pas du tout, le plus simple est celui o\`u  
la matrice de d\'epart est \'egale \`a $\In$: toutes les  
divisions se font par $1$~!  
Cette remarque d'apparence anodine est cependant la cl\'e  
de l'\elidz. En effet, il suffit de  
faire le changement de variable \,$F:=A-\I_n$\, et de  
d\'ecider d'\'evaluer le \cir pour l'entr\'ee  
\,$\I_n+zF$, o\`u $z$ est une nouvelle variable, dans  
l'anneau  $\aqo{\A[z]}{z^{n+1}}$:
\emph{l'anneau des \dlis \`a l'ordre \,$n$\, \`a  
\coes dans $\A$}, que nous noterons souvent  \,$\A_n$.
\index{anneau des  \dlis \`a l'ordre \,$n$}
  
Quelle que soit la matrice \,$F$\, \`a \coes dans  
\,$\A$\, prise en entr\'ee, chaque n{\oe}ud \,$v_{\ell}$\,  
intervenant dans une division est maintenant un  
\dli du type  
$$v_{\ell}=1+c_{\ell, 1}z+\cdots+c_{\ell, n}z^n,\qquad  
(c_{\ell, 1},\ldots,c_{\ell, n}\in\A )$$  
\cad un \elt inversible de \,$\A_n=\aqo{\A[z]}{z^{n+1}}$.  
A la fin du calcul on r\'ecup\`ere donc  
\,$\det{(\I_n+zF)}$\, dans \,$\A_n$, \cad en  
fait: on r\'ecup\`ere \,$\det(\I_n+zF)$\, dans  
$\A[z]$({\footnote{~Dans le cas pr\'esent, il serait  
donc plus astucieux, d'appliquer la \pcd avec la  
matrice \,$A$\, \`a la place de la matrice \,$F$, car on  
obtient ici \`a tr\`es peu pr\`es le \polcar de  
\,$F$. Dans le cas pr\'esent la \pcd d'\elid est donc tr\`es 
proche de la \mjb modifi\'ee. Cette derni\`ere est cependant un peu 
plus simple, car dans la \met de Strassen on manipule tr\`es 
rapidement des \pols de degr\'e \,$n$\, 
(l'ordre de la matrice). Pour terminer notons que c'est  
un fait d'exp\'erience assez curieux que les \pcds  
\gui{rapides} de calcul sans division du \deter passent toutes  
par le calcul du \polcarz.}}).  
Et il suffit de faire \,$z=1$\, pour obtenir  
\,$\det(A)$. 
En fait la division dans \,$\A_n$\, d'un \elt \,$a(z)$\, par un  
\elt inversible $b(z)=1-zu(z)$ ne n\'ecessite que des  
additions et \muls dans \,$\A_n$: on peut en effet
faire une division en puissances croissantes de \,$a(z)$\, par
\,$b(z)$\, jusqu'\`a l'ordre \,$n$\, en \,$z$. On peut
\'egalement invoquer la formule (valable dans \,$\A_n$)
\begin{equation} \label{eqInvDevLim}
(1-zu(z))^{-1}=(1-w)^{-1}=(1+w)(1+w^2)(1+w^4)  
\cdots(1+w^{2^k})
\end{equation}
si \,$2^{k+1}\geq n+1$\, (il suffit de prendre \,$k=\esup{\log{(n+1)}} 
-1$).  

\ss Ainsi toutes les  
affectations (correspondant \`a l'\eva du  
\cirz) dans \,$\A_n$\, se ram\`enent \`a des  
additions et \muls de \pols tronqu\'es, \cad  
encore \`a des additions et \muls dans $\A$.  
Le \tho suivant est maintenant clair:  
\begin{theorem}  
\label{thEliDiv}  
La \pcd d'\elid  de Strassen
\index{elimination@\'elimination des divisions (Strassen)}  
peut \^etre appliqu\'ee  
\`a tout \cari pourvu qu'on soit dans le cas suivant: on  
conna\^{\i}t un point \,$(\xi_1,\ldots,\xi_n)$\,  de  
\gui{l'espace des entr\'ees} tel que, lorsque le \cir est  
\'evalu\'e en ce point, toutes les divisions qui doivent  
\^etre ex\'ecut\'ees le sont par des \elts  
inversibles de l'anneau de base (on rajoute alors ces  
\elts et leurs inverses \`a l'ensemble des  
constantes $C$ du \cirz).\\
En particulier l'\elid est toujours possible si l'anneau de 
base est un corps infini.  
\end{theorem}  
\begin{defi}  
\label{def Eli Div}  
\emph{\'Eliminer les divisions (\`a la Strassen)
dans un \cari \`a partir du point  
\,$(\xi_1,$ $\ldots,\- \xi_n)$},
\index{eliminer@\'eliminer les divisions!\`a la Strassen dans un \cariz}  
 c'est lui appliquer la  
\pcd d'\elid de Strassen en utilisant  
\,$(\xi_1,\ldots,\xi_n)$\, comme point en lequel le \cir  
est \'evalu\'e sans divisions. Nous appellerons ce point le
\emph{centre d'\elidz}\index{centre d'\'elimination 
des divisions}. 
\end {defi}  
Sur un corps infini, l'existence d'un centre d'\elid 
pour un \cari r\'esulte du fait qu'on peut toujours \'eviter 
l'ensemble 
des z\'eros d'une \fam finie de \pols non (formellement) 
nuls: 
leur produit est un \poly non formellement nul et un tel \poly 
d\'efinit une fonction non identiquement nulle sur $\K^n$ ($n$\, 
est le nombre de variables) si $\K$ est infini.

\subsubsection*{Un exemple d'\elid}   

Donnons \`a titre d'exemple 
le r\'esultat de l'\elid pour l'\apg simplifi\'e,  
dans le cas \,$n=3$,  
pour une matrice \,$\I_3+zF$.
Le  circuit initial
 est donn\'e par le \prev \vref{ancien}.

\begin{proevaH}[Calcul du \deter  de la 
 \maca \,$\I_3+zF$\, par la \mpgz.]
 \label{ancien} 
\acl{ancien}{\Prodeva\, du \deter et de la \LU-\deco de la 
 \maca $\,\I_3+zF$}
 \Entree Les \coes \,$f_{ij}$\, de la matrice \,$F$\, dans un \acoma 
\,$\A$.
\Sortie Le \deter \,$d_{3}=\det(\I_3+zF)$. Le calcul est correct 
si on se situe dans un anneau \,$\B$\,
contenant  \,$z$\, et 
\,$\A$\, et dans lequel tous les \elts de la forme $1+zb$ 
sont inversibles.  Les \oparis de ce \prev sont effectu\'ees 
dans \,$\B$. Notez que les \coes de \,$\I_3+zF$\, sont les
\elts \,$zf_{ij}$\, pour \,$i\neq j$\, et les
\elts \,$(1+zf_{ii})$\, pour \,$i=j$.
\Debut 
\Pro{1}{Traitement du premier pivot}
{l_{21}:=z\,f_{\,21}/(1+zf_{\,11})~; \;  
 l_{31}:=z\,f_{\,31}/(1+zf_{\,11})}
\Pro{2}{} 
{b^{[1]}_{\,22}:=l_{21}\,z\,f_{\,12}\;;\;	
b^{[1]}_{\,23}:=l_{21}\,z\,f_{\,13} \;;\;}
\hsud
{b^{[1]}_{\,32}:=l_{31}\,z\,f_{\,12} \;;\; 
b^{[1]}_{\,33}:=l_{31}\,z\,f_{\,13}
}
\Pro{3}{} 
{f^{[1]}_{\,22}:=z\,f_{\,22}-b^{[1]}_{\,22}\;;\;
f^{[1]}_{\,23}:=z\,f_{\,23}-b^{[1]}_{\,23} \;;\;}
\hsud
{f^{[1]}_{\,32}:=z\,f_{\,32}-b^{[1]}_{\,32}\;;\;
f^{[1]}_{\,33}:=z\,f_{\,33}-b^{[1]}_{\,33}
}
\Pro{4}{Traitement du deuxi\`eme pivot}{
l_{32}:=f^{[1]}_{\,32}/(1+f^{[1]}_{\,22})\;;\;
d_2:=(1+z\,f_{\,11})\,(1+f^{[1]}_{\,22})}
\Pro{5}{}{b^{[2]}_{\,33}:=l_{32}\,f^{[1]}_{\,23}}
\Pro{6}{}{f^{[2]}_{\,33}:=f^{[1]}_{\,33}-b^{[2]}_{\,33}}
\Pro{7}{}{d_3:=d_2\,(1+f^{[2]}_{\,33})}
\fin
\end{proevaH}

Pour passer de l'ancien \cir (\prev \ref{ancien})  
au nouveau (\prev \vref{nouveau}), chaque porte  
\,$y_{ij}$\, ou \,$d_i$\, (sauf les portes d'entr\'ee)  
a \'et\'e remplac\'ee par les portes \,$y_{ijk}$\, ou  
\,$d_{ik}$, avec \,$k=0,\,\ldots,\,3$, qui donnent  
les quatre premiers \coes de la s\'erie formelle  
en \,$z$.

Dans l'\algo transform\'e \ref{nouveau}, 
nous n'avons pas \'ecrit les portes nulles  
(pour les bas degr\'es) et nous n'avons pas mentionn\'e les  
\,$y_{ijk}$\, qu'il est inutile d'\'evaluer pour obtenir le r\'esultat 
final.

Il faut remarquer que pour le \deter et  
m\^eme le \polcar des matrices  
\,$3\times 3$, les formules directes sont  
bien entendu pr\'ef\'erables.

\begin{proevaH}[Calcul du \deter  d'une 
 \maca \,$F$\, par la \mpg apr\`es \elid \`a la Strassen.]
 \label{nouveau} 
\acl{nouveau}{\Prodeva\, du \deter d'une \maca apr\`es \elid \`a la 
Strassen}
\Entree Les \coes \,$f_{ij}$\, de la matrice \,$F$\, dans un \acomaz.
\Sortie Le \deter \,$\det(F)$. En fait, on calcule
m\^{e}me
\,$d_{3}=\det(\I_3+zF)=1+d_{31}z+d_{32}z^2+d_{33}z^3$. 
L'\algo fonctionne \gui{en ligne droite}
et sans aucune hypoth\`ese restrictive. 
Les \oparis de ce \prev sont effectu\'ees 
dans \,$\A$.
\Debut \\
\textbf{Renommages :} $l_{211}=f_{21}$, $l_{311}=f_{31}$,
$f^{[1]}_{221}=f_{22}$, $f^{[1]}_{231}=f_{23}$, 
\hsu $l_{321}=f^{[1]}_{321}=f_{32}$, $f^{[1]}_{331}=f_{33}$.
\Pro{1}{Traitement du premier pivot}  
{l_{212}:=-f_{21}\, f_{11} \;;\;   
l_{312}:=-f_{31}\, f_{11}  \;;\; d_{21}:=f_{11}+f_{22} }
\Pro{2}{} 
{  
f^{[1]}_{222}:=-l_{211}\,f_{12} \;;\;   
f^{[1]}_{232}:=-l_{211}\,f_{13} \;;\;  
 }
\hsud
{  
f^{[1]}_{322}:=-l_{311}\,f_{12} \;;\;  
f^{[1]}_{332}:=-l_{311}\,f_{13} \;;\;  
f^{[1]}_{333}:=-l_{312}\,f_{13} }
\Pro{3}{} 
{d_{22}:=f_{11}\,f^{[1]}_{221}+f^{[1]}_{222}    }
\Pro{4}{Traitement du deuxi\`eme pivot}{   
l_{322}:=f^{[1]}_{322}-f^{[1]}_{221}\,f^{[1]}_{321}    
}
\Pro{6}{}{  
b^{[2]}_{332}:=l_{321}\,f^{[1]}_{231}\;;\;  
b^{[2]}_{333}:=l_{321}\,f^{[1]}_{232}+  
l_{322}\,f^{[1]}_{231}}
\Pro{7}{}{  
f^{[2]}_{331}:=f^{[1]}_{331} \;;\;  
f^{[2]}_{332}:=f^{[1]}_{332}-b^{[2]}_{332} \;;\;  
f^{[2]}_{333}:=f^{[1]}_{333}-b^{[2]}_{333}}
\Pro{8}{}{
 d_{31}:=d_{21}+f^{[2]}_{331} }
\Pro{9}{}{  
d_{32}:=d_{22}+d_{21}\,f^{[2]}_{331}+f^{[2]}_{332}  \;;\; 
d_{33}:=d_{22}\,f^{[2]}_{331}+d_{21}\,f^{[2]}_{332} +  
f^{[2]}_{333}}  
\fin
\end{proevaH}

\newpage
  
\subsection{Co\^ut de l'\elid}

Quel est le co\^{u}t de la transformation d'un \cir  
avec division en un \cir sans division, lorsque les  
sorties sont des \pols de degr\'e $\le n$ en les  
entr\'ees~?  
  
\ss  
Tout d'abord si on utilise les \algos \usls pour les \oparis   
dans \,$\A_n$, la taille du \cir sera en  
gros multipli\'ee par $n^2$ (ce qui fait qu'on reste dans  
le cadre des \cirs de taille \pollez). 
Par exemple le produit de deux \elts de \,$\A_n$\, r\'eclame
\,$(n+1)(n+2)/2$\, \muls et \,$n\,(n+1)/2$\, additions,
tandis que la division de \,$a(z)$\, par \,$1-zu(z)$\, 
n\'ecessite \,$n\,(n+1)/2$\, \muls et autant d'additions si
on effectue la division en puissance croissante.

Si on applique ces constatations 
dans le cas du calcul du \deter
(et du \polcar par la m\^{e}me occasion) d'une \maca 
par \elid dans l'\apg comme nous l'avons vu \`a la section
\ref{subsecelidivPG},
on trouve une taille de \cir \'equivalente \`a 
\,$\sum_{k=1}^{n-1}\,n^2\,(n-k)^2$\, 
\cad  un \,$\frac{1}{3}\,n^5+\O(n^4)$, \`a comparer au
 \,$\frac{1}{10}\,n^5+\O(n^4)$\, que nous avons obtenu pour
l'\ajb modifi\'e.

Notons aussi que la \mul dans \,$\A_n$\, par l'\algo \usl se fait 
naturellement en \prof \,$\O(\log n)$\, tandis que la division par 
puissance croissante est en \prof \,$\O(n\,\log n)$. 
On peut pallier ce dernier inconv\'enient en utilisant 
la \form (\ref{eqInvDevLim}) qui donne un \cir de taille  
\,$\O(n^2\,\log n)$\, et de \prof
\,$\O(\log^2 n)$.

 Il existe par ailleurs des \pcds de \mul rapide  
pour les \polsz: les \oparis  $+,-,\times$ et  
division par un \elt \,$f$\, v\'eri\-fiant  
\,$f(0)=1$\, dans \,$\A_n$\, peuvent \^etre  
ex\'ecut\'ees par des \cirs de taille  
\,$\O(n\,\log{n}\,\log{\log{n}})$\, et de \prof  
\,$\O(\log{n})$\,  (voir \cite{Cant} et, infra, le \tho \vref{thCaK}).   

Plus \gnlt nous utiliserons la notation suivante.
\begin{notation} 
\label{notamuP}
{\rm  Pour un anneau \,$\A$\, fix\'e par le contexte,
nous noterons \,$\mu_P(n)$\, \indexnota{muP@$\mu_P(n)$} le nombre 
d'\oparis \ncrs pour la \mul de deux \pols de degr\'e \,$n$\, en \prof 
\,$\O(\log{n})$.
} 
\end{notation}
 
Strassen obtient alors pr\'ecis\'ement le r\'esultat suivant: 

\begin{theorem}  
\label{prop Eli Div2}  
Lorsqu'on \'elimine les divisions \`a la
Strassen pour l'\eva d'une \fam de \polys de degr\'es 
\,$\leq n$\,  la \prof du
\cir est  multipli\'ee par  
$\O(\log{n})$ et sa taille par  
\,$\O(\mu_P(n))$.
\end{theorem}  

Notons aussi le r\'esultat suivant simple et int\'eressant 
concernant les \caris
qui \'evaluent des \fams de \pols du second degr\'e.

\begin{prop}  
\label{prop Eli Div3}  
Lorsqu'on \'elimine les divisions \`a la
Strassen pour l'\eva d'une \fam de \polys de degr\'e 
$\leq 2$, la longueur \muv du
\cari est  inchang\'ee. 
\end{prop}  

\preuve 
Lorsqu'on applique la \pcd d'\elidz, supposons qu'on ait 
\,$f=f_0+z f_1+z^2 f_2$,
et  \,$g=g_0+z g_1+z^2 g_2$\, dans 
l'anneau des 
\dlis \`a l'ordre $2$ en \,$z$\, sur 
\,$\A\,[(x_i)]$\, (ici on suppose \spdg que $(0,\ldots,0)$ est le centre 
d'\elid et donc que les \,$f_j$\, et \,$g_j$\, sont \hogs de degr\'e 
\,$j$\, 
en les entr\'ees \,$x_i$).
On obtient pour  le produit  \,$h=fg$\, modulo \,$z^3$, 
\,$h=h_0+z h_1+z^2 h_2$\, avec \,$h _0=f_0g_0$, 
\,$h_1=f _1g_0+f_0g _1$\, et 
\,$h_2=f _2g_0+f_1g _1+f _0g _2$\,
avec la seule \mul essentielle \,$f _1g _1$\, puisque
\,$f _0$\, et \,$g _0$\, sont des constantes. 
Et on a un calcul analogue pour \,$k=f/g$.
\,$k_0=f _0/g_0$, 
\,$k _1 = f_1/g_0 - g _1 (f_0/g_0^2)$\, et 
\,$k _2 = f_2/g_0 - k_1g _1/g_0 - g_2k_0/g_0$\, avec la seule \mul 
essentielle \,$k_1g _1$.
\qed  

\ms On pourrait \gnr  avec un \cari calculant une 
\fam de \pols de degr\'es $\leq d$.

\section[Calcul des d\'eriv\'ees partielles]{Calcul de toutes les 
d\'eriv\'ees partielles d'un \pol ou d'une fraction rationnelle} 
\label{secpartialder}

Nous donnons une \met pour transformer un \cari \,$\Gamma$\, qui calcule 
une fraction rationnelle \,$f$\, en un \cari \,$\Gamma'$\,  de taille 
comparable (la taille est multipli\'ee par au plus $5$), qui calcule \`a 
la fois la fonction \,$f$\, et toutes ses d\'eriv\'ees partielles.

Si le \cari \,$\Gamma$\, est sans division, il en est de m\^{e}me pour 
\,$\Gamma'$. La \met est due \`a Baur~\&~Strassen \cite{Baur}. Nous 
suivons l'expos\'e  simple et constructif que Morgenstern en fait  
dans \cite{Mor}.  

Une application importante de ce r\'esultat concerne le calcul de 
l'adjointe d'une matrice avec un co\^{u}t
voisin de celui de son \deterz.
En effet les \coes \,$b_{ij}$\, de  
la  matrice adjointe de \,$A$\, sont donn\'es par:  
$$b_{ij}=(-1)^{i+j}\det{(A\downarrow_{ji})}=  
\frac{\partial~\det{(A)}}{\partial~a_{ji}}$$  
o\`u \,$\det{(A\downarrow_{ji})}$\, est le  
mineur d'ordre \,$n-1$\, obtenu en supprimant la  
\,$j^{\grave{e}me}$\, ligne et la \,$i^{\grave{e}me}$\,  
colonne de la matrice \,$A$.  

\ms 
Nous montrons le r\'esultat par r\'ecurrence sur la longueur du \prev 
qui calcule la fonction.

Supposons donc par exemple qu'un \pol \,$f(x_1,\ldots ,x_n)$\, soit 
calcul\'e par un \prev \,$\Gamma$\, sans division de longueur \,$s$. On 
peut num\'eroter \,$x_{n+1},\ldots ,x_{n+s}$\, les variables du 
programme.
La variable \,$x_{n+1}$\, repr\'esente un \pol \,$g(x_1,\ldots ,x_n)$\, 
de l'un des 4 types suivants: 

\centerline {(1) $x_i+x_j$, \,\,\, (2) $x_i\,{\times}\,x_j$, \,\,\, (3) 
$c+x_i$, \,\,\, (4)  $c\,{\times}\,x_i$}

\noindent avec $1\leq i,j\leq n$ et $c$ une constante.
On a aussi 

\centerline {$f(x_1,\ldots ,x_n)=f_1(x_1,\ldots ,x_n,g(x_1,\ldots 
,x_n))$}

\noindent  o\`u le \pol \,$f_1(x_1,\ldots ,x_n,x_{n+1})$\, est calcul\'e 
par le \prev \'evident \,$\Gamma_1$\, \gui{extrait} de \,$\Gamma$\, et 
de longueur \,$s-1$.

Par \hdr \,$f_1$\, et les \,$n+1$\, d\'eriv\'ees partielles de \,$f_1$\, 
peuvent \^{e}tre calcul\'ees par un \prev \,$\Gamma'_1$\, de longueur 
\,$\ell_1\leq 5(s-1)$.

On consid\`ere alors les formules qui permettent de calculer les 
d\'eriv\'ees partielles de \,$f$\, \`a partir de celles de \,$f_1$\, 
dans les 4 cas envisag\'es pr\'ec\'edemment:
\begin{itemize}
\item  [$(1)$]  $\partial f/\partial x_h = \partial f_1/\partial 
x_h\quad $ si $h\neq i,j$,
\\ $\partial f/\partial x_i = \partial f_1/\partial x_i\,+\, \partial 
f_1/\partial x_{n+1}$,
\\ $\partial f/\partial x_j = \partial f_1/\partial x_j\,+\, \partial 
f_1/\partial x_{n+1}$.
\item  [$(2)$]   $\partial f/\partial x_h = \partial f_1/\partial 
x_h\quad $ si $h\neq i,j$, 
\\ $\partial f/\partial x_i = \partial f_1/\partial x_i\,+\, 
x_j\,{\times}\, \partial f_1/\partial x_{n+1}$,
\\ $\partial f/\partial x_j = \partial f_1/\partial x_j\,+\, 
x_i\,{\times}\,\partial f_1/\partial x_{n+1}$.
\item  [$(3)$]   $\partial f/\partial x_h = \partial f_1/\partial 
x_h\quad $ si $h\neq i$,
\\ $\partial f/\partial x_i = \partial f_1/\partial x_i\,+\, \partial 
f_1/\partial x_{n+1}$.
\item  [$(4)$]   $\partial f/\partial x_h = \partial f_1/\partial 
x_h\quad $ si $h\neq i$,
\\ $\partial f/\partial x_i = \partial f_1/\partial x_i\,+\, 
c\,{\times}\,\partial f_1/\partial x_{n+1}$.
\end{itemize}

C'est le deuxi\`eme cas qui consomme le plus d'instructions nouvelles: 4 
en tout (2 instructions pour calculer $\partial f/\partial x_i$ et 2 
pour $\partial f/\partial x_j$). Il faut par ailleurs rajouter 
l'instruction qui permet de calculer \,$x_{n+1}$\, en fonction des 
\,$x_i$\, pr\'ec\'edents. 

Ceci nous permet donc de construire \`a partir de \,$\Gamma'_1$\, un 
\prev
\,$\Gamma'$\, pour calculer \,$f$\, et ses \,$n$\, d\'eriv\'ees 
partielles. Ce \prev  \,$\Gamma'$\, a une longueur major\'ee par

\centerline {$ \ell_1+1+4\leq 5(s-1)+5=5s$.}

Par ailleurs l'initialisation de la r\'ecurrence est imm\'ediate.

\smallskip Le cas d'un \prev avec divisions se traite de la m\^{e}me 
mani\`ere et aboutit \`a la m\^{e}me majoration.

\newpage \thispagestyle{empty} 
 
\chapter{Notions de complexit\'e}
\label{ChapNoCo}
\minitoc

\subsubsection*{Introduction}   

Ce chapitre est consacr\'e aux notions de \cobi d'une part, directement 
issue de la mod\'elisation du travail des ordinateurs, et de \com \arit 
d'autre part, en relation avec le nombre d'\oparis ex\'ecut\'ees par un 
\algoz. 

Les deux premi\`eres sections sont consacr\'ees \`a la \cobi et 
constituent une pr\'esentation rapide en guise de \gui{rappels}. 

Les trois 
derni\`eres sections d\'ecrivent de mani\`ere pr\'ecise la \com \arit 
des familles de \carisz, elles servent donc de base de travail pour les 
calculs de \com d\'evelopp\'es dans tout le reste de l'ouvrage.

\ss Dans la section \ref{sec arit/algo} nous introduisons les 
classes importantes de \com \arith \,$\SD(f(n),g(n))$. 
Nous discutons le rapport entre \com \arith
(le nombre d'\oparis ex\'ecu\-t\'ees) et \com \bin (le temps 
d'ex\'ecution effectivement utilis\'e lorsqu'on travaille avec des 
entr\'ees repr\'esentant les \elts de l'anneau \,$\A$\, 
convenablement cod\'es).

\ss Ceci nous conduit \`a la notion de \famu de \caris
et aux classes $\NC^k$ qui sont discut\'ees dans la section 
\ref{sec unif}. 

\ss Enfin dans la section \ref{sec modpar} nous discutons 
bri\`evement un mod\`ele de machine \paral  (les PRAMs) 
correspondant aux \caris et assez proche de la pratique des 
architectures \paralsz.

\section{Machines de Turing et Machines \`a Acc\`es Direct}  
\label{secMTMAD} 
Nous donnons ici quelques indications succinctes sur les mod\`eles de 
calcul \algq dans lesquels est prise en compte la taille des objets \`a 
manipuler. Par exemple le temps utilis\'e pour additionner deux entiers 
\'ecrits en base $10$ est manifestement du m\^eme ordre de grandeur que 
la place occup\'ee par l'\'ecriture de ces deux entiers, tandis que 
l'\algo \usl pour la multiplication de deux entiers de tailles 
\,$k$\,  et \,$\ell$\, utilise un temps du m\^eme ordre de grandeur que 
\,$k\times \ell$.

Lorsque dans les ann\'ees 30 des math\'ematiciens et logiciens ont 
r\'efl\'e\-chi \`a la mani\`ere de d\'ecrire en termes pr\'ecis ce 
qu'est 
un calcul \algqz, ils ont abouti \`a des r\'esultats assez vari\'es 
quant \`a la forme, mais identiques quant au fond. Tous les mod\`eles 
\'elabor\'es ont abouti \`a la m\^eme notion de \gui{fonction calculable 
de 
$\NN$  vers $\NN$}. 

\subsubsection*{La \MT abstraite} 

Cependant, c'est Alan Turing qui a emport\'e la conviction par la 
simplicit\'e de son mod\`ele et par son caract\`ere vraiment 
m\'ecanique. 
Il est parti de l'id\'ee qu'un calcul doit pouvoir \^etre ex\'ecut\'e 
par une machine id\'eale qui, \`a l'instar d'un calculateur humain, 
dispose d'une feuille de papier et d'un crayon, et proc\`ede selon une 
suite d'\ops \elrs bien r\'epertori\'ees une fois pour toutes, 
ex\'ecut\'ees conform\'ement \`a un plan de travail d\'etaill\'e ne 
laissant place \`a aucune ambig\"uit\'e. 
Ce mod\`ele est bas\'e sur la notion d'\op \elrz. 
Une telle \op doit \^etre suffisamment simple pour ne consommer qu'une 
quantit\'e fixe de temps et d'\'energie. 
On imagine donc que la machine dispose d'un alphabet fini fix\'e une 
fois pour toutes, et qu'une \op \elr consiste \`a lire, \'ecrire ou 
effacer une lettre \`a un endroit pr\'ecis (la feuille de papier doit 
\^etre divis\'ee en cases, par exemple on prend du papier quadrill\'e), 
ou encore \`a se d\'eplacer vers une case voisine sur la feuille de 
papier. Naturellement on n'autorise qu'un nombre fini de lettres 
distinctes.
Dans le premier mod\`ele, Turing utilise une feuille de papier 
constitu\'ee d'une simple succession de cases sur une seule ligne 
potentiellement infinie: la bande de la machine de Turing. Par la 
suite, il a sembl\'e plus naturel d'utiliser pour mod\`ele une Machine 
de Turing qui utilise plusieurs bandes pour son travail. Quant au crayon 
(muni d'une gomme), il est repr\'esent\'e par ce qu'il est convenu 
d'appeler une t\^ete de lecture{\footnote{~Il serait plus correct mais 
plus lourd de parler d'une t\^ete de 
lecture/effa\-\c{c}age/\'ecriture.}} qui se d\'eplace le long de la 
bande. Il y a une 
t\^ete de lecture pour chacune des bandes. Au d\'epart, certaines bandes 
doivent contenir l'entr\'ee de l'\algo (convenablement cod\'ee), tandis 
que les autres sont enti\`erement vides. Lorsque la machine s'arr\^ete, 
on lit le r\'esultat \`a un endroit convenu. Une t\^ete de lecture est 
capable de reconna\^{\i}tre si la case lue est vide, d'y \'ecrire alors 
une lettre, si elle n'est pas vide de lire la lettre qui s'y trouve et 
\'eventuellement de l'effacer. Pour plus de d\'etails nous renvoyons \`a 
l'ouvrage \cite{TG} o\`u sont traduits et comment\'es les articles 
originaux de Turing, ainsi qu'aux ouvrages \cite{Ste} et \cite{Bal}.

Le caract\`ere tr\`es \elr du fonctionnement abstrait 
de la \MT en a fait un candidat naturel, non seulement 
pour les questions de calculabilit\'e th\'eorique, mais \egmt pour 
les questions de \comz, et en particulier pour la question de 
l'appr\'eciation du 
temps et de l'espace \ncrs \`a l'ex\'ecution d'un \algoz. 
Une fois  l'\algo traduit dans le mod\`ele de la \MTz, le \emph{temps 
d'ex\'ecution} est simplement mesur\'e par le nombre d'\ops \elrs qui 
sont effectu\'ees avant d'aboutir \`a l'arr\^et.
 L'\emph{espace \ncr \`a l'ex\'ecution} est repr\'esent\'e par 
le nombre de cases r\'eellement utilis\'ees.

\subsubsection*{Programmes \elrs} 

On peut donner un mod\`ele \'equivalent \`a la machine de Turing en 
termes de \prgs ex\'ecutables, sans doute plus parlant pour quiconque a 
d\'ej\`a \'ecrit un \prg \ifoz. On consid\`ere des \prgs de nature 
tr\`es simple. Ils sont \'ecrits en utilisant des variables enti\`eres 
\,$N_1,\ldots,\,N_r$\, (les entiers sont suppos\'es \'ecrits en binaire) 
ou \boles \,$B_1,\ldots,\,B_s$\, ($\in\{0,1\}$). Un \gui{\prg \elrz} 
est une suite finie d'instructions num\'erot\'ees de l'un des types 
suivants:
\begin{itemize}
\item [$(A)$] Affectations
\begin{itemize}
\item [$(1)$] $B_j\leftarrow N_i \; \mathrm{mod}\; 2$ 
\item [$(2)$] $N_i\leftarrow N_i \; \mathrm{div}\; 2$ 
\item [$(3)$] $N_i \leftarrow 2N_i + B_j$ 
\item [$(4)$] $B_j\leftarrow 0$ 
\item [$(5)$] $B_j\leftarrow 1$
\end{itemize}
%
\item [$(B)$] Branchements
\begin{itemize}
\item [$(1)$] Direct: aller \`a l'instruction \num $\ldots$ 
\item [$(2)$] Conditionnel bool\'een: si $B_j=0$ aller \`a 
l'instruction \num $\ldots$
\item [$(3)$] Conditionnel entier: si $N_i=0$ aller \`a l'instruction 
\num $\ldots$
\end{itemize}
%
\item [$(S)$] Arr\^et.
\end{itemize}
Les variables sont toutes initialis\'ees \`a $0$ sauf celles qui 
repr\'esentent les entr\'ees du \prgz.

Puisque les entiers sont \'ecrits en binaire, 
on voit que chaque affectation ou branchement peut correspondre \`a un 
travail r\'ealis\'e en consommant un temps et une \'energie 
ind\'ependantes de 
l'\'etat des variables. Le temps d'ex\'ecution est donc raisonnablement 
estim\'e comme \'etant le nombre d'instructions ex\'ecut\'ees avant 
d'aboutir \`a l'arr\^et.

\subsubsection*{Machines \`a Acc\`es Direct} 
Comme tout mod\`ele abstrait, la machine de Turing est une 
id\'ealisa\-tion. 
Le point le plus contestable est l'hypoth\`ese implicite selon 
laquelle une \op \elr est \'equivalente \`a une autre quel que soit 
l'\'etat de la bande (dans la version Machine) ou des variables (dans la 
version \prg \ifo \elrz). Une telle conception se heurte \`a des 
limitations physiques. Elle n'est en tout cas pas conforme \`a ce qui se 
passe concr\`etement dans les ordinateurs actuels.

Alan Turing participa \`a l'aventure des premiers ordinateurs. Les 
ordinateurs ont une conception globale qui diff\`ere sensiblement de la 
\MT abstraite. Les donn\'ees ne sont pas trait\'ees \gui{l\`a 
o\`u elles sont}, comme dans l'image du crayon qui se d\'eplace sur la 
feuille de papier, mais elles sont transf\'er\'ees depuis la 
p\'eriph\'erie (un disque dur par exemple) vers le centre o\`u elles 
sont trait\'ees, 
\cad vers un microprocesseur, avant d'\^etre renvoy\'ees vers la 
p\'eriph\'erie.  
Ces transferts permanents prennent d'autant plus de temps que les 
donn\'ees sont plus \'eloign\'ees et que l'espace \ncr \`a leur 
stockage est plus grand.  

Ceci a donn\'e lieu \`a un autre mod\`ele de 
calcul, le mod\`ele MAD des \emph{Machines \`a Acc\`es Direct} 
(RAM en version anglaise abr\'eg\'ee), avec de nombreuses variantes. 
Dans un mod\`ele MAD, on doit consid\'erer une infinit\'e potentielle de 
\gui{registres} 
(correspondant au stockage des donn\'ees en m\'emoire, ou aux cases 
d'une bande de \MTz). 
Il serait logique (mais ce n'est pas 
en \gnl l'option choisie), de consid\'erer que chaque registre ne 
contient qu'une information dont la taille est fix\'ee une fois pour 
toutes. Pour traiter le registre dont l'adresse est l'entier $n$, on 
consid\`ere que l'\op de transfert vers l'unit\'e centrale requiert un 
temps \'egal \`a la taille en binaire de l'entier \,$n$. Dans le 
mod\`ele de Turing, le temps correspondant peut \^etre nul mais aussi 
beaucoup plus grand que \,$\log\,n$, selon la position des t\^etes de 
lecture sur chaque bande. 

En fin de compte, selon l'\algo utilis\'e (et selon le mod\`ele MAD
choisi), les temps 
d'ex\'ecution \,$T$\,  et  \,$T'$\, obtenus dans le mod\`ele MT (\MT \`a 
plusieurs bandes) ou dans les mod\`eles MAD pour une entr\'ee de taille 
\,$t$\, sont soumis \`a des majorations respectives du type suivant
(voir un exemple pr\'ecis dans \cite{Ste} chapitre 2, sections 5.5 et 
5.6):
$$ a\,T'\leq T\leq b\,T'^2),\quad c\,T\leq T'\leq d\,T\,(T+t)^2\,.
$$

Signalons le terme d'\emph{accumulateur} qui dans le mod\`ele MAD 
d\'esigne le microprocesseur.   

\subsubsection*{L'espace de travail proprement dit} 

 Nous terminons cette section avec un commentaire et une d\'efinition 
plus pr\'ecise de \gui{l'espace de travail} utilis\'e dans les mod\`eles 
MT ou MAD. Dans le mod\`ele MT nous avons d\'efini \gui{l'espace \ncrz} 
comme le nombre total de cases effectivement utilis\'ees au cours de 
l'ex\'ecution de l'\algoz. En fait, si on veut \'etudier l'espace de 
travail proprement dit utilis\'e par un \algoz, il est judicieux  
d'op\'erer une distinction entre l'espace \ncr aux donn\'ees 
d'entr\'ee-sortie d'une part, et l'espace \ncr au travail 
proprement dit d'autre part. On convient dans ce cas que
les bandes contenant les entr\'ees sont utilis\'ees en lecture 
uniquement et qu'elles sont lues en une seule passe. De m\^eme, 
les bandes contenant les sorties sont utilis\'ees en \'ecriture 
uniquement, et elles sont \'ecrites en une seule passe.

Par exemple lorsqu'on veut faire la preuve par $9$ pour un produit 
\,$a\times b=c$\, o\`u \,$a$, \,$b$\, et \,$c$\,  sont consid\'er\'ees 
comme des entr\'ees \'ecrites en base $10$, il suffit de lire en une 
seule passe les donn\'ees et aucun stockage des r\'esultats 
\itmds n'est \ncrz. On donne \`a la fin le r\'esultat 
(oui, ou non) sans avoir utilis\'e aucun espace pour le travail 
proprement 
dit{\footnote{~De m\^eme si on veut additionner deux entiers il suffit 
de les lire en une seule passe et d'\'ecrire au fur et \`a mesure le 
r\'esultat sur la bande de sortie. Cependant les entr\'ees et la sortie 
ne sont pas \'ecrites dans le m\^eme sens. En effet, pour pouvoir 
encha\^{\i}ner des \algosz, la convention naturelle est que la t\^{e}te 
de lecture sur chaque entr\'ee doit \^{e}tre au d\'epart \`a 
l'extr\'emit\'e droite de l'entr\'ee, et la t\^{e}te de lecture sur 
chaque sortie doit \^{e}tre \`a la fin \`a l'extr\'emit\'e droite de la 
sortie}}. Si on \'ecrivait cela 
sous forme d'un \prg \ifo \elr du type  que nous avons d\'ecrit 
ci-dessus, cela signifierait que les variables de travail sont toutes 
\bolesz, que les variables repr\'esentant les entr\'ees sont 
seulement utilis\'ees en lecture (elles ne peuvent \^etre utilis\'ees 
que via les affectations $A1$ et $A2$) et les variables repr\'esentant 
la sortie sont seulement utilis\'ees en \'ecriture (elles ne peuvent 
\^etre utilis\'ees que via les affectations $A3$). 

Ainsi certains \algos utilisent un espace de travail nul (dans le cas 
optimal) ou nettement inf\'erieur \`a la taille des entr\'ees-sorties. 
Pour les \'etudes de \com d'\algos on est particuli\`erement 
int\'eress\'e par ceux qui n'utilisent aucun espace de travail d'une 
part, par ceux qui utilisent un espace de travail \lin par rapport 
\`a la taille de l'entr\'ee d'autre part et enfin par ceux qui utilisent 
un espace de travail de l'ordre de grandeur de \,$C\,\log(n)$\, o\`u
\,$C$\, est une constante et \,$n$\, est la taille de l'entr\'ee. On 
appelle ces derniers des \algos $\LOSP$.

\section{Complexit\'e binaire, les classes $\P$, $\NP$ et $\DieseP$}
\label{secPNP}
\subsection{Calculs faisables} 

Malgr\'e la grande abondance des mod\`eles de calcul propos\'es, un 
consensus a fini par s'\'etablir sur ce qu'est un \emph{calcul 
faisable}.
On dit qu'un calcul est faisable, ou encore qu'il est dans la classe 
$\P$ \indexnota{P@$\P$}
si on conna\^{\i}t un \algo qui dans les mod\`eles MT ou MAD 
n\'ecessite un temps \poll par rapport \`a la taille de l'entr\'ee. Plus 
pr\'ecis\'ement, on ne dit rien concernant tel calcul isol\'e (celui des 
100.000 premi\`eres d\'ecimales de $\pi$ par exemple), mais on dit 
quelque chose concernant un calcul \gnl correspondant \`a des 
entr\'ees de tailles variables et en tout cas arbitrairement grandes 
(celui de la \,$k\,$\eme d\'ecimale de $\pi$ par exemple). On demande 
que,
pour un certain \poly \`a \coes positifs ou nul \,$P$,  pour toute 
entr\'ee de taille inf\'erieure ou \'egale \`a \,$n$, l'\algo donne sa 
r\'eponse en un temps major\'e par $P(n)$. Les \algos $\LOSP$  sont dans 
la classe \,$\P$, et ils sont consid\'er\'es \`a juste titre comme bien 
meilleurs que les \algos qui travailleraient en temps \emph{et} espace 
\pollz. 

On voit que la notion d'\algo de classe \,$\P$\, est une notion 
asymptotique, qui peut \^{e}tre assez \'eloign\'ee de la r\'ealit\'e des 
calculs. Un \algo ayant un temps de calcul \gui{\linz} \'egal \`a 
\,$n+10^{100}$\,  correspond en pratique \`a quelque chose d'infaisable, 
tandis que si son temps de calcul est \gui{exponentiel} major\'e par 
\,$\sup(n,2^n/2^{10^{100}})$\,  il reste facile \`a ex\'ecuter pour 
toutes les entr\'ees concr\`etement 
envisa\-geables, alors m\^eme qu'il n'est pas dans la classe \,$\P$. 

De nombreux auteurs distinguent les \emph{\pbs faisables} (les 
entr\'ees sont des entiers, ou cod\'ees par des entiers, mais la sortie 
est du type oui/non, donc cod\'ee par un \bol dans \,$\{0,1\}$) des 
\emph{fonctions faisables} (la ou les sorties sont des entiers) et ils 
r\'eservent le symbole $\P$ pour les \pbs faisables. 
La classe des fonctions faisables
(calculables \etpz) est alors not\'ee  \,$\FP$. En fait une fonction 
\,$f:\NN\rightarrow \NN$\, est faisable \ssi d'une part la taille de la 
sortie est \polt major\'ee en fonction de la taille de l'entr\'ee, et 
d'autre part, le \pb \,$f(n)\leq p\;?$\, est dans la classe $\P$.  
Nous n'introduirons donc pas deux notations distinctes et nous ferons 
confiance au contexte pour lever les ambig\"uit\'es \'eventuelles.

\ss Citons des \pbs de base qui ont re\c{c}u dans le pass\'e une 
solution \algq satisfaisante, ce qui les mettait dans la classe $\P$ 
bien avant qu'elle ne f\^{u}t invent\'ee. 
La r\'esolution des \slis par la \met chinoise du pivot, 
appel\'ee en Occident \mpgz, 
donne un \algo de classe $\P$ lorsque les \coes et les inconnues sont 
des nombres rationnels. Le calcul du nombre de racines r\'eelles d'un 
\poly par la \met de Sturm, qui avait \'et\'e salu\'ee lors de sa 
d\'ecouverte pour sa clart\'e et son \'el\'egance, fournit un \algo 
\etp lorsque les \coes du \pol sont des nombres rationnels. Le calcul du 
\polcar d'une \maca par la \met de Leverrier est un autre exemple 
c\'el\`ebre.  
Le calcul int\'egral lui-m\^{e}me a un aspect \algq (pour 
le calcul automatique de certaines aires par exemple) qui frappa les 
contemporains de Leibniz et Newton et qui est devenu aujourd'hui une des 
branches du calcul formel.
\subsection[Quand les solutions sont faciles \`a tester]{Probl\`emes 
dont les solutions sont faciles \`a tester}  
 La conjecture $\P\neq \NP$ est apparue dans les ann\'ees 70 
(Cook, \cite{Cook}). 
Elle correspond \`a l'id\'ee intuitive suivante: il y a des 
\pbs dont les solutions sont faciles \`a tester mais qui sont 
difficiles \`a r\'esoudre.
On pourrait dire a priori que la plupart des \syses 
qu'on cherche \`a r\'esoudre correspondent \`a ce paradigme. 
Il est remarquable que cette id\'ee intuitive n'ait pu recevoir une 
forme math\'ematique pr\'ecise qu'avec l'av\`enement de la th\'eorie de 
la \com des \algosz. Tard venue dans le monde des conjectures
math\'ematiques, la conjecture $\P\neq \NP$ appara\^{\i}t aujourd'hui 
comme l'une des plus importantes, l'une dont la signification est la 
plus profonde. Elle a r\'esist\'e \`a toutes les tentatives d'en venir 
\`a bout, et beaucoup d'experts pensent qu'on ne dispose pas aujourd'hui 
des concepts \ncrs \`a sa solution, alors m\^eme qu'elle a 
quasiment la force d'une \'evidence.  
Nous allons en donner quelques commentaires relativement informels. Ils 
sont  \ncrs pour aborder dans les chapitres 
\ref{chap DeterUniv} et  \ref{chap Perma}
l'analogue en \coag de la conjecture  $\P\neq \NP$ en \cobiz.
Nous recommandons l\`a encore sur ce sujet les ouvrages \cite{Bal} et 
\cite{Ste}. 

Comme exemple de \pb dont les solutions sont faciles \`a tester 
mais qui sont difficiles \`a r\'esoudre, nous allons consid\'erer les 
\pbs de programmation \linz. Un tel \pb est donn\'e 
par une matrice $A$ (de type $n\times m$) et un vecteur colonne $b$ (de 
type $m\times 1$) \`a \coes r\'eels, et une solution du \pb est 
un vecteur colonne $x$ (de type $n\times 1$) tel que le vecteur $y=Ax-b$ 
ait toutes ses coordonn\'ees $\geq 0$~({\footnote{~Un \pb de 
programmation \lin est en \gnl \'enonc\'e sous forme d'un 
\pb d'optimisation. Nous en pr\'esentons ici une version 
\'equivalente plus facile \`a discuter pour notre propos actuel.}}).
Pour en faire un \pb dont la nature \algq est bien pr\'ecise, 
nous nous limitons aux matrices \,$A'=[A\vert b]$\, \`a \coes entiers 
cod\'es en binaire. Quant aux solutions, nous avons le choix. Si nous 
demandons des solutions en nombres rationnels, nous parlons de 
programmation \lin en rationnels, et si nous demandons des 
solutions en nombres entiers, nous parlons de programmation \lin 
en entiers. Pour chacun de ces deux \pbs une solution \,$x$\, 
\'eventuelle est facile \`a tester. Un \algo qui donne en \gnl 
une solution rapide (s'il en existe une) pour la programmation 
\lin en rationnels a \'et\'e mis au point dans les ann\'ees 50 
(\cf \cite{Dan}). Il est en \gnl tr\`es performant et il est 
encore aujourd'hui fr\'equemment utilis\'e, c'est l'agorithme de 
Dantzig. L'inconv\'enient est que pour certaines matrices \,$A'$, 
l'\algo a un mauvais comportement et son temps de calcul peut devenir 
exponentiel par rapport \`a la taille de \,$A'$. Dans les ann\'ees 70 on
a trouv\'e d'autres \algosz, qui dans la plupart des cas sont nettement 
plus lents que celui de Dantzig, mais qui tournent \etp pour n'importe 
quelles matrices \,$A'$\, (\cf \cite{Karm,Kha1,Kha2} et l'ouvrage 
\cite{Schri}).  
Depuis, on sait donc que la programmation \lin en rationnels est 
dans la classe $\P$. Par contre pour ce qui concerne la programmation 
\lin en entiers, on n'est toujours pas capable de r\'esoudre ce 
\pb par un \algo de la classe $\P$, m\^eme si on ne s'int\'eresse 
qu'aux solutions de petite taille. En fait, on pense qu'on en sera \`a 
tout jamais incapable, car une r\'eponse dans l'autre sens signifierait 
que la conjecture \,$\P\neq \NP$\, est fausse. 

Pour expliquer comment est d\'efinie la classe $\NP$, nous essayons 
d'examiner avec un peu de recul ce que signifierait en \gnl 
\gui{savoir r\'esoudre un \pb dont on sait tester facilement les 
solutions}. Nous commen\c{c}ons par remarquer que pour bien poser la 
question, il faut savoir donner le \pb sous une forme cod\'ee, 
qui puisse \^etre prise comme entr\'ee d'un \prg \ifo (ou d'une Machine 
de Turing). On peut donc toujours consid\'erer que l'on a une suite 
infinie de \pbs \,$P_n$\, o\`u \,$n$\, est justement la forme cod\'ee en 
binaire du \pb (les entiers \,$n$\, qui ne coderaient pas 
correctement une instance de notre \pb doivent pouvoir \^etre 
faciles \`a rep\'erer). Quant aux solutions, elles doivent \egmt 
pouvoir \^etre cod\'ees, donn\'ees comme telles en entr\'ee ou \`a la 
sortie d'un \prg \ifoz. Nous supposons donc \spdg que la solution 
\'eventuelle est elle aussi cod\'ee par un entier \,$x$. Maintenant 
consid\'erons la fonction 
\,$\varphi:\NN\times \NN\rightarrow \{0,1\}$\, qui est d\'efinie comme 
suit:     
\,$\varphi(n,x)=1$\, si \,$x$\,  est le code d'une solution du \pb 
\,$P_n$\,  
et  
\,$\varphi(n,x)=0$\, sinon.  Supposer qu'on sait tester facilement les 
solutions de notre famille \,$P_n$\, peut \^etre raisonnablement 
interpr\'et\'e comme signifiant que la fonction \,$\varphi$\, est dans 
la 
classe \,$\P$. Tandis que supposer que le \pb est intrins\`equement 
difficile \`a r\'esoudre
 peut \^etre raisonnablement interpr\'et\'e comme signifiant que la 
question
$$ 
\exists x\in\NN \qquad   \varphi(n,x)=1\; \; ?
$$
n'a pas de r\'eponse dans la classe $\P$. 
Maintenant, nous devons apporter une restriction. Il se peut que le 
\pb soit intrins\`equement difficile \`a r\'esoudre pour une trop 
bonne raison, \`a savoir que les solutions \'eventuelles sont de taille 
trop grande. Plus exactement que la taille de toute solution \,$x$\, du 
\pb \num$n$ croisse trop vite par rapport \`a celle de \,$n$. 
\emph{Nous notons dans la suite de cette section \,$|x|$\,  la taille de 
l'entier naturel \,$x$,} 
\cad la longueur de son \'ecriture en binaire.
Nous pouvons maintenant \'enoncer ce qu'est un \pb dans la classe 
$\NP$. C'est r\'epondre \`a une question du type suivant:\\
\emph{Existe-t-il une solution \,$x(n)$\, de taille raisonnable pour 
telle famille 
\,$P_n$\, de \pbs dont les solutions sont faciles \`a tester~?}\\
Plus pr\'ecis\'ement une famille de \pbs cod\'ee dans \,$\NN$\, est 
dite dans la classe $\NP$  si sa solution revient \`a r\'esoudre une 
question du type
\begin{equation} \label{eqNP1}
\exists x\in\NN \qquad  (|x|\leq a+|n|^k\;\;\mathrm{et} \;\;
\varphi(n,x)=1)\quad ?
\end{equation}
o\`u \,$a$\,  et \,$k$\,  sont deux entiers positifs donn\'es et o\`u 
\,$\varphi:\NN\times \NN\rightarrow \{0,1\}$\, est dans la classe $\P$.
Autrement dit si on pose
\begin{equation} \label{eqNP}
\psi(n)=\sup\left\{\varphi(n,x)\;;\;|x|\leq a+|n|^k\right\}
\end{equation}
et si la fonction \,$\varphi$\,  est dans la classe $\P$, alors la 
fonction 
\,$\psi$\, est dans la classe $\NP$. On peut d'ailleurs supposer \spdg 
que 
\,$\varphi(n,x)=0$\, si \,$|x|> a+|n|^k$.

Le ${\cal N}$  de $\NP$ est mis \indexnota{NP@$\NP$}
pour \emph{non d\'eterministe}\index{non d\'eterministe}. 
La raison en est la suivante. La fonction \,$\psi$\, ci-dessus pourrait 
\^etre calcul\'ee \etp par une machine dont le fontionnement serait 
\gui{non d\'eterministe}. Plus pr\'ecis\'ement, en utilisant nos \prgs 
\ifos \elrs ci-dessus, on admettrait des intructions de 
branchement non d\'eter\-mi\-nis\-te: aller \`a l'instruction \num 
$\ldots$  ou  $\ldots$  (selon l'humeur du moment). 
Le \prg peut alors aboutir \`a plusieurs r\'esultats diff\'erents selon 
le chemin choisi lors de son ex\'ecution. 
Et ce \prg serait r\'eput\'e calculer (pour une entr\'ee $x$ fix\'ee) la 
plus grande des valeurs qu'il peut d\'elivrer en sortie. 
L'acronyme $\NP$ vaut alors pour: calculable \etp par une machine \`a 
fonctionnement non d\'eterministe. 

Notez que si on avait $\P=\NP$ (ce que personne ne croit), on pourrait 
non seulement calculer, dans l'exemple ci-dessus, la fonction \,$\psi$\, 
\etpz,
mais \egmtz, dans le cas d'une r\'eponse positive \,$\psi(n)=1$, 
trouver une solution \,$x$\, pour \,$\varphi(n,x)=1$\,  \etpz. En effet,
on pourrait calculer un tel \,$x$\, par dichotomie \etp en posant un 
nombre \poll de fois la question  

\centerline{ \,$\exists x \quad   x\leq p\;\,\mathrm{et} \; 
\varphi(n,x)=1\quad ?$ 
}
 
\noi qui serait r\'esoluble \etp sur les entr\'ees \,$n,\,p$\, 
(on d\'e\-mar\-re\-rait avec \,$p=2^{a+|n|^k}$).

Certains \pbs qui peuvent sembler a priori \^{e}tre dans la 
classe $\NP$ sont ramen\'es dans la classe $\P$ lorsque quelqu'un 
d\'ecouvre un \algo rapide pour les r\'esoudre. 
Des succ\`es spectaculaires ont \'et\'e \`a la fin du 20\`eme si\`ecle 
la solution \etp des \syses \lins \`a \coes et 
inconnues enti\`eres, celle des \pbs de programmation \lin 
en rationnels et la d\'etermination de \gui{petits vecteurs} dans un 
r\'eseau (qui conduit notamment \`a la factorisation \etp des \polys sur 
\,$\QQ[X]$).  

\ss Cook a montr\'e (\cf \cite{Cook}) que certains \pbs de la 
classe $\NP$ sont \uvlsz: si on d\'emontre pour l'un d'entre eux qu'il 
est dans la classe $\P$, alors $\P=\NP$. 
Un tel \pb est dit \NPcz. \index{NPcom@\NPc} 
Par exemple la programmation \lin en entiers est un \pb 
\NPcz, m\^eme si on limite a priori la taille des solutions par 
un entier fixe.

Nous pouvons expliquer informellement pourquoi il existe des \pbs 
\NPcsz. 

Un ordinateur qui ne serait soumis \`a aucune limitation physique
de temps et d'espace serait une \emph{machine \uvlez} en ce sens qu'il 
est capable d'ex\'ecuter n'importe quel \prg qu'on lui soumet
(en faisant abs\-trac\-tion des limitations physiques).
Un des premiers \thos d'Alan Turing \'etait \emph{l'existence d'une \MT 
\uvlez.} 
Une cons\'equence importante de l'existence d'une \MT \uvle est, via le 
processus diagonal de Cantor, l'existence de \pbs bien pos\'es 
(pour les \MTsz) mais qui ne pourront \^etre r\'esolus par aucun 
proc\'ed\'e m\'ecanique du type \MTz: l'ensemble des (codes Turing de) 
fonctions m\'ecaniquement calculables de \,$\NN$\, vers \,$\NN$\, 
(au sens des \MTsz) n'est pas m\'ecaniquement calculable (au sens des 
\MTsz).  
L'existence de \pbs \NPcs est un r\'esultat de nature similaire.

Introduisons la notation \,$\gen{x_1,\ldots,x_k}\,$
\indexnota{a x@$\gen{ x_1,\ldots,x_k}$, (code d'un  \,$k$-uple)}
pour un \emph{code dans \,$\NN$\, d'un  \,$k$-uple 
d'entiers}\index{code d'un $k$-uple d'entiers}{\footnote{~On consid\`ere 
un codage naturel, de sorte que les fonctions de codage 
\,$(x_1,\ldots,x_k)\mapsto \gen{x_1,\ldots,x_k}$\,  et celles de 
d\'ecodage ($\gen{x_1,\ldots,x_k}\mapsto x_i$\,  et 
\,$\gen{x_1,\ldots,x_k}\mapsto k$) sont dans la classe $\P$.
On suppose aussi \spdg que \,$\gen{x_1,\ldots,x_k}\geq x_i$.}}.
En termes de \prgs \ifos \elrsz, l'existence d'une \MT \uvle signifie 
qu'on sait \'ecrire un \prg \elr \gui{\uvlz} en ce sens qu'il remplit le 
contrat suivant:
\begin{itemize}
\item Il prend en entr\'ee 2 entiers binaires \,$n,\,x$\,  et un entier 
b\^aton\index{entier b\^aton}  \,$t$\, ({\footnote{~Un entier b\^aton 
sert 
de compteur, il est cod\'e en binaire par \,$2^t-1$, \cad (si \,$t\geq 
1$) le mot form\'e de \,$t$\, fois la lettre $1$.  
Ici il est \ncr de prendre pour \,$t$\, un entier b\^aton parce 
qu'on veut que la fonction \uvle soit calculable \etp par rapport \`a la 
taille de ses entr\'ees.}}), o\`u \,$n$\, est un texte de \prg \elr 
\,$Q_n$\, cod\'e en binaire, \,$x$\, est un code pour la liste des 
entr\'ees pour \,$Q_n$\,  et \,$t$\, est le nombre d'\'etapes \elrs 
pendant lequel on d\'esire que soit ex\'ecut\'e \,$Q_n$. 
\item Il donne en sortie une \emph{description instantan\'ee 
de}\index{description!instantan\'ee}  (\cad un
 codage binaire \,$U(n,x,t)$\, qui d\'ecrit de mani\`ere exacte) 
\emph{l'\'etat o\`u se trouve la machine qui ex\'ecute le \prg \,$Q_n$\, 
apr\`es 
l'ex\'e\-cution de \,$t$\, \'etapes \elrs de calcul sur l'entr\'ee 
\,$x$}: la valeur de chacune des variables \,$x_i$\, du \prg d'une 
part, le num\'ero \,$h$\, de l'instruction en cours d'autre part 
(cod\'es par  \,$\gen{x_1,\ldots,x_\ell,h}$). 
\end{itemize}
Si le temps d'ex\'ecution est \,$t_0$\, on demande que pour \,$t>t_0$\, 
on ait  \,$U(n,x,t)=U(n,x,t_0)$. Nous supposons aussi \spdg que les 
variables de sortie sont en \'ecriture seulement, \cad ne sont 
utilis\'ees que via les affectations de type $A3$.

Il n'est pas tr\`es difficile de v\'erifier qu'un \prg \elr \uvl \'ecrit 
de mani\`ere naturelle calcule la fonction \uvle 
\,$U$\, \etpz.  

Comme cons\'equence on obtient quelque chose qui poura\^{\i}t \^etre 
compris comme une \'enum\'eration dans la classe $\P$ de tous les \prgs 
dans la classe $\P$ s'ex\'ecutant sur une entr\'ee de taille \polt 
major\'ee. Expliquons nous.

Tout d'abord notons \,$(n,x,t)\mapsto V(n,x,t)$\, la fonction (dans la 
classe $\P$) qui donne l'\'etat de la variable en sortie (ou, s'il y a 
plusieurs sorties pr\'evues, l'\'etat de la premi\`ere d'entre elles).

Soit maintenant \,$\varphi:\NN\times \NN\rightarrow \{0,1\}$\, une 
fonction dans la classe $\P$ qui v\'erifie \label{NPcomp}
\begin{equation} \label{eqNP2}
\forall x\quad |x|> a+|n|^k \;\Rightarrow \;\varphi(n,x)=0
\end{equation}
Alors la fonction \,$\psi(n)=\sup\left\{\varphi(n,x)\;;\;x\in 
\NN\right\}$\, r\'esoud un pobl\`eme dans la classe $\NP$. 

Vu que \,$\varphi$\,  est dans la classe $\P$ et vu le caract\`ere \uvl 
de \,$V$\, il existe un entier \,$m_0$\, et deux entiers \,$b,\,\ell$\, 
tels que 
$$ \varphi(n,x)=V(m_0,\gen{n,x},t)\quad \mathrm{avec}\quad t\leq 
b+|n|^\ell
\quad \mathrm{si}\quad |x|\leq a+|n|^k
$$
D\'efinissons par ailleurs (avec \,$z,\,t$\, des entiers b\^atons, et 
\,$n$\, et \,$x$\,  des entiers binaires,)
$$ \Phi(p,x)=
\formul{ll}{ \inf(1,V(m,\gen{n,x},t)) &\mathrm{si}\;  |x|\leq z 
\cr 0&  \mathrm{sinon}
}
\quad \mathrm{avec}\quad p=\gen{m,n,z,t}.
$$
C'est naturellement une fonction \,$\NN\times \NN\rightarrow \{0,1\}$\, 
dans la classe $\P$ pour laquelle on a 
$$ 
\forall x\quad |x|> |p| \;\Rightarrow \;\Phi(p,x)=0
$$
et \`a partir de laquelle on peut d\'efinir
\begin{equation} \label{eqNPc}
\Psi(p)=\sup\left\{\Phi(p,x)\;;\;|x|\leq |p|\right\}
\end{equation}
qui est dans la classe $\NP$.
Maintenant il est clair que si on pose
$$ 
\lambda(n)=\gen{m_0,n,a+|n|^k,b+|n|^\ell}
$$
alors la fonction \,$\lambda$\, est dans la classe \,$\P$\,  et 
$$ 
\psi(n)=\Psi(\lambda(n)).
$$
Ceci montre le caract\`ere \uvl de la fonction \,$\Psi\,$
au sens o\`u nous le souhaitions. 
En \'ecrivant cette preuve en d\'etail, on peut donner des pr\'ecisions 
\supts sur la mani\`ere dont le \pb \,$\NP$\, associ\'e 
\`a \,$\varphi$\, \`a \'et\'e r\'eduit \etp \`a celui associ\'e \`a 
\,$\Phi$. 
En particulier on peut fabriquer une variante o\`u la 
r\'eduction est dans la classe $\LOSP$.

\subsection{Probl\`emes de comptage}  
\label{comptage}
Si \,$E$\, est un ensemble fini, nous noterons \,$\#E$\, le nombre 
d'\elts de \,$E$. \indexnota{DieseE@$\#E$ nombre d'\elts de \,$E$}

Lorsqu'on a une \fam de \pbs dont les solutions sont faciles \`a 
tester et
de taille \polt major\'ees, on peut se poser non seulement la question 
de savoir si une solution existe, mais \egmt combien de solutions 
existent.  \indexnota{DieseP@$\DieseP$}

Pr\'ecis\'ement si  \,$\varphi:\NN\times \NN\rightarrow \{0,1\}$\, est 
une fonction dans la classe $\P$ qui v\'erifie
$$ \forall x\quad |x|> a+|n|^k \;\Rightarrow \;\varphi(n,x)=0
$$
alors la fonction
\begin{equation} \label{eqDieseP}
\theta(n)=\#\left\{x\; |\;x\in \NN,\;  \varphi(n,x)=1\;\right\}=
\sum_{|x|\leq a+|n|^k}{\varphi(n,x)}
\end{equation}
compte le nombre de solutions (pour la question cod\'ee par \,$n$).
A priori cette fonction est plus difficile \`a calculer que
la fonction d\'efinie par l'\'equation (\ref{eqNP}) 
\,$\psi(n)=\sup\left\{\varphi(n,x)\;;\;x\in \NN\right\}\,$
(qui est dans la classe \,$\NP$). La taille de \,$\theta(n)$\,  est 
\polt 
major\'ee en fonction de celle de \,$n$.  
Les fonctions \,$\theta$\, obtenues de cette mani\`ere d\'efinissent une 
nouvelle classe de complexit\'e, \emph{les fonctions de comptage pour 
les \pbs dont les solutions sont faciles \`a tester}, que l'on 
note $\DieseP$ (prononcer di\`ese P). Cette classe a \'et\'e introduite 
par Valiant dans \cite{Val2}. 

Si on veut que la classe \,$\DieseP$\,  soit une classe 
de \pbs plut\^ot qu'une classe de fonctions, on la d\'efinit 
comme la classe des \pbs du type \,$\theta(n)\leq p\; ?$.
En effet, puisque \,$\theta(n)\leq2^{a+|n|^k}$\, il est facile de 
calculer par dichotomie, \etpz, la fonction
\,$\theta$\, \`a partir des tests \,$\theta(n)\leq p\; ?$.

On conjecture que les deux 
inclusions 
$$ \P\subset \NP\subset \DieseP$$ 
 sont strictes. 

De m\^eme qu'il existe des \pbs \NPcsz, il existe des 
fonctions $\DieseP$-compl\`etes. En fait la r\'eduction que nous avons 
esquiss\'ee dans le cas $\NP$ ci-dessus fonctionne aussi pour les 
fonctions de comptage.
D\'efinissons en effet
\begin{equation} \label{eqDPc}
\Theta(p)=\#\left\{x\; |\;x\in \NN,\;  \Phi(p,x)=1\;\right\}=
\sum_{|x|\leq |p|}{\Phi(p,x)}
\end{equation}
alors, avec la m\^eme fonction \,$\l$\,  que ci-dessus, on obtient
\,$ \theta(n)=\Theta(\lambda(n)).$

\section[Complexit\'es \arith et \bin] 
{Complexit\'e \arith et \com \bin des circuits}
 
\label{sec arit/algo}

\subsection{Complexit\'e arithm\'etique}   
\index{complexit\'e!arithm\'etique}
La taille (en fait le nombre d'\oparisz)  
et la \prof d'un \cari ou d'un \prev
 sont les deux \parats qui 
 mesurent ce qu'on appelle la \emph{\com \arithz}  
de ce \cari ou de ce \prevz.  
  
Ce sont des fonctions de ce que nous avons appel\'e les  
\parats d'entr\'ee du \cariz. Comme on  
s'int\'eresse souvent \`a la \com asymptotique des  
\algos (\cad \`a leur comportement quand ces \parats  
tendent vers l'infini), nous allons utiliser les  
notations classiques \,$\O,\,o,\,\Theta,\,\Omega$\,  
d\'efinies de la mani\`ere suivante:  
  
\begin{nota}   
\'Etant donn\'ees deux fontions \,$f$\, et \,$g$\, de  
\,$\NN^*$\, dans \,$\RR_{+}^*$, on dit que:  
\indexnota{o(f)@$\o(f)$}\indexnota{O(f)@$\O(f)$}

$\bullet~$ \,$g\in \O(f)$\, et l'on \'ecrira  
\,$g(n)=\O(f(n))$\, s'il exis\-te une constante   
r\'eelle \,$c>0$\, telle que 
\,$\forall n,\,n\in \NN^*\impliq g(n) \leq c\,f(n)$.  
  
$\bullet~$ \,$g\in \o(f)$\, et l'on \'ecrira  
\,$g(n)=\o(f(n))$\, si pour tout r\'eel  
\,$\varepsilon > 0 $, il existe \,$k\in\NN^*$\, tel que  
\,$\forall n,\,(n\in \NN^* \mbox{ et } n>k) \impliq  
g(n) \leq \varepsilon\,f(n)$.  
  
$\bullet~$ \,$g\in \Omega (f)$\, et l'on \'ecrira  
\,$g(n)=\Omega(f(n))$\, si \,$f(n)=\O(g(n))$. 
\indexnota{Omega@$\Omega(f)$}

$\bullet~$ \,$g\in \Theta(f)$\, et l'on \'ecrira  
\,$g(n)=\Theta(f(n))$\, si \,$g(n)=\O(f(n))$\,  
et \,$f(n)=\O(g(n))$. On dit dans ce cas que \,$f$\,  
est du m\^eme ordre que \,$g$.
\indexnota{oTeta@$\Theta(f)$}
\end {nota}  

Remarquons que pour montrer que \,$g\in \O(f)$, il suffit  
de trouver une constante r\'eelle \,$K_0$\, et un entier  
\,$n_0 \in \NN^*$\, tels que \,$g(n) \leq K_0\,f(n)$\,  
pour tout \,$n \geq n_0$.  
Nous appellerons une telle constante \,$K_0$\, une
\emph{constante asymptotique (cach\'ee dans le grand~$\,\O$)}.
Dans la suite chaque fois que ce sera possible nous  
nous appliquerons \`a faire appara\^{\i}tre la constante
asymptotique cach\'ee dans le grand~$\,\O$\,  
dans l'\'etude de \com des \algosz. Et l'entier
 \,$n_0$\, sera parfois pr\'ecis\'e.

\begin{nota}\emph{(\comar d'une \fam de \cirsz)}\label{notaSD}~\\   
On \'ecrira qu'un \algo est (dans la classe)
\index{complexit\'e!arithm\'etique d'une \fam de \cirs}
 \,$\SD(f(n),g(n))$\, 
 \indexnota{SD@$\SD(f(n),g(n))$} pour  
dire qu'il correspond \`a une \fam de \caris  
de taille \,$t(n) = \O(f(n))$\, 
et de \prof \,$p(n) = \O(g(n))$.  
\end{nota}
  
Par exemple, l'\algo simplifi\'e du pivot de Gauss, tel  
qu'il a \'et\'e d\'evelopp\'e dans la section
\ref{subsec Gauss}, est $\SD(n^3,n)$.

Un \algo est dit \emph{optimal} lorsqu'il  
n'y a pas d'\algo asymptotiquement plus performant, 
du point de vue de la taille.  

Il y a des \pbs dont on conna\^{\i}t la \com  
\sqlez, \cad l'ordre asymptotique exact     
du nombre d'\oparis  \ncrs pour le  
r\'esoudre, comme par exemple le \pb de  
l'\eva d'un \pol \`a une \idtr  
sur un \acom quelconque\footnote{~L'\algo  
de Horner est optimal pour ce \pbz, \cf page \pageref{schemaHorner}.}.  
D'autres \pbsz, par contre, comme celui de la  
\mul des matrices, sont des \pbs dont on  
ignore la \com exacte \`a cause de l'\'ecart entre les  
bornes inf\'erieure et sup\'erieure asymptotiques que  
l'on conna\^{\i}t\footnote{~Le \pb de la  
\mul des matrices est \,$\Omega(n^2)$\, et  
\,$\O(n^{2,376})$.}.
  
\ss Il faut remarquer que le grand~$\O$  de la  
notation introduite ci-dessus pr\'esente
l'inconv\'enient majeur de \gui{cacher} la constante
asymptotique qui permet de le d\'efinir.  
Elle a pourtant une importance pratique consid\'erable puisque deux  
\algos permettant par exemple de r\'esoudre respectivement  
le m\^eme \pb avec \,$100n^3$\, et \,$10^9n^2$\,  
\oparis  sont tels que le second a une \com  
asymptotique nettement meilleure que le premier  
(il peut arriver qu'il soit aussi optimal) alors que le  
second, asymptotiquement moins performant, reste plus  
rapide tant que le nombre d'\ops \`a effectuer n'a pas  
atteint la borne astronomique de $10^{23}$.

\subsection{Complexit\'e binaire}   \index{complexit\'e!binaire} 
\label{subsecCobicaris}
On raconte que l'inventeur du jeu d'\'echec demanda comme  
r\'e\-com\-pen\-se un grain de bl\'e sur la premi\`ere case,  
deux sur la deuxi\`eme, quatre sur la troisi\`eme et ainsi  
de suite jusqu'\`a la soi\-xan\-te-quatri\`eme. Cela fait a  
priori un \cari de \prof 64. Mais pour calculer  
$ 2^{64}-1= 2^{2^6}-1= 18.446.744.073.709.551.615$ un \cari de taille 
(et de \profz)  
6+1 suffit:

\begin{samepage}
\sni \textsf{D\'ebut
\hsu $v_0:=x\quad \quad $\, (porte d'entr\'ee,  
on \'evaluera avec $x=2$)  
\hsu $v_1:=v_0\times v_0$  
\hsu $v_2:=v_1\times v_1$  
\hsu $v_3:=v_2\times v_2$  
\hsu $v_4:=v_3\times v_3$  
\hsu $v_5:=v_4\times v_4$  
\hsu $v_6:=v_5\times v_5$~~~ ($v_6=2^{2^6}$)  
\hsu $v_7:=v_6-1$  
\fin}\label{2p2n}
\end{samepage}
  
\ss De m\^{e}me, un \cari de taille 20 \'evalu\'e sur l'entr\'ee 2 
permet de calculer \,$2^{2^{20}}=2^{1.048.576}=6,7411\ldots  \, 
10^{315.652}$. 
Ceci montre clairement qu'il y a une diff\'erence  
consid\'erable entre la taille d'un \cir  et celle  
des objets qu'il peut produire lorsqu'on l'\'evalue sur  
des entiers cod\'es en binaire.  
  
La \emph{\com \bin d'un \cirz}  (ou d'une \fam de \cirsz)
\index{complexit\'e!binaire d'une \fam de \cirs}
\index{binaire!complexit\'e}
est par  d\'efinition la \com du calcul d'\eva  
qu'il produit lorsqu'on prend ses entr\'ees \emph{dans un anneau  
fix\'e avec un codage fix\'e}. L'exemple le plus simple  
et le plus important est l'anneau des entiers cod\'es en  
binaire.  
  
Naturellement, si on accepte de coder un entier par un  
\cari sans division ayant pour seules entr\'ees des  
constantes d\'etermi\-n\'ees a priori ($-1,0,1,2$ par  
exemple) et si on note \,$\ZZ_{\mathrm{preval}}$\, l'anneau  
des entiers ainsi cod\'e, on voit que l'\eva d'un  
\cari sans division dans \,$\ZZ_{\mathrm{preval}}$\, est en  
temps \lin (il suffit de mettre les \cirs bout  
\`a bout en changeant seulement certaines \profs et  
certains identificateurs). 
\indexnota{Zpreval@$\ZZ_{\mathrm{preval}}$} Le \pb avec  
\,$\ZZ_{\mathrm{preval}}$\, est alors report\'e du c\^ot\'e  
du test de signe, de la division euclidienne, ou de  
l'\eva des \cirs avec divisions exactes.  
  
Il est donc crucial de pr\'eciser \`a la fois l'anneau et  
le codage choisi pour cet anneau lorsqu'on veut parler de  
la \com \bin d'un \cariz.  

Signalons \`a ce sujet qu'en g\'eom\'etrie \agqz, 
la notion usuelle de degr\'e d'un \pol peut \^etre souvent 
remplac\'ee avantageusement par la notion de \prof d'un \preval 
\arith qui lui correspond. Il s'agit l\`a d'un sujet de 
recherche actif et prometteur (\cf \cite{Giu,Giu2}).

\subsubsection*{Un exemple: \com binaire de l'\apg}  

Elle est mesur\'ee par le nombre d'\ops \boles \ncrs pour  
ex\'ecuter l'\algo
avec des entr\'ees cod\'ees sous forme de suites de bits. 
Cette \com d\'epend de mani\`ere importante du corps 
\,$\K$\, et du codage choisi pour les \elts de \,$\K$. 

Si le corps \,$\K$\, est un corps fini, la \com 
\bin  est proportionnelle \`a la \comarz. 
C'est \gui{le bon cas} pour l'\algoz.

Appliqu\'e dans le cadre de calculs num\'eriques 
(ce qui constitue  aujourd'hui une partie importante 
du travail des ordinateurs), l'\algo est en \gnl
ex\'ecut\'e avec des nombres en virgule flottante, cod\'es par 
des suites de bits de
longueur fixe, et la \com  \bin est de nouveau 
proportionnelle \`a la \comarz. 
Mais naturellement, on ne travaille pas vraiment avec les
\elts du corps des r\'eels. 
D'o\`u la n\'ecessit\'e de garantir les r\'esultats
avec une pr\'ecision demand\'ee. L'analyse num\'erique 
matricielle remplit des rayons entiers de biblioth\`eques.

Dans cet ouvrage, nous ne prenons en compte que les calculs  
e\-xacts (en pr\'ecision infinie dit-on parfois), et nous ne 
ferons gu\`ere d'autre allusion aux
aspects proprement num\'eriques des \algos 
que nous commenterons (voir cependant page \pageref{caricunu}).

\ms  
La \mpg appliqu\'ee dans le corps des rationnels r\'eserve
quelques d\'esagr\'eables surprises. M\^eme si les entr\'ees 
sont des nombres entiers
(suppos\'es cod\'es en binaire de la mani\`ere usuelle), on doit 
\immt passer au corps des fractions. 
Un rationnel est alors cod\'e par un couple d'entiers, le
\nume avec un signe et le \deno strictement positif. 
Avec les rationnels
ainsi cod\'es (ce qui est le codage binaire naturel), on est 
alors devant l'alternative
suivante: simplifier les nouvelles entr\'ees de la matrice 
d\`es qu'elles sont
calcul\'ees, ou ne jamais simplifier. La deuxi\`eme solution est 
d\'esastreuse, car les
fractions successives voient en \gnl les tailles de leur \nume et 
\deno cro\^{\i}tre de mani\`ere exponentielle. 
La premi\`ere solution, quoique moins
d\'esastreuse, est n\'eanmoins co\^uteuse, car elle 
implique des calculs syst\'ematiques de pgcd. 
La derni\`ere formule donn\'ee dans la 
propri\'et\'e \ref{propri Gauss} permet d'exprimer \,$a_{\,ij}^{[p]}$\, 
comme quotient de deux
\deters extraits de la matrice de d\'epart (et elle se 
\gni au cas o\`u
des permutations de lignes ou de colonnes sont effectu\'ees). On 
a donc la garantie que
toutes les fractions qui sont calcul\'ees au cours de 
l'\algo restent de taille
raisonnable ($\O(n\,(t+\log\, n))$\, si on part d'une matrice 
\,$n\times n$\, \`a \coes  
entiers major\'es par \,$t$\, en  taille binaire, \ie 
major\'es par \,$2^t$\, en valeur absolue).
Le nombre d'\oparis dans \,$\ZZ$\, 
doit donc \^etre multipli\'e par un facteur \,$nt$\, 
pour tenir compte du calcul de simplification des 
fractions.  
La \com \binz, elle, a une majoration fort d\'ecevante 
en  $\,\O(n^5t^2)$\, (\`a des facteurs \logqs pr\`es)
si on utilise les \algos \usls pour la \mul ou la division de
deux entiers.

Appliqu\'ee avec le corps des fractions de \,$\ZZ[X]$\, ou 
\,$\ZZ[X,Y,Z]$\,  
la \mpg se heurte au m\^eme type de 
difficult\'es, mais  
tr\`es nettement aggrav\'ees, car les calculs de pgcd de 
\polsz, surtout en  
plusieurs variables, sont tr\`es co\^uteux.

\subsubsection*{Situations dans lesquelles la \com    
\bin d'un \cir est en rapport \'etroit avec sa \com \arith}

Nous signalerons ici trois situations de ce type.  
  
\ms  Le \emph{premier cas } 
est celui d'une \fam de  
\caris \'evalu\'es dans un anneau avec un codage pour  
lequel les \oparis  produisent des objets de taille  
bien contr\^ol\'ee, du fait m\^eme de la structure du  
\cariz.  

\begin{prop} \label{prop coalgo1}~  
Consid\'erons une \fam de \caris \,$\Gamma_n$\, de  
taille \,$\sigma_n$\, et de \prof \,$\pi_n$\, ($n$\,  
est un \parat contr\^olant le nombre d'entr\'ees  
du \cir  \,$\Gamma_n$.   
Supposons en outre que la production du \cir  \,$\Gamma_n$\, 
r\'eclame un temps \,$\tau_n$.   
Soit enfin \,$\A$\, un anneau donn\'e dans un codage pour lequel  
les \oparis  sont \etp \,$\O(N^k)$\, avec  
\,$k>1$\, et la taille \,$t(x)$\, des objets v\'erifie l'\ine  
\,$t(a\circ b)\le t(a)+t(b)$. \\  
Alors la production puis l'ex\'ecution de ce \cir   
r\'eclame, dans le mod\`ele MAD, un temps major\'e par  
\,$\tau_n+\sigma_n\cdot \O((2^{\pi_n}N)^k)$\, ($N>n$\,  
est la taille de la liste des entr\'ees).   
En particulier si \,$\sigma_n=\O(n^h)$, \,$\pi_n\le \ell\,\log\,n$\,  
et \,$\tau_n=\O(n^{c})$\, (pour des  
constantes convenables \,$h$, \,$\ell$\, et \,$c$) alors  
l'ex\'ecution de l'\algo correspondant \`a la \fam 
 \,$\Gamma_n$\, est (globalement) \etpz,
pr\'ecis\'ement en  
\,$\O(n^{c}+ n^{h+\ell k}N^k)$. 
\end{prop}  
  
\prv  
Dans le mod\`ele MAD, on peut utiliser un registre distinct  
pour chacune des variables du \prevz.  
La taille de tous les r\'esultats \itmds est  
major\'ee par \,$2^{\pi_n}N$\,  puisqu'elle double au 
maximum quand la \prof augmente d'une unit\'e.\\  
Les transferts entre les registres de travail et  
l'accumulateur repr\'esentent un temps de l'ordre de  
\,$\sigma_n\cdot (2^{\pi_n}N+\log(\sigma_n))$\, qui est  
n\'egligeable devant l'estimation du temps d'ex\'ecution des \oparis   
proprement dites: \,$\sigma_n\cdot\O((2^{\pi_n}N)^k)$. 
\qed

\begin{remark} 
\label{rem coalgo1}
\emph{Dans le mod\`ele des machines de Turing, on obtient  
les m\^emes majorations pour \,$n$\, fix\'e.  
Par contre, lorsque \,$n$\, varie, se pose le \pb  
de la gestion d'un nombre \emph{non fix\'e a priori} de  
variables de travail, alors qu'une telle machine n'a,  
quant \`a elle, qu'un nombre fix\'e a priori de bandes de  
travail.   
Les transferts de donn\'ees entre d'une part la  
bande o\`u est stock\'ee la liste des (contenus des)  
variables de travail et d'autre part les bandes o\`u sont  
ex\'ecut\'ees les \oparis  prennent normalement un  
temps de l'ordre de  
\,$(\sigma_n)^2 (2^{\pi_n}N+\log(\sigma_n))$\,  
car la bande de stockage doit \^etre relue pour chacune  
des \,$\sigma_n$\, \ops arithm\'eti\-ques, et sa taille est  
seulement major\'ee par  
\,$\sigma_n\cdot(2^{\pi_n}N+\log(\sigma_n))$.  
Il s'ensuit que la majoration en temps obtenue peut  
parfois \^etre un peu moins bonne que celle indiqu\'ee  
pour le mod\`ele MAD.  
}  
\end{remark}  
    
De nombreuses variantes de la situation pr\'ec\'edente  
peuvent \^etre utilis\'ees.  
Par exemple, pour l'\eva dans \,$\ZZ$, c'est  
seulement la \prof \muv qui doit \^etre en  
\,$\O(\log\,n)$\,  
pour qu'on ait un bon contr\^ole de la taille des objets  
produits, et donc de l'ensemble du calcul d'\evaz.  
    
Un \algo est dit \emph{bien parall\'elis\'e} lorsqu'il  
correspond \`a une \fam de \caris \,$(\Gamma_n)$\, dont  
la taille \,$\sigma_n$\, est  
optimale et dont  
la \prof est en \,$\O(\log^{\ell}(\sigma_n))$\, (pour un  
certain exposant \,$\ell >0$). Si la taille est \polle  
en \,$n$, la \prof est alors \pog\index{polylogarithmique},  
\cad en \,$\O(\log^{\ell}(n))$.  
En fait, nous utilisons dans cet ouvrage le terme \emph{bien 
\parasz}\index{bien parall\'elis\'e!algorithme}  
avec un sens un peu plus lib\'eral
pour le mot optimal. Pour les \algos \etp nous demandons seulement
que, en ce qui concerne la taille, l'exposant du  \,$n$\, ne soit pas 
tr\`es loin de celui du meilleur \algo \sql connu (la \prof 
\'etant, elle
\pogz). C'est en ce sens que nous consid\'erons que les \algos
de Csanky, de Chistov ou de Berkowitz sont bien \parassz.
 
\ms  
Le \emph{deuxi\`eme cas} 
est celui d'une \fam de \caris  
dont la \prof n'est pas \ncrt \logq  
et pour laquelle on a un argument de nature \agq  
qui permet de mieux majorer la taille des objets  
\itmds que l'argument de \profz.  
C'est par exemple le cas de l'\apg simplifi\'e  
(\'eventuellement modifi\'e par \elid \`a la Strassen) ou de l'\ajbz. 
M\^eme dans le cas d'un \algo  
bien \paras com\-me celui de Berkowitz,  
ex\'ecut\'e dans \,$\ZZ$, les majorations de taille ob\-te\-nues  
par un argument \agq direct sont meilleures que  
celles obtenues par l'argument de \profz.  
  
Signalons un calcul de majoration simple qui permet  
souvent un contr\^ole satisfaisant de la taille des objets  
\itmds dans le cas de l'\eva dans un  
anneau du style \,$\Mat_n(\ZZ[x,y])\,${(\footnote{~Ceci  
d\'esigne l'anneau des matrices \,$n\times n$\, \`a  
\coes dans $\ZZ[x,y]$.})} cod\'e en \rpn  
dense (voir la note \vref{labdense}), 
les entiers \'etant eux-m\^emes cod\'es en binaire.  
Si \,$A=(a_{ij})$\, est une matrice dans cet anneau,  
on note \,$d_A$\, le degr\'e maximum d'une entr\'ee  
\,$a_{ij}(x,y)$\, et \,$\ell_A:=\log(\Sigma_{ijhk}\vert  
a_{ijhk}\vert)$, o\`u \,$a_{ijhk}$\, est le \coe de  
\,$x^hy^k$\, dans \,$a_{ij}(x,y)$. On a alors la taille de  
\,$A$\, qui est major\'ee par \,$n^2d_A^2\ell_A$\, et les  
formules suivantes sont faciles \`a v\'erifier:  
$$
\begin{array}{ll} 
\ell_{A\pm B}\le 1+\max{(\ell_A,\ell_B)}& \ell_{AB}\le \ell_A+\ell_B    
\\ 
d_{A\pm B}\le 1+\max{(d_A,d_B)}& d_{AB}\le d_A+d_B     
\end{array}$$
Ceci signifie que ce type d'anneau se comporte comme  
\,$\ZZ$\, pour tous les calculs de majoration de taille des  
objets produits lors de l'\eva d'un \cariz. 
En particulier si la taille du \cir   
\num$n$\, est \polle en \,$n$\, et si sa \prof  
\muv est \logqz, alors la taille des  
objets est \polt major\'ee.  
La plupart des \algos que nous examinons dans cet ouvrage ont  
pour le type d'anneau que nous venons de signaler,  
une majoration \polle de la taille des objets  
\itmdsz. Signalons en revanche le mauvais comportement de l'\aghb  
pour la taille des objets \itmdsz. 
  
\ms  
Le \emph{troisi\`eme cas} 
est celui d'une \fam de  \caris (sans divisions) 
\'evalu\'es dans un cadre de calcul num\'erique \label{caricunu}
bien contr\^ol\'e.  
Lors de l'\eva du \cirz,
les entr\'ees sont des nombres dyadiques interpr\'et\'es  
comme des nombres r\'eels pris   
avec une pr\'ecision fix\'ee. Toutes les portes du  
\cir  sont elles-m\^emes \'evalu\'ees avec une pr\'ecision fix\'ee.  
Un calcul de majoration d'erreur est \ncr pour que  
le r\'esultat du calcul ait un sens math\'ematique pr\'ecis.  
Ce calcul dit une chose du genre suivant: sachant que vous  
d\'esirez les sorties avec une pr\'ecision absolue \,$p$\,  
(\cad de \,$p$\, digits apr\`es la virgule), et que les  
entr\'ees sont prises sur l'intervalle contr\^ol\'e  
par le \parat \,$n$, alors vous devez \'evaluer le  
\cir  \,$\Gamma_{n,p}$\, en effectuant tous les calculs  
\itmds avec la pr\'ecision \,$\vep(n,p)$\,  
(en particulier les entr\'ees doivent \^etre prises avec  
cette pr\'ecision). 
Par exemple, on pourra imaginer une \fam de \caris  
\'evaluant en ce sens la fonction \,$(x^2+1)/ln(1+x)$\,  
sur l'intervalle \,$]\,0,\infty\,[~$:  le \cari
 \,$\Gamma_{n,p}$\, doit permettre d'\'evaluer cette fonction
sur l'intervalle
\,$[2^{-n},2^n]$\, avec la pr\'ecision \,$p$, 
en ex\'ecutant tous
les calculs avec une pr\'ecision \,$\vep(n,p)$. 

Si la \fam peut \^etre produite \etpz, et si
la pr\'ecision requise  \,$\vep(n,p)$\, peut
\^etre major\'ee par un \pol en \,$(n,p)$\,  
alors la fonction r\'eelle ainsi calcul\'ee est dite  
calculable \etp (\cf \cite{Ko,Ho,KF,LLM}). 
Cela signifie que cette fonction  
peut \^etre \'evalu\'ee avec la pr\'ecision \,$p$\,  
sur n'importe quel r\'eel dans l'intervalle  
contr\^ol\'e par \,$n$\, en un temps qui d\'epend  
\polt de \,$ (n,p)$. Il s'agit donc d'analyse  
num\'erique enti\`erement s\^ure et parfaitement  
contr\^ol\'ee. 
 
Ce type d'\algos est en phase d'\^etre impl\'ement\'e sur  
machine, cela peut \^etre consid\'er\'e comme une des  
t\^aches importantes \`a r\'ealiser par le Calcul Formel.

\section[Familles uniformes de \cirsz]{Familles uniformes de \cirs 
\arits et \bols}       
\label{subsec uniforme}\label{sec unif}
 
Les \algos de calcul \agq usuels ont  
un nombre d'entr\'ees et de sorties qui d\'epend d'un ou  
plusieurs \parats entiers, comme par exemple  
\gui{la \mul de deux matrices} (3 \parats  
pour fixer les tailles des deux matrices)  
ou \gui{le produit d'une liste de matrices} (une liste  
d'entiers pour \paratsz) ou \gui{le \deter  
d'une matrice} (un \paratz). Nous avons appel\'e ces  
\parats des \emph{\parats d'en\-tr\'ee}.  
Comme nous l'avons d\'ej\`a dit,
ce n'est pas seulement la taille et  
la \prof du \cir  (en fonction des \parats  
d'entr\'ee) qui sont importantes, mais aussi son
co\^ut de production. Pour calculer le \deter  
d'une matrice \`a \coes entiers dans la situation  
la plus \gnle possible, par exemple, on doit  
d'abord produire le texte du \prev correspondant au \cir   
qu'on envisage, et ensuite ex\'ecuter ce \prev sur la liste  
d'entr\'ees voulue. Si le \cir  est de faible \prof et de  
faible taille mais que le co\^ut de la production du \prev  
correspondant cro\^{\i}t tr\`es vite lorsque le \parat  
d'entr\'ee augmente, on ne peut gu\`ere \^etre satisfait du  
r\'esultat.
  
C'est la raison pour laquelle on a introduit la notion de  
\emph{\famu de \carisz}\index{famille uniforme!de circuits 
arithm\'etiques}.  
On dit qu'une \fam de \caris (index\'ee par les \parats  
d'entr\'ee) est uniforme lorsque le co\^ut de production du  
\cir  (en tant que texte d'un \prevz)  
d\'epend \gui{de mani\`ere raisonnable} des \parats  
d'entr\'ee.   
Une premi\`ere notion d'uniformit\'e consiste \`a  
demander que le co\^ut de production du \cir  soit  
dans la classe \,$\P$, \cad \etpz.  
Une deuxi\`eme notion, plus forte,  consiste \`a  
demander qu'on soit dans la classe $\LOSP$ \cad  
que l'espace de travail \ncr \`a la production  
du \cir  soit \logqz.  
  
Ces notions d'uniformit\'e sont relativement satisfaisantes  
mais elles n\'eces\-si\-teraient d'\^etre mieux explicit\'ees
dans chaque cas concret.  
Il est clair qu'une \fam de \cirs d\'ependant d'un  
\parat d'entr\'ee \,$n$\, qui aurait une \prof en \,$\log\,{n}$,  
une taille en \,$n^2$\, et un co\^ut de production en \,$n^{2001}$\,  
ne serait pas un tr\`es bon cru pour l'ann\'ee 2001. Dans la  
litt\'erature sur le sujet r\`egne un silence discret. En fait  
tout le monde consid\`ere apparemment
qu'il est bien clair que le co\^ut de  
production du \cir n'a en \gnl pas un ordre de grandeur bien  
sup\'erieur \`a sa taille.  
  
Nous nous contenterons de confirmer cette impression par  
l'\'etude d'un cas d'\'ecole, la \mul rapide  
des matrices \`a la Strassen. 
Nous renvoyons pour cette \'etude au chapitre 
\ref{chap multimat}  section \ref{subsecStraUnif}
\thoz~\ref{thMuStra0}.

\subsubsection*{Classes de \com $\NC$}   

Pour d\'efinir les notions de taille et \prof en
 \comar \paral on a utilis\'e des \fams de
\caris sans exiger que ces \fams soient uniformes.
 
En \com \binz, les entr\'ees et les
sorties d'un \algo sont des mots \'ecrits sur un alphabet fix\'e,
par exemple l'alphabet \,$\{0,1\}$ (ou si on pr\'ef\`ere
des entiers \'ecrits en binaire). Il est alors naturel
d'utiliser les \fams de \cibos  
pour d\'efinir les notions de taille et \prof
d'un \algo \paralz. Dans un \ciboz, chaque entr\'ee est 
un \elt de  \,$\{0,1\}$, et les portes sont de
trois sortes: \,$\lor$, \,$\land$ (avec deux ant\'ec\'edents)
ou \,$\lnot$\, (\`a un seul ant\'ec\'edent). 
Pour chaque longueur de l'entr\'ee d'un \algo \paralz, 
cod\'ee comme une suite finie de \bolsz, le 
\cibo correspondant doit calculer la sortie,
cod\'ee de la m\^eme mani\`ere. 
Mais sans uniformit\'e de la \famz,
on aboutirait \`a des contre-sens intuitifs \'evidents,
puisque toute fonction \,$f$\, de \,$\NN$\, vers 
 \,$\{0,1\}$\, telle que \,$f(n)$\, ne d\'epend que 
de la longueur de \,$n$\, 
est r\'ealisable par une \fam 
non uniforme de  \cibos  de taille \,$n+1$\, 
et de \prof  \,$0\,$
(a vrai dire, l'entr\'ee du \cir \num$n$ 
ne sert \`a rien, et aucune \op \bole
n'est ex\'ecut\'ee). Or une telle fonction
peut ne pas \^etre calculable.

\ss Pour un entier naturel \,$k$\, donn\'e, on note  
\,$\NC^k$\,\indexnota{NCk@$\NC^k$!complexit\'e binaire}
la classe de toutes les fonctions   
qui peuvent \^etre calcul\'ees par une \famu
\index{famille uniforme!de circuits bool\'eens}  
de \cibos dans  \,$\SD(n^h,\log^k\,{n})$\, o\`u   
\,$h$\, est un entier positif (par hypoth\`ese,  
le \cir \,$C_n$\, a un nombre de portes d'entr\'ee
\polt reli\'e \`a \,$n$).  
L'uniformit\'e est prise ici au sens le plus  
fort que nous avons consid\'er\'e au d\'ebut de 
cette section. C'est la $\LOSP$   
uniformit\'e,  
c'est-\`a-dire, pour une \fam de \cirs  
\,$(C_n)_{n\in \N}$, l'existence d'une machine  
de Turing qui, pour l'entr\'ee \,$n$, donne en  
sortie le codage du \cir \,$C_n$\, en utilisant un  
espace m\'emoire en \,$\O(\log\,{n}).$  
  
On pose \,$\NC = \bigcup_{k\in\NN}\NC^k$\indexnota{NC@$\NC$}. Il s'agit
d'un acronyme pour Nick's Class du nom de Nicholas Peppinger  
qui a propos\'e cette classification des \algos \paralsz.  
  
Alors \,$\NC\subseteq \P$\, mais l'inclusion dans  
l'autre sens (\cad l'\egt des deux classes) est un  
\pb ouvert, et il est conjectur\'e que l'inclusion  
est stricte.  
  
\ss  On peut d\'efinir des notions analogues en \com \arith
(\cite{VonZur,Bur}). Il serait alors 
th\'eoriquement \ncr de
distinguer dans les notations
la classe \,$\NC$\,  au sens de la \com\arith  de celle
d\'efinie pr\'ec\'edemment. 
En outre, en  \com\arith 
on peut exiger ou ne pas exiger l'uniformit\'e de la \fam
de \cirsz, et on peut aussi vouloir indiquer sur
quel \acom on travaille. 

Dans le cadre de cet ouvrage, nous ne d\'esirons pas 
multiplier les notations et nous garderons la notation
\,$\NC^k\,$\indexnota{NCk@$\NC^k$!complexit\'e arithm\'etique}  
pour parler des 
\famus de \caris en  \,$\SD(n^h,\log^k\,{n})$, (o\`u \,$n$\, 
est la somme des \parats d'entr\'ee du \cir et  
\,$h$\, est un entier positif).
\emph{Nous demandons en outre que le degr\'e de tous les \pols 
\'evalu\'es aux noeuds du \cir soit major\'e par un \pol en \,$n$.} 
Enfin, nous prendrons l'uniformit\'e en un sens plus modeste:
la \fam des \cirs doit seulement \^etre construite \etpz.
La seule vraie preuve d'uniformit\'e 
que nous faisons est d'ailleurs
celle du \tho  \ref{thMuStra0}, et la construction
que nous donnons n'est pas $\LOSP$
(par contre, notre r\'esultat est plus pr\'ecis en ce qui concerne
 le temps de construction du \cariz).

La plupart des autres \algos d\'evelopp\'es dans cet ouvrage 
ont une preuve d'uniformit\'e plus simple, 
ou alors analogue \`a celle donn\'ee pour 
le \tho  \ref{thMuStra0}.

Dans le cas des \fams non \ncrt uniformes,
qui ont \'et\'e intensivemnt \'etudi\'ees par Valiant,
nous utiliserons les notations \,$\VNC^k$\,  et \,$\VNC$\,
en \comar et  \,$\BNC^k$\,  et \,$\BNC$\, en \com \bolez.
(voir chapitres \ref{chap DeterUniv} et  \ref{chap Perma}).  
\section{Machines parall\`eles \`a acc\`es direct}   
\label{sec modpar}
  
Nous pr\'esentons bri\`evement dans cette section 
quelques
mod\`eles de \gui{machines} susceptibles d'ex\'ecuter des 
\fams de \cirsz, \ariths ou \bolsz. 
Nous ne d\'evelopperons pas cependant les questions
de la programmation pour les machines \parals concr\`etes. 

Le principal objet de la conception d'\algos  
\parals est la r\'educ\-tion du temps de calcul  
permettant de r\'esoudre un \pb donn\'e  
moyennant un nombre suffisant mais  
raisonnable de processeurs.  
   
\subsection[Une id\'ealisation des calculs parall\`eles]
{Une id\'ealisation des calculs parall\`eles sur ordinateur} 
\label{secPRAM}

A d\'efaut de mod\`ele unique nous devons faire un choix. 
En \algq \sqle la Machine  
\`a Acc\`es Direct ou \gui{Random Access Machine}  
(RAM) est une abstraction de l'ordinateur \sql de Von Neumann. 
Nous consid\'erons ici  
le mod\`ele analogue en \algq \paralz,  
celui des \emph{machines \parals \`a acc\`es direct\/}   
(Parallel Random Access Machines) ou PRAM, qui constitue   
le mod\`ele \gui{standard} (\cf \cite{Cos,VZG,Kar}).

Une \emph{machine \paral \`a acc\`es direct\,} ou PRAM  
est une machine virtuelle (et un mod\`ele id\'eal abstrait)  
compos\'ee d'un nombre illimit\'e de processeurs partageant  
une m\'emoire commune, \emph{la m\'emoire globale\,},  
elle-m\^eme constitu\'ee d'un nombre illimit\'e de  
registres\footnote{~Le nombre de processeurs ainsi que le  
nombre de registres de m\'emoire partag\'ee sont habituellement  
fonctions de la taille du \pb \`a traiter.}, auxquels ils  
ont acc\`es pour y lire ou pour y \'ecrire des donn\'ees ou des  
r\'esultats de calcul.  
  
Chaque processeur a sa propre m\'emoire locale suppos\'ee  
\egmt de taille illimit\'ee, et inaccessible aux autres  
processeurs. Elle lui permet d'ex\'ecuter en une seule unit\'e  
de temps ou \emph{\'etape de calcul\,} la t\^ache, consid\'er\'ee  
comme \elrz, compos\'ee de la suite d'instructions  
suivantes:  
  
\noi -- chercher ses op\'erandes dans la m\'emoire globale;  
  
\noi -- effectuer l'une des \oparis \,$\{+,-,\times\}$\,  
(et \'eventuel\-lement la division quand elle est permise) 
sur ces op\'erandes;  
  
\noi -- \'ecrire le r\'esultat dans un registre de la m\'emoire  
commune (ou globale).  
  
\ss  
Faisant abstraction de tous les \pbs d'acc\`es \`a la  
m\'emoire glo\-ba\-le, de communication et d'interconnexion 
entre processeurs, \emph{une unit\'e de temps} ou \emph{\'etape  
de calcul \paralz} dans un tel mod\`ele abstrait  
correspond  
\`a l'ex\'ecution simultan\'ee de cette t\^ache par un certain  
nombre de processeurs, \emph{les processeurs actifs\,}, d'autres  
processeurs pouvant rester \emph{inactifs}.   
   
L'ex\'ecution des t\^aches par l'ensemble des processeurs actifs   
est synchronis\'ee: une \'etape d\'emarre d\`es que les   
op\'erandes sont disponibles, \cad au d\'emarrage du processus,  
quand chaque processeur sollicit\'e puise ses donn\'ees dans  
la m\'emoire globale, ou d\`es la fin d'une \'etape quand chaque  
processeur actif a livr\'e le r\'esultat de son calcul, le  
d\'eroulement de ce calcul \'etant li\'e aux contraintes de  
d\'ependance entre donn\'ees dans l'\algo consid\'er\'e.   
  
Il existe plusieurs variantes du mod\`ele PRAM selon le mode  
d'acc\`es \`a la m\'emoire globale, concurrent ou exclusif.   
  
Ce sera une PRAM-EREW\footnote{~EREW comme \gui{Exclusive Read,  
Exclusive Write}.} si la lecture ou l'\'ecriture dans un m\^eme  
registre n'est permise qu'\`a un seul processeur \`a la fois,  
une PRAM-CREW\footnote{~CREW comme \gui{Concurrent Read,  
Exclusive Write} etc.} si la lecture est concurrente et  
l'\'ecriture exclusive, une PRAM-ERCW si la lecture est exclusive  
et l'\'ecriture concurrente, et une PRAM-CRCW si la lecture et  
l'\'ecriture simultan\'ees dans un m\^eme registre de la m\'emoire  
globale sont permises pour plusieurs processeurs \`a la fois.  
Dans les deux derniers cas, il faut \'eviter que deux processeurs  
mettent simultan\'ement dans un m\^eme registre des r\'esultats  
dif\-f\'e\-rents, ce qui donne d'autres variantes 
 de machines PRAM selon le mode de gestion de la  
concurrence d'\'ecriture (mode prioritaire, arbitraire, etc.).  
  
M\^eme s'il existe une hi\'erarchie entre ces diff\'erentes  
variantes, de la \gui{moins puissante} (EREW) \`a la \gui{plus  
puissante} (CRCW prioritaire), ces mod\`eles PRAM sont en  
fait \'equivalents, pour la classe des \pbs qui nous  
int\'eressent, dans le sens o\`u ils se ram\`enent l'un \`a  
l'autre par des techniques de simulation  
(\cf \cite{Cos,VZG,Kar}).  
  
Nous utiliserons 
pour la description et  
l'analyse des \algos qui nous concernent, la variante  
{\bf PRAM-CREW} dont la conception est tr\`es proche de la  
notion de \cari ou de \prevz,  
puisqu'une PRAM-CREW peut \^etre repr\'esent\'ee par un  
\cari dans lequel les n{\oe}uds d'entr\'ee  
repr\'esentent les donn\'ees du \pbz, et chacun  
des autres n{\oe}uds (internes) repr\'esente aussi bien  
un processeur actif (et l'\op qu'il ex\'ecute)  
que le contenu d'un registre de la m\'emoire globale  
correspondant au r\'esultat de cette op\'eration. 
  
Enfin la \prof du \cari ou du \prev
 telle que nous l'avons d\'efinie  
pr\'ec\'edemment (section \ref{sec circprev}) correspond  
au nombre d'\'etapes du calcul \paralz.  
  
\subsection{PRAM-complexit\'e et Processeur-efficacit\'e}     
\label{secParcomp}  
Plusieurs \parats permettent de mesurer ce que  
nous appelle\-rons la PRAM-\com d'un \algo donn\'e.  
Ces \parats sont:  
  
\noi -- \emph{le temps \paralz\,} qui est \'egal au nombre  
d'\'etapes du calcul \paral et qui correspond au 
temps d'ex\'ecution de l'\algo \paralz; c'est aussi  
ce que l'on appelle \emph{la \com \paralz\,}  
ou la \emph{\profz} de l'\algoz;  
  
\noi -- \emph{le nombre de processeurs\,} \cad le nombre maximum  
de processeurs simultan\'ement actifs durant une \'etape  
quelconque du calcul, sachant qu'un processeur peut \^etre  
sollicit\'e durant une ou plusieurs \'etapes successives;  
  
\noi -- \emph{le temps \sqlz} de l'\algo \cad le nombre  
d'\oparis  qui interviennent dans le calcul ou, ce qui  
revient au m\^eme, le temps \paral si on ne disposait que  
d'un seul processeur, ou encore la somme des nombres de  
processeurs actifs durant toutes les \'etapes du calcul  
\paralz. C'est ce que l'on appelle aussi  
\emph{la taille} et parfois m\^eme \emph{la surface de calcul}  
\cite{Cos} ou \emph{la \com \sqlez} de l'\algoz;  
  
\noi -- \emph{le travail potentiel} ou \emph{la surface totale}  
de l'\algo qui est le produit du nombre de processeurs  
utilis\'es par le nombre d'\'etapes du calcul \paralz,  
\cad le temps \sql si tous les processeurs  
\'etaient actifs durant toutes les \'etapes du calcul.  
  
\ss On peut r\'esumer la parfaite analogie des \parats jusqu'ici  
d\'efinis entre PRAM-CREW, \cari et \prev
par le tableau suivant:

\begin{center}
\label{propram}
{\small
\begin{tabular}{|l|l|l|}
\hline PRAM-CREW & Programme d'Evaluation & Circuit Arithm\'etique \\  
\hline Temps \paral & Profondeur & Profondeur\\  
\hline Temps \sql & Longueur & Taille \\  
\hline Nombre de processeurs & Largeur & Largeur \\  
\hline  
\end{tabular}\\  }
  
\ss {\small {\bf Tableau \ref{propram}}}  
\end{center}  
  
\noindent  
Le nombre de processeurs dans une PRAM est l'\'equivalent  
de la largeur dans un \prevz,  
le temps \sql dans une PRAM est l'analogue de la  
longueur (ou la taille) d'un \prevz,  
et le temps \paral correspond \`a la \profz.  
  
\ss  
\emph{L'efficacit\'e\,} d'un \algo est alors d\'efinie  
comme le rapport entre le temps \sql et le travail  
potentiel de cet \algoz, ou encore le rapport entre  
surface de calcul et surface totale de l'\algo consid\'er\'e.  
  
\ss  
Pour revenir \`a l'exemple de l'\apg (voir page  
\pageref{F6}), la PRAM-CREW  
qui r\'ealise cet \algo peut \^etre repr\'esent\'ee  
par le tableau suivant (rectangle de gauche) dont les lignes  
corres\-pondent aux $7$ \'etapes successives du calcul  
et les colonnes aux processeurs (ceux marqu\'es d'une  
croix sont les processeurs actifs 
au cours d'une \'etape donn\'ee): 

\begin{minipage}[tbc]{2cm}\label{PRAM-Gauss} 
$$\hspace{2,8cm} 4~processeurs \hspace{2,9cm} 2~processeurs $$
$$\hspace*{2,7cm} \overbrace{\hspace{2,5cm}~} \hspace{3,1cm} 
\overbrace{\hspace{1,4cm}}~$$
$$
\begin{array}{c|c|c|c|c|c|c|c|}  
\cline{2-5} \cline{7-8} 
Etape~1\dans ~ & \times & \times & \, & \, & 
~~~~~ Etape~1\dans ~ & \times & \times \\  
\cline{2-5} \cline{7-8}   
Etape~2\dans ~ & \times & \times & \times & \times & 
~~~~~ Etape~2\dans ~ & \times & \times \\  
\cline{2-5} \cline{7-8}   
Etape~3\dans ~ & \times & \times & \times & \times & 
~~~~~ Etape~3\dans ~ & \times & \times \\  
\cline{2-5} \cline{7-8}   
Etape~4\dans ~ & \times & \times & \, & \, & 
~~~~~ Etape~4\dans ~ & \times & \times \\  
\cline{2-5} \cline{7-8}   
Etape~5\dans ~ & \times & \, & \, & \, & 
~~~~~ Etape~5\dans ~ & \times & \times \\  
\cline{2-5} \cline{7-8} 
Etape~6\dans ~ & \times & \, & \, & \, & 
~~~~~ Etape~6\dans ~ & \times & \times \\  
\cline{2-5} \cline{7-8}  
Etape~7\dans ~ & \times & \, & \, & \, & 
~~~~~ Etape~7\dans ~ & \times & \, \\  
\cline{2-5} \cline{7-8} 
-\cdot - ~~~~~ & \cdot & \cdot & \cdot & \cdot & 
~~~~~ Etape~8\dans ~ & \times & \, \\  
\cline{2-5} \cline{7-8}   
-\cdot - ~~~~~ & \cdot & \cdot & \cdot & \cdot & 
~~~~~ Etape~9\dans ~  & \times & \, \\  
\cline{2-5} \cline{7-8}   
\end{array} $$ 
\end{minipage} \\[7mm]

Le m\^eme \algo peut \^etre simul\'e par une PRAM \`a deux 
processeurs (rectangle de droite) au lieu de quatre, 
moyennant une augmentation du nombre d'\'etapes 
(\cad un \gui{ralentissement} des calculs) avec 9 \'etapes 
au lieu de 7. 

Pour chaque rectangle, la surface marqu\'ee repr\'esente  
la surface de calcul ou le temps \sqlz, la  
surface totale repr\'esentant le travail potentiel; 
la longueur et la largeur du rectangle repr\'esentent  
respectivement le temps \paral et le nombre de processeurs. 
L'efficacit\'e de cet \algo passe de \,$15/28$\, quand il 
est r\'ealis\'e par la PRAM initiale \`a \,$15/18$\, avec 
la PRAM modifi\'ee \cad de \,$54\,\%$ \`a \,$83\,\%$ environ. 

\medskip 
Nous introduisons maintenant la notation classique  
suivante pour la PRAM-complexit\'e qui sera utilis\'ee 
dans la suite.  

\begin{nota} 
On note \,$\PRAM(p(n),t(n))$\,  
la classe des \pbs de taille \,$n$\, r\'esolus  
par un \algo PRAM-CREW en \,$\O(t(n))$\, \'etapes,  
avec \,$\O(p(n))$\, processeurs. Tout \algo \,$P$\,  
qui, ex\'ecut\'e sur une telle machine, permet de  
r\'esoudre un \pb de cette classe, est lui-m\^eme  
consid\'er\'e, par abus de langage, comme appartenant  
\`a cette classe, et on dira que \,$P$\, est un \algo  
\,$\PRAM(p(n),t(n))$.  
\end{nota}  

\emph{La Processeur-efficacit\'e\,} d'un  \algo repr\'esent\'e 
par une PRAM-CREW est une notion rela\-tive \cite{KalPan,KalPan2,Sni} 
estim\'ee \`a partir du temps \sql d'un \algo   
choisi comme \algo de r\'ef\'erence: il s'agit en ce qui 
nous concerne, pour l'\agr \linz, de l'\algo de 
la multiplication des \macas d'ordre \,$n$\, 
suppos\'e \^etre r\'ealis\'e par une PRAM-CREW en 
\,$\log\,{n}$\, \'etapes, avec \,$\M(n)$\, processeurs. 
On peut \'evidemment supposer \,$\M(n)=\O(n^3)$\,
et \,$\M(n)=\Omega(n^2)$.

\begin{defi} 
Un \algo P est dit processeur-efficace (par rapport \`a 
un \algo de r\'ef\'erence de temps \sql \,$S(n)$) 
s'il existe $k,m\in \NN^*$\, tels que P soit dans 
\,$\PRAM(S(n)\log^m(n),\log^k(n))$. 
\end{defi}

Nous verrons plus loin des exemples d'\algos 
\gui{processeur-effi\-ca\-ces} (comme celui de l'inversion des 
matrices fortement r\'eguli\`eres, page \pageref{invtrian}) 
pour lesquels on prend comme \algo de r\'ef\'erence celui 
de la multiplication \usle (resp. rapide) des \macas \,$n\times n$\, 
r\'ealis\'ee par un \cari en 
\,$\SD(n^3,\log\,{n})$ (resp. \,$\SD(n^\alpha,\log\,{n})$\, 
pour $\alpha<3$).

\subsection{Le principe de Brent}  
Le principe de Brent affirme qu'on peut r\'epartir 
intelligemment le travail entre les diff\'erentes 
\'etapes d'un calcul \paralz, afin de diminuer de 
mani\`ere significative la proportion des processeurs 
inactifs (\cf \cite{Brent} lemme 2.4).
\begin{proposition}\label{brent}  
Un \algo \paral dont le temps \sql sur une PRAM 
est \'e\-gal \`a \,$s(n)$\, et dont le temps \paral est  
\'egal \`a \,$t(n)$\, peut \^etre simul\'e sur une PRAM  
utili\-sant \,$p$\, processeurs et 
\,$\einf{s(n)/p}   + t(n)$\, \'etapes de calcul
sans changer le temps \sqlz.  
\end{proposition}  
\preuve  
Supposons, en effet, qu'un calcul \paral peut \^etre  
effectu\'e en \,$t(n)$\, \'etapes \parals \`a raison  
de \,$m_i$\, \oparis  de base par \'etape.  
Si l'on impl\'emente directement ce calcul sur une  
PRAM pour \^etre ex\'ecut\'e en \,$t(n)$\, \'etapes,  
le nombre de processeurs utilis\'es sera alors \'egal \`a  
\,$m=\max{\{m_i\,|\,1\leq i \leq t(n)\}}$. En prenant  
\,$p$\, processeurs au lieu de \,$m$\, avec \,$p < m$\,  
(pour le cas \,$p\geq m$, la proposition est triviale)  
on peut ex\'ecuter le m\^eme calcul en faisant effectuer  
les \,$m_i$\, \ops de base de la \,$i\,$\eme \'etape par  
les \,$p$\, processeurs en \,$\esup{m_i/p}$\,  
\'etapes, et comme \,$\esup{ m_i/p } \leq \einf{  
m_i/p } + 1$\, le nombre total d'\'etapes avec une  
PRAM \`a \,$p$\, processeurs n'exc\`edera pas  
$$
\,\somm_{i=1}^{t(n)}(\einf{ m_i/p } + 1) \leq  
t(n)\,+\einf{ \somm_{i=1}^{t(n)} m_i/p\,} \leq  
t(n) + \einf{ s(n)/p } \,.
$$
\qed
  
\ms Ce principe est tr\`es utile lorsque le temps \paral  
\,$t(n)$\, est n\'egli\-gea\-ble (quand \,$n\tend \infty$)  
devant le temps \sql \,$s(n)$\, de l'\algo  
puisqu'on peut pratiquement diviser le nombre  
de processeurs par \,$t(n)$\, en doublant simplement  
le temps d'ex\'ecution \paral de l'\algoz:
on prend \,$p=\esup{ s(n)/t(n)}$.  
Par exemple, un \algo \,$\SD(n^\alpha,\log^k(n))$\,  
o\`u \,$\alpha$\, est un r\'eel positif et \,$k$\, un entier naturel  
quelconque, donne par application de ce principe de Brent un  
\algo \,$\PRAM(n^\omega/\log^k(n),\log^k(n))$.  
  
Cela permet dans la pratique, au prix d'un ralentissement relatif  
(\mul du temps de calcul par une petite constante),  
d'am\'eliorer l'efficacit\'e d'un \algo \paral en diminuant le  
temps d'inactivit\'e des processeurs par une r\'eduction du  
rapport entre le travail potentiel (\ie la surface totale)  
et le travail r\'eel (\ie la surface de calcul), et ceci par  
une r\'eorganisation des calculs dans le sens d'une meilleure  
r\'epartition des processeurs entre les \'etapes \paralsz.

Nous en d\'eduisons la propri\'et\'e suivante qui relie la 
complexit\'e des \caris \`a celle des PRAM.
  
\begin{proposition}\label{SD-PRAM}  
Un \algo \paral en \,$\SD(f(n),g(n))$\, est un \algo  
\,$\PRAM(f(n)/g(n),g(n))$. Inversement, tout \algo 
dans \,$\PRAM(p(n),t(n))$\, est un \algo en 
\,$\SD(p(n)t(n),t(n))$.
\end{proposition}   

\rem Dire qu'un \algo est processeur-efficace par rapport 
\`a un \algo de r\'ef\'erence de temps \sql \,$S(n)$\,
revient \`a dire qu'il est \,$\SD(S(n)\log^m(n),\log^k(n))$\, 
pour un couple \,$(m,k)\in \NN^{*}\times \NN^{*}$.

\newpage \thispagestyle{empty}  

\chapter{Diviser pour gagner}
\label{chap divpar}\label{sec divpar}
\minitoc
\acvide

\subsubsection*{Introduction}   

Dans ce chapitre, nous pr\'esentons une approche bien  
connue sous le nom de \gui{divide and conquer}  
que l'on peut traduire par \gui{diviser pour  
r\'egner} auquel nous pr\'ef\'erons le concept  
\dpg  
parce  que mieux adapt\'e, nous semble-t-il, au calcul  
\paralz.

Apr\`es en avoir donn\'e le principe \gnl
nous l'utilisons  pour \'etudier  deux \pbs classiques de  
l'\algq \paral que nous serons  
amen\'es \`a utiliser dans la suite:  
  
\begin{itemize}
\item  le calcul du produit de $n$ \elts  
d'un mono\"{\i}de;
\item  le \pb du calcul \paral  
des pr\'efixes  
(\gui{Parallel Prefix Al\-go\-rithm})
\end{itemize}

Pour ce dernier \pbz, nous d\'evelopperons,  
en plus de l'\algo classique, une  
\met r\'ecursive due \`a Ladner \& Fischer  
\cite{Lad} pour obtenir une famille  
de \cirs de taille \lin et de \prof  
\logqz. C'est le meilleur  
r\'esultat connu \`a l'heure actuelle. 

\ss Nous
appliquerons la strat\'egie \dpg  
en plusieurs autres occasions dans les chapitres suivants,
notamment pour les \muls rapides de matrices et 
de \pols et pour l'\agr \lin 
rapide sur les corps.

\section{Le principe g\'en\'eral}  
\label{subsec principegene}

L'approche \dpg    
s'applique pour r\'esoudre une famille de probl\`emes
\,$(P_n)_{n\in\N}$. Elle consiste \`a \gui{diviser} le \pb  
num\'ero \,$n$\, en  
$q~(q\geq 2)$ sous-\pbs du style \,$P_m$\, avec
\,$m<n$,  auxquels on peut appliquer,  
en \paral et de mani\`ere r\'ecursive,  
le m\^eme \algo que celui qui permet  
de r\'esoudre le \pb initial, pour  
r\'ecup\'erer ensuite le r\'esultat final  
\`a partir des solutions des sous-\pbsz.  
  
Le param\`etre entier \,$q$ repr\'esente le  
nombre des sous-\pbs qui seront trait\'es  
en \paralz. Lorsqu'il  
 ne d\'epend pas de \,$n$, il s'appelle le \emph{degr\'e de  
\parasmz} de l'\algo ainsi obtenu. \index{degr\'e de  
\parasmz} 
  
Une telle approche r\'ecursive de conception  
d'\algos permet souvent d'apporter  
une solution efficace \`a un \pb dans  
lequel les \,$q$\, sous-\pbs  \,$P_m$\, sont  
des copies r\'eduites du \pb initial,  
et avec \,$m$\, sensiblement \'egal  
(\`a $\lceil n/p \rceil$ par exemple, o\`u  
\,$p$\, est un entier donn\'e~$\,\geq 2$).  
  
Cette \met nous permet \egmt  
d'analyser la complexit\'e de l'\algo  
qu'elle produit et de calculer des majorants  
asymptotiques de la taille et de  
la \prof du \cir \arith correspondant,  
avec une estimation pr\'ecise  
de la constante cach\'ee du \gui{grand $\O$}.  
  
En effet, supposons que le \pb \`a  
traiter est le \pb  \,\num$ n = m_0\,p^\nu$\,  
($m_0,\,p,\,\nu\in \NN^*$) et qu'il peut  
\^etre scind\'e en \,$q$\, sous-\pbs  
 \,$P_m$\, avec \,$m=m_0\,p^{\nu -1}$,  
suceptibles d'\^etre trait\'es en \paralz.  
Remarquons tout de suite que \,$q$\, est un  
entier \,$\geq 2$\, d\'ependant  
\'eventuellement de \,$\nu$: c'est pourquoi  
on \'ecrira, dans le cas \gnlz,  
\,$q = q(\nu)$.  
  
Le co\^ut \,${\hat{\kappa}}(\nu) = (\tau(\nu),\,\pi (\nu))$\, 
de cet \algo o\`u  \,$\tau(\nu)$\,  
(resp. \,$\pi(\nu)$) d\'esigne la taille (resp. la \profz) du \cir  
correspondant, se calcule par \recu  
sur \,$\nu$\, \`a l'aide des formules suivantes:  
\begin{equation}\label{paral1}  
\formule{  
\,\tau(\nu ) &=& q(\nu)\,\tau(\nu-1) \,+ \,
         \tau'(\nu) \\[1mm]  
\,\pi(\nu ) &=& \pi(\nu-1)\, +\, \pi'(\nu)  
}  
\end{equation}

\sni  
o\`u \,$\tau'(\nu)$\, (resp.  
\,$\pi'(\nu)$) repr\'esente la  
taille (resp. la \profz)  
des \cirs correspondant \`a la double op\'eration  
de partitionnement du \pb  
et de r\'ecup\'eration de sa solution \`a  
partir des solutions partielles.  
  
L'absence du facteur \,$q$\, dans l'\egt  
exprimant la \prof  
\,$\pi$\, est due au fait que les  
\,$q$\, sous-\pbsz, de m\^eme  
taille, sont trait\'es \emph{en \paralz}  
avec des \cirs de \prof maximum  
\,$\pi(\nu-1)$.  
  
Si l'on se donne \,$\tau(0)$\,  
et \,$\pi(0)$\, le \sys (\ref{paral1})  
ci-dessus admet pour solution:  
\begin{equation}\label{paral2}  
\formule{  
\,\tau(\nu) &=&  
q(1)\,q(2)\,\cdots\,q(\nu)\,  
\tau(0) +  
\somm_{i = 1}^\nu\,
\left[\prod\nolimits_{j = i+1}^\nu q(j)\right]\, 
\tau'(i) \\ 
\,\pi(\nu) &= & 
\pi(0) + \somm_{i = 1}^\nu  
\,\pi'(i)  
}
\end{equation}  

Dans le cas o\`u \,$q = q(\nu)$\, est une  
constante, sachant que la \profz, ne  
d\'ependant pas de \,$q$, reste la m\^eme,  
le \sys  \,(\ref{paral2})\, devient  \,(\ref{paral3})\, ci-dessous.
Nous rappelons pr\'ecis\'ement les hypoth\`eses dans l'\'enonc\'e qui 
suit.  
\begin{proposition} 
\label{propTPDQ} 
Soient $m_0,\,p,\,q\in \NN^*$ fix\'es et \,$\nu\in \NN^*$ variable.
Nous supposons que le \pb \`a  
traiter est le \pb \,$P_n$\, avec \,$n = m_0\,p^\nu$\, et qu'il peut  
\^etre scind\'e en \,$q$\, sous-\pbs  
de type \,$P_m$\, avec   \,$m=m_0\,p^{\nu -1}$,  
suceptibles d'\^etre trait\'es en \paralz.
Nous notons  \,$\tau'(\nu)$\, (resp.  
\,$\pi'(\nu)$) la taille (resp. la \profz)  
des \cirs correspondant \`a la double op\'eration  
de partitionnement du \pb et de r\'ecup\'eration de sa solution \`a  
partir des solutions partielles.
Enfin \,$\tau_0$\, et \,$\pi_0$\, sont la taille et la \prof 
d'un \cir qui traite le \pb  \,$P_{m_0}$.
Alors la taille et \prof d'un \cir produit 
en utilisant la \met \dpg  sont:  
\begin{equation}\label{paral3}  
\formule{  
\,\tau(\nu) &=& q^\nu\, 
\tau_0 +  
\sum_{i = 1}^\nu\,q^{\nu - i}\,  
\tau'(i) \\ [2mm]  
\,\pi(\nu) &=&  
\pi_0 + \sum_{i = 1}^\nu  
\,\pi'(i)  
}
\end{equation}
\end{proposition}

En particulier si \,$\tau'(\nu)=\O(n^r)$\, avec
\,$r\neq \log{q}$\, et  \,$\pi'(\nu)=\O(\nu^\ell)$\, on obtient: 
\begin{equation}\label{paral4}  
\formule{  
\,\tau(\nu) &=&  \O(q^\nu)\;= \;\O(n^{\sup(r,\log{q})}) \\ [1mm]  
\,\pi(\nu)  &=&  \O(\nu^{\ell+1})\;=\;\O(\log^{\ell+1}{n})  
}
\end{equation}

\smallskip Donnons un aper\c{c}u rapide sur quelques cas particuliers 
significatifs que nous allons traiter dans la suite.

Dans le calcul \paral des pr\'efixes section \ref{subsec CPP}, 
nous avons de mani\`ere naturelle
\,$p=q=2$,  \,$r=1$\, et \,$\ell=0$\, ce qui conduit \`a une famille de
\cirs en \,$\SD(n\,\log{n},\log{n})$, et nous verrons qu'on peut encore 
tr\`es l\'eg\`erement am\'eliorer la borne sur la taille.

Dans la \mul des \pols \`a la Karatsuba section \ref{secKarat}, nous 
avons \,$p=2$, \,$q=3$,
\,$r=1$\, et \,$\ell=0$\, ce qui conduit \`a une famille de
\cirs en \,$\SD(n^{\log{3}},\log{n})$. 

Dans la \mul rapide des matrices \`a la Strassen section 
\ref{sec strass}, nous avons \,$p=2$, \,$q=7$,
\,$r=2$\, et \,$\ell=0$\, ce qui conduit \`a une famille de
\cirs en \,$\SD(n^{\log{7}},\log{n})$.

Enfin pour l'inversion des \matgs section \ref{invtrian},
 nous avons \,$p=2$, \,$q=2$,
\,$r=\alpha $\, et \,$\ell=1$\, ce qui conduit \`a une famille de
\cirs en \,$\SD(n^{\alpha},\log^2{n})$.  
\section{Circuit binaire \'equilibr\'e}  
\label{subsec CBI}

L'approche \dpgz,  
appliqu\'ee \`a ce premier \pbz,   
nous donne la construction d'un type
particulier de \caris de taille   
\lin et de \prof \,$\esup{\log {n}}$\,    
que l'on appelle les \cbes  
(\gui{Balanced Binary Trees}). 

Un \cbe est un \cari prenant en entr\'ee  
une liste \,$(x_1,x_2,\ldots,x_{n-1},x_n)$\, de  
\,$n$\, \elts d'un mono\"{\i}de \,${\cal M}$\,  
(loi associative not\'ee \,$*$\,  avec \elt neutre  
not\'e 1) et donnant en sortie le produit  
\,$\Pi = x_1*x_2*\cdots *x_{n-1}*x_n$.  
  
On peut supposer \,$n=2^\nu$\, o\`u  
\,$\nu\in \NN^*$\, quitte \`a compl\'eter la liste  
donn\'ee par \,$2^{\esup{\log  n}}-n$\,  
\elts \'egaux \`a 1, ce qui ne change pas  
le r\'esultat.  
  
Le \cir est d\'efini de mani\`ere r\'ecursive en divisant le \pb en  
deux sous-\pbs de taille  
\,$2^{\nu -1}$, qui correspondent \`a deux  
\gui{sous-\cirsz}  acceptant chacun en entr\'ee une  
liste de taille moiti\'e.  
  
Ces deux sous-\cirs calculent respectivement  
et en \paral les deux produits partiels  
\,$\Pi_1 = x_1*\cdots *x_{2^{\nu-1}}$\, et 
\,$\Pi_2 = x_{2^{\nu-1}+1}*\cdots *x_{2^\nu}$.  
On r\'ecup\`ere ensuite le produit \,$\Pi$\, en  
multipliant ces deux produits partiels.  
  
Ainsi un \cbe pour une entr\'ee de taille  
\,$2^\nu$\, est d\'efini par  
\recu sur $\nu$:  
  pour \,$\nu=0$\, c'est le \cir trivial  
\,${\cal C}_0$\, de taille \prof nulles.
Pour \,$\nu\geq 1$, le \cir  
\,${\cal C}_\nu$\, prend en entr\'ee une  
liste de longueur \,$2^\nu$, fait agir  
deux copies du \cir \,${\cal C}_{\nu-1}$\,  
pour calculer \,$\Pi_1$\,  et   
\,$\Pi_2$\, qu'il utilise  pour r\'ecup\'erer  
le r\'esultat final \,$\Pi=\Pi_1*\Pi_2$\,  
(comme l'indique la figure \vref{F4}).  
\begin{figure}[htbp] 
\begin{center}
\includegraphics*[width=10cm]{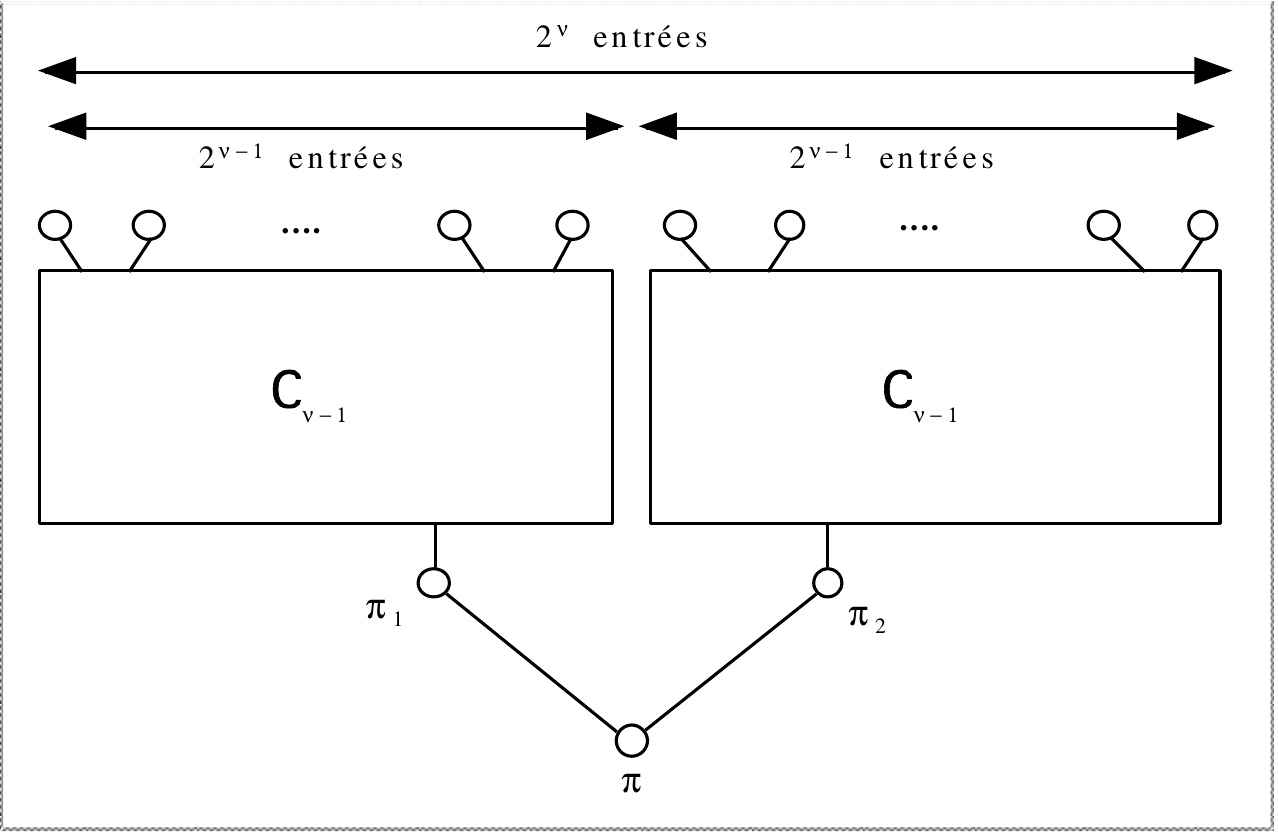}
\end{center}
\caption[Construction r\'ecursive des \cbes]{\label{F4}  
Construction r\'ecursive du \cbe \,${\cal C}_\nu\,$}  

\centerline{\hspace{1,5cm} (\ˆ partir du \cbe \,${\cal C}_{\nu-1}$)}  
\end{figure}
Si l'on note \,$\tau(\nu)$ et  
 \,$\pi(\nu)$\, la taille  
et la \prof du \cir  
\,${\cal C}_\nu$, on obtient les relations:
$$
\left\{\begin{array}{lcll} 
\tau(\nu) & = &  
2 \tau(\nu-1) \,+\, 1 & \mbox{  avec  }  
\tau(0) = 0 \\[1mm]  
\pi(\nu) & = & \pi(\nu-1) \,+\,1 & \mbox{  avec  }  
\pi(0) = 0 \end{array}\right.
$$  
qui admet la solution exacte:
\begin{equation} \label{eqCBIE}
\formule{  
\tau(\nu) &=&  
2^\nu-1 \\[1mm] \pi\,(\nu) &=& \nu\,.}
\end{equation}
\begin{proposition} 
\label{propCBIE} 
Un \cbe  qui prend en entr\'ee une liste  
une liste \,$(x_1,x_2,\ldots,x_{n-1},x_n)$\, dans un mono\"{\i}de 
\,${\cal M}$\,  et donne en sortie le produit  
\,$\Pi = x_1*x_2*\cdots *x_{n-1}*x_n$\, est un \cari   
de \prof \,$\esup{\log  n }$.   
Il est de taille \,$n-1$\, si \,$n$\, est une  
puissance de $2$, et cette taille est en tous cas  
major\'ee par \,$2n-3$\, lorsque  
\,$n$\, n'est pas une puissance de $2$. 
\end{proposition}
Notons qu'on peut trouver une  
majoration l\'eg\`erement meilleure de la taille pour \,$n>3$.

\section{Calcul \paral des pr\'efixes}     
\label{subsec CPP}
  
\'Etant donn\'ee une liste de \,$n$\, \elts  
\,$x_1,x_2,\ldots,x_n$\, (ou \,$n\,$-uplet)  
d'un mono\"{\i}de \,$({\cal M},*,1)$\, dont la loi  
(en \gnl non commutative) est not\'ee  
\muvt et dont l'\elt neutre est  
not\'e $1$, le \pb du calcul des pr\'efixes  
consiste \`a calculer les produits partiels
$$
\Pi_k=\prod\nolimits_{i=1}^kx_i~~~~\mathrm{pour}~~~(1\leq k\leq n).
$$  

La solution  na\"{\i}ve de ce \pb donne un \cir  
de taille $n-1$ (c'est la taille minimum) et de  
\prof $n-1$.

\subsubsection*{Premi\`ere \met de \paran}
Il est facile de voir que ce calcul peut \^{e}tre 
\paras pour obtenir
un circuit de \prof \,$\esup{\log{n}}$.

On peut toujours supposer \,$n=2^\nu$\, o\`u  
\,$\nu = \esup{\log {n}} \in \NN^*$, quitte \`a  
compl\'eter la liste donn\'ee par  
\,$2^{\esup{\log {n}}} - n$\, copies de  
l'\elt neutre $1$.  
  
Ce \pb peut se d\'ecomposer en deux sous-\pbs de taille  
\,$n/2 = 2^{\nu\!-1}$ qui seront trait\'es en \paralz:  
  
\noi $\bullet$\, le calcul des pr\'efixes pour la liste  
\,$\,x_1,x_2,\ldots,x_{2^{\nu\!-1}}$; 
  
\noi $\bullet$\, le calcul des pr\'efixes pour la liste  
\,$\,y_1,y_2,\ldots,y_{2^{\nu\!-1}}$\, o\`u  
\,$y_i = x_{2^{\nu\!-1} + i}$\,  
pour \,$1\leq i\leq 2^{\nu\!-1}$.  
  
 La solution du \pb principal est ensuite  
obtenue par \mul du produit  
\,$\Pi_{2^{\nu\!-1}}$, faisant partie de la solution  
du premier sous-\pbz, par les \,$2^{\nu\!-1}$\,  
produits partiels des \,$y_i$\,  
qui constituent la solution du second sous-\pbz.  
Cette derni\`ere \'etape de r\'ecup\'eration  
augmente par cons\'equent de  
\,$2^{\nu\!-1}$\, \muls la taille du \cir  
et de 1 sa \profz.  
  
\ss Pour le cas \,$n = 7$\, par exemple  
(on prend \,$n = 8$\, pour avoir une  
puissance de 2 et on fait \,$x_8=1$),  
on obtient le circuit \vref{F5} qui montre 
le d\'eroulement de cette \pcd pour  
le calcul des sept (ou huit) produits  
\,$\Pi_1 = x_1\,,~\Pi_2 = x_1*x_2\,,~\Pi_3 =  
x_1*x_2*x_3\,,~\ldots~,~\Pi_7 =  
x_1*x_2*\cdots *x_7$\,  
($\Pi_8 = \Pi_7$\, puisque \,$x_8=1$).
\setcounter{agh}{\arabic{algorithm}}
\stepcounter{algorithm}
\begin{agh}
\begin{center}
\includegraphics*[width=11cm]{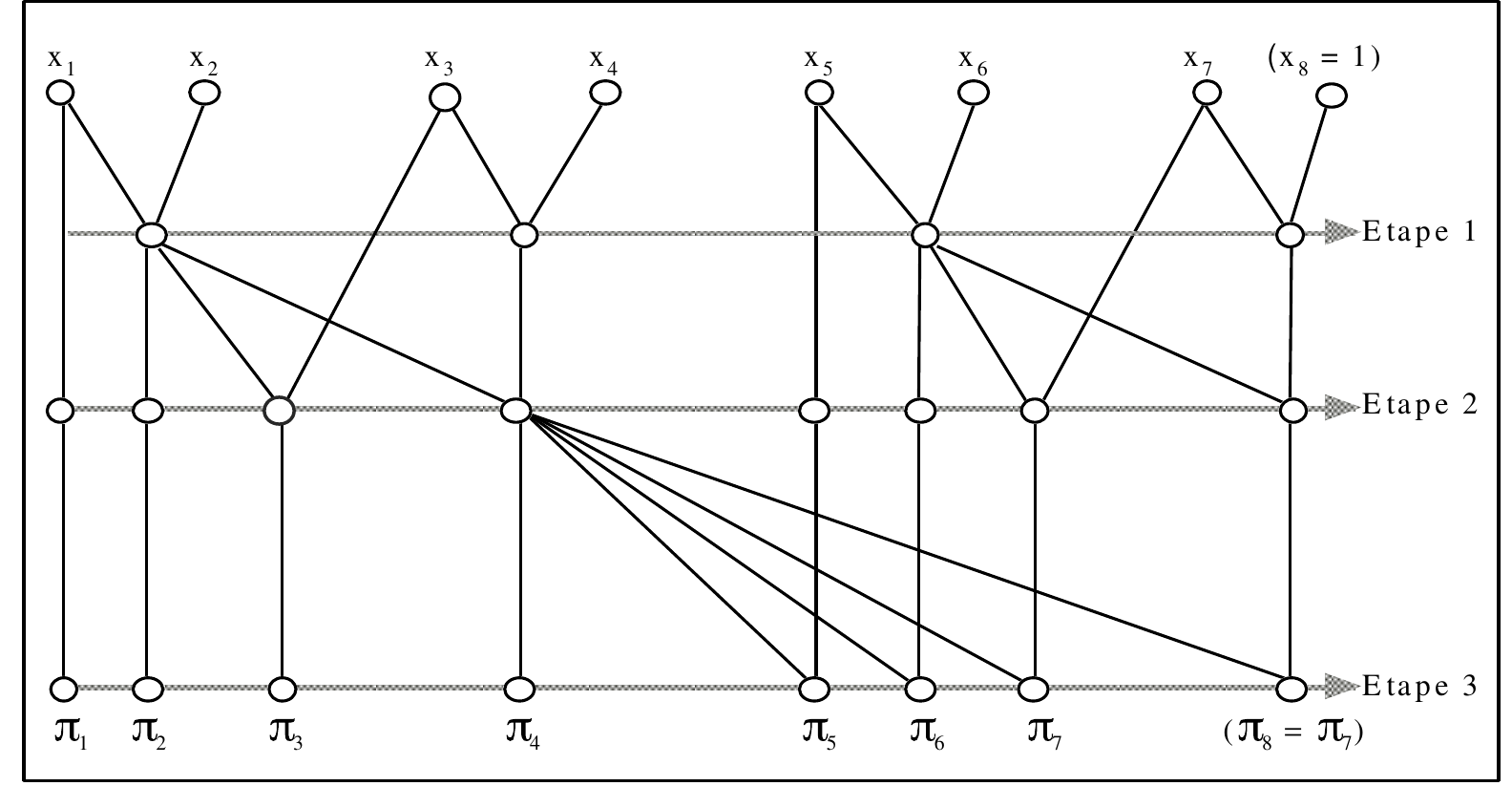}
\end{center}
\caption[Calcul Parall\`ele des Pr\'efixes ]{\label{F5}  
~~Calcul Parall\`ele des Pr\'efixes pour \,$n=7\,$} 
\end{agh}  
Appliqu\'ees \`a notre \pbz, les relations  
\,(\ref{paral3})\, donnent la taille et la \prof  
du \cir \arith correspondant au calcul  
\paral des pr\'efixes pour une liste donn\'ee de taille  
\,$2^\nu$. 

Il suffit en effet de faire \,$p=q=2$, 
\,$\tau'\,(i) = 2^{i-1}$, 
\,$\pi'\,(i) = 1$\, (pour   
\,$i\geq 1$) et 
\,$\tau(0) =  
\pi(0) = 0$\,  pour obtenir:  
$$\tau(\nu) = \sum_{i=1}^\nu 2^{\nu-i}2^{i-1} =  
\nu\,2^{\nu\!-1} = {n\over 2}\,\log {n},  
~~\mathrm{et}~~\pi(\nu)  
= \sum_{i=1}^\nu 1 = \nu = \log {n}.$$   
Ainsi le \pb du calcul des pr\'efixes pour une  
liste de \,$n$\, \elts 
se parall\'elise bien, et il admet  une  
solution en \,$\SD(n\,\log {n},\log {n})$\, ou encore,  
en utilisant le principe de Brent (proposition \ref{brent}), 
une solution qui est  
\,$\PRAM(n,\log {n})$.  
  
\ss  
Ladner \& Fischer \cite{Lad}  
obtiennent un meilleur r\'esultat en donnant une  
construction r\'ecursive d'un \cir en \,$\SD(n, \log {n})$. 
C'est ce que nous allons d\'evelopper au paragraphe suivant.  

\subsubsection*{Am\'elioration du calcul des pr\'efixes
(Ladner \& Fischer)}  
  
\'Etant donn\'es un  mono\"{\i}de \,$({\cal M},*,1)$,  
un entier \,$n\geq 2$, et \,$x_1,\ldots,x_n$\,  
dans \,$\cal M$, nous allons construire,  
\`a l'instar de Ladner \& Fischer \cite{Lad}  
deux familles de \cirs \,$({\cal P}_k\,(n))_{\,n\in \NN^*}$\, de  
tailles \,$S_k\,(n)$\, ($k\in\{0,1\}$) major\'ees  
respectivement par \,$4n$\, et \,$3n$\, et de  
\profs respectives  
\,$D_0\,(n)=\esup{\log {n}}$\,  
et \,$D_1\,(n)=\esup{\log {n}} + 1 $\,  
qui calculent les pr\'efixes  
\,$\Pi_1,\Pi_2,\dots,\Pi_{n}$\, du $n$-uplet  
\,$(x_1,x_2,\ldots ,x_n)$.  
  
Cette construction se fait conjointement et de  
mani\`ere r\'ecursive \`a partir du \cir trivial  
\,${\cal P}_0\,(1)={\cal P}_1\,(1)$\, r\'eduit \`a  
une seule porte (la porte d'entr\'ee).  
La figure \vref{F7}  montre le  
d\'eroulement de cette construction r\'ecursive  
conjointe des deux familles  
\,$({\cal P}_0\,(n))_{\,n\in \NN^*}$ et  
\,$({\cal P}_1\,(n))_{\,n\in \NN^*}$.
\begin{figure}[ht]  
\begin{center}
\includegraphics*[width=6.5cm]{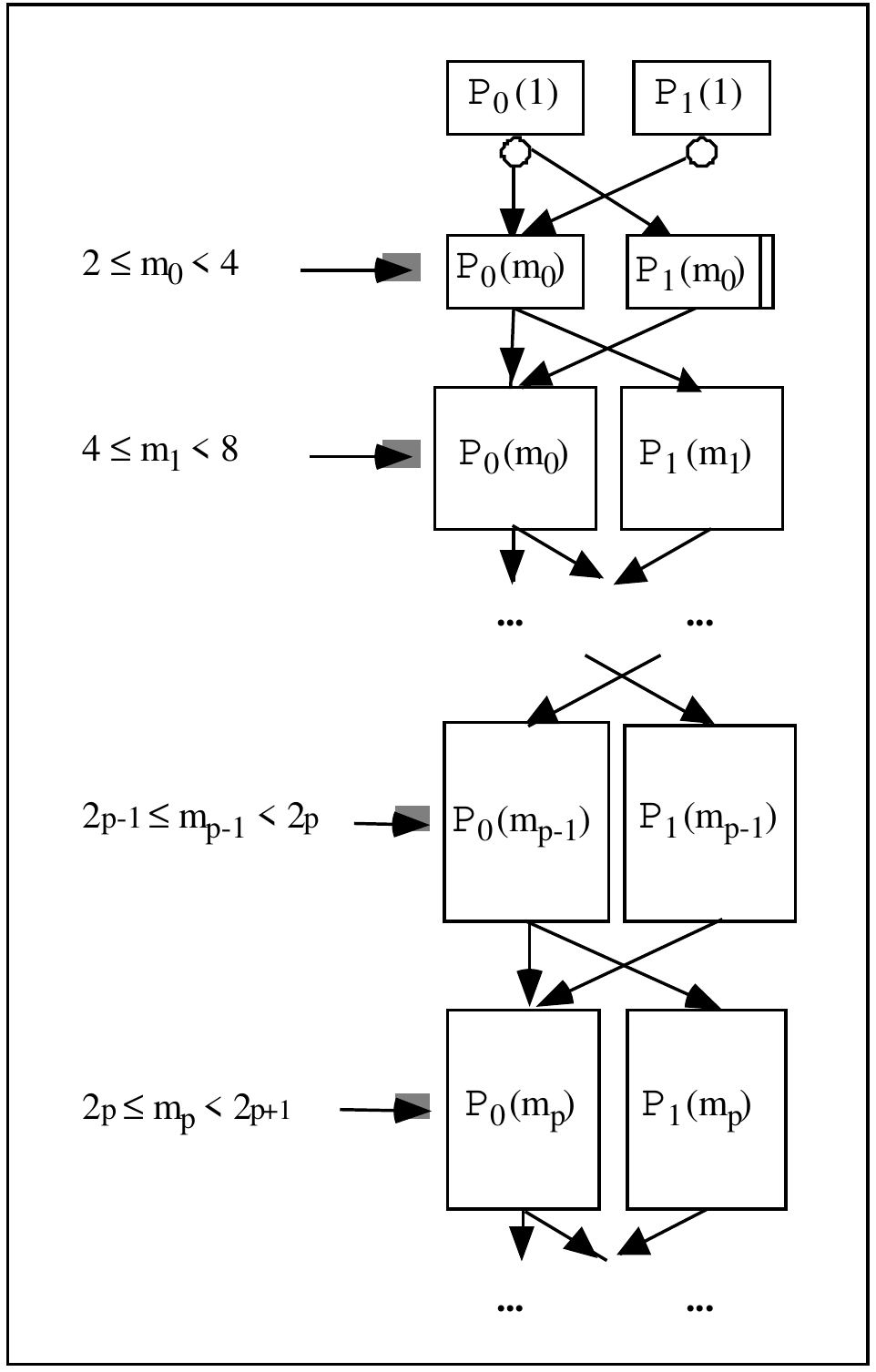}
\end{center}
\caption[Sch\'ema de la construction r\'ecursive  
des \cirs \,${\cal P}_k\,(n)$] {\label{F7} 
Sch\'ema de la construction r\'ecursive des \cirs 
$\,({\cal P}_k\,(n))_{\,n\in \NN}$}  
    
\centerline{$(k\in\{0,1\},\,(m_p = 2\,m_{p-1}\,  
\mbox{ pour } ~1\leq p\leq \esup{\log {n}})$}  
\end{figure}

\paragraph{Construction de la famille \,$({\cal P}_0\,(n))_{\,n\in 
\NN^*}$}~~

\sni On d\'efinit r\'ecursivement le \cir  
\,${\cal P}_0\,(n)$\, \`a partir des \cirs  
\,${\cal P}_{1}\,(\lfloor {n \over 2} \rfloor)$\,  
et \,${\cal P}_{0}\,(\lceil {n \over 2} \rceil)$\,   
appliqu\'es respectivement aux entr\'ees  
\,$\;\;  (x_1, \ldots , x_{\lfloor {n \over 2} \rfloor})\quad $\,
et
\,$ (x_{\lfloor {n \over 2} \rfloor + 1}, \ldots , x_n)$\,  
qui forment une partition de la liste donn\'ee  
\,$(x_1,\ldots ,x_n)$. 

Comme \,${\cal P}_{1}\,(\lfloor {n \over 2} \rfloor)$\,  
calcule  
\,$\Pi_1,\Pi_2,\dots,\Pi_{\lfloor {n \over 2}\rfloor}$,  
il suffit d'effectuer en \paral et en une seule 
\'etape les \,$\lceil {n \over 2} \rceil$\, \muls  
de \,$\Pi_{\lfloor {n \over 2} \rfloor}$\, par les  
\,$\lceil {n \over 2} \rceil$\, sorties de
\,${\cal P}_{0}\,(\lceil {n \over 2} \rceil)$\,  
pour avoir les pr\'efixes  \,$\Pi_{\lfloor  
{n \over 2} \rfloor + 1},\Pi_{\lfloor {n \over 2} \rfloor + 2},  
\ldots , \Pi_n$\, et par cons\'equent tous les pr\'efixes  
\,$\Pi_1,\Pi_2,\dots,\Pi_{n}$\, de la liste  
\,$(x_1,\ldots ,x_n)$.  
         
Partant du \cir trivial \,${\cal P}_0\,(1)=  
{\cal P}_1\,(1)$, la figure \vref{F8}  
 illustre cette construction.  
\begin{figure}[ht]  
\begin{center}
\includegraphics*[width=11cm]{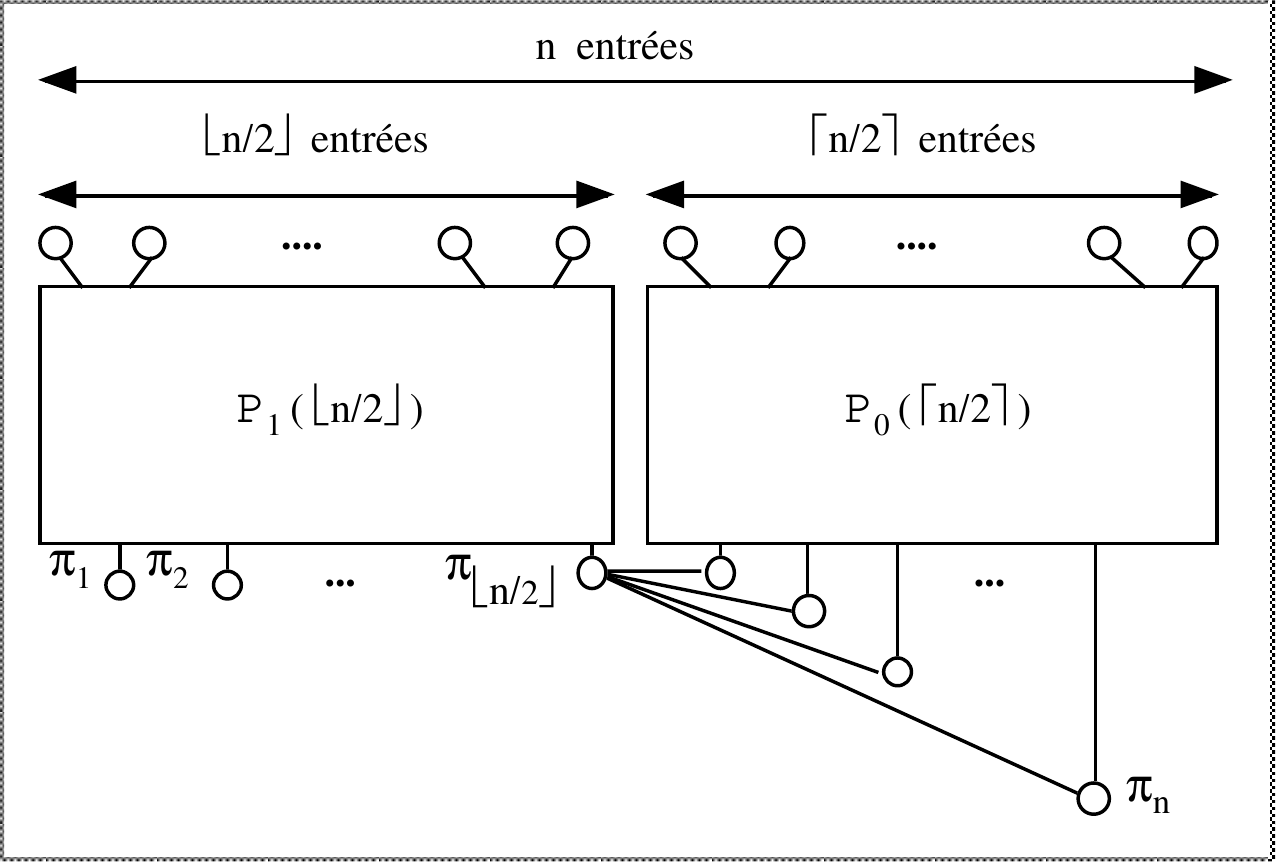}
\end{center}
\caption[Construction r\'ecursive des \cirs  
\,${\cal P}_0\,(n)\,$]{\label{F8}  
~Construction r\'ecursive des \cirs  
\,${\cal P}_0\,(n)$.}  
\end{figure}  

La construction du \cir \,$({\cal P}_1\,(n))$,  
quant \`a elle, se fait \`a partir du \cir  
\,${\cal P}_{0}\,(\lfloor {n \over 2} \rfloor)$,  
elle est  
illustr\'ee par la figure \vref{F9}.    
\begin{figure}[ht]   
\begin{center}
\includegraphics*[width=11cm]{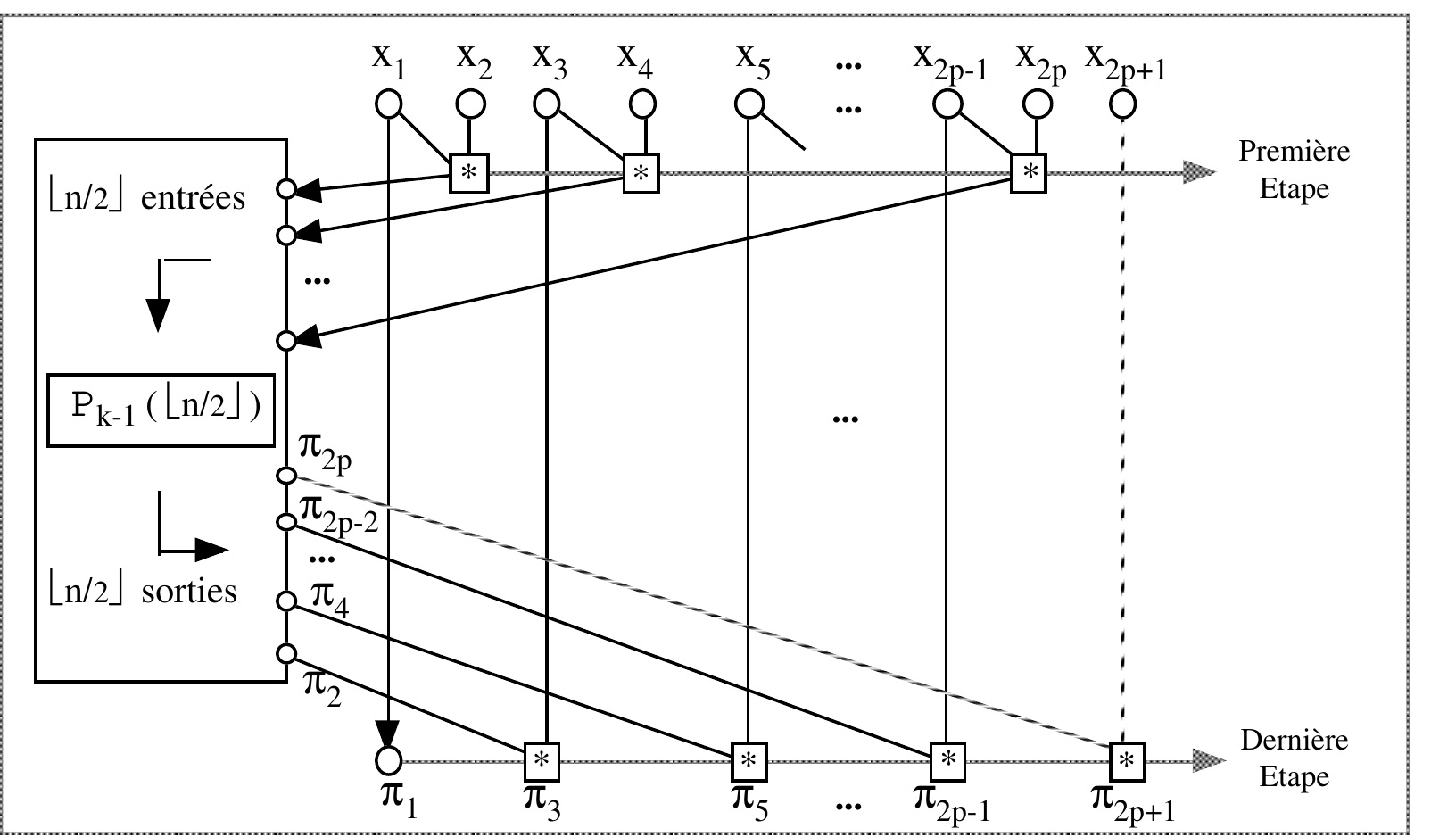}
\end{center}
\caption[Construction r\'ecursive des \cirs  
\,${\cal P}_1\,(n)\,$]{\label{F9} Construction  
du \cir \,${\cal P}_1\,(n)$\, \`a partir du  
\cir \,${\cal  P}_{0}\,(\lfloor {n \over 2}  
\rfloor)$}
  
\centerline{($p=\lfloor {n \over 2}  
\rfloor$\, et les 2 lignes en pointill\'e sont  
absentes si \,$n$\, est pair)}  
\end{figure}  

\paragraph{Construction de la famille \,$({\cal P}_1\,(n))_{\,n\in 
\NN^*}$} ~

\sni 
$\bullet~$ On commence par calculer en  
\paral (\cad en une seule \'etape) les  
produits $~x_1*x_2,\,x_3*x_4,\,\cdots  
,\, x_{2p-1}*x_{2p}~$ (o\`u \,$p = \lfloor  
{n\over 2}\rfloor$) d'un \elt de rang  
impair par l'\elt suivant (de rang pair)  
dans la liste donn\'ee \,$(x_1,\ldots ,  
,x_n)$\, ($n = 2p$\, si \,$n$\, est pair et  
\,$n = 2p + 1$\, si \,$n$\, est impair). 

\sni
$\bullet~$ \`A ce $p$-uplet on applique le \cir  
\,${\cal P}_{0}\,(\lfloor {n\over 2}\rfloor) =  
{\cal P}_{0}\,(p)$\, 
pour obtenir  
en sortie les \,$p$\, pr\'efixes de longueur paire:  
\,$\Pi_2,\,\Pi_4,\,\dots,\,\Pi_{2p}$. 

\sni 
$\bullet~$ On multiplie enfin, et en \paralz,  
les pr\'efixes \,$\Pi_2,\Pi_4,\dots,\Pi_{2p}$\,  
res\-pectivement par les entr\'ees \,$(x_3,x_5,\ldots,  
x_{2p-1}$\, (et \'even\-tuel\-lement \,$x_{2p+1}$\,  
si $n$ est impair) pour obtenir, en plus de  
\,$\Pi_1$\, ($\Pi_1=x_1$\, est d\'ej\`a donn\'e),  
les autres pr\'efixes de longueur impaire:  
\,$\Pi_3,\,\Pi_5,\,\dots,\,\Pi_{2p-1}$\,  
(et \'e\-ven\-tu\-el\-le\-ment \,$\Pi_{2p+1}$\, si  
\,$n$\, est impair).  
  
\ms On obtient ainsi le \cari \paral  
\,${\cal P}_{1}\,(n)$\, \`a partir du \cir  
\,${\cal P}_{0}\,(\lfloor {n\over 2}\rfloor)$\,  
en ajoutant au maximum deux \'etapes  
(\`a l'entr\'ee et \`a la sortie) comportant au  
total \,$n-1$\, \oparis  ($2p-1$\, si $n$ est pair  
et \,$2p$\, si $n$ est impair).  
  
Les circuits \vref{F8b} sont  
des exemples de \cirs \,${\cal P}_0\,(n)$\,  
et \,${\cal P}_1\,(n)$\, pour quelques valeurs  
de \,$n$.  
\setcounter{agh}{\arabic{algorithm}}
\stepcounter{algorithm}
\begin{agh}[ht]   
\begin{center}
\includegraphics*[width=11cm]{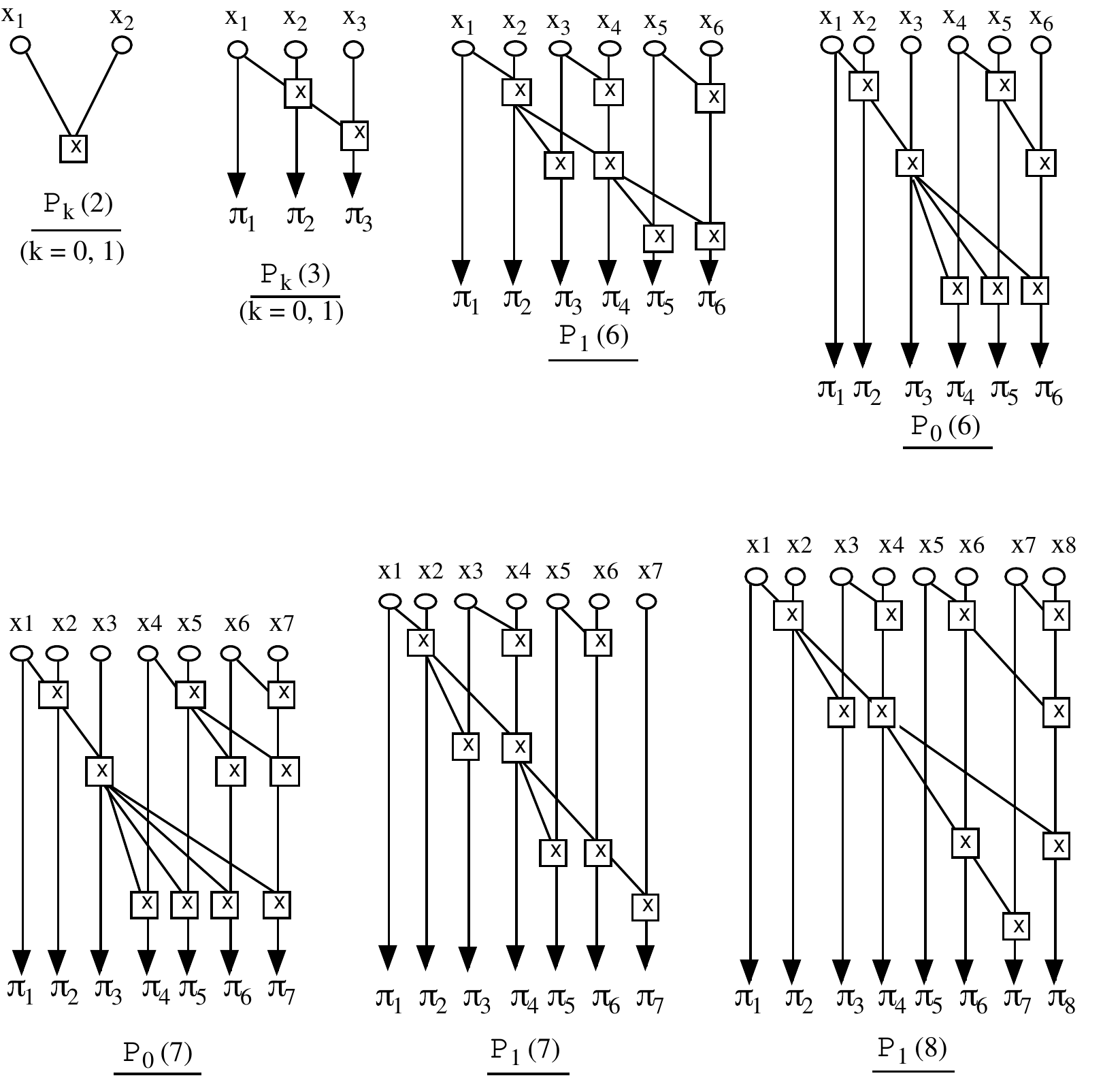}
\end{center}
\caption[Exemples de \cirs  
\,${\cal P}_k\,(n)~(k\in\{0,1\})\,$]{\label{F8b}  
~~Circuits \,${\cal P}_0\,(n)\,,~{\cal P}_1\,(n)~$  
\,pour quelques valeurs de \,$n$.}  
\end{agh} 

\paragraph{Analyse de la complexit\'e des circuits} ~  
  
\sni Si l'on note \,$S_k(n)$\, (resp. \,$D_k(n)$) la taille  
(resp. la \profz) du \cir \,${\cal P}_{k}\,(n)$\,   
pour \,$n\geq 2$\, et \,$k\in\{0,1\}$, cette construction  
r\'ecursive donne les  relations suivantes:  
\begin{itemize}  
\item Pour la taille:  
\begin{equation}\label{lad1}  
\formule{  
  S_1(n) &= &  
S_{0}(\lfloor {n\over 2}\rfloor) + n-1 \\[1mm]  
S_0(n) &= &  S_{1}(\lfloor {n\over 2}\rfloor) +  
S_{0}(\lceil {n\over 2}\rceil) + \lceil {n\over 2}\rceil  
}  
\end{equation}  
\item Pour la \profz:  
\begin{equation}\label{lad2}  
\formule{  
 D_1(n) &\leq &  
D_{0}(\lfloor {n\over 2}\rfloor) + 2 \\[1mm]  
D_0(n) &= & \max{\{D_{0}(\lfloor {1\over 2} \lfloor  
{n\over 2}\rfloor \rfloor)+1\,,  
\,D_{0}(\lceil {n\over 2} \rceil)}\}\,  
+ \,1 
}
\end{equation}
\end{itemize}  
$\mbox{avec } S_{k}(1) = D_{k}(1) = 0  
~\mbox{ pour tout } ~k\in\{0,1\}$.  
  
Il faut remarquer que l'\ine  
\,$D_1(n) \leq D_{0}(\lfloor  
{n\over 2}\rfloor) + 2$\, dans \,(\ref{lad2})\,  
peut \^etre stricte (voir par exemple le \cir  
\,${\cal P}_{1}\,(6)$\, dans les circuits \vref{F8b}  
pour s'en convaincre).  
  
La deuxi\`eme \'equation dans (\ref{lad2}) est  
justifi\'ee par le fait que,  
dans le \cir \,${\cal P}_{0}\,(n)$,  
le n{\oe}ud correspondant au produit  
\,$\Pi_{\lfloor {n \over 2} \rfloor}$\,  
-- dont on a besoin pour calculer en une  
\'etape \supt les autres pr\'efixes --  
se trouve exactement \`a la \prof  
\,$D_{0}\,(\lfloor {1\over 2} \lfloor  
{n\over 2}\rfloor \rfloor)+1$\,  
dans le sous-\cir \,${\cal P}_{1}\,(\lfloor  
{n\over 2}\rfloor)$\, de  
\,${\cal P}_{0}\,(n)$\, qui calcule ce produit.  
  
\ms  
Il est facile de voir, \`a partir des \'equations  
(\ref{lad2}), par une \recu  imm\'ediate sur  
\,$n$, que les \profs \,$D_k(n)$\, des \cirs  
\,${\cal P}_{k}\,(n)$\, pour \,$k\in\{0,1\}$\,   
v\'erifient:  
$$ D_0(n) = \esup{\log n}~~~
\mbox{ et }~~~D_1\,(n) \leq \esup{\log n} + 1\,.$$  
  
\ms  
Pour calculer les tailles des \cirs \`a partir  
des \'equations (\ref{lad1}),  nous allons d'abord  
consi\-d\'e\-rer le cas o\`u \,$n$\,  
est une puissance de 2 en faisant  
\,$n=2^\nu$\, o\`u \,$\nu=\esup{\log n}$.  
  
Posant \,$\tau_k(\nu) = S_k(2^\nu)\,\mbox{ avec  }  
\tau_{k}(0) = 0\, \mbox{ pour } k\in\{0,1\}$\, les  
\'equations \,(\ref{lad1})\, deviennent:  
  
\begin{equation}\label{lad3}  
\formule{    
\tau_0(\nu)& = &   
\tau_{0}(\nu-1) + \tau_{0}(\nu-2) + 2^\nu-1 \\[1mm]  
\tau_1(\nu)& = &  \tau_1(\nu-1) +  
\tau_1(\nu-2) + 3.2^{\nu-2}\,.  
}  
\end{equation}  

Posant \,$u_0(\nu) = 4.2^\nu +1 -  
\tau_0(\nu)$\, et  
\,$u_1(\nu) = 3.2^\nu - \tau_1(\nu)$,  
les relations \,(\ref{lad3})\, permettent de  
v\'erifier que  
\,$u_k(\nu + 2) = u_k(\nu + 1) + u_k(\nu)$\,  
($\nu \in\N$,  \,$k\in\{0,1\}$).  
Comme \,$u_0(0) = 5\,,$ \,$u_0(1) = 8\,,$
\,$u_1(0) =  3$\, et \,$u_1(1) = 5\,,$  
on en d\'eduit que: $$u_0(\nu) = F(\nu + 5)  
\,\mbox{ et } u_1(\nu) = F(\nu + 4)$$  
o\`u \,$(F(\nu))_{\nu\in \NN}$\, est la suite de  
Fibonacci\footnote{~La suite de  
Fibonacci est d\'efinie par \,$F(0)=0,\,F(1)=1$\,  
et la relation  
\,$F(\nu + 2) = F(\nu + 1) + F(\nu)$\, pour tout  
\,$\nu\in \NN$.}. Par cons\'equent:  

\begin{equation}\label{lad4}  
\formule{  
\tau_0(\nu)& = &  4.2^\nu +1 - F(\nu + 5) \\[1mm]  
\tau_1(\nu) &= &  3.2^\nu - F(\nu + 4)\,.  
}
\end{equation}  
qui donne, lorsque \,$n$\, est une puissance de 2,  
les majorations souhai\-t\'ees. 
Dans le cas contraire, il est facile --  
en utilisant directement les relations  
\,(\ref{lad2}) -- d'obtenir, par  
\recu sur \,$n$, les majorations suivantes  
vraies pour tout \,$n \geq 2$:  
$$S_0(n)\leq 4\,n - 7~~  
\mbox{ et }~~S_1(n)\leq 3\,n - 3\,.$$  
  
Ce qui donne le r\'esultat suivant de Ladner  
\& Fischer (\cite{Lad}) qui montre que le calcul des pr\'efixes  
est \,$\PRAM(n/\log{n},\log{n})$:  
  
\begin{theorem}\label{prefix}  \emph{(Ladner  \& Fischer)}
Le calcul des pr\'efixes d'une liste de $n$  
\elts dans un mono\"{\i}de  
(non \ncrt commutatif) se fait par  
un \cir \arith \paral de \prof  
\,$\esup{\log n}$\, et de taille  
inf\'erieure \`a \,$4n$\, et aussi par   
un \cir \arith \paral de \prof  
\,$1+\esup{\log n}$\, et de taille  
inf\'erieure \`a \,$3n$.  
\end{theorem} 

\newpage \thispagestyle{empty}

\chapter{Multiplication rapide des \pols} 
\label{sec multipoly}\label{chap multipoly}
\minitoc
\acvide

\subsubsection*{Introduction}   

Soit \,$\A$\, un \acom unitaire
et \,$\A[X]$\, l'anneau des \pols \`a une \idtr  
sur \,$\A$. 

Le produit de deux \pols \,$A=\sum_{i=0}^na_iX^i$\, et  
\,$B=\sum_{i=0}^mb_iX^i$\, est d\'efini par 
$$\,C=AB=\somm_{k=0}^{m+n}c_kX^k~~~\mbox{ avec }~~c_k=  
\somm_{i=0}^ka_ib_{k-i}\,~\mbox{ pour }~\,0\leq k\leq 
m+n\,.
$$  
L'\algo \usl pour le calcul des \coes du \poly  
\,$C$\, cor\-respond \`a un \cari de \prof  
\,$\O(\log{m})$\, (si l'on suppose \,$m\leq n$) et de  
taille \,$\O(mn)$\, avec pr\'ecis\'ement \,$(m+1)\,(n+1)$\,  
\muls et \,$mn$\, additions dans l'anneau de base  
\,$\A$. \, 
Pour \,$m=n$, cela donne un \algo en  
\,$\SD(n^2,\log{n}).$

\ms Dans les trois premi\`eres sections nous exposons  
deux fa\c{c}ons d'am\'e\-lio\-rer la \mul des \polsz.

\ss Dans la section \ref{secKarat} nous expliquons 
la \met de 
Karatsuba, facile \`a impl\'ementer pour n'importe
quel anneau commutatif, avec un r\'esultat en 
\,$\SD(n^{\log{3}},\log{n})$. 

\ss Un bien meilleur r\'esultat est obtenu en  
\,$\SD(n\log{n},\log{n})$\, gr\^ace \`a la transformation  
de Fourier discr\`ete (\cite{Aho,Knu}) pour un anneau  
auquel s'applique une telle transformation. Ceci fait 
l'objet des 
sections \ref{secTFDu} et \ref{subsecTFDfa}.

\ss Dans la section \ref{TFDmrGen} nous exposons une 
am\'elioration
 due \`a  
Cantor et Kaltofen \cite{Cant} qui ont \'etendu le   
r\'esultat \`a tout \acom unitaire en exhibant  
un \algo en \,$\SD(n\,{\log{n}}\,\log{\log{n}},  
\log{n})$\, (avec
le m\^eme nombre de \muls dans l'anneau de base, le  
facteur \,$\log{\log{n}}$\, \'etant d\^u \`a l'augmentation  
du nombre d'additions).
Pour r\'ealiser ce travail il a fallu
 l'adjonction de racines principales de l'unit\'e \`a 
l'anneau  
consid\'er\'e. 
On peut comparer la borne obtenue avec la meilleure borne  
inf\'erieure actuellement connue, qui est \,$\O(n)$. 
 
\ss Dans la section \ref{secMTT} nous donnons le lien
entre la \mul des \pols et celle des \mttisz.
Nous en d\'eduisons un r\'esultat de complexit\'e 
int\'eressant concernant le produit d'une \mto 
arbitraire par une matrice arbitraire.

\section{M\'ethode de Karatsuba}\label{secKarat}

Consid\'erons deux \pols arbitraires \,$A$\, et \,$B$\,
et leur produit \,$C$. Si les \pols  \,$A$\, et \,$B$\, 
sont de degr\'e $<d$ (d\'etermin\'es chacun par \,$d$\, 
\coesz), 
leur produit \,$C=AB$\, peut \^{e}tre calcul\'e en 
appliquant directement la formule qui le d\'efinit. 
Il y a alors \,$d^2$\, \muls et $(d-1)^2$ additions. Les 
\,$d^2$\, \muls peuvent \^{e}tre calcul\'ees 
en une seule \'etape de 
calcul \paral et les \,$2d-1$\, \coes de \,$C$\, sont 
ensuite calcul\'es en $\esup{\log  d}$ 
\'etapes \parals (le \coe r\'eclamant 
l'addition la plus longue est celui de degr\'e $d-1$).
  
Une premi\`ere fa\c{c}on d'am\'eliorer cette \mul
  est d'adop\-ter une d\'emarche r\'ecursive
bas\'ee sur le fait que le produit de deux  
\pols de de\-gr\'e~$1$ peut s'effectuer avec seulement  
$3$ \muls au lieu de $4$ (le nombre 
d'additions/soustractions 
passant de $1$ \`a $4$). En effet, on peut calculer  
\,$a+bX+cX^2=(a_1+a_2X)\,(b_1+b_2X)$\, en posant: 
\begin{equation} \label{EqKara1}
a=a_1b_1, \quad \quad c=a_2b_2,  
\quad \quad b=(a_1+a_2)\,(b_1+b_2)-(a+c),
\end{equation}
ce qui correspond \`a un \cari de \prof totale $3$, de 
largeur $4$ et de \prof \muv $1$.

Consid\'erons maintenant deux \pols arbitraires \,$A$\, et 
\,$B$ et leur produit  \,$C\,$. 
Ces \pols  s'\'ecrivent de mani\`ere unique, sous la forme:
$$
\left\{\begin{array}{lcl}  
A & = & A_1(X^2)+X\,A_2(X^2) \\  
B & = & B_1(X^2)+X\,B_2(X^2) \\  
C & = & C_1(X^2)+X\,C_2(X^2) 
\end{array}\right.  
$$
avec \,$C_1  =  A_1B_1+X\,A_2B_2$\, et 
\,$C_2  =  A_1B_2+A_2B_1.$  
Si  \,$A_1,\,B_1,\,A_2,\,B_2$\, sont de degr\'es $\leq k-1$
(avec $k$ \coesz)
alors  \,$A$\, et \,$B$\, sont de degr\'e $\leq  2k-1$
(avec $2k$ \coesz).

Supposons qu'un \prev $Kara^{(k)}$ calcule les \coes 
du produit de deux \pols arbitraires de degr\'e $\leq k-1$,
avec une \prof \muv \'egale \`a $\mu^{(k)}$, une \prof
totale  \'egale \`a \,$\pi^{(k)}$, une largeur
\'egale \`a \,$\lambda^{(k)}$, un nombre de \muls
\'egal \`a \,$m^{(k)}$, un nombre d'\adsos
\'egal \`a \,$a^{(k)}$, et donc avec pour nombre
total d'\oparis \,$s^{(k)}=a^{(k)}+m^{(k)}$. 
L'utilisation des \'equations (\ref{EqKara1})
donne un \cari $Kara^{(2k)}$ que nous avons d\'ecrit 
sch\'e\-ma\-tiquement dans le \prev \vref{prevKara}.

\begin{proeva}[$Kara^{(2k)}$] \label{prevKara}
    \acl{prevKara}{Produit de deux \pols \`a la Karatsuba}
\Entree Les $4k$ \coes dans $\A$ (un \acomaz) de deux \pols 
de degr\'e 
$<2k$~: $A(X)  =  A_1(X^2)+X\,A_2(X^2)$ et $B(X)  =  
B_1(X^2)+X\,B_2(X^2)$.
\Sortie Les \coes du produit des deux \polsz: $C(X)  =  
C_1(X^2)+X\,C_2(X^2)$.
\Debut
\Pro{1}{}
 {D_{1}:=A_1+A_2\;;\;D_{2}:=B_1+B_2}
\Pro{\pi^{(k)}}{}{D_3:=Kara^{(k)}(A_1,B_1)\;;\;D_4:=Kara^{(k
)}(A_2,B_2
)}
\Pro{\pi^{(k)}+1}{}{D_5:=Kara^{(k)}(D_1,D_2)\;;\;D_{6}:=D_3+
D_4\;;\;C_
1:
=D_3+XD_4 }
\Pro{\pi^{(k)}+2}{}{C_2:=D_5-D_6}
\fin
\end{proeva}

Notez que la ligne \'ecrite avec la \prof $\pi^{(k)}$ 
repr\'esente la derni\`ere ligne des deux \prevs 
\,$Kara^{(k)}(A_1,B_1)$\, et 
\,$Kara^{(k)}(A_2,\alb B_2)$,
qui ont d\'emarr\'e en \paral avec les deux affectations 
indiqu\'ees
sur la ligne de \prof $1$. Sur la ligne \'ecrite 
avec la \prof $\pi^{(k)}+1$, la premi\`ere affectation 
correspond
\`a la derni\`ere ligne du \prev \,$Kara^{(k)}(D_1,D_2)$\, 
qui a commenc\'e
\`a la \prof $2$, tandis que les deux autres affectations 
sont effectu\'ees \`a la \prof \,$\pi^{(k)}+1$.

On constate donc que lorsqu'on passe de \,$Kara^{(k)}$\, \`a 
\,$Kara^{(2k)}$\,
selon la m\'ethode d\'ecrite dans le \prev \ref{prevKara}:
\begin{itemize}
\item  la \prof  passe de \,$\pi^{(k)}$\, \`a 
\,$\pi^{(2k)}=\pi^{(k)}+2$,
\item  la \prof \muv n'a pas chang\'e  
($\mu^{(2k)}=\mu^{(k)}$),
\item  la largeur passe de \,$\lambda^{(k)}$\, \`a
\,$\lambda^{(2k)}=\sup(3\,\lambda^{(k)},\lambda^{(k)} + 4k-
2)$,
\item  le nombre de \muls est maintenant 
\,$m^{(2k)}=3m^{(k)}$,
\item   le nombre d'\adsos
(\footnote{~On ne compte pas les \ops de  
substitution de \,$X^2$\, \`a \,$X$\, ou vice-versa, ni les 
\muls par 
\,$X$\, ou  
par \,$X^2$, qui reviennent en fait \`a des d\'ecalages de  
\coesz.}) 
est \,$a^{(2k)}=k+k+3a^{(k)}+(2k-1)+(2k-2)+(2k-
1)=3a^{(k)}+8k-4$,
\item  et le nombre total
d'\oparis passe de  \,$s^{(k)}$\,  \`a 
\,$s^{(2k)}=3s^{(k)}+8k-4$.
\end{itemize}

En comparaison, pour la multiplication  \usle des \polsz, 
le nombre de \muls
passe de \,$\tilde{m}^{(k)}=k^2$\, \`a 
\,$\tilde{m}^{(2k)}=4k^2=4\tilde{m}^{(k)}$, 
le nombre d'\adsos de $\tilde{a}^{(k)}=(k-1)^2$ \`a 
\,$\tilde{a}^{(2k)}=(2k-1)^2=4\tilde{a}^{(k)}+4k-3$\, 
et le nombre total d'\oparis de \,$\tilde{s}^{(k)}$\, \`a 
\,$\tilde{s}^{(2k)}=4\tilde{s}^{(k)}+4k-3$.

Si on veut minimiser le nombre de \muls on initialisera
la processus r\'ecursif avec  \,$Kara^{(1)}$\, 
(le produit de deux constantes) et on mettra
en place les \caris successifs \,$Kara^{(2)}$, 
\,$Kara^{(4)}$, 
\,$Kara^{(8)}$, 
\ldots, \,$Kara^{(2^\nu)}$ selon la \pcd d\'ecrite
ci-dessus. Le circuit  \,$Kara^{(2^\nu)}=Kara_\nu$\, est 
ensuite utilis\'e
pour le produit de deux \pols de degr\'es \,$<n=2^\nu$\, et 
\,$\geq 
2^{\nu-1}$.
Pour deux \pols de degr\'e exactement $n-1$ on aura ainsi 
remplac\'e 
le \cari \usl qui utilise \,$4^\nu=n^2$\, \muls par un \cari
\,$Kara^{(n)}=Kara_\nu$\, qui utilise \,$3^\nu=n^{\log 
3}\simeq 
n^{1.585}$\, \mulsz{\footnote{~$\log 
3=1.58496250072115618145373894394$.}}. 
Le gain concernant le nombre total d'\oparis est du m\^{e}me 
style.
En notant \,$s_{\nu}$\, pour \,$s^{(n)}$,
on passe en effet de \,$s_{\nu}$\, \`a 
\,$s_{\nu+1}=3\,s_\nu+8.2^{\nu}-
4$. 
Les premi\`eres valeurs  de \,$s_{\nu}$\, sont \,$s_{0}=1$, 
\,$s_{1}=7$, 
\,$s_{2}=33$\, et la relation de \recu se r\'esoud avec 
l'aide de 
Maple en: 
$$s_\nu=7\cdot3^\nu-8\cdot2^\nu+2.
$$ 
En fait \,$s_{\nu}$\, devient meilleur
que \,$4^n+(2^n-1)^2$\, \`a partir de \,$\nu=4$\, (pour des 
\pols
de degr\'e \,$15$). Enfin, concernant la largeur 
\,$\lambda_{\nu}$\,
du \cari \,$Kara_\nu$, la 
r\'esolution de la \recu donne  
\,$\lambda_{\nu}=2\cdot3^\nu$\, pour \,$\nu\geq 2$.

\ss Nous pouvons conclure avec la proposition suivante.
\begin{prop}\label{prop Kara}
La \mul de deux  
\pols de degr\'e \,$\leq  n$\, par la \met de Karatsuba
se fait en  
\,$\SD(n^{\log{3}},\log{n})$.
Plus pr\'ecis\'e\-ment, le produit de deux \pols de degr\'es 
\,$<2^\nu=n$
peut \^{e}tre r\'ealis\'e par un \cari de \prof \muv $1$,
de \prof totale \,$1+2\,\nu$, de largeur 
\,$2\cdot3^\nu=2\,n^{\log{3}}$, 
avec \,$3^\nu=n^{\log{3}}$\, \muls et 
\,$6\cdot3^\nu-8\cdot2^\nu+2=6n^{\log{3}}-8n+2$\,
\adsosz.
\end{prop}  

Notons que pour deux \pols dont les degr\'es sont compris
entre \,$2^{\nu-1}$\, et \,$2^\nu$\,, on 
obtient seulement les majorations suivantes
en appelant \,$n$\, le plus grand degr\'e: \,$3+2\,\log 
n$\, 
pour la \profz, \,$6\,n^{\log{3}}$\, pour la largeur et 
\,$21\,n^{\log{3}}-8\,n+2$\, pour la taille du circuit.

\ss Remarquons qu'on aurait pu envisager une autre  
partition de \coes des \pols \,$A$\, et  
\,$B$\, pour une application r\'ecursive,  
\`a savoir 
$$A = A_1+X^kA_2~~~  
\mbox{ et }~~~B = B_1+X^kB_2\,.
$$   
avec \,$A_i$\, et \,$B_i$\, de degr\'es \,$\leq k-1$, 
\,$A$\, et \,$B$\, de degr\'es $\leq 2k-1$.
Alors 
\,$C=A\,B=A_1\,B_1+X^k(A_1\,B_{2}+A_2\,B_1)+X^{2k}\,A_2\,B_2
$.
Une \pcd r\'ecursive bas\'ee sur
cette partition produirait des \caris avec une  estimation 
analogue \`a la pr\'ec\'edente pour ce qui concerne 
la taille mais une \prof de \,$1+3\,\log n$\, au lieu 
de \,$1+2\,\log n$\, (pour le produit de deux 
\pols de degr\'e \,$n-1$\, lorsque \,$n=2^{\nu}$).

\section{Transformation de Fourier discr\`ete usuelle} 
\label{secTFDu}

Un bien meilleur r\'esultat, que nous exposons dans
cette section et la suivante, est obtenu en  
\,$\SD(n\log{n},\log{n})$\, gr\^ace \`a la \emph{transformation  
de Fourier discr\`ete} pour un anneau  
auquel s'applique une telle transformation.  
La transformation de Fourier discr\`ete, que nous 
d\'esignerons  
ici par le sigle TFD, est d\'efinie sur un \acom  
unitaire \,$\A$, pour un entier donn\'e \,$n\geq 2$, \`a  
condition de disposer dans \,$\A$\, 
d'une \emph{racine \,$n$\,-\,\`eme  
principale}
\index{racine de l'unit\'e!principale \,$n\,$-\,\`eme} de 
$1$, 
\cad d'un \elt \,$\xi \in\A$\,  
v\'erifiant: $$\,\xi\neq 1\,, \,\xi^n=1\,, \mathrm{\;et\;}  
\,\somm_{j=0}^{n-1}\xi^{ij}= 0\, \mathrm{\;pour\;} 
\,i=\,1,\ldots, n-
1\,.$$ 

Dans un anneau int\`egre, toute racine 
primitive\footnote{~C'est un \,$\xi$\, tel que \,$\xi^n=1$\, 
mais  \,$\xi^m\neq 1$\, si \,$1\leq m<n$.}  
\,$n$\,-\,\`eme de $1$ est
\index{racine de l'unit\'e!primitive \,$n\,$-\,\`eme} 
principale, mais ceci  
peut-\^etre mis en d\'efaut dans un anneau contenant  
des diviseurs de z\'ero. Dans un anneau int\`egre,  
s'il y a une racine primitive \,$n\,$-\,\`eme de $1$,  
il y en a \,$\varphi(n)$\, o\`u \,$\varphi$\,  
d\'esigne l'indicatrice d'Euler.
Dans \,$\CC$, les racines \,$n$\,-\,\`emes principales  
de $1$ sont les nombres complexes e$^{2ik\pi/n}$  
tels que \,$1\leq k<n$\, et \,$k$\, premier avec \,$n$. 

Il est clair que si \,$\xi$\, est une racine  
\,$n\,$-\,\`eme principale de $1$, alors il en est  
de m\^eme de \,$\xi^{-1}$.  
  
\begin{defi}  
La transformation de Fourier discr\`ete d'ordre  
\,$n$\, sur \,$\A$, associ\'ee \`a la racine  
principale \,$\xi$, est l'application \lin  
$$\TFD_{n,\xi} ~:~ \A^n \dans \A^n$$ d\'efinie,  
pour tout $~(a_0,\,a_1,\,\ldots,\,a_{n-1}) \in  
\A^n~$ par:  
$$\TFD_{n,\xi}\,(a_0,\,a_1,\,\ldots,\,a_{n-1})  
= (A(1),\,A(\xi),\,\ldots,\, A(\xi^{n-1}))$$  
o\`u \,$A$\, est le \poly  
\,$A(X)=a_0+a_1X+\ldots+a_{n-1}X^{n-1}\,.$
\end{defi}
  
Cette application peut aussi \^etre vue comme un  
homomorphisme de \,$\A$-\agrs  
$$\TFD_{n,\xi} ~:~\aqo{\A[X]}{X^n-1} \dans \A^n~$$  
qui \`a tout \poly \,$A$\, de degr\'e $\leq n-1$\,  
associe le vecteur form\'e des valeurs de \,$A$\, aux  
points \,$1,\xi,\ldots,\xi^{n-1}$.  
En effet, en notant \,$\odot$\, la loi produit  
(coordonn\'ee par coordonn\'ee) de l'\agr  
$\A^n$, il est imm\'ediat de v\'erifier que:  
$$\TFD_{n,\xi}(AB)=\TFD_{n,\xi}(A) \odot  
\TFD_{n,\xi}(B)\,.$$
  
En tant qu'\aliz,  
\,$\TFD_{n,\xi}$\, est repr\'esent\'ee  
dans les bases canoniques par la matrice  
de Vandermonde particuli\`ere:
$$
W_{n,\xi}=\cmatrix{  
1 & 1 & 1 & \ldots & 1 \cr  
1 & \xi & \xi^2 & \ldots & \xi^{n-1} \cr  
1 & \xi^2 & \xi^4 & \ldots & \xi^{2(n-1)} \cr  
\vdots & \vdots & \vdots & \ddots & \vdots \cr  
1 & \xi^{n-1} & \xi^{2(n-1)} & \ldots & \xi^{(n-1)^2}  
}\,.
$$
Si, de plus, \,$n\,1_\A$\, est inversible dans l'anneau  
\,$\A$\, (on d\'esignera par \,$n^{-1}$\, son inverse),  
alors la matrice \,$W_{n,\xi}$\, est inversible  
dans \,$\A^{n\times n}$\, et on v\'erifie qu'elle admet  
pour inverse la matrice  
$$
W_{n,\xi}^{-1}=n^{-1}
\cmatrix{ 
1 & 1 & 1 & \ldots & 1 \cr  
1 & \xi^{-1} & \xi^{-2} & \ldots & \xi^{1-n} \cr  
1 & \xi^{-2} & \xi^{-4} & \ldots & \xi^{2(1-n)} \cr  
\vdots & \vdots & \vdots & \ddots & \vdots \cr  
1 & \xi^{1-n} & \xi^{2(1-n)} & \ldots & \xi^{-(n-1)^2}  
} 
= n^{-1}\,W_{n,{\xi^{-1}}}\,.
$$  
Dans ce cas, modulo l'identification pr\'ec\'edente,  
l'application \,$\TFD_{n,\xi}$\, est un isomorphisme  
d'\agrs \,$\TFD_{n,\xi}\,:\,\aqo{\A[X]}{X^n-1}  
\dans \A^n$. 
Nous \'enon\c{c}ons ce r\'esultat.

\begin{prop}\label{prop TFD} 
Supposons que l'\acom  \,$\A$\, poss\`ede une racine 
\,$n$\,-\,\`eme principale de $1$, not\'ee \,$\xi$, et que  
\,$n\,1_\A$\, est inversible dans \,$\A$.  
Alors la transformation de Fourier discr\`ete  
\,$\TFD_{n,\xi}~:\A^n\dans\A^n\,$
est un isomorphisme de \Amosz, et 
$\TFD_{n,\xi}^{-1}=  
(n1_{\A})^{-1}\TFD_{n,\xi^{-1}}$.   
Par ailleurs, si on identifie le \Amo $\A^n$ 
source de l'\ali \,$\TFD_{n,\xi}$\, avec  
$\aqo{\A[X]}{X^n-1}$ (en choisissant le repr\'esentant
de degr\'e $\leq n$ et en l'exprimant sur la base des
mon\^omes) alors  \,$\TFD_{n,\xi}$\, d\'efinit  
un isomorphisme de l'\agr \,$\aqo{\A[X]}{X^n-1}$\,  
(munie de la \mul des \polsz) vers l'\agr   
\,$\A^n$\, (munie de la \mul \,$\odot$\,  
coordonn\'ee par coordonn\'ee).  
En bref, pour deux \pols de degr\'e $< n$, on a:  
$$
AB\,\equiv\, TFD_{n,\xi}^{-1}\,(TFD_{n,\xi}\,(A) \odot  
TFD_{n,\xi}\,(B)) \; \;\mathrm{modulo}\;(X^n-1).
$$  
\end{prop}

\noindent  
C'est la cl\'e de l'\algo de \mul rapide,  
que nous explici\-tons dans la section suivante.

\section{Transformation de Fourier discr\`ete rapide}  
\label{secTFDmr}  

\subsection{Cas favorable}\label{subsecTFDfa} 
Le r\'esultat \'enonc\'e dans la proposition  
\ref{prop TFD} pr\'ec\'edente peut \^etre appliqu\'e  
au calcul du produit \,$AB=\sum_{k=0}^{2n-1} c_kX^k$\,  
de deux \pols \,$A=\sum_{i=0}^{n-1}a_iX^i$\, et  
\,$B=\sum_{i=0}^{n-1}b_iX^i$\, \`a une \idtr  
sur \,$\A$, \`a condition que l'anneau \,$\A$\, s'y  
pr\^ete. Nous supposons qu'il poss\`ede une racine  
\,$2n\,$-\,\`eme principale \,$\omega$\, de $1$, et que  
\,$(2n)\,1_\A$\, est inversible dans \,$\A$, alors la  
proposition \ref{prop TFD} pour la TFD d'ordre \,$2n$\,  
sur l'anneau \,$\A$\, se traduit par:  
$$AB=TFD_{2n,\omega}^{-1}\,(TFD_{2n,\omega}\,(A) \odot  
TFD_{2n,\omega}\,(B))\,.$$  
car le calcul de \,$AB$\, modulo \,$X^{2n}-1$\, donne  
exactement \,$AB$.  
Le calcul du produit \,$AB$\, de deux \pols
de degr\'es inf\'erieur ou \'egal \`a \,$n$\, 
par la TFD est r\'esum\'e dans l'\algo \vref{algoTFD}. 
\begin{algor}[Multiplication des \pols via la Transformation 
de Fourier Discr\`ete.] 
\label{algoTFD}
    \acl{algoTFD}{Transformation de Fourier Discr\`ete}
\Entree Deux \pols \,$A$\, et \,$B$\, de degr\'es \,$<n$\, 
sur un anneau
\,$\A$\, convenable (voir proposition \ref{prop TFD}). 
\Sortie Le produit $A\,B$. 
\Debut 
\Etap{1}{}{Deux TFD d'ordre \,$2n$\,  
appliqu\'ees \`a \,$A$\, et \,$B$.}  
\Etap{2}{}{\'Evaluation de  \,$2n$\,  
\muls dans \,$\A$\, pour obtenir la  
transform\'ee}
\hsu  de Fourier discr\`ete de \,$A\,B$.
\Etap{3}{}{Calcul de l'inverse d'une TFD d'ordre \,$2n$\, 
pour  
obtenir \,$A\,B$.}              
\fin 
\end{algor}

Le lemme suivant nous permet tout d'abord de  
montrer comment une TFD d'ordre \,$2^{\nu}$\,  
peut \^etre effectu\'ee rapidement au moyen  
d'une strat\'egie \gui{diviser pour gagner}.  

\begin{lem}\label{lem TFD}  
Soit \,$n$\, un entier \,$\geq 2$\, et \,$\nu=\lceil  
\log{n}\rceil$. La transformation de Fourier discr\`ete  
d'ordre \,$n$\, et son inverse, dans un anneau poss\'edant  
une racine \,$2^\nu\,$-\,\`eme principale de $1$ et dans  
lequel \,$2_\A$\, est inversible, se font en  
\,$\SD(n\log{n},\log{n})$.  
Plus pr\'ecis\'ement, la taille \,$S(n)$\, et la \prof  
\,$D(n)$\, du \cari correspondant sont respectivement  
major\'ees par \,$ n\,(3\,\log{n}+3)$\, et  
\,$2\,\log{n}+2$\, pour la transformation directe et  
par \,$n\,(3\,\log{n}+4)$\, et \,$2\,\log{n}+3$\,  
pour la transformation inverse.
\end{lem}  
\prv  
Soit \,$A=\sum_{i=0}^{n-1}a_iX^i$\, un \poly de  
degr\'e $\leq n-1$\, \`a \coes dans \,$\A$,  
\,$\nu=\lceil \log{n}\rceil$\, (de sorte que  
\,$2^{\nu-1} < n \leq 2^\nu$) et \,$\omega$\,  
une racine \,$2^\nu\,$-\,\`eme principale de $1$.    
Il s'agit de calculer les valeurs de \,$A$\, aux points  
\,$1,\omega,\omega^2,\ldots,\omega^{2^\nu-1}$.  
Le \poly \,$A\,$
peut  \^etre mis sous la forme   
\,$A=A_1\,(X^2)+XA_2\,(X^2)$\, avec  
\,$\deg{A_1},\,\deg{A_2}\,\leq 2^{\nu-1}-1$.  
  
\ss  
Remarquons que \,$\xi=\omega^2$\, est une racine  
\,$2^{\nu-1}\,$-\,\`eme principale de $1$, que  
\,$\omega^{2^{\nu-1}}=-1$\, et que  
\,$A\,(\omega^{i})=A_1\,(\xi^i) + \omega^{i}\,  
A_2\,(\xi^i)$\, pour \,$0\leq i\leq 2^{\nu-1}-1$.  
Comme \,$\omega^{2^{\nu-1}+i}=-\omega^i$, on a aussi   
\,$A\,(\omega^{2^{\nu-1}+i})=A_1\,(\xi^i) - \omega^{i}\,  
A_2\,(\xi^i)$\, pour \,$0\leq i\leq 2^{\nu-1}-1$.  
  
\ss Ce qui donne toutes les valeurs recherch\'ees de \,$A$\,  
 et ram\`ene r\'ecursi\-ve\-ment l'\eva de  
\,$A$\, en les \,$2^\nu$\, points  
\,$\omega^i~~(0\leq i\leq {2^\nu-1})$,  
\cad la TFD d'ordre \,$2^\nu$, au calcul suivant:  
  
\ss  
$\bullet~$ deux TFD d'ordre \,$2^{\nu-1}$\, appliqu\'ees  
\`a \,$A_1$\, et \,$A_2\, $ et effectu\'ees en \paralz;  
  
$\bullet~$ \,$2^{\nu-1}$\, \muls (par les  
\,$\pm\,\omega^i$\, avec \,$0\leq i\leq 2^{\nu-1}-1)$\,  
effectu\'ees en \paral et en une seule \'etape  
de calcul, suivies de \,$2^\nu$\, additions dans l'anneau  
de base \,$\A$\, effectu\'ees \egmt en une seule  
\'etape \paralz.  
  
\ss  
Si \,$S$\, et \,$D$\, d\'esignent respectivement la  
taille et la \prof de l'\algo r\'ecursif ainsi d\'efini,  
on obtient les relations suivantes valables pour tout  
entier \,$\nu \geq 1$:
$$
\left\{\begin{array}{lcl}  
S(2^\nu) & \leq & 2\,S(2^{\nu-1}) + 3\,2^{\nu-1} \\  
D(2^\nu) & \leq & D(2^{\nu-1}) + 2 \,.  
\end{array}\right.
$$

Ce qui donne, par sommation, sachant que  
\,$S(1)=D(1)=0$:
$$
\left\{\begin{array}{lcl}  
S(2^\nu) & \leq & 3\,\nu\,2^{\nu-1} \\  
D(2^\nu) & \leq & 2\,\nu\,.  
\end{array}\right.
$$
Comme \,$2^{\nu-1}< n \leq 2^\nu$\, et par cons\'equent  
 \,$\nu-1< \log{n}$, on en  
d\'eduit que \,$S(n)< 3\,n\,(1+\log{n})$\, et que  
\,$D(n)< 2\,(1+\log{n})$.  
  
\ss  
Pour la TFD inverse d'ordre \,$n$, nous avons vu que  
\,$TFD_{n,\omega}^{-1}=(n\,1_\A)^{-1}\,TFD_{n,\omega^{-1}}$.   
Cela signifie que l'on peut r\'ecup\'erer les \coes  
du \poly \,$A$\, de degr\'e $\leq n-1$, \`a partir  
du vecteur \,$\vec{A}=(A(1),\,\alb A(\omega),\alb\ldots,  
\alb\,A(\omega^{n-1}))$\, form\'e des valeurs de ce \poly  
aux points \,$\omega^i$, en effectuant sur le vecteur  
\,$\vec{A}$\, la TFD d'ordre \,$n$\, associ\'ee \`a la  
racine principale \,$\omega^{-1}=\omega^{n-1}$\, et en  
multipliant ensuite ce vecteur par \,$(n\,1_\A)^{-1}$. 
Par cons\'equent, la TFD inverse d'ordre \,$n$\, peut se  
faire par un \cari de taille \,$S(n)+n$\,  
et de \prof \,$D(n)+1$. \qed  
  
\ss
Ce r\'esultat et l'\algo \ref{algoTFD} 
qui a introduit  
le lemme \ref{lem TFD} nous permettent d'estimer avec  
pr\'ecision la \com de l'\algo de la \mul  
rapide des \pols et d'\'enoncer le \tho suivant.  

\begin{theorem} \label{thTFD1} 
On consid\`ere un anneau \,$\A$\, poss\'edant  
une racine \,$2^{\nu+1}\,$-\,\`eme principale de $1$ et dans  
lequel \,$2_\A$\, est inversible.\\
Alors, en utilisant l'\algo \ref{algoTFD}
avec l'\eva r\'ecursive d\'ecrite dans la preuve du 
lemme  \ref{lem TFD},
la \mul de deux \pols de degr\'es \,$<n\leq 2^\nu$\,
\`a   \coes dans  \,$\A$\,  se fait \`a l'aide  
d'un \cari de taille $Ê\leq n\,(18\,\log{n}+44)$\,  
et de \prof $\leq 4\,\log{n}+10$. 
\end{theorem}  
  
\prv  
Supposons d'abord \,$n=2^\nu$. On ex\'ecute en \paral deux  
TFD d'ordre \,$2n$\, suivies d'une \'etape \paral  
avec \,$2n$\, \muls dans l'anneau de base,  
et on termine par une transformation inverse  
d'ordre \,$2n$.  La preuve du lemme  
 \ref{lem TFD} donne la majoration de la taille par
\,$9(\nu+1)\,2^\nu+4n=9\,n\,\log{n}+13\,n$\, et de la profondeur par
\,$4\,\nu+6=4\,\log{n}+6$. Dans le cas g\'en\'eral,
il faut remplacer \,$n$\, par \,$2n$\, et \,$\log{n}$\,
par \,$1+\log{n}$. \qed  

\medskip 

Rappelons que pour un anneau \,$\A$\, fix\'e par le contexte,
nous notons \,$\mu_P(n)$\,  le nombre d'\oparis \ncrs pour la \mul de 
deux \pols de degr\'e \,$n$\, en \prof \,$\O(\log{n})$. 
Le \tho pr\'ec\'edent nous dit donc qu'on a \,$\mu_P(n)=\O(n\,\log n)$\,
si  \,$2_\A$\, est inversible et si l'anneau poss\`ede des racines 
\,$2^{\nu}\,$-\,\`emes principales de l'unit\'e pour tout \,$\nu$.

\subsection[Cas d'un \acom arbitraire]  
{Algorithme de la TFD   
rapide pour un \acoma} 
\label{TFDmrGen}  

L'\algo que nous venons de d\'evelopper n'est  
pas valable lorsque  \,$2_{\A}$\, divise  
z\'ero dans l'anneau \,$\A$\, (puisque, dans un tel  
anneau, la division par $2$ ne peut pas \^etre  
d\'efinie de mani\`ere unique, m\^eme lorsqu'elle  
est possible). On peut essayer de contourner cette  
difficult\'e en rempla\c{c}ant $2$ par un entier  
\,$s\geq 2$\, tel que \,$s\,1_{\A}$\, ne divise pas  
z\'ero dans \,$\A$.  
Lorsqu'un tel entier \,$s\geq 2$\, existe, et \`a  
supposer qu'on dispose d'une racine principale  
\,$s\,$-\,\`eme de 1 dans \,$\A$, il faut encore disposer  
d'un \algo performant pour la division par \,$s$\,  
(quand elle est possible) pour pouvoir effectuer la  
transformation de Fourier inverse. En outre, un tel  
entier \,$s$\, n'existe pas \ncrtz.

Pour se d\'ebarrasser radicalement de ce probl\`eme,  
l'id\'ee de Cantor-Kaltofen dans \cite{Cant} est de  
calculer s\'epar\'ement \,$uAB$\, et \,$vAB$\, avec  
deux entiers \,$u$\, et \,$v$\, premiers entre eux,  
puis de r\'ecup\'erer \,$AB$\, en utilisant une  
relation de Bezout entre $u$ et $v$. Par exemple,  
on prend \,$u=2^\nu\geq 2n$\, et \,$v=3^\mu\geq 2n$.  
On calcule sans aucune division $2^\nu AB$ par la  
formule $$2^\nu AB=TFD_{2^\nu,\omega_{2,\nu}^{-1}}  
\,(TFD_{2^\nu,\omega_{2,\nu}}\,(A) \odot   
TFD_{2^\nu,\omega_{2,\nu}}\,(B))\,.$$
(o\`u $\omega_{2,\nu}$ est une racine  
\,$2^\nu\,$-\,\`eme principale de $1$). De m\^eme,  
on calcule \,$3^\mu AB$\, par la formule    
$$
3^\mu AB=TFD_{3^\mu,\omega_{3,\mu}^{-1}}\,  
(TFD_{3^\mu,\omega_{3,\mu}}\,(A) \odot   
TFD_{3^\mu,\omega_{3,\mu}}\,(B))
$$   
(o\`u $\omega_{3,\mu}$ est une racine  
\,$3^\mu\,$-\,\`eme principale de $1$). 

\ss Il reste n\'eanmoins un obstacle de taille,  
qui consiste en la n\'ecessit\'e de rajouter  
un substitut formel \`a \,$\omega_{2,\nu}$\, (et  
\,$\omega_{3,\mu}$) lorsqu'on ne les a pas sous  
la main dans l'anneau \,$\A$.  
Or l'id\'ee toute simple de faire les calculs  
dans l'anneau \,$\A[\lambda_{2,\nu}]$, o\`u  
\,$\lambda_{2,\nu}$\, est un susbstitut formel de  
\,$\omega_{2,\nu}$\, ne donne pas le r\'esultat  
souhait\'e. En effet, une \opari dans l'anneau  
\,$\A\,[\lambda_{2,\nu}]$\, correspond a priori  
\`a grosso modo \,$n$\, \oparis dans \,$\A$, ce  
qui annule le b\'en\'efice de la transformation  
de Fourier discr\`ete.

L'id\'ee de Cantor et Kaltofen pour r\'esoudre  
ce deuxi\`eme probl\`eme est d'appliquer une  
strat\'egie \gui{diviser pour gagner}, un  
peu semblable \`a celle du lemme \ref{lem TFD}.  
  
\ss  
La d\'efinition pr\'ecise de l'anneau  
\,$\A\,[\lambda_{2,\nu}]$\,
et la description de l'\algo font appel aux \emph{\pols  
cyclotomiques}\index{polyn\^ome!cyclotomique}, 
dont nous rappelons  
maintenant quelques propri\'et\'es.
 
Le \,$n\,$-\,\`eme \poly cyclotomique est d\'efini  
\`a partir d'une racine \,$n\,$-\,\`eme primitive
de 1, \cad un g\'en\'erateur \,$\omega_n$\, du  
groupe multiplicatif (cyclique) des racines  
\,$n\,$-\,\`emes de 1 dans une cl\^oture \agq  
de \,$\QQ$, par exemple dans \,$\CC$\, avec  
\,$\omega_n=\mathrm{e}^{i\,2\pi/n}$. 
    
Le \emph{$\,n\,$-\,\`eme \poly cyclotomique} est, par   
d\'efinition, le \poly
$$\,\Phi_n\,(X)\,=  
\prod_{\hbox{\footnotesize ${\begin{array}{c} 1\leq h<n \\  
(h,n)=1\end{array}}$}}\!
{\left(X-\omega_n^h\right)}\,.$$   
C'est un \poly unitaire \`a \coes entiers  
dont les z\'eros sont les racines \,$n\,$-\,\`emes  
primitives de 1 et dont le degr\'e est \'egal  
\`a \,$\varphi(n)$. C'est aussi un \poly  
r\'eciproque: \,$X^{\varphi\,(n)}\,\Phi_n\,(1/X)  
= \Phi_n\,(X)$. Les \pols cyclotomiques poss\`edent
en outre les  
propri\'et\'es suivantes:\\[2mm]    
$~\bullet~$ $\Phi_n\,(X)=\prod_{d|n} \Phi_d\,(X)$; \\  
\hspace*{6mm}($d|n$\, signifie que \,$d$\, est  
un diviseur positif de \,$n$) \\[1mm]  
$~\bullet~$ $\Phi_p\,(X)=X^{p-1}+\cdots+X+1$\,  
pour tout nombre premier \,$p$;\\[1mm]  
$~\bullet~$ $\Phi_{ms^k}\,(X)=  
\Phi_{ms}\,(X^{s^{k-1}})~$ si \,$k\geq 2$;\\[1mm]  
$~\bullet~$ $\Phi_m\,(X)\,\Phi_{mp}\,(X)=  
\Phi_{m}\,(X^{p})~$ si \,$p$\, premier
ne divise pas \,$m$;\\[1mm]  
$~\bullet~$ $\Phi_{2n}(X)=\Phi_n\,(-X)~$ si \,$n$\,  
est impair \,$\geq 3$.  
   
\ss   
On en d\'eduit, en particulier, que:\\[2mm]   
$~\bullet~$ $\Phi_n\,(1) \, =\,  
\formule{  
\,p~~\mbox{ si \,$n$\, est une puissance d'un  
nombre premier \,$p\,$}\\ 
\,1~~\mbox{ sinon.}  
}$ 

\ss Rajouter formellement une racine primitive  
\,$s^q\,$-\,\`eme de 1 dans \,$\A$\, revient \`a  
consid\'erer l'anneau \,$\aqo{\A[Y]}{\Phi_{s^q}(Y)}  
= \A[\lambda_{s,q}]$. Dans cet anneau, une  
addition \'equivaut \`a \,$\varphi(s^q)$\,  
additions dans \,$\A$. Pour une \mulz,  
on peut travailler dans \,$\A[Y]$\, modulo  
\,$(Y^{s^q}-1)$\, puis r\'eduire le r\'esultat  
obtenu modulo \,$\Phi_{s^q}(Y)$. Cette derni\`ere  
op\'eration est relativement peu co\^uteuse car  
\,$\Phi_{s^q}(Y)=\Phi_s(Y^{s^{q-1}})$\, est un  
\poly unitaire qui a tr\`es peu de coefficents  
non nuls. Cette remarque permet de voir que les  
\muls dans \,$\A[\lambda_{s,q}]$\, ne  
sont pas tellement plus co\^uteuses que les  
additions. Elle donne une id\'ee de comment pourra  
\^etre appliqu\'ee une strat\'egie diviser pour  
gagner, de mani\`ere \`a rendre peu  
co\^uteux les calculs dans l'anneau  
\,$\A[\lambda_{s,q}]$.
L'\algo de Cantor-Kaltofen donne alors 
le r\'esultat suivant:  

\begin{theorem}\label{thCaK} 
Il existe une \famu de \caris
 de \prof  \,$\O(\log{n})$\, qui calculent  
le produit de deux \pols de degr\'e \,$<n$\,  
\`a \coes dans  un \acoma \,$\A$\, avec  
\,$\O(n\,\log{n})$\, \muls  
et \,$\mu_P(n)=O(n\,\log{n}\,\log{\log{n}})$\,  
\adsosz.  
\end{theorem}  

\begin{remark} \label{rem CaK}  
\emph{ 
L'\algo de Cantor-Kaltofen prend en entr\'ee  
deux \pols \,$A$\, et \,$B$\, de degr\'e  
$< n$\, et donne en sortie \,$C=AB$.  
Il calcule tout d'abord \,$s_1^{q_1}\,C$\, et  
\,$s_2^{q_2}\,C$, o\`u \,$s_1$\, et \,$s_2$\,  
sont deux petits entiers premiers entre eux, et  
\,$s_1^{q_1}$\, et \,$s_2^{q_2}$\, ne sont pas  
trop grands par rapport \`a \,$n$.   
La constante cach\'ee du \gui{grand \,$\O\,$} dans  
l'estimation \,$\O(n\,\log{n}\,\log{\log{n}})$\,  
de la taille du \cir calculant $s^{q}\,C$ est  
de l'ordre de \,$4s^2\,(3s+1)$\, si \,$s$\, est  
premier. Il s'ensuit qu'en utilisant les deux  
valeurs optimales \,$s_1=2$\, et \,$s_2=3$,  
l'\algo de Cantor-Kaltofen ne devient plus  
performant que l'\algo en \,$\O(n^{\log{3}})$\,  
que pour les valeurs de \,$n$\, qui sont de  
l'ordre de \,$6\,10^4$. 
}
\end{remark}

\begin{remark} \label{rem DFTapprox} 
 \emph{ 
La \mul rapide des \pols est en fait couramment
uti\-lis\'ee en analyse num\'erique, en prenant des
approximations num\'eriques des racines de l'unit\'e dans 
\,$\CC$. Cela laisse supposer qu'une impl\'ementation efficace
de cette \mul rapide est \egmt possible en calcul formel
avec des anneaux tels que \,$\ZZ$\, ou un anneau de \pols 
sur \,$\ZZ$. \emph{Il suffit} en effet de faire le calcul 
num\'erique approch\'e avec une pr\'ecision
suffisante pour que le r\'esultat du calcul soit garanti 
avec une pr\'ecision meilleure que \,$1/2$.
Une autre solution voisine, mais o\`u la pr\'ecision est 
plus facile \`a contr\^{o}ler, serait de faire un
calcul num\'erique approch\'e non dans
\,$\CC$\, mais dans un anneau d'entiers
\,$p$\,-\,adiques (voir par exemple \cite{Serre}): un tel
anneau contient une racine primitive
\,$(p-1)$\,-\,\`eme de l'unit\'e, et  
\,$(p-1)$\, y est inversible.
}
\end{remark}

\section{Produits de \mtos}\label{secMTT}

Nous signalons ici une interpr\'etation matricielle du 
produit de deux \pols  \,$A$\, et \,$B$\, de degr\'es \,$m$\, et \,$n$.
On consid\`ere le \Amo libre \,$P_{m+n+1}\simeq\A^{m+n+1}$\, 
des \pols de degr\'e \,$\leq m+n$\, muni de la base canonique
des mon\^{o}mes \,$X^k$. La multiplication par \,$A$\, 
(resp.
\,$B$, resp. \,$AB$) tronqu\'ee au degr\'e \,$m+n$\, est 
repr\'esent\'ee sur cette base par une \mtt \,$T_A$\,
(resp. \,$T_B$, resp. \,$T_{AB}$) et on a 
\,$T_A\,T_B=T_{AB}$.
Par exemple avec \,$m=3,\,n=2$\, on obtient le produit
$$
\left[ 
{\begin{array}{cccccc}
a_0 & 0 & 0 & 0 & 0 & 0 \\
a_1 & a_0 & 0 & 0 & 0 & 0 \\
a_2 & a_1 & a_0 & 0 & 0 & 0 \\
a_3 & a_2 & a_1 & a_0 & 0 & 0 \\
0 & a_3 & a_2 & a_1 & a_0 & 0 \\ 
0 & 0 & a_3 & a_2 & a_1 & a_0  
\end{array}}
 \right] 
 \,\times\, 
 \left[ 
{\begin{array}{cccccc}
b_0 & 0 & 0 & 0 & 0 & 0  \\
b_1 & b_0 & 0 & 0 & 0 & 0 \\
b_2 & b_1 & b_0 & 0 & 0& 0 \\
0 & b_2 & b_1 & b_0 & 0& 0 \\
0 & 0 & b_2 & b_1 & b_0& 0 \\
0 & 0 & 0 & b_2 & b_1 & b_0 
\end{array}}
 \right] 
$$
qui est \'egal \`a la \mtti dont la premi\`ere colonne est 
donn\'ee
par les \coes du produit \,$AB$:
$$
 \left[ 
{\begin{array}{cccccc}
a_0 & 0 & 0 & 0 & 0 & 0 \\
a_1 & a_0 & 0 & 0 & 0 & 0 \\
a_2 & a_1 & a_0 & 0 & 0 & 0 \\
a_3 & a_2 & a_1 & a_0 & 0 & 0 \\
0 & a_3 & a_2 & a_1 & a_0 & 0 \\ 
0 &  0 &a_3 & a_2 & a_1 & a_0  
\end{array}}
 \right] 
 \,\times\, 
 \left[ 
{\begin{array}{c}
b_0  \\
b_1  \\
b_2  \\
0   \\
0   \\
0 
\end{array}}
 \right] 
 \,=\,
 \left[ 
{\begin{array}{ccccc}
a_0\,b_0  \\
a_1\,b_0 + a_0\,b_1  \\
a_2\,b_0 + a_1\,b_1 + a_0\,b_2  \\
a_3\,b_0 + a_2\,b_1 + a_1\,b_2  \\
a_3\,b_1 + a_2\,b_2  \\
a_3\,b_2 
\end{array}}
 \right]
$$

Inversement, le produit de deux \mttis dans \,$\A^{n\times 
n}$\,
peut s'interpr\'eter
comme le produit de deux \pols
de degr\'es $\leq n-1$, tronqu\'e au degr\'e \,$n-1$
(\cad encore comme le produit dans l'anneau des \dlis
\,$\A_{n-1}=\aqo{\A[X]}{X^{n}}$). Par exemple
$$
 \left[ 
{\begin{array}{ccccc}
a_0 & 0 & 0 & 0  \\
a_1 & a_0 & 0 & 0  \\
a_2 & a_1 & a_0 & 0 \\
a_3 & a_2 & a_1 & a_0
\end{array}}
 \right] 
 \,\times\, 
 \left[ 
{\begin{array}{c}
b_0  \\
b_1  \\
b_2  \\
b_3
\end{array}}
 \right] 
 \,=\,
 \left[ 
{\begin{array}{ccccc}
a_0\,b_0  \\
a_1\,b_0 + a_0\,b_1  \\
a_2\,b_0 + a_1\,b_1 + a_0\,b_2 \\ 
a_3\,b_0 + a_2\,b_1 + a_1\,b_2 + a_0\,b_3  
\end{array}}
 \right]\,.
$$

En bref il n'y a pas de diff\'erence significative
entre le produit de 2 \polsz, le produit de 2 \mttis
carr\'ees et le produit d'une \mtti par un vecteur.

\ms Voyons maintenant la question du produit d'une 
\mto arbitraire par un vecteur.
Par exemple
$$
 \left[ 
{\begin{array}{ccccc}
a_3 & a_2 & a_1 & a_0\\
a_4 & a_3 & a_2 & a_1\\
a_5 & a_4 & a_3 & a_2\\
a_6 & a_5 & a_4 & a_3\\
a_7 & a_6 & a_5 & a_4\\
a_8 & a_7 & a_6 & a_5
\end{array}}
 \right] 
 \,\times\, 
 \left[ 
{\begin{array}{c}
b_0  \\
b_1  \\
b_2  \\
b_3
\end{array}}
 \right] 
 \,=\,
 \left[ 
{\begin{array}{c}
c_3\\
c_4\\
c_5\\
c_6\\
c_7\\
c_8
\end{array}}
 \right]\,.
$$
Il suffit d'ins\'erer
la premi\`ere matrice dans la matrice de la
multiplication par le \pol \,$A=\sum_{i=0}^8a_iX^i$,
tronqu\'ee au degr\'e \,$11$,
dans le \Amo libre des \pols de degr\'es $\leq 11$:
$$
 \left[ 
{\begin{array}{ccccc}
a_0 & 0   & 0   & 0   \\
a_1 & a_0 & 0   & 0   \\
a_2 & a_1 & a_0 & 0   \\
a_3 & a_2 & a_1 & a_0 \\
a_4 & a_3 & a_2 & a_1 \\
a_5 & a_4 & a_3 & a_2 \\
a_6 & a_5 & a_4 & a_3 \\
a_7 & a_6 & a_5 & a_4 \\
a_8 & a_7 & a_6 & a_5 \\
0   & a_8 & a_7 & a_6 \\
0   & 0   & a_8 & a_7 \\
0   & 0   & 0   & a_8
\end{array}}
 \right] 
 \,\times\, 
 \left[ 
{\begin{array}{c}
b_0  \\
b_1  \\
b_2  \\
b_3
\end{array}}
 \right] 
 \,=\,
 \left[ 
{\begin{array}{c}
a_0\,b_0  \\
a_1\,b_0 + a_0\,b_1  \\
c_2\\
c_3\\
c_4\\
c_5\\
c_6\\
c_7\\
c_8\\
c_{9}\\
a_7\,b_3 + a_8\,b_2\\
a_8\,b_3
\end{array}}
 \right]\,.
$$
On voit alors que le calcul se ram\`ene au produit du
\pol \,$A$\,  par le \pol \,$B=\sum_{i=0}^3b_iX^i$.
On en d\'eduit le r\'esultat important suivant
o\`u l'on voit que le produit par une \mto n'est
gu\`ere plus cher que le produit par une matrice creuse.

\begin{prop}\label{prop mto}
Le produit d'une \mto et d'une matrice arbitraire, toutes deux carr\'ees 
d'ordre \,$n$\, peut se faire par une   \fam de \caris en 
\,$\SD(n\,\mu_P(n),\alb\log n ))$.
\end{prop}

\rem Plus pr\'ecis\'ement supposons que dans l'\acom  
\,$\A$\, la \mul d'un \pol de degr\'e $\leq n$
par un \pol de   degr\'e $\leq m$ soit en 
\,$\SD(\mu(n,m),\alb\lambda(n,m))$.
Alors le produit \,$T\,B$\, d'une \mto  
\,$T\in\A^{n\times m}$\,
par une matrice   \,$B\in\A^{m\times p}$\, est en 
\,$\SD(p\,\mu(n+m,m),\lambda(n+m,m))$.
Ceci n'est qu'un exemple des r\'esultats de \com \arit
 concernant les \mtosz.
Nous renvoyons le lecteur int\'eress\'e par le sujet \`a l'ouvrage
\cite{BP}.

\newpage \thispagestyle{empty}

\chapter{Multiplication rapide des matrices}  
\label{chap multimat}
\minitoc
\acvide

\subsubsection*{Introduction}

La \mul des matrices \`a
\coes dans un \acom unitaire
\,$\A$\, a fait l'objet de multiples
investigations durant les trente derni\`eres
ann\'ees en vue de r\'eduire le nombre
d'\oparis (dans \,$\A$)
\ncrs au calcul du produit d'une matrice
\,$m\times n$\, par une matrice \,$n\times p$,
et d'am\'eliorer la borne sup\'erieure
asymptotique de ce nombre. Il s'est
av\'er\'e que c'est le nombre de \muls essentielles
qui contr\^ole la \com
asymptotique de la \mul des \macasz, comme 
nous allons le voir
tout d'abord \`a travers l'\algo de la \mul rapide
de Strassen.

L'\algo conventionnel (dit \uslz) pour le
calcul du produit
\,$C=(c_{ij})\in \A^{m\times p}$\, d'une matrice
\,$A=(a_{ij})\in \A^{m\times n}$\, par une
matrice \,$B=(b_{ij})\in \A^{n\times p}$\, se
fait par \,$mnp$\, \muls et
\,$mp\,(n-1)$\, additions en calculant en \paral
(en une seule \'etape) les \,$mnp$\, produits
\,$a_{ik}b_{kj}$\, et en calculant ensuite, en
\paral et en \,$\esup{\log{n}}$\,
\'etapes, les \,$mp$\, sommes \,$c_{ij}$\,
intervenant dans les formules
$$c_{ij} =  \som_{k=1}^n a_{ik}b_{kj}~\mbox{ pour
}~1\leq
i\leq m~\mbox{ et }~1\leq j\leq p\,.$$

En particulier pour la \mul de deux
\macas d'ordre
\,$n$, cet \algo correspond \`a un
\cari de taille
\,$n^2(2n-1)$\, et de \prof
\,$\esup{\log{n}} + 1$\,
avec \,$n^3$\, \muls et \,$n^2(n-1)$\,
additions.

Dans un premier temps, les investigations
portaient sur la diminution du nombre de
\muls en essayant d'y r\'eduire le
coef\-ficient de \,$n^3$\, sans s'occuper de
l'exposant de \,$n$, et c'est Winograd
qui r\'eussit le premier \`a r\'eduire ce
\coe de moiti\'e, mais en doublant presque
le nombre d'additions, ce qui constitue, malgr\'e
ce prix, un progr\`es dans la \com
asymptotique si l'on sait que dans une large classe
d'anneaux la \mul est beaucoup plus
co\^uteuse que l'addition\footnote{~Ce qui n'est
pas vrai par exemple dans le corps des fractions
rationnelles \,$\QQ\,(X)$\, o\`u
l'addition est plus co\^uteuse.}.
   Beaucoup pensaient
que ce r\'esultat de Winograd
serait optimal au sens que \,${1\over 2}n^3$\,
\muls seraient \ncrs pour le
calcul du produit de deux matrices
\,$n\times n$\, (voir \cite{Knu}, page 481).

Mais une ann\'ee plus tard (1969), Strassen montra
que l'on pouvait multiplier
deux matrices \,$n\times n$\, en utilisant seulement
\,$\O(n^{2,8})$\, \mulsz. Ce r\'esultat
\'etait bas\'e sur le fait tr\`es simple
que le produit de
deux matrices \,$2\times 2$\, \`a
\coes dans un anneau \emph{non \ncrt
commutatif\/} pouvait \^etre calcul\'e avec seulement
7 \muls au lieu de 8, le nombre
d'additions passant de 4 \`a 18, et il donna les
relations prouvant ce fait dans son fameux article
\,\emph{Gaussian elimination is not optimal\,}
\cite{Str}. Winograd donna un peu plus tard
\cite{Winog} une variante de la \mul
rapide de Strassen avec seulement 15 additions.

Comme ces relations n'utilisent pas la commutativit\'e
de la \mulz, elles s'appliquent
r\'ecursivement au calcul du produit de deux matrices
quelconques \`a \coes dans \,$\A$\, selon
la strat\'egie \gui{diviser pour gagner}.

\ms La section \ref{sec strass} est consacr\'ee \`a
une
analyse d\'etaill\'ee de la \mul rapide
des matrices dans la version Strassen-Winograd.
Nous \'etudions \egmt l'uniformit\'e de la
construction 
de la famille de \caris qui correspond \`a la version 
originale de Strassen, comme annonc\'e dans la section
\ref{subsec uniforme}.

\ss Dans la section \ref{invtrian}  nous montrons que 
l'inversion des \matgs  \fregs peut \^etre
r\'ealis\'ee 
par des \caris avec une taille de m\^eme ordre que les
circuits de la \mul des \macas et une \prof d'ordre 
$\,\O(\log^2{n})\,$ au lieu de $\,\O(\log{n})\,$.

\ss Dans la section \ref{sec combili}  nous
introduisons
les notions de \com\bilz, de \com\muv et de \rgtez.
Nous montrons le r\^{o}le central jou\'e par la notion
de \rgte
dans la \com asymptotique de la \mul des \macas (\tho
\ref{thRgExp} d\^u \`a Strassen).
Nous montrons \egmt le r\'esultat de Sch\"onhage,
qui dit que l'exposant de la \mul des
\macas ne d\'epend que de la \cara du corps de base
(on conjecture en fait que cet exposant est le m\^eme
pour tous les
corps et pour l'anneau des entiers relatifs.)

\ss Dans la section \ref{sec bini} nous nous attaquons
\`a des
\algos nettement plus sophistiqu\'es qui s'appuient
sur la notion
de \cabaz, introduite par Bini.
Malgr\'e leurs performances asymptotiques,
aucun des \algos de cette section ne semble devoir
\^etre impl\'ement\'e sur machine dans un proche
avenir.
Il nous a pourtant sembl\'e que ce serait un crime
contre la beaut\'e
que de ne pas d\'evoiler au moins en partie les
id\'ees fascinantes
qui y sont \`a l'{\oe}uvre.
Nous n'avons cependant pas expos\'e la \gui{\met du
laser}
due \`a Strassen (cf. \cite{Bur,Stra2}), car nous
n'avons pas
vu comment en donner une id\'ee assez exacte en termes
suffisamment simples.
Cette \met a conduit \`a la meilleure borne connue
pour
l'exposant de la \mul des \macasz. L'estimation
actuelle
de cet exposant \,$\omega$\, est de $2,376$: Winograd
\& Coppersmith,
1987 (\cite{Win,Win2}). 

\section{Analyse de la m\'ethode de Strassen} \label{sec strass}
\acvide 

\subsection[La \met et sa \com]{La \met de Strassen (version Winograd) 
et sa \com}
On consid\`ere dans un anneau \,$\B$\, (non \ncrt commutatif) deux  
matrices \,$A$\, et \,$B$:  
$$A = \cmatrix{  
a_{11}    &    a_{12}       \cr  
a_{21}    &    a_{22}  
}\quad  
B=\cmatrix{  
b_{11}    &    b_{12}       \cr  
b_{21}    &    b_{22}  
}\quad \mathrm{ avec~~} C=AB=\cmatrix{  
c_{11}    &    c_{12}       \cr  
c_{21}    &    c_{22}  
}  
\,.$$  
Alors la matrice \,$C$\, peut \^etre obtenue par le  
calcul suivant:  
\begin{agh}
$$\begin{array}{ll}  
m_1:=a_{11}\,b_{11}\;\;&  
m_2:=a_{12}\,b_{21}\;\; \\    
m_3:=(a_{11}-a_{21})\,(b_{22}-b_{12}) &  
m_4:=(a_{21}+a_{22})\,(b_{12}-b_{11}) \\  
m_5:=(a_{21}+a_{22}-a_{11})\,(b_{22}-b_{12}+b_{11}) \\  
m_6:=(a_{11}+a_{12}-a_{21}-a_{22})\,b_{22} \\  
m_7:=a_{22}\,(b_{22}-b_{12}+b_{11}-b_{21}) & \;\; \\ 
c_{11}:=m_1+m_2 &  
c_{12}:=m_1 + m_5 + m_4 + m_6 \\ 
c_{21}:=m_1 + m_3 + m_5 -m_7& 
c_{22}:=m_1 + m_3 + m_4+ m_5   
\end{array}$$
\end{agh}

Ces relations de Strassen (version  
Winograd), appliqu\'ees \`a l'anneau des \macas
d'ordre $2k$,  ram\`enent le calcul du produit de   
deux matrices \,$2k \times 2k~(k \in \NN^*)$\,  
\`a celui de sept produits de matrices  
\,$k \times k$\,  
et de $15$ sommes de matrices de m\^eme type. 

L'analyse de \com faite \`a la section \ref{subsec principegene}
montre que ce passage de 8 \`a 7 \muls est un 
avantage d\'ecisif, ind\'ependamment du nombre des 
additions utilis\'ees par ailleurs. Cela tient \`a ce que
7 est le degr\'e de \parasm dans la proc\'edure \dpg tandis que le 
nombre d'additions n'intervient que dans la constante du \,$\O(n^2)$\,
\oparis \ncrs pour, partant  du \pb initial \,$P_n$, d'une part cr\'eer 
les 7 sous-\pbs de
type \,$P_{\esup{n/2}}$, et d'autre part r\'ecup\'erer la solution du 
\pb initial
\`a partir des solutions, calcul\'ees en \paralz, des 7 sous-\pbs  (cf. 
proposition \vref{propTPDQ}).
  
Posant:   
$$\cmatrix{  
A_{11} & A_{12} \cr A_{21} & A_{22} }\,   
\cmatrix{ B_{11} & B_{12} \cr B_{21}  
& B_{22} }
\,=\,
\cmatrix{
C_{11} & C_{12} \cr C_{21} &  
C_{22}   
}$$
o\`u les  
\,$A_{ij},\,B_{ij},\,C_{ij}$\, $(1\leq i,j\leq 2)$  
sont des matrices \,$k \times k$, on a un sch\'ema
de \prev comportant les instructions  suivantes 
dans lesquelles les affectations des variables  
\,$M_i~(1\leq i\leq 7)$\, correspondent aux  
7 \muls et celles des variables  
\,$N_i~(1\leq i\leq 11)$\, et   
\,$C_{ij}~~(1\leq i,j\leq 2)$\, correspondent aux 15  
\adsosz({\footnote{~Dans la suite, nous dirons simplement
additions, en sous-entendant \adsosz.}}), 
avec indication des \'etapes du calcul \paralz:   
\begin{algor}[Multiplication de matrices par blocs, \`a la 
Stras\-sen-Winograd]
    \label{MuSW}
    \acl{MuSW}{Multiplication par blocs, \`a la Strassen-Winograd}
\Debut
\Eta{1}{}{N_1 := A_{11}-A_{21}\;;\;N_2 := A_{21}+A_{22}\;;\; 
}\hsud{
N_3 := B_{12}-B_{11}\;;\; N_4 := B_{22}-B_{12}}
\Eta{2}{}{N_5 := N_2-A_{11}\;;\; N_6 := B_{22}-N_3}
\Eta{3}{}{N_7 := A_{12}-N_{5}\;;\; N_8 := N_6-B_{21}}
\Eta{4}{Les 7 multiplications}
{M_1 := A_{11}B_{11}\;;\; M_2 := A_{12}B_{21}\;;\;
M_3 := N_1N_4\;;\;M_4 := N_2N_3\;;\;}
\hsud{M_5 := N_5N_6\;;\;M_6 := N_7B_{22}\;;\;M_7 := A_{22}N_8}
\Eta{5}{}{C_{11} := M_1 + M_2\;;\;N_9 := M_1 + M_5\;;\; 
N_{10} := M_4 + M_6}
\Eta{6}{}{N_{11} := M_3 + N_9 \;;\;C_{12} := N_9 + N_{10} }
\Eta{7}{}{C_{21} := N_{11}-M_7\;;\;C_{22} := M_4+N_{11} }
\fin
\end{algor}

Appliqu\'e r\'ecursivement \`a une matrice
\,$m2^\nu \times m2^\nu~(m\in \NN^*,\nu\in \NN)$\,
ce programme donne un  \cari \paral de
taille \,$S(m2^\nu)$\, et de \prof \,$D(m2^\nu)$\,
dans l'anneau \,$\A$, v\'erifiant les relations de
\recuz\footnote{~Signalons que
pour la version originale de Strassen avec 18
additions
(cf. page \pageref{OrigStrass}), la \prof
v\'erifie la relation
\,$D(m2^\nu) = D(m2^{\nu-1}) + 3$.}:
\begin{equation}\label{strasswin}
\formule{
S(m2^\nu) = & 7\,S(m2^{\nu-1}) + 15\cdot
m^2\,4^{\nu-1} \\ D(m2^\nu) =
& D(m2^{\nu-1}) + 6 \,\cdot
}
\end{equation}
   
La derni\`ere \'equation est justifi\'ee par le fait  
que les \'etapes o\`u il n'y a que des additions de  
matrices \,$m2^{\nu - 1} \times m2^{\nu - 1}$\, ont   
une \prof \'egale \`a 1 (les \,$m^2\,4^{\nu-1}$\,  
additions correspondantes dans \,$\A$\, se faisant en  
\paralz) alors que l'\'etape comprenant les   
\muls de matrices (Etape~$4$) est de \prof  
\,$D(m2^{\nu-1})$.    
   
Utilisant l'\algo usuel pour la \mul de deux  
matrices \,$m\times m$, on peut \'ecrire  
\,$S(m)=m^2\,(2m-1)$\, et \,$D(m) = \esup{ \log{m}} + 1 $.   
    
Ce qui donne \,$D(n) = \,D(m2^\nu) = D(m)+6\nu =  
6 \esup{ \log{n} } + \esup{ \log{m} } +1 $\,  
comme r\'esultat pour la \prof du \cir   
\arith correspondant au calcul 
du produit de deux matrices \,$n\times n$\, si l'on  
prend \,$n=m2^{\nu}$   (la version originale de Strassen donne \,$D(n) = 3 \esup{  
\log{n} } + \esup{ \log{m} } + 1 $). 
  
Concernant la taille, la premi\`ere \'equation dans  
(\ref{strasswin}) donne successivement:   
   
   
%
\begin{center}   
\begin{tabular}{|r|rrl|}   
\hline   
$\;$ & $\;$ & $\;$ & $\;$ \\   
$1 \,\times $ & \,$S(2^{\nu}m)~=$ &  
$7\,S(2^{\nu-1}m)$ & $+ ~~15\cdot 4^{\nu-1}m^2$\, \\   
$7 \,\times $ & \,$S(2^{\nu-1}m)~=$ &  
$7\,S(2^{\nu-2}m)$ & $+ ~~15\cdot 4^{\nu-2}m^2$\, \\   
$\vdots~$ & $\vdots$\, & $\;$ & $~\vdots$ \\   
$7^{\nu-1} \,\times $ & \,$S(2m)~=$ &  
$7\,S(m)$ & $+ ~~15\cdot m^2$\, \\   
$7^{\nu} \,\times $ & \,$S(m)~=$ &  
$m^2\,(2m-1)$\, & $\;$ \\   
$\;$ & $\;$ & $\;$ & $\;$ \\   
\hline   
$\;$ & $\;$ & $\;$ & $\;$ \\   
$\dans$ & $S(2^{\nu}m)~=$ & $7^{\nu}m^2\,(2m-1)$ &  
$+~~5\,m^2\,(7^{\nu}- 4^{\nu})~~$ \\   
$\;$ & $\;$ & $\;$ & $\;$ \\   
\hline   
\end{tabular}   
\end{center}   
   
\ss Ce qui donne comme r\'esultat \,$S(m2^\nu) =  
7^{\nu}m^2\,(2m+4) - 5m^24^\nu$\, pour la taille du  
\cari correspondant au calcul par la \met  
de Strassen (variante Winograd)\footnote{~La version  
originale de Strassen donne \,$S(m2^\nu) =  
7^{\nu}m^2\,(2m+5) - 6m^24^\nu$.} du produit    
de deux matrices \,$m2^\nu \times m2^\nu$. Ainsi:   
\begin{equation}\label{strassres}   
\formule{  
 S(m2^\nu) = &  
2m^2\,(m+2)\,7^{\nu} - 5m^24^\nu \\ D(m2^\nu) = &  
6 \esup{ \log{n} } + \esup{ \log{m} } +1    
}   
\end{equation}   
   
En particulier, si \,$n$\, est une puissance de 2  
(\cad \,$m=1$), et comme \,$7^{\nu}= 2^{\nu\log{7}}$:   
$$S(n) = 6\,n^{\log{7}} - 5\,n^2~~\mathrm{et}~~ D(n) = 6  
\esup{ \log{n} }$$ (on obtient \,$3 \esup{\log{n}} +1 $\, 
seulement pour la version  originale de Strassen).   

\ms Mais le \coe de \,$n^{\log{7}}\simeq 
n^{2.807}\,$({\footnote{~$\log{7}\simeq 
2.8073549220576041074$.}}) 
dans \,$S(n)$\, peut \^etre ramen\'e \`a \,$4,15$\, lorsque
\,$n$\, est une puissance de 2. En effet, si \,$n=32$\, on 
peut 
v\'erifier directement que le nombre d'\oparis dans la \mul 
usuelle des matrices \,$n^2(2n-1)$\, ne d\'epasse gu\`ere 
\,$3,\!9\,n^{\log{7}}$\, et pour \,$n\geq 32$, on pose 
\,$\log{n} 
=\nu + 5 \geq 5$\,, de sorte que \,$n = 32.\,2^\nu\, 
~(m=32)\,$. 

\ss La premi\`ere des \'equations \,(\ref{strassres})\, donne 
alors: $$\begin{array}{llll} S(n) & = & S(m2^\nu) & \; \\ \; & 
< & 
2m^2\,(m+2)\,7^{\nu} & \; \\ \; & \leq & 2^{11}*34*7^{\nu} \\ 
\; & 
\leq & 2^{11}*34*(1/7)^5*7^{\log{n}} & ~(\mbox{puisque }~7^\nu 
= 
7^{\log{n}-5})\\ \; & \leq & 4,15\,n^{\log{7}}\,. & ~(\mbox{on 
remplace $7^{\log{n}}$ par $n^{\log{7}}$}) 
\end{array} $$

Ceci conduit donc au r\'esultat suivant d\^u \`a Strassen, 
mais dans lequel nous int\'egrons la version (avec 15 
additions) de Winograd: 
\begin{theorem}\label{thStrassMuma} 
La \mul de deux matrices \,$n\times n$\, \`a \coes  dans un 
anneau 
arbitraire \,$\A$\, est dans la classe $\SD\,(n^{\log{7}},\log 
n)$. Plus pr\'ecis\'ement, lorsque \,$n$\, est une puissance 
de 2, elle se fait soit avec un \cari dont la taille et la 
\prof sont respectivement major\'ees par 
\,$4,15\,n^{\log{7}}\,$ et \,$6\esup{\log{n}}$, soit par un 
\cir dont la taille et la \prof sont respectivement major\'ees 
par 
\,$4,61\,n^{\log{7}}$ et \,$3\esup{\log{n}}$. 
\end{theorem}

Notez aussi que la \prof \muv de ces \cirs est \'egale \`a $1$.
En fait la conclusion dans le \tho pr\'ec\'edent est 
non seulement qu'il existe une famille de \caris dans la classe 
\,$\SD(n^{\log{7}},\log n)$\, qui r\'ealise la \mul des \macasz, 
mais qu'on sait construire explicitement une famille \emph{uniforme} 
de tels \carisz. Ceci est l'objet du paragaphe
qui suit avec le \thoz~\ref{thMuStra0}. 

\subsection[Une famille uniforme de \caris]{Un exemple de construction 
uniforme d'une \fam de \caris}\label{subsecStraUnif}

Nous allons maintenant tenir une promesse que nous avions faite
dans la section \ref{sec unif}. Celle d'analyser un exemple de 
construction r\'ecursive uniforme typique d'une \fam de \caris 
pour laquelle le co\^ut de production d'un \cir de la famille 
n'a pas un ordre de grandeur bien sup\'erieur \`a sa taille.
Nous utiliserons pour cet exemple la \mul rapide des matrices  
originale  
de Strassen \cite{Str}  
qui repose sur le calcul suivant. On consid\`ere dans  
un anneau \,$\B$\, (non \ncrt commutatif) deux  
matrices \,$A$\, et \,$B\,$:  
$$A = \cmatrix{  
a_{11}    &    a_{12}       \cr  
a_{21}    &    a_{22}  
}\quad  
B=\cmatrix{  
b_{11}    &    b_{12}       \cr  
b_{21}    &    b_{22}  
}\quad \mathrm{ avec~~} C=AB=\cmatrix{  
c_{11}    &    c_{12}       \cr  
c_{21}    &    c_{22}  
}  
\,.$$  
Alors la matrice \,$C$\, peut \^etre obtenue par le  
calcul suivant, qui n\'eces\-si\-te 18 \refstepcounter{bidon}
\label{OrigStrass}
\adsos et 7 \mulsz:  
\begin{agh}
$$\begin{array}{ll}  
m_1:=(a_{12}-a_{22})\,(b_{21}+b_{22})\qquad&  
m_2:=(a_{11}+a_{22})\,(b_{11}+b_{22}) \\    
m_3:=(a_{11}-a_{21})\,(b_{11}+b_{12}) &  
m_4:=(a_{11}+a_{12})\,b_{22} \\  
m_5:=a_{11}\,(b_{12}-b_{22}) &  
m_6:=a_{22}\,(b_{21}-b_{11}) \\  
m_7:=(a_{21}+a_{22})\,b_{11} & \;\; \\ 
c_{11}:=m_1+m_2-m_4+m_6 &  
c_{12}:=m_4+m_5 \\ 
c_{21}:=m_6+m_7 &  c_{22}:=m_2-m_3+m_5-m_7   
\end{array}$$  
\end{agh}
  
Ceci peut \^etre r\'e\'ecrit sous forme d'un \cari  
de \prof 4. Concernant les variables en entr\'ee,  
on note \,$x_{0,i,j}$\, pour \,$a_{ij}$\, et  
\,$x_{0,2+i,2+j}$\, pour \,$b_{ij}$.   
On obtient le \prev    \vref{EvaMumaStrass}, que nous appelons 
\,$P_1$\,.  

\begin{proeva}[$P_1$ : produit de deux \macas d'ordre $2$ sur
un anneau non \ncrt commutatif, 
\`a la Strassen.]\label{EvaMumaStrass}
\acl{EvaMumaStrass}{Produit \`a la Strassen de deux \macas d'ordre 
$2$}
\Entree Les $8$ \coes $x_{0,ij}$ dans $\A$ (un anneau arbitraire) 
de deux \macas \,$A$\, et \,$B$\, d'ordre 
$2$.
\Sortie Les \coes $x_{4,ij}$ du produit: $C  =  A\,B$.
\Debut
\Pro{1}{}{x_{1,1}:=x_{0,12}-x_{0,22} \;;\;  
x_{1,2}:=x_{0,11}+x_{0,22} \;;\;}
\hsud{    
x_{1,3}:=x_{0,11}-x_{0,21} \;;\;  
x_{1,4}:=x_{0,11}+x_{0,12} \;;\;}
\hsud{  
x_{1,5}:=x_{0,11} \;;\;  
x_{1,6}:=x_{0,22} \;;\; 
x_{1,7}:=x_{0,21}+x_{0,22} \;;\;}
\hsud{  
x_{1,8}:=x_{0,43}+x_{0,44} \;;\;  
x_{1,9}:=x_{0,33}+x_{0,44} \;;\;}
\hsud{    
x_{1,10}:=x_{0,33}+x_{0,34} \;;\;   
x_{1,11}:=x_{0,44} \;;\;}
\hsud{  
x_{1,12}:=x_{0,34}-x_{0,44} \;;\;  
x_{1,13}:=x_{0,43}-x_{0,33} \;;\;  
x_{1,14}:=x_{0,33}  }
\Pro{2}{Les 7 multiplications}{ 
x_{2,1}:=x_{1,1}\,x_{1,8}    \;;\; x_{2,2}:=x_{1,2}\,x_{1,9}  \;;\;  
x_{2,3}:=x_{1,3}\,x_{1,10}   \;;\; }
\hsud{x_{2,4}:=x_{1,4}\,x_{1,11}\;;\;  
x_{2,5}:=x_{1,5}\,x_{1,12}   \;;\;  
x_{2,6}:=x_{1,6}\,x_{1,13}   \;;\; }
\hsud{x_{2,7}:=x_{1,7}\,x_{1,14}   } 
\Pro{3}{}{ x_{3,1} := x_{2,1}+x_{2,2} \;;\; x_{3,2} := x_{2,4}-x_{2,6} 
\;;\;}
\hsud{  
x_{3,3} := x_{2,2}-x_{2,3} \;;\; x_{3,4} := x_{2,5}-x_{2,7}}
\Pro{4}{}{
x_{4,11}:=x_{3,1}-x_{3,2} \;;\;  x_{4,12}:=x_{2,4}+x_{2,5} \;;\;}
\hsud{  
x_{4,21}:=x_{2,6}+x_{2,7} \;;\;  x_{4,22}:=x_{3,3}+x_{3,4}
 }
\fin
\end{proeva}

La \met de Strassen consiste \`a utiliser ces formules  
de mani\`ere r\'ecursive. Si on doit multiplier des matrices  
carr\'ees \`a \,$m=2^n$\, lignes et colonnes, on les partitionne  
chacune en 4 \macas  \`a \,$2^{n-1}$\,  lignes et  
colonnes, qui jouent le r\^ole des \,$a_{ij}$\, et \,$b_{ij}$\, dans  
les formules pr\'ec\'edentes. On obtient en d\'efinitive un  
\cari de \prof $3n+1=3\log(m)+1$ comportant \,$6m^2$\,  
additions/soustractions et \,$7^n=m^{\log(7)}$\, \muls  
(la \met \usle donne un \cir de \prof \,$1+n$\,  
comportant $m^3$ additions/soustractions et \,$8^n=m^3$\,  
\mulsz). Notre \pb est de d\'eterminer  
la \com en temps pour l'\'ecriture du \prev  
correspondant.

Supposons qu'on ait \'ecrit le \prev \,$P_n$\, pour la  
\mul de deux \macas \`a \,$m=2^n$\,  
lignes et colonnes, avec les entr\'ees \,$x_{0,ij}$\, et  
\,$x_{0,2^n+i,2^n+j}$\, avec \,$1\le i,j\le 2^n$, et les sorties  
\,$x_{3n+1,ij}$ ($1\le i,j\le 2^n$). 
  
Comment \'ecrit-on le \prev \,$P_{n+1}$~? 

Les entr\'ees sont maintenant $x_{0,i,j}$ et  
$x_{0,2^{n+1}+i,2^{n+1}+j}$ avec  
$1\le i,j\le 2^{n+1}$. 
Notons $X_{0,uv}$ ($1\le u,v \le 2$ ou $3\le u,v \le 4$)   
les matrices extraites \`a \,$2^n$\, lignes et colonnes,  
avec \,$m=2^n$\, et \,$m^2=2^{n+1}$:  
$$\begin{array}{llll}  
X_{0,11}[i,j]:=x_{0,i,j} &\;\;& X_{0,12}[i,j]:=x_{0,i,m+j} \\  
X_{0,21}[i,j]:=x_{0,m+i,j} && X_{0,22}[i,j]:=x_{0,m+i,m+j}\\  
X_{0,33}[i,j]:=x_{0,m^2+i,m^2+j} && 
X_{0,34}[i,j]:=x_{0,m^2+i,m^2+m+j}\\  
X_{0,43}[i,j]:=x_{0,m^2+m+i,m^2+j} && 
X_{0,44}[i,j]:=x_{0,m^2+m+i,m^2+m+j}\\  
\end{array}$$
On commence par cr\'eer (conform\'ement au programme  
$P_1$ appliqu\'e aux matrices $X_{0,uv}$) les \gui{matrices  
$X_{1,k}$}  pour $1\le k\le 14$, au moyen des affectations  
matricielles:  
$$\begin{array}{lll}  
X_{1,1}:=X_{0,12}-X_{0,22} \quad \quad &\quad \quad & 
X_{1,2}:=X_{0,11}+X_{0,22} \\  
X_{1,3}:=X_{0,11}-X_{0,21} & & X_{1,4}:=X_{0,11}+X_{0,12}\\  
X_{1,5}:=X_{0,11} & & X_{1,6}:=X_{0,22} \\  
X_{1,7}:=X_{0,21}+X_{0,12} & & X_{1,8}:=X_{0,34}+X_{0,44} \\  
X_{1,9}:=X_{0,33}+X_{0,44} && X_{1,10}:= X_{0,33}+X_{0,34} \\  
X_{1,11}:=X_{0,44} && X_{1,12}:=X_{0,34}-X_{0,44} \\  
X_{1,13}:=X_{0,43}-X_{0,33} && X_{1,14}:=X_{0,33}\\  
\end{array}$$  

\ms Cela signifie pr\'ecis\'ement dans l'anneau de base \,$\B$,   
avec \,$X_{1,k}[i,j]=\alb x_{1,k,i,j}$\, pour \,$1\le i,j\le m=2^n$:   

\noi \textsf{
\Pro{1}{}{x_{1,1,i,j}:=x_{0,i,m+j}-x_{0,m+i,m+j}}
\hsud{x_{1,2,i,j}:=x_{0,i,j}+x_{0,m+i,m+j} }
\hsud{x_{1,3,i,j}:=x_{0,i,j}-x_{0,m+i,j} }
\hsud{x_{1,4,i,j}:=x_{0,i,j}+x_{0,i,m+j}}
\hsud{x_{1,5,i,j}:=x_{0,i,j} }
\hsud{x_{1,6,i,j}:=x_{0,m+i,m+j} }
\hsud{x_{1,7,i,j}:=x_{0,m+i,j}+x_{0,i,m+j} }
\hsud{x_{1,8,i,j}:=x_{0,m^2+i,m^2+m+j}+x_{0,m^2+m+i,m^2+m+j} }
\hsud{x_{1,9,i,j}:=x_{0,m^2+i,m^2+j}+x_{0,m^2+m+i,m^2+m+j} }
\hsud{x_{1,10,i,j}:= x_{0,m^2+i,m^2+j}+x_{0,m^2+i,m^2+m+j}}
\hsud{x_{1,11,i,j}:=x_{0,m^2+m+i,m^2+m+j} }
\hsud{x_{1,12,i,j}:=x_{0,m^2+i,m^2+m+j}-x_{0,m^2+m+i,m^2+m+j} }
\hsud{x_{1,13,i,j}:=x_{0,m^2+m+i,m^2+j}-x_{0,m^2+i,m^2+j} }
\hsud{x_{1,14,i,j}:=x_{0,m^2+i,m^2+j}}
}

\ms Ensuite on cr\'ee les \gui{matrices \,$X_{2,k}\,$} pour  
\,$1\le k\le 7$. Pour cela il s'agit d'\'ecrire 7 fois,  
avec \`a chaque fois une renum\'erotation convenable,  
le programme \,$P_n$.  
Pour \,$k$\, de 1 \`a 7, on r\'e\'ecrit \,$P_n$\,  
avec les transformations suivantes:
\begin{itemize}
\item  les variables d'entr\'ees \,$x_{0,i,j}$\,  
($1\le i,j \le m$) sont remplac\'ees par les  
variables $x_{1,k,i,j}$,
\item les variables d'entr\'ees \,$x_{0,m+i,m+j}$\,  
($1\le i,j \le m$) sont remplac\'ees par les  
variables $x_{1,7+k,i,j}$,
\item  toute variable \,$x_{p,u}$\, dans \,$P_n$\, avec  
une \prof \,$p\ge 1$\, est remplac\'ee par  
la va\-ria\-ble \,$x_{p+1,k,u}$.
\end{itemize}

En particulier, on obtient en sortie les variables,
de \prof \,$3n+2$,  
$x_{3n+2,k,i,j}$ ($1\le i,j \le m$) qui sont les \coes  
des matrices $X_{2,k}$ ($1\le k\le 7$).

Il reste enfin \`a r\'ealiser les affectations  matricielles:
\textsf{
\Pro{3}{}{  
X_{3,1} := X_{2,1}+X_{2,2};\;\;  X_{3,2} := X_{2,4}-X_{2,6} }
\hsud{X_{3,3} := X_{2,2}-X_{2,3} ;\;\; X_{3,4} := X_{2,5}-X_{2,7} }
\Pro{4}{}{  
X_{4,11}:=X_{3,1}-X_{3,2} ;\;\; X_{4,12}:=X_{2,4}+X_{2,5} }
\hsud{X_{4,21}:=X_{2,1}+X_{2,7} ;\;\; X_{4,22}:=X_{3,3}+X_{3,4}
}}

\ms Cela signifie pr\'ecis\'ement, avec \,$1\le i,j \le m=2^n$:  
 \textsf{
\Pro{3n+3}{}{  
      x_{3n+3,1,i,j} := x_{3n+2,1,i,j}+x_{3n+2,2,i,j}  ;\;}
\hsud{x_{3n+3,2,i,j}:= x_{3n+2,4,i,j}-x_{3n+2,6,i,j} ;\;}
\hsud{x_{3n+3,3,i,j} := x_{3n+2,2,i,j}-x_{3n+2,3,i,j} ;\;}
\hsud{x_{3n+3,4,i,j}:= x_{3n+2,5,i,j}-x_{3n+2,7,i,j}}  
\Pro{3n+4=3(n+1)+1}{}{  
      x_{3n+4,11,i,j}:=x_{3n+3,1,i,j}-x_{3n+3,2,i,j} ;\;}
\hsud{x_{3n+4,12,i,j}:=x_{3n+2,4,i,j}+x_{3n+2,5,i,j} ;\;}
\hsud{x_{3n+4,21,i,j}:=x_{3n+2,1,i,j}+x_{3n+2,7,i,j} ;\;}
\hsud{x_{3n+4,22,i,j}:=x_{3n+3,3,i,j}+x_{3n+3,4,i,j}  
}}

\ms Le \prg qui, pour l'entr\'ee \,$n$\, donne  
en sortie le texte du \prev \,$P_n$\, est un  
\prg du type \gui{loop program} (ou \emph{ 
\prg \`a boucles pour}: \gui{\textsf{\pour{u}{1}{r}}\ldots }) 
de structure simple.  
Lorsqu'on le r\'ealise sous forme d'une machine  
de Turing \'ecrivant le texte \,$P_n$, la gestion des  
boucles occupe un temps n\'egligeable par rapport aux  
instructions qui permettent d'\'ecrire successivement  
\,$P_1$, \,$P_2,\,\ldots,\,P_n$. Il faut pr\'evoir que,  
\`a la fin de l'\'etape \num$\,i$, le texte \,$P_i$\,  
doit \^etre recopi\'e sur une bande o\`u il sera lu  
pendant l'\'etape $i+1$, car durant cette \'etape,  
la premi\`ere bande o\`u a \'et\'e \'ecrite \,$P_i$\,  
sera effac\'ee par l'\'ecriture de \,$P_{i+1}$.  
Si \,$t(n)$\, est le temps d'ex\'ecution pour  
l'\'ecriture de \,$P_n$\, et \,$s(n)$\, la taille  
de \,$P_n$, on obtient les formules r\'ecurrentes suivantes,
o\`u les \,$c_{i}$\, sont des constantes:  
$$\begin{array}{rclll} 
s(n+1)&  \le & c_0\ n\ m^2+ 7\  s(n)   \quad\quad   \mathrm{ et}  
\\[1mm] 
t(n+1)& \le  & c_1\ n\ m^2+t(n)+c_2\ s(n)+c_3\ s(n+1)\,,  &   &   \\
\end{array}$$
d'o\`u, puisque $nm^2=n2^{2n}=n4^n$ est n\'egligeable  
devant $7s(n)\ge 7^n$,  
$$s(n)=\O(7^n) \quad \mathrm{ et} \quad t(n)=\O(7^n)\,.$$  

Nous pouvons r\'esumer comme suit (rappelons que nous notons
\,$\log k$\, pour \,$\max(\log_2{k}, 1)$). 

\begin{theorem} 
\label{thMuStra0} Lorsqu'on utilise la \met 
r\'ecursive de Strassen pour construire une \fam de
\caris pour la \mul des \macas d'ordre \,$m=2^n$,
on peut construire une machine de Turing qui \'ecrit le  
code du \prev $Q_m=P_n$  en un temps   
du m\^eme ordre de grandeur que la taille de sa sortie:
 \,$\O(m^{\log 7})$.
\end{theorem}

Naturellement, comme d'habitude le
r\'esultat sur le temps de calcul \,$\O(m^{\log 7})$\, est encore
valable lorsque \,$m$\, n'est pas une puissance de \,$2$, 
en compl\'etant les matrices dans
\,$\A^{m\times m}$\, par des lignes et colonnes de $0$.

\section{Inversion des matrices triangulaires} 
\label{invtrian}

Les notations que nous pr\'ecisons maintenant
concernant
la \mul des \macas seront utilis\'ees dans
toute la suite de l'ouvrage quand nous aurons \`a
faire des
calculs de \comz.

\begin{nota} \label{IDNConst}
Nous supposerons que le calcul du produit
de deux matrices \,$n\times n$\, se fait
par un  \cari de taille
\,$\muM(n)=\Ca\,n^{\alpha}\,$ de \prof
\,$\gaM(n)=\Ka\log{n}$\,
et de largeur
\,$\laM(n) = \La\,n^\alpha\,/\,\log{n}$\,
o\`u \,$2< \alpha \leq 3$,
  \,$\Ka$\, et
\,$\La$\, sont des constantes r\'eelles
positives~$\,\geq 1$
et \,$\Ca\geq 3$~{\rm(}{\footnote{~Certains calculs de
\com dans la
suite de l'ouvrage conduiraient \`a des formules
l\'eg\`erement
diff\'erentes pour les cas \,$\alpha>2$\, et \,$\alpha
=2$. C'est la
raison pour laquelle nous avons pr\'ef\'er\'e exclure
cette derni\`ere
valeur, qui n'est de toute mani\`ere pas
d'actualit\'e.
L'hypoth\`ese \,$\Ca\geq 3$\, qui est v\'erifi\'ee
pour la
\mul rapide de Strassen et pour toutes les autres
\muls rapides connues,
n'est pas non plus restrictive et simplifie quelques
calculs.}}{\rm)}.
\indexnota{ga@$\gaM(n)=\Ka\log{n}$}\indexnota{Ka@$\Ka$}
\indexnota{muM@$\muM(n)=\Ca\,n^{\alpha}$}
\indexnota{Ca@$\Ca$}
\indexnota{laM@$\laM(n) = \La\,n^\alpha\,/\,\log{n}$}
\indexnota{La@$\La$}
\end{nota}

L'approche \dpg
donne un \algo qui montre que le probl\`eme de
l'inversion d'une
\matg inversible (autrement dit, \fregz)   admet une
solution
en \,$\SD(n^{\alpha},\alb\log^2{n})$\, avec une
constante asymptotique de l'ordre de \,$4\,\Ca$\,
pour la taille et de
l'ordre de \,$\Ka$\, pour la \prof du \cirz.

\begin{prop}\label{divtri}
Soit \,$\A$\, un anneau arbitraire, \,$n$\, un entier
$\geq 2$ et
\,$A\in \A^{n\times n}$\, une \matg  inversible. \\
Alors
l'inverse de \,$A$\, peut \^etre calcul\'ee par une
\famu de
\caris de taille \,$\tau(n)$\,  et de \prof
\,$\pi(n)$\,
v\'erifiant: \,$\tau(n) \leq 4\,\Ca\,n^\alpha$\,
 et \,$\pi(n) \leq \Ka\,\log^2{n}+\O(\log{n})\,.$
 \indexnota{Cap@$\Ca'$}
\end{prop}

\preuve
On peut toujours supposer
\,$A\in \A^{2^{\nu}\times 2^{\nu}}$\, o\`u
\,$\nu=\esup{\log{n}}$\, (\ie $2^{\nu-1} < n
\leq 2^{\nu}$) quitte \`a rajouter
\,$2^\nu\!-n$\, lignes et \,$2^\nu\!-n$\,
colonnes de z\'eros
\`a la matrice \,$A$, et remplir la partie
\gui{sud-est} restante par la matrice unit\'e
\,$\I_{2^\nu\!- n}$, ce qui revient \`a remplacer
la matrice \,$A$\, par la matrice
\,$A'=\clmatrix{ A & 0_{2^\nu\!-n,n} \\
0_{n,2^\nu\!-n} & \I_{2^\nu\!-n}
} \in \A^{2^{\nu}\times 2^{\nu}}$\,
o\`u \,$0_{p,q}$\, d\'esigne, pour tous entiers
naturels \,$p$\, et \,$q$, la matrice nulle \`a
\,$p$\, lignes et \,$q$\, colonnes.

Le calcul de \,$A^{-1}$, si \,$A$\, est inversible,
se ram\`ene \'evidemment \`a celui de \,$A'^{-1}$\,
puisque dans ce cas \,$A'$\, est inversible et
$$A'^{-1}=
\clmatrix{ A & 0_{2^\nu\!-n,n} \\
0_{n,2^\nu\!-n} &
\I_{2^\nu\!-n}
}^{-1}=
\clmatrix{
A^{-1} & 0_{2^\nu\!-n,n} \\
0_{n,2^\nu\!-n} & \I_{2^\nu\!-n}
}.
$$
Ainsi, remplac\'ee par \,$A'$, la matrice \,$A$\, peut
\^etre consid\'er\'ee comme une matrice
\,$2^{\nu}\times 2^{\nu}$\,
et s'\'ecrire (si elle est \tgiz): \\
$A = \clmatrix{
A_1 & 0 \\ A_3 & A_2
} ~~\mathrm{ o\grave{u}}
~~A_1,\,A_2,\,A_3\,\in
\A^{2^{\nu\!-1}\times 2^{\nu\!-1}}$~ avec
$~A_1,\,A_2~$ \tgisz. Donc:

\noindent \centerline{ \,$A$\,  est \freg
\,$\equiva A_1 \mbox{ et } A_2$\,
 sont \fregsz.}

\ms
De plus: $~~A^{-1} = \cmatrix{ A_1^{-1} &
0_{2^{\nu\!-1},2^{\nu\!-1}} \cr
- A_2^{-1}A_3A_1^{-1} & A_2^{-1} }.$

\ms
Le calcul de \,$A_1^{-1}$\, et \,$A_2^{-1}$\, se fait
en \paral avec un  \cari de taille
\,$\tau(2^{\nu -1})$\, et de \prof
\,$\pi(2^{\nu -1})$. On r\'ecup\`ere ensuite
le r\'esultat, \cad la matrice
\,$A^{-1}$, \`a partir de \,$A_1^{-1}$\, et
\,$A_2^{-1}$, en calculant le
produit \,$A_2^{-1}A_3A_1^{-1}$\, de trois matrices
\,$2^{\nu\!-1}\times 2^{\nu\!-1}$.

Ce qui donne les relations de \recu vraies
pour tout $\nu\geq 1$ avec
\,$\tau(1) = \pi(1) = 1$:
$$
\left\{\begin{array}{lcl} \tau(2^{\nu}) &
= & 2\,\tau(2^{\nu -1}) +
2\,\Ca\,2^{(\nu\!-1)\alpha} \\[1mm]
\pi(2^{\nu}) & = & \pi(2^{\nu -1}) + 2\,\Ka\,\nu \,.
\end{array}\right.
$$

On obtient par sommation lorsque \,$n=2^\nu$, avec
\,$a=2^{\alpha-1}$:
$$\formule{\tau(n) &=&
\left({\Ca\,n^\alpha\,-\,(\Ca+1-a)\,n}\right)\,/\,({a-1})
\leq \frac{1}{a-1}\;\Ca \,n^\alpha\\[1mm]
\pi(n) &
= &\Ka\,(\log{n}+1)\log{n} \,+\,1\,.}
$$
(ici on a utilis\'e sur la premi\`ere ligne
l'hypoth\`ese
\,$\Ca\geq 3$).

Pour le cas g\'en\'eral, on remplace \,$n$\, par
\,$2^{\esup{\log\,n}}< 2\,n$, \,$\log\,n$\, par
\,$1+\log\,n$\, et
on obtient les majorations ($2<a\leq 4$): $$ \formule{
\tau(n)&\leq & \frac{2\,a}{a-1}\;\Ca\,n^\alpha\; \leq
\;4\,\Ca\,n^\alpha\\[1mm] \pi(n)  &\leq
&\Ka\,(\log^2{n}+3\log{n}+2)+1\,. ~~~~~~~~~~~\Box}$$

  
\section{Complexit\'e bilin\'eaire} 
\label{sec combili} 

Soit un corps $\K$ et trois 
\Kevs \,$E$, \,$F$, \,$G$\, de dimensions finies. 
Rappelons 
qu'une  \emph{\abi \,$\psi:(x,y)\mapsto \psi(x,y)$\,  de \,$E\times 
F$\,  
vers \,$G\,$}\index{application!bilin\'eaire}  est une application qui 
est s\'epar\'ement \lin en \,$x$\, et en   \,$y$.
\index{bilin\'eaire!application}

\paragraph{Retour sur les \'egalit\'es de Strassen-Winograd}~
\label{paraRetSW}

\noindent  R\'e\'ecrivons les \egts de Strassen-Winograd 
donn\'ees dans la section \ref{sec strass}, sous une forme o\`u nous 
isolons les \mulsz:

$  \begin{array}{ll}
\alpha_1 := \alpha_{11}  & \beta_1 := \beta_{11}  \\
\alpha_2 := \alpha_{12}  & \beta_2 := \beta_{21}  \\
\alpha_3 := \alpha_{11}-\alpha_{21} & \beta_3 := \beta_{22}-\beta_{12} 
\\
\alpha_4 := \alpha_{21}+\alpha_{22} & \beta_4 := \beta_{12}-\beta_{11}  
\\
\alpha_5 := \alpha_{21}+\alpha_{22}-\alpha_{11} &  
\beta_5 := \beta_{22}-\beta_{12}+\beta_{11} \\   
\alpha_6 := \alpha_{12}-\alpha_{21}-\alpha_{22}+\alpha_{11}\quad  &  
\beta_6 := \beta_{22} \\
\alpha_7 := \alpha_{22}  & \beta_7 := \beta_{22}-
\beta_{12}+\beta_{11}-
\beta_{21}\\
&  \\ 
\mu_i := \alpha_i\beta_i\quad (i=1,\ldots,7) 
&  \gamma_{11} := \mu_1 + \mu_2 \\
&  \gamma_{12} := \mu_1 + \mu_4  + \mu_5 + \mu_6 \\
&  \gamma_{21} := \mu_1 + \mu_3 + \mu_5 - \mu_7 \\
&  \gamma_{22} := \mu_1 + \mu_3 + \mu_5 + \mu_4\\
\end{array} $    

\ms Ici nous avons consid\'er\'e avec trois matrices \,$A,\; B,\; C$\, 
les entr\'ees de \,$A$\, comme des formes \lins \,$\alpha_{ij}$, 
celles de 
\,$B$\, comme des formes \lins \,$\beta_{jk}$, celles de 
\,$C$\, comme des formes \lins \,$\gamma_{ik}$. Ces formes 
\lins sont d\'efinies sur l'espace des \macas d'ordre 2 sur un 
anneau $\A$.  Les 8 affectations qui d\'efinissent le produit 
\,$C:=AB$, 
$$ \gamma_{ik}:=\alpha_{i1}\beta_{1k}+\alpha_{i2}\beta_{2k}
$$
ont \'et\'e remplac\'ees par d'autres affectations, avec l'avantage de 
n'avoir que 7 \mulsz.

L'analyse de \com nous a montr\'e que ce passage de 8 \`a 7 \'etait un 
avantage d\'ecisif, ind\'ependamment du nombre des 
additions utilis\'ees par ailleurs. 

Nous avons utilis\'e 7 formes \lins \,$\alpha_\ell$\, sur l'espace 
o\`u vit la matrice \,$A$,  7 formes \lins \,$\beta_\ell$\, sur 
l'espace o\`u vit la matrice \,$B$, effectu\'e les 7 produits  
\,$\mu_\ell=\alpha_\ell\beta_\ell$\,  et r\'ecup\'er\'e les  
\,$\gamma_{ik}$\, comme \colis des \,$\mu_\ell$.

Si nous appelons \,$(c_{11},c_{12},c_{21},c_{22})$\, la base canonique 
de  l'espace o\`u vit la matrice \,$C$, nous pouvons \'ecrire
$$ C:= \mu_1 c_1+ \mu_2 c_2+\cdots+\mu_7 c_7\quad i.e., \quad 
C:= \som_{\ell=1}^7 {\; \alpha_\ell\cdot\beta_\ell\cdot c_\ell} 
$$
o\`u les \,$c_\ell$\, sont des \colis suivantes des \,$c_{ij}$:
$$\begin{array}{llll} 
c_1:=c_{11}+c_{12}+c_{21}+c_{22} & c_2:=c_{11}  &  
c_3:=c_{21}+c_{22}& c_4:=c_{12}+c_{22}   \\ 
c_5:=c_{12}+c_{21}+c_{22}  & c_6:=c_{12} & c_7:=-c_{21}
\end{array}$$

Bilan des courses: 7 formes \lins \gui{en \,$A$}, 7 formes 
\lins \gui{en \,$B$} et 7 vecteurs \gui{en \,$C$}. En math\'ematiques un 
peu plus savantes on r\'e\'ecrit ceci en utilisant la notation 
tensorielle. L'\abi \,$(A,B)\mapsto C=AB$\, correspond au tenseur 
suivant, (le premier membre de l'\egt provient directement de la 
d\'efinition)
$$ 
\som_{i,j,k\in\{1,2\}} {\alpha_{ij}\otimes\beta_{jk}\otimes 
c_{ik}}\quad 
= \quad 
\som_{\ell=1}^7 {\; \alpha_\ell\otimes \beta_\ell\otimes c_\ell}
$$
On peut consid\'erer, au choix, que ces tenseurs appartiennent \`a un 
espace tensoriel abstrait construit \`a partir des trois espaces 
\,$E$, 
\,$F$, \,$G$\, o\`u vivent les matrices \,$A$, \,$B$, \,$C$, ou bien 
qu'ils sont dans l'espace des \abis de \,$E\times F$\,  vers \,$G$. 
Dans 
ce dernier cas un \emph{tenseur 
\elrz}\index{tenseur!element@\'el\'ementaire} 
\,$\alpha\otimes \beta\otimes c$\, est \'egal par d\'efinition \`a 
l'\abi 
$$(A,B)\longmapsto \alpha(A)\cdot\beta(B)\cdot c$$  

\subsection{Rang tensoriel d'une application bilin\'eaire}
\label{paraRangTens}

Consid\'erons plus \gnlt  un corps $\K$, trois 
\Kevs \,$E$, \,$F$, \,$G$\, de dimensions finies. 
Soient
$(e_i)_{i\in I}$, $(f_j)_{j\in J}$, $(g_\ell)_{\ell\in L}$ des bases 
de 
 \,$E$, \,$F$, \,$G$\,  et notons $(e_i\sta)_{i\in I}$, 
$(f_j\sta)_{j\in J}$, $(g_\ell\sta)_{\ell\in L}$ les bases duales.
Toute \abi \,$\psi$\,  de \,$E\times F$\,  vers \,$G$\,  est alors une 
somme de tenseurs \elrsz: \,$\psi$\, est compl\`etement 
d\'etermin\'ee 
par les images qu'elle donne pour les vecteurs des bases canoniques de 
\,$E$\,  et \,$F$, et si 
\,$\psi(e_{i},f_{j})=\sum_\ell {\,\gamma_{ij\ell}\,g_\ell}$\,  on 
obtient ipso facto
$$ \psi=\som_{i,j,\ell} {\; \gamma_{ij\ell}\; e_i\sta\otimes 
f_j\sta\otimes g_\ell}
$$
 Les \,$\gamma_{ij\ell}$\, peuvent \^etre appel\'ees les 
\emph{coordonn\'ees de \,$\psi$\, sur les trois 
bases}\index{tenseur!coordonn\'ees d'un} \,$(e_i)_{i\in I}$, 
\,$(f_j)_{j\in J}$, \,$(g_\ell)_{\ell\in L}$.

L'important du point de vue du calcul sont les r\`egles de 
manipulation 
des tenseurs, qui disent qu'on a le droit d'utiliser $\otimes$ comme 
\gui{n'importe quel} produit (en utilisant la lin\'earit\'e par rapport 
\`a 
chacune des entr\'ees, l'associativit\'e, mais pas la 
commutativit\'e).

On peut par exemple supprimer les symboles $\otimes$ et calculer avec
des variables formelles \,$x_i,\; y_j,\; z_\ell$\, \`a la place des
\,$e_i\sta,\;  f_j\sta,\;  g_\ell$\,  \`a condition de ne pas 
autoriser 
la commutation de deux variables entre elles (par contre, elles 
commutent avec les \elts de $\K$). L'objet abstrait correspondant 
\`a ce calcul s'appelle \emph{l'anneau des \pols non commutatifs \`a 
\coes dans $\K$}\index{anneau des \pols non commutatifs}.

\begin{definition} 
\label{defRGtens} \emph{(Rang tensoriel d'une \abiz)}
Soient $\K$ un corps, \,$E$, \,$F$, \,$G$\,  trois \Kevs de dimension 
finie. On note \,$\Bil(E,F;G)$\, l'espace des 
\abis de \,$E\times F$\, vers \,$G$.  Soit \,$\psi\in\Bil(E,F;G)$. 
On appelle \emph{\rgte de \,$\psi\,$}\index{rang tensoriel!d'une 
application bi\-li\-n\'e\-ai\-re} le plus petit entier \,$r$\, tel que 
\,$\psi$\, puisse s'\'ecrire sous forme
$$ \som_{\ell=1}^r {\; \vep_\ell\otimes \varphi_\ell\otimes g_\ell}
$$
o\`u les \,$\vep_\ell$\, sont dans \,$E\sta$, les  \,$\varphi_\ell$\, 
sont dans \,$F\sta$\, et les \,$g_\ell$\, sont dans \,$G$. Autrement 
dit 
encore c'est le plus petit entier \,$r$\, tel que \,$\psi$\, puisse 
s'\'ecrire comme compos\'ee de trois applications selon le format 
suivant
$$ E\times F\surver{\vep\times \varphi} \K^r\times \K^r
\surver{\mu_r}\K^r\surver{g}G
$$
o\`u $\; \vep:E\rightarrow \K^r$,  $\; \varphi:F\rightarrow \K^r$  et
 $\; g:\K^r\rightarrow G$  sont des \alis et $\;\mu_r:\K^r\times 
\K^r\rightarrow \K^r$ est le produit coordonn\'ee
par coordonn\'ee.\\   
Le \prev \arith correspondant s'appelle un {\em \cabi de 
$\,\psi$}\index{bilin\'eaire!calcul}\index{calcul!bilin\'eaire}.
Le \rgte de $\,\psi\,$ est encore appel\'e la {\em \com 
\bilz}\index{complexit\'e!bilin\'eaire}\index{bilin\'eaire!complexit\'e} 
de $\,\psi$. Nous le noterons $\,R(\psi)$\indexnota{rang@$R(\psi)$, 
$R_\K(\psi)$}, ou s'il y a ambiguit\'e $\,R_\K(\psi)$.

\end{definition}
L'importance du \rgte dans les questions de \com 
\agq a \'et\'e soulign\'ee par Gastinel (\cite{Gas2}) et 
Strassen (\cite{Stra,Stra2,Stras}).
\begin{remark} 
\label{remRGtens1}
\emph{Nous laissons libre choix pour l'interpr\'etation du 
tenseur \,$\vep_\ell\otimes \varphi_\ell\otimes g_\ell$. 
Pour les gens savants cet objet vit dans un espace tensoriel abstrait 
\,$E\sta\otimes F\sta\otimes G$, 
canoniquement isomorphe \`a l'espace des \abis 
\,$E\times F\rightarrow G$. 
Mais on peut consid\'erer aussi que cet objet est \'egal 
par d\'efinition \`a l'\abi 
\,$(x,y)\mapsto\vep_\ell(x)\cdot\varphi_\ell(y)\cdot g_\ell$.   
} 
\end{remark}
\begin{remark} 
\label{remRGtens2}
\emph{Le lecteur ou la lectrice peut donner la d\'efinition analogue 
pour le \rgte d'une \aliz, ou celui d'une forme \bilz, et v\'erifier 
qu'on retrouve la notion usuelle de rang pour ces objets.   
} 
\end{remark}
\begin{remark} 
\label{remRGtens3}
\emph{ 
Nous pourrions remplacer dans la d\'efinition \ref{defRGtens} le corps 
$\K$ par un \acoma $\A$, \`a condition de consid\'erer des espaces 
convenables analogues aux \evcsz. Une possibilit\'e est de 
consid\'erer que \,$E$, \,$F$\, et \,$G$\, doivent \^etre des 
\emph{\Amos libres}, \cad des modules 
(isomorphes \`a) \,$\A^e$, \,$\A^f$\, et \,$\A^g$. Par exemple 
pour le produit matriciel, le cadre le plus naturel serait de choisir de 
travailler sans aucune hypoth\`ese pr\'ecise, \cad sur l'anneau $\ZZ$.   
} 
\end{remark}
\begin{remark} 
\label{remRGtens4}
\emph{Contrairement au rang d'une \aliz, le \rgte d'une \abi est en 
g\'en\'eral difficile \`a d\'eter\-mi\-ner. Il ne semble pas qu'on 
connaisse d'\algo qui r\'ealise ce travail, sauf pour quelques classes 
de corps particuliers (les corps finis ou les corps \agqt 
clos par exemple). Mais m\^eme dans ces cas, les \algos sont 
impraticables. La d\'etermination du \rgte des \abis sur un corps fini 
fix\'e est un probl\`eme 
$\NP$-complet (cf. \cite{Ha}): le probl\`eme, en prenant pour 
entr\'ees 
un entier \,$r$\, et une \abi \,$\psi$\, donn\'ee par ses 
coordonn\'ees 
sur trois bases, est de d\'eterminer si le \rgte de \,$\psi$\, est 
$\leq r$\, ou non.   
} 
\end{remark}

\subsubsection*{Rang tensoriel de la \mul des matrices}
\label{paraRgTMuMa}
  
La notation tensorielle n'a pas seulement l'avantage de 
l'\'el\'egance. 
Elle a aussi le m\'erite de nous aider \`a r\'efl\'echir sur les 
calculs  mis en {\oe}uvre. 
La meilleure synth\`ese de l'id\'ee de Strassen est 
peut-\^etre de dire que le miracle s'est produit quand il a pu \'ecrire 
le  tenseur de la \mul des \macas d'ordre 2 comme somme de 7 tenseurs 
\elrsz.

Chaque fois qu'on arrive \`a \'ecrire la \mul des \macas d'ordre 
\,$k$\, 
comme une somme de \,$h$\, tenseurs \elrsz, avec une bonne valeur de 
\,$\log h/\log k$\,  on obtient imm\'ediatement que la \mul des 
matrices  tombe dans la classe
\,$\SD(n^{\log h/\log k},\,\log{n})$\, car le calcul de \com de la 
section \ref{sec strass} pourra fonctionner \`a l'identique.

Plus pr\'ecis\'ement, la possibilit\'e d'effectuer des produits de 
\macas d'ordre \,$k^\ell$\, par blocs de taille  \,$k^{\ell-1}$\,  
implique r\'ecursivement que ce produit matriciel est repr\'esent\'e 
par 
une somme de \,$h^\ell$\, tenseurs \elrsz. En outre si la \prof du 
\prev correspondant au produit des matrices d'ordre 
\,$k$\, est un entier \,$K$\, alors celle du \prev 
correspondant au produit des matrices d'ordre \,$k^\ell$\, est \'egale 
\`a \,$\ell\, K$, et sa  \prof \muv (\cad la \prof mesur\'ee 
en ne tenant compte que des \muls essentielles, cf. d\'efinition 
\ref{defLongueurs} page \pageref{profmul}) est \'egale \`a \,$\ell$.

Depuis la d\'ecouverte de Strassen, un nouveau sport a \'et\'e 
cr\'e\'e, 
au\-quel ont particip\'e quelques grands noms de la \coagz: 
faire diminuer  \,$\log h/\log k$\,  en \'elaborant des identit\'es 
\agqs in\'edites, pour des valeurs de  \,$k$\, de plus en plus 
grandes.

\ms Un aspect fascinant de la notation tensorielle pour les \abis est 
qu'elle \'etablit une sym\'etrie\refstepcounter{bidon}\label{symrgte} 
entre les trois espaces 
\,$E$, \,$F$, \,$G$\, en jeu (rappelons qu'il s'agit ici d'espaces de 
matrices). Sym\'etrie qui n'est pas directement visible sur la 
d\'efinition. En fait, il n'y aurait vraiment sym\'etrie que si nous 
consid\'erions notre tenseur comme repr\'esentant l'application 
tri\lin 
$$ (A,B,C)\; \longmapsto\; \som_{i,j,k}{\; 
\alpha_{ij}\beta_{jk}\gamma_{ki}}
$$
L'\'ecriture tensorielle permet de traiter des arguments de dualit\'e 
sous forme scripturale. Pour montrer que ce jeu d'\'ecriture est bien 
plus qu'un jeu, 
prenons de nouveau les \egts de Strassen-Winograd que nous 
r\'e\'ecrivons avec des tenseurs o\`u nous ne marquons pas la 
diff\'erence  entre formes \lins et vecteurs. Cela donne alors 
pour le produit matriciel l'\egt
$$ \som_{i,j,k\in\{1,2\}} {a_{ij}\otimes b_{jk}\otimes c_{ik}}\quad = 
\quad 
\som_{\ell=1}^7 {\; a_\ell\otimes b_\ell\otimes c_\ell}
$$
avec
{\footnotesize {$$\begin{array}{lll} 
a_1 := a_{11}  & b_1 := b_{11}   &   c_1:=c_{11}+c_{12}+c_{21}+c_{22}  
\\ 
a_2 := a_{12}  & b_2 := b_{21}  & c_2:=c_{11}   \\
a_3 := a_{11}-a_{21} & b_3 := b_{22}-b_{12}  & c_3:=c_{21}+c_{22}  \\ 
a_4 := a_{21}+a_{22} & b_4 := b_{12}-b_{11}  & c_4:=c_{12}+c_{22}   \\
a_5 := a_{21}+a_{22}-a_{11} &b_5 := b_{22}-b_{12}+b_{11}  & 
c_5:=c_{12}+c_{21}+c_{22}  \\ 
a_6 := a_{12}-a_{21}-a_{22}+a_{11} & b_6 := b_{22} &  c_6:=c_{12}    
\\
a_7 := a_{22}  & b_7 := b_{22}-b_{12}+b_{11}-b_{21}  & c_7:=-c_{21}   
\end{array}$$}}

\noindent Pour r\'ealiser une sym\'etrie par permutation circulaire 
dans 
la d\'efinition, \'echangeons les indices \,$i$\, et \,$k$\, dans les 
\,$c_{ik}$. Alors, vu l'invariance  par permutation circulaire nous 
pouvons remplacer partout \,$a$, \,$b$\,  et \,$c$\, par \,$b$, 
\,$c$\,  
et \,$a$. Et finalement nous permutons \`a nouveau les indices \,$i$\, 
et \,$k$\, dans les (nouveaux) \,$c_{ik}$\, pour revenir \`a la 
d\'efinition. Ceci nous donne d'autres \egtsz, qui peuvent tout 
aussi bien servir que les premi\`eres:
{\footnotesize {$$\begin{array}{lll} 
a_{1}:= a_{11}+a_{21}+a_{12}+a_{22} & b_{1}:= b_{11} & c_{1}:= c_{11} 
\\
a_{2}:= a_{11} & b_{2}:= b_{12} & c_{2}:= c_{12} \\
a_{3}:= a_{12}+a_{22} & b_{3}:= b_{11}-b_{21} & c_{3}:= c_{22}-c_{21} 
\\
a_{4}:= a_{21}+a_{22} & b_{4}:= b_{21}+b_{22} & c_{4}:= c_{21}-c_{11} 
\\
a_{5}:= a_{21}+a_{12}+a_{22} & b_{5}:= b_{21}+b_{22}-b_{11} & c_{5}:= 
c_{22}-c_{21}+c_{11} \\
a_{6}:= a_{21} & b_{6}:= b_{12}-b_{21}-b_{22}+b_{11} & c_{6}:= c_{22} 
\\
a_{7}:=-a_{12} & b_{7}:= b_{22} &  c_{7}:= c_{22}-c_{21}+c_{11}-c_{12}
\end{array}$$}}

\noindent Naturellement, dans le cas pr\'esent, on obtient seulement
sans fatigue un nouveau \sys d'\idas pour 
traiter le m\^eme produit matriciel. Mais si nous \'etions parti d'un 
produit de \marcs non carr\'ees, la permutation 
circulaire deviendrait un outil vraiment efficace, produisant des 
identit\'es correspondant \`a un cas de figure vraiment nouveau. Cette 
remarque importante remonte \`a 1972 (cf. \cite{HM,Pan1}).

\ss Notez que nous avons une situation famili\`ere analogue si nous 
consid\'erons le cas des \alis de \,$E$\, vers \,$F$. La dualit\'e 
nous 
dit que le passage de \,$\varphi$\, \`a \,$\tra{\varphi}$\, est un 
\isoz. En termes de matrices c'est une banale transposition. En termes 
d'\'ecriture tensorielle, c'est un jeu d'\'ecriture. Les tenseurs 
remplacent les matrices lorsqu'il y a plus que deux espaces en cause.
\begin{notation} 
\label{notaRGtens}   \emph{(Rang tensoriel de la \mul des matrices)} 
\\
Soient \,$m,\; n,\; p$\, trois entiers $>0$ et $\K$ un corps. On note
\,$\gen{m,n,p}_{\K}\,$\indexnota{a mnpK@$\gen{m,n,p}_{\K}$ (tenseur)} 
(ou \,$\gen{m,n,p}$)
\indexnota{a mnp@$\gen{m,n,p}$ (tenseur)} l'\abi 
$$ (A,B)\longmapsto AB\quad \mathrm{o\grave u}\quad 
A\in\K^{m\times n},\;  B\in\K^{n\times p} \; \mathrm{et}\;  
AB\in\K^{m\times p}
$$
On note donc le \rgte par \,$R\,\gen{m,n,p}$, ou si on doit pr\'eciser 
\,$R_\K\,\gen{m,n,p}$. 
\end{notation}
\begin{proposition} 
\label{propRGtens} \emph{(Rang tensoriel de la \mul des matrices)}
\begin{itemize}
\item [$(1)$] Si \,$m\leq m',\;n\leq n'$\, et \,$p\leq p'$\, alors 
\,$R\,\gen{m,n,p}\leq R\,\gen{m',n',p'}$
\item [$(2)$] $R\,\gen{mm',nn',pp'}\le R\,\gen{m,n,p}\cdot 
R\,\gen{m',n',p'}$ 
\item [$(3)$] $R\,\gen{m_1+m_2,n_1+n_2,p_1+p_2}\le\sum_{i,j,k\in 
\{1,2\}} R\,\gen{m_i,n_j,p_k}$ 
\item [$(4)$] $R\,\gen{n^\ell,n^\ell,n^\ell} \le 
(R\,\gen{n,n,n})^\ell$
\item [$(5)$] $R\,\gen{m,n,p}$\,  est invariant par permutation des 
entiers \,$m,\,n,\,p$.
\item [$(6)$] $R\,\gen{1,n,1}=n$\, et $\, R\,\gen{m,n,1}=mn$
\end{itemize}
\end{proposition}
\preuve
Le point $(1)$ est facile: on peut compl\'eter des matrices 
correspondant au format $(m,n,p)$ par des $0$ pour en faire des 
matrices 
au format $(m',n',p')$.\\  
Les points (2), (3) r\'esultent de la possibilit\'e de faire des 
produits de matrices par blocs. Et (4) r\'esulte de (2).\\
Le point (5) a \'et\'e expliqu\'e avant la proposition (voir page 
\pageref{symrgte}). On peut redire \`a peu pr\`es la m\^eme chose sous 
la forme suivante un peu plus abstraite, qui d\'ecrit 
peut-\^etre mieux l'essence du r\'esultat. Si \,$E_1$, \,$E_2$, 
\,$E_3$\, sont trois \Kevs de dimensions finies 
\,$p,\; n,\; m$, alors la \mul des matrices correspondante 
\,$\gen{m,n,p}:(A,B)\mapsto AB$\, est un \elt canonique 
\,$\theta$\, de l'espace
$$ \Bil(\Hom(E_2,E_3),\Hom(E_1,E_2);\Hom(E_1,E_3))=\Bil(E,F;G)
$$
Si nous notons \,$\mathrm{L}_3(E,F,G)$\, l'espace des formes 
tri\lins sur 
\,$E\times F\times G$, nous avons un \iso canonique
$$ \Bil(E,F;G)\simeq \mathrm{L}_3(E,F,G\sta)\qquad \psi\longmapsto 
\varphi=
((x,y,\gamma)\mapsto \gamma(\psi(x,y))
$$
Dans la situation pr\'esente, il y a aussi une dualit\'e canonique
entre \,$\Hom(E_1,E_3)=G$\, et \,$\Hom(E_3,E_1)$\, donn\'ee sous forme 
matricielle
par \,$(C,D)\mapsto \Tr(CD)$, ce qui fournit un \iso canonique entre
 \,$\Hom(E_3,E_1)$\, et \,$G\sta$. Une fois mis bout \`a bout tous ces 
\isos canoniques, on voit que l'\elt canonique
$$ \theta\in\Bil(\Hom(E_2,E_3),\Hom(E_1,E_2);\Hom(E_1,E_3))
$$
correspond \`a l'\elt canonique
$$ \theta'\in \mathrm{L}_3(\Hom(E_2,E_3),\Hom(E_1,E_2),\Hom(E_3,E_1))
$$
donn\'e sous forme matricielle par \,$(A,B,D)\mapsto \Tr(ABD)$.
Maintenant il est bien connu que \,$\Tr(ABD)=\Tr(BDA)=\Tr(DAB)$\, et 
\,$\Tr(ABD)=\Tr(\tra{(ABD)})=\Tr(\tra{D}\tra{B}\tra{A})$. 
Ceci \'etablit les sym\'etries demand\'ees.\\
Voyons maintenant le point $(6)$. Nous reprenons les notations 
pr\'ec\'edentes avec \,$E_1=\K$, donc on identifie \,$\Hom(E_1,E_i)$\, 
\`a \,$E_i$\, ($i=2,3$). Regardons l'espace  \,$\Bil(E,F;G)$\, sous la 
forme (canoniquement \'equivalente) \,$\Hom(E,\Hom(F,G))$. Si on 
\'ecrit 
le produit d'une matrice par un vecteur colonne \,$(A,X)\mapsto AX$\, 
sous forme
$$ \som_{i,j=1}^{i=m,j=n} {\alpha_{ij}\otimes \beta_{j1}\otimes 
c_{i1}}\quad = \quad 
\som_{\ell=1}^r {\; \alpha_\ell\otimes \beta_\ell\otimes c_\ell}
$$
on voit que l'\ali correspondante de \,$E$\, vers \,$\Hom(F,G)$\, est 
nulle sur \,$\Inter_{l=1}^r{\Ker(\alpha_\ell)}$. 
Mais dans le cas pr\'esent, modulo les identifications 
pr\'ec\'edentes, 
cette \ali n'est autre que l'application identique de \,$E$. Son noyau 
est donc r\'eduit \`a \,$\{0\}$\, et \,$r$\, est au moins \'egal \`a 
la 
dimension de \,$E$\, \cad \`a \,$mn$.
\qed 

\subsection{Exposant de la \mul des \macas} 
\label{paraExpMu}

\begin{definition} 
\label{defExpMu} 
On dit que \emph{$\alpha$ est un exposant acceptable pour la \mul des 
\macasz}\index{exposant!acceptable pour la \mul des \macas} 
si celle-ci peut \^etre r\'ealis\'ee en \,$\SD(n^\alpha,\alb\log n)$.
La borne inf\'erieure des exposants acceptables est appel\'ee 
\emph{l'exposant de la \mul des \macasz}\index{exposant!de la \mul des 
\macas} et elle est not\'ee \,$\omega$.
\end{definition}
A priori on devrait mettre en indice le corps $\K$ pour les exposants 
$\alpha$ et $\omega$. Les r\'esultats concernant ces exposants dont 
nous 
rendons compte sont cependant ind\'ependants du corps $\K$ 
consid\'er\'e. 
\begin{theorem} 
\label{thRGtens} \emph{(Rang tensoriel et exposant de la \mul des 
\macasz)} 
\begin{itemize}
\item [$(1)$] S'il existe \,$n$\, et \,$r$\, tels que 
\,$R\,\gen{n,n,n}=r$\, alors l'exposant \,$\alpha= {\log r\over\log 
n}$\, est acceptable pour la \mul des \macasz.
\item [$(2)$] S'il existe \,$m$, \,$n$, \,$p$ et \,$r$ tels que 
\,$R\,\gen{m,n,p}=r$\, alors l'exposant  
\,$\alpha= 3\,{\log r\over \log mnp}$\, est acceptable pour la \mul 
des 
\macasz.
\end{itemize}
\end{theorem}
\preuve
Comme nous l'avons d\'ej\`a remarqu\'e le point (1) r\'esulte du 
m\^eme 
calcul de \com que celui fait dans la section \ref{sec strass}.\\
Le point (2) r\'esulte du point (1) puisque d'apr\`es les points (2) 
et 
(5) de la proposition \ref{propRGtens} on a
(avec \,$N=mnp$) 

\smallskip \centerline{$R\,\gen{N,N,N}\leq 
R\,\gen{m,n,p}\,R\,\gen{n,p,m}\,R\,\gen{p,m,n}=(R\,\gen{m,n,p})^3 
$}

\smallskip \qed

\ss En fait la conclusion dans le \tho pr\'ec\'edent est non 
seulement qu'il existe une famille de \caris dans la classe 
\,$\SD(n^\alpha,\log n)$\, qui r\'ealise la \mul des \macasz, mais 
qu'on 
sait construire explicitement une telle \famu de \carisz. En 
outre le temps de construction du \cari num\'ero \,$n$\, est 
proportionnel \`a sa taille, selon les lignes de la preuve du 
\tho \ref{thMuStra0}.

Le point (1) du \tho ci-dessus peut \^etre pr\'ecis\'e comme 
suit
(par le m\^eme calcul qu'\`a la section \ref{sec strass}).

\begin{proposition} 
\label{propthRGtens} \emph{(Pr\'ecision pour le \tho 
\ref{thRGtens}~(1))}
Supposons que l'\abi \,$\gen{n,n,n}$\, puisse \^etre calcul\'ee par
un \cari de \prof  \,$\ell$\, contenant  \,$r$\, \muls essentielles et  
\,$s$\, autres \oparis (addition, soustraction, \mul par une 
constante), 
avec
\,$r>n^2$. Alors 
l'\abi \,$\gen{n^\nu,n^\nu,n^\nu}$\, peut \^etre calcul\'ee par
un \cari de \prof  \,$\nu\,\ell$\, contenant  \,$r^\nu$\, \muls 
essentielles et  \,$s\,{r^\nu-n^{2\nu}\over r-n^2} $\, autres 
\oparisz.
\end{proposition}

\subsection[Complexit\'es \bil et multiplicative]
{Complexit\'e \bil versus \com multiplicative} 
\label{paraComBiMu}
  
Soient $\K$ un corps, \,$H$, \,$G$\,  deux \Kevs de 
dimension finie. Une \emph{application quadratique}
\index{application!quadratique}\index{quadratique!application} de 
\,$H$\, vers \,$G$\, est 
par d\'efinition une application de la forme 
\,$\Psi:x\mapsto \psi(x,x)$\, o\`u 
\,$\psi\in \Bil(H,H;G)$. 
Si \,$(h_i)_{i\in I}$\, et \,$(g_\ell)_{\ell\in L}$\, sont des bases 
de 
\,$H$\, et \,$G$\, il revient au m\^eme de dire que chaque 
coordonn\'ee 
de \,$\Psi(x)$\, est une forme quadratique en \,$x$, \cad un \pol 
homog\`ene du second degr\'e en les coordonn\'ees de \,$x$. Si les 
coordonn\'ees de \,$x$\, sont prises comme variables, on peut alors 
consid\'erer les \prevs \ariths sans division qui permettent de 
calculer 
les coordonn\'ees de \,$\Psi(x)$. La \emph{\com \muv de 
\,$\Psi\,$}\index{complexit\'e!multiplicative} est alors d\'efinie 
comme 
la plus petite longueur \muv d'un tel \prevz. Nous la noterons 
\,$M(\Psi)$. Comme les changements de base ne co\^utent rien en 
longueur 
\muv cette d\'efinition ne d\'epend pas du choix des bases.
Le lemme suivant est une paraphrase de la proposition 
 \ref{prop Eli Div4} dans le cas d'une application quadratique.
\begin{lemma} 
\label{lemCOMU} 
Avec les notations pr\'ec\'edentes la \com \muv d'une application 
quadratique $\Psi$ est aussi \'egale au plus petit entier \,$r$\, tel 
que \,$\Psi$\, puisse s'\'ecrire comme compos\'ee de trois 
applications 
selon le format suivant
$$ H\surver{(\eta,\zeta)} \K^r\times \K^r
\surver{\mu_r}\K^r\surver{g}G
$$
o\`u $\;\eta:H\rightarrow \K^r$,  $\;\zeta:H\rightarrow \K^r$\,  et
 $\;g:\K^r\rightarrow G$\,  sont des \alis et \,$\mu_r:\K^r\times 
\K^r\rightarrow \K^r$\, est le produit coordonn\'ee
par coordonn\'ee. Un \prev \arith correspondant \`a cette 
\deco s'appelle un \emph{calcul quadratique de 
\,$\Psi$}\index{calcul!quadratique}\index{quadratique!calcul}. 
\end{lemma}

\begin{remark} \label{remlemCOMU}
\emph{ Si on consid\'erait des \caris avec division on ne pourrait pas 
diminuer pour autant la longueur \muv pour \'evaluer une application 
quadratique, au moins dans le cas d'un corps infini, d'apr\`es le 
\tho \ref{thEliDiv} et la proposition \ref{prop Eli Div3}. 
} 
\end{remark}

\begin{proposition} 
\label{propCOMU} Soient $\K$ un corps, \,$E$, \,$F$, \,$G$  trois 
\Kevs de dimension finie. Soit \,$\psi\in\Bil(E,F;G)$\, et 
\,$H=E\times F$.
Alors \,$\psi$\, est une application quadratique de \,$H$\, vers 
\,$G$. 
Sa \com \bil \,$R(\psi)$\, et sa \com \muv \,$M(\psi)$\, sont 
reli\'ees 
par
$$ M(\psi)\leq R(\psi)\leq 2M(\psi)
$$
\end{proposition}
\preuve
La premi\`ere \ine est \'evidente. Pour la seconde 
consid\'erons
un programme quadratique comme dans le lemme \ref{lemCOMU} qui calcule
\,$\psi(u,v)$\, avec \,$m=M(\psi)$\, \muls essenteielles. On a donc
$$ \psi(u,v)=\som_{\ell=1}^m{\; \alpha_\ell(u,v)\cdot\beta_\ell(u,v) 
\cdot g_{\ell}}
$$
o\`u les \,$\alpha_\ell$\, et \,$\beta_\ell$\,  sont dans \,$H\sta$. 
Remarquons qu'on a
$$\begin{array}{rcl} 
\alpha_\ell(u,v)\cdot\beta_\ell(u,v)& =  &  
\alpha_\ell(u,0)\cdot\beta_\ell(u,0)+
\alpha_\ell(0,v)\cdot\beta_\ell(0,v)+    \\ 
&   & \alpha_\ell(u,0)\cdot\beta_\ell(0,v)+
\alpha_\ell(0,v)\cdot\beta_\ell(u,0)
\end{array}$$
Puisque \,$\psi(u,v)$\, est \bil on a \,$\psi(u,0)=0$\, et on peut 
supprimer les termes \,$\alpha_\ell(u,0)\cdot\beta_\ell(u,0)\cdot 
g_\ell$\, dont la somme est nulle. M\^eme chose avec \,$\psi(0,v)=0$\, 
et finalement on obtient 
$$ \psi(u,v)=\som_{\ell=1}^m{\; 
(\alpha_\ell(u,0)\cdot\beta_\ell(0,v)+ 
\beta_\ell(u,0)\cdot\alpha_\ell(0,v))\cdot g_{\ell}}
$$
et ceci montre que \,$R(\psi)\leq 2r$.
\qed

\ss On en d\'eduit le r\'esultat suivant qui relie le \rgte et 
l'exposant de la \mul des \macasz.
\begin{theorem} 
\label{thRgExp} 
L'exposant \,$\omega$\, de la \mul des \macas est \'egal \`a la borne 
inf\'erieure des exposants \,$\alpha$\, qui v\'erifient, pour au moins 
un entier \,$n$, l'\ine 
\,$R\,\gen{n,n,n}\leq n^\alpha$. On a aussi
$$ \omega =\lim_{n\rightarrow \infty} {\log R\,\gen{n,n,n} \over \log 
n}
=\lim_{n\rightarrow \infty} {\log M(\gen{n,n,n}) \over \log n}
$$
o\`u chaque suite converge vers sa borne inf\'erieure.
\end{theorem}
\preuve
Il est clair d'apr\`es la proposition \ref{propCOMU} que les deux 
suites 
consid\'er\'ees ont la m\^eme borne inf\'erieure \,$\beta$. \\
On a le r\'esultat direct plus pr\'ecis dans le \tho 
\ref{thRGtens}: tout exposant \,$\alpha$\, strictement sup\'erieur 
\`a 
la borne inf\'erieure des \,${\log R\,\gen{n,n,n}\over \log n}$\, est 
acceptable  pour la \mul des \macasz. \\
Pour la r\'eciproque on consid\`ere un \,$\alpha>\beta$. 
Pour un \,$n_0$\, assez grand on a un \prev sans
division de longueur $< n_0^\alpha/2$\, qui calcule 
l'application quadratique \,$\gen{n_0,n_0,n_0}$.  
A fortiori sa longueur \muv est $< n_0^\alpha/2$\, et on a 
\,$R\,\gen{n_0,n_0,n_0}\leq 2M(\gen{n_0,n_0,n_0})< n_0^\alpha$. 
\qed

\ss Malgr\'e la relation tr\`es \'etroite entre \,$M(\psi)$\,  et 
\,$R(\psi)$,
c'est seulement la consid\'eration du \rgte qui permet de d\'emontrer 
les r\'esultats de base concernant l'exposant de la 
\mul des matrices. Cela tient \`a ce que la proposition 
\ref{propRGtens} ne serait pas vraie en rempla\c{c}ant 
le \rgte par la longueur multiplicative. 
Le fait d'interdire la commutation dans les tenseurs
est ce qui permet de traiter correctement le produit 
des matrices par blocs.

\subsection{Extension du corps de base}
\label{paraChangeCorps}

Soient $\K$ un corps, \,$E$, \,$F$, \,$G$\,  trois \Kevs de dimension 
finie.  Soit $\psi\in\Bil(E,F;G)$. Si 
$\L$ est une extension de $\K$ on peut \emph{\'etendre $\psi$ \`a 
$\L$} 
de mani\`ere naturelle. Nous nous en tiendrons ici \`a un point de vue 
pragmatique et purement calculatoire. Si
\,$(e_i)_{i\in I}$, \,$(f_j)_{j\in J}$, \,$(g_\ell)_{\ell\in L}$\, 
sont 
des bases de 
 \,$E$, \,$F$, \,$G$  et si 
$$ \psi\quad =\quad \som_{i,j,\ell} {\; \gamma_{ij\ell}\; 
e_i\sta\otimes 
f_j\sta\otimes g_\ell}
$$
nous consid\'erons trois \,$\L\,$--\,espaces vectoriels   \,$E_\L$, 
\,$F_\L$, 
\,$G_\L$  ayant les m\^emes bases et l'extension $\psi_\L$ de $\psi$ 
est 
d\'efinie par la m\^eme \egtz. 
Comme tout \cabi dans \,$\K$\, est aussi un \cabi dans \,$\L$\, on a 
\ncrt l'\ine \,$ R_\L(\psi_L)\leq R_\K(\psi)$, mais 
il se peut que l'u\-ti\-li\-sa\-tion de constantes dans \,$\L$\, puisse 
faciliter le calcul de \,$\psi$\, et l'\ine peut \^etre 
stricte.

Nous allons cependant voir dans ce paragraphe que l'exposant de la 
\mul des matrices ne peut pas changer lorsqu'on passe d'un corps 
\,$\K$\, \`a une extension \,$\L$. 
\begin{lemma} 
\label{lemChangCor} 
Avec les notations pr\'ec\'edentes \\
1) Si $\L$ est une extension finie de $\K$ de degr\'e \,$n$\, il 
existe 
un entier \,$m\leq n^3$\, tel que pour toute \abi \,$\psi$\, on a 
\,$ R_\K(\psi)\leq m\,R_\L(\psi_L)$. 
En particulier l'exposant de la \mul des matrices ne change pas 
lorsqu'on passe de \,$\K$\, \`a \,$\L$.\\
2) Si $\L= \K(t)$ (corps des fractions rationnelles 
en \,$t$) et si $\K$ 
est infini, on a l'\egt \,$ R_\L(\psi_L)= R_\K(\psi)$. 
\end{lemma}
\preuve
Dans le cas 2) la famille finie des constantes \,$c_s(t)$\, dans 
\,$\L$\, utilis\'ees par le \cari peut \^etre remplac\'ee par des 
constantes  \,$c_s(a)$\,  o\`u \,$a\in\K$\, est choisi de mani\`ere \`a 
n'annuler aucun des d\'enominateurs.\\
Dans le cas 1) consid\'erons une base \,$b=(1,b_2,\ldots,b_n)$\, de 
\,$\L$\, 
lorsqu'on le voit comme \Kevz. La \mul dans \,$\L$\, 
repr\'esente une \abi sur \,$\K$\, lorsqu'elle est traduite dans les 
coordon\'ees sur la base \,$b$. 
Cette \abi  \,$\L\times \L\rightarrow \L$\, 
peut \^etre r\'ealis\'ee par  \,$m\le n^3$\,  \muls essentielles dans 
$\K$. 
En fait la constante \,$m$\, peut \^etre prise \'egale au \rgte de 
cette 
\abiz, qui  est en g\'en\'eral not\'e  \,$R_\K(\L)$.\\
Tout \cabi dans \,$\L$\, peut  alors \^etre mim\'e par un \cabi dans 
\,$\K$\, de 
la mani\`ere suivante. Chaque variable \,$x_i$\, sur \,$\L$\,  est 
remplac\'ee par \,$n$\, variables \,$x_{i,n}$\, sur \,$\K$\, qui 
repr\'esentent les coordonn\'ees de \,$x_i$\, sur la base \,$b$.
Seules les \muls essentielles du calcul dans \,$\L$\, produisent des 
\muls 
essentielles dans \,$\K$. Dans cette simulation, le nombre de \muls 
essentielles est multipli\'e par la constante \,$m$.\\
Si maintenant \,$ \alpha>\omega_\L $\, il existe un entier \,$n$\, tel 
que
\,$ R_\L\,\gen{n,n,n}< n^\alpha$, donc pour une puissance convenable 
\,$N=n^\ell$\, 
on a $\; R_\K\,\gen{N,N,N}\leq  R_\K(\L)\cdot R_\L\,\gen{N,N,N}< 
N^\alpha$,
donc  \,$ \alpha\geq \omega_\K$.
\qed 
 
\ss On en d\'eduit le r\'esultat suivant d\^u \`a Sch\"onhage 
(\cite{Scho2}).

\begin{proposition} 
\label{propChangCor} 
L'exposant de la \mul des \macas sur un corps $\K$ ne d\'e\-pend que 
de 
la \cara de $\K$.
\end{proposition}
\preuve
Il suffit de prouver que l'exposant ne change pas lorsqu'on passe
d'un corps premier $\K$  ($\QQ$  ou l'un des $ \FF_p$) \`a l'une de 
ses 
extensions $\L$. Supposons qu'on ait \,$R_\L(\gen{n,n,n})\leq r\leq 
n^\alpha$, \cad qu'on ait sur le corps $\L$ une \egt
$$ 
\som_{i,j,k\in\{1,\ldots,n\}} {\alpha_{ij}\otimes\beta_{jk}\otimes 
c_{ik}}\quad = \quad 
\som_{\ell=1}^r {\; \alpha_\ell\otimes \beta_\ell\otimes c_\ell}
$$
On peut consid\'erer les coordonn\'ees des 
\,$ \alpha_\ell,\;  \beta_\ell,\;  c_\ell$\, sur les bases 
\,$\alpha_{ij},\; \beta_{ij}$, \,$c_{ij}$\, comme \,$3rn^2$\, 
\idtrs \,$z_s$.
L'\egt des deux tenseurs ci-dessus signifie que ces 
\idtrs v\'erifient un \sys de \,$n^6$ \'equations 
\polles de degr\'e $3$ dont tous les \coes sont \'egaux \`a 1 ou 0.
Maintenant, le lecteur ou la lectrice conna\^{\i}t peut-\^etre le beau 
r\'esultat suivant{\footnote{~Ce r\'esultat admet de nombreuses 
preuves, 
dont certaines tout \`a fait explicites. Essentiellement c'est un 
r\'esultat de la th\'eorie de l'\'elimination. On peut le faire 
d\'ecouler du 
Null\-stellen\-satz de Hilbert, du lemme de normalisation de Noether 
ou 
encore de la th\'eorie des bases de Gr\"obner. Il se trouve dans les 
bons livres d'\agrz.}}: lorsqu'un \syse 
\polles sur un corps $\K$  admet une solution dans une extension 
$\L$, alors il admet une solution dans une extension finie de $\K$.  
On 
est donc ramen\'e au premier cas du lemme pr\'ec\'edent.  
\qed

\begin{remark}\label{remExpMu}
 \emph{ L'exposant \,$\omega_\K$\, peut donc \^etre \'etudi\'e en 
prenant pour $\K$ la cl\^oture \agq de $\QQ$ ou de $\FF_p$. 
Lorsqu'on a affaire \`a un corps \agqt clos $\K$ le \rgte 
\,$R_\K(\gen{n,n,n})$\, est calculable en principe (sinon en pratique) 
car  savoir si \,$R_\K(\gen{n,n,n})\le m$\, revient \`a d\'eterminer 
si un \syse \agqs admet ou non une solution 
(comme dans la preuve de la proposition \ref{propChangCor}). 
Et on sait, en principe, r\'epondre \`a ce genre 
de questions par un \algo d'\'elimination. 
On ne sait cependant pas grand chose concernant 
\,$\omega_\K$. Cet exposant mythique est un nombre compris entre $2$  
et $2,38$. Mais on ne sait apparemment toujours rien sur la vitesse avec 
laquelle la suite 
\,$\log R_\K(\gen{n,n,n})/\log n$\, (cf. \tho \ref{thRgExp}) 
converge vers \,$\omega_\K$. Il se pourrait que la vitesse de 
convergence soit si lente que le nombre  \,$\omega_\K$ serait 
d\'efinitivement impossible \`a calculer. 
} 
\end{remark}

\section[Calculs bilin\'eaires approximatifs]{Acc\'el\'eration par la 
\met des \cabas}  
\label{sec bini}
 
\subsection{M\'ethode de Bini} 
\label{paraCaBiAp}

La \met des \cabas est inspir\'ee des \mets 
num\'eriques approch\'ees, elle a \'et\'e invent\'ee par Bini et elle 
a 
quelque parent\'e avec l'\elidiv de Strassen. Elle a d\'ebloqu\'e la 
situation pour l'exposant \,$\omega$\, de la \mul des matrices.

Un exemple est le produit de deux matrices \,$A,B\in\A^{2\times 2}$\, 
avec la premi\`ere qui a son \coe \,$a_{2,2}$\, nul. On sch\'ematise 
ce 
produit sous la forme suivante (figure \ref{fbini}).
\begin{figure}[ht]   
\begin{center}
\includegraphics*[width=3cm]{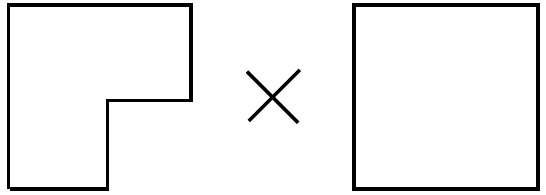}
\end{center}
\caption[Produit matriciel \`a trou de Bini]
{\label{fbini} Produit matriciel \`a trou de Bini }  
\end{figure}  
Ce \emph{produit matriciel \`a trous}\index{produit matriciel \`a 
trous}  
correspond \`a un tenseur de rang $6$ qui s'\'ecrit avec la notation 
des 
\pols non commutatifs
$$ \psi=a_{11}b_{11}c_{11} + a_{12}b_{21}c_{11}+ a_{21}b_{11}c_{21}+ 
a_{11}b_{12}c_{12}+ a_{12}b_{22}c_{12}+ a_{21}b_{12}c_{22}
$$
On introduit (pour les m\^emes variables) un tenseur de rang $5$ 
perturb\'e par des \,$\vep\, x_{ij}$\, (avec \,$x=a,\; b$\, ou \,$c$).
$$\begin{array}{rcl} 
\varphi(\vep)& = &  (a_{12}+\vep a_{11})\,(b_{12}+\vep b_{22})\,c_{21} 
+ 
(a_{21}+\vep a_{11})\, b_{11}\,(c_{11}+\vep c_{12})  \\ 
& &  - a_{12}\,b_{12}\,(c_{11}+c_{21}+\vep c_{22})  
- a_{12}\, (b_{11}+b_{12}+\vep b_{21})\, c_{11} \\
& & + (a_{12}+a_{21})\,(b_{12}+\vep b_{21})\,(c_{11}+\vep c_{22})
\end{array}$$
Lorsqu'on d\'eveloppe on obtient
$$ \varphi(\vep) = \vep\, \psi + \vep^2\, \theta(\vep)
$$
Num\'eriquement on a donc lorsque \,$\vep$\, est suffisamment petit 
\,$\psi\simeq \varphi(\vep)/\vep$. On dit que \,$\varphi$\, constitue 
une approximation d'ordre $1$  de \,$\psi$. On peut transformer ceci 
en 
un calcul purement formel dans l'anneau des d\'eve\-lop\-pe\-ments 
limit\'es \`a l'ordre $1$  en \,$\vep$\, comme lorsqu'on \'elimine les 
divisions \`a la Strassen. Naturellement, il n'y a pas de miracle, 
cela 
ne donne pas une \'ecriture de \,$\psi$\, comme somme de $5$  tenseurs 
\elrsz. Mais il y a n\'eanmoins quelque chose \`a gagner en prenant un 
peu de recul et en analysant en d\'etail ce qui se passe.
Tout d'abord en appliquant le sch\'ema suivant
\begin{figure}[ht]   
\begin{center}
\includegraphics*[width=10cm]{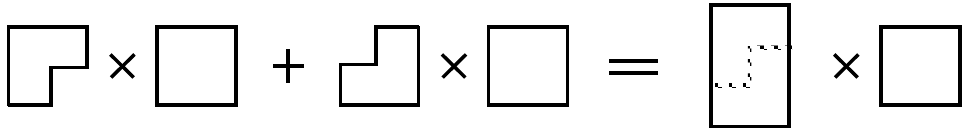}
\end{center}
\caption[Produit matriciel plein]
{\label{fbini2} Produit matriciel plein }  
\end{figure}  
(o\`u on fait un produit par blocs non rectangulaires~!) on constate 
que 
le produit matriciel \,$\gen{3,2,2}$\, peut \^etre r\'ealis\'e de 
mani\`ere approximative (\`a l'ordre 1) par une somme de 10 tenseurs 
\elrs au lieu de 12. Gr\^ace au produit de matrices par blocs 
rectangulaires on pourra alors r\'ealiser \,$\gen{12,12,12}$\,  de 
mani\`ere approximative (\`a un ordre convenable, nous allons 
d\'efinir 
cela un peu plus loin) comme somme de \,$10^3$\, tenseurs
\elrs au lieu des \,$12^3$\, \ncrs dans la \met \uslez. 
Enfin il reste \`a r\'ealiser que lorsqu'on passe (gr\^ace au produit 
par blocs) \`a  \,$\gen{12^n,12^n,12^n}$\, l'ordre d'approximation ne 
cro\^{\i}t pas trop vite et que le co\^ut du d\'ecryptage d'un \caba 
en  un \cabi exact devient n\'egligeable devant \,$10^{\vep n}$\, pour 
n'importe quel \,$\vep>0$. Tout ceci ram\`ene l'exposant \,$\omega$\, 
\`a 
$$ 3\,\log(10)/\log(12)<2,78<3\,\log(7)/\log(8)=2,807
$$

Nous devons maintenant donner des d\'efinitions et \'enonc\'es plus 
pr\'ecis pour v\'erifier que ce plan de travail fonctionne bien.

\begin{definition} 
\label{defCabas} 
Soient $\K$ un corps, \,$E$, \,$F$, \,$G$\,  trois \Kevs de dimensions 
finies. Soit $\psi\in\Bil(E,F;G)$. Soit 
$\L=\K[\vep]$\, l'anneau des \pols en la variable \,$\vep$\, sur $\K$. 
Un \elt
\,$\varphi(\vep)$\, de \,$\Bil(E_\L,F_\L;\alb G_\L)$\, est appel\'e une 
\emph{approximation d'ordre \,$q$\, de \,$\psi\,$}\index{approximation 
d'ordre \,$q$\, d'une \abiz} si on a
$$ \varphi(\vep)\equiv\vep^{q}\,\psi \mathrm{\qquad  modulo\quad  } 
\vep^{q+1}
$$
Un \cabi de \,$\varphi(\vep)$\, est appel\'e un \emph{\caba de 
\,$\psi$\, \`a l'ordre $q$}\index{calcul!bilin\'eaire approximatif}. 
On appelle \emph{\rgtm de  \,$\psi$\, \`a l'ordre \,$q$} le plus petit 
rang possible pour un \caba de  \,$\psi$\, \`a l'ordre \,$q$. On le 
note 
\,$\Ra(\psi,q)$. Enfin, le \emph{\rgtm de  $\psi$}\index{rang 
tensoriel!marginal} est le plus petit des  \,$\Ra(\psi,q)$\, et il est 
not\'e \,$\Ra(\psi)$\indexnota{Ra@$\Ra(\psi)$}. Nous dirons aussi plus 
simplement le \emph{rang  marginal de  \,$\psi$}\index{rang marginal}.
\end{definition}
\begin{remark} \label{remCabas1}
\emph{Nous utilisons ici des \cabis sur un anneau, comme il \'etait 
indiqu\'e dans la remarque \ref{remRGtens3}. De m\^eme l'extension de 
\,$\psi$\, \`a l'anneau $\K[\vep]$\,  se fait comme dans le cas d'une 
extension du corps de base (cf. page \pageref{paraChangeCorps}).
} 
\end{remark}

\begin{remark} \label{remCabas2}
\emph{Il est clair que lorsque \,$q$\, augmente, le \rgtm \`a l'ordre 
\,$q$\, d'une \abi ne peut que diminuer, autrement dit 
$$ R(\psi)\; =\; \Ra(\psi,0)\; \geq \; \Ra(\psi,1)\; \geq \;\cdots
\; \geq \; \Ra(\psi,q)\; \geq \; \Ra(\psi,q+1)\; \cdots
$$
Le \rgtm \`a l'ordre \,$q$\, d'une \abi est a priori nettement plus 
difficile \`a calculer que son \rgtez. Le rang marginal est encore plus 
difficile \`a \'etablir. En fait on est satisfait quand on a \'etabli 
une bonne majoration du rang marginal en explicitant un \cabaz. 
} 
\end{remark}

Quel est le co\^ut d'un d\'ecryptage d'un \caba \`a l'ordre \,$q$~? 
Nous avons d\'ej\`a fait un calcul analogue dans la preuve de la 
proposition 
\ref{prop Eli Div4} qui concernait la possibilit\'e d'une mise en 
forme 
simplifi\'ee des \caris sans division. 

On commence par consid\'erer que le \cabi de \,$\varphi(\vep)$\, se 
passe non pas sur l'anneau \,$\L=\K[\vep]$\, mais sur l'anneau des 
\dlis \`a l'ordre \,$q$. Ensuite on simule toute 
variable \,$Z$, qui repr\'esente un \elt de  
\,$\K[\vep]$\, modulo \,$\vep^q$, par \,$q+1$\, variables 
\,$Z^{[k]}$\,  
dans $\K$ ($0\leq k\leq  q$) qui repr\'esentent les \coes de  \,$Z$\, 
en 
dessous du degr\'e \,$q$.
Quand on doit calculer le \coe de \,$\vep^{q}$\, dans un tenseur 
\,$X(\vep)\,Y(\vep)\,Z(\vep)$\, on voit qu'on doit faire la somme des 
\,$X^{[i]}\,Y^{[j]}\,Z^{[k]}$\, pour tous les triplets $(i,j,k)$ dont 
la 
somme vaut \,$q$. Il y a au plus \,$(q+1)(q+2)/2$\, triplets de ce 
type. 
En termes de \rgtez, cela signifie donc que le \rgte de \,$\psi$\, est 
major\'e par  \,$(q+1)(q+2)/2$\, fois son rang marginal \`a l'ordre 
\,$q$. Nous avons donc \'etabli le lemme suivant. 
\begin{lemma} 
\label{lemCabas} 
Soit \,$\psi$\, une \abi d\'efinie sur un corps $\K$. Si un \caba \`a 
l'ordre \,$q$\,  de \,$\psi$\, a une \com \bil \,$\ell$, on en 
d\'eduit 
par d\'ecryptage un \cabi de \,$\psi$\, de \com \bil $\leq \ell\cdot 
(q+1)(q+2)/2$. En bref 
$$ 
R(\psi)\, \leq \, {(q+1)(q+2)\over2}\;\; \Ra(\psi,q)\quad\quad   
\left(\, \leq \, q^2\; \Ra(\psi,q)\quad  \mathrm{si\quad  }q\geq 
4\right)
$$
\end{lemma}
Maintenant nous devons examiner comment se comporte le rang marginal 
du 
produit matriciel lorsqu'on utilise des produits par blocs.
\begin{proposition} 
\label{propRGmarg} \emph{(Rang tensoriel marginal de la \mul des 
matrices)}
\begin{itemize}
\item [$(1)$] Le rang marginal \,$\Ra(\gen{m,n,p},q)$\,  est une 
fonction croissante de chacun des entiers \,$m$, \,$n$, \,$p$\, et 
d\'ecroissante de l'entier \,$q$.
\item [$(2)$] $ \Ra(\gen{mm',nn',pp'},q+q') \le 
\Ra(\gen{m,n,p},q)\cdot 
\Ra(\gen{m',n',p'},q')$. En particulier
avec $N=mnp$ on a\\
\centerline{$ \Ra(\gen{N,N,N},3q) \le 
\left(\Ra(\gen{m,n,p},q)\right)^3$  
}
\item [$(3)$] $ \Ra(\gen{m_1+m_2,n_1+n_2,p_1+p_2},q) \le 
\sum_{i,j,k\in 
\{1,2\}} \Ra(\gen{m_i,n_j,p_k},q)$ 
\item [$(4)$] $ \Ra\left(\gen{n^\ell,n^\ell,n^\ell},\ell q\right) \le 
\Ra(\gen{n,n,n},q)^\ell$
\item [$(5)$] $\Ra(\gen{m,n,p},q)$\,  est invariant par permutation 
des 
entiers $m,n,p$.
\end{itemize}
\end{proposition}
\preuve
Tout se passe comme avec le \rgte usuel dans la preuve de la 
proposition 
\ref{propRGtens}. Le seul point qui demande un peu d'attention est le 
point $(2)$. La meilleure mani\`ere de le comprendre est (encore une 
fois) de prendre du recul. Il faut prendre du recul sur ce que 
repr\'esente le tenseur \,$\gen{m_1m_2,n_1n_2,p_1p_2}$\, par rapport 
aux 
tenseurs \,$\gen{m_1,n_1,p_1}$\, et \,$\gen{m_2,n_2,p_2}$. Lorsque 
nous 
voyons une matrice \,$A$\, de type
\,$m_1m_2\times n_1n_2$\, comme une matrice de type 
\,$m_1\times n_1$\,  ayant pour entr\'ees des matrices \,$A_{ij}$\, de 
type $m_2\times n_2$ nous rep\'erons une entr\'ee de la grosse matrice 
par deux paires d'indices $((i_1,j_1),(i_2,j_2))$ correspondant au 
couple d'indices 
$(i_1(m_2-1)+i_2,j_1(n_2-1)+j_2)$ comme dans l'exemple 
d\'ecrit par la figure \vref{fblocs} avec 
$(m_1,n_1)=(5,6)$, $(m_2,n_2)=(3,4)$, $(i_1,i_2)=(3,1)$ et 
$(j_1,j_2)=(4,2)$.

\begin{figure}[h]   
\begin{center}\label{fblocs}
{\setlength{\unitlength}{1pt}%
\begin{picture}(340,170)(-58,0) 
\thicklines
\multiput(0,0)(40,0){7}{\line(0,1){150}}
\multiput(-0.4,0)(0,30){6}{\line(1,0){240.6}}
\multiput(10,0)(40,0){6}{\qbezier[50](0,1)(0,75)(0,149)}
\multiput(20,0)(40,0){6}{\qbezier[50](0,1)(0,75)(0,149)}
\multiput(30,0)(40,0){6}{\qbezier[50](0,1)(0,75)(0,149)}
\multiput(0,10)(0,30){5}{\qbezier[80](1,0)(120,0)(239,0)}
\multiput(0,20)(0,30){5}{\qbezier[80](1,0)(120,0)(239,0)}
\put(-30,80){$(3,1)$}
\put(121,158){$(4,2)$}
\put(135.5,85){\circle*{4}}
\end{picture}}
\end{center}
\caption[Num\'erotation par blocs]
 {Num\'erotation par blocs }  
\end{figure}
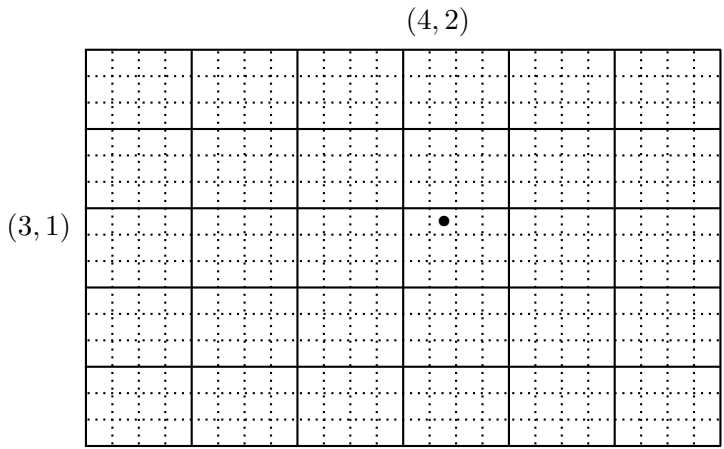

Cependant la mise en ligne de la paire $(i_1,i_2)$ sous forme 
\,$i_1(m_2-1)+i_2$, si elle est indispensable au dessin et \`a une 
premi\`ere compr\'ehension des choses, est plut\^ot un obstacle pour 
ce 
qui concerne la compr\'ehension du calcul \gui{emboit\'e} que 
repr\'esente 
un produit par blocs.
Prenons en effet les indices dans la grande matrice sous forme des 
couples $(i,j)=((i_1,i_2),(j_1,j_2))$  comme dans la figure 
\ref{fblocs}
(et non pas $(i_1(m_2-1)+i_2,j_1(n_2-1)+j_2)$ ni non plus 
$((i_1,j_1),(i_2,j_2))$). 
Nous obtenons, en notation de \pols non commutatifs:
$$ \gen{m_1m_2,n_1n_2,p_1p_2}\; =\; \sum_{i,j,k}
{\; a_{i,j}\;  b_{j,k} \; c_{i,k}}
$$
(o\`u la somme est prise sur \,$i\in I_1\times I_2$, 
\,$j\in J_1\times J_2$, \,$k\in K_1\times K_2$\, pour des ensembles 
d'indices de cardinalit\'es convenables). Alors nous avons 
l'\egtz:
$$ \gen{m_1m_2,n_1n_2,p_1p_2}\; =\; \gen{m_1,n_1,p_1}\; 
\gen{m_2,n_2,p_2}
$$
o\`u
$$\begin{array}{rcl} 
\gen{m_1,n_1,p_1}&  =  & \sum_{i_1,j_1,k_1} 
{\; a_{i_1,j_1}\;  b_{j_1,k_1}\;  c_{i_1,k_1}}  \; \mathrm{et}  \\ 
\gen{m_2,n_2,p_2}&  =  & \sum_{i_2,j_2,k_2} 
{\; a'_{i_2,j_2}\;  b'_{j_2,k_2}\;  c'_{i_2,k_2}} 
\end{array}$$
(nous avons mis des $'$ pour le cas o\`u les ensembles d'indices
dans le premier tenseur ne seraient pas disjoints de ceux du second)
\emph{\`a condition de respecter les r\`egles de calcul suivantes}
$$\begin{array}{rcll} 
x_{(i_1,i_2),(j_1,j_2)}& = &  x_{i_1,j_1}\;  x'_{i_2,j_2}  & \; 
(x \mathrm{\; vaut \; pour\; }a,\; b\; \mathrm{ou\; }c)\\ 
x_{i_1,j_1}\;  y'_{i_2,j_2}& =  & y'_{i_2,j_2}\; x_{i_1,j_1}   & \; 
(\mathrm{idem\; avec\; }x\; \mathrm{et}\; y\neq x)    \\
\end{array}$$

Une fois ceci constat\'e, nous n'avons m\^eme plus besoin de penser 
au calcul emboit\'e que repr\'esente le produit par blocs. Nos 
nouvelles 
r\`egles de calcul fonctionnent 
toutes seules et produisent automatiquement le r\'esultat (2) 
aussi bien dans la proposition
\ref{propRGtens} que dans la proposition \ref{propRGmarg}.
Nous sommes en effet ramen\'es maintenant 
\`a la constatation banale suivante
concernant les \dlisz: le premier terme 
non nul du produit de deux \dlis est le 
produit des premiers termes non nuls  de chacun des deux 
\dlisz. Et les ordres des deux \dlis s'ajoutent.
\qed

\ss Le raisonnement fait au d\'ebut de ce paragraphe, en tenant 
compte de la proposition \ref{propRGmarg} et du lemme \ref{lemCabas} 
donne alors le r\'esultat de Bini.
\begin{theorem} 
\label{thRGmarg} 
\begin{itemize}
\item [$(1)$] S'il existe \,$n$\,  et \,$r$\, tels que 
\,$\Ra\,\gen{n,n,n}\leq r$\, alors \,$n^\omega\leq r$, \cad 
\,$\omega\leq  {\log r\over\log n}$.
\item [$(2)$] S'il existe \,$m$, \,$n$, \,$p$\, tels que 
\,$\Ra\,\gen{m,n,p}\leq r$\, alors  \,$(mnp)^{\omega/3}\leq r$, \cad
\,$\omega\leq 3\, {\log r\over \log mnp}$. 
\end{itemize}
\end{theorem}

En bref, pour ce qui concerne l'exposant \,$\omega$, 
une \ine
\,$\Ra\, \gen{m, n, p}\alb\leq r$\,
 donne le m\^eme r\'esultat qu'une
\ine \,$R\,\gen{m,n,p}\leq r.$ 

On pourra remarquer que la preuve du \tho \ref{thRGmarg} est 
tout \`a fait explicite. Si on connait un \caba \`a l'ordre \,$q$\, 
qui 
utilise \,$r$\, \muls essentielles pour le produit matriciel 
$\gen{m,n,p}$ et si
\,$\alpha > 3\,{\log r \over \log mnp}$\, alors on sait construire un 
entier \,$N$\, et un \cabi pour le produit matriciel  $\gen{N,N,N}$ 
qui 
utilise moins de \,$N^\alpha$\, \muls essentielles.
\begin{cor} 
\label{corRGmarg} 
On a pour l'exposant de la \mul des \macas  
\,$\omega\leq 3\,\log 10 /\log 12 <2.7799$.
\end{cor}

\subsection[Une am\'elioration d\'ecisive de Sch\"onhage]{Une 
premi\`ere 
am\'elioration d\'ecisive de Sch\"onhage}
\label{paraSchon1}

La \met de Bini n'a pas donn\'e dans un premier temps 
une am\'elio\-ra\-tion tr\`es importante de l'exposant \,$\omega$\, 
mais 
elle a ouvert la voie aux am\'e\-lio\-ra\-tions suivantes, beaucoup 
plus substantielles.
  
Dans la \met de Strassen on remplace pour calculer le produit 
matriciel  $\gen{2,2,2}$  le \cabi avec $8$ \muls (correspondant \`a 
la d\'efinition du produit) par un \cabi avec seulement $7$ \muls 
essentielles, et on obtient \,$\omega\leq 3\,\log7/\log8$. 
Dans la \met de Bini, 
on utilise un produit de matrices \`a trous dans lesquel les $6$ \muls 
qui interviennent dans la d\'efinition du produit matriciel peuvent 
\^etre remplac\'ees (dans un \cabaz) par seulement $5$ \muls 
essentielles. Cependant au lieu d'aboutir \`a
\,$\omega\leq 3\,\log5/\log6$\, comme dans  la \met de Strassen, on 
a abouti \`a \,$\omega\leq 3\,\log10/\log12$.  Sch\"onhage a pens\'e 
qu'il y avait l\`a quelque chose d'immoral et il a obtenu dans un 
travail m\'emorable (voir \cite{Scho2}) l'am\'e\-lio\-ra\-tion 
d\'ecisive suivante.

\begin{theorem} 
\label{thSchBi} Si dans un produit de matrices \`a trous, on est 
capable 
de remplacer, dans un \cabaz, les \,$\rho$\, \muls qui interviennent 
dans la d\'efinition du produit matriciel par seulement \,$\theta$\, 
\muls essentielles, alors \,$\omega\leq 3\,{\log \theta \over \log 
\rho}$. 
En particulier \,$\omega\leq 3\,{\log 5 \over \log 6}\leq 2,695$. 
\end{theorem}

Le reste de ce paragraphe est consacr\'e \`a la preuve de ce 
\thoz, selon les lignes de \cite{Scho2}. 
La preuve est faite sur un corps \,$\K$\, infini, ce qui est 
l\'egitime 
d'apr\`es la proposition \ref{propChangCor}.  Le plus simple est de
commencer sur un exemple. Nous allons voir directement sur l'exemple 
de 
Bini quelle est la machinerie mise en {\oe}uvre par Sch\"onhage. 
La \met it\'erative de Strassen donne des produits matriciels \`a 
trous successifs du type suivant (figures \ref{fbinisch1} et 
\ref{fbinisch2}).
\begin{figure}[h]   
\begin{center}
\includegraphics*[width=6cm]{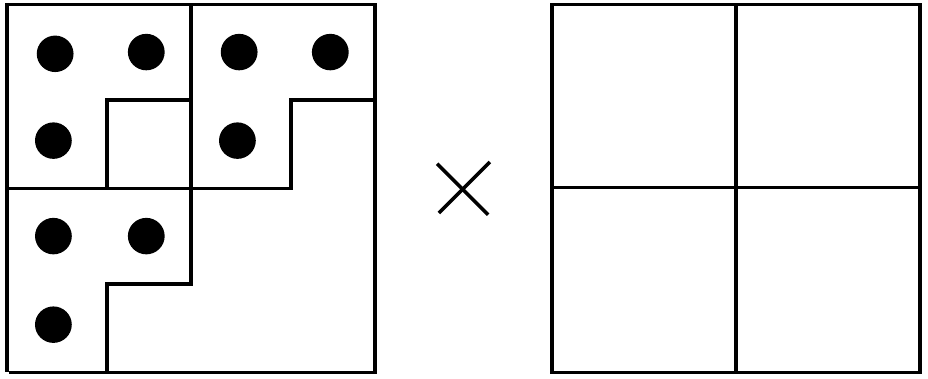}
\end{center}
\caption[Bini it\'er\'e une fois]
{\label{fbinisch1}  Bini it\'er\'e une fois}  
\end{figure}  
Le produit matriciel \`a trous it\'er\'e une fois de la figure 
\ref{fbinisch1} peut \^etre obtenu par un \caba d'ordre $2$ et de rang 
$5^2$ (au lieu de $6^2$). Ceci se d\'emontre comme le point $(2)$ dans 
la proposition \ref{propRGmarg}.
\begin{figure}[h]   
\begin{center}
\includegraphics*[width=10cm]{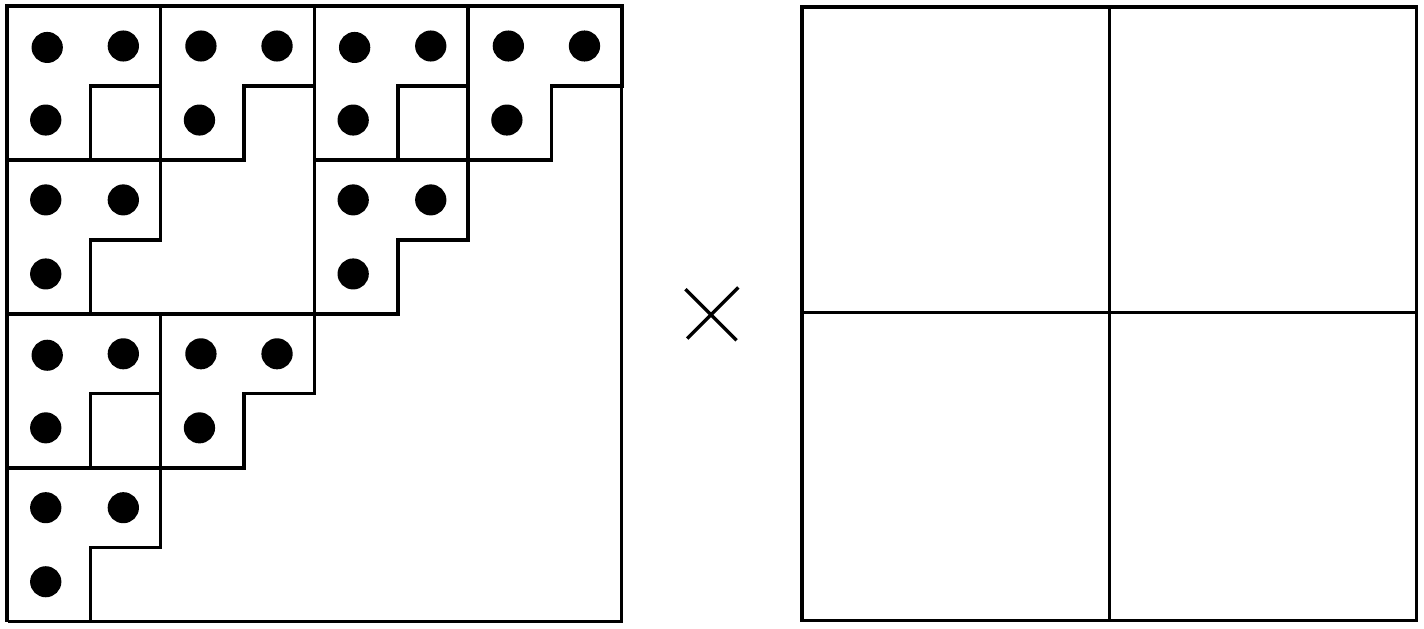}
\end{center}
\caption[Bini it\'er\'e deux fois]
{\label{fbinisch2} Bini it\'er\'e deux fois }  
\end{figure}  
Le produit matriciel \`a trous it\'er\'e deux fois (figure 
\ref{fbinisch2}) peut \^etre obtenu par un \caba d'ordre $3$ et de 
rang 
$5^3$. Ceci se d\'emontre aussi comme le point $(2)$ dans la 
proposition 
\ref{propRGmarg}.

Plus \gnltz, la m\^eme preuve donne.
\begin{proposition} 
\label{propRGmarg2} 
Notons \,$\varphi$\, l'\abi qui correspond \`a un produit matriciel 
\`a 
trou
\,$(A,B)\mapsto A\,B$\, o\`u certaines entr\'ees fix\'ees de \,$A$\, 
et  
\,$B$\, sont nulles et les autres sont des variables in\'ependantes.
Notons \,$\varphi^{\otimes k}$\, le produit matriciel \`a trou obtenu 
en it\'erant \,$k-1$\, fois le produit  \,$\varphi$.
Alors on a
$$ R\left(\varphi^{\otimes k}\right) \le \left(R(\varphi)\right)^k
$$
et
$$ \Ra\left(\varphi^{\otimes k},kq \right) \le 
\left(\Ra(\varphi,q) \right)^k
$$
De m\^eme si on embo\^{\i}te dans un produit par blocs deux produits 
matriciels \`a trous  \,$\varphi$\, et  \,$\psi$\, et qu'on note
 \,$\varphi\otimes \psi$\, le produit matriciel \`a trous que l'on 
obtient,
on a les \ines
$$ R\left(\varphi \otimes \psi\right) \le R(\varphi)\,R(\psi)
$$
et
$$ \Ra\left(\varphi \otimes \psi,q+q'\right) \le 
\Ra(\varphi,q)\,\Ra(\psi,q')
$$
\end{proposition}

Revenons \`a notre exemple. Dans le produit de la figure 
\ref{fbinisch2} 
nous pouvons s\'electionner les $3$ colonnes $2$, $3$, $5$ de la 
premi\`ere matrice, qui contiennent chacune $4$ v\'eritables entr\'ees 
et les lignes $2$, $3$, $5$ de la deuxi\`eme, qui contiennent $8$ 
entr\'ees. On obtient le produit \`a trous 
\,$U\times V=W$\, suivant (figure \ref{fbinisch3}). 
Du point de vue du \cabaz, cette extraction de lignes et de colonnes 
revient simplement \`a remplacer des variables par des $0$ et donc ne 
peut que le simplifier. 
\begin{figure}[ht]   
\begin{center}
\includegraphics*[width=8cm]{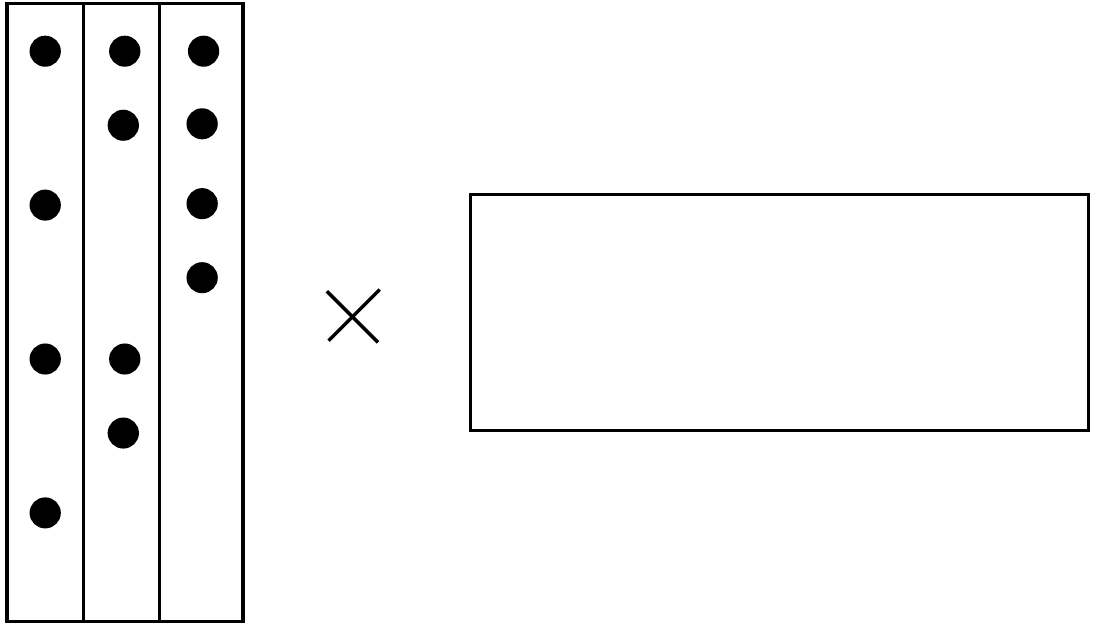}
\end{center}
\caption[Produit \`a trous extrait de \gui{Bini it\'er\'e deux fois}]
{\label{fbinisch3} Produit \`a trous extrait de \gui{Bini it\'er\'e deux 
fois} }  
\end{figure}  
Si nous consid\'erons maintenant une matrice fix\'ee 
\,$G\in\K^{4\times 
8}$\, l'\ali \,$\mu_G:U\mapsto GU$\, est en fait une \ali entre deux 
\evcs de dimension $12$. Admettons un moment que les \coes 
de \,$G$\,  peuvent \^etre choisis de mani\`ere que \,$\mu_G\,$
soit un \iso (lemme de compression). 
En posant \,$GU=U'$\, on voit que le produit 
matriciel \emph{sans trou} \,$\gen{4,3,8}$\, est r\'ealis\'e sous 
forme
\,$(U',V)\mapsto U'\times V$\,  par un \caba d'ordre $3$  et de rang 
$5^3$: d\'ecrypter \,$U'$\, pour obtenir \,$U$\, (sans aucune \mul 
essentielle) puis calculer \,$U\times V$.

De mani\`ere plus g\'en\'erale, nous pouvons consid\'erer le produit 
\`a 
trous \,$A_k\times B_k$\, obtenu en it\'erant \,$k-1$\, fois le 
processus de Bini. La matrice \,$A_k$\,  est de plus en plus creuse. 
Dans chaque colonne, le nombre d'entr\'ees v\'eritables est \'egal \`a 
une puissance de $2$. En s\'electionnant les colonnes ayant un m\^eme 
nombre d'entr\'ees (disons \,$m_k$\, colonnes avec \,$2^{h_k}$\, 
entr\'ees non nulles), on obtient un produit \`a trous  \,$U_k\times 
V_k=W_k$\, 
 \`a l'int\'erieur du format \,$\gen{2^k,m_k,2^k}$. Chaque colonne de
\,$U_k$\, a exactement  \,$2^{h_k}$\, entr\'ees v\'eritables.
En appliquant le lemme de compression, nous choisissons une matrice 
convenable
\,$G_k\in\K^{2^{h_k}\times 2^k}$\, nous rempla\c{c}ons \,$U_k$\, par 
\,$U'_k=G_k\,U_k$\, et nous obtenons un produit matriciel sans trou 
\,$\gen{2^{h_k},m_k,2^k}$\,  sous la forme \,$(U'_k,V_k)\mapsto 
U'_k\times V_k$\,  par un \caba d'ordre \,$k$\,  et de rang \,$5^k$. 

Quel est le comportement asymptotique de ce calcul? 
On peut facilement se convaincre que le produit \,$m_k\,2^{h_k}$\, est 
obtenu comme l'un des termes du \dev de \,$(1+2)^k$\, selon la formule 
du bin\^ome. 
Cela tient \`a ce que la matrice \`a trous initiale 
poss\`ede  une colonne \`a deux entr\'ees et une autre \`a une entr\'ee.
 Comme la formule du bin\^ome est une somme de \,$(k+1)$\, termes, le 
plus grand de ces termes est certainement sup\'erieur \`a
\,$3^k/(k+1)$. 
Donc par un choix optimal de \,$h_k$\, nous obtenons 
$$ N_k=2^{h_k}\cdot m_k\cdot 2^k \ge {6^k \over k+1}
$$
Donc en appliquant la proposition \ref{propRGmarg2}
$$ \Ra\left(\gen{N_k,N_k,N_k},3k\right)
\le 5^{3k}
$$
d'o\`u en appliquant le lemme \ref{lemCabas}
$$ R\gen{{N_k},{N_k},{N_k}}\le 9k^2\,5^{3k}
$$
ce qui donne bien par passage \`a la limite 
$\omega\leq 3\,{\log 5 \over \log 6}$.

Avant de passer \`a la preuve dans le cas g\'en\'eral, nous montrons 
le lemme de compression.
\begin{lemma} 
\label{lemCompression} \emph{(lemme de compression)} 
Soit \,$A=(a_{ij})$\, une \gui{matrice \`a trous} de format 
\,$m\times n$\, 
dont les entr\'ees sont ou bien nulles ou bien des variables 
ind\'ependantes. Nous supposons que la matrice poss\`ede \,$p$\, 
variables et \,$m-p$\, entr\'ees nulles dans chaque colonne. 
Si on sp\'ecialise les variables dans le corps \,$\K$\, on obtient un 
\evc \,$E_A$\, de dimension \,$n\,p$. 
On suppose le corps \,$\K$\, infini. 
Alors il existe une matrice $G\in \K^{p\times m}$ telle que l'\ali
$$ 
\mu_{G,A}~:~E_A\longrightarrow \K^{p\times n}\; ,\quad U\longmapsto 
G\,U
$$
soit bijective.
\end{lemma}
\preuve Les colonnes de \,$U$\, sont transform\'ees ind\'ependamment 
les 
unes des autres. Chaque colonne \,$U_j$\, de la matrice \,$U$, en ne 
gardant que les entr\'ees non nulles, est transform\'ee selon le 
sch\'ema
 $$\,U_j\,\longmapsto\, G_j\,U_j\,$$
o\`u \,$G_j$\, est une \maca extraite de \,$G$\, en ne gardant que 
\,$p$\, colonnes de \,$G$. Il s'ensuit que l'\ali \,$\mu_{G,A}$\, est 
bijective \ssi les matrices \,$G_j$\, sont inversibles. Pour cela, il 
suffit que \,$p$\, colonnes distinctes de \,$G$\, soient toujours 
ind\'ependantes. 
Ce probl\`eme de g\'eom\'etrie combinatoire admet toujours une 
solution 
sur un corps ayant suffisamment d'\eltsz. Si on a d\'ej\`a construit 
une matrice convenable \,$G$\, avec  \,$\ell\geq p$\, colonnes, pour 
rajouter une colonne,
il faut choisir un vecteur en dehors des hyperplans d\'efinis par 
n'importe quel \sys de  \,$p-1$\, colonnes extraites de  \,$G$. 
\qed

\ss Passons maintenant \`a la preuve du cas g\'en\'eral.
Nous supposons que nous avons un produit de matrice \`a trous 
$ \,A\times B\,$
par exemple du style suivant (figure \ref{fbinisch4})
\begin{figure}[ht]   
\begin{center}
\includegraphics*[width=6cm]{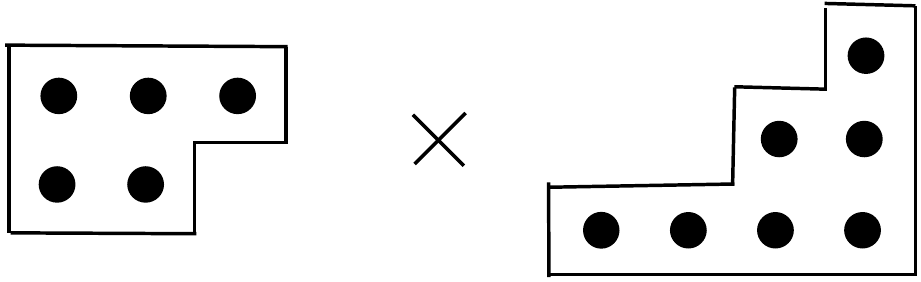}
\end{center}
\caption[Un exemple arbitraire de produit matriciel \`a trous]
{\label{fbinisch4} Un exemple arbitraire de produit matriciel \`a 
trous 
}  
\end{figure}  
qui peut \^etre r\'ealis\'e par un \caba de mani\`ere \'economique.
Supposons que les colonnes successives de \,$A$, au nombre de \,$t$\, 
contiennent respectivement \,$m_1$, \,$m_2,\ldots $, \,$m_t$\, 
entr\'ees 
v\'eritables. Supposons que les \,$t$\, lignes successives de \,$B$\, 
contiennent respectivement \,$n_1$, \,$n_2,\ldots $, \,$n_t$\, 
entr\'ees 
v\'eritables. 
A priori ce produit \`a trous r\'eclame 
$$ \rho=m_1n_1+\cdots+m_tn_t
$$
\muls : le tenseur qui correspond \`a sa d\'efinition est une somme de 
\,$\rho$\, tenseurs \elrs (dans l'exemple ci-dessus, \,$t=3$, 
$(m_1,m_2,\alb m_3)\alb =(2,2,1)$, $(n_1,n_2,n_3)=(1,2,4)$ et 
\,$\rho=10$).
 Supposons qu'un \caba \`a l'ordre \,$\ell$\, et de rang \,$\theta$\, 
permette de r\'ealiser ce produit \`a trous. Si on it\`ere $k-1$ fois 
ce 
\cabaz, on obtient un nouveau produit de matrices 
\`a trous \,$A_k\times B_k$. 
Par exemple pour $A_2\times B_2$, on obtient le produit \`a trous 
suivant (figure \ref{fbinisch5}).
\begin{figure}[ht]   
\begin{center}
\includegraphics*[width=12cm]{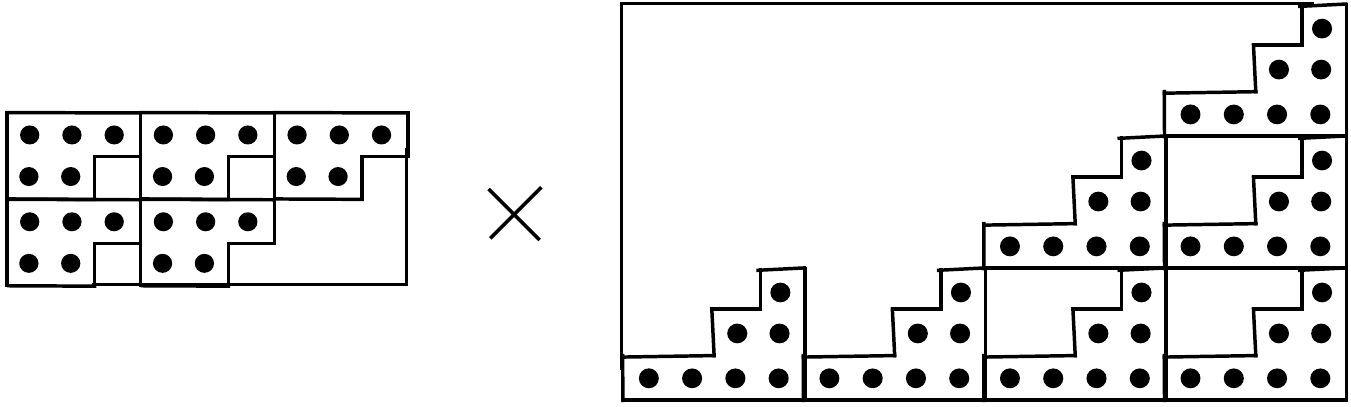}
\end{center}
\caption[L'exemple pr\'ec\'edent it\'er\'e une fois]
{\label{fbinisch5} L'exemple pr\'ec\'edent it\'er\'e une fois }  
\end{figure}  
Une colonne de \,$A_k$\, doit \^etre index\'ee par un \,$k$-uple 
\label{IndicesDesCol} \,$j=(j_1,\ldots,j_k)$\, d'\elts de 
\,$\{1,\ldots,t\}$. Une telle colonne contient alors
$$ m_{j_1}\cdots\, m_{j_k}=m_1^{u_1}\cdots\, m_t^{u_t}
$$
entr\'ees non nulles, o\`u chaque \,$u_i$\, est \'egal au nombre des
\,$j_s$\, \'egaux \`a \,$i$.
De m\^eme une ligne de \,$B_k$\, doit \^etre index\'ee par un 
\,$k$-uple \,$j=(j_1,\ldots,j_k)$\, d'\elts de 
\,$\{1,\ldots,t\}$\, et elle contient \,$n_{j_1}\cdots\, n_{j_k}$\, 
entr\'ees non nulles.
Parmi toutes les colonnes de \,$A_k$\,  on d\'ecide de 
s\'elec\-tion\-ner toutes celles qui fournissent une certaine liste 
d'exposants $(u_1,\ldots,u_t)$. 
En particulier elles ont toutes le m\^eme  nombre 
\,$\mu_k=m_1^{u_1}\cdots m_r^{u_t}$\, d'entr\'ees non nulles (avec
\,$u_1+\cdots+u_t=k$). Le nombre des colonnes en question est \'egal 
au 
\coe multinomial
$$ \lambda_k = {k \choose u_1,\ldots,u_t}= 
{k\,!\over u_1!\,\cdots\, u_t!}
$$
De m\^eme, nous s\'electionnons parmi les lignes de \,$B_k$\, toutes 
celles correspondant aux m\^emes indices (qui sont des $k$-uples 
\,$j=(j_1,\ldots,j_k)$). Elles ont toutes le m\^eme nombre d'entr\'ees 
non nulles \,$\nu_k=n_1^{u_1}\cdots\, n_t^{u_t}$.
Nous obtenons de cette mani\`ere un produit de matrice \`a trous
\,$U_k\times V_k$. Comme les colonnes de \,$U_k$\, ont toutes le 
m\^eme 
nombre \,$\mu_k$\, d'entr\'ees non nulles, on peut utiliser le lemme 
de 
compression. M\^eme chose pour \,$V_k$\, en tenant compte du fait que 
toutes les lignes ont le m\^eme nombre \,$\nu_k$ d'entr\'ees non 
nulles. 
En d\'efinitive nous obtenons un produit matriciel sans trou de type 
\,$\gen{\mu_k,\lambda_k,\nu_k }$\,  qui est r\'ealis\'e par un \caba 
d'ordre \,$k\ell$\, et de rang $\leq \theta^k$.\\
A quoi est \'egal  \,$\mu_k\cdot\lambda_k\cdot\nu_k \,$~? C'est l'un 
des 
termes du \dev multinomial de 
\,$\left(m_1n_1+\cdots+m_tn_t\right)^k=\rho^k$. Si on choisit le terme 
le plus grand dans cette somme on obtient donc
$$ M_k=\mu_k\cdot\lambda_k\cdot\nu_k \geq 
{\rho^k\over {k+t \choose t-1}}
\geq {\rho^k\over {(k+1)^{t-1}}} 
$$
car il y a \,${k+t \choose t-1}$\, termes dans cette somme.
On termine comme dans le cas particulier examin\'e au d\'ebut:
$$ \Ra\left(\gen{M_k,M_k,M_k},3k\ell\right)
\le \theta^{3k},\qquad R\,\gen{M_k,M_k,M_k}
\le 9\,k^2\,\ell^2\,\theta^{3k}
$$
et par passage \`a la limite en appliquant le \tho \ref{thRgExp},
\,$\omega\leq 3\,{\log \theta \over \log \rho}$.

\begin{remark} \label{remBinSch}
\emph{ Dans \cite{Scho2} Sch\"onhage indique des produits matriciels 
\`a 
trous avec un rang marginal plus avantageux que celui de Bini, ce qui 
donne \,$\omega\leq 2,6087$. Mais ce dernier r\'esultat est surpass\'e 
par la formule asymptotique qu'il obtient ensuite et que nous exposons 
dans le paragraphe suivant.
} 
\end{remark}

\label{remBinSch2}
\begin{remark} 
\emph{Dans le lemme \ref{lemCabas} il est possible de remplacer 
\,$(q+1)(q+2)/2$\, par \,$1+6q$. M\^eme avec cette am\'elioration, 
c'est 
uniquement pour des entiers \,$N$\, tr\`es grands que la \pcd de 
Bini aussi bien que celle de Sch\"onhage fournissent un meilleur \cabi  
pour \,$\gen{N,N,N}$\,  que celui qui d\'ecoule de la \pcd 
originale de Strassen. Ces \mets ne sont donc pas impl\'ement\'ees sur 
machine. 
} 
\end{remark}

\subsection{Sommes directes d'\abis}
Approfondissant son analyse des produits de matrices \`a 
trous, Sch\"on\-ha\-ge a remarqu\'e que certains produits du type 
ci-dessous (figure \ref{fbinisch6}) permettent de construire \`a 
partir 
d'un \caba des \cabis
exacts donnant un meilleur exposant pour la \mul des \macas que celui 
\'etabli dans le \tho \ref{thSchBi}.  
\begin{figure}[ht]   
\begin{center}
\includegraphics*[width=6cm]{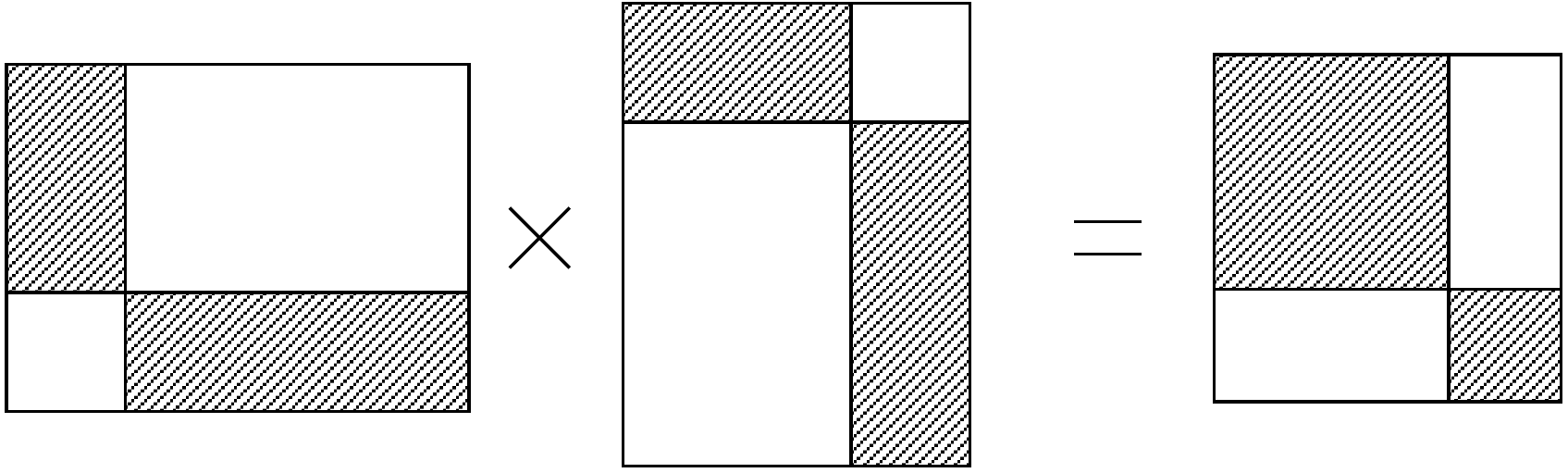}
\end{center}
\caption[Somme directe de deux produits matriciels]
{\label{fbinisch6} Somme directe de deux produits matriciels }  
\end{figure}  
L'exemple de la figure \ref{fbinisch6}  correspond \`a la 
\emph{somme disjointe}\index{somme disjointe (de deux \abisz)} (on 
peut 
dire aussi \emph{somme directe}
ou encore \emph{juxtaposition}) des
deux \abis \,$\gen{1,2,1}$\, et \,$\gen{3,1,3}$.
De mani\`ere g\'en\'erale, la somme directe de deux \abis\index{somme 
directe (de deux \abisz)}
$ \,\varphi_1:E_1\times F_1\rightarrow G_1$\,  et
$ \,\varphi_2:E_2\times F_2\rightarrow G_2$\,  est l'\abi
$$\,\varphi:(E_1\oplus E_2)\times (F_1\oplus F_2)\,\longrightarrow\, 
(G_1\oplus G_2)\,$$ 
d\'efinie par 
\,$\varphi((x_1,y_1),(x_2,y_2))=(\varphi_1(x_1,y_1),
\varphi_2(x_2,y_2))$. Du point de vue des \cabisz,
un \cabi possible pour la somme disjointe consiste \`a faire seulement
les deux calculs en \paral avec toutes les variables distinctes.
%
\begin{notation} 
\label{notaSomAbi}
On note \,$\varphi_1\oplus\varphi_2$\, la somme directe des \abis $ 
\,\varphi_1$\, et  $ \,\varphi_2$.
On note  \,$\ell\odot\varphi$\, pour la somme directe de
$ \,\ell$\, exemplaires de $ \,\varphi$. 
\end{notation}
On fait alors les remarques suivantes. Le premier lemme est \`a la 
fois 
simple et crucial.
\begin{lemma} 
\label{lemSdir1}~ Supposons \,$R\,\gen{f,f,f}\leq s$\, et 
\,$R(s\odot\gen{m,n,p})\leq r$. Alors
\,$R\,\gen{fm,fn,fp}\leq r$.   
\end{lemma}
\prv
L'\abi \,$\gen{fm,fn,fp}$\, peut \^etre r\'ealis\'ee comme un produit 
par blocs, chacune des deux matrices  \,$A$\, et \,$B$\,  
qu'on multiplie \'etant d\'ecoup\'ee en  \,$f^2$\, blocs de m\^eme 
format. 
Les  \,$f^3$\, \muls correspondantes de type \,$\gen{m,n,p}$\, qui 
sont 
a priori \ncrs pour ce produit par blocs peuvent \^etre 
remplac\'ees par seulement
\,$s$\, produits (entre \colis convenables des blocs), selon le 
sch\'ema 
fourni par le \cabi qui montre \,$R\,\gen{f,f,f}\leq s.$
\qed
\begin{lemma} 
\label{lemSdir2} ~\\ 
$(1)$ $ R(ss'\odot\gen{mm',nn',pp'}) \le R(s\odot\gen{m,n,p})\cdot 
R(s'\odot\gen{m',n',p'})$.  En particulier
avec $N=mnp$ on a \\
\centerline{$R(s^3\odot\gen{N,N,N}) \le 
\left(R(s\odot\gen{m,n,p})\right)^3$  }
\\
$(2)$\\
 $\Ra(ss'\odot\gen{mm',nn',pp'},q+q') \le 
\Ra(s\odot\gen{m,n,p},q)\cdot 
\Ra(s'\odot\gen{m',n',p'},q')$. 
En particulier avec $N=mnp$ on a\\ 
\centerline{ $\Ra(s^3\odot\gen{N,N,N},3q) \le 
\left(\Ra(s\odot\gen{m,n,p},q)\right)^3$ }
\end{lemma}
\prv C'est toujours la \met du produit par blocs,
appliqu\'ee avec les produits matriciels \`a trous correspondants.
On peut consid\'erer qu'il s'agit d'un cas particulier de la 
proposition 
\ref{propRGmarg2}.
\qed

\ss
On en d\'eduit la proposition suivante qui \gni le 
\tho \ref{thRGmarg}.
\begin{proposition} 
\label{propSdir} 
S'il existe \,$s$, \,$m$, \,$n$,  \,$p$,  \,$r$\, tels que 
\,$\Ra(s\odot\gen{m,n,p})\leq r$\, alors  
\,$s\,(mnp)^{\omega/3}\leq r$, \cad
\,$\omega\leq 3\, {\log (r/s)\over \log mnp}$. 
\end{proposition}
Autrement dit, pour ce qui concerne l'exposant \,$\omega$, 
l'\ine
\,$\Ra(s\odot\gen{m,n,p})\leq r$\, donne le m\^eme r\'esultat qu'une
\ine \,$R\,\gen{m,n,p}\leq r/s.$ 

\ss\prv  
Si \,$\Ra(s\odot\gen{m,n,p},q)\leq r$, en appliquant le lemme 
\ref{lemSdir2} on obtient avec \,$N=mnp$\, 
$$\Ra\left(s^3\odot\gen{N,N,N},3q\right)\le r^3$$ 
puis aussi 
$$ \Ra\left(s^{3\ell}\odot\gen{N^\ell,N^\ell,N^\ell},3\ell q\right)\le 
r^{3\ell}$$  et donc 
\begin{equation} \label{eqSdir1}
R\left(s^{3\ell}\odot\gen{N^\ell,N^\ell,N^\ell}\right)\le 
9\,\ell^2\,q^2 
r^{3\ell}\,.
\end{equation}
Par passage \`a la limite, cela nous ram\`ene 
au cas o\`u on conna\^{\i}t des entiers  \,$s$, \,$m$, \,$r\,$
tels que  \,$R(s\odot\gen{m,m,m})\le r$. 
On veut alors montrer \,$\omega\leq  {\log (r/s)\over \log m}$. 
Posons \,$\l=\, {\log (r/s)\over \log m}$.\\
Supposons tout d'abord qu'on connaisse un \cabi qui montre que
 \,$R\,\gen{f,f,f}\leq s$\, et posons 
\,$\alpha_0=\log s/\log f$\, (${\alpha_0}$\, est un exposant 
acceptable). 
Si \,${\alpha_0}\leq \l$\, on n'a rien \`a faire. 
Si  \,${\alpha_0}> \l$\, 
le lemme \ref{lemSdir2} nous dit que \,$R\,\gen{fm,fm,fm}\leq r$.
Donc l'exposant 
$${\alpha_1}= \log r/\log fm= 
{\alpha_0}\log r/(\log s+{\alpha_0}\log m)$$  
est acceptable pour la \mul des \macasz. Un calcul simple montre
alors que \,$\l<{\alpha_1}<{\alpha_0}$. Donc nous avons am\'elior\'e 
la 
situation en passant de \,${\alpha_0}$\, \`a \,${\alpha_1}$.\\
Nous voyons maintenant le travail qui nous reste \`a faire.\\
Primo, montrer que si on a  
\,$R\,\gen{f,f,f}\leq  s'=f^{\alpha_0}$\, avec un entier \,$s'\neq s$, 
cela n'est pas trop grave, car on peut utiliser
 \,$R\,s^\ell\odot\gen{m^\ell,m^\ell,m^\ell}\le r^\ell$\,  et
\,$R\,\gen{f^k,f^k,f^k}\leq  
s'^k=\left(f^k\right)^{\alpha_0}$\,  avec
\,$s'^k\leq s^\ell$\, et le rapport de leurs logarithmes
aussi proche qu'on veut de 1. Donc par le lemme \ref{lemSdir2} 
\,$R\,\gen{f^km^\ell,f^km^\ell,f^km^\ell}\leq r^\ell$,
ce qui conduit \`a un exposant acceptable 
$${\alpha'_1}= \log r\ell/\log f^km^\ell= 
{\alpha'_0}\log r/(\log s+{\alpha'_0}\log m)$$
avec \,${\alpha'_0}$\, aussi proche qu'on veut de \,${\alpha_0}$. \\
Secundo, montrer que si on recommence, les exposants successifs
\,${\alpha_n}\,$
qu'on obtient convergent bien vers \,$\l$.\\
Nous ne ferons pas ce travail, car les d\'etails techniques deviennent 
vraiment trop lourds.
\qed

\subsubsection*{La conjecture additive de Strassen} 
On a \'evidemment
$$\begin{array}{rcll} 
R(\varphi_1\oplus\varphi_2)& \leq   & R(\varphi_1)+R(\varphi_2)  &   
\quad \quad \mathrm{et}   \\[2mm] 
\Ra(\varphi_1\oplus\varphi_2,q)& \leq  & 
\Ra(\varphi_1,q)+\Ra(\varphi_2,q)\,.  
\end{array}$$
Une conjecture de Strassen est que la premi\`ere \ine est en 
fait toujours une \egtz. On appelle cette conjecture \emph{la 
conjecture additive (pour le \rgte des \abisz)}\index{conjecture 
addititve}.
Bien que plausible, cette conjecture a \'et\'e \'ebranl\'ee par 
Sch\"onhage qui a montr\'e que la variante avec \gui{\rgtmz} \`a la 
palce 
de \gui{\rgtez} est fausse, d'apr\`es le r\'esultat du lemme 
\ref{lemSdirApp}.
Si la conjecture additive est vraie, ou m\^eme si seulement 
\,$R(s\odot\gen{m,m,m})= s\, R\,\gen{m,m,m}$  pour \,$s$\, et  \,$m$\, 
arbitraires, la preuve de la proposition \ref{propSdir}  est beaucoup 
simplifi\'ee, car on d\'eduit de l'\'equation \ref{eqSdir1}  
directement
$ R\left(\gen{N^\ell,N^\ell,N^\ell}\right)\le 9\,\ell^2\,q^2 
(r/s)^{3\ell}$.
Mais cela ne fournirait les \cabis  demand\'es que si on \'etait 
capable
de trouver un \cabi de rang con\-ve\-na\-ble pour \,$\gen{m,m,m}$\, 
\`a 
partir d'un \cabi pour \,$s\odot\gen{m,m,m}.$

\begin{lemma} 
\label{lemSdirApp} 
Pour \,$k>1$\, et \,$m=(k-1)^2$:
$$ \Ra(\gen{k,1,k}\oplus\gen{1,m,1})=k^2+1<  
k^2+m=\Ra\,\gen{k,1,k}+\Ra\,\gen{1,m,1}
$$
\end{lemma}
\prv Nous montrons seulement \,$\Ra(\gen{k,1,k}\oplus\gen{1,m,1})\le 
k^2+1$.
Nous repr\'esentons le produit \,$\gen{k,1,k}$\,  par le \pol non 
commutatif
\,$\sum_{i=1}^k {a_ib_jc_{j,i}}$\, et le produit \,$\gen{1,m,1}$\,  
par 
le \pol non commutatif
\,$\sum_{\ell=1}^m {u_\ell v_\ell w}$. Pour simplifier les \'ecritures 
qui suivent, nous prenons \,$\ell=(i,j)$\, avec \,$1\le i,j\le k-1$.
Nous introduisons en outre les notations 
$$ u_{i,k}=v_{k,j}=0,\quad u_{k,j}=-\sum_{i=1}^{k-1} u_{i,j}
,\quad v_{i,k}=-\sum_{j=1}^{k-1} v_{i,j}
$$
de sorte que
$$ \sum_{\ell=1}^m {u_\ell v_\ell}=\sum_{i=1,j=1}^{k-1}u_{i,j}v_{i,j}=
\sum_{i=1,j=1}^{k}u_{i,j}v_{i,j}
$$
On consid\`ere alors le \pol non commutatif 
suivant (qui correspond \`a un \caba
avec \,$k^2+1$\, \muls essentielles)
$$ \sum_{i=1,j=1}^k {(a_i+\vep u_{i,j})\,(b_j+\vep v_{i,j})\,
(\vep^2 c_{j,i}+w)} - \left(\sum_{i=1}^k a_i \right) 
\left(\sum_{j=1}^k b_j \right)w
$$
qui, une fois d\'evelopp\'e donne
$$ \vep^2\,\left(\sum_{i=1,j=1}^k 
\left(a_i b_j c_{j,i} +u_{i,j}v_{i,j}w\right)\right) + \vep^3 Q\,.
$$
\qed

\subsection{L'in\'egalit\'e asymptotique de Sch\"onhage}

\index{Sch\"onhage!in\'egalit\'e asymptotique}
Revenons au produit \`a trou de la figure  \ref{fbinisch6} 
qui repr\'esente la juxtaposition \,$\gen{2,1,2}\oplus\gen{1,3,1}$. Si 
nous it\'erons 
une fois (\`a la Strassen) ce produit \`a trous, 
nous obtenons un nouveau
produit \`a trou correspondant \`a la figure \ref{fbinisch7},  
\begin{figure}[ht]   
\begin{center}
\includegraphics*[width=10cm]{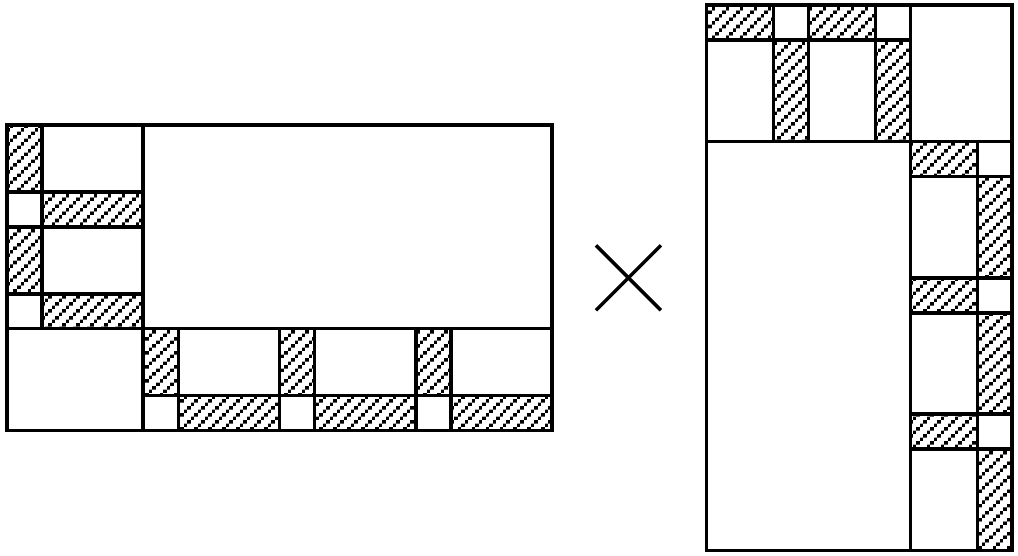}
\end{center}
\caption[Somme directe it\'er\'ee]
{\label{fbinisch7} Somme directe,  it\'er\'ee une fois, de deux 
produits 
matriciels }  
\end{figure}  
qui peut \^etre r\'eorganis\'e, par changement de num\'erotation des 
lignes et colonnes, en le produit \`a trou  qui correspond \`a la 
figure   
\ref{fbinisch8},
\begin{figure}[ht]   
\begin{center}
\includegraphics*[width=10cm]{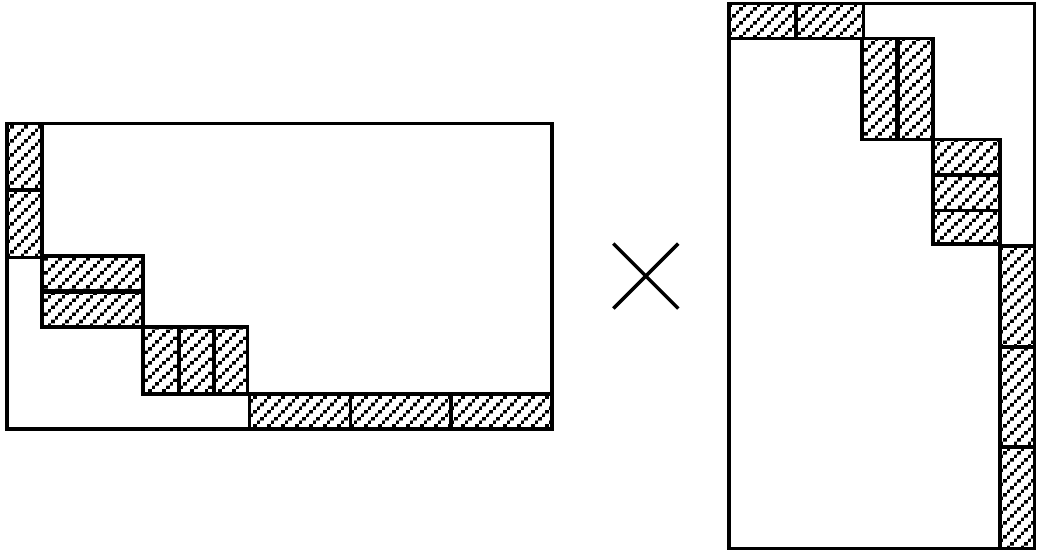}
\end{center}
\caption[Somme directe it\'er\'ee, variante]
{\label{fbinisch8} Somme directe, it\'er\'ee une fois et 
r\'eorganis\'ee 
}  
\end{figure}  
et nous voyons clairement que cela signifie
$$ \left(\gen{2,1,2}\oplus\gen{1,3,1}\right)^{\otimes2}\simeq
\gen{4,1,4}\,\oplus\,2\odot\gen{2,3,2}\,\oplus\,\gen{1,9,1}
$$
Le lecteur ou la lectrice est invit\'ee \`a r\'ealiser par 
elle-m\^eme l'it\'eration une deuxi\`eme fois, et \`a v\'erifier que
$$ \left(\gen{2,1,2}\oplus\gen{1,3,1}\right)^{\otimes3}\simeq
\gen{8,1,8}\oplus3\odot\gen{4,3,4}\oplus3\odot\gen{2,9,2}
\oplus\gen{1,27,1}
$$
avec la parent\'e \'evidente avec la formule du bin\^ome.
Cette parent\'e n'est pas un hasard. 
C'est bien la m\^eme machinerie combinatoire qui est \`a 
l'{\oe}uvre dans les deux cas.
En it\'erant \,$k-1$\, fois on obtiendra
$$ \left(\gen{2,1,2}\oplus\gen{1,3,1}\right)^{\otimes k}\simeq
\sum_{i=1}^k{k \choose i}\odot\gen{2^i,3^{k-i},2^i} \
$$
o\`u le \,$\sum$\,  indique une somme disjointe
d'\abisz.
En fait nous avons une formule du multin\^ome
\gnlez, o\`u les sommes indiquent des sommes disjointes
d'\abisz, et o\`u l'\iso correspond \`a une organisation convenable 
des 
lignes et colonnes du produit matriciel \`a trous correspondant
au premier membre
$$ \left(\sum_{i=1}^t\gen{m_i,n_i,p_i}\right)^{\otimes k}\simeq
\sum_{(u_1,\ldots,u_t)}{k\choose u_1,\ldots,u_t}\odot
\gen{\prod_{i=1}^tm_i^{u_i},
\prod_{i=1}^tn_i^{u_i},\prod_{i=1}^tp_i^{u_i}}
$$
(la deuxi\`eme somme est prise sur tous les \,$t$-uples 
\,$(u_1,\ldots,u_t)$\, tels que \,$\sum u_i=k$).
On en d\'eduit la formule asymtotique suivante.
\begin{theorem} 
\label{thFormAsSch} \emph{(formule asymptotique de Sch\"onhage)} 
Supposons qu'on ait 
$$ \Ra\left(\bigoplus_{i=1}^t\gen{m_i,n_i,p_i}\right)\le r
{\rm\quad et\quad   } 
\sum_{i=1}^t{(m_i\,n_i\,p_i)^\beta}= r
$$
Alors on obtient pour l'exposant de la \mul des \macas 
\,$\omega\le 3\,\beta.$
\end{theorem}
\prv  Notons d'abord que le \tho \ref{thSchBi} donne
$$\,2\le\omega\le {\log r \over \log \left(\sum_{i=1}^t 
m_in_ip_i)\right)}\,$$ 
donc \,$\log r \ge 2 \log t.$
En appliquant  la formule du multin\^ome et l'\ine
de la proposition \ref{propRGmarg2}  on obtient 
$$ \Ra\left({k\choose u_1,\ldots,u_t}\odot
\gen{\prod_{i=1}^tm_i^{u_i},
\prod_{i=1}^tn_i^{u_i},\prod_{i=1}^tp_i^{u_i}}\right)\leq r^k
$$
Pour un choix particulier de \,$u_1,\ldots,u_t$, nous notons ceci
sous la forme
$$ \Ra(S_k\odot\gen{M_k,N_k,P_k})\leq r^k
$$
ce qui nous donne, d'apr\`es la proposition \ref{propSdir}
$$\,\omega\leq 3\, {\log (r^k/S_k)\over \log M_kN_kP_k}.$$
Quel est le choix optimal de  \,$u_1,\ldots,u_t$~?
Nous consid\'erons l'\egt
$$ \left(\sum_{i=1}^t{(m_i\,n_i\,p_i)^\beta} \right)^k = r^k=
\sum_{(u_1,\ldots,u_t)}{{k\choose u_1,\ldots,u_t}\,
\left(\prod_{i=1}^tm_i^{u_i}\,n_i^{u_i}\,p_i^{u_i}\right)^\beta}
$$
La somme de droite a 
\,${k+t-1 \choose t-1}\le {r^k\over (k+1)^{t-1}}$\, 
termes et donc pour le plus grand d'entre eux on obtient
$$ {k\choose u_1,\ldots,u_t}\,
\left(\prod_{i=1}^tm_i^{u_i}\,
n_i^{u_i}\,p_i^{u_i}\right)^\beta=
S_k\,\left(M_kN_kP_k\right)^\beta 
\ge {r^k\over (k+1)^{t-1}} 
$$
ce qui donne
$$ \omega\leq 3\, {\log \left(r^k/(S_k(k+1)^{t-1})\right)\over 
\log(M_kN_kP_k)}\le
3\,\beta + 3\,{ (t-1) \log(k+1)\over \log(M_kN_kP_k) }
$$
D'o\`u le r\'esultat par passage \`a la limite: 
$\beta\,\log(M_kN_kP_k)\,$
est \'equivalent \`a \,$\log r^k/S_k$, et on a 
\,$S_k<t^k$\, et \,$r\ge t^2$.
\qed

\begin{corollary} 
\label{corFAS} 
L'exposant de la \mul des \macas v\'erifie
\,$\omega\le 2,5479.$
\end{corollary}
\prv
On applique la formule asymptotique avec la somme disjointe
\,$\gen{4,1,4}\oplus\gen{1,9,1}$\,  du lemme \ref{lemSdirApp}.
\qed


\newpage \thispagestyle{empty}


\chapter{Alg\`ebre lin\'eaire \sqle rapide}
\label{chap AlinCorps}
\minitoc
\acvide
 
\subsubsection*{Introduction}
 
Une cons\'equence importante de la \mul rapide des
matrices est la recherche de \mets de calcul permettant
de ramener les \pbs classiques d'\agr \lin
\`a une \coag du m\^eme ordre que celle
de la \mul des matrices.
 
\ss Bien que nous utilisions syst\'ematiquement
la \mul rapide des matrices, qui est
obtenue par un \algo tr\`es bien \parasz,
les \algos obtenus dans ce chapitre
ne sont pas eux-m\^emes bien \parassz.
Leur \prof est en \gnl en $\O(n)$, ce qui
explique le titre du chapitre (\agr \lin \emph{\sqlez} rapide)
 
\ss Nous avons d\'ej\`a vu \`a la section \ref{invtrian}
 que  l'inverse
d'une \matg d'ordre \,$n$\, peut se calculer par une
\famu de \caris de taille
\,$\O(n^\alpha)$\, et de \prof \,$\O(\log^2{n})$.
 
\ss Nous allons dans ce chapitre montrer que, pour autant
qu'on travaille sur un corps et qu'on ait droit \`a la
division\footnote{~Ceci est l\'egitime si la division n'est pas
trop co\^uteuse en termes de \cobiz.},
des \fams de \caris ayant des tailles voisines
peuvent \^etre construites pour r\'esoudre les principaux
\pbs de l'\agr \lin sur un corps. Mais dans tous les
\algos que nous exhiberons, le temps \paral
(la \prof du \cirz) n'est plus \pogz.
 
En outre, comme ce sont des \cirs avec divisions, ils ne peuvent
pas \^{e}tre ex\'ecut\'es sur toutes les entr\'ees, et nous donnerons
en g\'en\'eral une version sous la forme d'un \algo \gui{avec
branchements} (les branchements sont gouvern\'es par des tests
d'\egt \`a 0 dans le corps). Dans ces \algos (qui ne correspondent
plus \`a des \caris proprement dits), nous aurons pour le temps
\sql et le temps \paral des estimations tr\`es voisines de celles
obtenues pour les \caris avec divisions.
 
Par exemple le calcul du \deter et de l'inverse d'une matrice
carr\'ee (si elle est inversible) peuvent \^etre r\'ealis\'es
par une \famu de \caris avec divisions de taille
$\,\O(n^{\alpha})$\, (voir section \ref{sec detinv}).
Ceci est une cons\'e\-quence de
l'\algo \sql de Bunch \& Hopcroft pour la \LUP-\deco que
nous d\'eveloppons dans la section \ref{sec BUHO}.
Cet \algo se pr\'esente naturellement sous la forme d'un
\algo avec branchements.
 
En ce qui concerne le calcul du \polcar plusieurs
\mets d'acc\'el\'eration de l'\algo de Frobenius
(section \ref{subsec.Melr}) assez
sophistiqu\'ees ont \'et\'e mises au point par Keller-Gehrig.
L'\algo avec branchements,
qui utilise un temps \sql en \,$\O(n^\alpha\,\log{n})$\, n\'ecessite au
pr\'ealable une \met rapide
pour la mise en forme \gui{\'eche\-lon\-n\'ee en lignes} d'une matrice
arbitraire.
Ceci est expliqu\'e dans les sections \ref{sec FEL} et
\ref{sec kellerseq}.
 
Dans la derni\`ere section \ref{sec kalto}, nous quittons le cadre
de l'\agr \lin sur les corps, mais nous restons dans
celui de l'\agr \lin \sqle acc\'el\'er\'ee
gr\^ace \`a la \mul rapide des matrices.
Nous d\'ecrivons la \met de
Kaltofen, inspir\'ee de l'\algo probabiliste de Wiedemann,
tr\`es efficace pour les matrices creuses sur des corps finis.
Elle donne le meilleur temps \sql actuellement connu pour
le calcul du d\'eterminant, du \polcar et de l'adjointe d'une
\maca \emph{sur un \acomaz}.
L'\algo utilise la \mul rapide des \pols et celle des matrices.
Contrairement \`a l'\algo de Wiedemann, celui de Kaltofen n'a
cependant pas encore fait l'objet d'une impl\'ementation satisfaisante.

\ss Si les \algos d\'evelopp\'es dans ce chapitre
 sont th\'eoriquement plus rapides
que les \algos \gui{\uslsz} donn\'es au chapitre
\ref{chap BasicAlgoAlin}, il y a encore
malheureusement loin de la th\'eorie \`a la pratique.
En fait seule la premi\`ere forme
de la \mul rapide des matrices (celle de Strassen,
correspondant \`a \,$\alpha = \log{7}\simeq 2,807$) commence \`a
\^etre impl\'ement\'ee.
Outre la difficult\'e pratique d'impl\'ementer
d'autres \algos de \mul rapide des matrices,
les \coes \,$\Ca$\, pour de meilleures valeurs
de \,$\alpha\,$ sont trop grands.
Leur impl\'ementation  ne se r\'ev\`elerait efficace
que pour des matrices de tailles astronomiques.

\section[L'Algorithme de Bunch \& Hopcroft]{L'Algorithme de Bunch \&
Hopcroft \\ pour la \LUP-\deco \\ des matrices surjectives}
\label{sec BUHO}
 
Dans la section \ref{subsec Gauss}
nous avons pr\'esent\'e l'\algo \vref{LUPdeco}
qui est l'\algo \sql \usl
(par la \mpgz) pour la \LUP-\deco  des matrices
surjectives.
La \pcd de la \LUP-\deco que
nous allons d\'evelopper ici fait appel \`a la
\mul rapide des matrices.
Cette \pcd, que nous noterons {\bf lup},
est due \`a Bunch et Hopcroft \cite{lup}.
 
L'\algo de Bunch \& Hopcroft prend en entr\'ee une matrice
 de rang  \,$n$, \,$A\in\K^{n\times p}$
$(1\leq n\leq p)$ et donne en sortie un triplet
$\,(L,U,P)$\, tel que \,$L$\, est une matrice \utg
\infeez, \,$U$\, une
\matgu \freg
et \,$P$\, une matrice de permutation.
On \'ecrira: ${\bf lup}(A,n,p)={\bf lup}_{n,p}(A)=(L,U,P)$.
 
\ss Pour \,$n=1$, \,$A$\, est une matrice ligne de rang 1:
il existe donc un \elt non nul de A occupant la
\,$i\,$-\,\`eme place de cette ligne ($1\leq i\leq p$).
Il suffit de prendre \,$L=[1]$\, et
\,$U=AP$\, o\`u \,$P$\, est la matrice
de permutation d'ordre \,$p$\, correspondant \`a
l'\'echange des colonnes 1 et \,$i$. On a donc
\,${\bf lup}(A,1,p)=([1],\,A\,P,\,P)$\,
pour la matrice \,$P$\, ainsi d\'efinie.
 
\ss Supposant la propri\'et\'e vraie pour tout entier
\,$n$\, compris entre 1 et \,$2^{\nu-1}$,
on la d\'emontre pour
\,$2^{\nu-1}< n \leq 2^\nu$\,  ($\nu = \esup{ \log{n} }$).
On pose \,$n_0=2^{\nu-1}$, \,$n_1=n-n_0$\,
et \,$p_1=p-n_0$.
Pour obtenir \,${\bf lup}(A,n,p)$\, avec
\,$A\in \K^{n\times p}$, on consid\`ere
la partition suivante de la matrice \,$A$:
\begin{equation} \label{EqBU}
A=\left[\begin{array}{c} A_1 \\ A_2
\end{array} \right] ,\quad A_1\in \K^{n_0\times p},
\, A_2\in \K^{n_1\times p}
\end{equation}
Si \,$A$\, est une matrice surjective, \,$A_{1}$\, et
\,$A_{2}$\, le sont \egmtz.
On commence par appeler
\,${\bf lup}(A_1,n_0,p)$\, qui donne une \LUP-\deco
\,$(L_1,\,U_1,\,P_1)$\, de \,$A_1$.
On consid\`ere alors les partitions suivantes des matrices
\,$U_{1}$\, et \,$A_{2}\,P_{1}^{-1}$:
\begin{equation} \label{EqBU2}
\left\{
\begin{array}{l}
\,U_{1}=\left[\, V_1 \,|\, B  \,\right]\in \K^{n_0\times p}\quad
\mathrm{et} \quad
A_{2}\,P_{1}^{-1} = \left[\, C \,|\, D \,\right]\in
\K^{{n_1}\times p}\\[2mm]
\,V_1 \in \K^{n_0\times n_0},\,
B \in \K^{n_0 \times p_1},\, C\in \K^{n_1\times n_0},\,
D \in \K^{n_1 \times p_1}
\end{array}
\right.
\end{equation}
$V_1$\, \'etant \tgu et
inversible (puisque \,$U_{1}$\, est \fregz). Posant
\,$C_1=C\,V_{1}^{-1}$\, et \,$E=D-C_{1}\,B,$
on v\'erifie que:
\begin{equation} \label{EqBU3}
A=\left[\begin{array}{cc}
L_1 & 0 \\ C_{1} &
\mathrm{I}_{n_1} \end{array} \right]\,
\left[\begin{array}{cc} V_1 & B \\
0 & E \end{array} \right]\,P_1\,.
\end{equation}
Comme la matrice \,$E$\, satisfait \`a l'\hdr
(elle est surjective puisque \,$A$\, l'est),
on peut appliquer la \pcd \,${\bf lup}(E,n_1,p_1)$\,
qui donne la \LUP-\deco \,$E=L_{2}\,U_{2}\,P_{2}$\,
dans laquelle \,$U_{2}$\, est une matrice
\,$n_1 \times p_1$\, \tgu\fregz. \\
Il suffit de poser
 \,$Q=\left[\begin{array}{cc}
\mathrm{I}_{n_0} & 0 \\
0 & P_{2} \end{array} \right]$\, et \,$B_2:=B\,P_{2}^{-1}$\,
pour obtenir la
\decoz:
\begin{equation} \label{EqBU4}
A=\left[\begin{array}{cc}
L_1 & 0 \\ C_{1} &
L_2 \end{array} \right]\,
\left[\begin{array}{cc} V_1 &
B_{2} \\ 0 & U_{2}
\end{array} \right]\,Q\,P_1\,.
\end{equation}
Ce qui donne \,${\bf lup}(A,n,p)=(L,U,P)$\, avec:
$$
L=\left[\begin{array}{cc}
L_1 & 0 \\ C_{1} &
L_2 \end{array} \right]~,~~U=
\left[\begin{array}{cc} V_1 &
B_2 \\ 0 & U_{2}\end{array}
\right]~,~~\mathrm{et~~}P=Q\,P_1\,.
$$
En r\'esum\'e on obtient le sch\'ema r\'ecursif de l'\algo
 \vref{prevBU}.

\begin{algor}[: ${\bf lup}_{(n,p)}$, \LUP-\deco \`a la
    Bunch \& Hopcroft pour une
    matrice surjective.]\label{prevBU}
    \acl{prevBU}{\LUP-\deco \`a la Bunch \& Hopcroft}
\Entree Une matrice surjective \,$A\in\K^{n\times p}$\,
($\K$ est un corps).
\Sortie Les matrices \,$L,\,U,\,P$\, de la \LUP-\deco de \,$A$\,.
\Debut  On utilise la partition donn\'ee en (\ref{EqBU})
\Etap{1}{\recu avec  $A_{1}\in\K^{n_0 \times p}$,
\,$\nu=\esup{\log{n}}$,
  \,$n_0=2^{\nu-1}$.}
 {$(L_{1},U_{1},P_{1}):={\bf lup}_{(n_0,p)}(A_{1})$ }
\Etap{2}{pas d'\opari ici}
{$B_{2}:= A_{2}\,P_{1}^{-1}\quad $  avec $A_2\in\K^{n_1 \times p}$,
$n_1=n-n_0$.}
\Etap{3}{inversion d'une \matgu \reg}{$V_{2}:=V_{1}^{-1}\quad $ avec
$V_1 \in \K^{n_0\times n_0}$\,\,\, (cf. la partition  (\ref{EqBU2}))}
\Etas{4,\, 5,\, 6}{}{C_1:=C\,V_2\,;\,\,F:=C_1\,B\,;\,\,
E:=D-F\,.}
\Etap{7}{\recu avec \,$E\in\K^{n_1\times p_1}$\,}
   {$(L_{2},U_{2},P_{2}):={\bf lup}_{(n_1,p_1)} (E)$}
\Eta{8}{pas d'\opari ici}
{B_2:=B\,P_{2}^{-1},\,\,P:=Q\,P_1.}
\fin
\end{algor}
 
L'\algo obtenu est un \algo avec branchements. Ceci est
in\'evitable puisque la sortie \,$P$\, d\'epend de mani\`ere
discontinue de l'entr\'ee \,$A$. Les branchements sont tous
command\'es par le test d'\egt \`a 0 dans le corps \,$\K$. Notons
\,$\tau(n,p)$\, le nombre d'\oparis ex\'ecut\'ees par cet \algo
pour les matrices \,$A\in\K^{n\times p}$\, et \,$\pi(n,p)$\, son
\emph{temps \paral \aritz}, \cad sa \prof si on ne prend pas en
compte les \'etapes de recherche d'\'el\'ements non nuls ni les
produits d'une matrice par une matrice de permutation.
 
On a alors en suivant le sch\'ema r\'ecursif
\vref{prevBU} les \ines suivantes.
 
Tout d'abord concernant le nombre d'\oparisz:
\begin{equation} \label{EqBU5}
\tau(n,p)\leq\left\{
\begin{array}{ccccl}
 \tau(n_0,p)&+& 4\,\Ca\,{n_0}^\alpha & + &
~\Ca\,{n_0}^\alpha~~~ +\\ \esup{p_1/n_0}\,\Ca\,{n_0}^\alpha & + &
p_1\,n_1 &+& ~\tau(n_1,p_1)
\end{array}
\right.
\end{equation}
 
Le terme \,$p_1\,n_1$\, correspond \`a la soustraction
\,$E:=D-F$\, et le terme
\,$\Ca{n_0}^\alpha\esup{p_1/n_0}$\, correspond au
calcul du produit \,$C_1\,B$\, dans lequel \,$B \in \K^{n_0 \times
p_1}$\, et
\,$C_1\in \K^{n_1\times n_0}$: on peut toujours compl\'eter \,$C_1$\,
par des lignes de 0 pour en faire une matrice carr\'ee et on d\'ecoupe
\,$B$\, en \,$\esup{p_1/n_0}$\, blocs carr\'es (apr\`es lui avoir
\'eventuellement rajout\'e des colonnes de 0); on effectue alors
en \paral \,$\esup{p_1/n_0}$\, \muls dans \,$\K^{n_0 \times n_0}$.
 
\ss Ensuite, concernant le temps \paral \aritz, on obtient de la
m\^{e}me mani\`ere en utilisant le r\'esultat de l'inversion des
matrices triangulaires (section \ref{invtrian})~:
\begin{equation} \label{EqBU6}
\pi(n,p)\leq \pi(n_0,p)+\Ka\,[(\nu-1)^2+
5\,(\nu-1)+2]+2+\pi(n_1,p_1)
\end{equation}
 
On en d\'eduit pr\'ecis\'ement~:
 
\begin{theorem} \label{thBU}
La \LUP-\deco d'une matrice surjective de type \,$(n,p)$\, sur
\,$\K$\, peut \^etre effectu\'ee par un \algo (avec branchements)
qui ex\'ecute un nombre d'\oparis \'egal \`a \,$\tau(n,p)\,$
 en temps \paral (\aritz) \,$\pi(n,p)$\, major\'es par
$$ \tau(n,p)\leq {1\over 2}\,\gamma_{\alpha}\,(\esup{p\over
n}+1)\, n^{\alpha}+\frac{1}{2}\,\esup{p\over n}\,{n}^{2}\log{n}
~~\mathrm{et}~~ \pi(n,p)\leq 4\,(5\,\Ka+1)\,n $$
$$\mathrm{o\grave{u}}~~\gamma_{\alpha}=
\Ca\,\max{\,\left(4\,,\frac{1}{2^{\alpha-2}-1}\right)}\,\cdot$$
\indexnota{Caz@$\gamma_{\alpha}$}
\end{theorem}
 
\noindent Notons que pour \,$p=n$\,,  la taille du
circuit correspondant \`a l'\algo de Bunch \& Hopcroft est
exactement major\'ee par \,$\gamma_{\alpha}\, n^{\alpha}+
\frac{1}{2}\,{n}^{2}\log{n}\,.$
 
 
\ms\prv \\
Le calcul
de \,${\bf lup}(A,n,p)$\, se fait de mani\`ere r\'ecursive.
Nous donnons les majorations pour le cas
o\`u \,$n=2^\nu$, et il est clair que
si \,$n<2^\nu$, le calcul ne peut \^{e}tre que plus rapide.
 
Pour le temps \paral \arit on a \,$\pi(1,p)=0$\, donc,
vue la \recu (\ref{EqBU6}), le r\'esultat
ne d\'epend pas de \,$p$\, et
$$ \pi(n,p)=\pi(n)=\pi(2^{\nu})\leq 2\,\pi(2^{\nu-1})
+\Ka\,[(\nu-1)^2+ 5\,(\nu-1)+2]+2\,. $$
La relation de r\'ecurrence
\,$f(\nu)=2f(\nu-1)+c\,[(\nu-1)^2+5\,(\nu-1)]+2\,(c+1)$\, avec
\,$f(0)=0$\, est r\'esolue par {\sc Maple} en $$f(\nu)=\left(
10\,c+2 \right) \left( {2}^{\nu}-1 \right) -c{\nu}^{2}-7\,c\,\nu$$
major\'e par \,$\left( 10\,c+2 \right) {2}^{\nu}= 2\left( 5\,c+1
\right)n \,$, ce qui donne le r\'esultat.
 
\ms Pour calculer le nombre d'\oparis on pose
\,$r_\nu=r=\esup{p/n}$\, et on suppose \spdg que \,$p=r n$. L'\ine
(\ref{EqBU5}) se r\'e\'ecrit, puisque \,$p_1=(2r-1)n_0$:
$$
\tau(2^{\nu},p)\;\leq\;\left\{
\begin{array}{ll}
\tau(2^{\nu-1},p)\,+\,\left(4\Ca+
\Ca+(2r_\nu-1)\,\Ca\right)\,2^{(\nu-1)\alpha}\\[1mm]
+(2\,r-1)\,2^{2(\nu-1)}+ \,\tau(2^{\nu-1},p_1)
\end{array}
\right.$$ ce qui donne~:
$$ \tau(2^{\nu},p)\;\leq\; 2\,\tau(2^{\nu-1},p)\,+
(4\Ca+2\,r\,\Ca)\,2^{(\nu-1)\alpha}+\,2^{2\nu- 1}\,r\,. $$
Dans le d\'eroulement r\'ecursif de l'\algoz, lorsqu'on traite les
matrices de type\footnote{~La majoration vaut aussi pour les
matrices de type \,$2^\kappa\times p'$\, avec \,$p'\leq p$.}~
$2^\kappa\times p$\,, on a \,$r_\kappa=2^{\nu-\kappa}r_\nu\,$. Et
donc en ramenant \`a \,$r=r_\nu$\, on obtient les \inesz~:
$$ \tau(2^{\kappa},p)\;\leq\; 2\,\tau(2^{\kappa-1},p)\,+
\,2^{(\kappa-1)\alpha}\,(4\Ca\,+\,2^{\nu-\kappa+1}\,r\,\Ca)\,+
\,2^{2\kappa-1}\,2^{\nu-\kappa}\,r\,.$$
 
Sachant que \,$\tau(1,p)=0$, on obtient par sommation
(et simplifica\-tion de la solution d'une relation de
r\'ecurrence) la majoration suivante~:
$$\tau(n,p)\leq \frac{1}{2}\,{n}^{2}r\log{n}
+2\,{\frac {\Ca\,{n}^{a}r}{{2}^{a}-4}}
+4\,{\frac {\Ca\,{n}^{a}}{{2}^{a}-2}}\cdot$$
Ce qui donne le r\'esultat annonc\'e. \qed
 
\section[Calcul du d\'eterminant et de l'inverse]{Calcul du
d\'eterminant et de l'inverse \\ d'une matrice carr\'ee}
\label{sec detinv}

La \LUP-\deco pr\'ec\'edente permet un calcul \sql rapide
du \deter et de l'inverse d'une \maca
(inversible) en ramenant ces \pbs \`a la
\mul rapide des \macas d'ordre \,$n$.
 
\ss En effet, si l'on passe par la \LUP-\decoz, le calcul du
\deter d'une matrice \,$A\in \K^{n\times n}$\, s'effectue avec le
m\^eme ordre de \com \sqle que la \mul des matrices \,${n\times
n}$\, puisque si \,$A=L\,UP$, alors \,$\det\,{A}=\epsilon\,
\det{U}$\, ce qui revient \`a calculer le produit des \,$n$\,
\elts diagonaux de la \matg \,$U$\, ($\epsilon =\pm1$ est la
signature de la permutation repr\'esent\'ee par la matrice \,$P$).
Il y a donc, apr\`es la \LUP-\decoz, un calcul \supt en
\,$\SD(n,\log{n})$\, (par un circuit binaire \'equilibr\'e).
 
 
Il en est de m\^eme pour
le calcul de l'inverse de \,$A$, quand elle est inversible,
puisque \,$A^{-1}=P^{-1}U^{-1}L^{-1}$\, ce qui revient, en plus
de la \LUP-\decoz, \`a inverser deux matrices triangulaires
($\,U$\, et \,$L\,$) et \`a effectuer un produit de matrices
\,$n\times n$.
 
\ss A priori les \algos de calcul du \deter et de l'inverse
tels que nous venons de les d\'ecrire sont
des \algos avec branchements.
 
Dans cette perspective, le co\^ut de la recherche des \elts
non nuls comme celui des permutations de lignes
ou de colonnes, \cad des \muls \`a gauche ou \`a droite
par une matrice de permutation, n'est pas pris en consid\'eration
dans les comptes  d'\oparis aussi bien du
point de vue de leur nombre total
que de celui de leur profondeur.
 
N\'eanmoins, on peut aussi prendre le point de vue selon lequel
nous avons construit des \famus de \caris avec divisions, qui
calculent des fractions rationnelles formelles en les \coes
de la matrice donn\'ee au d\'epart. Il n'y a alors pas de
\LUP-\deco mais seulement une \LU-\decoz, sans
aucun branchement. Naturellement la contrepartie est
que l'\algo ne peut pas \^{e}tre ex\'ecut\'e concr\`etement
sur un corps avec une matrice arbitraire. C'est seulement
pour une \gui{matrice g\'en\'erique} que le \cari fonctionne: une telle
matrice est une matrice qui, lorsqu'on lui applique l'\algo avec
branchements,
subit tous les tests \,$x=0~?$\, en donnant une r\'eponse n\'egative.
 
\noindent Dans nos \'enonc\'es nous adoptons de pr\'ef\'erence ce
second point de vue.

\begin{prop}\label{propBU}
Le calcul du d\'eterminant d'une \maca d'or\-dre \,$n$\, sur un
corps \,$\K$\, est r\'ealis\'e par une \famu de \cirs \ariths avec
divisions en \,$\SD(n^\alpha,n)$\,. \\ Les constantes
asymptotiques sont respectivement major\'ees par
$\,\gamma_\alpha\,$ pour la taille $\,\O(n^\alpha)\,$ et par
$\,4\,(5\,\Ka+1)\,$ pour la \prof $\,\O(n)\,$ (les m\^emes
majorations des constantes que celles don\-n\'ees au th\'eor\`eme
\ref{thBU}).
\end{prop}
 
\begin{prop}\label{invcar}
L'inversion d'une \maca d'ordre \,$n$\, sur un corps \,$\K$\, est
un \pb r\'esolu par une \famu de \cirs \ariths avec divisions en
\,$\SD(n^\alpha,n)$\, avec la m\^eme estimation que celle de la
proposition \ref{propBU} pour la constante asymptotique de la
\prof $\,\O(n)\,$, et une constante asymptotique major\'ee par
\,$\zeta_\alpha =\gamma_\alpha+9\,\Ca\,$ pour la taille
$\,\O(n^\alpha)\,$.
\indexnota{za@$\zeta_\alpha$}
\end{prop}
 
\noindent Dans la constante \,$\zeta_\alpha\,$ de la proposition
ci-dessus, le terme $\,\gamma_\alpha\,$ correspond \`a la
\LUP-\deco et le terme $\,9\,\Ca\,$ \`a l'inversion de deux
matrices triangulaires suivie de la \mul de deux \macasz.
 
\section{Forme r\'eduite \'echelonn\'ee en lignes}
\label{sec FEL}

Dans cette section nous donnons un aper\c{c}u sur
une \met r\'ecursive permettant de r\'eduire les
matrices \`a \coes dans un corps commutatif
\,$\K$, \`a la forme \'echelonn\'ee en
lignes avec une \com \sqle du m\^eme
ordre que celle de la \mul des matrices.
 
\'Etant donn\'ee une matrice \,$A$\, de type \,$(n,p)\,$
sur \,$\K$, la r\'eduction de \,$A$\, \`a la forme
\'echelonn\'ee en lignes consiste \`a transformer \,$A$,
en ayant exclusivement recours \`a des transformations
\elrs \emph{unimodulaires} sur les
lignes\footnote{~Rappelons (cf. page \pageref{unimod})
qu'il s'agit d'une part de la transformation
qui consiste \`a ajouter \`a une ligne une combinaison
\lin des autres et d'autre part
des \emph{\'echanges sign\'es de lignes}
du type \,$(L_{i},L_{j})\leftarrow (L_{j},- L_{i})$.},
en une matrice de de type \,$(n,p)$\, sur
\,$\K$\, avec un nombre de z\'eros strictement croissant
apparaissant \`a gauche des lignes successives de la matrice
r\'eduite. Si l'on note \,$E$\, la matrice unimodulaire
correspondant \`a ces transformations\footnote{~C'est-\`a-dire
la matrice obtenue en faisant subir \`a la matrice unit\'e
d'ordre \,$n$\, les m\^emes transformations.}, cela revient \`a
multiplier la matrice \,$A$\, \`a gauche par la matrice \,$E$.
 
\ss
Prenons par exemple la \maca d'ordre 6
{\footnotesize $$A=
\left [\begin {array}{cccccc}
1&2&3&4&5&6\\
1&2&3&6&7&8\\
1&2&3&2&3&1\\
1&1&2&1&1&1\\
1&3&4&7&9&6\\
3&6&9&10&11&20
\end {array}\right]\,.
$$}
 
On peut la r\'eduire \`a la forme \'echelonn\'ee en lignes
en effectuant des transformations du style pivot de Gauss sur
les lignes. Ces m\^emes transformations, effectu\'ees sur les
lignes de la matrice unit\'e d'ordre 6, donnent la matrice
unimodulaire \,$E$\, qui r\'esume ces transformations:
{\footnotesize
$$
E=\left [\begin {array}{cccccc} 1&0&0&0&0&0\\
-1&0&0&1&0&0\\
-1&0&1&0&0&0\\
-2&0&-1&0&0&1\\
-2&0&0&1&1&0\\
-4/5&1&1&-3/5&-3/5&0\end {array}\right]\,.
$$}
 
La matrice r\'eduite \'echelonn\'ee en lignes est alors
donn\'ee par le produit:
{\footnotesize $$
EA=\left [\begin {array}{cccccc} 1&2&3&4&5&6\\
0&-1&-1&-3&-4&-5\\
0&0&0&-2&-2&-5\\
0&0&0&0&-2&7\\
0&0&0&0&0&-5\\
0&0&0&0&0&0
\end {array}\right]\,.
$$}
 
\sni
Comme nous l'avons fait pour la \LUP-\decoz, il
s'agit ici
de d\'ecrire une version rapide de la \mpg
sur les lignes. Mais contrairement \`a la \LUP-\decoz,
aucune hypoth\`ese \supt n'est faite sur la matrice
\,$A$\, et aucune permutation de colonnes n'est permise.
En contrepartie, dans la \deco \,$A=FU$\, qui
r\'esulte
de cette \met de r\'eduction ($\,F=E^{-1}\,$), la matrice
\,$F$\, poss\`ede seulement la propri\'et\'e d'\^etre
unimodulaire.
 
La forme \'echelonn\'ee en lignes trouve sa justification et son
application dans des \pbs comme la r\'esolution des \syses \lins
ou la d\'etermination d'une base pour un sous-espace de
\,$\K^{n}$\, d\'efini par un \sys g\'en\'erateur. Elle sera aussi
utilis\'ee dans la section \ref{sec kellerseq} pour le calcul
rapide du \polcar sur un corps. La \met que nous allons exposer
ci-dessous est due \`a Keller-Gehrig \cite{KeGe} et elle est
reprise dans \cite{Bur}.
 
\subsubsection*{Description de la \pcd rapide}
 
On consid\`ere une matrice \,$A\in \K^{n\times p}$.
Pour la r\'eduire \`a la forme \'echelonn\'ee en lignes, on
peut supposer \spdg que
\,$n=p=2^\nu$\,
quitte \`a compl\'eter la matrice \,$A$\, avec suffisamment de
lignes et/ou colonnes de z\'eros\footnote{\,Les lignes et les
colonnes ajout\'ees ne feront l'objet d'aucune manipulation.}.
 
\sni
La \pcd principale que nous noterons \,{\bf Fel}\,
utilise les \pcds auxiliaires \,${\bf Fel}_{1}$,
\,${\bf Fel}_{2}$\, et \,${\bf Fel}_{3}$\, suivantes.
 
\mni {\bf Proc\'edure} \,${\bf Fel}_{1}$:
 
\noi C'est une \pcd r\'ecursive qui transforme
une \maca \,$A\in \K^{2n\times n}\,$
($\,n=2^\nu\,$) dont la moiti\'e \infee est \tgu
en une \matguz. \\
Plus pr\'ecis\'ement, si \,$A=\left[\begin{array}{cc} A_{1} \\
A_{2}\end{array}\right]$\, avec \,$A_{1},A_{2} \in
\K^{n\times n}$\, et \,$A_{2} $\, \tguz,
la \pcd \,${\bf Fel}_{1}$ calcule une matrice
unimodulaire
\,$E\in \mathrm{SL}_{2n}\,(\mathrm{K})$\, et une matrice
\,$T\in \K^{n\times n}$\, \tgu
telles que \,$E\,A=\left[\begin{array}{cc} T\\0\end{array}
\right]$.
Utilisant l'approche \gui{diviser pour gagner}  on divise la matrice
\,$A$\, donn\'ee en huit blocs
\,$2^{\nu-1}\times 2^{\nu-1}$\,
(si \,$\nu=0$, le traitement de la matrice \,$A$\, est
\immz)
et on applique de mani\`ere r\'ecursive la \pcd
\,${\bf Fel}_{1}$\, aux blocs \,$2^{\nu-1}\times 2^\nu\,$
qui
poss\`edent la m\^eme propri\'et\'e que \,$A$.
On obtient, avec des notations \'evidentes, le d\'eroulement
suivant de la \pcdz:
 
$$
\left[\begin{array}{cc}
A_{11} & A_{12}\\ A_{13} & A_{14} \\
A_{21} & A_{22}\\ 0 & A_{24} \end{array} \right]
~\stackrel{E_{1}}{\dans}~
\left[\begin{array}{cc}
A_{11} & A_{12}\\ A'_{13} & A'_{14} \\
0 & A'_{22} \\ 0 & A_{24} \end{array} \right]
\stackrel{E_{2}}{\dans}
\left[\begin{array}{cc}
A'_{11} & A'_{12} \\ 0 & A''_{14} \\
0 & A'_{22}\\  0 & A_{24} \end{array} \right]
~\stackrel{E_{3}}{\dans}~
$$
$$
\left[\begin{array}{cc}
A'_{11} & A'_{12}\\ 0 & A''_{14} \\
0 & A''_{22}\\  0 & 0 \end{array} \right]
~\stackrel{E_{4}}{\dans}~
\left[\begin{array}{cc}
A'_{11} & A'_{12}\\ 0 & A'''_{14} \\
0 & 0 \\  0 & 0 \end{array} \right] ~~~~~~\mathrm{avec~:}
$$
$$
E_{1}\left[\begin{array}{cc} A_{13} \\
A_{21}\end{array}\right] =
\left[\begin{array}{cc} A'_{13} \\
0 \end{array}\right]~~\mathrm{et}
~~\left[\begin{array}{cc} A'_{14} \\
A'_{22}\end{array}\right] =
E_{1}\left[\begin{array}{cc} A_{14} \\
A_{22} \end{array}\right]~;
$$
$$
E_{2}\left[\begin{array}{cc} A_{11} \\
A'_{13}\end{array}\right] =
\left[\begin{array}{cc} A'_{11} \\
0 \end{array}\right]~~\mathrm{et}
~~\left[\begin{array}{cc} A'_{12} \\
A''_{14}\end{array}\right] =
E_{2}\left[\begin{array}{cc} A_{12} \\
A'_{14} \end{array}\right]~;
$$
$$
E_{3}\left[\begin{array}{cc} A'_{22} \\
A_{24}\end{array}\right] =
\left[\begin{array}{cc} A''_{22} \\
0 \end{array}\right]~~\mathrm{et}
~~~E_{4}\left[\begin{array}{cc} A''_{14} \\
A''_{22} \end{array}\right]=
\left[\begin{array}{cc} A'''_{14} \\
0\end{array}\right]\,.
$$
 
$$
\mathrm{Posant}~~~E=\left[\begin{array}{ccc}
\I_{2^{\nu-1}} & 0 & 0 \\ 0 & E_{4} & 0 \\
0 & 0 & \I_{2^{\nu-1}}
\end{array}\right]
\left[\begin{array}{cc} E_{2} & 0 \\
0 & E_{3} \end{array}\right]
\left[\begin{array}{ccc}
\I_{2^{\nu-1}} & 0 & 0 \\ 0 & E_{1} & 0 \\
0 & 0 & \I_{2^{\nu-1}}
\end{array}\right]\,,
$$
on a bien \,$EA=\left[\begin{array}{cc} T \\
0 \end{array}\right]$\, o\`u \,$T=
\left[\begin{array}{cc} A'_{11} &  A'_{12} \\
0 &  A'''_{14} \end{array}\right]$\, est une
\matgu \,$A'_{11}$\, et
\,$A'''_{14}$\, le sont.
 
\mni {\bf Proc\'edure} \,${\bf Fel}_{2}$:
 
\noi Elle prend en entr\'ee une \maca
\,$A\in \K^{n\times n}$\, ($\,n=2^\nu\,$) et retourne
une matrice unimodulaire \,$E\in \mathrm{SL}_{n}\,(\mathrm{K})\,$
et une \matgu \,$T$\, v\'erifiant
\,$EA=T$.\\ L\`a encore, on obtient avec l'approche
\gui{diviser pour gagner} et des notations analogues
\`a celles utilis\'ees pr\'ec\'edemment, le d\'eroulement
suivant de la \pcdz:
$$
~~~~~~~\left[\begin{array}{cc}
A_{11} & A_{12}\\ A_{21} & A_{22} \end{array} \right]
~\stackrel{E_{1}}{\dans}~
\left[\begin{array}{cc}
A_{11} & A_{12}\\ A'_{21} & A'_{22} \end{array} \right]
$$
$$
\stackrel{E_{2}}{\dans}~
\left[\begin{array}{cc}
A'_{11} & A'_{12} \\ 0 & A''_{22} \end{array} \right]
~\stackrel{E_{3}}{\dans}~
\left[\begin{array}{cc} A'_{11} & A'_{12} \\ 0 & A'''_{22}
\end{array} \right]
$$
o\`u \,$E_{1}$\, est la matrice unimodulaire correspondant
\`a l'\algo \,${\bf Fel}_{2}$\, appliqu\'e de mani\`ere
r\'ecursive \`a la matrice \,$A_{21}\,$
($\,E_{1}A_{21}=A'_{21}$\, est donc une \matguz,
et l'on pose \,$A'_{22}=E_{1}A_{22}\,$) alors
que les matrices \,$E_{2}$\, et \,$E_{3}$\, corres\-pondent
\`a l'application respective de l'\algo
\,${\bf Fel}_{1}$\, \`a la matrice
\,$\left[\begin{array}{cc} A_{11} \\ A'_{21}
\end{array}\right]\,$
qui est de type \,$(2^\nu,2^{\nu-1})$\, et de l'\algo
\,${\bf Fel}_{2}$\, \`a la matrice \,$A''_{22}$\, qui est
carr\'ee
d'ordre \,$2^{\nu-1}$. Cela se traduit par le fait que
\,$E_{3}A''_{22}=A'''_{22}$\, est \tgu
et que:
$$
E_{2}\left[\begin{array}{cc} A_{11} \\
A'_{21} \end{array}\right] =
\left[\begin{array}{cc} A'_{11} \\
0 \end{array}\right]~~\mathrm{avec}
~~\left[\begin{array}{cc} A'_{12} \\
A''_{22}\end{array}\right] =
E_{2}\left[\begin{array}{cc} A_{12} \\
A'_{22} \end{array}\right]\,.
$$
Posant $~E=\left[\begin{array}{cc}
\I_{2^{\nu-1}} & 0 \\ 0 & E_{3} \end{array}\right]
 E_{2} \left[\begin{array}{cc}
\I_{2^{\nu-1}} & 0 \\ 0 & E_{1} \end{array}\right]\,,$
on a bien \,$EA=\left[\begin{array}{cc} T \\
0 \end{array}\right]$\, o\`u \,$T=
\left[\begin{array}{cc} A'_{11} &  A'_{12} \\
0 &  A'''_{22} \end{array}\right]$\, est une \matguz.
 
\mni {\bf Proc\'edure} \,${\bf Fel}_{3}$:
 
\noi Elle prend en entr\'ee une \matgu
\,$A\in \K^{n\times n}$\, (avec \,$n=2^\nu\,$) et donne
en
sortie une matrice unimodulaire \,$E\in \mathrm{SL}_{n}({\rm
K})\,$
et une matrice \,$S$\, sous forme \'echelonn\'ee en lignes
v\'erifiant \,$EA=S$. \\
On consid\`ere la partition \,$A=\left[\!\begin{array}{cc}
A_{11} & A_{12}\\ 0 & A_{22} \end{array} \right]$\, en
blocs
\,$2^\nu\times 2^\nu$\, de la matrice \,$A$\, donn\'ee
($A_{11}$\, et \,$A_{22}$\, sont des \matgusz).
 
\sni
Le d\'eroulement de la \pcd est alors
illustr\'e par le sch\'ema sui\-vant dans lequel c'est d'abord
l'\algo \,${\bf Fel}_{3}$\, qui est appliqu\'e \`a la matrice
\,$A_{11}$\, pour donner la matrice
\,$\left[\!\begin{array}{c}
S_{11}\\ 0 \end{array} \right]$\, o\`u \,$S_{11}$\, est une
matrice surjective\footnote{\,Le nombre \,$r$\, de ses lignes
est \'egal \`a son rang qui est aussi celui de \,$A_{11}$.}
\'echelonn\'ee en lignes; c'est ensuite \,${\bf Fel}_{1}\,$
qui est appliqu\'e \`a la matrice \,$\left[\begin{array}{c}
A_{23}\\ A_{24}\end{array} \right]$\, pour donner la matrice
\,$\left[\begin{array}{c} A'_{23} \\ 0 \end{array} \right]\,$
o\`u \,$A'_{23}$\, est \tguz;
et c'est enfin \,${\bf Fel}_{3}$\, qui, appliqu\'e \`a la
matrice $\, A'_{23}$, donne la matrice
\'echelonn\'ee en lignes \,$S_{23}$:
$$
\left[\!\begin{array}{cc}
A_{11} & A_{12}\\ 0 & A_{22} \end{array} \right]
\stackrel{E_{1}}{\dans}
\left[\!\begin{array}{cc} S_{11} & A'_{12}\\ 0 & A_{23} \\
0 & A_{24} \end{array} \right] \stackrel{E_{2}}{\dans}
\left[\!\begin{array}{cc} S_{11} & A'_{12}\\ 0 & A'_{23}
\\ 0 & 0 \end{array} \right]
\stackrel{E_{3}}{\dans}
\left[\!\begin{array}{cc} S_{11} & A'_{12} \\
0 & S_{23} \\ 0 & 0 \end{array} \right]
$$
avec \,$E_{1}A_{12}=\left[\!\begin{array}{c}
A'_{12}\\ A_{23} \end{array} \right]$.
Si maintenant on pose:
$$
E=\left[\! \begin{array}{ccc}
\I_{r} & 0 & 0 \\ 0 & E_{3} & 0 \\
0 & 0 & \I_{r} \end{array}\right]
\left[\! \begin{array}{cc}\I_{r} & 0 \\
0 & E_{2} \end{array}\right]
\left[\! \begin{array}{cc}
E_{1} & 0  \\ 0 & \I_{2^{\nu-1}}
\end{array}\right]\,
$$
o\`u \,$r$\, est le rang de \,$A_{11}$, alors \,$EA=S\,$
avec \,$S=\left[\! \begin{array}{cc}
S_{11} & A'_{12} \\ 0 & S_{23} \\
0 & 0 \end{array} \right]$\, qui est bien une matrice
\'echelonn\'ee en lignes puisque \,$S_{11}$\, et
\,$S_{23}$\, le sont.
 
\mni {\bf Proc\'edure principale} \,${\bf Fel}$:
 
\noi
Elle prend en entr\'ee une \maca
\,$A\in \K^{n\times n}$\, ($\,n=2^\nu\,$) et retourne
une
matrice unimodulaire \,$E\in \mathrm{SL}_{n}(\mathrm{K})$\, et
une
matrice \,$S$\, sous forme \'echelonn\'ee en lignes v\'erifiant
\,$EA=S$.\\
Le cas \,$\nu=0$\, est trivial. Pour \,$\nu\geq 1$, on
applique la \pcd auxiliaire \,${\bf Fel}_{2}$\, pour
transformer la matrice \,$A$\, en une \matgu
\,$T$\, puis la \pcd \,${\bf Fel}_{3}\,$
pour transformer \,$T$\, en une matrice \'echelonn\'ee en
lignes.
 
\subsubsection*{Analyse de complexit\'e}
 
L'\'etude de complexi\'e de la \pcd principale \,${\bf Fel}\,$
passe par celle des trois \algos auxiliaires \,${\bf Fel}_{1}$,
\,${\bf Fel}_{2}$\, et \,${\bf Fel}_{3}$. Si l'on d\'esigne par
\,$\tau_{1}$, \,$\tau_{2}$\, et \,$\tau_{3}$\, les tailles et par
\,$\pi_{1}$, \,$\pi_{2}$\, et \,$\pi_{3}$\, les \profs respectives
de ces trois \algosv on a les majorations suivantes dans
lesquelles les \coes \,$\Ca$\, et \,$\Ka$\, sont les constantes
intervenant dans la taille et la \prof des complexit\'es
arithm\'etiques de la \mul des matrices.
 
\sni 
\s{pour les tailles:}
\begin{itemize}
\item [] $\tau_{1}(2^\nu)\leq 4\,\tau_{1}(2^{\nu-1}) +
2\,\Ca\,2^{\nu\alpha}$\,
\item [] $\tau_{2}(2^\nu)\leq \tau_{1}(2^{\nu-1}) +
2\,\tau_{2}(2^{\nu-1}) + (2^\alpha + 1)\,\Ca\,
2^{(\nu-1)\alpha}$\,
\item [] $\tau_{3}(2^\nu)\leq \tau_{1}(2^{\nu-1}) +
2\,\tau_{3}(2^{\nu-1}) + \Ca\,2^{(\nu-1)\alpha}\,$
\end{itemize}
 
\sni 
\s{pour les profondeurs:}
\begin{itemize}
\item []  $\pi_{1}(2^\nu)\leq 3\,\pi_{1}(2^{\nu-1}) +
2\,\Ka\,\nu$\,
\item [] $\pi_{2}(2^\nu)\leq \pi_{1}(2^{\nu-1}) +
2\,\pi_{2}(2^{\nu-1}) + \Ka\,(2\nu-1)$\,
\item [] $\pi_{3}(2^\nu)\leq \pi_{1}(2^{\nu-1}) +
2\,\pi_{3}(2^{\nu-1}) + \Ka\,(\nu-1)$\,
\end{itemize}
 
Il faut remarquer que dans la \pcd \,${\bf Fel}_{1}\,$ les
\'etapes \,$\stackrel{E_{2}}{\dans}$\, et
\,$\stackrel{E_{3}}{\dans}$\, peuvent \^etre ex\'ecut\'ees en
\paralz, ce qui explique la diminution du \coe (de 4 \`a 3) entre
\,$\tau_{1}$\, et \,$\pi_{1}$.
 
\mni Utilisant les \ines ci-dessus et le fait que:
$$ \tau(2^\nu)=\tau_{2}(2^\nu) + \tau_{3}(2^\nu)~~~~
\mathrm{et}~~~~\pi(2^\nu)=\pi_{2}(2^\nu) + \pi_{3}(2^\nu)\,, $$
nous allons montrer le r\'esultat suivant concernant la \com du
\pb de la r\'eduction \`a la forme \'echelonn\'ee en lignes.
 
\begin{prop}\label{fel}
La r\'eduction \`a la forme \'echelonn\'ee en lignes d'une \maca
d'ordre \,$n$\, sur un corps commutatif \,$\K$\, est r\'ealis\'ee
par une \famu de \caris de taille \,$\tau(n)$\, et de \prof
\,$\pi(n)\,$ avec les majorations suivantes:
$$ \tau(n)\leq
\frac{21\,\Ca}{2^{\alpha-2}-1}\,n^\alpha~~~\mathrm{et} ~~~\pi(n)
\leq (3\,K_{{\alpha}}+2)\,n^{\log{3}}\,. $$
\end{prop}
 
\prv Les sommations des relations \,$\tau_{1}(2^k)\leq
4\,\tau_{1}(2^{k-1}) + 2\,\Ca\,2^{k\alpha}$\, d'une part et des
relations \,$\pi_{1}(2^k)\leq 3\,\pi_{1}(2^{k-1}) + 2\,\Ka\,k$\,
d'autre part pour \,$k$\, allant de 0 \`a \,$\nu$\, (avec
\,$\tau_{1}(1)=\pi_{1}(1) =1\,$) donnent les majorations suivantes
pour la taille et la \prof du \cari correspondant \`a la \pcd
\,${\bf Fel}_{1}$:
$$ \tau_{1}(2^\nu)< \frac{2^{\alpha+1}\Ca}
{2^{\alpha-2}-1}\,2^{\nu\alpha}~~\mathrm{et}~~ \pi_{1}(2^\nu) <
\left({3\over 2}\,\Ka+1\right)\,{3}^\nu -\Ka\,\nu-{3\over
2}\,\Ka\,. $$
Tenant compte de ces relations et du fait que
\,$\tau_{2}(1)=\pi_{2}(1)=1$, les sommations pour \,$k$\, allant
de 0 \`a \,$\nu$\, des \ines relatives \`a la taille et la \prof
du \cari correspondant \`a la \pcd \,${\bf Fel}_{2}$\, nous
donnent la majoration:
$$ \tau_{2}(2^\nu) < {1\over 2}\,E_{\alpha}\,
\frac{2^{\nu\alpha}}{2^{\alpha-1}-1} < {1\over
2}\,E_{\alpha}\,2^{\nu\alpha}~~~ \mathrm{dans~laquelle} $$
$$ E_{\alpha}\,=\,\frac{2^{\alpha+1}}
{2^{\alpha-2}-1}\,\Ca+\left(2^{\alpha}+1\right)\,\Ca\,<\,
\frac{25\,\Ca}{2^{\alpha-2}-1} $$
(avec \,$2<\alpha \leq 3$) et la majoration:
$$ \pi_{2}(2^\nu)\, <\, \left({3\over 2}\,\Ka+1\right)\,
{3}^\nu-\Ka\,\left({2}^\nu+\nu +{1\over 2}\right)\,. $$
 
\sni On obtient, par des calculs analogues, les majorations
suivantes pour la taille et la \prof du \cari correspondant \`a la
\pcd \,${\bf Fel}_{3}$:
$$ \tau_{3}(2^\nu) < {1\over 2}\,F_{\alpha}\,
\frac{2^{\nu\alpha}}{2^{\alpha-1}-1} < {1\over
2}\,F_{\alpha}\,2^{\nu\alpha}~~~~\mathrm{et} $$
$$ \pi_{3}(2^\nu) < \left({3\over 2}\,\Ka+1\right)\, {3}^\nu -
{3\over 2}\,\Ka\,(2^{\nu+1}-1) $$
$$ \mathrm{o\grave{u}}~~~F_{\alpha}=\left(\frac{2^{\alpha+1}}
{2^{\alpha-2}-1}+ 1\right)\,\Ca<\frac{17\,\Ca}{2^{\alpha-2}-1} $$
 
\sni Le r\'esultat annonc\'e d\'ecoule des majorations ci-dessus
et du fait que l'on a: \,$\tau(2^\nu)=
\tau_{2}(2^\nu)+\tau_{3}(2^\nu)$, \,$\pi(2^\nu)=\pi_{2}(2^\nu)+
\pi_{3}(2^\nu)$\, et \,$n=2^\nu$. \qed
 
\mni \rem Le fait de consid\'erer des \macas dont le nombre de
lignes (et de colonnes) est une puissance de 2 n'est pas une
hypoth\`ese restrictive. On peut en effet plonger toute matrice
\,$A\in \K^{n\times p}\,$ dans une \maca d'ordre \,$2^\nu$\, en
prenant \,$\nu=\max{(\esup{ \log{n} }, \esup{ \log{p} })}\,$ et en
compl\'etant la matrice donn\'ee par \,$2^\nu-n$\, lignes et
\,$2^\nu-p$\, colonnes nulles. Les rang\'ees ajout\'ees, form\'ees
de z\'eros, ne subissent aucune transformation au cours du
d\'eroulement de la \pcd d\'ecrite et le r\'esultat \'enonc\'e
dans la proposition \ref{fel} reste valable \`a condition de
remplacer \,$n$\, par \,$\max{(n,p)}$.

\section{M\'ethode de Keller-Gehrig}
\label{sec kellerseq}
 
Les \algos de Keller-Gehrig \cite{KeGe} sont des versions
acc\'el\'er\'ees de l'\algo de Frobenius que nous avons  d\'ecrit
\`a la section \ref{sec.TriFro}.
 
Dans la section pr\'esente nous ne d\'ecrirons en d\'etail que le
plus simple de ces \algosz. Nous reprenons les  notations de la
section~\ref{sec.TriFro}.
 
La matrice \,$A\in\K^{n\times n}$ d\'efinit l'\endom \,$h_A$\, de
\,$\K^n$. Nous appelons \,$a=(e_{1},\ldots,e_{n})$\, la base
canonique de \,$\K^{n}$.
 
\subsubsection*{Acc\'el\'eration dans le cas simple}
Nous examinons ici le cas le plus simple (et le plus fr\'equent)
o\`u \,$k_{1}=n$\, \cad le cas o\`u
\,$b=(e_{1},Ae_{1},\ldots,A^{n-1}e_{1})$\, est une base
de~$\,\K^{n}$.
 
Nous d\'esignons par \,$[S']_{S}$\, la matrice d'un \sys de
vecteurs (ou d'un vecteur ou d'un \endoz) \,$S'$\, dans une base
\,$S$. Alors \,$U=[b]_{a}$\, est la matrice de passage de \,$a$\,
\`a \,$b$, et on a:
$$ [h_{A}]_{b}= U^{-1}AU=[(e_{1},Ae_{1},\ldots,A^{n-1}e_{1})]_{b}
= \left[\matrix{
  0       &  \ldots   &  0   &   a_{0} \cr
  1       &  \ddots   &  \vdots    &   a_{1}      \cr
  \vdots  &  \ddots   &  0    &    \vdots     \cr
  0       &  \ldots   &  1   &   a_{n-1}
}\right] $$
o\`u \,$a_0,a_1, \ldots , a_{n-1}$\, sont les \coes (dans
\,$\K\,$) de la relation de d\'epen\-dance
\,$A^ne_{1}=a_{n-1}A^{n-1}e_{1}+ \ldots + a_1Ae_{1} + a_0e_1$.
Ceci prouve que \,$A$\, est semblable \`a une matrice de Frobenius
et que son \polcar est
\,$\rPA(X)=(-1)^n\,(X^n-(a_{n-1}X^{n-1}+\cdots+a_{1}X+a_0))$.
 
\ss L'\algo de Keller-Gehrig, dans ce cas le plus simple, consiste
\`a calculer la matrice \,$U$\, puis le produit \,$U^{-1}A\,U$\,
pour obtenir par simple lecture de la derni\`ere colonne les \coes
du \polcar de \,$A$. Prenant \,$\nu=\esup{ \log{n} }$, le calcul
de \,$U$\, se fait en \,$\nu$\, \'etapes. L'\'etape \num $k$
($\,1\leq k\leq \nu\,$) consiste \`a:
 
\sni --- calculer la matrice \,$A^{2^{k}}$\, (\'el\'evation au
carr\'e de la matrice \,$A^{2^{k-1}}$\, d\'ej\`a calcul\'ee \`a
l'\'etape pr\'ec\'edente);
 
\sni --- calculer la matrice \,$A^{2^{\nu-1}}\,
[e_{1}\,|\,Ae_{1}\,|\,\ldots\,|\,A^{2^{k-1}-1}e_{1}]$\, \`a partir
de la matrice \,$[e_{1}\,|\,Ae_{1}\,|\,\ldots\,
|\,A^{2^{k-1}-1}e_{1}]$\, calcul\'ee \`a l'\'etape \,$k-1$\, pour
obtenir la matrice \,$[e_{1}\,|\,Ae_{1}\,|\,\ldots\,
|\,A^{2^{k}-1}e_{1}]$\, de l'\'etape \,$k$.
 
\ss \`A la fin de ces \,$\nu$\, \'etapes, on obtient la matrice
 
\ss
\centerline{$[e_{1}\,|\,Ae_{1}\,|\,\ldots\,|\,A^{2^\nu-1}e_{1}]
\in \K^{n\times 2^\nu}$}
 
\sni qui admet comme sous-matrice la matrice recherch\'ee
 
\ss \centerline{$U=[e_{1}\,|\,Ae_{1}\,|\,\ldots\,|\,A^{n-1}e_{1}]
\in \K^{n\times n}$}
 
\sni puisque \,$2^\nu-1\geq n-1$.
 
\ss On calcule ensuite la derni\`ere colonne de \,$U^{-1}\,A\,U$\,
en commen\c{c}ant par inverser la matrice \,$U$\, (en passant par
sa \LUP-\decoz). Enfin on calcule la derni\`ere colonne \,$V$\, de
\,$A\,U$\, en multipliant \,$A$\, par la derni\`ere colonne de
\,$U$, puis on calcule \,$U^{-1}\,V$.

\ms L'analyse de \com dans ce cas simple nous donne donc:
 
\begin{proposition}
\label{propKeGe1} On peut calculer le \polcar d'une \maca d'ordre
\,$n$\, \`a \coes dans un corps \,$\K$\, au moyen d'un \cari avec
divisions en \,$\SD(n^\alpha \log{n}, n\,\log{n})$, de taille
major\'ee  plus pr\'ecis\'ement par $$\,2\, \Ca\,n^\alpha\,\esup{
\log{n} }+ \zeta_{\alpha}\,n^\alpha +\O(n^2)\, $$ o\`u  \,$\Ca$\,
et \,$\zeta_{\alpha}$\, sont les constantes intervenant dans les
\coms \sqles de la \mul des matrices et de l'inversion des \macas
 (voir proposition \vref{invcar}).
\end{proposition}

\subsubsection*{Le cas g\'en\'eral}
L'\algo pr\'ec\'edent fournit d\'ej\`a une \famu de \caris avec
divisions qui calcule le \polcar d'une matrice sur un corps, au
sens des circuits avec divisions. Autrement dit, le \cari \'evalue
correctement le \polcar en tant que fraction rationnelle: en tant
qu'\elt du corps \,$\K((a_{ij}))$\, o\`u les \coes \,$a_{ij}$\, de
la \maca sont pris comme des \idtrsz.
 
Mais il \'echoue \`a calculer le \polcar de toute matrice qui n'a
pas un \polmin de m\^{e}me degr\'e que le \polcarz.
 
On est donc dans une situation pire que pour le calcul du \deter
\`a la Bunch \& Hopcroft, car dans ce dernier cas, il suffit de
multiplier \`a droite et \`a gauche la matrice par des matrices
unimodulaires (\`a petits \coes entiers) prises au hasard pour
obtenir une matrice qui poss\`ede une \LU-\deco avec une tr\`es
grande probabilit\'e\footnote{~Si cette \met
est seulement probabiliste en th\'eorie, elle fonctionne toujours
en pratique.}.
Et ceci m\^{e}me si son \deter est nul (cf. l'\algo \vref{apgs2}).
L'\algo de Bunch \& Hopcroft sans branchement, avec le
preprocessing que nous venons d'indiquer n'\'echouera que dans le
cas d'une matrice \,$n\times n$\, dont le rang est strictement
\infe \`a \,$n-1$.

\smallskip C'est donc en produisant un
\algo avec branchements qui fonctionne dans tous les cas que
Keller-Gehrig r\'ealise son v\'eritable tour de force. Et pour
cela il lui fallait d'abord d\'evelopper sa m\'ethode de
r\'eduction rapide d'une matrice \`a la forme \'echelonn\'ee en
lignes (sur un corps). Dans cette r\'eduction nous avons vu que la
profondeur de l'\algo (avec branchements) est  un
\,$\O(n^{\log\,3})$. Keller-Gehrig obtient pr\'ecis\'ement le
r\'esultat suivant:
 
\begin{theorem}
\label{thKG} Le \polcar d'une \maca d'ordre \,$n$\, sur un corps
\,$\K$\, peut \^etre calcul\'e par un \algo avec branchements qui
a pour taille un \,$\O(n^{\alpha}\log\,n).$
\end{theorem}

\subsubsection*{Une version plus rapide pour les cas favorables}
Notons que  Keller-Gehrig propose une version plus rapide
pour un \algo avec divisions mais sans branchements, qui calcule
le \polcar dans les m\^{e}mes conditions qu'\`a 
la proposition~\ref{propKeGe1}:

\begin{proposition}
\label{propKeGe2} On peut calculer le \polcar d'une \maca d'ordre
\,$n$\, \`a \coes dans un corps \,$\K$\, au moyen d'un \cari avec
divisions qui a pour taille un $\O(n^\alpha)$.
\end{proposition}

\subsubsection*{Une version parall\`ele}

Signalons enfin qu'une \paran de l'\algo de Keller-Gehrig
a \'et\'e obtenue par Giesbrecht \cite{Gie94,Gie95}.

\section{M\'ethode de Kaltofen-Wiedemann}
\label{sec kalto}

Pour \gnr l'\algo de Wiedemann (section \ref{subsec.wied}) \`a un
\acoma \,$\A$\, en \'evitant les divisions qu'il contient et le
d\'ebarrasser en m\^eme temps de son aspect al\'eatoire, l'id\'ee
de Kaltofen \cite{Kal} est de lui appliquer la \met de l'\elid de
Strassen (cf. le \tho \vref{thEliDiv}). Il doit pour cela exhiber
une matrice particuli\`ere \,$C\in \K^{n\times n}$\, et un couple
de vecteurs \,$u,v\in\A^{n\times 1}$\, pour lesquels l'\algo de
Wiede\-mann s'effectue sans divisions et tels que le \polgmin de
la \srl \,$(\tra{u}\,C^i\,v)_{\,i\in\NN}$, qui est donn\'e par
l'\algo de Berlekamp/Massey \cite{Dor}, est de degr\'e \,$n$\, (et
n'est autre, par cons\'e\-quent, que le \polmin \,$\rP^C$\, et,
\`a un signe pr\`es, le \polcar \,$\rP_C$\, de \,$C\,$).

\ss Kaltofen  consid\`ere la suite de nombres entiers \,$(a_i)\in
\NN^\NN$\, d\'efinie par $$ a_i={i \choose \einf{ i/2 }},
\;\;\mathrm{avec}\;\; \left\{\begin{array}{ll} a_{n+1}=2\,a_{n} &
\mathrm{si} \;n\; \mathrm{est \;impair\; et} \\[1mm]
a_{n+1}=2\,\frac{n+1}{n+2}\,a_{n}& \mathrm{si} \;n\; \mathrm{est
pair}.
\end{array}\right.
$$
 
Les premiers termes sont \,$1$, $ 1$, $ 2$, $ 3$, $6$, $10$, $20$,
$35$, $70$, $ 126$, $252$, $462$, $924$, $1716$, $3432$, $6435$,
$12870$, $24310$, $48620$, $92378$, \ldots
 
Il applique l'\algo de Berlekamp/Massey aux \,$2n$\, premiers
termes:
$$a_0=1,\,a_1=1,\,a_2=2,\,\ldots,\,a_{2n-1} = {2n-1 \choose n-1}
$$
Il constate que les restes successifs dans l'\algo d'Euclide
\'etendu, jusqu'au \,$(n-1)\,$-\,\`eme, ont un \coe dominant
\'egal \`a \,$\pm 1$, avec un degr\'e ne diminuant que d'une seule
unit\'e \`a chaque pas (\cad que \,$\dg R_{i} = 2n-1-i$\, pour
\,$1\leq i\leq n-1\,$). Ce qui garantit le fait que les \pols
\,$R_i,\,Q_i,\,U_i,\,V_i~(1\leq i\leq n)$\, appartiennent \`a
\,$\ZZ[X]$\, et que \,$\dg R_n = n-1$. Il constate \egmt que les
multiplicateurs \,$V_i~(1\leq i\leq n)$\, ont un \coe dominant et
un terme constant \'egaux \`a \,$\pm 1$\, et que, par
cons\'equent, dans la derni\`ere \egt obtenue:
$$\,U_{n}\,X^{2n} + V_n\somm_{i=0}^{2n-1}a_iX^i=R_n ~~~~ (\mbox{
avec }~\dg R_n = n-1)\,, $$
\,$V_n$\, est un \poly de degr\'e \,$n$\, qui, \`a un signe
pr\`es, s'\'ecrit:
$$\pm V_n=X^n - (c_{n-1}X^{n-1} + \cdots + c_1X + c_0)\,. $$
(avec $c_0 = \pm 1,\,c_{n-1}=1$) Kaltofen  montre m\^eme, \`a
partir de l'\algo qui calcule les \coes de \,$V_n$, que ces
derniers sont en fait donn\'es par la formule:
$$c_i = (-1)^{\einf{ \frac{n-i-1}{2}}} {\einf{\frac{n+i}{2}
}\choose i} ~~~~~\mbox{ pour } ~~~0\leq i\leq n-1\,. $$
C'est donc le \poly \,$f(X) = X^n - c_{n-1}X^{n-1} - \ldots - c_1X
- c_0$\, ainsi obtenu qui est le \polmin de la \srl
\,$(a'_i)_{\,i\in\NN}$\, dont les \,$2n$\, premiers termes
co\"{\i}ncident avec les \,$2n$\, premiers termes
\,$a_0,a_1,\ldots,a_{2n-1}$\, de la suite \,$(a_i)$.
 
\ms Il consid\'ere alors la matrice \,$C$\, transpos\'ee de la
matrice compagnon du \poly \,$f(X)$: $$C=\left[
\begin{array}{ccccc} 0 & 1 & 0 & \ldots & 0 \\ 0 & 0 & 1 & \ldots
& 0 \\ \vdots & \vdots & \vdots & \ddots & \vdots \\ 0 & 0 & 0 &
\ldots & 1 \\ c_{0} & c_{1} & c_{2} & \ldots & c_{n-1}
\end{array}\right]\,,$$
Par exemple, pour \,$n=7$\, on obtient $$C := \left[
\begin{array}{ccccccc} 0 &  1 &  0 &  0 &  0 &  0 &  0\cr 0 &  0 &
1 &  0 &  0 &  0 &  0\cr 0 &  0 &  0 &  1 &  0 &  0 &  0\cr 0 &  0
&  0 &  0 &  1 &  0 &  0\cr 0 &  0 &  0 &  0 &  0 &  1 &  0\cr 0 &
0 &  0 &  0 &  0 &  0 &  1\cr -1 &  4 &  6 &  -10 &  -5 &  6 &  1
\end{array}\right]
$$ Le \polcar de \,$C$\, n'est autre que \,$\rP_C=(-1)^nf(X)$. Il
consid\`ere enfin les deux vecteurs: $$V=\left[\begin{array}{c}
a_0 \\ a_1 \\ \ldots \\ a_{n-1} \end{array}\right] ~~~\mbox{ et
}~~~ E_1 = \left[\begin{array}{c} 1 \\ 0 \\ \vdots \\ 0
\end{array}\right]~~ \mbox{ de }~~\A^{n\times 1}\,.$$
On v\'erifie \immt que les \srls \,$(a'_i)_{\,i\in \NN}$\, et
\,$(\tra{E_1}\,C^i\,V)_{\,i\in \NN}$, qui admettent un \polg
unitaire commun de degr\'e \,$n$, sont telles que \,$a'_i =
\,\tra{E_1}\,C^i\,V~(= a_i)$\, pour tout \,$i$\, compris entre
\,$0$\, et \,$2n-1$. On en d\'eduit que \,$a'_i =
\,\tra{E_1}\,C^i\,V$\, pour tout \,$i\in \NN$.
 
\ms Ainsi, par construction m\^eme de \,$C$, l'\algo de Wiedemann,
prenant \,$C$\, en entr\'ee avec les deux vecteurs \,$E_1$\, et
\,$V$, s'effectue avec les seules \ops d'addition et de \mul dans
\,$\ZZ$\, pour donner en sortie le \polmin de la \srl
\,$(\tra{E_1}\,C^i\,V)_{\,i\in \NN}$, et par cons\'equent le
\polcar de \,$C$.
 
\ms Soit maintenant \,$A=(a_{ij})$\, une \maca d'ordre \,$n$\, \`a
\coes dans \,$\A$. Il s'agit de calculer le \polcar de \,$A$\, en
n'utilisant que les \oparis de \,$\A$. Cela se fait par \elid dans
l'\algo de Wiedemann pour la matrice \,$A$\, en prenant comme
centre d'\elid le point form\'e par la matrice \,$C$\, et les deux
vecteurs auxilaires \,$E_1$\, et \,$V$. Comme les \coes du \polcar
de \,$A$\, (les sorties de l'\algo de Wiedemann) sont des \pols de
degr\'e $\leq n$ en les \coes \,$(a_{ij})$, on utilise l'\elid de
Strassen en degr\'e \,$n$.
 
On consid\`ere donc une \idtr \,$z$\, sur \,$\A$\,.
 
On pose \,$F=A-C$, et on applique l'\algo de Wiedemann dans
l'anneau \,$\A_n=\aqo{\A[z]}{z^{n+1}}$\, \`a la matrice
\,$B=C+zF$\, avec les vecteurs auxiliaires \,$E_1$\, et \,$V$. On
r\'ecup\`ere le \polcar de \,$A$\, en rempla\c{c}ant \,$z$\, par
$1$ dans les sorties.
 
Cet \algo calcule le \polgmin \,$g_z(X)\in \A_n[X]$\, de la \srl
\,$(\tra{E_1}\,B^{i}\,V)_{\,i\in \NN}$. Comme les seules divisions
se font par des \pols en \,$z$\, de terme constant \'egal \`a
\,$\pm 1$, l'ensemble du calcul se fait uniquement avec des
additions et \muls dans \,$\A$.
 
\ms D'o\`u l'\algo \vref{algokalto} de Kaltofen pour le calcul  du
\polcar d'une \maca \,$A\in\A^{n\times n}$.
 
\begin{algor}[Algorithme de Kaltofen-Wiedemann]
             \label{algokalto}
         \acl{algokalto}{Algorithme de Kaltofen-Wiedemann}
\Entree Un entier $n\geq 2$ et une matrice
        $A=(a_{ij})\in\A^{n\times n}$.
\Sortie Le \polcar $\rPA(X)$ de $A\,$. \Debut ~~~~~(on pose
$\A_n=\aqo{\A[z]}{z^{n+1}}$) \Varloc $i$, $k\in\N$ ;
        $V=(v_i)\in\ZZ^{n\times 1}$ (vecteur du centre d'\elidz);
        $C=(c_{ij})\in \ZZ^{n\times n}$
        (matrice du centre d'\elidz);
        $B\in (\A_n)^{n\times n}$;
    $(r_k)_{k=0..2n-1}\in (\A_n)^{2n}$.
\Eta{1}{Calcul du centre d'\elidz, et initialisation.}
     {C:= 0\in\A^{n\times n}} ;
\hsu \pour{i}{1}{n} \hsd $\ds v_i:={i-1 \choose \einf{
\frac{i-1}{2} }}$;
 $\ds c_{n,i}:=(-1)^{\einf{\frac{n-i}{2}}}{\einf{\frac{n+i-1}{2}
}\choose i-1}$ \hsu \finpour; \hsu \pour{i}{1}{n-1}
$c_{i,i+1}:=1\;$ \finpour; \hsu $B:=C+z\times(A-C)$;
\Etap{2}{Calcul de la \srl}
   {\pour{k}{0}{2n-1}}
\hsd $r_k:=$ premi\`ere coordonn\'ee de $B^k\times V$ dans $\A_n$
\hsu \finpour; \Etap{3}{Berlekamp-Massey}
  {Appliquer la \pcd de Berlekamp-Massey \`a la suite
  $(r_k)_{k=0..2n-1}$}
\hsu puis remplacer $z$ par $1$ dans le \polgmin trouv\'e. \fin
\end{algor}

\subsubsection*{Complexit\'e de l'algorithme}
 
\mni On utilise comme d'habitude la notation \vref{notamuP} ainsi
que la notation \vref{IDNConst}. L'\'etude de \com donne le
r\'esultat suivant d\^u \`a Kaltofen \cite{Kal}:

\begin{theorem}\label{thKalWied}
Le calcul du \deterz, du \polcar et de l'adjointe d'une \maca
d'ordre \,$n$\, sur un \acoma   \,$\A$\, se fait \`a l'aide d'une
\famu de \caris en 
\,$\SD(n\,\log{n},n^{{{\alpha+3}\over 2}}\,\mu_P(\esup{\sqrt n\,}))$. 
\end{theorem}
Si on utilise une \mul rapide des \pols en 
\,$\O(n\,\log{n})$\,
ou en \,$\O(n\,\log{n}\,\log \log n)$\, \oparis
(selon l'anneau consid\'er\'e), cela fait donc,
\,$\O(n^{{\alpha\over 2}+2}\,\log{n})$\, ou
\,$\O(n^{{\alpha\over 2}+2}\,\log{n}\,\log \log n)$\,  \oparis pour 
l'\algo de Kaltofen. 
Nous verrons au chapitre \ref{chap PolCarAnn}
que les \algos \parals en \prof \,$\log^2n$\, font moins bien
dans le cas d'un anneau vraiment arbitraire (ils utilisent
\,$\O(n^{\alpha +1}\,\log{n})$\, \oparisz) mais un peu mieux 
($\O(n^{\alpha +\frac{1}{2}})$\, \oparisz) dans le cas d'un
anneau o\`u les entiers \,$\leq n$\, sont non diviseurs de z\'ero.

Dans le cours de la preuve qui suit nous ferons \egmt l'analyse de 
\com de la version \elr de l'\algo de Kaltofen. Nous obtenons
le r\'esultat suivant.

\begin{prop}\label{propKalWied}
Dans la version \sqle simple de l'\algo de Kaltofen,
le calcul du \deterz, du \polcar et de l'adjointe d'une \maca
d'ordre \,$n$\, sur un \acoma   \,$\A$\, se fait \`a l'aide d'une
\famu de \caris de taille 
\,$\O(n^4)$\, et plus pr\'ecis\'ement avec un nombre
de \muls \'egal \`a \,$4n^4+\O(n^3)$\, et un nombre d'additions
du m\^{e}me ordre de grandeur. Le nombre de \muls essentielles est 
de \,$2n^4+\O(n^3)$.
\end{prop}

\smallskip \prv On remarque tout d'abord que le co\^ut de l'\'etape 1 est
n\'egli\-geable. Les entiers qu'elle calcule sont des constantes
du \cir dis\-ponibles une fois pour toutes et leur calcul ne doit
pas \^{e}tre pris en compte (ils sont de toute fa\c{c}on
calculables en \,$\O(n^2)$\, \oparisz). Quant \`a l'affectation
\,$B:=C+z\times(A-C)$\, dans \,$(\A_n)^{n\times n}$\, elle
signifie du point de vue des \oparis dans \,$\A$\, qu'on effectue
\,$2n-1$\, soustractions qui peuvent \^{e}tre effectu\'ees en une
seule \'etape \paralz.
 
\ss L'\'etape 3 est pour l'essentiel un \algo d'Euclide \'etendu.
Elle se fait avec un \cari de \prof \,$\O(n\,\log{n})$\, et de
taille \,$\O(n^2\,\mu_P(n))$\, o\`u $\,\mu_P(n)\,$ est le nombre
d'\oparis \ncrs pour la \mul de deux \pols de degr\'e \,$n$\, dans
\,$\A[z]$\, en \prof \,$\O(\log{n})$.
Cela est d\^u au
fait que l'\algo d'Euclide \'etendu utilis\'e comporte \,$\O(n)$\,
\'etapes avec chacune \,$\O(n)$\, \oparis dans l'anneau des \dlis
\,$\A_n$\, (certaines de ces \ops sont des divisions par des \elts
inversibles).

\ss Pour obtenir le r\'esultat \'enonc\'e, reste l'\'etape 2, la
plus co\^uteuse en nombre d'\oparisz.

\ss Voyons tout d'abord la version \elrz.
On calcule successivement les \,$V_k=B^kV$\, pour \,$k=1,\ldots ,2n-1$
\,par \,$V_{k+1}=BV_{k}$.
Cela fait en tout \,$2n^3-n^2$\, \muls et  \,$n(n-1)(2n-1)$\, 
additions dans \,$\A_{n}$.
Chacune des  \,$2n^3-n^2$\,  \muls est le produit d'une entr\'ee de \,$B$\, par une
coordonn\'ee de l'un des \,$V_{k}$. 
Or les entr\'ees de \,$B$\, sont des
\elts de la forme \,$c+bz$\, o\`u \,$c$\, est une
constante (une des entr\'ees non nulles de \,$C$) et \,$b$\, est une
entr\'ee de \,$A-C$. Un tel produit consomme donc \,$n$\, \muls
essentielles, \,$n+1$\, \muls du type \gui{produit d'un \elt de \,$\A$\,
par une constante} et \,$n$\, additions. 
En r\'esum\'e, l'\'etape 2 dans la version \sqle \elr consomme 
\,$2n^4-n^3$\,
\muls essentielles, \,$2n^4+\O(n^3)$\, \muls non essentielles
et \,$4n^4+\O(n^3)$\, additions.

\ss Voyons maintenant la version acc\'el\'er\'ee. On subdivise 
l'\'etape 2 en quatre
sous-\'etapes qui sont les suivantes, num\'erot\'ees de $2.1$ \`a
$2.4$, dans lesquelles on pose \,$r=\esup{ \sqrt{n}\,}$\,,
\,$s=\esup{2n/r }-1$, \,$U_0=E_{1}$\,  et \,$V_0=V$:
 
\mni\textsf{\textbf{\'Etape} $2.1$:
 \por{j}{1}{r-1} \textbf{calculer }
 $V_j :=B^{\,j}\,V_0$
\hsz \textbf{\'Etape} $2.2$: Calculer la matrice $B^{\,r}$ 
\hsz \textbf{\'Etape} $2.3$:
 \por{k}{1}{s} \textbf{calculer }
 $U_k := \,(\tra{B^r})^{k}\,E_1$
\hsz \textbf{\'Etape} $2.4$:
 \por{j}{0}{r-1}  et \por{k}{0}{s} \textbf{calculer }
\hsept $b_{kr+j}(z):=\,\tra{U_k}(z)\,V_j(z)\,.$}
 
\ss Notez que \,$r(s+1)\geq 2n$\,
 si bien que que les entiers \,$kr+j$\,
parcourent tout l'intervalle \,$[0,2n-1]$.

\ms Au cours des sous-\'etapes \,$2.1$\, et \,$2.2$, les \coes
calcul\'es sont des \pols en \,$z$\, de degr\'e \,$\leq r$\, (dans
\,$B^j\,V$, ils sont de degr\'e \,$\leq j\,$), \cad que chaque
\mul de deux \coes correspond \`a un \cari de \prof
\,$\O(\log{r})$\, avec \,$\mu_P(r)$\, \ops de base dans \,$\A$.
Cela donne l'analyse suivante pour les diff\'erentes
sous-\'etapes.
 
\ss $\bullet~$ {\bf Sous-\'etape\/} \,$2.1$: Pour obtenir tous les
vecteurs \,$B^j\,V$\, pour \,$1\leq j \leq r-1$\, on peut
proc\'eder en \,$\einf{ \log{r} }$\, \'etapes
\parals o\`u chaque \'etape \,$i$\, ($\,i=1,\dots,\einf{
\log{r} }\,$) consiste \`a \'elever au carr\'e la matrice
\,$B^{2^{i-1}}$\, puis \`a la multiplier \`a droite par la matrice
\,$[\,V~|~B\,V~|~\ldots~|~B^{2^{i-1}-1}\,V\,]$\, qui est une
matrice \,$n\times 2^{i-1}$\, pour obtenir la matrice
\,$[\,V~|~B\,V~|~\ldots~|~B^{2^{i}-1}\,V\,]$\, qui est une matrice
\,$n\times 2^{i}$\, dont les \coes sont des \pols de \,$\A[z]$\,
de degr\'e \,$< 2^{i}\leq r$.
 
Chacune de ces \,$\einf{ \log{r} }$\, \'etapes correspond
donc \`a un \cari de \prof \,$\O(\log{n}\log{r})$\, et de taille
\,$\O(n^{\alpha}\,\mu_P(r))$, ce qui donne au total, pour la
sous-\'etape \,$2.1$,  
un \cari de \prof \,$\O(\log^3{n})$\, et de taille
\,$\O(n^{\alpha}\,\mu_P(r)\,\log{n})\,$.
 
\ss $\bullet~$ {\bf Sous-\'etape\/} \,$2.2$: Si \,$r$\, est une 
puissance de 2, le calcul de
\,$B^{\,r}$\, se fait en \'elevant au carr\'e
la matrice \,$B^{{r}/{2}}$\, d\'ej\`a calcul\'ee. 
Sinon il faut 
faire le produit de certaines des matrices \,$B^{2^{i-1}}$: par exemple
si \,$r=39=32+4+2+1$, on a \,$B^{39}=B^{2^5}B^{2^2}B^2B$.
Pour chaque 
produit les \coes des matrices sont de
degr\'e \,$\leq r/2$\, dans \,$\A[z]$.
Ceci correspond de nouveau \`a un \cari de \prof \,
$\O(\log^3{n})$\, et de taille
\,$\O(n^{\alpha}\,\mu_P(r)\,\log{n})$. 
Pour la suite nous posons \,$B_1=\,\tra{\!B^r}$\,
 
\ss $\bullet~$ {\bf Sous-\'etape\/} \,$2.3$: 
Nous ne pouvons plus utiliser la technique de
l'\'etape 2.1 qui ici donnerait a priori
une \famu de \caris dans 
\,$\SD(\log^3{n},n^{\alpha}\,\mu_P(n)\,\log{n})$.

Partant du vecteur
\,$U_0 = E_1$, la sous-\'etape \,$2.3$\, de notre \algo consiste
\`a calcu\-ler, pour \,$k$\, allant de $1$ \`a \,$s$, le vecteur
\,$U_k(z) = B_1\,U_{k-1}(z)$. 
Posons  
\,$s_{1}=\esup{(n+1)/r }$. Notons que \,$U_{k-1}=
B_1^{k-1}\,E_1$\, se r\'e\'ecrit dans \,$\A_{n}$\,
 sous la
forme \,$U_{k-1}(z) = \sum_{\ell=0}^{s_{1}-1}\,z^{r \ell}\,U_{k-1,\ell}$\,
o\`u chacun des \,$U_{k-1,\ell}$\, est un vecteur dont les
composantes sont des \pols en \,$z$\, de degr\'e \,$< r$. On peut
donc identifier \,$U_{k-1}(z)$\, avec la matrice \,$n\times s_{1}\,$: 
$$W_{k}(z)=[\,U_{k-1,0}~|~U_{k-1,1}~|~\ldots~|~U_{k-1,s_1-1}\,]\,.$$

Le calcul du vecteur \,$U_k(z)$\, \`a \,$n$\,
lignes et \,$s_1$\, colonnes se fait comme suit. On calcule
la matrice \,$B_1W_{k-1}(z)$\, dont les entr\'ees sont des 
\pols de degr\'e \,$\leq 2r$, puis on
r\'eorganise les sommes correspondantes pour obtenir \,$U_{k}(z)$
(ce qui n\'ecessite au plus \,$n^2$\, additions dans \,$\A$).
Le produit \,$B_1W_{k-1}(z)$\, est celui d'une matrice
\,$n\times n$\, par une matrice \,$n\times s_1$, toutes 
les entr\'ees \'etant de degr\'e \,$\leq r$.
Ceci peut se
faire avec \,$r^2$\, \muls \parals de blocs \,$s_1 \times
s_1 $. Chaque \mul de blocs se fait en \,$\O(s_1^\alpha)$\,
\oparis sur des \pols de degr\'e \,$\leq r$. On obtient
donc chaque  \,$U_k(z)$\, en 
\,$\SD(\log^2r,r^{2+\alpha}\,\mu_P(r))$.
 
Cela donne au total, pour la sous-\'etape \,$2.3$, une
\famu de \caris dans \,$\SD(s\,\log^2n,s\,r^{2+\alpha}\mu_P(r))$\,
\cad encore dans
\,$\SD(n^{\frac{1}{2}}\,\log^2n,n^{\frac{3+\alpha}{2}}\mu_P(\esup{\sqrt{n}\,})$.
 
\ss $\bullet~$ {\bf Sous-\'etape\/} \,$2.4$: Cette \'etape peut
\^etre \egmt ramen\'ee \`a la \mul d'une matrice \,$(s+1) \times
n$\, par une matrice \,$n\times r$: $$\left[\begin{array}{c}
\,\tra{U_0}(z) \\ \,\tra{U_1}(z) \\ \vdots \\ \,\tra{U_s}(z)
\end{array}\right]~\times~
\left[\begin{array}{c|c|c|c} V_0(z) & V_1(z) & \cdots &
V_{r-1}(z)\,
\end{array}\right]~=~~~~~~~~~~~~~~~~~~~~$$
$$~~~~~~~~~~~~~~~~~~~~~~~~~~~~~~ \left[\begin{array}{ccc}
\,\tra{U_0}(z)\,V_0(z) & \,\cdots\, & \,\tra{U_0}(z)\,V_{r-1}(z)
\\ \vdots & \ddots & \vdots \\ \,\tra{U_s}(z)\,V_0(z) & \,\cdots\,
& \,\tra{U_s}(z)\,V_{r-1}(z)
\end{array}\right]$$
dont l'\elt en position \,$(k+1,j+1)$\, pour \,$0\leq k\leq s$\,
et \,$0\leq j\leq r-1$\, n'est autre que le \coe recherch\'e:
$~\,\tra{U_k}(z)\,V_j(z) = b_{kr+j}(z).$
 
Utilisant \`a nouveau la \mul par blocs \,$(s+1)\times
(s+1)$, nous concluons que la sous-\'etape \,$2.4$\,
correspond \`a un \cari de \prof \,$\O(\log^2{n})$\, et de taille
\,$\O(n^{{\alpha+2}\over 2}\,\mu_P{(\esup{\sqrt{n}\,})})$. 
 
\ss On peut r\'esumer le calcul de \com dans le tableau
(\ref{cokal}) suivant qui donne, pour chaque \'etape, la \com
\arith du \cir correspondant, en m\^eme temps que le r\'esultat
\gnlz. Nous avons \'egalement indiqu\'e la taille lorsqu'on
ex\'ecute l'\algo avec une
\mul acc\'el\'er\'ee des \pols mais sans \mul rapide des matrices,
sur les lignes \gui{avec \,$\alpha =3$}.
 
\ms
\begin{center}
\label{cokal} \setlength{\extrarowheight}{1.5mm}
\begin{tabular}{|l|c|c|}
\hline Etape~ & ~Profondeur~ & ~Taille~ \\[1mm] 
\hline Etape~1 &
$\O(1)$ & n\'egligeable \\[1mm] \hline 
Etape~2 & $\O(n^{1\over
2}\,\log^2{n})$ & 
\,$\O(n^{{\alpha+3}\over 2}\,\mu_P{(\esup{\sqrt{n}\,})})$\, \\ 
~avec \,$\alpha =3$ & $\cdots $ & 
\,$\O(n^3\,\mu_P{(\esup{\sqrt{n}\,})\,\log\,n})$\, 
\\[1mm] 
\hline Etape~3 &
\,$\O(n\,\log{n})$\, & \,$\O(n^2\,\mu_P{(n)})$\, \\[1mm] 
\hline
Total & \,$\O(n\,\log{n})$\, & \,$\O(n^{{\alpha+3}\over
2}\,\mu_P{(\esup{\sqrt{n}\,})})$\, \\ 
~avec \,$\alpha =3$ & $\cdots $ & 
\,$\O(n^3\,\mu_P{(\esup{\sqrt{n}\,})\,\log\,n})$\, 
\\[1mm] 
\hline
\end{tabular} \\
 
\ms {\small {\bf Tableau \ref{cokal}}} \\ {\small Complexit\'e de
l'\algo de Kaltofen-Wiedemann} \\
 
\end{center}
\qed 

\ms Dans notre preuve c'est l'\'etape 3 qui d\'etermine la \prof
du \cari correspondant \`a l'\algo de Kaltofen-Wiedemann. Mais on
peut r\'eduire la \prof de l'\'etape 3 par diverses \metsz.
 
Une premi\`ere  est de ne pas utiliser l'\algo de
Berlekamp/Massey pour le calcul du \polmin d'une \srlz. Une telle
\metz, d\'evelopp\'ee dans \cite{KalPan} (voir aussi \cite{BP})
ram\`ene ce calcul \`a la r\'esolution d'un \sli qui a la forme de
Toeplitz, en utilisant le calcul du \polcar de sa matrice par la
\met de Le Verrier am\'elior\'ee par Csanky (cf. section \ref{sec
Csanky}). On obtient un \cari de \prof de \,$\O(\log^3{n})$ et de
m\^eme taille, \cad \,$\O(n^2\,\mu_P{(n)})$\,.
L'inconv\'enient de cette am\'elioration est qu'elle s'applique
uniquement lorsque \,$n!$\, ne divise pas z\'ero dans l'anneau
\,$\A$.
 
Une deuxi\`eme \metz, qui ne se heurte pas \`a l'obstacle
pr\'ec\'edent, consiste \`a utiliser une version
\parase de l'\algo d'Euclide \'e\-ten\-du.
Voir \cite{Moe,LiRo} et \cite{VonZurbook} corollaire 11.6
page~304.
 
\ss Cependant, il ne suffit pas de r\'eduire la \prof de l'\'etape
3 pour obtenir une \prof \pogz. Il faudrait le faire \egmt pour
l'\'etape 2 et plus pr\'ecis\'ement la sous-\'etape 2.3.

On a
donc \`a l'heure actuelle un \pb ouvert: peut-on obtenir un \cirz,
de la taille de cet \algo et de \prof \pogz, permettant de
calculer le \polcar sur un \acomaz?
 
\ms L'\algo de Kaltofen-Wiedemann obtient le r\'esultat
asymptotique ci-dessus, \`a savoir 
\,$ \O(n^{{\alpha+3} \over 2}\,\mu_P(\esup{\sqrt n\,}))$, 
le meilleur temps \sql de tous
les \algos connus pour le calcul du \polcar sur un \acomaz,
gr\^ace \`a la \mul rapide des matrices, bien s\^ur, mais aussi
gr\^ace \`a la \mul rapide des \polsz. Et pour les \pols la \mul 
rapide est d\'esormais couramment impl\'ement\'ee sur machine. 

Ainsi lorsqu'on ne dispose pas d'une \mul rapide des matrices,
on obtient un temps \sql asymptotiquement meilleur que tous les 
autres \algos fonctionnant sur un \acomaz, d\`es qu'on acc\'el\`ere la
\mul des \polsz, ne serait-ce que par la \met de 
Karatsuba.

Notons que sur un anneau commutatif qui ne poss\`ede pas de 
racines principales de l'unit\'e, la \met qui utilise la
transformation de Fourier rapide est en 
\,$\O(n\,\log\,n\,\log\log\,n)$\,
et elle ne devient plus performante que la \met de Karatsuba en
\,$\O(n^{\log{3}})$\, que pour \,$n$\, tr\`es grand, de l'ordre
de plusieurs milliers (cf. section \ref{TFDmrGen} 
et notamment la remarque \vref{rem CaK}). 

Un vaste
champ d'exp\'erimentation s'ouvre donc, maintenant que
diff\'erentes \muls rapides commencent \`a avoir une r\'eelle
port\'ee pratique en calcul formel.

\subsection*{Conclusion}

Nous terminons ce chapitre en renvoyant le lecteur \`a deux surveys
r\'ecents d'Erich Katofen et Gilles Villard \cite{KaV,KaV2} concernant 
la complexit\'e aussi bien \agq que binaire du calcul des \deters
(nous nous int\'eressons plut\^ot au calcul du \polcar
dans cet ouvrage).

Ils montrent \`a quel point 
l'\agr \lin est un sujet de recherche actif en calcul formel
et l'importance des m\'ethodes modulaires et seminu-m\'eriques
pour le traitement des \pbs concrets.

\newpage \thispagestyle{empty}
 
\chapter{Parall\'elisations de la \met de Leverrier}
\minitoc
\acvide
\label{chap par Lever}

\subsubsection*{Introduction}   

Csanky \cite{Csan} fut le premier \`a  
prouver que les probl\`emes du calcul  
des \detersz, de l'inversion des  
matrices, de la r\'esolution des \syses \lins et du calcul  
du \polcarv dans le cas d'un 
anneau contenant le corps des  
rationnels, sont dans la classe  
\,$\NC$, \cad dans la classe des probl\`emes qui  
peuvent \^etre r\'esolus en temps  
\paral \pog avec  
un nombre \poll de processeurs
par une \famu de \carisz. 
  
\ss Il montre, en effet, que tous ces probl\`emes  
se ram\`enent au calcul du \polcar  
et que ce dernier se calcule en  
\,$\SD(n^{\alpha +1},\log^2{n})$.  
En particulier ils sont dans la classe  
\,$\NC^2$.
  
\ms Nous pr\'esentons le travail de Csanky
dans la section \ref{sec Csanky}. 
Dans la section suivante nous donnons
l'am\'elioration  
due \`a \PrSa \cite{PrSa}
qui montre que le calcul du \polcar peut \^etre
r\'ealis\'e dans \,$\SD(n^{\alpha +1/2},\log^2{n})$.
Dans la section \ref{sec GaliPan} nous donnons une  
meilleure estimation de la \com th\'eorique  
de l'\algo pr\'ec\'edent, l\'eg\`erement am\'elior\'e, 
 due \`a \GaPa \cite{GalPan}.
 
 \ms
Dans le chapitre \ref{chap PolCarAnn}, nous examinerons des
\algos qui r\'esolvent les m\^emes probl\`emes sur un anneau 
commutatif arbitraire.  

\section{Algorithme de Csanky}
  \label{sec Csanky}
  
Pour calculer le \polcar, Csanky utilise la  
\met de Le Verrier en la  
\paralt de la mani\`ere suivante.  
  
On se donne un entier \,$n$, un corps  
\,${\cal K}$\, (ou plus \gnlt un  
anneau dans lequel \,$n!$\, est inversible) et une  
matrice \,$A\in {\cal K}^{n\times n}$\, de \polcar: \\  
$$P(X)=\det{(A-X\In)}=(-1)^n[X^n -  
c_1X^{n-1} -\dots - c_{n-1}X -c_n\,]\,.$$  
On pose \,$s_k=\Tr(A^k)$\, pour \,$k=1,2,\dots,n$.  
  
\ss  
La \met de Le Verrier consiste \`a r\'esoudre  
l'\'equation   
\begin{equation} \label{LeV}
    S\,\vec{c}=\vec{s}
\end{equation}
 o\`u
$$
\vec{c}=\left[
\begin{array}{c} c_1\\ \vdots\\ \vdots\\ c_n 
\end{array}
          \right]~,  
~\vec{s}=\left[
\begin{array}{c} s_1\\ \vdots \\ \vdots\\ s_n 
\end{array}
\right]  
~~\mathrm{et}~~S = \left[
\begin{array}{ccccccc}  
1 & 0 & \cdots & \cdots  & 0 \\  
s_1 & 2 & \ddots &  & \vdots  \\[1mm]  
\vdots & \ddots & \ddots & \ddots & \vdots \\  
s_{n-2} &       & \ddots & \ddots  & 0 \\  
s_{n-1} & s_{n-2} & \cdots & s_{1} & n  
\end{array}
\right].
$$

Cette \'equation admet la solution unique  
\,$\vec{c}=S^{-1}\vec{s}$\, qui donne les  
\coes du \polcarz.  
  
Ceci donne l'\algo de Csanky \vref{algoCsan}  
en quatre grandes \'etapes.  
  
\begin{algor}[\Algo de Csanky, principe \gnlz.] 
\label{algoCsan}
    \acl{algoCsan}{\Algo de Csanky}
\Entree Un entier \,$n\in\N$ et une matrice $A\in\A^{n\times n}$. 
L'anneau $\A$ contient le corps $\QQ$.
\Sortie Les \coes du \polcar $\rPA(X)$ de $A$. 
\Debut 
\Etap{1}{}{Calculer en \paral les  
puissances \,$A^2,A^3,\dots,A^n$;}  
\Etap{2}{}{Calculer en \paral les  
traces \,$s_1,s_2,\dots,s_n$\, des  
matrices} 
\hsu$A,A^2,\dots,A^n$. 
\Etap{3}{}{Cr\'eer et inverser la \matg \,$S$\, (\'equation
\ref{LeV}).}
\Etap{4}{}{Calculer le produit  
\,$S^{-1}\vec{s}=\vec{c}$.}
\fin 
\end{algor}
  
L'analyse de \com pour cet  
\algo utilise les r\'esultats de  
\com de la technique  \gui{diviser pour gagner} 
et notamment son application au  
calcul \paral de l'inverse d'une  
\matg que nous avons  
d\'ecrite au \S~\ref{invtrian}. 

\subsubsection*{La complexit\'e de l'algorithme}  
  
\noi  
$\bullet~$  
Le calcul en \paral des puissances  
\,$A^2,\dots,A^n$\, de la matrice  
\,$A$\, se ram\`ene \`a un \algo de calcul  
\paral des pr\'efixes  repr\'esent\'e par un  
\cari \paral de \prof  
\,$\O(\log{n})$\, et de taille major\'ee  
par \,$4n$\, 
(\tho \vref{prefix}), 
mais dont les n{\oe}uds internes  repr\'esentent eux-m\^emes  
des \cirs de \mul de matrices \,$n\times n$, \cad des \cirs  
de taille \,$\O(n^\alpha)$\, et de \prof  
\,$\O(\log{n})$. Ce qui donne au total,  
pour r\'ealiser l'\'etape~1 un \cari 
en \,$\SD(n^{\alpha+1},\log^2{n})$.  
  
\sni  
$\bullet~$  
On calcule ensuite les traces des matrices  
\,$A,A^2,\dots,A^n$, \cad les  
coefficients \,$s_k=\Tr{(A^k)}$\, qui  
forment la \matg \,$S$.\\  
Ce sont des sommes de \,$n$\, \elts  
de \,$\A$\, que l'on calcule en \paral pour  
\,$1\leq k\leq n$\, en \,$\SD(n^{2},\log{n})$.  
  
\sni  
$\bullet~$  
Le calcul de \,$S^{-1}$\, se fait comme indiqu\'e  
au \S~\ref{invtrian}.  
La matrice \,$S$\, est en effet triangulaire et \fregz. 
D'apr\`es la proposition  
\ref{divtri}, le calcul de la matrice  
\,$S^{-1}$\, se fait par un \cari  
\paral en \,$\SD(n^{\alpha},\log^2{n})$.  
  
\sni  
$\bullet~$  
Enfin, le calcul de \,$\vec{c} = S^{-1}\vec{s}$,  
qui est le produit d'une \matg par  
un vecteur, se fait en \paral par un \cari  
de taille \,$n^2$\, et de \prof  
\,$\esup{ \log{n} }$, la \prof \'etant  
essentiellement due aux additions. 

En fait, on a un tout petit peu mieux.

\begin{theorem}  \emph{(Csanky)}\label{thCsa}\\
Soit \,$\A$\, un anneau v\'erifiant les hypoth\`eses
pour l'\algo de Le Verrier: la division par \,$n!$,  
quand elle est possible, est unique  et explicite.
Le calcul du \polcarz, de l'adjointe et l'inverse d'une \maca
d'ordre \,$n$\, 
 est en \,$\SD(n^{\alpha+1},log^2{n})$. \\  
\end{theorem}
  
\prv  
 Une l\'eg\`ere modification de l'\algo de Csanky  
pour le \polcar d'une \maca d'ordre  
\,$n$\, montre que l'hypoth\`ese d'un anneau dans lequel \,$n!$\, est  
inversible, peut \^etre remplac\'ee par  
l'hypoth\`ese pour l'\algo de Le Verrier.  
En effet soit \,$A \in \A^{n\times n}$\, et  
\,$S$\, la matrice utilis\'ee dans l'\algo de  
Csanky pour le calcul du \polcar.  
  
Au lieu de calculer \,$S^{-1}$\, (ce qui n'est  
possible que si \,$n!$\, est inversible dans  
\,$\A\,$), on calcule \,$n!\,S^{-1}$.  
Il suffit pour cela de d\'evelopper le \polcar  
de \,$S$\, en calculant le produit  
\,$(X-1)(X-2)\cdots (X-n)$, ce qui revient  
\`a calculer les valeurs des \pols  
sym\'etriques \elrs  
\,$\sigma_1,\sigma_2,\dots,\sigma_n$\, de \,$n$\,  
variables au point \,$(1,2,\dots,n\,$).  
  
Le \tho de Cayley-Hamilton permet alors  
d'\'ecrire:  
$$(-1)^{n+1}\,n!\,S^{-1} = S^{n-1} +  
\somm_{k=1}^{n-1}(-1)^k\sigma_kS^{n-k-1}$$  
ce qui ram\`ene le calcul de \,$n!\,S^{-1}$\,  
\`a celui des puissances \,$S^2,S^3,\dots,S^{n-1}$.  
  
Or ce calcul se fait en \paralz, d'apr\`es le  
calcul des pr\'efixes par exemple (proposition  
\ref{prefix})\, en \,$\SD(n^{\alpha +1},\log^2{n})$. 

Nous laissons le lecteur
ou la lectrice terminer pour ce qui concerne les calculs de  l'adjointe 
et de l'inverse.\qed 

\subsubsection*{Variante de Sch\"onhage}
Signalons  qu'il existe  une variante de   
la \met de Csanky/Le Verrier due \`a Sch\"onhage \cite{Schon}
qui donne une \famu de \caris avec divisions calculant
le \polcar avec une faible \prof sur un corps de
\cara finie.

Sch\"onhage  utilise le r\'esultat suivant concernant les sommes  
de Newton (\S~\ref{subsec Newton}) connu sous le nom de  
\emph{crit\`ere de Kakeya\/} \cite{Kak}:  
\index{Kakeya!crit\`ere de}
\index{Sch\"onhage!variante de la \met de Csanky}

\begin{prop} Soit \,$J$\, une partie finie \`a \,$n$\,  
\elts de \,$\NN$\, et \,$(s_j)_{j\in J}$\, le  
\sys correspondant de \,$n$\, sommes de Newton
\`a \,$n$\, \idtrs sur un corps \,$\K$\, de \cara nulle. 
Alors \,$(s_j)_{j\in J}$\, est un \sys fondamental  
de \pols \syms sur \,$\K$\, 
(\cf d\'efinition \ref{defsfps}) si et seulement
si \,$\NN\,\backslash J$\, est stable pour l'addition dans  
\,$\NN$.  
\end{prop}  

Par exemple, pour tout entier \,$p$\, positif, la  
partie \,$J(p,n)\subset \NN\setminus p\,\NN$\, constitu\'ee  
des \,$n$\, premiers entiers naturels qui ne sont  
pas des multiples de \,$p$,
satisfait ce crit\`ere, et Sch\šnage \cite{Schon} l'utilise  
pour adapter la \met de Le Verrier au calcul du
\polcar sur un corps de  
\cara \,$p>0$.

Notez qu'en \cara \,$p$\, l'\'egalit\'e \,$(x+y)^p=x^p+y^p$\, implique 
que les
sommes de Newton v\'erifient  les
\egts \,$s_{kp}={s_k}^p$.

Prenons maintenant un exemple. Le \pol \gnl de degr\'e
8 est \,$P(X)=X^8-\sum_{i=1}^8 a_iX^{8-i}$.
Si nous sommes sur un corps de \cara 3,  nous consid\'erons les
8 premi\`eres relations de Newton qui donnent les
sommes \,$s_j$\, pour  \,$j\in\NN\setminus 3\,\NN$\,
 (cf. l'\'egalit\'e (\ref{Newt3}) page \pageref{Newt3}):
\begin{equation} \label{EqKS} \left[ 
{\begin{array}{cccccccc}
1 & 0 & 0 & 0 & 0 & 0 & 0 & 0 \\
s_1 & 2 & 0 & 0 & 0 & 0 & 0 & 0 \\
s_3 & s_2 & s_1 & 1 & 0 & 0 & 0 & 0 \\
s_4 & s_3 & s_2 & s_1 & 2 & 0 & 0
 & 0 \\
s_6 & s_5 & s_4 & s_3 & 
s_2 & s_1 & 1 & 0 \\
s_7 & s_6 & s_5 & s_4 & 
s_3 & s_2 & s_1 & 2 \\
s_9 & s_8 & s_7 & s_6 & 
s_5 & s_4 & s_3 & s_2 \\
s_{10} & s_9 & s_8 & s_7 & 
s_6 & s_5 & s_4 & s_3
\end{array}}
 \right] ~
\cmatrix{a_1\cr a_2\cr a_3\cr a_4\cr a_5\cr a_6\cr a_7\cr a_8 }~=~
\cmatrix{s_1\cr s_2\cr s_4\cr s_5\cr s_7\cr s_8\cr s_{10}\cr s_{11} }
\end{equation}
Compte tenu des relations
\begin{equation} \label{EqKS2}
    s_3={s_1}^3, \,s_6={s_2}^3, \,s_9={s_1}^9,
\end{equation}
le \deter de la \maca est \'egal \`a
$$
\begin{array}{c}
\begin{array}{rcl}
d &=& -s_1\,s_2^3\,s_5 +s_2^2\,s_4^2+s_1^3\,s_4\,s_5 +
s_1^5\,s_2\,s_5+ s_1^3\,s_2^2\,s_5- s_1^2\,s_2^3\,s_4+\\ &&
s_1^2\,s_5^2+ s_1^4\,s_2^4+ s_1^6\,s_2^3- s_1^5\,s_7-
s_1^4\,s_8-
s_4^3- s_1\,s_7\,s_2^2+ s_1^{12}\\&&
-s_1^4\,s_4^2+
s_1^2\,s_2^5+ s_1^{10}\,s_2+ s_2^6+ s_2^2\,s_8- s_1^8\,s_2^2+ 
s_1^8\,s_4-s_1^2\,s_2\,s_8\\ &&
+s_4\,s_8- s_5\,s_7+  s_1^4\,s_2^2\,s_4+ s_2\,s_5^2-
s_1^6\,s_2\,s_4+
s_1\,s_2\,s_4\,s_5
\end{array}\end{array}
$$

Un point non trivial est que \,$d$\, n'est pas
une fonction identiquement nulle (si le corps de base est infini).
En fait dans le cas g\'en\'erique, \cad si on consid\`ere les
\,$a_i$\, comme des \idtrs et les \,$s_i$\, 
($i=1,\,2,\,4,\alb\,5,\alb\,7,\,8,\,10,\,11$) comme
donn\'es par les relations (\ref{EqKS}) et (\ref{EqKS2}), 
les \'el\'ements \,$s_i$\, 
sont alg\'ebri\-quement ind\'ependants.
Cela implique alors que les \,$a_{i}$\, ($i=1,\ldots ,\alb8$) peuvent 
s'exprimer
comme fractions rationnelles en les \,$s_{i}$\, 
($i=1,\,2,\,4,\alb\,5,\alb\,7,\,8,\,10,\,11$)
avec \,$d$\, pour \denoz.

Un autre point non trivial consiste \`a r\'esoudre les 
\slis du type (\ref{EqKS}) 
(lorsque  le \deter correspondant
est non nul) par un \algo (avec divisions) bien \parasz.

L'\algo de Sch\šnage \cite{Schon}  
correspond \`a  une \fam de \caris (avec divisions) dans  
\,$\SD(n^{\alpha +1},\log^2{n})$\, 
(voir aussi le livre \cite{BP} Annexe C pages 372--377). 
  
\section{Am\'elioration de Preparata et Sarwate}
\label{sec PrepaSar}

\subsubsection*{Principe g\'en\'eral}
Consid\'erons un anneau \,$\A$\,  
 v\'erifiant les hypoth\`eses
pour l'\algo de Le Verrier, et  
une matrice  carr\'ee \,$A\in\A^{n\times n}$.  
L'am\'e\-lio\-ration apport\'ee par \PrSa  
\cite{PrSa} \`a l'\algo  de Csanky provient du fait 
que pour calculer les  
traces \,$s_k=\Tr(A^k)~(1\leq k\leq n)$, on n'a pas  
besoin de calculer toutes les puissances de \,$A$.  
  
\ms  
Il suffit en effet, si l'on pose  
\,$p=\esup{\sqrt{n}\,}$, de disposer  
des \,$2p$\,  matrices  
\,$I_n,\,A,\,\ldots,\,A^{p-1}$\, et  
\,$A^p=B,\,B^2,\,\ldots,\,B^p=A^{p^2}$,  
ce qui revient \`a calculer  
\,$2\,\esup{\sqrt{n}\,}-2$\,  
puissances de matrices \,$n\times n$\,  
au lieu des \,$n-1$\, puissances de \,$A$.  
Il est fait appel pour cela \`a deux \pcds  
r\'ecursives not\'ees \,\textsf{Powers}$(A,r)$\, et  
\textsf{Superpowers}$(A,r)$\, permettant de calculer  
les puissances successives d'une matrice  
carr\'ee \,$A$\, jusqu'\`a l'ordre $r$.  
  
Les traces des puissances de \,$A$\, seront  
alors obtenues en consid\'erant les matrices  
\,$U_j~(1\leq j\leq n)$\, d\'efinies de la  
mani\`ere suivante: \\  
\,$U_j=L_jC_j$\, o\`u \,$L_j\in \A^{p\times n}$\,  
est la matrice form\'ee uniquement des \,$j\,$-\,\`emes  
lignes ($\,1\leq j\leq n\,$) des $p$ matrices  
\,$I_n,\,A,\,\ldots,\,A^{p-1}$\, et o\`u  
\,$C_j\in \A^{n\times p}$\, est la matrice form\'ee  
des \,$j\,$-\,\`emes colonnes ($\,1\leq j\leq n\,$) des autres  
matrices \,$A^p,\,A^{2p},\,\ldots,\,A^{p^2}$.  
  
Les matrices \,$U_j~(1\leq j\leq n)$\, sont des matrices  
carr\'ees  d'ordre \,$p$\, dont les \,$p^2$\, \coes  
ne sont autres que les \,$j\,$-\,\`emes \,\emph{\elts  
diagonaux\/} des matrices  
\,$A^p,\,A^{p+1},\,\ldots,\,A^{p^2+p-1}$.  
  
Plus pr\'ecis\'ement, l'\elt \,$u_{kl}^{[j]}$\,  
qui est position \,$(k,l)$\, dans la matrice \,$U_j$\,  
et qui est obtenu par \mul de la \,$j\,$-\,\`eme  
ligne de la matrice \,$A^{k-1}$\, par la \,$j\,$-\,\`eme  
colonne de la matrice \,$A^{pl}$\, est donc le  
\,$j\,$-\,\`eme \elt de la diagonale du produit  
\,$A^{k-1}A^{pl}=A^{pl+k-1}$, \cad que   
\,$u_{kl}^{[j]}=a_{\,jj}^{[pl+k-1]}$\, pour   
\,$1\leq k,l \leq p$, si l'on d\'esigne par   
\,$a_{\,rs}^{[m]}$\, l'\elt en  position  
\,$(r,s)$\, de la matrice \,$A^{m}$.

Posant \,$m=pl+k-1$\, ($\,m$\, prend toutes les  
valeurs comprises entre \,$p$\, et \,$p^2+p-1$\,  
quand \,$k$\, et \,$l$\, varient de 1 \`a \,$p\,$)   
on obtient, avec les notations ci-dessus, et pour  
\,$p\leq m\leq p^2+p-1$:
$$\Tr (A^m) = \sum_{j=1}^n a_{\,jj}^{[m]} =  
\sum_{j=1}^n u_{kl}^{[j]}
$$
(o\`u \,$l$\, et \,$k-1$\, sont respectivement  
le quotient et le reste euclidiens de  
\,$m$\, par \,$p\,$).

Comme les matrices \,$A,\,\ldots,\,A^{p-1}$\,  
sont d\'ej\`a disponibles, cela nous donne donc  
les traces de toutes les puissances  
\,$A,\,\ldots,\,A^{p},\,\ldots,\,A^{p^2}$\,  
donc celles de toutes les matrices  
\,$A,\,\ldots,\,A^n$\, puisque  
\,$p^2+p-1 \geq n + \sqrt{n} -1 \geq n$.  
  
\ms  
D'o\`u l'\algo de \PrSa qui comprend deux  
parties, la premi\`ere pour le calcul du \polcar de la  
matrice donn\'ee \,$A\in\A^{n\times n}$, et la  
deuxi\`eme pour le calcul de l'adjointe et de l'inverse  
de cette matrice.  

\subsubsection*{Calcul du \polcar}

Avant de donner l'\algo \vref{algoPrSa}, voyons tout d'abord  
les sous-\pcds utilis\'ees dans cet \algoz.  
Il s'agit essentiellement de la \pcd  
\textsf{Superpowers} qui est d\'efinie de mani\`ere r\'ecursive  
\`a partir de la \pcd \textsf{Powers} (elle-m\^eme  
d\'efinie de mani\`ere r\'ecursive) en vue d'acc\'el\'erer  
le calcul des puissances d'une \maca donn\'ee  
(dans notre cas, c'est la matrice \,$A\in\A^{n\times n}\,$).  
  
Chacune de ces deux sous-\pcds prend donc en  
entr\'ee \,$A$\, et un entier \,$p > 1$\, et donne en  
sortie la \marc \,$n\times np$\, form\'ee  
des \,$p$\, puissances de \,$A$:  

\centerline{\textsf{Powers}$(A,p)~=~$  
\textsf{Superpowers}$(A,p)~=~[\,A~|~A^2~|~\ldots~|~A^p\,]$.}  
  
\mni \textsf{Powers$(A,s)$ \\  
$\bullet~$ $m:= \esup{ s/2 }$;\\  
$\bullet~$ $[\,A\,|\,\ldots\,|\,A^m\,] :=$ Powers$(A,m)$;\\  
$\bullet~$ \pour{i}{m+1}{s}  
$A^i:= A^{\einf{ i/2 }}\,A^{\esup{ i/2 }}$}.  

\mni  \textsf{Superpowers$(A,p)$ \\  
$\bullet~$ $r:= \esup{ \log{p} }$; \\  
$\bullet~$ $s:= \einf{ p/r }~;~~q:= p-rs$; \\  
$\bullet~$ $[\,A\,|\,\ldots\,|\,A^s\,] :=$ Powers$(A,s)$;\\  
$\bullet~$ \pour{k}{1}{r-1} \,$A^{sk}\,\times$ Powers\,$(A,s)$;  
\hsu (cela donne toutes les puissances  
de \,$A$\, jusqu'\`a l'ordre \,$rs$) \\  
$\bullet~$ \pour{i}{1}{q}
 \,$A^{sr}\,\times [\,A ~|~ \ldots ~|~ A^q\,]$;  
\hsu (pour avoir les \,$q=p-rs$\, puissances  
restantes de \,$A$).}  
\begin{algor}[\Algo de \PrSa] 
\label{algoPrSa}
    \acl{algoPrSa}{\Algo de \PrSa}
\Entree Un entier $n$ et une matrice $A\in\A^{n\times n}$. 
\Sortie Le vecteur \,$\vec{c}$\,  
des coefficients du \polcar de \,$A$. 
\\[1mm]\textbf{Les \'etapes du calcul} \,\,\,
(avec \,$p = \esup{\sqrt{n}\,}$),
\begin{itemize}  
\item[1.~] Calculer les puissances  
\,$A,\,\ldots,\,A^{p}$\,  
en appelant Superpowers$(A,p)$;  
\item[2.~] Calculer les puissances  
\,$A^p,\,\ldots,\,A^{p^2}$\, en faisant  
Superpowers$(A^p,p)$;  
\item[3.~] Calculer en \paral les \,$n$\, produits
\,$U_j=L_jC_j$\, ($\,1\leq j\leq n\,$);
\item[4.~] Former le vecteur \,$\vec{s}$\, et  
la \matg \,$S$\, en calculant en  
\paralz, \`a partir des matrices  
\,$U_j=(u_{kl}^{[j]})$\, obte\-nues \`a l'\'etape  
pr\'ec\'edente, les \,$n$\, traces  
$~s_m=\sum_{j=1}^n u_{kl}^{[j]}~~(1\leq m\leq n)$.  
~On prendra, pour chaque valeur de \,$m$,  
$~l=\einf{ m/p }~$ et $~k=m+1-lp$;  
\item[5.~] Calculer \,$S^{-1}$\,  (en utilisant  
l'approche \gui{diviser pour gagner});  
\item[6.~] Calculer le produit  
\,$S^{-1}\vec{s}=\vec{c}$.  
\end{itemize}  
\end{algor}

\subsubsection* { La complexit\'e de l'\algo}  
      
Nous utilisons comme d'habitude les notations
\vref{IDNConst}.
Nous allons d\'eterminer  
les \parats de \com de la famille  
de \caris \parals repr\'esentant  
l'\algo de \PrSa en commen\c{c}ant  
par la \com des sous-\pcds qu'il  
utilise.  
  
\ss  
Les \parats de \com pour l'\algo  
principal \ref{algoPrSa}  (re\-pr\'esent\'e par la  
colonne PS$(A,n)\,$) et les \pcds  
auxiliaires \textsf{Powers} (colonne PW$(A,p)\,$) et 
\textsf{Superpowers} 
(colonne SPW$(A,p)\,$) seront d\'esi\-gn\'es,  
conform\'ement au tableau suivant, respectivement par:  
  
\ms\begin{center}  
\begin{tabular}{|l|c|c|c|}  
\hline $\downarrow$ Param\`etre\,/ Proc\'edure  
$\rightarrow$ & PS$(A,n)$\, & PW$(A,p)$\, &  
SPW$(A,p)$\, \\  
\hline ~Taille~ & \,$\tau(n)$\, &  
\,$\tau_1(p)$\, & \,$\tau_2(p)$\, \\  
\hline ~Profondeur & \,$\pi(n)$\, &  
\,$\pi_1(p)$\, & \,$\pi_2(p)$\, \\ 
\hline ~Largeur & \,$\lambda(n)$\, &  
\,$\lambda_1(p)$\, & \,$\lambda_2(p)$\, \\  
\hline  
\end{tabular}  
\end{center} 
    
\ms  
La d\'efinition de la \pcd \textsf{Powers}  
nous donne les relations de \recuz:  
  
\ms  
$\left\{\begin{array}{lll}  
\,\tau_1(p)\, &  =  & \,\tau_1(\esup{ p/2 })  
+ \einf{ p/2 } \muM(n)\, \\  
\,\pi_1(p)\, &  =  & \,\pi_1(\esup{ p/2 }) +  
\gaM (n)\, \\  
\,\lambda_1(p)\, & = & \,\max{\{\,\lambda_1(\esup{ p/2 }),  
\einf{ p/2 } \laM(n)\,\}}  
\end{array}\right.$

\mni on en d\'eduit pour \,$p\geq 2$\, par sommation de $1$  
\`a \,$r=\esup{ \log{p} }\,$:  

\ms  
$\left\{\begin{array}{lll}
\,\tau_1(p)\, & \leq & \,(2p-3)\,\Ca\,n^\alpha \\  
\,\pi_1(p)\, & \leq & \,\Ka\,\esup{ \log{p} }  
\log{n} \\  
\,\lambda_1(p) & = & \einf{ p/2 } \laM(n) \leq  
{1\over 2}\,\La\,p\,\log{n}\,.  
\end{array}\right.$  
  
\ms  
La d\'efinition de la \pcd \textsf{Superpowers}
dans laquelle  
\,$r=\esup{ \log{p} }$,  
\,$s=\einf{ p/r }$\, et \,$q:= p-rs$, permet  
d'\'ecrire\footnote{~Le calcul  
pr\'eliminaire des entiers \,$r,\,s,\,q$\,
n'intervient pas: il fait partie de la construction
du \cari correspondant.}:   
  
\ms  
$\left\{\begin{array}{lll}  
\,\tau_2(p)\, & = & \,\tau_1(s) + (p-s)\,\muM(n) \\  
\,\pi_2(p)\, & = & \,\pi_1(s) + r\,\gaM(n) \\ 
\,\lambda_2(p)\, & = & \,\max{\{\lambda_1(s), s\,\laM(n)\}}  
\end{array}\right.$  

\mni  
qui donnent, avec les majorations pr\'ec\'edentes: 
  
\ms
$\left\{\begin{array}{llll}  
\,\tau_2(p)\, & \leq & \,(p+s-3)\,\Ca\,n^\alpha  
& \leq  [\,p+(p/\log{p})-3\,]\,\Ca\,n^\alpha \\  
\,\pi_2(p)\, & \leq & \,2\Ka\,(\esup{ \log{p} }  
\log{n})\, &  \\ 
\,\lambda_2(p)\, & = & \,\einf{ p/r } \laM(n)  
& \leq (\La\,p\,n^\alpha)\,/\,(\log{p}\log{n})\,.  
\end{array}\right.$ \\  
  
\ms 
L'\algo utilise en plus des \pcds  
ci-dessus une \pcd d'inversion de  
\matgz.   
Nous avons vu (proposition \ref{divtri}) 
que l'inversion d'une \matg 
\freg se fait par  
un \cari \paral de taille major\'ee par  
\,$\Ca\,(2n-1)^\alpha$\, donc par  
\,$8\,\Ca\,n^\alpha$, de \prof au plus  
\'egale \`a $$\Ka\,(\log^2{(n)}+3\log{(n)}+2)  
+1\,\leq \,2\,\Ka\,(\log{(n)}+1)^2\,.$$ 
Et sa  
largeur est \,$\O(n^\alpha/\log^2{n})$\,  
si on applique le principe de Brent.  
  
\ms   
Ceci permet d'\'etablir la \com de  
la premi\`ere partie de l'\algo principal.   
Compte tenu du fait que  
\,$p=\esup{ \sqrt{n}\,}$\,  
et que \,$2\leq \alpha \leq 3$, le tableau  
\vref{PSW}   
indique le r\'esultat des majorations pour  
la taille et la \prof et pour chaque \'etape. 
  
\medskip 
{\begin{minipage}[t]{10,6cm}  
{  
\centerline{{\bf Complexit\'e de l'Algorithme de  
\PrSa}}\label{PSW}
\vspace*{2mm}   
\begin{center}
\begin{tabular}{|l|c|c|
    } 
\hline  
Etapes & Taille &  
Profondeur 
\\  
\hline  
Etape 1 & $[\,p+(p/\log{p})-3\,]\,n^\alpha$ &  
$2\,\Ka\,\esup{ \log{p} }\,\log{n}$ 
\\  
\hline  
Etape 2  & $[\,p+(p/\log{p})-3\,]\,n^\alpha$ &  
$2\,\Ka\,\esup{ \log{p} }\,\log{n}$ 
\\  
\hline 
Etape 3 & \,$n\,[p\muM(p) + (p-1)p^2]$\,  
& \,$\Ka\,\log{n}$\,  
\\  
\hline  
Etape 4 & \,$n\,(n-1)$\, &  
\,$\esup{ \log{n} }$\,  
\\  
\hline  
Etape 5 & \,$8\,\Ca\,n^\alpha$\, &  
$2\,\Ka\,(\log{(n)}+1)^2$  
\\  
\hline  
Etape 6 & \,$n^2$\, &  
\,$\esup{ \log{n} } + 1$\,  
\\  
\hline  
Total & ~$\,\tau(n)=\O(n^{\alpha+{1\over 2}})$\, &  
~$\,\pi(n)=\O(\log^2{n})\,$~  
\\  
\hline  
\end{tabular}   
\end{center}

\noindent \centerline{
{\bf Tableau \ref{PSW}}
}     
}\end{minipage}}
  
\bs On en d\'eduit le r\'esultat  
suivant de \PrSaz, dans lequel nous avons  
\egmt int\'egr\'e, le calcul de l'adjointe et de  
l'inverse qui constitue la deuxi\`eme partie de cet \algo:  

\begin{theorem}
Soit \,$\A$\, un anneau v\'erifiant les hypoth\`eses
pour l'\algo de Le Verrier. 
Le \polcarz, le \deterz, l'ad\-jointe  
et l'inverse (s'il existe) d'une \maca  
\,$A\in\A^{n\times n}$\, se fait par un \cari  
de taille \,$\tau(n)$, de \prof \,$\pi(n)$\, et de  
largeur \,$\lambda(n)$\, major\'ees respectivement  
par: $$\left \{ \begin{array}{lll}  
\tau(n) & \leq & 4\,\Ca\,n^{\alpha+{1\over 2}}  
+ \o(n^{\alpha+{1\over 2}}) \\  
\pi(n) & \leq & 5\,\Ka\,\log^2{n} + \O(\log{n}) \\  
\lambda(n) & \leq &  
(2\,\La\,n^{\alpha+{1\over 2}})\,/\,(\log^2{n}) \\  
\end{array}\right. $$  
o\`u \,$\Ca,\,\Ka,\,\La$\, d\'esignent  
les constantes asymptotiques de la  
\mul \paral des matrices en  
\,$\SD(n^\alpha,\log{n})$.
\end{theorem}  

\subsubsection*{Calcul de l'adjointe et de l'inverse}

L'\algo de \PrSa ne calcule pas toutes les  
puissances de la matrice \,$A$. Par cons\'equent le  
calcul de l'adjointe de \,$A$\, \`a partir de la formule  
de Cayley-Hamilton doit se faire en n'utilisant  
que les \,$2\,\esup{ \sqrt{n}\,}$\, puissances de  
\,$A$\, d\'ej\`a calcul\'ees, avec en plus les  
\coes \,$c_1,c_2,\ldots,c_n$\, du \polcar et les  
matrices \,$L_j$\, form\'ees des lignes des premi\`eres  
puissances de \,$A$\, \egmt disponibles.  
  
L'astuce est de consid\'erer les \,$p$\, matrices  

\ss
\centerline{$~B_{i-1}=\sum_{j=0}^{p-1}c_{n-p\,(i-1)-j-1}\,A^j~$}

\sni   
($1\leq i\leq p$) form\'ees avec les \coes du  
\polcarz, avec la convention \,$c_0=-1$\, et  
\,$c_k=0$\, si \,$k<0$\, (rappelons que \,$n\leq p^2\,$).  
On calcule ensuite la somme  

\ss \centerline{$~\sum_{k=0}^{p-1}B_{k}\,A^{pk}=\sum_{k=0}^{p-1}\,  
\sum_{j=0}^{p-1}c_{n-pk-j-1}\,A^{pk+j}~$}

\sni en r\'epartissant les calculs sur  
\,$\esup{ \log{p} }$\, \'etapes \parals avec au  
maximum \,$p\,/\,\log{p}$\, \muls de matrices  
\,$n\times n$\, (\ie des produits du type  
\,$B_{k}\times A^{pk}\,$) par \'etape.  
  
Or cette somme est \'egale \`a 
\,$\sum_{\ell=1}^{n}\,c_{n-\ell}\,A^{\ell-1} = \Adj A\,,$  
puisque d'une part \,$c_{n-\ell}=0$\, si \,$\ell>n$\, et que   
d'autre part, si \,$\ell$\, est compris entre 1 et \,$n$,  
\,$\ell$\, correspond de mani\`ere unique \`a un couple  
\,$(k,j)$\, tel que \,$1\leq j,\,k \leq p-1$\, et  
\,$\ell-1=p\,k+j$\, (division euclidienne de \,$\ell-1$\,  
par \,$p$). Ce qui donne l'adjointe puis l'inverse.  
  
\ss Ainsi la deuxi\`eme partie de l'\algo de \PrSa  
pour le calcul de l'adjointe et de l'inverse de  
\,$A$\, peut \^etre d\'etaill\'ee comme suit.  
  
\mni
\textsf{\textbf{Entr\'ees}: \\[1mm]  
--- Les puissances \,$A,\,\ldots,\,A^{p}$\, de la  
matrice \,$A$, ainsi que les puissances  
\,$A^{2p},\,\ldots,\,A^{p^2}$\, de la  
matrice \,$A^p$, toutes disponibles \`a l'issue des  
deux pre\-mi\`eres \'etapes de l'\algo principal  
\ref{algoPrSa};\\[0,5mm]   
--- La matrice \,$L=[\,L_1\,|\,L_2\,|\,\cdots\,|\,L_n\,]\,  
\in \A^{p\times n^2}$\, form\'ee des \,$n$\, matrices  
\,$L_k$\, ($\,1\leq k\leq n\,$) d\'ej\`a calcul\'ees; \\[0,5mm]      
--- Enfin la matrice \,$C$\, form\'ee \ˆ partir des  
coefficients \,$c_1,c_2,\ldots,c_n$\, du \polcar  
$\rPA(X)= (-1)^n\,\left(X^n - \sum_{i=1}^nc_iX^{n-i}\right)$:
$$ 
C = \left[\begin{array}{cccc}  
c_{n-1} & c_{n-2} & \cdots & c_{n-p}  \\  
c_{n-p-1} & c_{n-p-2} & \cdots & c_{n-p-p}  \\  
\vdots & \vdots & \vdots & \vdots \\   
c_{n-p\,(p-1)-1} & c_{n-p\,(p-1)-2} & \cdots &  
c_{n-p\,(p-1)-p}  
\end{array}\right]~\in\A^{p\times p}\,\,.
$$  
On a alors $~B_{i-1}= \sum_{j=0}^{p-1} c_{ij}\,A^j~$ et  
il est facile de voir que la \,$k\,$-\,\`eme ligne de cette  
matrice n'est autre que la \,$i\,$-\,\`eme ligne de la matrice  
\,$T_k:=C\,L_k$\, o\`u \,$L_k$, rappelons-le, est  
la matrice form\'ee des \,$k\,$-\,\`emes lignes des matrices  
\,$\In,\,A,\,\ldots,\,A^{p-1}$.\\[2mm]
\textbf{Sortie}: \\[1mm]   
L'adjointe et l'inverse de \,$A$, \cad les matrices \\  
\,$\Adj A = A^{n-1}-c_1A^{n-2}-\dots-c_{n-1}A-c_n\In~~$  
et $~~\,A^{-1} = {1\over c_n}\,\Adj A$.\\[2mm]
\textbf{Les \'etapes du calcul}:  \\[1mm]  
Faisant suite aux \'etapes (1 \`a 6) qui calculent  
le \polcarv elles seront num\'erot\'ees de 7 \`a 10.    
On pose \,$r = \esup{\log{p}}$\, et \,$D_0=0_{nn}$\,  
(la \maca d'ordre \,$n$\, nulle)  
et \,$s = \einf{ p/r }$:
\begin{itemize}  
\item[7.~] Calculer le produit  
\,$T=C\,L=[\,CL_1\,|\,CL_2\,|\,\cdots\,|\,CL_n\,]$\, 
(ce qui revient \`a calculer en \paral les produits  
de la matrice \,$C$\, qui est une matrice \,$p\times p$\,  
par les \,$n$\, matrices \,$CL_k$\, qui sont des matrices  
\,$p\times n\,$).\\ Cette \'etape permet d'\'ecrire les  
matrices \,$B_{i-1}$\, ($\,1 \leq i \leq p\,$).  
\item[8.~] \pour{k}{1}{r} $D_k := D_{k-1} +  
\sum_{k=(i-1)s}^{is-1} B_{k}\,A^{pk}$;  
\item[9.~] Calculer \,$\Adj\, A :=-(D_r + \sum_{k=rs}^{p-1}  
B_{k}\,A^{pk})$;  
\item[10.~] Calculer \,$A^{-1}={1\over{c_n}}\,\Adj\, A$.  
\end{itemize} 
} 

\ms La \com de cette deuxi\`eme partie de l'\algo  
de \PrSa poss\`ede les m\^emes bornes  
que l'\algo principal du \polcarz.

L'\'etape 7  se fait en \,$\Ka\,\log{n}$\, \'etapes  
comportant au total \,$n\,p\,\,\Ca\,p^\alpha$\, \ops  
\ariths dans l'anneau de base, utilisant au maximum  
\,$n\,p\,\,\La\,(p^\alpha\,/\log{p})$\, processeurs.

Les \'etapes 8 et 9 sont les plus co\^uteuses. Elles  
correspondent \`a un total de \,$r+1$\, \'etapes  
\parals comportant \,$p$\, \muls de \macas
d'ordre \,$n$\, \cad \,$p\,\,\Ca\,n^\alpha$\,  
\ops \ariths de base, ce \`a quoi il
faut rajouter des additions de  
matrices \,$n\times n$.  
Cela fait un \cir de \prof  
$$(r+1)\,\La\,\log{n} + \O(\log{n})\leq  
{1\over 2}\,(\log{n}+1)\,\La\,\log{n} +  
\O(\log{n})\,.$$  
Le nombre de processeurs utilis\'es au  
cours de ces \,$r+1$\, \'etapes \parals est \'egal  
\`a \,$s\,\La\,(n\alpha\,/\log{n}) \leq  
(\La\,p\,n\alpha)\,/\log{n}\log{p}$\, puisque  
\,$s = \einf{ p/r } \leq p\,/\log{p}$.

\section{Am\'elioration de Galil et Pan} \label{sec GaliPan} 
  
\GaPa \cite{GalPan} r\'eduisent les \'etapes  
les plus co\^uteuses de l'\algo pr\'ec\'edent \`a quatre  
\muls de \marcsz.  
  
Il s'agit plus pr\'ecis\'ement des l'\'etapes 1, 2 et 3 de  
l'\algo principal (calcul du \polcarz) d'une part et des \'etapes 8 et  
9  du calcul de l'adjointe d'autre part.  

Par une r\'eorganisation des \'etapes 1 et 2 de l'\algo  
principal qui font intervenir les \pcds  
r\'ecursives \textsf{Powers} et \textsf{Superpowers}, on remplace  
l'appel \`a ces \pcds par l'appel r\'ecursif \`a  
une \pcd unique permettant de calculer  
les matrices 

\ss \centerline{\,$[\,A\,|\,A^2\,|\,\ldots\,|\,A^{p-1}\,]~$ et  
$~[\,A^p\,|\,A^{2p}\,|\,\ldots\,|\,A^{p(p-1)}\,]$}
  
\sni \`a partir des matrices 

\ss \centerline{$[\,A\,|\,A^2\,|\,\ldots\,|\,A^{s-1}\,]~$ et  
$~[\,A^{s}\,|\,A^{2s}\,|\,\ldots\,|\,A^{(s-1)\,s}\,]$}
  
\sni  
o\`u \,$s= \esup{\sqrt{p}\,}$.  
Cela se fait en effectuant le produit d'une \marc
 \,$n\,(s-1)\times n$\, par une \marc
 \,$n\times ns$\, qui donne les puissances  
restantes:
$$
\left[\begin{array}{c}  
A\\A^2\\ \vdots \\ A^{s-1}\end{array}\right]~\times~  
\left[\begin{array}{c|c|c|c}  
A^s\, & \,A^{2s}\, & \,\cdots\, & \,A^{s^2}  
\end{array}\right]\,.
$$  
  
L'\'etape 3 de l'\algo principal calcule les \,$n$\,  
produits \,$U_j=L_j\,C_j \in \A^{p\times p}$\, pour  
en d\'eduire les traces des puissances de \,$A$.  
Il est possible de r\'eduire cette \'etape au calcul  
d'un seul produit de deux \marcs de  
types respectifs \,$p\times n^2$\, et \,$n^2\times p$\,  
o\`u \,$p=\esup{ \sqrt{n}\, }$.  
En effet, si on \'ecrit les \elts de chaque  
matrice \,$A^k$\, (pour \,$0\leq k\leq p-1\,$) sur une  
seule ligne, de mani\`ere \`a la repr\'esenter par la  
suite ordonn\'ee de ses \,$n$\, lignes, \cad par  
\,$v_{k} \in \A^{1\times n^2}$, et si l'on fait de  
m\^eme avec les matrices \,$A^{pk}$\,  
($\,1\leq k\leq n\,$), mais en d\'eroulant cette fois  
chacune d'elles sur une seule colonne ($\,A^{pk}$\, sera  
donc repr\'esent\'ee, dans l'ordre de ses colonnes,  
par \,$w_k\in \A^{n^2\times 1}\,$), le calcul des traces  
\,$s_k$\, revient alors \`a calculer  
le produit des deux \marcsz:  
$$\left[\begin{array}{c}  
v_0 \\ v_1 \\ \cdots \\ v_{p-1}  
\end{array}\right]~\times~  
\left[\begin{array}{c|c|c|c}  
w_1\, & \,w_2\, & \,\cdots\, & \,w_p\,  
\end{array}\right] =$$ 
$$\left[\begin{array}{cccc}  
v_0\,w_1 & v_0\,w_2 & \,\cdots\, & v_0\,w_p \\  
v_1\,w_1 & v_1\,w_2 & \,\cdots\, & v_1\,w_p \\  
\vdots & \vdots & \vdots & \vdots \\  
v_{p-1}\,w_1 & v_{p-1}\,w_2 & \,\cdots\, & v_{p-1}\,w_p
\end{array}\right]\,.$$ 
  
Il est clair que l'\elt \,$v_{i-1}\,w_j$\, de la  
\,$i\,$-\,\`eme ligne et \,$j\,$-\,\`eme colonne de cette matrice est  
\'egal \`a \,$s_{pj+i-1}=\Tr A^{pj+i-1}~~(1\leq i,j\leq p)$.  
  
\ss  
On modifie enfin les \'etapes 8 et 9 du  
calcul de l'adjointe de \,$A$\, en prenant  
\,$q=\esup{ \sqrt[3]{n}\,  }$, \,$t=\einf{ (n+1)\,/\,  
q }$\, de mani\`ere \`a avoir \,$qt \leq n+1 <  
q(t+1)$, et on change les dimensions de la matrice  
\,$C$\, en la rempla\c{c}ant par \,$C^*=(c_{ij})\in\A^{(t+1)  
\times q}$\, (avec les m\^emes notations et la m\^eme  
convention pour les \,$c_{ij}\,$) ainsi que les dimensions  
des matrices \,$L_1,\,L_2,\,\ldots,\,L_n$\, en les  
rempla\c{c}ant par des matrices  
\,$L_1^*,\,L_2^*,\,\cdots,\,L_n^*$\,  
d\'efinies exactement de la m\^eme fa\c{c}on mais \`a partir  
des lignes des matrices  
\,$\In,\,A,\,\cdots,\,A^{q-1}$, ce qui  
fait qu'elles sont de type \,$q\times n$\, au lieu d'\^etre  
de type \,$p\times n$.

On calcule alors la matrice \,$T^*\in \A^{(t+1)\times n^2}$\,  
en effectuant le produit d'une matrice  
\,$(t+1)\times q$\, par une matrice \,$q\times n^2\,$:
$$T^*=C^*\,[\,L_1^*\,|\,L_2^*\,|\,\cdots\,|\,L_n^*\,] =  
[\,C^*L_1^*\,|\,C^*L_2^*\,|\,\cdots\,|\,C^*L_n^*\,]
$$
en tenant compte du fait que la \,$(i+1)\,$-\,\`eme ligne du  
bloc \,$C^*L_j^*$\, n'est autre que la \,$j\,$-\,\`eme ligne  
de la matrice 
$$B_i = \somm_{j=0}^{q-1}c_{n-1-qi-j}\,A^{j}  
= \somm_{j=0}^{q-1}c_{i+1,j+1}A^{qi+j}~~~~(\mathrm{ici}~~  
0\leq i\leq t)\,.$$
  
\ss  
Avec ces modifications, les \'etapes 8 et 9 se ram\`enent  
donc, comme on peut le constater, au calcul du produit  
de deux \marcs (avec les m\^emes  
notations que ci-dessus) qui est un produit d'une matrice  
\,$n\times n\,(t+1)$\, par une matrice  
\,$n\,(t+1)\times n\,$:
$$ \left[\begin{array}{c|c|c|c|c}  
B_0 & B_1 & B_2 & \cdots & B_{t} \end{array}\right]~  
\times~\left[\begin{array}{c} \In\, \\ A^q \\ A^{2q}  
\\ \vdots \\ A^{tq} \end{array}\right]\,.
$$

Posant \,$\ell=qi+j$, ce dernier produit est en effet  
\'egal \`a 
$$\somm_{i=0}^t B_i\,A^{qi}=\somm_{i=0}^{t}  
\,\somm_{j=0}^{q-1}c_{n-1-qi-j}\,A^{qi+j} =  
\somm_{\ell=0}^{(t+1)\,q-1}\,c_{n-1-\ell}\,A^{\ell}\,.
$$  
Comme \,$(t+1)\,q-1 > n$\, (d'apr\`es la d\'efinition  
m\^eme de \,$t$\, et de \,$q\,$) et que \,$c_{n-1-\ell}=0$\,  
pour \,$\ell\geq n$, la matrice ainsi obtenue est  
exactement l'oppos\'ee de l'adjointe de \,$A$:  
$\Adj{A} = -\sum_{l=1}^{n}\,c_{n-\ell}\,A^{\ell-1}.$  
  
Les calculs de ces quatre produits de \marcs 
auxquels \GaPa r\'eduisent l'\algo  
de \PrSaz, et qui sont des \muls  
d'ordres respectifs donn\'es par le tableau suivant o\`u  
\,$p=\esup{\sqrt{n}\,},\,s=\esup{\sqrt{p}\,}\,  
,\,q=\esup{ \sqrt[3]{n}\,},\,t=\einf{ (n+1)\,/\,  q }$,

\vspace*{1mm}
\begin{center}  
\begin{tabular}{|l|c|c|}  
\hline  
Multiplication & 1er facteur & 2\`eme facteur \\  
\hline 1\`ere \mul & \,$n\,(s-1)\times n$\, &  
\,$n\times ns$\, \\  
\hline 2\`eme \mul & \,$p\times n^2$\,  &  
\,$n^{2}\times p$\, \\ 
\hline 3\`eme \mul & \,$(t+1)\times q$\, &  
\,$q\times n^2$\, \\  
\hline 4\`eme \mul & \,$n\times n\,(t+1)$\, &  
\,$n\,(t+1)\times n$\, \\  
\hline  
\end{tabular}  
\end{center} 
  
\sni  
s'effectuent en \,$\O(\log^2{n})$\, \'etapes \paralsz.  
  
\ss  
On fait d'autre part appel aux r\'esultats concernant  
les notions d'\algo bi\lin et de \rgte  
(voir la section \ref{sec combili}), 
pour am\'e\-lio\-rer la \com th\'eorique  
de l'\algo de \PrSa ainsi remani\'e, en  
faisant passer l'exposant de \,$n$\, dans cette  
\com (en taille et en nombre de processeurs) de  
$2,876$ \ˆ $2,851$ (si on prend le \,$\alpha \approx  
2,376$\, de Winograd \& Coppersmith \cite{Win}).  
  
\ss  
Rappelons (voir la section \ref{paraRangTens}) 
que le \rgte de l'\abi
$$f\,:\, \A^{m\times n}\times  
\A^{n\times p} \dans \A^{m\times p}
$$ 
associ\'ee \`a la \mul des matrices \,$m\times n$\, par des  
matrices \,$n\times p$\, \`a coefficients dans \,$\A$\,  
(on note \,$\langle m,n,p \rangle_\A$\, cette
\abiz) est d\'efini  
comme le rang de l'\algo \bil ou du tenseur  
d\'efinissant \,$\gen{m,n,p}_\A$, 
\cad le nombre minimum de  
\muls essentielles \ncrs au \cabi correspondant.  
Ce rang est not\'e  
\,$R\,\gen{m,n,p}\,$
(nous omettons \,$\A$\, en indice dans la mesure o\`u
tous les r\'esultats cit\'es s'appliquent \`a n'importe quel anneau).

Outre les propri\'et\'es \'etablies dans la section \ref{sec combili},
il y a un r\'esultat d\^u \`a Coppersmith \cite{Cop} 
pour le cas des \marcs qui nous occupe ici. 
Il est utilis\'e par \GaPa pour \'etablir 
 qu'il existe une constante positive  
$\beta$\, estim\'ee dans un premier temps \`a  
\,$\beta=(2\log{2})\,/\,(5\log{5})\approx 0,172$\, puis  
\`a \,$\beta\approx 0,197$\, qui v\'erifie la propri\'et\'e  
$$\,R\,\langle m,m^\beta,m \rangle = \O(m^{2+\epsilon})~~~  
\mbox{ pour tout }~~~\epsilon > 0\,.$$  
  
Les modifications des \'etapes les plus co\^uteuses  
aboutissent \`a des \muls de \marcs  
de rangs respectifs:

\vspace*{1mm}
\begin{center}  
\begin{tabular}{|l|l|}  
\hline  
Multiplication & Rang tensoriel \\  
\hline 1\`ere \mul &  
\,$R\,\langle n^{5/4},n,n^{5/4}\rangle$\, \\  
\hline 2\`eme \mul &  
\,$R\,\langle n^{1/2},n^2,n^{1/2}\rangle$\, \\  
\hline 3\`eme \mul &  
\,$R\,\langle t+1,q,n^2\rangle$\, \\  
\hline 4\`eme \mul &  
\,$R\,\langle n,n\,(t+1),n \rangle$\, \\  
\hline  
\end{tabular}  
\end{center}

\sni o\`u \,$q\asymp n^{1/3}$\, et \,$t\asymp n^{2/3}$\,  
v\'erifient aussi \,$qt \leq n+1 < q\,(t+1)$.

On a alors:  
\begin{theorem}  \emph{(Galil \& Pan)}\label{thGaPa}\\
Le calcul du \polcarz, de l'adjointe et l'inverse d'une \maca
 d'ordre \,$n$\, 
 est en \,$\SD(n^{\alpha+{1\over 2}-\delta},log^2{n})$\,
o\`u \,$\delta$\, est un r\'eel strictement positif  
d\'ependant de \,$\alpha$. \\  
En particulier, pour \,$\alpha \approx 2,376$\, la taille  
du  \cari est un \,$\O(n^{2,851})$.  
\end{theorem}

Il suffit en effet, pour \'etablir ce r\'esultat,  
d'\'evaluer les quatre rangs tensoriels indiqu\'es dans le tableau  
ci-dessus en utilisant la constante \,$\beta$\,  
de la \mul des \marcs  
($\,\beta<1\,$).   
Pour cela, on pose \,$m=n^{1\,/\,(4-4\beta)}$\, et  
\,$r=n^{1\,/\,(4-4\beta)}$, ce qui donne les  
estimations 

\ss \centerline{$R\,\langle m,m^\beta,m\rangle =  
\O(n^{(2+\epsilon)\,/\,(4-4\beta)})$} 

\sni et 

\ss \centerline{$R\,\langle r,r,r \rangle =  
\O(n^{\alpha(4-5\beta)\,/\,(4-4\beta)})$} 

\sni   
qui, multipli\'ees entre elles, donnent  
$$R\,\langle n^{5/4},n,n^{5/4}\rangle =  
\O(n^\rho)~~~\mbox{ o\`u }~~~\rho = \alpha +  
\frac{1}{2} + \delta_1 ~~~\mbox{ et }~~~   
\delta_1 = \frac{\epsilon -  
\beta(\alpha - 2)}{4-4\beta}\,\cdot$$  
Comme \,$\alpha > 2$\, et \,$\beta <1$, on peut prendre  
\,$0< \epsilon < \beta(\alpha - 2)$\, et \,$\delta_1 >0$,  
ce qui \'etablit le r\'esultat pour la premi\`ere  
\mulz.  
  
\ss  
Pour les trois autres \mulsz, on remarque que:  
  
\noi $\bullet$~ d'une mani\`ere \gnle \,$R\,\langle  
m,m^4,m\rangle = \O(m^{\alpha + 3})$\, (\mul par  
blocs \,$m\times m\,$) et que, par cons\'equent,  
\,$R\,\langle n^{1/2},n^2,n^{1/2}\rangle =  
\O(n^{\frac{\alpha + 3}{2}}) = \O(n^{\alpha + \frac{1}{2} +  
\delta_2})$\, avec \,$\delta_2=\frac{\alpha-2 }{2} > 0$.\\  
$\bullet$~ \,$R\,\langle t+1,q,n^2\rangle =  
\O(n^{(\alpha-3)\,\eta+3}) = \O(n^{\alpha + \frac{1}{2} +  
\delta_3})$\, si l'on prend  
\,$t=n^\eta$, \,$q=n^{1-\eta}$, et \,$\delta_3 =  
(1-\eta)(3-\alpha) +{1\over 2}$\, avec \,$0<\eta<1$.\\  
$\bullet$~ \,$R\,\langle n,n\,(t+1),n \rangle =  
\O(n^{\alpha +\eta}) = \O(n^{\alpha +\frac{1}{2}+\delta_4})$\,  
avec \,$\delta_4 = \frac{1}{2} -\eta$\, pour le m\^eme  
\,$\eta$.  
  
\ss  
Prenant \,$0<\eta<{1\over 2}$\, et  
\,$(1-\eta)(\alpha-3)<{1\over 2}$, ce qui correspond au  
cas concret \,$\eta = {1\over 3}$, cela donne bien  
\,$\inf{(\delta_1,\delta_2,\delta_3,\delta_4)}>0$\, et  
\'etablit le r\'esultat \,$\o(n^{\alpha+1/2-\delta})$\,  
pour n'importe \,$\delta >0$\, strictement inf\'erieur \`a  
\,$\inf{(\delta_1,\delta_2,\delta_3,\delta_4)}$. 
Le r\'esultat num\'erique en d\'ecoule pour  
\,$\alpha<2,376$.
 
En fin de compte l'exposant de \,$n$\, dans la  
\com asymptotique pour le calcul du \polcar et de  
l'adjointe par la \met de \PrSa est  
de $2,876$ au lieu de $2,851$ de \GaPa pour 
\,$\alpha=2,376$.  

\subsubsection*{Conclusion}  
Les \algos de Csanky, de \PrSaz, de \GaPa ne sont en 
fait que des variantes \parases  
de la \met de Le Verrier (1840) mais elles ont le  
m\'erite d'avoir ing\'enieusement r\'eduit, et  
de mani\`ere spectaculaire, la \com des  
\caris permettant de r\'esoudre ces probl\`emes dans le  
cas d'un anneau commutatif autorisant les divisions exactes  
par les entiers. Les estimations de ces \algos \parals dans  
le cas de tels anneaux restent les meilleures connues \`a  
l'heure actuelle.  
  
\acvide
\newpage \thispagestyle{empty}

 
\chapter[\Polcar sur un anneau arbitraire]{Calcul du \polcar sur un
\acoma}
\minitoc
\label{chap PolCarAnn}

\subsubsection*{Introduction}
 
Dans ce chapitre, nous pr\'esentons des \algos  bien
\parass de calcul du \polcar sur un \acomaz.
 
Le premier r\'esultat de cette sorte, expos\'e dans la section
\ref{sec vonzur}, a \'et\'e obtenu en 1982.
L'estimation de son  temps \sql est pessimiste, mais il reste
d'un grand int\'er\^et th\'eorique
 
Dans les sections suivantes
nous expliquons les \algos de Chistov et de Berkowitz (am\'elior\'e)
qui sont dans \,$\SD(n^{\alpha +1}\log{n},\log^2{n})$.
 
On notera que le r\'esultat est cependant moins bon en temps
\sql que pour l'\algo de \PrSa
(qui r\'eclame la division par un entier arbitraire)
ou celui de Kaltofen (qui n'est pas bien \parasz).
 
\section{M\'ethode \gnle de parall\'elisation}
\label{sec vonzur}
 
Tout   \prev (donc tout \cariz) sans division \`a \,$n$\,
\idtrs \,$(x_i)_{i=1..n}$\, sur un anneau
\,$\A$\, calcule un \poly de \,$\A\,[x_1,\ldots ,x_n]$.
 Valiant, Skyum, Berkowitz et Rackoff
\cite{Vasbra} d\'emontrent le r\'esultat important suivant.
La preuve, d\'elicate, est bien expliqu\'ee dans \cite{Bur2}.
 
\begin{theorem}\label{vaskybra}
Soit \,$\Gamma$\, un \cari sans division,
de taille \,$\ell$, qui calcule un \poly \,$f$\, de
degr\'e \,$d$\, en \,$n$\, variables
sur un anneau $\A$. Alors il existe
un \carih \,$\Gamma'$\, de taille
\,$\O(\ell^3d^6)$\, et de \prof
\,$\O(\log({\ell\,d})\log{d})$\, qui calcule
(les composantes \hogs de) \,$f$\,{\rm (}{\footnote{~On trouve dans
\cite{Bur2} la majoration \,$\O(\log({\ell\,d})\log{d}+\log{n})$\,
pour la \profz.
Le terme \,$\log({n})$\, \supt est \ncr lorsque
\,$d=1$\, si on a \,$\log{1}=0$. Mais la convention de notation
\ref{notaLog} que nous avons choisie pour \,$\log{d}$, conforme \`a la
longueur du code binaire de \,$d$, nous donne \,$\log{1}=1$.}}{\rm )}.
En outre la construction de \,$\Gamma'$\, \`a partir de \,$\Gamma$\, 
est $\LOSP$.
\end{theorem}
 
En particulier:
 
\begin{cor} Toute famille \,$(Q_{\ell})$\, de
\pols de degr\'es  \,$d=\O(\ell^k)$\, qui peut \^etre
calcul\'ee au moyen d'une \famu de
\caris peut aussi \^etre calcul\'ee dans la classe
$\NC^2$.
\end{cor}

En appliquant le \tho \ref{vaskybra}
\`a l'\apg auquel on fait subir la \pcd
d'\elid \`a la Strassen, et vu que le \deter qu'il
calcule est un \poly de degr\'e \,$n$, on obtient
le r\'esultat suivant d\^u \`a Borodin, Hopcroft et
Von zur Ga\-then~\cite{Bor}:
\begin{prop} \label{propBHV}
Le \deter d'une matrice \,$n\times n$\,  est calcul\'e
par un \prev de taille
\,$\O\left(n^{18}\log^3{n}\,\log^3{\log{n}}\right)$\, et
de \prof \,$\O(\log^2{n})$.
\end{prop}
 
Dans la construction correspondant au \tho
\ref{vaskybra} est utilis\'ee la \mul rapide des \polsz.
Avec la \mul \usle des \polsz, la proposition \ref{propBHV}
 donne \,$\O\left(n^{21}\right)$\, \oparis
dans l'anneau de base.

\section{Algorithme de Berkowitz am\'elior\'e}
\label{sec Berkopar}

\subsubsection*{Introduction}
 
Utilisant la \met de partitionnement \cite{Gas,Fad},
attribu\'ee \`a Sa\-muelson (\cite{Sam}),
Berkowitz \cite{Ber} a pu exhiber un \cari \paral de taille
\,$\O(n^{\alpha+1+\epsilon})$\, et de \prof
\,$\O(\log^{2}n)$, o\`u \,$\epsilon$\,
est un r\'eel positif quelconque.
 
\ss
Il a ainsi am\'elior\'e de mani\`ere
d\'ecisive
la \com asymptotique du calcul des \detersz,
\polcarsz, et adjointes de matrices
\`a \coes dans un anneau commutatif quelconque
\,$\A$.
 
\ss
Nous allons donner une version l\'eg\`erement am\'elior\'ee
de l'\agbz, due \`a Eberly \cite{Eber}, qui ram\`ene sa taille \`a
\,$\O(n^{\alpha+1}\log{n})$\, sans en changer la \profz.
Pour cela nous donnons une version plus simple de la \recu
utilis\'ee pour le calcul des \coes du \polcarz.
Nous donnons
\egmt une estimation pr\'ecise de la constante qui intervient
dans le \gui{$grand~\O$} de la \com
\sqle (cf. \cite{These}).
 
\ss
Soit \,$A=(a_{ij})\in\A^{n\times n}$\, une \maca d'ordre \,$n \geq 2$\,
sur un
\acoma  \,$\A$.
Conform\'ement aux notations introduites dans la section
\ref{sec AlgLin},  pour tout entier
\,$r~(1 \leq r \leq n)$, on d\'esigne par \,$A_r$\,
la \smpd d'ordre \,$r$\, de \,$A$.
On notera ici $\,R_r$\, la matrice
\,$A_{r+1,1..r}\in\A^{1 \times r}$\, et \,$S_r$\,
la matrice  $\,A_{1..r,r+1}\in\A^{r \times 1}$.
Rappe\-lons  la formule de Samuelson (\ref{EqSam}) vue \`a la section
\ref{sec berkoseq}.
$$
P_{r+1}= \left\{
\begin{array}{lll}
(a_{r+1,r+1} - X)\, P_r(X)\, +\\[1mm]
\sum_{k=2}^{r+1} \;\left[(R_rA_r^{k-2}S_r)\, p_0 + \cdots
+ (R_r S_r)\, p_{k-2}\right]\; X^{r+1-k}
\end{array}
\right.
$$
o\`u \,$P_r(X)=\sum_{i=0}^r p_{r-i}X^i$.
Notons \,$Q_{r+1}$\, le \pol
$$
-X^{r+1}+a_{r+1,r+1}\,X^r+R_rS_r\,X^{r-1}+
R_rA_rS_r\,X^{r-2}+\cdots+R_rA_r^{r-1}S_r\,.
$$
On peut aussi \'ecrire  la formule de Samuelson
sous la forme (\ref{sam}):
$$
\overrightarrow{P_{r+1}}=\Toep(Q_{r+1})\times
\overrightarrow{P_r}
$$
o\`u \,$\overrightarrow{P_r}$\, est le vecteur colonne
\,$\tra{(p_0,p_{1},\ldots ,p_r)}$\, des \coes du \poly
\,$P$\, et \,$\Toep(Q_{r+1})\in\A^{(r+2)\times(r+1)}$\,
est la \mto suivante d\'efinie \`a partir
du \pol \,$Q_{r+1}$:
$$
\Toep(Q_{r+1}) =
\left[\begin{array}{ccccc}
-1          & 0 & \cdots & \cdots & 0 \\[1mm]
a_{r+1,r+1} & -1 & \ddots &   & \vdots \\[1mm]
R_rS_r & \ddots & \ddots & \ddots & \vdots \\[2mm]
\vdots & \ddots & \ddots & \ddots & 0 \\[2mm]
R_rA_r^{r-2}S_r &  &  \ddots
 & \ddots & -1 \\[2mm]
R_rA_r^{r-1}S_r &  R_rA_r^{r-2}S_r & \dots & R_rS_r & a_{r+1,r+1}
\end{array}\right]$$
 
Le calcul du \polcar consiste donc:
\begin{itemize}
\item
\`a calculer d'abord les \coes de la matrice
\,$\mathrm{Toep}(Q_{r+1})$\, -- qui interviennent dans
l'\egtz~(\ref{sam}) -- ou, ce qui revient au m\^eme,
la famille \,$T=\{RM^iS\}_{i=0}^{r-1}$\, lorsque
\,$R,\,M,\,S$\, sont respectivement des matrices
\,$1 \times r$, \,$r \times r$, et
\,$r\times 1$, et lorsque r est un entier tel que
\,$2 \leq r < n$\, ($M=A_r$, \,$R=R_r$, \,$S=S_r$);
\item
\`a calculer ensuite le \pol \,$P_n$\, dont le
vecteur des \coesz, compte tenu de (\ref{sam}),
est donn\'e par:
\begin{equation}\label{top}
\overrightarrow{P_n}=
\mathrm{Toep}(Q_n) \times \mathrm{Toep}(Q_{n-1})\times\cdots
\times \mathrm{Toep}(Q_1)
\end{equation}
\end{itemize}
Dans son papier original \cite{Ber}, Berkowitz d\'emontre que les
familles $$U=\{RM^i\}_{i=0}^{n^{1/2}}~~\mathrm{et}
~~V=\{M^{j\,n^{1/2}}S\}_{j=0}^{n^{1/2}}$$ peuvent \^etre
calcul\'ees par un \cari \paral en \,$\SD(n^{\alpha +
\epsilon},\alb\log^2{n})$\, pour en d\'eduire que le calcul du \polcar
se fait en \,$\SD(n^{\alpha+1+\epsilon}, \log^2{n})$\,.
 
\subsection*{La version parall\`ele am\'elior\'ee
et sa complexit\'e}
 
Nous utilisons comme d'habitude la notation \vref{IDNConst}.

\ms 
\begin{prop} \label{RMiS} 
On consid\`ere un entier \,$r\geq 2\,$ et des matrices 
\,$R\in\A^{1\times r}$, \,$M\in\A^{r\times r}$, 
\,$S\in\A^{r\times 
1}$. 
$$
\mathit{La~famille}~~~~ T=\{R\,M^iS\}_{i=0}^{r-1}
$$ 
peut \^etre 
cal\-cul\'ee par un \cari dont la taille et la \prof sont 
major\'ees respectivement par 
$$\Ca\, r^{\alpha}\,\log{r}+\O(r^{\alpha}) 
~~~\mathrm{et}~~~\Ka\,\log^2{r}+\O(\log{r})\,.
$$ 
\end{prop} 
 
\preuve \\ Soit \,$r\geq2$. On utilisera, pour l'analyse de \com
des \algosz, les entiers \,$\nu= \esup{ \log_4r } = \esup{
\frac{1}{2}\,\log{r} }$\, et \,$\eta= \esup{ \log r }$\, qui
v\'erifient les in\'egalit\'es~: \,$2^{2\nu - 2}<r\leq 2^{2\nu}$\,
et \,$2^{\eta - 1}<r\leq 2^{\eta}$\, (on a aussi \,$1\leq \nu \leq
2\nu -1 \leq \eta\leq 2\nu$).
 
\noindent Toute \maca \,$A$\, d'ordre \,$r$\, sera plong\'ee,
selon le cas, soit dans une matrice \,$ \cmatrix{ A  & 0  \cr  0 &
0 }$\, carr\'ee d'ordre \,$2^{\eta}$\, soit dans une matrice \,$
\cmatrix{ A  & 0  \cr  0  & 0 }$\, carr\'ee d'ordre \,$2^{2\nu}$\,
(chacun des \,$0$\, d\'esignant ici une matrice nulle de
dimensions convenables).
 
\noindent
Il faut cependant remarquer que, dans les deux cas, l'\'el\'evation au
carr\'e de la matrice \,$
\cmatrix{ A  & 0  \cr  0  & 0 }$\, se fait \`a l'aide d'un \cari de
taille \,$\Ca r^{\alpha}$\, et de profondeur \,$\Ka \log{r}$\, puisque
\,$\cmatrix{ A  & 0  \cr  0  & 0 }^{2}=\cmatrix{ A^{2} & 0 \cr 0 & 0
}$\,.
 
\noindent De m\^eme, le produit d'une matrice \,$2^{k} \times
2^{2\nu}~ (k=1,\ldots ,\nu)$\, par une matrice du type
\,$\cmatrix{ A  & 0  \cr  0 & 0 } \in \A^{2^{2\nu} \times
2^{2\nu}}$\, avec \,$A\in \A^{r\times r}$\, peut \^etre obtenu par
le calcul du produit d'une matrice \,$2^{k} \times 2^{\eta}$\, par
une matrice du type \,$\cmatrix{ A  & 0  \cr  0 & 0 } \in
\A^{2^{\eta} \times 2^{\eta}}$\, \`a cause du fait que, dans ces
deux produits, les \,$r$\, premi\`eres colonnes sont les m\^emes
alors que les colonnes restantes sont nulles. Ce qui fait que le
produit en question peut \^etre obtenu par \,$2^{2(\eta - k)}$\,
multiplications en \paral de blocs \,$2^{k} \times 2^{k}$\, et de
\,$2^{\eta - k}\, (2^{\eta - k}-1)$\, additions en \paral des
blocs produits obtenus, \cad par un \cari de taille \,$2^{\eta -
k}\,[2^{\eta - k}\Ca 2^{k\alpha} + (2^{\eta -
k}-1)2^{2k}]$\,(\footnote{~Le premier terme du crochet provient
des \muls de blocs \,$2^{k} \times 2^{k}$\,, et le second terme
indique le nombre d'additions dues aux additions des blocs.})\, et
de profondeur \,$(\Ka +1)\,k$\,.
 
\smallskip
Consid\'erons \`a pr\'esent, pour \,$k=1,\ldots ,\nu$\,, la
matrice \,$U_k$\, dont les lignes sont les \'el\'ements de la
famille $\,\{RM^i\}_{i=0}^{2^k-1}$\, consid\'er\'ee comme une
matrice \,$2^{k}\times 2^{\eta}$\, et la matrice \,$V_k$\, dont
les colonnes sont les \'el\'ements de la famille
\,$\{M^{j2^\nu}S\}_{j=0}^{2^k-1}$\, consid\'er\'ee comme une
matrice \,$2^{\eta}\times 2^{k}$. La famille \,$T$\, s'obtient
alors en calculant la matrice \,$U_\nu\in\A^{2^\nu \times
2^{\eta}}$\, puis la matrice \,$V_\nu\in\A^{2^{\eta}\times
2^\nu}$\, et enfin le produit matriciel \,$W_\nu=U_\nu\,V_\nu$. La
famille \,$T=\{R\,M^iS\}_{i=0}^{r-1}$\, est enti\`erement
d\'etermin\'ee par la donn\'ee de la matrice
\,$W_\nu=(w_{ij})=U_\nu\,V_\nu \in\A^{2^\nu \times
2^\nu}$\,puisque~: \,$R\,M^kS=w_{ij}$\, ($\,0 \leq k \leq 2^{2\nu
-1}\,$) si et seulement si $\,k=(i-1)+(j-1)2^\nu$\,  i.e. si et
seulement si \,$j=\einf{{k \over {2^\nu}}}+1$\, et
\,$i=k+1-\einf{{k \over {2^\nu}}} \,2^\nu$.
 
\ss Le calcul de $T$ se fait donc en deux phases~: une premi\`ere
phase de calcul des matrices  \,$U_\nu$\, et \,$V_\nu$\, et une
deuxi\`eme phase de calcul du produit \,$W_\nu=U_\nu\,V_\nu$.
 
\mni {\bf $\bullet$~ Co\^ut de la phase 1:}
 
\ms Le calcul de \,$U_\nu$\, et de \,$V_\nu$\, se fait de proche
en proche \`a partir de \,$U_0=\left[\begin{array}{cc} R & 0
\end{array}\right] \in \A^{1\times 2^\eta}\,$,
\,$V_0=\left[\begin{array}{c} S \\ 0 \end{array}\right] \in
\A^{2^\eta \times 1}\,$\, et des puissances de \,$M$\, obtenues
par \'el\'evations successives au carr\'e, c'est-\`a-dire les
matrices \,$M^{2^s}$\, ($\,1 \leq s\leq 2\nu-1\,$).
 
\ss On a en effet, pour $k =1,\ldots ,\nu$: $$\begin{array}{lcl}
\{RM^i\}_{i=0}^{2^k-1} & = & \{RM^i\}_{i=0}^{2^{k-1}-1} \cup
\{RM^{i+2^{k-1}}\}_{i=0}^{2^{k-1}-1}~~~~ \mathrm{et}
\\[1mm] \{M^{j2^\nu}S\}_{j=0}^{2^k-1} & = &
\{M^{j2^\nu}S\}_{j=0}^{2^{k-1}-1} \cup
\{M^{2^\nu(j+2^{k-1})}S\}_{j=0}^{2^{k-1}-1} \\
\end{array}$$
Ce qui donne, de mani\`ere plus pr\'ecise, les relations
matricielles suivantes (si on pose $\,N=M^{2^\nu}\,)$~: 
$$ 
U_k =
\left[\begin{array}{c} U_{k-1} \\ 
\widetilde{U_{k-1}}
\end{array}\right] \in \A^{2^k\times 2^{\eta}}~~\mathrm{et}~~~
V_k =\left[\begin{array}{cc} V_{k-1} & \widetilde{V_{k-1}}
\end{array}\right] \in \A^{2^{\eta}\times 2^k}
$$
$$\mathrm{avec}~~~\widetilde{U_{k-1}}= U_{k-1}\times
M^{2^{k-1}}~~\mathrm{et}~~~ \widetilde{V_{k-1}}=N^{2^{k-1}}\times
V_{k-1}~.$$
 
\ss D'o\`u l'\algo suivant pour le calcul de \,$U_\nu$\, et
\,$V_\nu$\, (comportant \,$2\nu$\, \'etapes successives) \`a
partir des donn\'ees initiales \,$U_0$\, , \,$V_0\,$ (\cad $R\,$,
$S$)~:
 
\begin{enumerate}
\item
L'\'etape $k$ ($1\leq k\leq \nu$) consiste \`a calculer \,$U_k$\,
et \,$M^{2^k}$; pour cela deux \ops seront ex\'ecut\'ees en \paral
sur \,$M^{2^{k-1}}$\, qui est une matrice \,$2^{\eta} \times
2^{\eta}$\, (d\'ej\`a calcul\'ee \`a l'\'etape \,$k-1\,$):
l'\'elever au carr\'e et la multiplier \`a gauche par
\,$U_{k-1}$\, qui est une matrice \,$2^{k-1}\times 2^{\eta}$.
 
A la fin de ces \,$\nu$\, \'etapes, on obtient la matrice
\,$U_\nu$\, et la matrice \,$N=M^{2^\nu}$\,
(figure \ref{F2a}).
\begin{figure}[hbtp]   
\begin{center}
\includegraphics*[width=11cm]{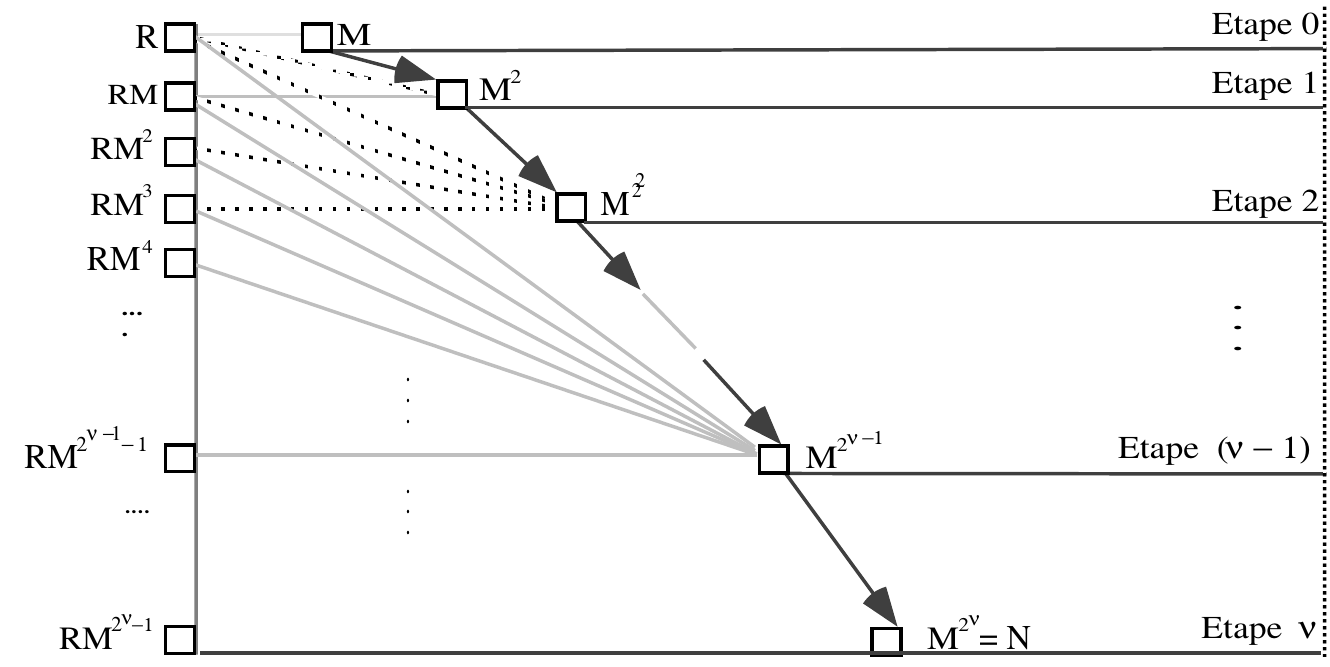}
\end{center}
\caption[Calcul de \,$U_\nu\,$]
{\label{F2a} Calcul de \,$U_\nu$\, }
 
\centerline{{\footnotesize les liens en trait
pointill\'e indiquent les \muls}}
 
\centerline{{\footnotesize \`a effectuer au cours
d'une \'etape pour passer \`a l'\'etape suivante}}
\end{figure}
 
\item
L'\'etape \,$\nu+k$\, ($1\leq k\leq \nu$) consiste \`a calculer
\,$V_k$\, et \,$N^{2^k}$: l\`a encore, il s'agit d'\'elever au
carr\'e une matrice \,$2^{\eta} \times 2^{\eta}$\, et de la
multiplier \`a droite par \,$V_{k-1}\,$ qui est une matrice
\,$2^{\eta} \times 2^{k-1}$.
 
A l'issue de ces \,$\nu$\, nouvelles \'etapes, on obtient
\,$V_\nu$\, et la matrice \,$N^{2^\nu}$\,
(figure \ref{F2b}).
\begin{figure}[hbtp]   
\begin{center}
\includegraphics*[width=11cm]{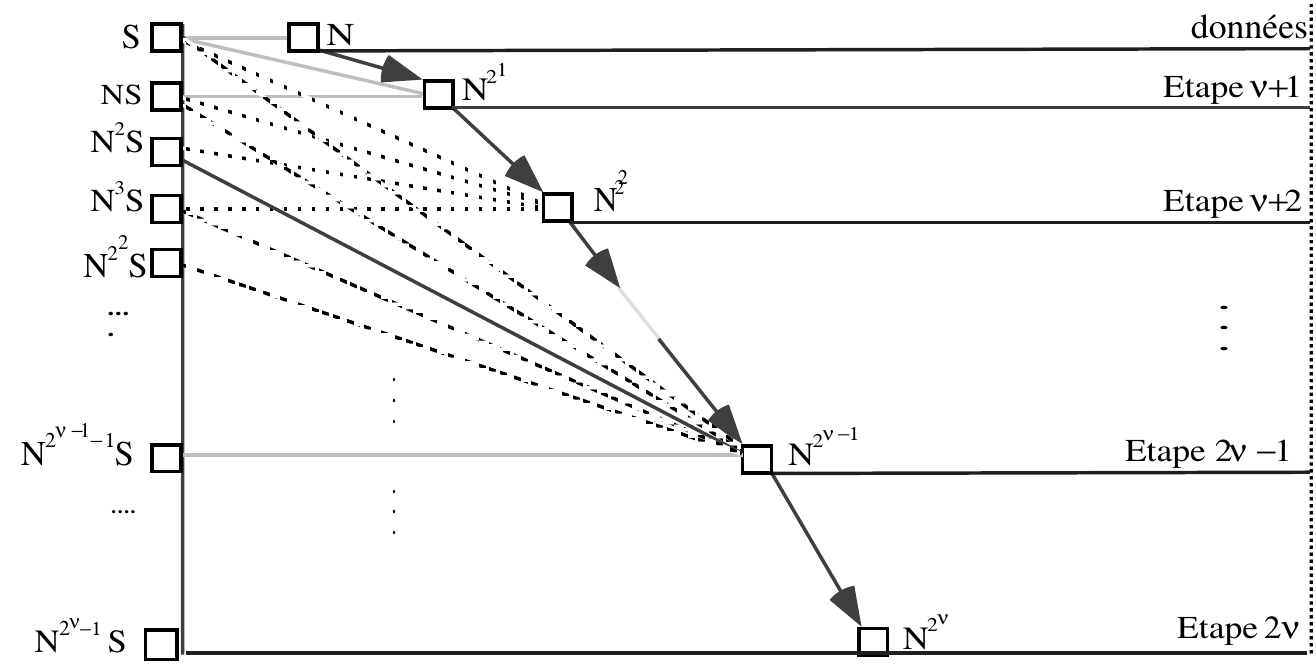}
\end{center}
\caption[Calcul de \,$V_\nu\,$]
{\label{F2b} Calcul de \,$V_\nu$\, }
 
\centerline{{\footnotesize les liens en trait
pointill\'e indiquent les \muls}}
 
\centerline{{\footnotesize \`a effectuer au
cours d'une \'etape pour passer \`a
l'\'etape suivante}}
\end{figure}
\end{enumerate}
 
Si l'on utilise les \muls par blocs \,$2^{k-1}\times 2^{k-1}\,$
(ils sont ici au nombre de \,$2^{\eta-k+1}$\, blocs ), l'\'etape
\,$k$\, (resp. \,$\nu + k$) ci-dessus est r\'ealis\'ee par un \cir
de taille~:
 
\noindent\centerline{$\,\Ca r^{\alpha} +
2^{\eta-k+1}\times[~2^{\eta-k+1}\Ca 2^{(k-1)\alpha} +
(2^{\eta-k+1}-1) 2^{2(k-1)}~]$.}
\\ 
et de \prof \'egale \`a $~\max\,{\{\Ka \,\eta ,(k-1)\,(\Ka -1) 
+\eta\}} = \Ka \,\eta\,$ (puisque \,$\Ka$\, est suppos\'e 
\,$\geq 1$).
 
Tenant compte du fait qu'il y a \,$2\nu$\, \'etapes et que \,$\eta
\leq 2\nu < \log{r}+2\,$, l'\algo calculant \,$U_\nu$\, et
\,$V_\nu$\, est donc r\'ealis\'e par un \cir de \prof \,$2\nu\Ka
\,\eta \leq \Ka\,(\log{r}+1)(\log{r}+2)$\, et de taille major\'ee
par~:\\[1mm]
$\,\Ca \,r^\alpha (\log{r}+2) + 2^{\eta+1}\,\sum_{k=1}^\nu (\Ca
\,2^{(\alpha-2)(k-1)}+1)$\, et donc par~: \\[1mm]
$(\Ca \,r^\alpha+2\,r)(\log{r}+2) + 4\,\Ca
\,r\sum_{k=1}^\nu2^{(\alpha-2)(k-1)}$
qui est \'egal \`a:
$$(\Ca \, r^\alpha+2\,r)(\log{r}+2) + 4\,\Ca \,
\frac{2^{(\alpha-2)\nu}-1}{2^{(\alpha-2)}-1}\,r\,.$$
Cette taille est donc major\'ee par~:
$$(\Ca\, r^\alpha+2\,r)(\log{r}+2) +
\frac{8\,\Ca}{2^{(\alpha-2)}-1}\,r^{\frac{\alpha}{2}} = \Ca\,
r^\alpha\log{r}+2\,\Ca\,r^\alpha+\O\,(r^{\alpha \over 2})\,. $$
qui est clairement \,$\O(r^\alpha\log{r})$\,.

\mni {\bf $\bullet$~ Co\^ut de la phase 2:} 

\ms Cette phase consiste \`a calculer le produit 
\,$U_\nu\times V_\nu$\, qui peut s'effectuer par des \muls de 
blocs \,$2^\nu \times 2^\nu$, en \paral et en 2 grandes 
\'etapes. 

\ss -- La premi\`ere \'etape consiste \`a calculer en \paral 
\,$\,2^{\eta-\nu}$\, produits de blocs \,$2^\nu \times 2^\nu$, 
avec 
\,$2^{\eta-\nu}\,\Ca\,2^{\nu\alpha}$\, \oparis dans l'anneau 
de base, ce qui donne une \prof totale de \,$\Ka\, \nu$; 

\ss -- Il s'agit dans la deuxi\`eme \'etape de calculer en 
\paral la somme des $2^{\eta-\nu}$ produits obtenus 
pr\'ec\'edemment, faisant intervenir 
\hbox{$(2^{\eta-\nu}-1)\,2^{2\nu}$} additions dans l'anneau de 
base, 
\`a l'aide d'une famille de \cirs binaires \'equilibr\'es de 
\prof \,$\eta-\nu\,$. 

\ss Le nombre total d'\oparisz, dans l'anneau de base, qui 
interviennent dans ces deux grandes \'etapes  du calcul de T, 
correspondant \`a une \prof totale de \,$(\Ka-1)\, \nu + 
\eta$, est donc major\'e par: 
$$ \,2^{\eta-\nu}\,\Ca \,2^{\alpha\nu}+ 
2^{\eta+\nu} \leq \Ca \,2^{(\alpha+1)\nu}+2^{3\nu} \leq \Ca 
\,2^{\alpha+1}\,r^{\frac{\alpha+1}{2}}+ 
8\,r^{\frac{3}{2}}\,.$$ 
 puisque \,$\eta \leq 2\nu$\,, \,$\alpha \leq 3$\, et \,$\Ca 
\geq 1\,$. 
Ce qui fait aussi \,$\O(r^{\frac{\alpha+1}{2}})$\, avec une 
constante asymptotique \'egale \`a \,$2^{\alpha +1}\Ca$\,. 
\ss Ainsi, le calcul de \,$T$\, \`a partir de \,$U_\nu$\, et 
\,$V_\nu$\, se fait par un \cir \paral de taille 
\,$\O(r^{\frac{\alpha+1}{2}})$\, et de \prof 
$$ \,(\Ka-1)\, \nu + \eta \leq {1 \over 2}\,(\Ka+1)\,\log{r} + 
\Ka\, 
\leq \Ka\, (1+\log{r})\, $$ 
puisque \,$\Ka \geq 1$, \,$\nu \leq {1 \over 
2}\,(2+\log{r})$\, et \,$r\geq 2$. 

\medskip  Nous r\'esumons dans le tableau ci-dessous l'analyse de 
\com 
qui vient d'\^etre faite et qui \'etablit le r\'esultat 
annonc\'e. 
\ss 
\begin{center} 
\label{coberk} 
\begin{tabular}{|l|c|c|} 
\hline Etapes~ & ~Profondeur~ & ~Taille~ \\ \hline 1\`ere 
phase & 
$\Ka\,(\log{r}+1)(\log{r}+2)$ & 
$\Ca\,r^\alpha\log{r}+\O\,(r^\alpha)$ 
\\ \hline 2\`eme 
phase & $\Ka\,(\log{r}+1)$ & $\O(r^{\frac{\alpha+1}{2}})$ \\ 
\hline Total & $\Ka\,(\log{r}+1)(\log{r}+3)$ & 
$\Ca\,r^\alpha\log{r}+\O\,(r^\alpha)$ \\ \hline 
\end{tabular}  \\ 
\end{center} 


\ms La diff\'erence essentielle avec l'\agb \cite{Ber} r\'e\-side
dans la simplification de la \recu permettant de calculer de
proche en proche les matrices \,$U_\nu$\, et \,$V_\nu$: \`a chaque
pas, la \mul par une seule matrice (au lieu de \,$\esup{
n^\epsilon }$\, matrices), avec recours \`a la \mul par blocs, a
permis de r\'eduire le nombre d'\oparis
 dans l'anneau de base, en \'eliminant
le facteur \,$n^\epsilon\,$.

\ss La d\'emonstration de la proposition \ref{RMiS} a permis de
donner une estimation pr\'ecise de la constante asymptotique:
cette constante est en effet \'egale \`a \,$\Ca $; elle est la
m\^eme que la constante asymptotique de la \mul des matrices.

\begin{theorem}\label{thBerko1}
Les \coes du \polcar d'une \maca d'or\-dre \,$n$\, peuvent \^etre 
calcul\'es par un \cari dont la taille et la \prof sont 
respectivement major\'ees par \,$\frac{1}{\alpha +1}\Ca\, 
n^{\alpha +1}\,\log{n} + \O(n^{\alpha +1})$\, et par 
\,$2\,\Ka\,\log^2{n} + \O(\log{n})\,.$ 

\end{theorem}
 
\prv
Le \polcar
de la matrice \,$A=(a_{ij})$\,
n'est autre que le \pol \,$P_n$\,
 donn\'e par la formule (\ref{top}):
$$ \overrightarrow{P_n}= \mathrm{Toep}(Q_n) \times
\mathrm{Toep}(Q_{n-1}) \times \cdots \times \mathrm{Toep}(Q_1). $$
 
\smallskip
\noindent Le calcul des \coes (de la forme \,$R\,M^iS\,$) du \pol
\,$Q_{k+1}$\, (pour $\,1 \leq k \leq n-1\,$) se fait, d'apr\`es la
proposition \ref{RMiS}, en \,$\O(k^\alpha\log{k})$\,. \\ De
mani\`ere plus pr\'ecise, le calcul de la totalit\'e des matrices
$\mathrm{Toep}(Q_{k+1})$ se fait donc 
avec une \prof major\'ee par \,$\Ka\,(\log{n} + 1)(\log{n} + 3)\,$ et
une taille major\'ee
par~: $$\Ca\,\sum_{k=1}^{n-1}[k^{\alpha}\log{k} +
2k^{\alpha}+\O(k^{\alpha \over 2})]~;$$ \cad par~: $$
\frac{\Ca}{\alpha +1}\, n^{\alpha +1}\,\log{n}+\frac{2\Ca}{\alpha
+1}\,n^{\alpha +1}+ \O(n^{{\alpha\over 2}+1})$$ \`a cause du
fait~: $$\sum_{k=1}^{n-1}k^{\alpha}\,\log{k}<
(\sum_{k=1}^{n-1}k^{\alpha})\,\log{n} ~~~\mathrm{et}~~~
\sum_{k=1}^{n-1}k^{\alpha}<n^{\alpha+1}\int_0^1x^{\alpha}dx=
\frac{n^{\alpha +1}}{\alpha +1}\,\cdot$$
 
\smallskip
D'autre part, le produit (\ref{top}) peut \^etre calcul\'e \`a
l'aide d'un \cir binaire \'equilibr\'e avec \,$\O(n^{\alpha+1})$\,
\oparis de base et une \prof major\'ee par \,$\Ka\,\log^2{n}$.
 
\ms Cela donne en fin de compte, dans l'anneau de base, un nombre
total d'\oparis major\'e par $$ \frac{\Ca}{\alpha +1}\,
n^{\alpha+1}\log{n} +  \frac{2\,\Ca}{\alpha +1}\,n^{\alpha+1} +
\O(n^{{\alpha\over 2}+1}) = \frac{\Ca}{\alpha +1}\,
n^{\alpha+1}\log{n} + \O(n^{\alpha+1}),$$ avec un \cari de \prof
major\'ee par 

\smallskip \centerline{
$
2\,\Ka\,(\log{n} + 1)^2\, + \Ka = 2\,\Ka\,\log^2{n} + 
\O(\log{n}).$} 
\begin{prop}\label{thBerko2}
Les \coes des \polcars de toutes les \smpds d'une matrice carr\'ee
d'ordre \,$n$\, peuvent \^etre calcul\'es en
\,$\SD(n^{\alpha+1}\log{n},\log^2{n})$\, (avec les m\^e\-mes
estimations que celles du \tho \emph{\ref{thBerko1}} pour les
constantes asymptotiques).
\end{prop}
 
\prv
En effet,
on a \,$A=A_n$, et les \coes du \polcar
\,$P_r$\, de la \smpd \,$A_r$\, de \,$A$\, ($\,1 \leq r \leq n\,$)
sont donn\'es par les vecteurs:
$$\overrightarrow{P_r}=\mathrm{Toep}(Q_r) \times
\mathrm{Toep}(Q_{r-1}) \times \cdots
\times \mathrm{Toep}(Q_1)$$
Ces vecteurs ne sont autres que les troncatures
successives (pour \,$r$\, allant de 2 \`a \,$n\,$) du
second membre de~(\ref{top}): ils peuvent donc \^etre
calcul\'es par un \algo \paral des pr\'efixes.
Le \cir que nous avons repr\'esent\'e
(figure \ref{F3}) correspond \`a l'une des solutions
du \gui{Calcul \paral des pr\'efixes} \cite{Lad}, que nous
avons pr\'esent\'ees dans la section \ref{subsec CPP}.
Il s'agit d'un \cir \paral de \prof
\,$\esup{ \log{n}} + 1$\, et de taille major\'ee par
\,$3n$.
\begin{figure}[hbtp]   
\begin{center}
\includegraphics*[width=11cm]{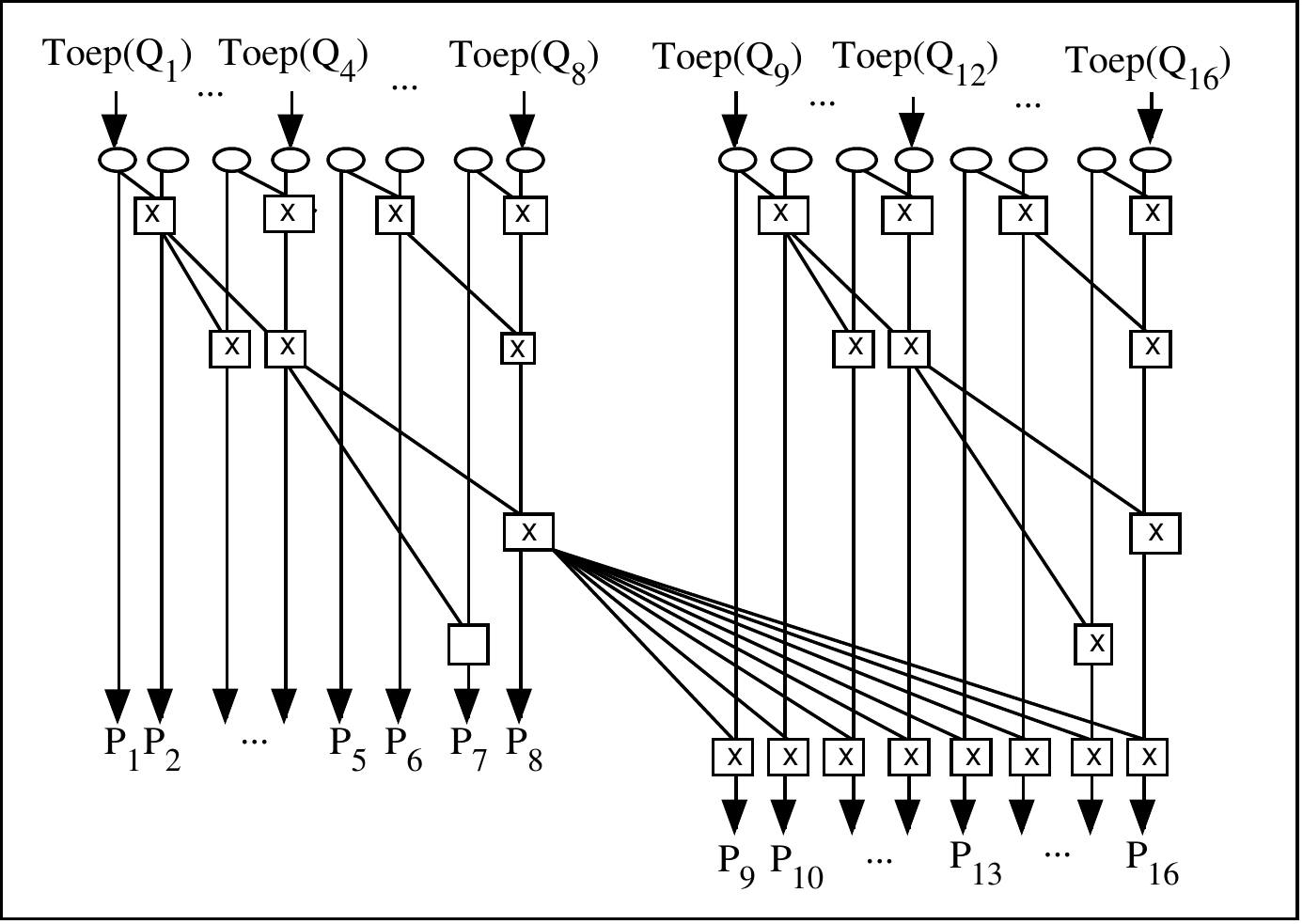}
\end{center}
\caption[Calcul des produits partiels]
{\label{F3} Calcul des \,$P_r$\, pour \,$1 \leq
r \leq n$\, (ici \,$n=2^4\,$) }
\end{figure}
Comme il s'agit de \muls matricielles, chaque
n{\oe}ud interne du \cir (repr\'esent\'e par une croix
dans la figure)
correspond \`a un \cir de \mul de matrices
de \prof \,$\Ka\,\log{n}$\, avec \,$\O(n^{\alpha})$\,
\oparis  dans l'anneau de base.
 
Le calcul des \,$P_r$\, ($\,2 \leq r \leq n\,$) \`a partir des
matrices \,$\mathrm{Toep}(Q_k)$\, se fait donc par un \cir  de
taille \,$\O(n^{\alpha +1})$\, et de \prof major\'ee par
\,$\Ka\,(\log{n}+1)\,(\log{n}+2)$. On conclut de la m\^eme fa\c
con que le \tho \ref{thBerko1} pour le produit des \mtosz. \qed
 
\begin{cor}
Le \deter et l'adjointe d'une \maca d'or\-dre \,$n$\, se calculent
en \,$\SD(n^{\alpha+1}\log{n},\log^2{n})$\, (avec les m\^emes
bornes que celles du \tho \emph{\ref{thBerko1}} pour les
constantes asymptotiques).
\end{cor}
 
\preuve
Le \deter de \,$A$\, n'est autre que
\,$\rPA(0)$. D'autre part, la matrice  adjointe
de \,$A$\, est donn\'ee par la formule:
$$\,\Adj(A)= -( p_0A^{n-1}+p_1A^{n-2}+
\cdots +p_{n-2}A+p_{n-1}\In)$$
o\`u  \,$\rPA(X)=\sum_{i=1}^np_{n-i}X^i$.
 
\mni \PrSa (voir section \ref{sec PrepaSar}) donnent un \algo
r\'ecursif (Powers\,$(A,n)$) pour calculer les \,$n$\, premi\`eres
puissances de \,$A$\, avec un \cari
\paral de \prof \,$\Ka\,
\log{n}\,(\log{n} + 1)$\, et de taille
major\'ee par \,$(2n-3)\,\Ca
\,n^{\alpha}$.(\footnote{~Le \gui{parallel
prefix algorithm} (section \ref{subsec CPP}) donne
le m\^eme r\'esultat pour la \prof mais
une taille major\'ee par
\,$3n\,\Ca \,n^{\alpha}$.})
 
Le r\'esultat est alors obtenu en remarquant que
\,$\mathrm{Adj\,}(A)$\, se calcule \`a partir des
puissances de \,$A$\, en \,$1+\esup{ \log{n}}$\,
\'etapes avec \,$\O(n^3)$\, \ops
arithm\'etiques de base.
\qed
 
\begin{remark}
\emph{
La \met de  Baur~\&~Strassen \cite{Baur} pour le calcul des d\'eriv\'ees
partielles (cf. section \ref{secpartialder})
montre que le calcul de l'adjointe d'une matrice a
toujours un co\^{u}t voisin de celui de son \deterz.
La construction originale ne se pr\'eoccupe pas de la \profz,
mais le r\'esultat a \'et\'e am\'elior\'e par Kaltofen et Singer
\cite{KalSin}: tout \cari de taille
\,$\tau$\, et de \prof \,$\pi$\, calculant une
fonction \polle (sur un anneau) ou une
fonction rationnelle (sur un corps) donne un
\cir de taille \,$4\tau$\, et de \prof \,$\O(\pi)$\,
qui calcule la fonction et toutes ses d\'eriv\'ees
partielles, et ceci ind\'ependamment du nombre de
variables d'entr\'ees du \cirz.
}
\end{remark}
\begin{remark}
\emph{Les \thos de \com que nous
venons d'\'e\-ta\-blir ne citent que deux param\`etres
de \comz: la taille et la \prof des \cirsz.
Mais une analyse minutieuse des \algos \'etudi\'es
nous permet \egmt d'avoir le nombre de
processeurs utilis\'es par ces \algos dans le
mod\`ele PRAM, \cad la largeur du \cari
correspondant. Ce troisi\`eme
param\`etre  peut \^etre
exprim\'e en fonction de la largeur d'un \cari
(de \prof \,$\Ka\,\log{n}$\, et de taille
\,$\Ca\,n^\alpha$) qui calcule le produit de deux \macas d'ordre
\,$n$.
Il est facile de v\'erifier que le r\'esultat trouv\'e
est le m\^eme que celui obtenu par application directe
du principe de Brent \`a cet \algo \paralz,
\cad un nombre de processeurs de l'ordre de
\,$\O(n^{\alpha+1}/\log{n})$.}
\end{remark}

\begin{remark}
\emph{Concernant les questions d'uniformit\'e et de co\^ut
de construction des \cirsz, ainsi que la taille des
\coes \itmdsz, le travail de
Matera \& Turull Torres \cite{Mater} donne, dans le
cas de l'anneau \,$\ZZ$\, des entiers relatifs, une
construction effective, avec une taille bien
contr\^ol\'ee des \coesz, des \cirs de base
qui interviennent dans l'\algo de Berkowitz.}
 
\noi \emph{Traduisant les \oparis  de \,$\ZZ$\, (addition et
\mulz) par des \cibos de \prof \,$\O(\log{b})$\, o\`u \,$b$\, est
la taille maximum de la \rpn binaire des \coes de la matrice
donn\'ee \,$A\in \ZZ^{n\times n}$, ils obtiennent~:} \\[1mm]
\emph{-- pour la \mul de deux matrices \,$n\times n$\, sur
\,$\ZZ$\, un \cibo de taille \,$\O(n^3\,b^2)$\, et de \prof
\,$\O(\log{(bn)})$~; \\ -- pour la taille des \coes \itmds
calcul\'es, une majoration de l'ordre de \,$\O(n(b+\log{n}))$~;
\\ -- pour l'\agbz, une \famu de \cibos de \prof
\,$\O(\log{(n)}\log{(bn)})$\, et de taille
\,$\O(n^6\,b^2\,\log^2{(n)})$\,.} \\[1mm] \emph{Cette
construction, appliqu\'ee \`a l'\algo am\'elior\'e que nous avons
pr\'esent\'e, donne une \famu de \cibos de m\^eme \profz, avec la
m\^eme majoration pour la taille des \coes \itmdsz, mais de taille
r\'eduite \`a \,$\O(n^{5}\,b^2\,\log^2{n})$. \\ Le facteur \,$n$\,
ainsi \'economis\'e provient essentiellement des \'etapes
correspondant aux figures \vref{F2a} et \vref{F2b} de notre
\algoz. }
\end{remark}
 
\section{M\'ethode de Chistov}
\label{secChistovPar}
 
\subsubsection*{Introduction}
 
On consid\`ere une \maca \,$A\in\A^{n\times n}$\,  on pose
\,$B=\In-XA\in\A[X]^{n\times n}$,
\,$B_r$\, est la \smpd d'ordre
\,$r$\, de \,$B$\, et \,$Q(X)=\det\,{B}$.
 
L'\algo est bas\'e sur les formules suivantes (ce sont les
\'equations (\ref{EqChis1}) et  (\ref{EqChis3}) \'etablies \`a la
section \ref{secChistovGene}) valables dans l'anneau des \dlis \`a
l'ordre \,$n$, \,$\A_n=\aqo{\A[X]}{X^{n+1}}$~:
\begin{equation} \label{EqChis01}
Q(X)^{-1}=\left[\,\det{(\In-XA)}\,\right]^{-1}=
\prod_{r=1}^n\left(B_{\,r}^{-1} \right)_{\\ r,r}.
\end{equation}
et, en notant \,$\Er$\, la \,$r$\,-\`eme colonne de $\Ir$:
 
\begin{equation} \label{EqChis03}
\left(B_{r}^{-1} \right)_{r,r} \ \mathrm{mod}\,X^{n+1}=
1+\sum_{k=1}^n\left(\tra{\Er}\,(A_r)^k\,\,\Er\right)\,X^k\,.
\end{equation}
 
Rappelons alors le principe \gnl de l'\algo \ref{algoChisG}
donn\'e en section
\ref{secChistovGene}.

\begin{agh}\textsf{\textbf{\Algo de Chistov, principe g\'en\'eral}
\Entree la matrice $A\in\A^{n\times n}$. \Sortie le \polcar $P(X)$
de $A$. 
\Debut 
\Etap{1}{}{Calculer  pour \,$r,k\in
\{1,\ldots,n\}$\, les produits \,$\tra{\Er}\,(A_r)^k\,\,\Er$,}
\hsu ce qui donne les \pols \,$\left(B_{\,r}^{-1} \right)_{\\
r,r}$ (formule $(\ref{EqChis03})$). 
\Etap{2}{}{Calculer le produit
des $n$ \pols pr\'ec\'edents modulo \,$X^{n+1}$,} \hsu ce qui
donne \,$Q(X)^{-1}\ \mathrm{mod}\,X^{n+1}$ (formule
$(\ref{EqChis01})$). 
\Etap{3}{}{Inverser modulo \,$X^{n+1}$\, le
\poly pr\'ec\'edent: on obtient \,$Q(X)$.} 
\Etap{4}{}{Prendre le
\poly r\'eciproque \`a l'ordre \,$n$\, du \poly \,$Q(X)$.} 
\hsu On
obtient \,$P(X)$\, en multipliant par $(-1)^{n}$. \fin }
\end{agh}

\subsection*{La version parall\`ele et sa complexit\'e}
 
\'Etudions pour chacune des \'etapes de cet \algoz, la taille et
la \prof d'un \cari correspondant qui tire le meilleur parti
de la \mul rapide des matrices et permet d'obtenir un temps \paral
en \,$\O(\log^{2} n)$.
 
\mni {\bf $\bullet$~ Co\^ut de l'\'etape 1:}
 
\sni Chacun des \elts \,$\tra{\Er}\,(A_r)^k\,\Er$\, qu'il s'agit
de calculer est obtenu en prenant la \,$r\,$-\,\`eme composante du
vecteur-colonne \,$(A_r)^k\,\Er$. On doit donc calculer
simultan\'ement, pour tous \,$r$\, compris entre \,$1$\, et
\,$n$\, et pour chaque \,$r$\,, les \,$n$\, produits
(matrice$\,\times\,$vecteur) \,$(A_r)^k\,\Er~~(1\leq k\leq n)\,$.
 
\ss Pour \'evaluer la \com de ce calcul, consid\'erons l'entier
\,$\nu\in \NN$\, tel que \,$2^{\nu -1} < r \leq 2^\nu$, \cad
\,$\nu =\esup{ \log{r} }$\,, et ramenons la matrice \,$A_r$\, \`a
une matrice $\, 2^\nu \times  2^\nu$\, en  remplissant de z\'eros
les rang\'ees \suptsz. Ainsi, toutes nos matrices \,$(A_r)^k$\,
seront consid\'er\'ees comme des matrices \,$2^\nu \times 2^\nu$,
et \,$\Er$\, comme une matrice $\, 2^\nu \times 1\,$: cela ne
change pas les produits \,$\tra{\Er}\,(A_r)^k\,\Er$\,
recherch\'es. Consid\'erons d'autre part l'entier \,$\eta\in
\NN$\, v\'erifiant \,$2^{\eta-1} \leq n < 2^{\eta}$\, \cad
\,$\eta=\einf{\log{n}}+1$\,.
 
\ss \noindent On proc\`ede alors en \,$\eta$\, sous-\'etapes
successives (num\'erot\'ees de 0 \`a $\eta-1$), chacune utilisant
le r\'esultat de la pr\'ec\'edente.
 
\ss \`A l'\'etape \,$j$\, $(0\leq j\leq \eta -1)$, on \'el\`eve au
carr\'e la matrice \,$(A_r)^{2^{j}}$\, puis on la multiplie \`a
droite par la matrice $$(\Er|A_r\,\Er|\dots|(A_r)^{2^j-1}\Er)\,\in
\, \A^{2^\nu\times 2^{j}}$$ pour obtenir la matrice
\,$(A_r)^{2^{j+1}}$\, et la  matrice
$$(\Er|A_r\,\Er|\dots|(A_r)^{2^j-1}\Er|(A_r)^{2^j}\Er|
\dots|(A_r)^{2^{j+1}-1}\Er)\,\in \, \A^{2^\nu\times 2^{j+1}}\,.$$
 
\`A la fin de ces \,$\eta$\, \'etapes (faisant \,$j=\eta -1\,$),
on obtient la matrice
$$(\Er|A_r\,\Er|\dots|(A_r)^{2^\eta-1}\Er)\,\in \, \A^{2^\nu
\times 2^{\eta}}$$ dont les \elts de la \,$r\,$-\,\`eme ligne,
plus pr\'ecis\'ement les \,$n\,$ premiers (on a \,$n <
{2^\eta}\,$), ne sont autres que les \elts
\,$\tra{\Er}\,(A_r)^k\,\Er$\, recher\-ch\'es.
 
Pour chaque \,$r~(1\leq r\leq n\leq 2^\eta-1)\,$, on a ainsi
\,$\eta$\, sous-\'etapes, chacune d'elles comportant
l'\'el\'evation au carr\'e d'une matrice \,$2^\nu \times
2^{\nu}$\, (en fait d'une matrice \,$r\times r\,$), et la \mul
d'une matrice \,$2^\nu \times 2^{\nu}$\, par une matrice \,$2^\nu
\times 2^j$. Utilisant pour cette derni\`ere \op les \muls par
blocs \,$2^{j} \times 2^{j}$\, (quitte \`a plonger la matrice
\,$2^\nu \times 2^{\nu}$\, dans une matrice \,$2^\eta \times
2^{\eta}$\, et la matrice \,$2^\nu \times 2^{j}$\, dans une
matrice \,$2^\eta \times 2^{j}$\,, on obtient pour chacune des
$\,\eta\,$ sous-\'etapes consid\'er\'ees un nombre d'\oparis
major\'e par~:
$$\Ca\,r^\alpha+2^{2\eta}\,[\Ca\,2^{(\alpha-2)j}+1]\leq
\Ca\,r^\alpha+4n^2\,[\Ca\,2^{(\alpha-2)j}+1].$$ Cela est d\^u au
fait que $\,\eta-1\leq \log{n}\,$ (et $\,2^{\eta}\leq 2n\,$) qui
permet d'obtenir les majorations suivantes~: $$\sum_{j=0}^{\eta
-1} \Ca\,r^\alpha \leq
\Ca\,r^\alpha\,\log{n}~~~\mathrm{et}~~~\sum_{j=0}^{\eta
-1}2^{(\alpha-2)j}\leq
\frac{2\,n^{\alpha-2}}{2^{\alpha-2}-1}\,\cdot$$ D'o\`u la
majoration du nombre d'\oparis intervenant (pour chaque valeur de
$r$) dans le calcul des $\,n\,$ produits
\,$\tra{\Er}\,(A_r)^k\,\Er\,$ ($\,1\leq k\leq n$)~: $$
\Ca\,r^\alpha\,\log{n}+
4n^2\,[\,\log{n}+\frac{2\,\Ca\,n^{\alpha-2}}{2^{\alpha-2}-1}\,]\,.$$
Comme $\,r\,$ varie de 1 \`a $\,n\,$, le calcul de l'\'etape~1
s'effectue \`a l'aide d'un \cari de taille \,$\O(n^{\alpha
+1}\log{n})$\, et de \prof \,$\O(\log^2{n})$.\\[1mm] Plus
pr\'ecis\'ement, la taille est major\'ee par~: $$
[\,4n^3+\Ca\,\sum_{r=1}^nr^\alpha\,]\,\log{n}+
\frac{8\,\Ca}{2^{\alpha-2}-1}\,\,n^{\alpha+1} \leq
\frac{\Ca}{\alpha+1}\,n^{\alpha+1}\,\log{n}+\O(n^{\alpha+1})\,.$$
et la profondeur par~:
$$\Ka\,\eta\,\log{n}+\sum_{j=0}^{\eta-1}(\Ka\,j+\eta -j) \leq
\frac{3\,\Ka+1}{2}\,\log^2{n}+\O(\log{n}).$$
 
\bs On peut remarquer qu'avec la \mul \usle des matrices
($\alpha=3\,$), l'\'etape~1 correspond \`a un \cari \paral de
\prof \,$\O(\log^2{n})$\, et de taille \,$\O(n^4\log{n}) +
\O(n^4)\,=\,\O(n^4\log{n})$\, puisque \,$\somm_{j=0}^{\eta
-1}\,n^22^j = n^2\,(2^\eta -1)< 2n^4 $.
 
\mni {\bf $\bullet$~ Co\^ut de l'\'etape 2:}
 
\sni
On doit calculer le produit tronqu\'e \`a l'ordre
\,$n$\, des \,$n$\, \pols de degr\'e $\leq n\,$
calcul\'es \`a l'\'etape pr\'ec\'edente.
Ce calcul se fait \`a l'aide d'un \cir
 binaire \'equilibr\'e en
\,$\SD(n^3,\log^2{n})$.
 
\mni {\bf $\bullet$~ Co\^ut de l'\'etape 3:}
 
\sni Il s'agit d'inverser modulo \,$X^{n+1}$\, le \poly
\,$\tilde{Q}$\, de degr\'e $\leq n$ obtenu \`a l'\'etape
pr\'ec\'edente. Ce \poly est de la forme \,$\tilde{Q}=1 - XR$\,
o\`u \,$R$\, est un \poly de degr\'e $n$ en $X$. Inverser
\,$\tilde{Q}$\, modulo \,$X^{n+1}$\, revient \`a calculer le
produit $$Q(X) = (1 +XR)(1 + X^{2}R^{2})\dots (1
+X^{2^\nu}R^{2^\nu})\ \mathrm{mod}\ X^{n+1}\,.$$ Cela s'effectue
en \,$\SD(n^2\log{n},\log{n}\log{\log{n}})\,$ \`a l'aide d'un \cir
binaire \'equi\-li\-br\'e.
 
\ms On peut acc\'el\'erer le calcul des deux \'etapes
pr\'ec\'edentes en utilisant une \mul rapide des \pols (mais cela
n'am\'eliore pas sensiblement le r\'esultat final).
 
\mni {\bf $\bullet$~ Co\^ut de l'\'etape 4:}
 
\sni Cette \'etape, de \prof $1$, n'intervient pas dans la \com de
l'\algoz.
 
\ms Nous donnons ci-dessous un tableau r\'esumant l'analyse
qui vient d'\^etre faite pour la \com de l'\algo de Chistov, montrant
que ce dernier est \,$\SD(n^{\alpha+1}\log{n},\log^2{n})$\, si
l'on utilise la \mul rapide des matrices ($\alpha < 3$) avec une
estimation pr\'ecise des constantes asymptotiques pour la taille
et pour la \profz.
 
\begin{center}
\label{cochis}
\begin{tabular}{|l|c|c|}
\hline Etape~ & ~Profondeur~ & ~Taille~ \\ \hline  ~ & ~ & ~ \\
Etape~1 & $\frac{3\Ka+1}{2}\,\log^2{n}+\O(\log{n})$ &
$\frac{\Ca}{\alpha+1}\,n^{\alpha+1}\,\log{n}+\O(n^{\alpha+1})$ \\
  ~ & ~ & ~ \\ \hline  ~ & ~ & ~ \\
Etape~2 & $\log^2{n}+\O(\log{n})$ & $\O(n^{3})$
\\   ~ & ~ & ~ \\ \hline   ~ & ~ & ~ \\
Etape~3 & $\O(\log{n}\log{\log{n}})$ & $\O(n^{2}\log{n})$ \\ ~ & ~ & ~ \\
\hline  Etape~4 & $1$ & n\'egligeable
\\ \hline
\end{tabular}  \\
 
\ms {\small {\bf Tableau \ref{cochis}}} \\ {\small Complexit\'e de
la version \paral de l'\algo de Chistov}
\end{center}
 
Si l'on utilise la \mul \usle ($\alpha = 3$), cela donne un \algo en
\,$\SD(n^4\log{n},\log^2{n})$\,. Dans ce dernier cas, l'\algo \sql \elr donn\'e \`a
la section \ref{secChisSeq} est donc pr\'ef\'erable (sur une machine \sqlez).
 
\begin{theorem}\label{thChistov}
~~L'\algo de Chistov calcule les \coes du \polcar d'une \maca d'or\-dre \,$n$\,
par un \cari  de \prof \,$\O(\log^2{n})$\, et de taille \,$\O(n^{\alpha+1}\log{n})$
avec des constantes asymptotiques estim\'ees respectivement \`a
$\,\frac{3}{2}\,(\Ka+1)\,$ pour la \prof et $\,\frac{1}{\alpha+1}\,\Ca\,$
pour la taille.
\end{theorem}

Enfin, comme dans le cas de l'\agbz, on obtient facilement
le r\'esultat compl\'ementaire suivant.
 
\begin{prop}\label{thchisto2}
Les \coes des \polcars
de toutes les \smpds d'une matrice
carr\'ee d'ordre \,$n$\, peuvent \^etre calcul\'es en
 \,$\SD(n^{\alpha+1}\log{n},\log^2{n})$\, par un \algo
directement d\'eriv\'e de celui correspondant au
th\'eor\`eme \ref{thChistov}.
\end{prop}
 
\begin{remark}
\emph{Remarquons que dans l'estimation de la taille des \caris
construits \`a partir des \algos de Chistov et de Berkowitz
am\'elior\'e, les termes en \,$n^{\alpha+1}\,\log{n}$\, sont les
m\^emes pour les deux \algos alors que les termes en
\,$n^{\alpha+1}\,$ sont respectivement estim\'es \`a
\,$\frac{8\,\Ca}{2^{\alpha-2}-1}\,n^{\alpha+1}\,$ pour Chistov et
\`a seulement $\,2\frac{\Ca}{\alpha+1}\,n^{\alpha+1}\,$ pour
Berkowitz am\'elior\'e (le rapport du premier \coe au second
\'etant strictement sup\'erieur \`a $16$). }
\end{remark}
\section[Applications des \algos]{Applications des \algos
\`a des \acoms}

\subsection*{Application  en \eva dynamique}
\label{conber}
 
Le calcul des \deters et des \polcars de toutes les
sous--matrices principales d'une matrice donn\'ee
trouve une application int\'eressante en
\emph{\eva dynamique}.
 
\ss Lorsqu'on travaille dans la cl\^oture
\agq dynamique \cite{D5} d'un corps
\,${\cal K}$, on se trouve dans la situation
standard suivante: on a des variables
\,$x_1,\ldots,x_n$\, qui repr\'esentent des
\elts \,$\xi_1,\ldots,\xi_n$\,
\agqs sur \,${\cal K}$.  On sait que ces
\elts v\'erifient un \sys \trg
d'\'equations  \agqsz.
 
\ss
De sorte que le corps \,${\K}[\xi_1,\ldots,\xi_n]$\,
est un quotient d'une \,${\K}\,$--\,\agr de dimension
finie
$$\A_{P_1, \ldots,P_n}=\aqo{\K[x_1,\ldots,
x_n]}{P_1(x_1),P_2(x_1,x_2),P_n(x_1,\ldots,x_n)}.$$
Chaque \,$P_i$\, est unitaire en \,$x_i$\, et cela
donne la structure de l'\agr de mani\`ere explicite.
 
\ss N\'eanmoins, cette \agr peut contenir
des diviseurs de z\'ero, ce qui signifie que plusieurs
situations diff\'erentes sont repr\'esent\'ees par
un seul calcul  dans \,$\A_{P_1,\ldots,P_n}$.
 
Lorsqu'on pose la question
\gui{ $Q(x_1,\ldots,x_n)=0\ ?$ }, le programme doit
calculer les \coes sous-r\'esultants de
\,$P_n$\, et \,$Q$\, par rapport \`a la variable
\,$x_n$\, (une discussion \gui{cas par cas} s'ensuit).
 
\ss \noindent Une solution est de calculer ces \coes dans l'\agr
${\cal K}[x_1,\alb\ldots,\alb x_{n-1}]$\, en utilisant l'\algo des
sous-r\'esultants \cite{GVLR,LRS,Lo} qui n\'ecessite des divisions
exactes et se situe naturellement dans le cadre d'un anneau
int\`egre, puis de les r\'eduire modulo l'id\'eal \,$\gen{P_1,\alb
\ldots,\alb P_{n-1}}$.
 
\ss
D\`es qu'on a trois $x_{i}$ le calcul s'av\`ere
tr\`es lourd. Une \'etude de \com montre
qu'on a un bien meilleur contr\^ole de la taille des
objets manipul\'es si on fait tous les calculs dans
l'\agr \,$\A_{P_1,\ldots,P_{n-1}}$.
 
\ss
Malheureusement l'\algo des
sous-r\'esultants ne peut plus s'ap\-pli\-quer.
En effet, des divisions requises par l'\algo
peuvent s'av\'erer impossibles, et m\^eme si l'\agr
est un corps, la division peut demander un effort
disproportionn\'e par rapport aux \mulsz.
 
\ss
Aussi semble-t-il que l'\agb (ou celui de Chistov),
appli\-qu\'e \`a la matrice de
Sylvester des \pols \,$P_n$\, et \,$Q$\, offre
la meilleure solution (en l'\'etat de l'art actuel)
pour calculer ces \coes sous-r\'esultants.
 
\ss
Il faut noter \`a cet \'egard que l'\algo
de Le Verrier-Fadeev-Csanky etc.
(en \cara nulle) ou celui propos\'e par
Kaltofen (cf. section \ref{sec kalto})
en \cara arbitraire  n'ont des performances \supees \`a l'\agb
que pour le calcul d'un \deter isol\'e,
mais non pour le calcul de tous les
\mips dominants d'une matrice donn\'ee.

\ss Signalons aussi que dans le cas o\`u on utilise l'\eva
dynamique pour la cl\^oture r\'eelle d'un corps ordonn\'e,
certaines discussions \gui{cas par cas} font appel aux signes de
tous les \coes sous-r\'esultants (\cf~\cite{tapasreels}).
 
\ss
Une autre application de l'\algo du calcul du \polcar en \eva
dynamique est la d\'etermination de la
signature d'une forme quadratique donn\'ee par une matrice
sym\'etrique arbitraire \,$S$.
Dans ce cas, la seule connaissance des signes des
\mips dominants
de la matrice \,$S$\, ne suffit pas toujours pour certifier
le rang et la signature{\footnote{Cela suffit dans le cas d'une matrice
\fregz.}}.
On pourra consulter \`a ce sujet le livre \cite{Gan}.
Mais il n'est pas difficile de voir que la connaissance des signes
des \coes du \polcar de la
matrice \,$S$\, permet de calculer et de certifier le rang
de \,$S$\, et la signature de la forme quadratique qui lui
est associ\'ee.
 
\subsection*{Cas des matrices creuses}
\label{creuses}
 
Signalons pour terminer que l'\agb et celui de Chistov sont
particu\-li\`erement bien adapt\'es au cas des matrices creuses,
notamment en version \sqle \elr o\`u le nombre d'\ops passe de
\,$\O(n^{4})$\, \`a \,$\O(n^{3})$\, lorsque seulement \,$\O(n)\,$
\coes de la matrice sont non nuls.

Parmi les autres \algos \'etudi\'es,
celui de Kaltofen-Wiede\-mann peut
\egmt \^etre adapt\'e au cas
des matrices creuses, avec une diminution similaire du
nombre d'\oparisz.

\newpage \thispagestyle{empty}

 
 
\chapter{R\'esultats exp\'erimentaux}
\label{chap experim}

\section{Tableaux r\'ecapitulatifs des complexit\'es}
\label{recap}
 
Dans cette section, nous donnons les tableaux r\'ecapitulatifs
des complexit\'es arithm\'etiques th\'eoriques pour les diff\'erents
\algos \'etudi\'es.
 
Y figure notamment le tableau des \coags
des \algos en version \sqle \elr (\cad n'utilisant que la \mul \usle des
matrices, des \pols et des entiers) que nous avons
exp\'eriment\'es.
 
\subsubsection*{Abbr\'eviations utilis\'ees}
\label{subsubsec}

Le mot {\sc C}te signifie \gui{constante asymptotique}
(pour les estimations de taille des \cirsz),
et {\sc Val.} signifie \gui{Domaine de validit\'e}:
\begin{itemize}
\item A.C.A. signifie \gui{anneau commutatif arbitraire},
\item A.I.A.D. signifie \gui{anneau int\`egre poss\'edant un \algo pour
les divisions exactes},
\item A.I.C. signifie \gui{anneau int\`egre et int\'egralement clos
poss\'edant un \algo pour les divisions exactes},
\item D. $n!$\, signifie \gui{la division par \,$n!$\, quand elle est
possible, est unique et explicite}.
\item Prob. signifie \gui{\algo de nature probabiliste}, il s'agit de
l'\algo de Wiedemann, qui fonctionne sur les corps,
avec des variantes possibles
dans le cas A.I.A.D.
\end{itemize}
 
Les sigles M.R.P. et M.U.P. d\'esignent
respectivement
la \mul rapide et la \mul \usle des
\polsz.
 
Rappelons que nous notons \,$\mu_P(n)$\, le nombre d'\oparis
dans la multiplication de deux \pols de degr\'e \,$n$\, en \prof
\,$\O(\log n)$. En M.U.P.
\,$\mu_P(n)=\O(n^2)$, avec la \met de Karatsuba \,$\mu_P(n)=\O(n^{\log
3})$, et en  M.R.P.
\,$\mu_P(n)=\O(n\log{n}\log{\log{n}})$\, ou \,$\O(n\log{n})$\, selon les
anneaux.
 
Les initiales  G et JB d\'esignent les \algos de Gauss
(sur un corps) et de
\JB (sur un anneau int\`egre poss\'edant un \algo pour les divisions
exactes) pour le calcul des d\'eterminants.
Rappelons que l'\ajbz, qui consomme un peu plus d'\oparisz,
pr\'esente des avantages significatifs
par rapport \`a l'\apgz, dans le nombreux anneaux commutatifs
(comme par exemple les anneaux de \pols \`a \coes entiers).

\label{T1}
\vspace{5mm}
\begin{center}
{\Large\textbf{A: Calcul des D\'eterminants}} \\[5mm]
{\bf M\'ethodes \sqles simples} \\[3mm]
\noi
\setlength{\extrarowheight}{1mm}
\begin{tabular}{|m{4.7cm}|c|c|c|c|}
\hline {\sc Algorithme} & {\sc Taille} & Cte &
{\sc Val.} & {\sc Date} \\[1mm]
\hline
{\sc Gauss} & $\O(n^3)$ &
$2/3$ & Corps &
$<1900$ \\[1mm]
\hline
{\sc Jordan-Bareiss}
& $\O(n^3)$ &
$4/3$ & A.I.A.D. &
$<1900$
\\[1mm]
\hline
{\sc Gauss avec \'{e}limina\-tion des divisions}
& $\O(n^5)$ & $1/3$ &
A.C.A. & $ \,\, 1973 $
\\[1mm]
\hline
{\sc Jordan-Bareiss mo\-difi\'{e}}
& $\O(n^5)$ & $1/10$ &
A.C.A. & $\,\,1982 $ \\[1mm]
\hline
\end{tabular}\\[5mm]

{\bf M\'ethodes rapides en \prof $\O(n)$} \\[3mm]
\begin{tabular}{|m{4.7cm}|c|c|c|c|}
\hline
{\sc Bunch\&Hopcroft}
 &  $\O(n^\alpha)$ & $\gamma_\alpha $(*)&
 ~Corps~  &  ~$1974$~ \\[1mm]
\hline
\end{tabular}
\end{center}

\mni(*) Voir \tho \ref{thBU} et proposition \ref{propBU}.

\medskip
Dans le premier tableau ci-apr\`es nous avons rajout\'e la colonne {\sc Cr.}
pour le traitement des
matrices creuses:
si une matrice \,$C\in\A^{n{\times}n}$\, a environ \,$k\cdot n$\, \coes
non nuls,
certains \algos sont acc\'el\'er\'es et leur temps d'ex\'ecution \sql
divis\'e par \,$n/k$. Nous avons indiqu\'e cette possibilit\'e
d'acc\'el\'e\-ra\-tion par un \gui{oui} dans la colonne {\sc Cr.}

\newpage
\begin{center} 
{\Large{\bf B: Calcul du Polyn\^ome Caract\'eristique}}\\[1mm] \label{T4} {\bf
Versions \sqles simples} \\[3mm]
\setlength{\extrarowheight}{1mm}
\begin{tabular}{|m{3.8cm}|c|c|c|c|}
\hline
{\sc Algorithme} & {\sc Taille}  &
 Cte & {\sc Val.}& {\sc Cr.} \\[1mm]
\hline
{\sc Wiedemann}  &
$\O(n^{3})$ & $2$ & Prob. & oui
\\[1mm]
\hline
{\sc Hessenberg} &
$\O(n^{3})$ & $2$ & Corps&
\\[1mm]
\hline
{\sc Frobenius} &
$\O(n^{3})$ & $10/3 $ & A.I.C. &
\\[1mm]
\hline
{\sc \PrSaz} &
$\O(n^{3,5})$ & $2 $ & D. $n!$ &
\\[1mm]
\hline
{\sc Berkowitz  }
&    $\O(n^{4})$  &    $1/2$  &
 A.C.A. & oui
\\[1mm]
\hline
{\sc Chistov} & $\O(n^{4})$ & $2/3$ & A.C.A. & oui
\\[1mm]
\hline
{\sc Faddeev-Souriau- Frame}~ (Le Verrier)
& $\O(n^{4})$ & $2$ & D. $n!$ & oui
\\[1mm]
\hline
{\sc Interpolation} ~~~~~~~ (Lagrange)  &
$\O(n^{4})$ &
\begin{tabular}{c}  $2/3$ (G)\\  $4/3$ (JB)\end{tabular} &
\begin{tabular}{c} Corps \\ A.I.A.D. \end{tabular}
&\\[1mm]
\hline
{\sc Kaltofen-Wiedemann}  &
$\O(n^{4})$ & $8$ & A.C.A. &oui
\\[1mm]
\hline
{\sc Jordan-Bareiss modifi\'e}
& $\O(n^{5})$ & $1/10$ & A.C.A.&
\\[1mm]
\hline
{\sc Gauss avec \elidz}
& $\;\O(n^{5})\;$ & $1/3$ & A.C.A.
&\\
\hline
\end{tabular} \\[5mm]

{\bf Taille avec \mul rapide des \polsz} \\[3mm]
\begin{tabular}{|l|c|}
\hline
{\sc Kaltofen-Wiedemann}  &
$\O(n^{3}\mu_P(\esup{\sqrt n\,})\log n)$ 
\\[1mm]
\hline
{\sc Jordan-Bareiss modifi\'e}
& $\O(n^{3}\mu_P(n))$ 
\\[1mm]
\hline
{\sc Gauss avec \elid}
& $\O(n^{3}\mu_P(n))$ 
\\
\hline
\end{tabular}

\end{center}

\newpage
\begin{center}
{\Large{\bf C: Calcul du Polyn\^ome
Caract\'eristique}}\\[1mm]    \label{T2}
{\bf M\'ethodes s\'equentielles rapides} \\[3mm]
\setlength{\extrarowheight}{1mm}
\begin{tabular}{|m{3,5cm}|c|c|c|}
\hline{\sc Algorithme}& {\sc Taille}&
Cte & {\sc Val.} \\[1mm] \hline
{\sc Keller-Gehrig} &
$\O(n^{\alpha}\log\,n)$ & $\cdots$ &  Corps
\\[1mm]
\hline
{\sc Interpolation} ~~~~ (de Lagrange)
    & $\O(n^{\alpha +1})$ &
    $\gamma_\alpha$
&
    Corps
\\[1mm] \hline
{\sc Faddeev-Souriau-Frame} &
$\O(n^{\alpha +1})$ & $\Ca$ &  D. $n!$
\\[1mm] \hline
{\sc Kaltofen-Wiedemann} &
$\O(n^{\frac{\alpha+3}{2}}\mu_P(\esup{\sqrt n\,}))$
&
    $\cdots$
&
A.C.A.
\\[1mm]
\hline
\end{tabular}
\end{center}

\vspace*{0,3cm}
 
\begin{center}
{\Large  {\bf D: Calcul du Polyn\^ome Caract\'eristique}} \\[1mm]
{\bf M\'ethodes parall\`eles en \prof $\O(\log^2{n})$} \\[3mm]
\label{T3}
\setlength{\extrarowheight}{1,5mm}
\begin{tabular}{|m{2.8cm}|m{2,5cm}|c|c|c|c|}
\hline
{\sc Algorithme} & {\sc Taille}  &
Cte &
{\sc K} & {\sc Val.}
\\[1mm]
\hline
{\sc Csanky \, 1976} & $\O(n^{\alpha +1})$ & $4\,\Ca$ &
$\Ka$ & D.  $n!$\,
\\[1mm]
\hline
{\sc Preparata~\&\,\,\,\, 
Sarwate \, 1978}  &
$\O(n^{\alpha + {1\over 2}})$ & $4\,\Ca$ &
$5\,\Ka$ & D.  $n!$\,
\\[1,5mm]
\hline
{\sc Galil \&  Pan \,\, 1989}  &
 $\O(n^{\alpha + {1\over 2} -\delta(\alpha)})$
 & $\cdots$ &
$\cdots$ & D.  $n!$\,
\\[1,5mm]
\hline
 {\sc B.H.G. (\dag) \,\,\,\,\,\,\, 1982}  &
$\O(n^{18+\epsilon })$~~(*) \hspace{1cm}
$\,\O(n^{21})$~~~~~~(**)  &
 $\cdots$   &
$\cdots$ & A.C.A.
\\[1,5mm]
\hline
 {\sc Berkowitz \, 1984}   &
 $\O(n^{\alpha + 1 + \epsilon})$  & $\cdots$ &
$\cdots$ & A.C.A.
\\[1,5mm]
\hline {\sc Chistov \, 1985} & $\O(n^{\alpha + 1}\,\log{n})$ &
$\frac{1}{\alpha+1}\,\Ca$ &
$\,\frac{3}{2}\,(\Ka+1)\,$ & A.C.A.
\\[1,5mm]
\hline {\sc Berkowitz am\'elior\'e \,1985}
& $\O(n^{\alpha +
1}\,\log{n})$ & $\,\frac{1}{\alpha+1}\,\Ca\,$ &
$3\,\Ka$ & A.C.A.
\\
\hline
\end{tabular}\\[5mm]
\end{center}
\noi (\dag) Borodin, Hopcroft \& v. z. Gathen. (*) M.R.P. (**) M.U.P.\\
 La colonne {\sc K} donne la constante asymptotique
 du  temps \paral en  $\O(\log^2{n})$.
 Le nombre $\delta(\alpha) > 0 $ d\'epend de \,$\alpha$.
 Enfin $\epsilon$ est
 positif arbitrairement petit.

\newpage
 
 
\section{Pr\'esentation des tests}
Les \algos consid\'er\'es dans les tableaux
de comparaison que nous pr\'e\-sen\-tons ci-dessous
ont \'et\'e exp\'eriment\'es \`a l'aide du logiciel
de Calcul Formel {\sc Maple} et \'ecrits dans le
langage de programmation qui lui est rattach\'e\footnote{Les
programmes ont tourn\'e avec la version
Maple~V~Release~3.}.
 
Les \algos sont ceux du tableau B, \cad les versions
\sqles simples pour le calcul du \polcarz.
Nous avons indiqu\'e \'egalement dans la colonne
\guig linalpoly \guid\, les performances de
l'\algo donn\'e par {\sc Maple} dans la version V.
Les versions plus r\'ecentes du logiciel utilisent
d\'esormais l'\agbz.
 
Chacun des tests de comparaison entre les diff\'erents
\algos a \'et\'e effectu\'e sur une m\^eme machine,
avec le m\^eme \'echantillon de matrices.
 
Les matrices utilis\'ees font partie de l'un des groupes
suivants, selon le type de l'anneau de base choisi:\\[3mm]
$\bullet~${\bf Groupe~1:}
les matrices \,$randmatrix(n,n)$\,
qui sont des \macas d'ordre \,$n$\, \`a
coefficients pris au hasard (entre -99 et +99) dans
l'anneau \,$\ZZ$\, des entiers relatifs; \\[2mm]
$\bullet~${\bf Groupe~2:}
les matrices \,$Mathard(n,x,y)$\, dont les
\elts sont des \pols en \,$[x,y]$\, de
degr\'e total \,$\leq 5$. Les coefficients de ces
\pols de \,$\ZZ[x,y]$\, sont aussi des entiers
compris entre -99 et +99;\\[2mm]
$\bullet~${\bf Groupe~3:}
les matrices \,$Matmod(n,lisvar,Ideal,p)$\, qui
sont des \macas d'ordre \,$n$\, dont les
coefficients sont des \elts choisis au hasard
dans l'anneau-quotient $$\aqo{\ZZ_p[lisvar]}{Ideal}$$
o\`u \,$p$\, est un entier positif (on le prendra
premier), \,$lisvar$\, une liste donn\'ee de
variables et \,$Ideal$\, une liste donn\'ee de
\pols en \,$lisvar$\, \`a coefficients dans
$\,\ZZ$.
L'anneau de base est donc ici, sauf exception, un
anneau dans lequel la division n'est pas permise. \\[2mm]
$\bullet~${\bf Groupe~4:}
les matrices \,$Jou(n,x)$, carr\'ees d'ordre
$\,n$, \`a coefficients dans \,$\ZZ[x]$, dont
les coefficients sont donn\'es par:
$$[Jou]_{\,ij} = x+x^2(x-ij)^2+(x^2+j)(x+i)^2~~~~
\mbox{ pour }~~~1\leq i,j\leq~n\,.$$ Quelle que soit
la valeur de \,$n$, le rang de la matrice
$\,Jou(n,x)$\, ne d\'epasse pas 3: c'est ce qui
explique la sup\'eriorit\'e, dans ce cas, des \algos
de Souriau-Faddeev et de \JB (nettement
plus performants pour les matrices de rang petit). \\[2mm]
$\bullet~${\bf Groupe~5:}~
ce sont des matrices creuses \`a coefficients entiers
choisis au hasard entre -99 et +99.
Elles sont donn\'ees par la proc\'edure {\sc Maple}
{\tt {\bf randmatrix}(n,n,{\bf sparse})\,.}\\[2mm]
Quant aux machines utilis\'ees, il s'agit
essentiellement d'un DEC~Alpha-600 \`a 175 Mhz et
320 Mo de m\'emoire centrale.\footnote{~Gr\^ace
notamment \`a l'hospitalit\'e du Laboratoire GAGE
(Ecole Polytechnique).}\\[3mm]
Les matrices intervenant dans les comparaisons
sont g\'en\'er\'ees par des codes {\sc maple}: une
proc\'edure {\sf Matmod} par exemple cr\'ee une
matrice du Groupe~3 \`a partir de la donn\'ee de deux
entiers positifs \,$n$\, (la taille
de la matrice) et \,$p$\, (on calcule modulo
\,$p$), d'une liste de variables \,$lisvar$,
et d'une liste \,$Ideal$\, de \pols en \,$lisvar$\,
comprenant autant de \pols \,$P_i$\, que de
variables \,$x_i$, chacun des \,$P_i$\, \'etant un
\pol en \,$[x_1, \ldots ,x_i]$, unitaire en
$\,x_i$. Ceci afin d'illustrer le genre d'application
de l'\agb lorsqu'on se place dans
l'\agr \,$\aqo{\ZZ_p[lisvar]}{Ideal}$, et la situation
indiqu\'ee dans la section \ref{conber}.

La proc\'edure {\sf Matmod} utilise comme sous-proc\'edure la
proc\'edure {\sf polmod} (donn\'ee dans l'annexe) qui prend en entr\'ee
un nombre entier \,$p$, un \pol \,$P$\, de
$\,\ZZ[lisvar]$, et donne en sortie un repr\'esentant simple
de l'image canonique de \,$P$\, dans l'anneau-quotient
\,$\aqo{\ZZ_p[lisvar]}{Ideal}$.

 
\section{Tableaux de Comparaison} \label{TabCompar}

Nous donnons dans les trois pages qui suivent les tableaux
correspondant aux cinq groupes de matrices que nous avons
pr\'ec\'edemment indiqu\'es.

Il ne s'agit que de quelques exemples, mais ils sont significatifs.

La comparaison entre le comportement pratique des \algos
montre un bon accord avec les calculs th\'eoriques de 
complexit\'e, surtout si on prend en compte la taille des objets 
\itmds cr\'e\'es par les diff\'erents \algosz. 
Sauf exception l'\agb est le plus performant, suivi de pr\`es 
par celui de Chistov.

Les performances a priori meilleures
pour les \algos de Hessenberg, Frobenius et Wiedemann
ne se r\'ev\`elent qu'avec des tests portant sur
des matrices \`a \coes dans des corps finis. 
En effet l'avantage en nombre d'\oparis est
contrebalanc\'e par la plus mauvaise taille des objets 
\itmds manipul\'es, par exemple d\`es que l'anneau des \coes contient
\,$\ZZ$. Il aurait fallu cr\'eer un autre groupe de matrices
pour mettre en \'evidence cet avantage.

Il serait \egmt int\'eressant d'\'elargir l'exp\'erimentation
en impl\'e\-mentant la version \sqle simple de l'\algo de \PrSaz.

\medskip Dans le groupe 1 nous avons pris des \macas
d'ordre \,$n$\, \`a coefficients dans \,$\ZZ$\, pour 10
valeurs de \,$n$\, comprises entre 16 et 128.
 
Dans le groupe 2, ce sont des \macas d'ordre \,$n$\,
\`a coefficients dans \,$\ZZ[x,y]$\, pour \,$n\in\{10,12,20\}$\,
et des matrices \`a coefficients dans \,$\ZZ[x]$\, pour
$\,n\in\{10,15,20,25\}$.
 
Parmi les matrices du groupe 3, nous avons pris des matrices
carr\'ees d'ordre $n\in\{8,10,12,16\}$ \`a coefficients dans
\,$\aqo{\ZZ_7[x]}{x^3-1}$\, (pour lesquelles Faddeev ne s'applique
pas) et des matrices \`a \coes dans $\,\aqo{\ZZ_{17}[x,y]}{H,L}$\,
(pour lesquelles Faddeev s'applique). \\ Ici \,$\gen{H,L}$\, est
l'id\'eal engendr\'e par les deux \pols \,$H=x^5-5xy+1$\, et
\,$L=y^3-2y+1$.
 
Dans le groupe 4, ce sont des \macas d'ordre
$\,n\in\{10,15,\alb 20,\alb 25\}$\, \`a coefficients dans \,$\ZZ[x]$,
mais de rang petit \,$\leq 3$.
 
Enfin les matrices du groupe 5 (des matrices creuses \`a \coes
entiers choisis au hasard entre -99 et +99) ont \'et\'e prises parmi
les matrices \,{\bf randmatrix}(n,n,{\bf sparse})\,
telles que \,$n\in\{32,50,64,128,200\}$.

\newpage
 
\vspace*{1cm}
 
\rotatebox{90}{ {\begin{minipage}[t]{17cm}
{
{\Large {\bf Premier Groupe~:~Matrices denses \,$n\times n\,$}}\\
{\large (Matrices \`a coefficients entiers)}\\[5mm]
\begin{tabular}{|p{3.5em}|p{5em}|@{}c@{}|@{}c@{}|@{}c@{}|@{}c@{}|@{}c@
{} |@{}c@{}|} \hline \multicolumn{2}{|c|}{matrice} &
\multicolumn{1}{|c|} {linalpoly \dag} & berkosam & Chistov &
Faddeev & barmodif \ddag& Hessenberg \\
\hline $n = 16$ &
\begin{tabular}{@{}l}CPU~Time \\ Mem. All. \end{tabular} &
\begin{tabular}{l} 0' 05" \\ 1,~244 Mb  \end{tabular} &
\begin{tabular}{l} 0' 02 " \\ 3,~014 Mb \end{tabular} &
\begin{tabular}{l} 0' 05" \\ 2,~003 Mb  \end{tabular} &
\begin{tabular}{l} 0' 12" \\ 1,~376 Mb   \end{tabular} &
\begin{tabular}{l} 0' 09" \\ 1,~244 Mb  \end{tabular} &
\begin{tabular}{l} 0' 03" \\ 1,~969 Mb  \end{tabular}
\\
\hline $n = 20$  &
\begin{tabular}{@{}l}CPU~Time \\ Mem. All. \end{tabular} &
\begin{tabular}{l} 0' 16" \\ 1,~834 Mb  \end{tabular} &
\begin{tabular}{l} 0' 04" \\ 3,~538 Mb  \end{tabular} &
\begin{tabular}{l} 0' 07" \\ 1,~834 Mb  \end{tabular} &
\begin{tabular}{l} 0' 28" \\ 2,~096 Mb  \end{tabular} &
\begin{tabular}{l} 0' 24" \\ 2,~031 Mb  \end{tabular} &
\begin{tabular}{l}  0' 20" \\ 5,~254 Mb  \end{tabular} \\
\hline $n = 25$ & \begin{tabular}{@{}l}CPU~Time \\ Mem. All.
\end{tabular} &
\begin{tabular}{l} 0' 40" \\ 1,~900 Mb  \end{tabular} &
\begin{tabular}{l} 0' 17" \\ 1,~900 Mb \end{tabular} &
\begin{tabular}{l} 0' 17" \\ 1,~900 Mb  \end{tabular} &
\begin{tabular}{l} 1' 07" \\ 2,~489 Mb  \end{tabular} &
\begin{tabular}{l} 1' 04" \\ 1,~834 Mb  \end{tabular} &
\begin{tabular}{l} 1' 28" \\ 14,~468 Mb  \end{tabular} \\
\hline $n=32$ & \begin{tabular}{@{}l}CPU~Time \\ Mem. All.
\end{tabular} &
\begin{tabular}{l} 2' 09" \\ 1,~507 Mb  \end{tabular} &
\begin{tabular}{l} 0' 30" \\ 1,~965 Mb  \end{tabular} &
\begin{tabular}{l} 0' 44 " \\ 1,~965 Mb  \end{tabular} &
\begin{tabular}{l} 3' 16" \\ 3,~014 Mb  \end{tabular} &
\begin{tabular}{l} 3' 25" \\ 2,~817 Mb  \end{tabular} &
\begin{tabular}{c} $\approx$ 2H \\ OUT$^*$ \end{tabular} \\
\hline $n=40$  & \begin{tabular}{@{}l}CPU~Time \\ Mem. All.
\end{tabular} &
\begin{tabular}{l} 6' 08" \\ 2,~096 Mb \end{tabular} &
\begin{tabular}{l} 1' 13" \\ 4,~062 Mb  \end{tabular} &
\begin{tabular}{l} 1' 48" \\ 2,~031 Mb  \end{tabular} &
\begin{tabular}{l} 17' 58" \\ 5,~241 Mb \end{tabular} &
\begin{tabular}{l} 11' 28" \\ 5,~307 Mb \end{tabular} &
\begin{tabular}{c}  -- \\ --  \end{tabular} \\
\hline $n=50$  &
\begin{tabular}{@{}l}CPU~Time \\ Mem. All. \end{tabular} &
\begin{tabular}{l} 19' 42" \\ 2,~424 Mb \end{tabular} &
\begin{tabular}{l} 3' 04" \\ 4,~193 Mb  \end{tabular} &
\begin{tabular}{l} 4' 37" \\ 2,~096 Mb  \end{tabular} &
\begin{tabular}{l} 2H 38' 50" \\ 9,~500 Mb  \end{tabular} &
\begin{tabular}{l} 42' 38" \\ 8,~124 Mb \end{tabular} &
\begin{tabular}{c}  -- \\ --  \end{tabular} \\
\hline $n = 64$  &
\begin{tabular}{@{}l}CPU~Time \\ Mem. All. \end{tabular} &
\begin{tabular}{l} 1H 17' 17" \\ 1,~900 Mb  \end{tabular} &
\begin{tabular}{l} 8' 47" \\ 4,~324 Mb  \end{tabular} &
\begin{tabular}{l} 13' 23" \\ 2,~162 Mb \end{tabular} &
\begin{tabular}{l} 13H 15' 45" \\ 18,~870 Mb  \end{tabular} &
\begin{tabular}{l} 3h 54' 46" \\ 14,~284 Mb  \end{tabular} &
\begin{tabular}{c}  -- \\ --  \end{tabular} \\
\hline $n = 128$ &
\begin{tabular}{@{}l}CPU~Time \\ Mem. All. \end{tabular} &
\begin{tabular}{c} 71H 33' 31" \\ -- \end{tabular} &
\begin{tabular}{c} 8H 17' 38" \\ 6,~552 Mb  \end{tabular} &
\begin{tabular}{c}  11H 22' 53" \\ 6,~814 Mb \end{tabular} &
\begin{tabular}{c} $\approx$ 170H  \\ OUT$^*$ \end{tabular} &
\begin{tabular}{c} $\approx$ 7H  \\ OUT$^*$ \end{tabular} &
\begin{tabular}{c}  -- \\ --  \end{tabular} \\
\hline
\end{tabular} \\[1mm]
{\small \dag ~~La proc\'edure Maple linalg[charpoly].
\\ \ddag ~~La proc\'edure
correspondant \`a la M\'ethode de Jordan-Bareiss modifi\'ee.
\\ $^*$ ~OUT signifie
qu'apr\`es un temps (\gui{CPU~Time}) plus ou moins long de calcul, un message
d'erreur \\ \gui{Out of memory} appara\^{\i}t lorsque la m\'emoire allou\'ee
(\gui{Mem. All.}) a d\'epass\'e le seuil de 350 Mbytes.}
 }\end{minipage}}
}
 
\newpage
 
\vspace*{2cm}
 
\rotatebox{90}{ {\begin{minipage}[t]{17cm}
{
{\Large {\bf Deuxi\`eme~Groupe~:~Mathard$(n,x,y)$}}\\[4mm]
\begin{tabular}{|p{3.5 em}|p{5 em}|c|c|c|c|c|}
\hline \multicolumn{2}{|c|}{matrice} & \multicolumn{1}{|c|} {linalpoly} & berkosam &
Chistov & Faddeev & barmodif  \\ \hline $n=10$ & \begin{tabular}{@{}l}CPU~Time \\
Mem. All.
\end{tabular} &
\begin{tabular}{l} 48' 35" \\ 44,~359 Mb \end{tabular} &
\begin{tabular}{l} 4' 04" \\ 8,~059 Mb \end{tabular} &
\begin{tabular}{l} 13' 17" \\ 10,~352 Mb \end{tabular} &
\begin{tabular}{l} 53' 31" \\ 47,~898 Mb \end{tabular} &
\begin{tabular}{c} $\approx$ 17' \\ OUT \end{tabular} \\
\hline $n=12$ & \begin{tabular}{@{}l}CPU~Time \\ Mem. All.
\end{tabular} &
\begin{tabular}{c} $\approx$ 1H 25' \\ OUT \end{tabular} &
\begin{tabular}{l} 15' 30" \\ 13,~891 Mb \end{tabular} &
\begin{tabular}{l} 54' 18" \\ 13,~563 Mb \end{tabular} &
\begin{tabular}{c} $\approx$ 1H 10' \\ OUT \end{tabular} &
\begin{tabular}{c}  -- \\ -- \end{tabular} \\
\hline $n=15$ & \begin{tabular}{@{}l}CPU~Time \\ Mem. All.
\end{tabular} &
\begin{tabular}{c} -- \\ -- \end{tabular} &
\begin{tabular}{l} 1H 23' 24" \\ 30,~403 Mb \end{tabular} &
\begin{tabular}{c} 5H 28' 28" \\ 26,~733 Mb \end{tabular} &
\begin{tabular}{c} -- \\ -- \end{tabular} &
\begin{tabular}{c}  -- \\ -- \end{tabular}  \\
\hline \hline \raisebox{1ex}{$\,n=10$} \raisebox{1ex}{$\& \, y=1$}
&
\begin{tabular}{@{}l} CPU~Time \\ Mem. All. \end{tabular} &
\begin{tabular}{l} 0' 47" \\ 4,~062 Mb \end{tabular} &
\begin{tabular}{l} 0' 05" \\ 2,~031 Mb \end{tabular} &
\begin{tabular}{l} 0' 10" \\ 2,~227 Mb \end{tabular} &
\begin{tabular}{l} 0' 27" \\ 3,~276 Mb \end{tabular} &
\begin{tabular}{l} 1' 13" \\ 2,~162 Mb \end{tabular} \\
\hline \raisebox{1ex}{$\,n=15$} \raisebox{1ex}{$\& \, y=1$} &
\begin{tabular}{@{}l} CPU~Time \\ Mem. All. \end{tabular} &
\begin{tabular}{c} $\approx$ 2H 10' \\ OUT  \end{tabular} &
\begin{tabular}{l} 0' 42" \\ 3,~341 Mb \end{tabular} &
\begin{tabular}{l} 1' 33" \\ 3,~603 Mb \end{tabular} &
\begin{tabular}{l} 5' 33" \\ 14,~480 Mb \end{tabular} &
\begin{tabular}{l} 22' 06" \\ 3,~800 Mb \end{tabular}  \\
\hline \raisebox{1ex}{$\,n=20$} \raisebox{1ex}{$\& \, y=1$} &
\begin{tabular}{@{}l} CPU~Time \\ Mem. All. \end{tabular} &
\begin{tabular}{c} -- \\ --  \end{tabular} &
\begin{tabular}{l} 3' 26" \\ 5,~110 Mb \end{tabular} &
\begin{tabular}{l} 7' 14" \\ 4,~979 Mb  \end{tabular} &
\begin{tabular}{l} 44' 30" \\ 47,~111 Mb \end{tabular} &
\begin{tabular}{l} 3H 24' 47" \\ 5,~831 Mb \end{tabular} \\
\hline \raisebox{1ex}{$\,n=25$} \raisebox{1ex}{$\& \, y=1$} &
\begin{tabular}{@{}l} CPU~Time \\ Mem. All. \end{tabular} &
\begin{tabular}{c} -- \\ -- \end{tabular} &
\begin{tabular}{l} 11' 56"\\ 8,~583 Mb \end{tabular} &
\begin{tabular}{c} 24' 52" \\ 10,~025 Mb \end{tabular} &
\begin{tabular}{c} $\approx$ 1H 10' \\ OUT \end{tabular} &
\begin{tabular}{c} 23H 54' 55"  \\ 9,~107 Mb \end{tabular} \\
\hline
\end{tabular}
}\end{minipage}}
}
 
 \newpage
 
\vspace*{1.5cm}
 
\rotatebox{90}{ {\begin{minipage}[t]{17cm}
{
{\Large {\bf Troisi\`eme~Groupe~:~Matmod$(n,lisvar,Ideal,p)$}}
\\[3mm]
{\Large Pour $~lisvar=[x],~~Ideal=[x^3-1]~~et~~p=7~:$}
\\[2mm]
\begin{tabular}{|p{4 em}|p{5 em}|c|c|c|c|c|}
\hline \multicolumn{2}{|c|}{matrice} & \multicolumn{1}{|c|}
 {linalpoly} & berkomod & Chistov & Faddeev & barmodif  \\
\hline $n=10$ &
\begin{tabular}{@{}l}CPU~Time \\ Mem. All. \end{tabular} &
\begin{tabular}{l} 0' 14" \\ 2,~162 Mb \end{tabular} &
\begin{tabular}{l} 0' 02" \\ 1,~703 Mb \end{tabular} &
\begin{tabular}{l} 0' 04" \\ 1,~769 Mb \end{tabular} &
\begin{tabular}{l} * \\ * \end{tabular} &
\begin{tabular}{l} 0' 43" \\ 1,~900 Mb \end{tabular} \\
\hline $n=12$ &
\begin{tabular}{@{}l}CPU~Time \\ Mem. All. \end{tabular} &
\begin{tabular}{l} 2' 03" \\ 8,~059 Mb \end{tabular} &
\begin{tabular}{l} 0' 04" \\ 1,~769 Mb \end{tabular} &
\begin{tabular}{l} 0' 07" \\ 1,~769 Mb \end{tabular} &
\begin{tabular}{l} * \\ *  \end{tabular} &
\begin{tabular}{l} 1' 52" \\ 2,~031 Mb \end{tabular} \\
\hline $n=16$ &
\begin{tabular}{@{}l}CPU~Time \\ Mem. All. \end{tabular} &
\begin{tabular}{c} $\approx$ 1H 30' \\ OUT \end{tabular} &
\begin{tabular}{l} 0' 12" \\ 1,~834 Mb \end{tabular}  &
\begin{tabular}{l} 0' 19" \\ 1,~900 Mb \end{tabular} &
\begin{tabular}{l} * \\ * \end{tabular} &
\begin{tabular}{l} 8' 31" \\ 2,~424 Mb \end{tabular} \\
\hline
\end{tabular}
\\[5mm] {\Large Pour $~lisvar=[x,y]~;~Ideal
= [H,L]\,^\dagger~;~et~~p = 17~:$} \\[2mm]
\begin{tabular}{|p{4 em}|p{5 em}|c|c|c|c|c|}
\hline \multicolumn{2}{|c|}{matrice} & \multicolumn{1}{|c|}
 {linalpoly} & berkomod & Chistov & Faddeev & barmodif  \\
\hline $n=10$ &
\begin{tabular}{@{}l}CPU~Time \\ Mem. All. \end{tabular} &
\begin{tabular}{l} 38' 13" \\ 37,~872 Mb \end{tabular} &
\begin{tabular}{l} 0' 28" \\ 2,~162 Mb \end{tabular} &
\begin{tabular}{l} 0' 50" \\ 3,~014 Mb \end{tabular} &
\begin{tabular}{l} 1' 25" \\ 1,~638 Mb \end{tabular} &
\begin{tabular}{l} 1H 06' 01" \\ 9,~369 Mb \end{tabular} \\
\hline $n=12$ &
\begin{tabular}{@{}l}CPU~Time \\ Mem. All. \end{tabular} &
\begin{tabular}{c} $\approx$ 2H \\ OUT \end{tabular} &
\begin{tabular}{l} 0' 50" \\ 2,~293 Mb \end{tabular} &
\begin{tabular}{l} 1' 33" \\ 3,~996 Mb \end{tabular} &
\begin{tabular}{l} 2' 43" \\ 2,~096 Mb \end{tabular} &
\begin{tabular}{l} 5H 59' 23" \\ 14,~087 Mb \end{tabular} \\
\hline $n=16$ &
\begin{tabular}{@{}l}CPU~Time \\ Mem. All. \end{tabular} &
\begin{tabular}{c} -- \\ -- \end{tabular} &
\begin{tabular}{l} 2' 11" \\ 2,~424 Mb \end{tabular} &
\begin{tabular}{c} 4' 25" \\ 7,~142 Mb \end{tabular} &
\begin{tabular}{l} 8' 03" \\ 3,~276 Mb \end{tabular} &
\begin{tabular}{c} $\approx$ 84H 40' \\  OUT  \end{tabular} \\
\hline
\end{tabular} \\[3mm]
\dag~~on a : \,$H=x^5-5xy+1$\, et \,$L=y^3-3y+1$.\\ $^*$~signifie que
Faddeev n'est pas applicable dans ce cas (puisque \,$p<n\,$). }
\end{minipage}} }
 
\newpage
 
\vspace*{1cm}
 
\rotatebox{90}{
{\begin{minipage}[t]{17cm}
{\vspace*{1cm}
{\Large {\bf Quatri\`eme~Groupe~:~Jou$(n,x)\,^\dagger$}} \\
{\large (Matrices de rang petit)} \\[5mm]
\begin{tabular}{|p{4 em}|p{5 em}|c|c|c|c|c|}
\hline \multicolumn{2}{|c|}{matrice} & \multicolumn{1}{|c|}{linalpoly} & berkosam &
Chistov & Faddeev & barmodif  \\ \hline $n=10$ &
\begin{tabular}{@{}l}CPU~Time \\ Mem. All. \end{tabular} &
\begin{tabular}{l} 0' 09" \\ 1,~310 Mb \end{tabular} &
\begin{tabular}{l} 0' 07" \\ 2,~489 Mb \end{tabular} &
\begin{tabular}{l} 0' 14" \\ 2,~555 Mb \end{tabular} &
\begin{tabular}{l} 0' 05" \\ 1,~244 Mb \end{tabular} &
\begin{tabular}{l} 0' 06" \\ 1,~310 Mb \end{tabular} \\
\hline
$n=15$ &
\begin{tabular}{@{}l}CPU~Time \\ Mem. All. \end{tabular} &
\begin{tabular}{l} 5' 10" \\ 5,~677 Mb \end{tabular} &
\begin{tabular}{l} 1' 00" \\ 3,~603 Mb \end{tabular} &
\begin{tabular}{l} 2' 08" \\ 3,~669 Mb \end{tabular} &
\begin{tabular}{l} 0' 17" \\ 1,~376 Mb \end{tabular} &
\begin{tabular}{l} 0' 24" \\ 1,~703 Mb \end{tabular} \\
\hline
$n=20$ &
\begin{tabular}{@{}l}CPU~Time \\ Mem. All. \end{tabular} &
\begin{tabular}{c} $\approx$  3H 30' \\  OUT \end{tabular} &
\begin{tabular}{l} 5' 09" \\ 5,~700 Mb \end{tabular} &
\begin{tabular}{l} 10' 55" \\ 5,~176 Mb  \end{tabular} &
\begin{tabular}{l} 0' 43" \\ 1,~834 Mb \end{tabular} &
\begin{tabular}{l} 1' 07" \\ 2,~293 Mb \end{tabular} \\
\hline
$n=25$ &
\begin{tabular}{@{}l}CPU~Time \\ Mem. All. \end{tabular} &
\begin{tabular}{c} $\approx$ 3H 30' \\ OUT \end{tabular} &
\begin{tabular}{l}  19' 01" \\ 8,~714 Mb  \end{tabular} &
\begin{tabular}{l}  44' 21" \\  9,~238 Mb \end{tabular} &
\begin{tabular}{l} 1' 39" \\ 2,~620 Mb \end{tabular} &
\begin{tabular}{l} 2' 32" \\ 2,~948 Mb \end{tabular} \\
\hline
\end{tabular} \\[5mm]
$^\dagger~{\rm {Jou}}(n,x)$ est une matrice \,$n\times n$\,
dont les coefficients \,$J_{ij}(x)\in {\bf Z}[x]$\,
sont donn\'es par la formule~: \\
$J_{ij}(x)=x^2(x-ij)^2+(x^2+j)(x+i)^2+x
~~{\rm pour}~1\leq i,j \leq n\, $. \\
Son rang est $\leq 3$ pour tout $x$ et pour
tout entier positif \,$n\,$. On remarque la sup\'eriorit\'e
des algorithmes de Faddeev et de Jordan-Bareiss modifi\'e
dans ce cas exceptionnel.
}\end{minipage}}
}

\newpage
 
\rotatebox{90}{
{\begin{minipage}[t]{17cm}
{\vspace*{15mm}
{\Large {\bf Cinqui\`eme Groupe~:~Matrices creuses
$\,n\times n\,$}} \\
{\large (Matrices \`a coefficients entiers)}\\[5mm]
\begin{tabular}{|p{4 em}|p{5 em}|c|c|c|}
\hline \multicolumn{2}{|c|}{matrice} & \multicolumn{1}{|c|} {linalpoly}
 & berkosam (creux) & Chistov (creux)  \\
\hline $n=32$ &
\begin{tabular}{@{}l}CPU~Time \\ Mem. All. \end{tabular} &
\begin{tabular}{l} 0' 59" \\ 2,~424 Mb \end{tabular} &
\begin{tabular}{l} 0' 05" \\ 1,~834 Mb \end{tabular} &
\begin{tabular}{l} 0' 35" \\ 1,~769 Mb \end{tabular} \\
\hline
$n=50$ &
\begin{tabular}{@{}l}CPU~Time \\ Mem. All. \end{tabular} &
\begin{tabular}{l} 5' 32" \\  2,~555 Mb \end{tabular} &
\begin{tabular}{l} 0' 21" \\ 1,~965 Mb \end{tabular} &
\begin{tabular}{l} 3' 55" \\ 1,~834 Mb \end{tabular} \\
\hline
$n=64$ &
\begin{tabular}{@{}l}CPU~Time \\ Mem. All. \end{tabular} &
\begin{tabular}{l} 14' 46" \\ 2,~620 Mb \end{tabular} &
\begin{tabular}{l} 0' 42" \\ 2,~031 Mb \end{tabular} &
\begin{tabular}{l} 12' 00" \\ 2,~031 Mb \end{tabular} \\
\hline
$n=128$ &
\begin{tabular}{@{}l}CPU~Time \\ Mem. All. \end{tabular} &
\begin{tabular}{l}  4H 29' 43" \\ 9,~173 Mb \end{tabular} &
\begin{tabular}{l} 6' 27" \\  2,~489 Mb  \end{tabular} &
\begin{tabular}{l} 6H 27' 19" \\ 3,~603 Mb \end{tabular} \\
\hline
$n=200$ &
\begin{tabular}{@{}l}CPU~Time \\ Mem. All. \end{tabular} &
\begin{tabular}{l} $> 16$ heures~: \\ Calcul stopp\'e  \end{tabular} &
\begin{tabular}{l} 29' 20" \\ 2,~555 Mb \end{tabular} &
\begin{tabular}{l}  Calcul stopp\'e \\ Calcul stopp\'e \end{tabular} \\
\hline
\end{tabular}
}\end{minipage}}
}
 \newpage \thispagestyle{empty}


\chapter[Le \deter et les \exas]{Le \deter et les \exas}
\acvide
\minitoc
\label{chap DeterUniv}

\subsubsection*{Introduction}   

Ce chapitre et le suivant donnent quelques aper\c{c}us sur le travail de 
Valiant (notamment \cite{Val1,Val2,Val3}) dans lequel il d\'ecrit un 
analogue \agq de la conjecture $\P\neq \NP$.
Notre expos\'e doit beaucoup au survey de von zur Gathen \cite{Gathen} 
et au livre de B\"{u}rgisser \cite{Bur2}. Une autre r\'ef\'erence
classique est le livre de B\"{u}rgisser, Clausen et  Shokrollahi 
\cite{Bur}.

\ss
Dans la section \ref{secForCirDes} nous discutons  diff\'erents 
codages possibles pour un \pol sur un anneau. La section \ref{secBrent} 
est consacr\'ee pour l'essentiel \`a la \met de Brent pour la 
\paran des \exasz. Dans la section \ref{secPlupDif} 
nous montrons pourquoi la plupart des \pols sont difficiles \`a 
calculer.
Enfin la section \ref{secDeterUniv} expose le r\'esultat de Valiant sur 
le caract\`ere universel du \deterz.

\section{Expressions, circuits et descriptions}
\label{secForCirDes}
Nous nous int\'eressons dans cette section \`a diff\'erentes approches 
concernant le codage d'un \pol arbitraire sur un \acom 
\,$\A$ (le codage des \elts de \,$\A$ est suppos\'e fix\'e). 

\ss Une premi\`ere mani\`ere de coder un \pol est de donner son degr\'e
total, les noms de ses variables et la liste de ses \coesz, 
dans un ordre convenu. 
C'est ce que nous avons appel\'e la \emph{\rpn 
dense}  des \polsz. Il est raisonnable de 
penser que pour l'immense majorit\'e des \pols il n'y a rien de mieux 
\`a faire, et nous donnerons un r\'esultat dans cette direction (voir 
le \tho \ref{thPlupDif}).

Certains \pols tr\`es utilis\'es  ont relativement peu de \coes non 
nuls. On peut choisir pour leur codage une
\emph{\rpn creuse}, dans laquelle on 
donne la liste des couples (\coe non nul, mon\^ome) effectivement 
pr\'esents dans le \polz, chaque mon\^ome \'etant cod\'e lui-m\^eme par 
la liste des exposants de chaque variable, \'ecrits en binaire.
Par exemple le \pol \,$\alpha X^{64}Y+\beta XY^{33}Z^4$\,
sera cod\'e par \,$[[a,[1000000,1,0]],[b,[1,100001,100]]]$, o\`u
\,$a$\,  et \,$b$\,  d\'esignent des codes pour \,$\alpha $\, et
\,$\beta$.

La \emph{taille \bolez} d'une \rpn creuse ou dense est la 
lon\-gueur du mot qui code le \polz. La taille peut \egmt \^etre 
ap\-pr\'e\-ci\'ee d'un point de vue purement \agqz, auquel cas chaque 
constante et chaque variable a 
conven\-tion\-nel\-lement la longueur $1$. Le point faible de la \rpn 
creuse est  que le produit d'un petit nombre de \pols creux est un \pol 
dense comme le montre l'exemple classique suivant:
\begin{equation} \label{Eq12.1}
(1+X)\times (1+X^{2})\times (1+X^{4})\times 
\cdots\times(1+X^{2^n})=\sum_{k=0}^{2^{n+1}-1}X^k
\end{equation}

\ss Un autre codage naturel est l'utilisation des \emph{\exasz}. Une 
\exa \index{expression!arithm\'etique} est un mot bien form\'e qui 
utilise comme ingr\'edients de base les \elts de \,$\A$ et les 
symboles de variables d'une part, les symboles \,$+,\; \times$\, d'autre 
part, et enfin les parenth\`eses ouvrante et fermante. 
D'un point de vue un peu plus abstrait, une \exn est vue comme \emph{un 
arbre \'etiquet\'e}. 
Aux feuilles de l'arbre, il y a des \elts de \,$\A$ 
(les constantes) et des symboles de variables,  
chaque n{\oe}ud est \'etiquet\'e par \,$+$\, ou \,$\times$.
En outre  deux branches partent exactement de chaque n{\oe}ud.  
La racine de l'arbre repr\'esente l'\exaz.
\begin{figure}[hbtp] 
\begin{center}
\includegraphics*[width=10cm]{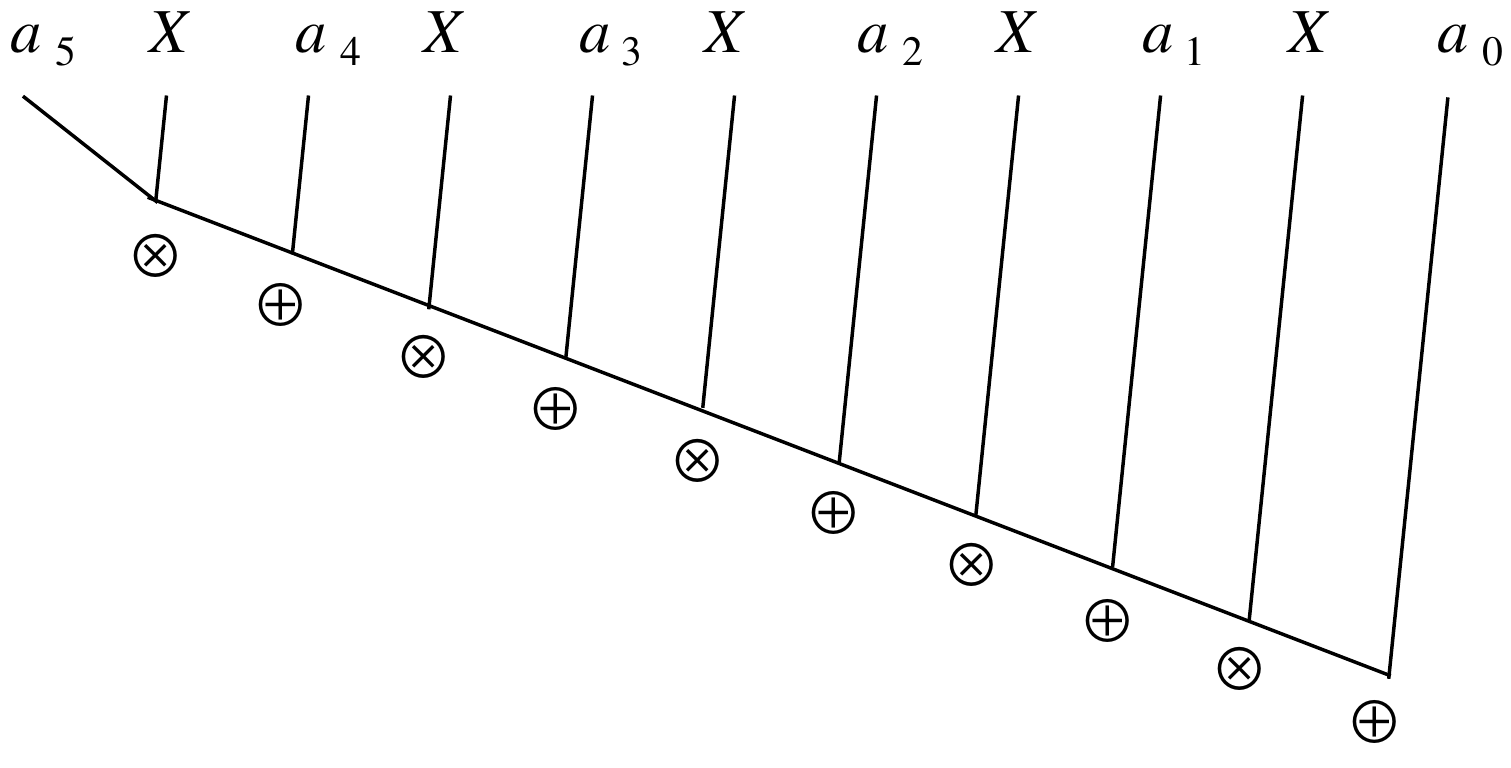}
\end{center}
\caption[L'arbre de Horner]
{\label{ffhorner}  L'arbre de l'\exn de Horner}  
\end{figure}  
 
La \emph{taille d'une \exnz}\index{taille!d'une expression 
artithm\'etique}  
peut \^etre appr\'eci\'ee d'un point de vue purement \agqz, on 
prend alors le nombre de n{\oe}uds dans l'arbre, sans compter les 
feuilles (la taille est alors \'egale au nombre de feuilles moins $1$). 
Si on adopte un point de 
vue proprement informatique,  il faut prendre en compte pour la 
\emph{taille \bolez} la longueur de l'\'ecriture explicite de  
l'\exn dans un langage pr\'ecis, 
o\`u les constantes et les variables ont des codes. 
M\^eme si on ne travaille qu'avec un nombre fini de constantes, la 
taille \bole de l'\exn ne peut \^etre consid\'er\'ee comme 
simplement proportionnelle \`a sa taille \agqz, ceci parce que 
l'ensemble des variables n'est pas born\'e a priori.

La \rpn dense peut naturellement \^etre vue comme une \rpn par 
\exns\index{repr\'esentation!par expressions} dans laquelle seules sont 
autoris\'ees des \'e\-cri\-tu\-res canoniques. Le nombre de \coes d'un 
\pol de degr\'e $d$ en $n$ variables est \'egal \`a 
$d+n \choose n$. 
La \rpn par \exa permet d'exprimer certains \pols (une petite 
minorit\'e, mais ce sont les \pols les plus utilis\'es) sous forme plus 
compacte, et plus efficace en ce qui concerne leur \evaz. Donnons 
en trois exemples. 

Le premier est celui de la \rpn \`a la Horner d'un \pol en une variable.
Dans les deux \'ecritures ci-dessous
\begin{equation} \label{Eq12.2}
\formule{a_5X^5+\cdots+a_1X+a_0 =\\
~~~~~~~~~~a_0+X \,(a_1+X \,(a_2+X \,(a_3+X \,(a_4+X \, a_5))))}
\end{equation}
l'expression dense r\'eclame pour son \eva $15$ \muls et  
l'\exn de Horner (dans le second membre) en
r\'eclame seulement $5$. \index{Horner!expression de}
En degr\'e \,$d$\, on obtient 
respectivement ${d+1 \choose 2}+d$ et \,$2d$\, \oparis 
respectivement pour l'\exn 
d\'evelopp\'ee et l'\exn de Horner.

Le deuxi\`eme exemple est celui d'un produit it\'er\'e. L'\exn
ci-dessous, qui est de taille $2n-1$ 
\begin{equation} \label{Eq12.3}
(X_1+Y_1)\times (X_2+Y_2)\times\cdots \times (X_n+Y_n)
\end{equation}
s'\'ecrit comme une somme de $2^n$ mon\^{o}mes, et a une taille de 
l'ordre de $n\,2^n$  en \rpn creuse (et plus grande encore en \rpn 
dense).

Le troisi\`eme exemple, sur lequel nous reviendrons plus en d\'etail est 
celui du \deter d'une \maca dont les entr\'ees sont $n^2$  variables 
ind\'ependantes. On ne sait pas si cette \fam de \pols peut \^etre ou 
non repr\'esent\'ee par une \fam d'\exns 
\emph{de taille \pollez}, \cad dont la taille  serait major\'ee par un 
$C\,n^k\leq C\,2^{k\,\log{n}}$ (avec $C$  et $k$ fix\'es). On 
conjecture 
que c'est faux. Par contre nous verrons que le \deter peut \^etre  
repr\'esent\'e par une \exn \emph{de taille 
\qplez}\index{quasi-polynomial}, \cad  major\'ee par 
un 
\,$C\,2^{(\log{n})^k}$ (avec \,$C$\,  et \,$k$\, fix\'es). Il est clair 
qu'en \rpn 
dense comme en \rpn creuse, le \deter a une taille 
$\geq n!\geq 2^n$ (pour $n\geq 5$) donc asymptotiquement beaucoup plus 
grande que \,$C\,2^{(\log{n})^k}$.

Notez par contre que la \fam de \pols de l'exemple $(\ref{Eq12.1})$ 
occupe une 
taille exponentielle en \rpn par \exasz, \`a cause du $X^{2^n}$: le 
degr\'e d'un \pol ne peut pas \^etre plus grand que la taille d'une 
\exa qui l'exprime. 

\ss Un troisi\`eme codage naturel est celui que nous avons retenu pour 
l'ensemble de cet ouvrage, 
le codage par les \prevs \ariths ou, ce qui revient au m\^eme, par les 
\carisz. 

Une \exa peut \^etre vue comme un cas particulier de \cariz.
Sa taille en tant qu'\exa est la m\^{e}me que celle du \cari qui lui 
correspond, \cad est \'egale au nombre d'\oparis lors de l'ex\'ecution 
du \cirz. 
La \rpn creuse peut \egmt \^etre simul\'ee efficacement par un 
\cirz. 

Pour un \cirz, les param\`etres pertinents sont \`a la fois la taille et 
la \profz. Un \pol calcul\'e par un \cari de \prof \,$p$\, a un degr\'e 
major\'e par \,$2^p$\, et on est particul\`erement int\'eress\'e par les 
\fams de \pols \,$(P_n)$\, qui peuvent \^etre \'evalu\'es par des \fams 
de \cirs dont la \prof est un \,$\O(\log(\deg(P_n)))$. 
Il semble cependant tr\`es improbable que le \deter (comme \pol de 
degr\'e \,$n$\, \`a \,$n^2$\, variables), 
qui est dans la classe  $\SD(n^4,\log^2 n)$  puisse \^etre 
r\'ealis\'e dans une classe $\SD(n^k,\log{n})$ (pour un entier \,$k$). 

\begin{convention} 
\label{convFora} 
Dans les chapitres  \ref{chap DeterUniv}  et  \ref{chap Perma}
les \cirs et les \exas que nous consid\'ererons seront toujours 
\emph{sans division et sans soustraction}.  
Rappelons que l'\elid \`a la Strassen 
montre qu'il ne s'agit pas d'une restriction importante (surtout dans le 
cas des corps,
voir \thos \ref{thEliDiv} et   
\ref{prop Eli Div2}). La soustraction, quant \`a elle, est simul\'ee en 
deux op\'erations par $x-y=x+(-1)\times y$.
\end{convention}

\ss Un dernier codage naturel que nous envisagerons est celui dans 
lequel un \pol \,$P(x_1,\ldots,x_k)$\,  est obtenu sous la forme 
\begin{equation} \label{Eq12.4}
P(x_1,\ldots,x_k)\;=\;\sum_{e_1,\ldots,e_\ell\in 
\left\{0,1\right\}}R(x_1,\ldots,x_k,e_1,\ldots,e_\ell)
\end{equation}
le \pol $R$ \'etant lui-m\^eme donn\'e par un \cir ou une \exaz. 
Ceci peut sembler a priori artificiel, mais nous verrons dans le 
chapitre \ref{chap Perma} que cette \'ecriture condens\'ee des \pols est 
en rapport assez 
\'etroit avec la conjecture \,$\P\neq \NP$.

\ss Nous donnons maintenant quelques d\'efinitions qui 
r\'e\-sul\-tent de la discussion pr\'ec\'edente.

\begin{definition} 
\label{defForEvDes} 
Soit \,$(P_n)$\, une \fam de \pols (index\'ee par \,$n\in \N$\, ou 
\,$\N^\ell$) \`a \coes dans un anneau commutatif \,$\A$. 
Notons \,$v_n$\, et \,$d_n$\, le nombre de variables et le degr\'e de 
\,$P_n$.
\begin{itemize}
\item  Nous disons que la \fam \,$(P_n)$\, est 
\emph{\pborz}\index{p-bornee@$p$-born\'ee}
\index{p-bornee@$p$-born\'ee!famille de polyn\^omes} si \,$v_n$\, et 
\,$d_n$\, sont 
major\'es par un \pol en \,$n$. On dit encore qu'il s'agit d'une 
\emph{\pfa de \polsz.}\index{p-famille@$p$-famille!de polyn\^omes}
\item  Nous disons qu'une \fam \,$(\varphi_n)$\, d'\exas est 
\emph{\pborz}
\index{p-bornee@$p$-born\'ee!famille d'expressions arithm\'etiques} si 
la taille de  \,$\varphi_n$\, est major\'ee par un \pol en \,$n$.
\item  Nous disons que la \fam \,$(P_n)$\, est 
\emph{\pforz}\index{p-exprimable@$p$-exprimable} 
\index{p-exprimable@$p$-exprimable!famille de polyn\^{o}mes}
si elle est r\'ea\-li\-sable par une 
\fam $p$-born\'ee d'\exas \,$\varphi_n$\, (en particulier,  \,$(P_n)$\,  
est \pborz).  
\item  Nous disons qu'une \fam \,$(\gamma_n)$\, de \caris est 
\emph{\pbor en taille} si la taille de  \,$\gamma_n$\, 
est major\'ee par un \pol en \,$n$,
\index{p-bornee@$p$-born\'ee!famille de circuits arithm\'etiques} 
\emph{\pbor en degr\'es}  si les \pols 
\'evalu\'es \`a tous les noeuds de  \,$\gamma_n$\, sont major\'es par un 
\pol en \,$n$, \emph{\pborz} si la \fam est \pbor en taille et en 
degr\'es.
\item  Nous disons que la \fam \,$(P_n)$\, est 
\emph{\,$p$-\'evaluable}\index{p-evaluable@$p$-\'evaluable} (ou encore 
\emph{\pcalz})
\index{p-calculable@$p$-calculable} 
\index{p-calculable@$p$-calculable!famille de polyn\^{o}mes} 
si elle est r\'ea\-li\-sa\-ble par une \fam \pbor de \caris
(en particulier,  \,$(P_n)$\,  est \pborz).  
\item  Nous disons que la \fam \,$(P_n)$\, est 
\emph{$\,qp$-exprimable}\index{qp-exprimable@$qp$-exprimable} 
si c'est une \pfa  
r\'eali\-sa\-ble par une \fam d'\exas  dont la taille est 
\qple en \,$n$\, 
(\cad major\'ee par un \,$C\,2^{(\log{n})^k}$\, avec \,$C$\,  
et \,$k$\, fix\'es).
\item  Nous disons que la \fam \,$(P_n)$\, est 
\emph{$\,qp$-\'evaluable}\index{qp-evaluable@$qp$-\'evaluable} (ou 
encore \emph{\qcalz}\index{qp-calculable@$qp$-calculable})
si c'est une \pfa  r\'ealisable par une \fam de \caris
dont la taille est \qple en \,$n$.
\item  Nous disons qu'un \pol \,$R$\, en les variables 
\,$x_1,\ldots,x_k,\,y_{1},\alb\ldots,\alb y_{\ell}$\, \emph{est une 
description du \pol \,$P\,$}\index{description!d'un poly@d'un \polz}  
en les variables 
\,$x_1,\alb\ldots,\alb x_k\,$
si
\begin{equation} \label{eqDescr1}
P(\s{x})\;=\;\sum\nolimits_{\s{y}\in 
\left\{0,1\right\}^\ell}\;R(\s{x},\s{y})
\end{equation}
\item  Nous disons que la \fam \,$(P_n)$\, est 
\emph{\pdbz}\index{p-descriptible@$p$-descriptible}  
 s'il existe une \fam \pcal de \pols \,$(R_n)$, telle que chaque   
\,$R_n$\, est 
une description de \,$P_n$. 
\item  Nous disons que la \fam \,$(P_n)$\, est 
\emph{\pdb en expressions}\index{p-descriptible@$p$-descriptible!en 
expressions}  
 s'il existe une \fam \pfor de \pols \,$(R_n)$, telle que chaque   
\,$R_n$\, est une description de \,$P_n$.
\end{itemize}
\end{definition}

Il faut souligner que toutes les notions introduites ici sont \emph{non 
uniformes}, \cad qu'on ne demande pas que les \fams d'\exns ou de 
\cirs soient des \famus 
(cf. section \ref{sec unif}).

Nous utiliserons les notations suivantes pour d\'ecrire les classes (de 
\fams de \polsz) correspondant aux d\'efinitions pr\'ec\'edentes. Le 
\,${\cal V}$\, est mis pour Valiant, qui a \'etabli la plupart des 
concepts et des r\'esultats des chapitres 
\ref{chap DeterUniv} et \ref{chap Perma}.
\begin{notation} 
\label{notaForEvDes}~
\begin{itemize}
\item 
La classe des \fams de \pols  \pfors est not\'ee \,$\VPF$, celle des 
\fams  $qp$-exprimables 
\,$\VQPF$. \indexnota{VPf@$\VPF$}\indexnota{VQPf@$\VQPF$}
\item 
La classe des \fams de \pols  \pcals est not\'ee \,$\VP$, celle des 
\fams  \,$qp$-calculables 
\,$\VQP$. \indexnota{VQP@$\VQP$} \indexnota{VP@$\VP$}
\item 
La classe des \fams de \pols  \pdbs est not\'ee 
\,$\VNP$, celle des \fams  \pdbs en \exns 
\,$\VNPF$.\indexnota{VNP@$\VNP$}\indexnota{VNPf@$\VNPF$}
\item
La classe des \fams de \pols 
\'evaluables par des \fams \pbors
de \caris  de \prof  \,$\O(\log^k(n))$\, est not\'ee
$\VNC^k$. La r\'eunion des $\VNC^k$ est not\'ee 
$\VNC$. \indexnota{VNCk@$\VNC^k$}\indexnota{VNC@$\VNC$}
\end{itemize}
\end{notation}

Ces classes sont d\'efinies relativement \`a un anneau commutatif 
fix\'e $\A$. Si on a besoin de pr\'eciser l'anneau
on notera  \,$\VPF(\A)$,  \,$\VP(\A)$, etc\ldots 
La plupart des r\'esultats sont cependant ind\'ependants de
l'anneau. Les conjectures sont \'enonc\'ees en \gnl pour des 
corps. 

\begin{remark} 
\label{remfamcirpb}
\emph{Vue la proposition \ref{prop Eli Div4}, 
s'il existe une \fam \pbor en taille de \caris 
qui calcule une \pfa de \polsz, alors il existe aussi
une \fam \pbor de \caris qui calcule la m\^{e}me \fam de \polsz.
Pour la m\^{e}me raison nous aurions pu demander, 
pour d\'efinir la classe
$\VQP$, que la \fam de \caris soit non seulement \,$qp$-born\'ee 
en taille mais aussi \pbor en degr\'es.
} 
\end{remark}

\section{Parall\'elisation des \exns et des circuits} 
\label{secBrent}

\subsubsection*{Parall\'elisation des \exns}

Des \exns comme celles de Horner, qui sont optimales 
quant \`a leur taille (i.e. pour le temps
\sql d'\evaz), pr\'esentent un d\'efaut de 
\parasm criant. 
Brent a d\'e\-cou\-vert que n'importe quelle \exa peut \^etre 
remplac\'ee par un \cir ou par une \exn dont la \prof est 
\logq en la taille de l'\exn initiale.
 
\begin{theorem} 
\label{thBrent} \emph{(Brent \cite{Brent})} 
Pour tout \pol \,$P$\, la \prof \,$\pi$\,  du meilleur \cir et la taille  
\,$\tau$\, de la 
meilleure \exn sont reli\'es par
\begin{equation} \label{eqParaBrent}
\log (\tau+1)\leq \pi\leq\, {2\over \log \,3/2} \,\log(\tau +1)
\end{equation}
\end{theorem}
NB: On a 
${2\over \log\,3/2}=3,4190\ldots $. 
Un calcul plus pr\'ecis (th\'eor\`eme 21.35 dans \cite{Bur}) donne 
\,$\pi\leq\, {2\over\log\,\phi} \,\log(\tau)\alb +1$\, o\`u 
\,$\phi$\, est le nombre d'or \,${1+\sqrt{5}\over 2}$\, et \,${2\over 
\log \,\phi}=2,8808\ldots $

\sni \prv  Dans cette preuve nous notons \,$t(\varphi)$\, le nombre de 
feuilles de l'arbre correspondant \`a l'\exn \,$\varphi$\, (c'est la 
taille de l'\exn $+1$) et \,$\pi(\gamma)$\, la \prof d'un \cir ou 
d'une \exn \,$\gamma$.

 La premi\`ere \ine est facile. Si \,$P$\, est une variable ou 
une constante la \prof et la taille sont nulles. 
Sinon lorsque \,$P$\, est \'evalu\'e par un \cir \,$\gamma$\, on a 
\,$\gamma=\gamma_1\circ \gamma_2$\, (o\`u $\circ$ repr\'esente $+$ ou 
$\times$) et si on suppose avoir d\'ej\`a r\'e\'ecrit \,$\gamma_1$\, et 
\,$\gamma_2$\,  avec des \exns \,$\varphi_1$\, et \,$\varphi_2$\, on 
obtient \,$t(\varphi)=t(\varphi_1)+t(\varphi_2)$\, et 
\,$\pi(\gamma)=1+ \max(\pi(\gamma_1),\pi(\gamma_2))$. 
En fait ce calcul correspond \`a une \pcd qui d\'eploie le  \cir 
en une \exn de m\^eme \profz. Or l'arbre d'une \exn de \prof \,$p$\, a 
au plus \,$2^p$\, feuilles. 

La deuxi\`eme \ine est nettement plus subtile.
L'id\'ee est la suivante. Appelons \,$x_1,\ldots,x_m$\, les variables de 
l'\exn \,$\varphi$\, qui repr\'esente \,$P$. Nous voyons cette \exn 
comme un arbre. 
Si on consid\`ere un n{\oe}ud  \,$N$\, de l'arbre,  
il repr\'esente une sous-\exn \,$\alpha$ (voir figure 
\vref{fbrent}).

En rempla\c{c}ant cette sous-\exn (ce sous-arbre) par une nouvelle 
variable \,$y$\, (\cad une feuille), on obtient une \exn (un arbre) 
\,$\beta$\, qui repr\'esente un \pol \,$b_0+b_1y$\, 
\refstepcounter{bidon}\label{PreuveBrent} avec
\,$b_0,b_1\in \A\,[x_1,\ldots,x_m]$. 
Les \pols \,$b_0$\, et
 \,$b_1$\, correspondent \`a des arbres \,$\beta_0$\, et
 \,$\beta_1$\,  qu'il est facile de
construire \`a partir de l'arbre \,$\beta$.

\begin{figure}[hbtp] 
\begin{center}
\includegraphics*[width=10cm]{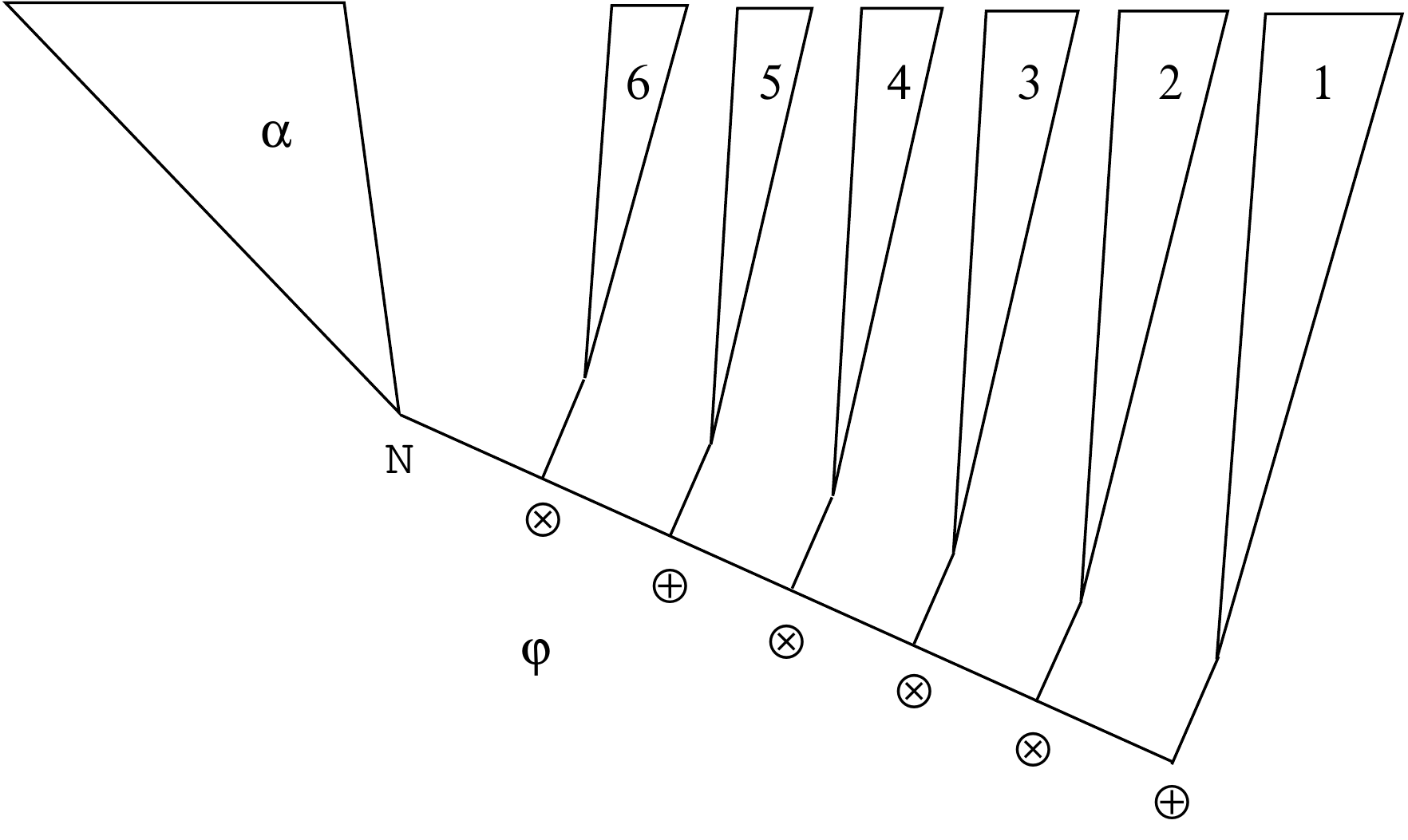}
\end{center}
\caption[Parall\'elisation d'une expression, 1]
{\label{fbrent} Parall\'elisation d'une expression, \`a la Brent.}  
\end{figure} 
\begin{figure}[hbtp] 
\begin{center}
\includegraphics*[width=12cm]{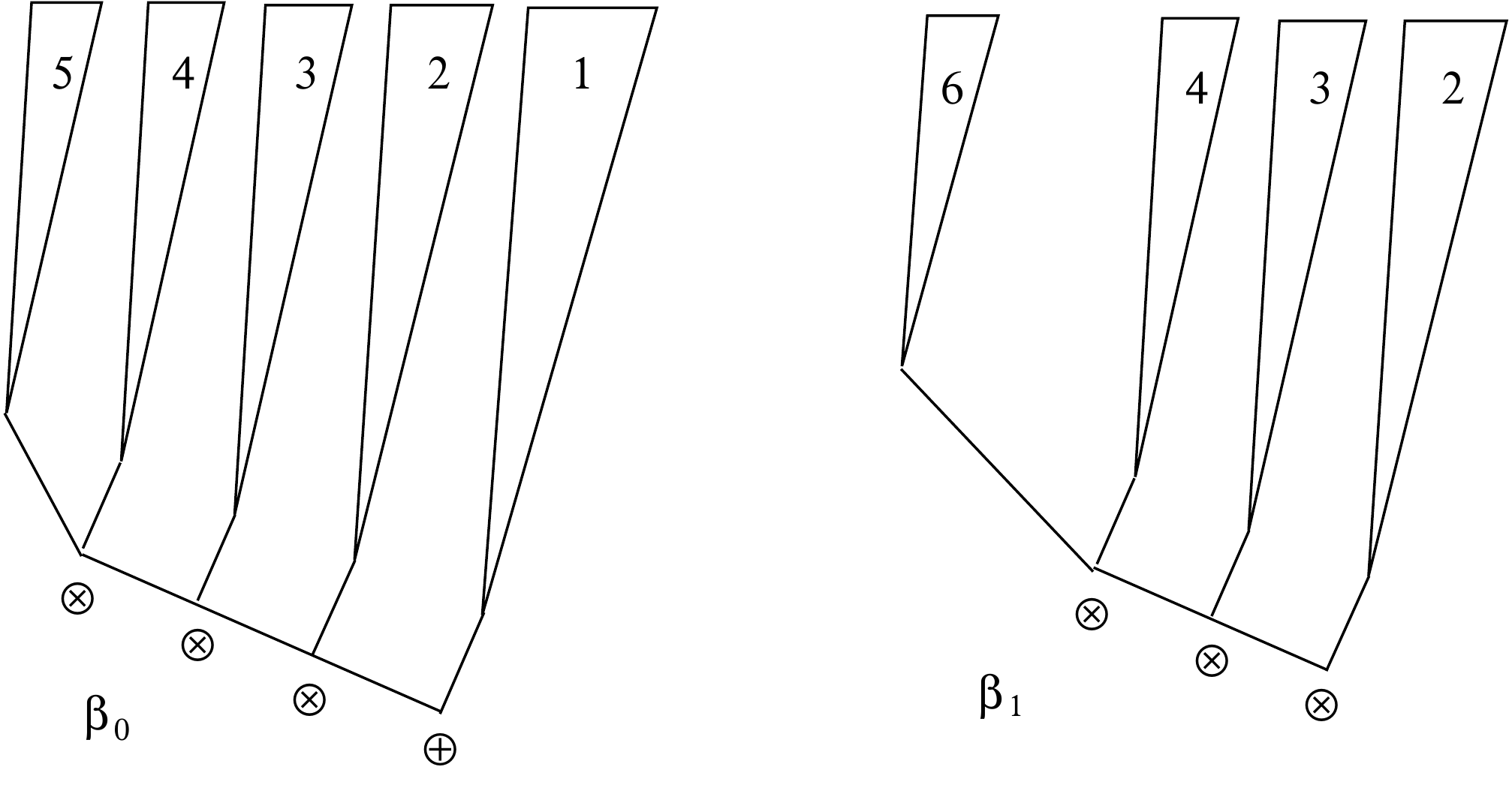}
\end{center}
\caption[Parall\'elisation d'une \exnz, 2]
{\label{fbrent2} Parall\'elisation d'une \exnz, \`a la Brent  
($\,\beta_0$\, et
 \,$\beta_1).$ }  
\end{figure}  
En effet (voir figure \vref{fbrent2}) pour \,$\beta_0$\, on substitue 
$0$ \`a \,$y$\, dans \,$\beta$, et 
on simplifie.
Pour \,$\beta_1$\, on part de la racine de \,$\beta$\, on suit le chemin 
jusqu'\`a  \,$y$\, et on supprime les n{\oe}uds \'etiquet\'es $+$, (et 
avec 
eux, la branche qui ne va pas \`a \,$y$).

On peut alors construire une \exn dans laquelle on met d'abord en 
\paral les \exns \,$\alpha$,  \,$\beta_0$\, et
 \,$\beta_1$\,  et o\`u on termine en calculant 
\,$\beta_0+(\beta_1\times \alpha)$.

La \prof \,$\pi(\gamma)$\, de cette \exn \,$\gamma$\, est 
major\'ee par  
\,$2+\max(\pi(\beta_0),\alb\pi(\beta_1),\alb\pi(\alpha))$.

Pour que cela soit efficace, il faut bien choisir 
le n{\oe}ud \,$N$\, (de mani\`ere 
que les tailles des trois \exns aient baiss\'e dans une proportion 
suffisante) et proc\'eder de mani\`ere r\'ecursive, 
\cad que chacune des $3$ \exns 
 est ensuite soumise de nouveau au m\^{e}me traitement (et ainsi de 
suite, cela va sans dire). 
Le choix du n{\oe}ud \,$N$\, se fait comme suit. Soit 
\,$t_0=t(\varphi)$. 
Si \,$t_0\le 4$\, on ne fait rien. Sinon on part de la racine de l'arbre 
et on choisit \`a chaque n{\oe}ud la branche la plus lourde 
\,$\varphi_k$. 
Si \,$t_k=t(\varphi_k)$\, on a donc \,$t_{k+1}\ge t_k/2$. On s'arr\^ete 
la derni\`ere fois que 
\,$t_k >(1/3) t_0$ 
(on aura fait un pas de trop lorsqu'on s'apercevra que le seuil a \'ete 
franchi, et il faudra retourner un cran en arri\`ere).
\begin{figure}[hbtp] 
\begin{center}
\includegraphics*[width=6cm]{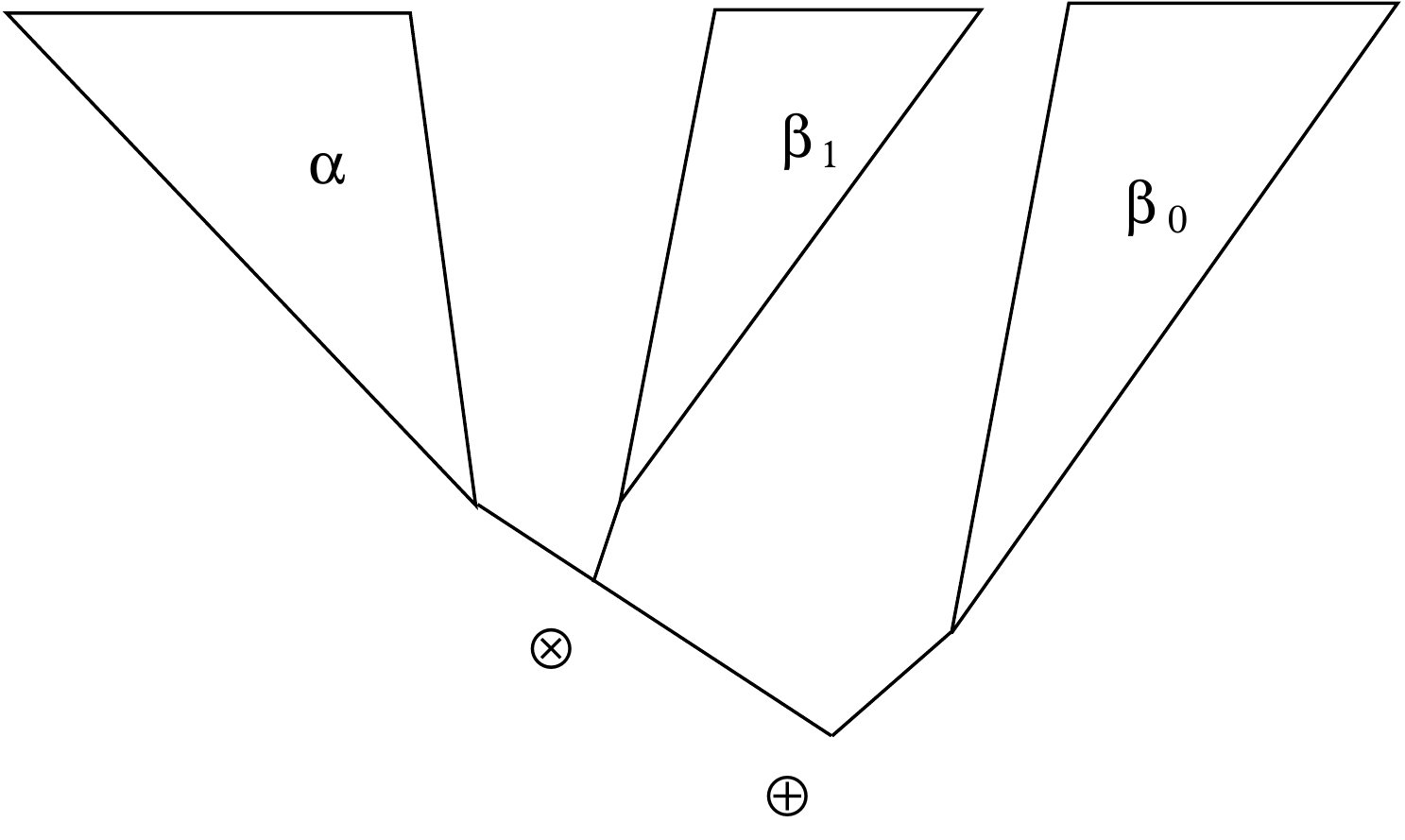}
\end{center}
\caption[Parall\'elisation d'une expression, 3]
{\label{fbrent3} Parall\'elisation d'une expression, \`a la Brent, 
\,$\beta_0+\beta_1\times \alpha$. }  
\end{figure}  

On a donc  \,$t_k>(1/3) t_0\ge t_{k+1}\ge t_k/2$, 
on en d\'eduit que $t_k$ et $t_0-t_k$ sont tous deux $\le (2/3) t_0$. 
Et on a \,$t(\alpha)=t_k$\, et 
\,$t(\beta_0),t(\beta_1)\le t_0-t_k$.
L'\ine voulue est donc \'etablie par \recu en 
v\'erifiant qu'elle fonctionne pour les \exns de taille $\le 3$.   
\qed

\ss Notez que les \pcds d\'ecrites sont uniformes.

\begin{remark} 
\label{remBrent}
\emph{Dans la premi\`ere \pcdz, on transforme un \cir
en une \exn de m\^eme \prof mais de taille peut-\^etre beaucoup plus 
grande. 
La seconde \pcd transforme toute \exn \gui{mal \'equilibr\'ee}
 de taille \,$\tau$\, en une
\exn \gui{bien \'equilibr\'ee} dont la taille \,$\tau'$\, n'a pas trop 
augment\'e{\footnote{~Le \tho donne 
\,$(1+\tau')\le (1+\tau)^{2\over \log (3/2)}$. En fait, lors d'une 
\'etape de \paran on a \,$t(\beta_0+\beta_1\times \alpha)\le 
(5/3)\,t(\varphi)$\, et cela conduit plus pr\'ecis\'ement \`a  
\,$(1+\tau')\le (1+\tau)^{\log(5/2)\over \log(3/2)}\le 
(1+\tau)^{2.26}$.}} et dont la \prof est devenue \logqz. 
Autrement dit la partie difficile du \tho de Brent fonctionne 
enti\`erement au niveau des \exnsz. 
} 
\end{remark}
\begin{figure}[hbtp] 
\begin{center}
\includegraphics*[width=10cm]{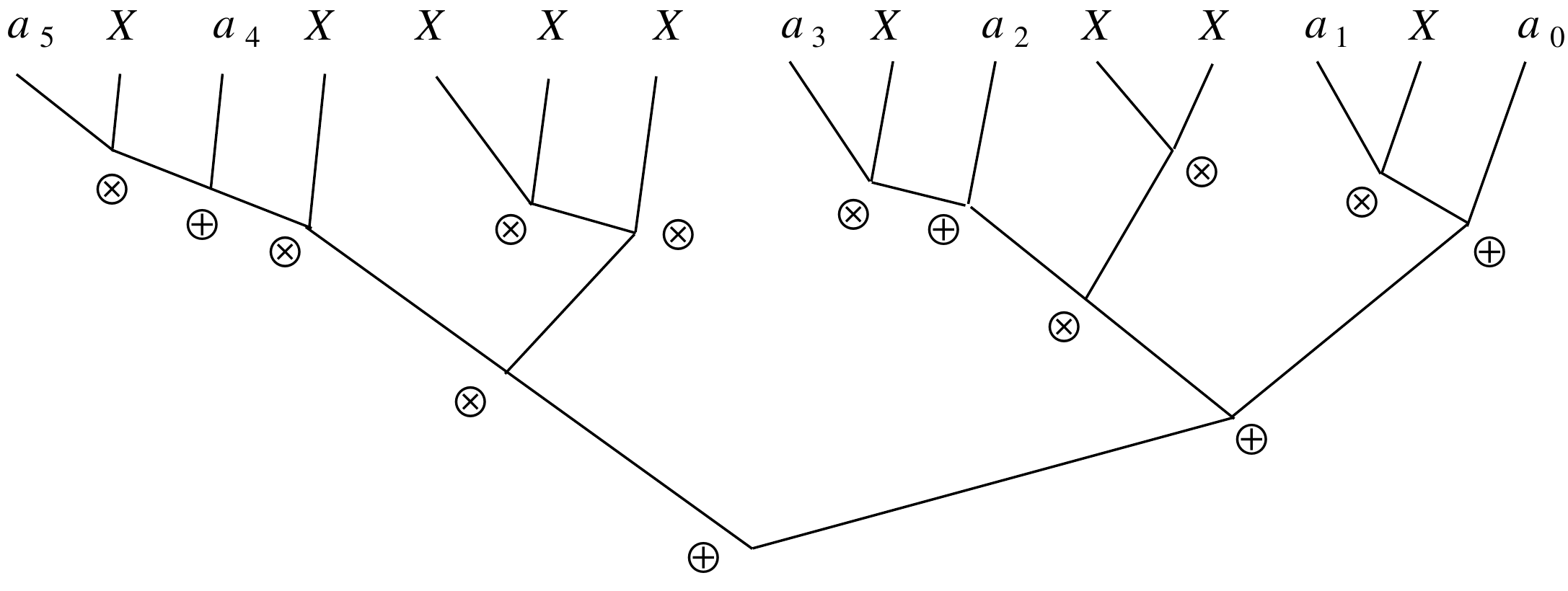}
\end{center}
\caption[L'arbre de Horner \paras]
{\label{ffhorner2} L'arbre de Horner, \paras \`a la Brent }
\index{Horner!parall\'elisation de l'\exnz}  
\end{figure}  
 
\ss Du \tho de Brent, on d\'eduit le corollaire suivant (la 
version uniforme serait \egmt valable):
\begin{corollary} 
\label{corBrent} On a $\VPF=\VNC^1$. 
\end{corollary}

On conjecture que par contre le \deter n'est pas r\'ealisable par une 
\exn de taille \pollez, et donc que \,$\VP\neq \VPF$.

La partie facile du \tho de Brent et l'\algo \paral de Berkovitz 
montrent \egmtz:
\begin{fact} 
\label{factdetfqp} 
Le \deter \,$n\times n$\, est r\'ealisable par une \exn de taille 
\qplez, en \,$2^{\O(\log^2n)}$.
\end{fact}

\subsubsection*{Parall\'elisation des \caris} 

Rappelons maintenant  le \tho \ref{vaskybra} (section
\ref{sec vonzur}) de
 Valiant \emph{et al.} \cite{Vasbra}  (voir aussi \cite{Hya}).
C'est en quelque sorte l'analogue pour les \caris du \tho
de Brent pour les \exasz.
Il donne une \pcd pour \parar n'importe quel \cari
\`a condition qu'il calcule un \pol de degr\'e raisonnable 
(le \cari purement \sql de taille \,$n$\, 
qui calcule \,$x^{2^n}$\, ne peut pas \^etre \parasz, 
mais c'est \`a cause de son degr\'e trop \'elev\'e).

\mni {\bf  
\Tho \ref{vaskybra}}    
\emph{Soit \,$\Gamma$\, un \cari sans division,  
de taille \,$\ell$, qui calcule un \poly \,$f$\, de  
degr\'e \,$d$\, en \,$n$\, variables
sur un anneau $\A$. Alors il existe  
un \carih \,$\Gamma'$\, de taille   
\,$\O(\ell^3d^6)$\, et de \prof  
\,$\O(\log({\ell\,d})\log{d})$\,  qui calcule \,$f$.}

\ms Ce \tho  implique \immt que
$\VP=\VNC^2$. Avec la partie facile du \tho de Brent il implique 
aussi qu'une \fam \,$qp$-born\'ee de \caris  peut \^etre \parase en une 
\fam de \cirs de \prof  \pog et donc en une \fam d'\exas 
de taille \qplez. En bref:  
\begin{corollary} 
\label{corBrentVaskybra} 
On a $\VP=\VNC^2=\VNC$ et $\VQP=\VQPF$. 
\end{corollary}

\begin{remark} 
\label{remVQP}
\emph{ 1) Ainsi $\VQP$ est la classe des
\fams de \pols r\'ealisables par des  \caris dont
le nombre de variables et les degr\'es sont \pbos et la
\prof est \pog (ce qui ne signifie pas pour autant qu'ils
soient dans $\VNC$).  Pour $\VQPF$ cela r\'esultait d\'ej\`a 
du \tho de Brent.\\
2) On conjecture a contrario que les inclusions  $\NC^2\subset 
\NC\subset \P$ sont strictes.
} 
\end{remark}

\section{La plupart des \pols sont difficiles \`a \'evaluer}  
\label{secPlupDif}
Pour \'etablir sous forme pr\'ecise ce qui est annonc\'e dans le titre 
de cette section 
nous avons besoin d'un r\'esultat de th\'eorie de l'\'elimination, dont 
la signification est intuitivement \'evidente. 
Si vous param\'etrez \gui{un 
objet g\'eom\'e\-tri\-que \,$S$} dans l'espace de dimension $3$ en 
donnant 
les $3$ cooordonn\'ees \,$x,y,z$\, comme fonctions \pollez s de deux 
\parats $u$ et $v$, l'objet que vous obtenez est en \gnl une 
surface, exceptionnellement une courbe ou plus exceptionnellement encore 
un point, mais jamais l'objet g\'eom\'e\-tri\-que ainsi cr\'e\'e ne 
remplira l'espace. Plus pr\'ecis\'ement, \`a partir des trois \pols 
\,$X(u,v)$, \,$Y(u,v)$, \,$Z(u,v)$\, qui param\`etrent l'objet 
\,$S$\, il est possible de calculer un \pol \,$Q$\, \`a trois variables 
non identiquement nul tel que \,$Q(X(u,v),Y(u,v),Z(u,v))$\, soit 
identiquement nul. 
Autrement dit tous les points de \,$S$\, sont sur la surface \agq  $S_1$ 
d'\'equation \,$Q(x,y,z)=0$. 
Donc, si le corps de base $\K$ est infini \gui{la plupart} des points de 
$\K^3$ sont en dehors de \,$S$. 
Pr\'ecis\'ement, consid\'erons un point \,$M$\, en dehors de \,$S_1$ 
(puisque le corps est infini, il en existe s\^urement), alors toute 
droite passant par \,$M$\, ne coupe \,$S_1$\, qu'en un nombre fini de 
points, major\'e par le degr\'e de \,$Q$.  Si le corps de base est celui 
des r\'eels ou celui des complexes, on en d\'eduit que le 
compl\'ementaire de \,$S_1$\, est un ouvert dense, ce qui donne encore 
une signification intuitive plus claire au
terme \gui{la plupart} utilis\'e dans la phrase 
ci-dessus. 

On peut montrer l'existence du \pol \,$Q$\, comme suit. Supposons  
les degr\'es de  \,$X$,  \,$Y$,  \,$Z$\, major\'es par   \,$d$. Pour   
\,$m$\, fix\'e, les \pols   \,$X^{m_1}Y^{m_2}Z^{m_3}$\, avec
\,$m_1+m_2+m_3\le m$\,  sont au nombre de \,${m+3 \choose 3}$, leur 
degr\'e est major\'e par \,$dm$, donc ils sont dans l'espace des \pols 
de degr\'e \,$\le dm$\,  en \,$2$\, variables, qui est de dimension 
\,${dm+2 \choose 2}$. Pour \,$m$\, assez grand, 
\,${m+3 \choose 3}>{dm+2 \choose 2}$, d'o\`u une relation de 
d\'ependance \lin non triviale entre les  
\,$X^{m_1}Y^{m_2}Z^{m_3}$, ce qui donne le \pol  \,$Q$. 

\noi Nous \'enon\c{c}ons maintenant le r\'esultat \gnlz, qui peut 
se d\'e\-mon\-trer de la m\^eme mani\`ere.

\begin{proposition} 
\label{propVarDim} Soit $\K$ un corps et $(P_i)_{1\le i\le n}$ une \fam 
de \pols en \,$m$\, variables \,$y_1,\ldots,y_m$\, avec \,$m<n$. Alors 
il existe un \pol non identiquement nul \,$Q(x_1,\ldots,x_n)$\, tel que 
\,$Q(P_1,\ldots,P_n)$\, est identiquement nul. 
En termes plus g\'eom\'etriques, l'image d'un espace $\K^m$ dans un 
espace $\K^n$ (avec \,$m<n$) par une application \polle est toujours 
contenue dans une hypersurface \agqz.
\end{proposition}

La proposition pr\'ec\'edente a la signification 
intuitive que, au moins en g\'eom\'etrie \agqz,
$\infty ^n>\infty ^m$ lorsque $n>m$.

Nous en d\'eduisons notre \thoz, dans lequel l'expresssion
\gui{la plupart} doit \^etre comprise au sens de la discussion qui 
pr\'ec\'edait la proposition \ref{propVarDim}.
\begin{theorem} 
\label{thPlupDif} Soit $\K$ un corps infini, \,$n$\, et \,$d$\, des 
entiers fix\'es. L'en\-sem\-ble des \pols de degr\'e \,$\le d$\, en 
\,$n$\, variables \,$x_1,\ldots,x_n$\, est un \evc  
\,$E(n,d)$\, sur $\K$ de dimension 
\,${n+d \choose d}$. 
Soit \,$t$\, une constante arbitraire fix\'ee. 
Notons \,$A(n,d,t)$\, la \fam  de tous les \caris qui repr\'esentent des 
\pols en \,$x_1,\ldots,x_n$\, de degr\'e 
\,$\le d$, avec au plus   \,${n+d \choose d}-1$\, 
constantes aux portes d'entr\'ees, et dont la taille est major\'ee par 
\,$t$. Alors \gui{la 
plupart} des \elts de  \,$E(n,d)$\, ne sont pas repr\'esent\'es 
par un \cari dans  \,$A(n,d,t)$.\\
En particulier, pour la plupart des \,$P\in E(n,d)$\, la taille
\,$\tau$\, du meilleur circuit admet la minoration 
\,$\tau\ge {n+d \choose d}$.
\end{theorem}

\sni \prv Chaque \cir dans  \,$A(n,d,t)$\, peut \^etre interpr\'et\'e 
comme calculant un \pol 
\,$P(c_1,\ldots,c_s;x_1,\ldots,x_n)\,$
dans lequel les \,$c_i$\,  sont les constantes du \cirz.
Le \pol \,$P$\,  en  \,$s+n$\,  variables correspondant fournit
(lorsqu'on fait varier les constantes)
une application de $\K^s$ vers \,$E(n,d)$\,  et chaque coordonn\'ee de 
cette application est une fonction \pollez.
Le fait de majorer la taille du \cir par  \,$t$\, implique que les \pols 
correspondants (en \,$s+n$\,  variables) sont en nombre fini. Finalement 
les \elts de  \,$E(n,d)$\, repr\'esent\'es par un \cari dans  
\,$A(n,d,t)$ sont contenus dans une r\'eunion finie d'hypersurfaces 
\agqs (d'apr\`es la proposition \ref{propVarDim}), qui est encore une 
hypersurface \agqz.
\qed

\ms On trouvera dans \cite{HS80,Stra74} des r\'esultats plus pr\'ecis 
sur ce sujet.
   
\begin{remark} 
\label{remThPlupDif}
\emph{ Notez a contrario, que le \cari qui exprime un \pol de
 \,$E(n,d)$\, directement comme somme de ses mon\^{o}mes 
utilise \,${n+d \choose d}$\,  constantes, et qu'il peut \^{e}tre
\'ecrit avec une 
taille  \,$\le 3\,{n+d \choose d}$. Tous les \,$x_1^{\mu_1}\cdots 
x_n^{\mu_n}$\, de degr\'e \,$\leq d$, qui sont au nombre de \,${n+d 
\choose d}$, peuvent en effet \^{e}tre calcul\'es en \,${n+d \choose d}-
n$, \'etapes (un produit de degr\'e \,$k>1$\, est calcul\'e en 
multipliant un   produit de degr\'e \,$k-1$, d\'ej\`a calcul\'e, par une 
variable). Il reste ensuite \`a multiplier chaque produit par une 
constante convenable, puis \`a faire l'addition.
} 
\end{remark}

\section{Le caract\`ere universel du d\'eterminant} 
\label{secDeterUniv}

Le but de cette section est de montrer que toute \exa  peut \^etre vue 
comme un cas particulier de 
l'\exn \gui{\deterz} dans laquelle les entr\'ees de la matrice ont 
simplement \'et\'e remplac\'ees par 
une des constantes ou une des variables de
l'\exnz, avec en outre le fait que le nombre de lignes de la \maca est
du m\^eme ordre de grandeur que la taille de l'\exnz.

Ceci n'est pas tr\`es surprenant, au vu de l'exemple classique 
ci-dessous, (inspir\'e de la matrice compagnon d'un \polz) dans lequel 
nous n'avons pas marqu\'e  les entr\'ees nulles:
$$ 
\det\left[\matrix{ 
 x   &      &      &     &   a_4   \cr 
 -1  &  x   &      &     &   a_3  \cr 
     &  -1  &  x   &     &   a_2   \cr 
     &      & -1   & x   &   a_1   \cr 
     &      &      & -1  &   a_0    
}\right]
= a_4+a_3x+a_2x^2+a_1x^3+a_0x^4\,.
$$

\subsubsection*{Projections} 

Nous introduisons maintenant formellement une notion pr\'ecise appel\'ee 
\emph{projection}\index{projection} pour d\'ecrire le processus de 
substitution extr\^emement limit\'e auquel nous allons avoir recours 
dans la suite.
\begin{definition} 
\label{defProjVal} Soit $\A$ un anneau commutatif fix\'e.
Soient \,$P\in \A\,[x_1,\alb\ldots,\alb x_k]$\, et
\,$Q\in \A\,[y_1,\ldots,y_\ell]$. Soient aussi \,$(P_n)$\, et
 \,$(Q_m)$\,  des \pfas de \pols \`a \coes dans $\A$.
\begin{itemize}
\item [$(1)$] On dit que \emph{$\,Q$\, est une 
\pro\index{projection!d'une expression} de 
\,$P\,$} 
si \,$Q$\, est obtenu \`a partir de  \,$P$\, en substituant \`a chaque 
\,$x_i$\, un \,$y_j$\, ou un \elt de $\A$.
\item [$(2)$] On dit que \emph{la \fam \,$(Q_m)$\, est une 
\ppro\index{p-projection@$p$-projection} de la \fam \,$(P_n)\,$}  s'il 
existe une 
fonction
\polt major\'ee \,$m\mapsto \varphi(m)$\, telle que, pour chaque \,$m$,  
\,$Q_m$\, est une \pro de \,$P_{\varphi(m)}\,$
\item [$(3)$] On dit que \emph{la \fam \,$(Q_m)$\, est une 
\qpro\index{qp-projection@$qp$-projection} de la \fam \,$(P_n)\,$}  s'il 
existe une 
fonction
quasi-poly\-nomia\-lement major\'ee \,$m\mapsto \varphi(m)$\, telle que, 
pour chaque \,$m$,  \,$Q_m$\, est une \pro de \,$P_{\varphi(m)}$.
\end{itemize}
\end{definition}

La proposition suivante est facile.
\begin{proposition} \label{propProj1}~
\begin{itemize}
\item [$(1)$] La compos\'ee de deux \pros est une \proz. M\^eme chose 
pour les \pprosz, ou pour les \qprosz.
\item [$(2)$] Les classes $\VP$ et $\VPF$ sont stables par 
\pproz.  
\item [$(3)$] La classe $\VQP=\VQPF$ est stable par 
\qproz.
\end{itemize}
\end{proposition}

\subsubsection*{R\'e\'ecriture d'une expression comme \deter} 

Dans le \tho de Valiant qui suit, la difficult\'e est de 
produire une matrice ayant pour \deter la somme des \deters de deux 
autres matrices. L'id\'ee est de faire cette construction non pas pour 
n'importe quelles matrices, mais en respectant un certain format.
C'est l'objet du lemme crucial qui suit.
Le format des matrices qui interviennent dans ce lemme  est visualis\'e 
ci-dessous sur un exemple avec \,$d=4.$
$$ 
\left[\matrix{ 
  a_{1} &  a_{2}   & a_{3}  &  a_{4}   &   0   \cr 
  1   &  x_{1,2}   & x_{1,3} &  x_{1,4}   &  b_1    \cr 
  0   &  1   & x_{2,3} &  x_{2,4}   &  b_2    \cr 
  0   &  0   & 1  &  x_{3,4}   &  b_3    \cr 
  0   &  0   & 0  &   1  &     b_4  
}\right]=
\left[\matrix{ 
 \alpha    &  0     \cr 
  T  &  \beta    
}\right]
$$
($\alpha$\, est un vecteur ligne, \,$\beta$\, un vecteur colonne et 
\,$T$\, est carr\'ee uni\tguz).
Notez que lorsqu'on d\'eveloppe le \deter d'une telle matrice sous la 
forme \,$\sum_{1\le i,j\le 4} a_i\,b_j\, c_{ij}$, le \pol 
\,$ c_{ij}$\, en facteur de \,$a_ib_j$\, est nul si la colonne de $a_i$ 
et la ligne de \,$b_j$\, se coupent dans la partie strictement 
sup\'erieure 
de la matrice \,$T$, puisque le cofacteur correspondant de \,$T$\, est 
nul.   
\begin{lemma} 
\label{lemValSom} Soient (pour \,$i=1,2$) deux entiers \,$d_i$, deux 
\macas
\,$T_i\in \A^{d_i\times d_i}$\, uni\tgusz, soient deux 
vecteurs lignes \,$\alpha_i\in\A^{1\times d_i}$\, 
et deux vecteurs colonnes 
\,$\beta_i\in\A^{d_i\times 1}$. 
Consid\'erons les trois matrices suivantes
$$ M_1 = 
\left[\matrix{ 
 \alpha_1   &  0     \cr 
  T_1  &  \beta_1   
}\right],\quad
M_2 = 
\left[\matrix{ 
 \alpha_2   &  0     \cr 
  T_2  &  \beta_2   
}\right],\quad 
 M = 
\left[\matrix{ 
  \alpha_1  & \alpha_2    &   0   \cr 
  T_1  &  0   & \beta_1     \cr
  0  & T_2    & \beta_2
}\right]\,.
$$
Alors on a 
$$ \det \,M = (-1)^{d_2}\, \det \,M_1 + (-1)^{d_1}\, \det\, M_2
$$
\end{lemma}
\prv Nous donnons seulement l'id\'ee directrice de cette preuve un peu 
technique.
Lorsqu'on d\'eveloppe compl\`etement le \deter comme indiqu\'e avant le 
lemme, le \pol en facteur d'un produit
\,$\alpha_{2,i}\beta_{1,j}$\, est nul car la ligne et la colonne 
correspondante se coupent dans la partie strictement sup\'erieure
de la matrice
\,$\left[\matrix{  T_1  &  0  \cr  0  & T_2 }\right]$\, 
(cf. le commentaire 
juste avant le lemme).  Pour voir que le \pol en facteur d'un produit 
\,$\alpha_{1,i}\,\beta_{2,j}$\, est nul \egmtz, 
il suffit de consid\'erer la matrice
$$ M'= \left[\matrix{ 
  \alpha_2  & \alpha_1    &   0   \cr 
  T_2  &  0   & \beta_2     \cr
  0  & T_1    & \beta_1
}\right]\,.
$$
Son \deter est identique (en tant qu'\exn d\'evelopp\'ee) \`a celui de 
\,$M$\, et l'argument pr\'ec\'edent s'applique. Il reste \`a 
consid\'erer, dans le \dev complet du \deter en somme de 
produits, les produits contenant un facteur 
\,$\alpha_{1,i}\,\beta_{1,j}$\, 
(et ceux contenant un facteur  \,$\alpha_{2,i}\beta_{2,j}$). 
Un exa\-men 
attentif montre que les seuls produits non nuls de ce type sont ceux qui 
empruntent la diagonale de \,$T_2$, donc on retrouve exactement les 
facteurs pr\'esents dans \,$\det\, M_1$\, au signe pr\`es. 
Ce signe correspond \`a une permutation 
circulaire des \,$d_2+1$\, derni\res colonnes de \,$M$.  
\qed

\begin{theorem} 
\label{thDetUni1} Toute \exn de taille \,$n$\, est la \pro
du \deter d'une matrice d'ordre inf\'erieur ou \'egal \`a \,$2n+2$.

\sni \emph{\underline{Pr\'ecisions} : la matrice est dans le format 
d\'ecrit au lemme pr\'ec\'edent, ses entr\'ees sont soit une constante 
de l'\exnz,
soit une variable de l'\exnz, soit $0$, $1$ ou $-1$, la derni\`ere 
colonne (resp. la derni\`ere ligne) ne contient que $0$, $1$ ou $-1$, 
une colonne quelconque  
contient au plus une variable ou une constante de l'\exnz.}
\end{theorem}
\begin{corollary} 
\label{corthDetUni1}~ 
\begin{itemize}
\item [$(1)$] Toute \fam \pfor est une \ppro de la famille \gui{\deterz} 
($\det_n$ est le \deter d'une \maca d'ordre $n$ donc un \pol de degr\'e 
$n$ en $n^2$ variables).
\item [$(2)$] $\VQP=\VQPF$ co\"{\i}ncide avec la classe des \fams qui 
sont des \qpros de la famille \deterz.
\end{itemize}
\end{corollary}
\prv On construit la matrice en suivant l'arbre de l'\exnz.
Pour une feuille \,$f$\, de l'arbre (constante ou variable) on prend la 
matrice  $2\times 2$ 
$$ \left[\matrix{ 
 f   &  0     \cr 
  1  &  1   
}\right]
$$
qui r\'epond bien aux sp\'ecifications souhait\'ees.
Supposons qu'on a construit les matrices \,$A_i$\, ($i=1,2$) qui ont 
pour \deters les \pols  \,$P_i$. Voyons d'abord la matrice pour 
\,$P_1+P_2$. Quitte \`a changer la derni\`ere colonne \,$\beta_i$\, en 
\,$-\beta_i$\, on peut aussi avoir les \deters oppos\'es, et on peut 
donc dans tous les cas appliquer le lemme \ref{lemValSom}. 
Donnons enfin la matrice \,$N$\, pour \,$P_1\times P_2$
$$ N=
\left[\matrix{ 
 A_1   &     0     \cr 
 J   &     A_2
}\right]
\quad \mathrm{avec }\quad  J=
\left[\matrix{ 
 0     & \cdots    &     0      & 1   \cr 
 0   &  \cdots         &     0      & 0   \cr
 \vdots&           &     \vdots  & \vdots   \cr
 0      &   \cdots&   0  & 0   \cr
}\right]\,.
$$
Cette matrice \,$N$\, r\'epond aux sp\'ecifications voulues, et comme 
elle est triangulaire par blocs, son \deter est \'egal \`a
\,$\det A_1\cdot\det A_2$.
\qed

\medskip Donnons par exemple la matrice construite comme indiqu\'e dans 
la preuve ci-dessus pour obtenir le d\'eterminant  \,$x+(2+y)z$\,
(nous n'avons pas mis les $0$):
$$ \left[\matrix{ 
 x   &  2  &  y  & . & . & .     \cr 
 1   &  .  &  .  & . & . & 1     \cr 
 .   &  1  &  .  &-1 & . & .     \cr 
 .   &  .  &  1  &-1 & . & .     \cr 
 .   &  .  &  .  & 1 & z & .     \cr 
 .   &  .  &  .  & . & 1 & -1      
}\right]\,.
$$

\section*{Conclusion}\label{sec Conclusion}

Dans le corollaire \ref{corthDetUni1} on a vu 
que toute \fam \pfor est une \ppro du \deter
et que toute \fam \qcal est une \qpro du \deterz. Cette derni\`ere 
propri\'et\'e s'\'enonce sous la forme suivante, qui ressemble \`a la 
$\NP$\,-\,compl\'etude.

\emph{La \fam \,$(\det_n)$\, est universelle pour $\VQP$ et les 
\qprosz.} 

Cependant le \deter lui-m\^eme n'est probablement pas 
\pfor et il est \pcalz, donc mieux que \qcalz.

Il se pose donc la question l\'egitime de trouver une \fam
\pfor qui soit universelle dans $\VPF$ par rapport aux 
\ppros et celle de trouver une \fam
\pcal qui soit universelle dans $\VP$ par rapport aux 
\ppros (le \deter serait un candidat naturel, mais pour le moment on ne 
conna\^{\i}t pas la r\'eponse \`a son sujet).
Le premier de ces deux probl\`emes a \'et\'e r\'esolu positivement par
Fich, von zur Gathen et Rackoff dans~\cite{FGR}. Le deuxi\`eme par
B\"urgisser dans \cite{Bur3}.

\ms La premi\`ere question admet une r\'eponse assez facile une fois 
connu le \tho de Brent. 
En effet toute \exn peut \^etre 
obtenue comme \pro d'une \exn de \prof comparable extr\^emement \parase 
qui combine syst\'ematiquement additions et \mulsz. Par exemple l'\exn 
(de \prof $2$) 
\,$\varphi_3=x_1+(x_2\times x_3)\,$
donne par \proz, au choix, l'une des deux \exns (de \prof $1$) 
\,$x_1+x_2$\, ou \,$x_2\times x_3$. Si maintenant on remplace chacun
des \,$x_i$\, par l'\exn \,$x_{i,1}+(x_{i,2}\times x_{i,3})\,$
on obtient une \exn \,$\varphi_9$\, \`a $9$ variables de \prof $4$ et
on voit que toute \exn de \prof $2$ est une \pro de \,$\varphi_9$. 
\begin{figure}[hbtp]   
\begin{center}
\includegraphics*[width=12cm]{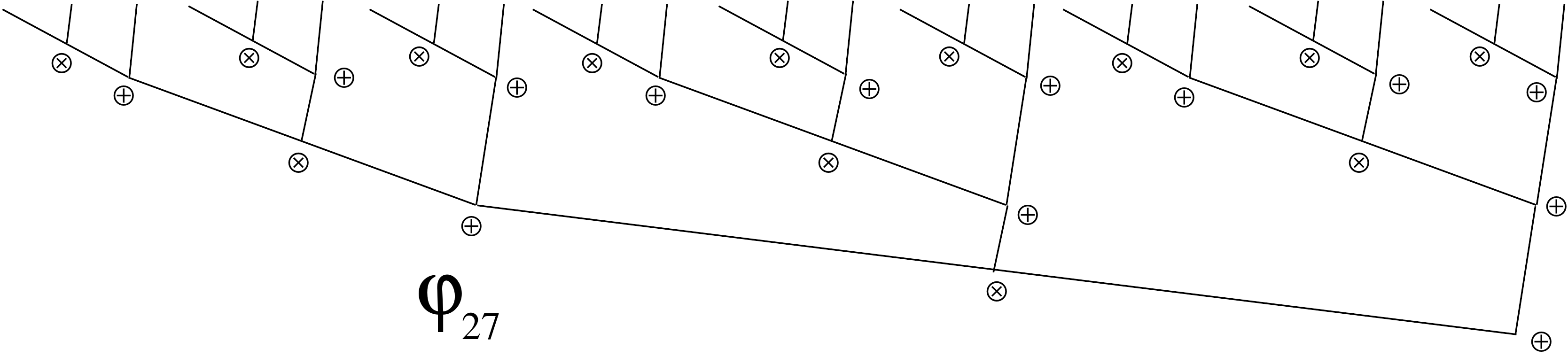}
\end{center}
\caption[Une \exn $p$-universelle]
{\label{fforuniv} Une \fam d'\exns $p$-universelle dans $\VPF$ }  
\end{figure}  
 
\noi En it\'erant le processus, toute \exn de \prof \,$n$\, est une 
\pro de l'\exn \,$\varphi_{3^n}$, qui est 
elle-m\^eme de \prof \,$2n$. 

Donc, apr\`es \paran \`a la Brent d'une \fam dans
$\VPF$, la \fam \parase est 
clairement une \ppro de
la \fam  \,$(\varphi_{3^n})$. 

Enfin la \fam \,$(\varphi_{3^n})$\, est elle-m\^eme 
 dans $\VPF$ car
\,$\varphi_{3^n}$\, est de taille \,$3^n-1$\, (on peut, si on a des 
scrupules, d\'efinir \,$\varphi_k$\, pour tout entier \,$k$\,
en posant \,$\varphi_k=\varphi_{3^\kappa}$\, o\`u 
\,$3^{\kappa-1}< k\leq 3^\kappa$).

\ms La deuxi\`eme question (trouver une \fam
\pcal qui soit universelle dans $\VP$ par rapport aux 
\pprosz) admet une r\'eponse po\-si\-ti\-ve du m\^eme style
(cf. \cite{Bur3}), mais 
nettement plus embrouill\'ee.

\newpage \thispagestyle{empty}

\chapter{Le permanent et la conjecture \texorpdfstring{$\P\not= 
\NP$}{P diff\'erent de NP}}  
\label{chap Perma}
\minitoc

\subsubsection*{Introduction } 

Ce chapitre est d\'edi\'e \`a la conjecture de Valiant.
Nous ne d\'emon\-tre\-rons que les r\'esultats les plus simples et nous 
souhaitons faire sentir l'importance des enjeux.

\smallskip Dans la section \ref{secFaforb} nous faisons une \'etude 
rapide des classes
de \com \bolez, qui constituent une variante non uniforme de
la \cobiz. 

Dans la section \ref{secAgBo} nous mettons en \'evidence 
quelques liens 
\'etroits et simples entre \fonbs et \polsz, et entre
\com \bole et \coagz.

Dans la section \ref{secCobiCobo} nous faisons le lien entre
\cobi et \com \bolez.
Dans la section \ref{secPermUniv} nous donnons quelques r\'esultats
sur le permanent.
Dans la section finale, nous rappelons la conjecture de Valiant
et discutons bri\`evement sa port\'ee.

Parmi les r\'ef\'erences utiles pour ce chapitre,
il faut citer le livre \cite{Weg} et l'article \cite{SkVa}, non encore 
signal\'es. 

\section{Familles d'\exns et de \cibosz}
\label{secFaforb}

\subsubsection*{Expressions, circuits et descriptions}  

L'analogue  \bol de l'anneau de \pols $\A\,[x_1,\ldots,x_n]$\, est 
l'\agr de Boole 
$$\BB\,[x_1,\ldots,x_n]\simeq 
\aqo{\FF_2[x_1,\ldots,x_n]}{x_1^2-x_1,\ldots,x_n^2-x_n}\,$$
avec dans \,$\BB\,[x_1,\ldots,x_n]$\, les \egts 
\,$a\land b=ab$, 
\,$a\lor b=a+b+ab$, \,$\lnot a=1+a$ et 
\,$a+b=(\lnot a\land b)\lor(a\land\lnot b)$. Cette interpr\'etation 
de l'\agr de Boole librement engendr\'ee par $n$ \elts comme quotient 
d'un anneau de \pols \`a  $n$ variables sur le corps  
$\FF_2$ laisse penser que les \mets \agqs sont a priori 
pertinentes pour r\'esoudre les \pbs \bolsz.

L'analogue \bol d'une fonction \pol \`a \,$n$\, variables est une 
\emph{\fonbz}\index{fonction bool\'eenne} 
 \,$f:\left\{0,1\right\}^n\rightarrow \left\{0,1\right\}$.
Nous aurons aussi \`a consid\'erer des 
\emph{\apbsz}\index{application!bool\'eenne} 
\,$g:\left\{0,1\right\}^n\rightarrow \left\{0,1\right\}^m$.

L'\agr de Boole \,$\BB\,[x_1,\ldots,x_n]$\, est isomorphe
\`a l'\agr des \fonbs \,$f:\left\{0,1\right\}^n\rightarrow 
\left\{0,1\right\}$. L'\iso fait correspondre \`a l'\elt \,$x_i$\,
de \,$\BB\,[\s{x}]$\, la \,$i$\,-\`eme fonction coordonn\'ee:
\,$(a_1,\ldots ,a_n)\mapsto a_i$.

Rappelons que si \,$p_1,\ldots,p_n$\, sont les variables \boles 
pr\'esentes dans une \exbz, on appelle 
\emph{litt\'eral}\index{litt\'eral} l'une des expressions
\,$p_1,\lnot p_1,\ldots,\alb p_n,\alb\lnot p_n.$  Et qu'une \exn
est dite \emph{en forme normale conjonctive} (resp. \emph{en forme 
normale disjonctive}) si elle est une conjonction de disjonctions de 
litt\'eraux (resp. une disjonction de conjonctions de litt\'eraux).
\index{forme normale!conjonctive}
\index{forme normale!disjonctive}

Il y a plusieurs types d'\'ecritures canoniques pour une  \fonbz,
en forme normale conjonctive, en forme normale disjonctive ou
sous forme d'un \pol en repr\'esentation creuse (chaque variable 
intervenant avec un degr\'e $\le 1$ dans chaque mon\^ome).
On peut aussi exprimer une \fonb au moyen d'une  \exb ou d'un \ciboz.

\begin{convention} 
\label{convForbs} 
Nous adopterons la convention qu'une \exb ou un 
\cibo\index{circuit!bool\'een} n'utilisent que les connecteurs 
\,$\land$, \,$\lor$\, et \,$\lnot$. 
En outre dans le cas d'une \exn 
\index{expression!bool\'eenne}l'usage du connecteur 
\,$\lnot$\, sera seulement 
implicite: on utilisera les 
litt\'eraux comme variables (aux feuilles de l'arbre), et nulle part 
ailleurs n'appara\^{\i}tra le connecteur  \,$\lnot$. La taille 
\index{taille!d'une expression bool\'eenne} et
la \prof d'une \exb \index{profondeur!d'une expression bool\'eenne} ne 
prendront en compte que les connecteurs
\,$\land$\,  et  \,$\lor$\, 
(les litt\'eraux sont tous consid\'er\'es comme de \prof nulle).
\end{convention}

Cette convention n'a pas de cons\'equence importante en ce qui concerne 
les \cirs car autoriser d'autres connecteurs ne ferait  diminuer la 
taille et/ou la \prof que d'un facteur constant. 
Par contre, en ce qui 
concerne les \exbsz, il s'agit d'une restriction significative 
de leur pouvoir d'expression: par exemple si on admet en plus le 
connecteur \,$a\oplus b\eqdef (a \land \lnot b)\lor(\lnot a \land b)$\, 
l'\exn \,$p_1\oplus p_2\oplus \cdots \oplus p_n$\, r\'eclamera tr\`es 
probablement une \'ecriture  nettement plus longue sans l'utilisation de 
\,$\oplus$.  
\subsubsection*{Classes de \com \bole}  
 Nous sommes particuli\`erement interess\'es ici par les analogues 
\bols des classes $\VNC$, $\VP$, $\VPF$, $\VNP$ et $\VNPF$ (cf. section 
\ref{secForCirDes}).

\begin{definition} 
\label{defForBoo} 
Soit \,$f_n:\left\{0,1\right\}^{v_n}\rightarrow 
\left\{0,1\right\}$\, une \fam de \fonbs (index\'ee par \,$n\in \N$\, ou 
\,$\N^\ell$).
\begin{itemize}
\item  Nous disons que la \fam \,$(f_n)$\, est 
\emph{\pborz}\index{p-bornee@$p$-born\'ee} 
\index{p-bornee@$p$-born\'ee!famille de fonctions \bolesz}si \,$v_n$\, 
est major\'e par 
un \pol en \,$n$. On dit encore qu'il s'agit d'une \emph{\pfa de 
\fonbsz.}\index{p-famille@$p$-famille!de fonctions bool\'eennes} 
\item  Nous disons qu'une \fam  d'\exbs  \,$(\varphi_n)$\, est 
\emph{\pborz} si la taille de \,$\varphi_n$\, est 
major\'ee par un \pol en \,$n$.
\index{p-bornee@$p$-born\'ee!famille d'expressions bool\'eennes}
\item  Nous disons que la \fam \,$(f_n)$\, est 
\emph{\pforz}
\index{p-exprimable@$p$-exprimable!famille de fonctions \bolesz} 
\index{p-exprimable@$p$-exprimable}\indexnota{BPE@$\BPF$}
si elle est r\'eali\-sa\-ble par une \fam \pbor d'\exbsz. 
La classe des \fams de \fonbs  \pfors 
est not\'ee \,$\BPF$. 
\item  Nous disons qu'une \fam de \cibos  \,$(\gamma_n)$\, est 
\emph{\pborz} si la taille de \,$\gamma_n$\, est 
major\'ee par un \pol en \,$n$.
\index{p-bornee@$p$-born\'ee!famille de circuits bool\'eens}
\item  Nous disons que la \fam \,$(f_n)$\, est 
\emph{$\,p$-\'evaluable}  (ou encore 
\emph{\pcalz}\index{p-calculable@$p$-calculable}) 
\index{p-calculable@$p$-calculable!famille de fonctions \bolesz}
si elle est r\'ealisable par une \fam \pbor de \cibosz. La classe des 
\fams de \fonbs  \pcals est not\'ee 
$\BP=\Ppoly$. \indexnota{P/poly@$\Ppoly$} \indexnota{BP@$\BP=\Ppoly$}
\item Nous notons  
$\BNC^k=\NCkpoly$\indexnota{BNCk@$\BNC^k=\NCkpoly$}\indexnota{BNC@$\BNC=
\NCpoly$} la classe des \fams de 
\fonbs r\'eali\-sa\-bles par une \fam de \cibos 
de taille \polle et de \prof 
en \,$\O(\log^k n)$, et $\BNC=\NCpoly$ d\'enote la r\'eunion des  
$\BNC^k$.
\item  Nous disons qu'une \fonb \,$g$\, en les variables 
\,$p_1,\ldots,\alb p_k,\,\alb r_{1},\alb\ldots,r_{\ell}$\, \emph{est une 
description de la \fonb 
\,$f\,$}\index{description!d'une fonction bool\'eenne}  en les variables 
\,$p_1,\alb\ldots,\alb p_k\,$
si 
\begin{equation} \label{eqDescr2}
f(\s{p})\;=\;\bigvee\nolimits_{\s{r}\in 
\left\{0,1\right\}^\ell}\;g(\s{p},\s{r})
\end{equation}
\item  Nous disons que la \fam \,$(f_n)$\, est 
\emph{\pdbz}\index{p-descriptible@$p$-descriptible} 
 s'il existe une \fam \pcal de \fonbs \,$(g_n)$, telle que 
chaque   \,$g_n$\, est une description de \,$f_n$. 
La classe des \fams de \fonbs  \pdbs est not\'ee 
\,$\BNP=\NPpoly$. \indexnota{NP/poly@$\NPpoly$} 
\indexnota{BNP/poly@$\BNP=\NPpoly$}
\item  Nous disons que la \fam \,$(f_n)$\, est 
\emph{\pdb 
en \exnsz}\index{p-descriptible@$p$-descriptible!en 
expressions}  s'il existe une \fam \pfor de \fonbs \,$(g_n)$, 
telle que chaque   \,$g_n$\, est une description de \,$f_n$. 
 La classe des \fams de \fonbs \pdbs en \exns est not\'ee
 \,$\BNPF$. \indexnota{BNPF@$\BNPF$}
\end{itemize}
\end{definition}

Il faut souligner que toutes les notions introduites ici sont \emph{non 
uniformes}, comme dans le cas \agqz.

La classe  $\Ppoly$ est clairement l'analogue \bol de la classe
$\VP$ en \coagz. C'est aussi un analogue non uniforme de la classe 
$\P$. Ce dernier point sera plus clair apr\`es le \tho 
\vref{thFamCib1}. De m\^{e}me nous verrons que la classe $\NPpoly$ est
un analogue non uniforme de la classe $\NP$.

Si on compare les d\'efinitions des \emph{descriptions}
dans le cas \agq et dans le cas \bolz, on voit qu'on utilise 
maintenant une disjonction \`a la place d'une somme (\forms 
\ref{eqDescr1} et \ref{eqDescr2}).

La notation $\Ppoly$ (voir par exemple  \cite{Bal} ou \cite{Weg}) 
s'explique comme suit: une
\fam \,$(f_n)$\, dans $\Ppoly$  peut \^etre calcul\'ee 
\etp si on a droit 
\`a \gui{une aide} (sous forme d'une \fam de \cibos 
\,$\gamma_n$\, qui calculent les fonctions \,$f_n$) qui n'est 
peut-\^etre pas uniforme mais qui est de taille \polle en  \,$n.$

Signalons que Karp et Lipton, qui introduisent la classe
$\Ppoly$ dans \cite{KaL} donnent une d\'efinition \gnle
pour une variante non uniforme $\cC/\mathit{poly}$  en \com
\bole d'une classe de \cobi arbitraire $\cC$. 
Leur d\'efinition justifie aussi les \'egalit\'es 
$\BNP=\NPpoly$ et
$\BNC^k=\NCkpoly$. 
Enfin la d\'efinition de Karp et Lipton
ne semble rien donner pour  $\P_{\fo}/\mathit{poly}$ par absence de
la classe  $\P_{\fo}$ en \cobiz.

\subsubsection*{La \com \bole des \oparis dans \,$\ZZ$\,} 

Le livre \cite{Weg} de Wegener contient une \'etude pr\'ecise et tr\`es
compl\`ete de la \com des \fams  de \fonbsz.
On y trouve notamment les r\'esultats donn\'es dans le \tho qui suit 
concernant la
 \com \bole des \oparis dans \,$\NN$.
 En fait les r\'esultats sont
uniformes et ils s'\'etendent \immt \`a \,$\ZZ$. 
\begin{theorem} 
\label{thWeg} (\thos 1.3, 2.4 et 2.8 du chapitre 3 de \cite{Weg})\\
L'addition et la multiplication dans \,$\NN$\, sont r\'ealisables
par des \fams de \cibos dans $\BNCU$.
Plus pr\'ecis\'ement: 
\begin{enumerate}
\item L'addition de deux entiers de taille \,$n\geq 3$\, est 
r\'ealisable par un \cibo de taille \,$\leq 9\,n$\, et de \prof 
\,$2\,\esup{\log{n}}+8$.
\item Le produit de deux entiers de taille \,$n$\, est r\'ealisable 
\`a la Karatsuba par un \cibo de taille \,$\O(n^{\log{3}})$\, et de 
\prof \,$\O(\log{n})$, ou m\^{e}me,
en suivant Sch\"onage et Strassen qui adaptent la transformation de 
Fourier discr\`ete rapide des \pols au cas des entiers,
 par un \cibo de taille \,$\O(n\,\log{n}\,\log{\log{n}})$\, et de \prof 
\,$\O(\log{n})$.
\end{enumerate}
\end{theorem}

Concernant la multiplication des entiers on lira aussi avec 
int\'er\^{e}t l'expos\'e de Knuth dans \cite{Knu}.

\subsubsection*{Parall\'elisation des \exbs} 

 Nous avons pour les \exbs un r\'esultat analogue \`a la \paran \`a la 
Brent des \exas  (voir \cite{Spi} ou \cite{Sav} \tho 2.3.3).
\begin{theorem} 
\label{BBrent}  
Pour toute \fonb \,$f:\left\{0,1\right\}^n\rightarrow 
\left\{0,1\right\}$\, la \prof \,$\pi$\, du meilleur \cibo et la taille 
\,$\tau$\, de la meilleure \exb sont reli\'es par
$$ 
\log (\tau +1)\leq \pi\leq\, {2\over \log\,3/2} \,\log(\tau +1)
$$
\end{theorem}
\prv Cela marche de la m\^eme mani\`ere que 
la \paran \`a la Brent des \exasz. 
Le \pol \,$b_0y+b_1$\, dans le cas \agq 
(cf. la preuve du \tho \ref{thBrent}, page
\pageref{PreuveBrent}) doit \^etre remplac\'e par une \exn 
\,$(y\land b_0)\lor(\lnot y\land b_1)$\, dans le cas \bolz.
\qed
 
\medskip 
L'analogue \bol du corollaire \ref{corBrent} ($\VPF=\VNC^1$) est:
\begin{corollary} 
\label{corBBrent} On a $\BPF=\BNCU$. 
\end{corollary}

\subsubsection*{Description des \cibos par des \exbs} 

Le lemme suivant est facile et utile.
\begin{lemma} 
\label{lemSCF} 
\'Etant donn\'e un \cibo \,$\gamma$\, de taille \,$\tau$\, avec les 
portes d'en\-tr\'ee \,$p_1,\ldots,p_n$,  
les portes internes \,$r_1,\ldots,r_\ell$  et une seule porte de sortie
(donc \,$\tau=\ell+1$), on peut construire une \exb en forme normale 
conjonctive \,$\varphi(\s{p},\s{r})\,$
de taille $\leq 7\,\tau$\, et de \prof  $\leq 2+\log(3\,\tau)$\, telle 
que, pour tous \,$\s{p}\in \left\{0,1\right\}^n$\,  on ait 
l'\'equivalence: 
\begin{equation} \label{eqlemSCF}
 \gamma(\s{p}) =1\qquad  
\Longleftrightarrow  \qquad 
\bigvee\nolimits_{\s{r}\in \left\{0,1\right\}^\ell}\; 
\varphi(\s{p},\s{r}) = 1\,.
\end{equation}
En outre dans le second membre il y a une seule affectation des \bols 
\,$r_1,\ldots,r_\ell$\, qui rend l'\exn vraie 
(lorsque $ \gamma(\s{p})=1  $). 
\end{lemma}
\prv  On remplace chaque affectation du \prev d\'efini par le \cibo 
par une \exb qui est vraie \ssi la valeur du \bol affect\'e est 
correcte. La conjonction de toutes ces \exbs  donne l'\exn 
\,$\varphi$. 
Une affectation \,$c:=\lnot a$\, est traduite par 
\,$(c\lor a)\land(\lnot c\lor \lnot a)$. 
Une affectation \,$c:= a\lor b$\, est traduite par 
\,$(\lnot c\lor a\lor b)\land(c\lor \lnot a)\land(c\lor \lnot b)$. 
Une affectation \,$c:= a\land b$\, est traduite par 
\,$(c\lor\lnot a\lor\lnot b)\land(\lnot c\lor a)\land
(\lnot c\lor b)$. \qed

\ms 
On en d\'eduit \immtz.
\begin{proposition} 
\label{propSCF} On a l'inclusion \,$\BP\subset \BNPF$\, et 
l'\egt  \,$\BNP= \BNPF$.
\end{proposition}

Signalons aussi le r\'esultat important de Valiant (pour une
preuve voir \cite{Bur2}).

\begin{theorem} 
\label{thVNP=VNPF} 
On a pour tout corps \,$\VNP=\VNPF$.
\end{theorem}

\subsubsection*{Expressions, circuits et descriptions: le cas des \apbs}  

 Nous pouvons reprendre avec les \fams d'\apbs les d\'efinitions 
donn\'ees au d\'ebut de cette section pour les \fams de \fonbsz. 
Notre objectif est surtout ici de d\'efinir l'analogue non uniforme de 
la classe $\DieseP$.

\begin{definition} 
\label{defAPB} 
Soit \,$f_n:\{0,1\}^{v_n}\, \longrightarrow \, \{0,1\}^{w_n}$\, une \fam 
d'\apbs (index\'ee par \,$n\in \N$\, ou 
\,$\N^\ell$). 
Soit  \,$(f_{n,k})$\, la \fam double de \fonbs 
$ f_{n,k}:\{0,1\}^{v_n}\; \rightarrow \; \{0,1\}$
qui donne la \,$k\,$- \`eme coordonn\'ee de  \,$f_n$\, si
\,$k\leq w_n$.
\begin{itemize}
\item  Nous disons que la \fam \,$(f_n)$\, est 
\emph{\pborz}\index{p-bornee@$p$-born\'ee}
\index{p-bornee@$p$-born\'ee!famille d'applications \bolesz}si \,$v_n$\, 
et \,$w_n$\, sont 
major\'es par un \pol en \,$n$. On dit encore qu'il s'agit d'une 
\emph{\pfa d'\apbsz.}\index{p-famille@$p$-famille!d'applications 
bool\'eennes}
\item  Nous disons que la \fam \,$(f_n)$\, est 
\emph{\pforz} si elle est \pbor et si la \fam double 
\,$(f_{n,k})$\, correspondante est \pforz.
D\'efinition analogue pour une \fam \emph{\pcalz}, 
\emph{\pdbz}, ou
\emph{\pdb en \exnsz}.
\item  Nous disons qu'une \fam \,$g_n:\{0,1\}^{v_n}\, \longrightarrow 
\,\NN$\,
est dans la classe 
\,$\DieseBP=\DiesePpoly$, ou encore qu'\emph{elle compte les solutions 
d'une
\fam \pcal de \fonbsz} si elle v\'erifie:
\begin{equation} \label{eqDieseBP}
\forall \s{p}\quad  g_n(\s{p})\,
\,=\,\#\left\{\; \s{q}\in\left\{0,1\right\}^{\ell_n}\; |\; 
h_n(\s{p},\s{q})=1\; \right\}
\end{equation}
o\`u \,$h_n:\{0,1\}^{v_n+\ell_n}\, \rightarrow \, \{0,1\}$\, est une 
\fam \pcal de  
\fonbsz.\indexnota{DieseBP@$\DieseBP=\DiesePpoly$} 
\noindent Si la \fam \,$(h_n)$\, est
\pfor on dira que la \fam \,$(g_n)$\, est dans la classe 
\,$\DieseBPF$.
\indexnota{DieseBPF@$\DieseBPF$}
\\
Si \,$g_n(\s{p})=\sum_{k=1}^{w_n} f_{n,k}(\s{p})\,2^{k-1}$\,
on dira que la \fam \,$(f_n)$\, est dans $\DieseBP$  
(resp. dans $\DieseBPF$) lorsque la \fam  \,$(g_n)$\, est dans 
$\DieseBP$  
(resp. dans $\DieseBPF$).
\end{itemize}
\end{definition}

Une cons\'equence \imme de la description des \cibos par les 
\exbs (lemme \ref{lemSCF}) est la proposition suivante, 
analogue \`a la 
proposition \ref{propSCF}.
\begin{proposition} 
\label{propDieseCF} ~
 On a  l'\egt  \,$\DieseBP= \DieseBPF$.
\end{proposition}

\begin{remark} 
\label{remDiesePB}
\emph{Il n'y a pas de diff\'erence de principe entre une \fam d'\apbs et 
une \fam de \fonbsz, puisque donner une  \fam d'\apbs revient \`a donner
une \fam double de \fonbsz.
Si on veut d\'efinir directement la classe $ \DieseBP=\DiesePpoly$ comme 
une classe de \fonbsz, on pourra dire que le \pb dans $ \DiesePpoly$ 
associ\'e
\`a la \fam \,$(h_n)\in\Ppoly$\, est le \pb suivant 
portant sur le couple  \,$(\s{p},m)$\, 
(o\`u \,$m$\, est cod\'e en binaire): 
$$ \#\left\{\; \s{q}\in\left\{0,1\right\}^{\ell_n}\; |\; 
h_n(\s{p},\s{q})=1\; \right\}\leq m~?$$ 
} 
\end{remark}

\subsubsection*{La plupart des \fonbs sont difficiles \`a \'evaluer}
On a aussi l'analogue suivant du \tho \vref{thPlupDif}: ici on trouve 
qu'une famille de \cibos de taille \qple ne peut calculer qune infime 
partie de toutes les \fonbsz. 
\begin{prop} 
\label{BPlupDif} Soit \,$\mathit{VQPB}(k)$\, l'ensemble des \fams 
\,$(f_n)$\, de 
\fonbs \`a \,$n$\, variables r\'ealisables par une \fam de \cibos de 
taille \,$2^{\,\log^k n}$. Soit \,$\vep>0$. Pour \,$n$\, assez grand 
seulement
une  proportion $<\vep$\, de \fonbs \`a \,$n$\, variables est dans 
\,$\mathit{VQPB}(k)$.
\end{prop}
\prv  Faisons les comptes. Le nombre total de \fonbs \`a \,$n$\, 
variables est \'egal \`a \,$2^{2^n}$. Le nombre total de \cibos
 \`a \,$n$\, variables  et de taille \,$t+1$\, est major\'e
par \,$N(t+1)=2\,N(t)\,(t+2n)^2$: en effet un \prev de taille
\,$t+1$\, est obtenu en rajoutant une instruction \`a un \prev de taille
\,$t$, instruction de la forme \,$x_{t+1}\aff y\circ z$\, avec 
\,$\land$\,  ou \,$\lor$\, pour \,$\circ\,$, et \,$y,\,z$\, sont \`a 
choisir parmi les litt\'eraux ou parmi les \,$x_i$\, ($1\leq i\leq t$). 
Cette majoration
conduit \`a \,$N(1)=2\,(2n)^2$, \,$N(2)= 2^2\,(2n)^2\,(2n+1)^2$, \ldots, 
\,$N(t)< 2^t\,((2n+t)!)^2= 2^{\O((n+t)\,\log (n+t))}$.
Donc si \,$t=2^{\,\log^k n}$, \,$\log N(t)=\O(2^{\log^{k+1} n})$\,
qui devient n\'egligeable devant \,$2^n$\, pour \,$n$\, grand.
\qed 

\medskip On trouvera des r\'esultats du m\^{e}me style mais nettement 
plus pr\'ecis dans le chapitre 4 du livre de Wegener \cite{Weg}.

\section[Bool\'een versus alg\'ebrique]{Bool\'een versus alg\'ebrique 
(non uniforme)} 
\label{secAgBo}

\subsection{\'Evaluation \bole des \caris}     

Rappelons ici le \pbz, d\'ej\`a \'evoqu\'e \`a la section 
\ref{subsecCobicaris}, de l'\eva d'un \cari sur un anneau \,$\A$\,
dont les \elts sont cod\'es en binaire.
Si l'anneau \,$\A$\, est fini, le temps \paral ou \sql du calcul \bol 
correspondant \`a l'\'ex\'ecution d'un \cari est simplement
proportionnel \`a la \prof ou \`a la taille du \cariz. Par ailleurs
rappelons que $\VP=\VNC=\VNC^2$,  $\VNC^1=\VPF$ et
 $\BNC^1=\BPF$.  On obtient donc:
\begin{lemma} 
\label{lemAnFini} 
Si une \pfa de \pols sur un anneau fini \,$\A$\, est dans la classe 
$\VP$ (resp. $\VPF$,  $\VNP$) son \eva \bole est donn\'ee
par une \fam  dans la classe \bole 
$\BNC^2\subset\BP$ (resp.  $\BPF$,  $\DieseBP$).
\end{lemma}

Dans le cas d'un anneau infini, l'\eva \bole d'un \cari
peut r\'eserver quelques mauvaises surprises (voir l'exemple de 
l'inventeur du jeu d'\'echec page \pageref{2p2n}). 
Il faudrait bannir toute constante (m\^{e}me 1~!)
d'un \cari sur \,$\ZZ$\, 
si on veut que l'\eva \bole (avec le codage naturel
binaire de \,$\ZZ$) ne produise pas d'explosion 
(et \,$\ZZ$\, est l'anneau infini le plus simple).

Une solution serait de coder les \elts de l'anneau par des 
\caris n'ayant que des constantes en entr\'ees{\footnote{~Du point de 
vue des calculs \etp on peut remarquer que le codage binaire usuel de
 \,$\ZZ$\, est \'equivalent \`a un codage par des \exas n'ayant que les
 constantes \,$0$, \,$1$, ou \,$-1$ aux feuilles de l'arbre. Il n'est 
donc pas artificiel de proposer un codage de
 \,$\ZZ$\,  par des \caris n'ayant que les
 constantes \,$0$, \,$1$, ou \,$-1$ aux portes d'entr\'ee, ce que nous 
avions not\'e  \,$\ZZ_{\mathrm{preval}}$.}}. 
Mais le test d'\egtz, le test de signe et bien d'autres op\'erations 
simples sur \,$\ZZ$\, semblent alors 
sortir de la classe \,$\P$.

Une autre solution serait d'apporter une restriction
plus s\'ev\`ere aux \fams \pbors de \carisz.
Avant d'y introduire la moindre constante, m\^{e}me 1,
la \fam devrait \^{e}tre \pbor (en taille et en degr\'es). 
Ensuite seulement on remplacerait certaines variables 
par des constantes. 

De mani\`ere \gnle il faut avoir une majoration convenable de
la taille des objets \`a calculer.

\begin{definition} 
\label{def-ffbpbor} 
Une \fam de fonctions \,$f_n:\ZZ^{v_n}\rightarrow \ZZ$\, est
dite \emph{\pbor en taille} si \,$v_n$\, est major\'e par un \pol 
en \,$n$\, et la taille de \,$f_n(x_1,\ldots ,x_{v_n})$\, est major\'ee 
par un \pol en la taille
de l'entr\'ee \,$x_1,\ldots ,x_{v_n}$\, (en utilisant les codages 
binaires \uslsz).
\end{definition}

On a alors l'extension importante
suivante du lemme \ref{lemAnFini}
\`a l'anneau \,$\ZZ$, sous une condition restrictive \suptz, qui est 
d'ailleurs in\'evitable.

\begin{lemma} 
\label{lemEvbZQ} 
On consid\`ere une \pfa \,$(P_n)$\, de \pols sur  \,$\ZZ$. On suppose 
que la \fam de fonctions \,$f_n:\ZZ^{v_n}\rightarrow \ZZ$\, d\'efinie 
par  \,$(P_n)$\,
est \pbor en taille.
Alors si \,$(P_n)$\,  est dans  $\VP(\ZZ)=\VNC^2(\ZZ)$ (resp.   
$\VPF(\ZZ)=\VNC^1(\ZZ)$, $\VNP(\ZZ)$) son \eva \bole est donn\'ee
par une \fam  de \cibos dans 
$\BNC^3\subset\BP$ (resp. $\BNC^2$,  $\DieseBP$).
\end{lemma}
\prv
Supposons que \,$(P_n)$\,  est dans la classe $\VP(\ZZ)$ et soit
\,$(\Gamma_n)$\, une \fam \pbor de \caris correspondant \`a
\,$(P_n)$.
 Pour tous \,$m,\,n$\, entiers positifs on veut 
construire un \cibo \,$\gamma_{n,m}$\, qui calcule (le code de)
\,$f_n(x_1,\ldots ,x_{v_n})$\, \`a partir des (codes des) 
\,$x_i$\, lorsqu'ils sont de taille \,$\leq m$. 
On sait que la taille de la sortie \,$y$\, est major\'ee
par un entier \,$p \leq C\,(n+m)^k$.
Il suffit alors de prendre les constantes de \,$\Gamma_n$\, modulo  
\,$2^{2p}$\, et d'ex\'ecuter les calculs indiqu\'es par le circuit
\,$\Gamma_n$\, modulo \,$2^{2p}$\,
pour r\'ecup\'erer \,$y$\, comme \elt de \,$\ZZ$\, \`a la fin du calcul. 
La taille du \cibo  \,$\gamma_{n,m}$\,  correspondant est bien
\polt major\'ee. Quant \`a sa \profz, compar\'ee \`a celle de
\,$\Gamma_n$,  elle a \'et\'e multipli\'ee par un 
\,$\O(\log{p})=\O(\log(m+n))$
(cf. le \thoz~\ref{thWeg})\\
Le r\'esultat pour  \,$(P_n)$\,  dans la classe $\VNP(\ZZ)$ se d\'eduit 
\immt du r\'esultat pour  \,$(P_n)$\,  dans la classe $\VP(\ZZ)$.
\qed

\medskip Notez que  si \,$(P_n)$\, est dans $\VNC^1$, l'hypoth\`ese 
que \,$(f_n)$\, est \pbor en taille 
est automatiquement v\'erifi\'ee si les constantes du \cir 
\,$\Gamma_n$\, ont une taille major\'ee
par un  \,$C\,n^\ell$.
Dans la section suivante, tous les \caris qui simulent des \cibos
utilisent les seules constantes $0$, $1$ et $-1$.
\subsection{Simulation \agq des \cirs et \exbs}  
\label{secSimBP}
Nous nous int\'eressons dans cette section \`a la possibilit\'e de 
\emph{simuler \agqtz} une \fonb 
\,$f:\{0,1\}^n\rightarrow \{0,1\}$, ou une application 
\,$g:\{0,1\}^n\rightarrow\ZZ$\, 
(par exemple cod\'ee par un \cibo \,$\{0,1\}^n\rightarrow\{0,1\}^m$).

Nous disons que \emph{le \pol \,$P$\, simule la \fonb  
\,$f$\, en \eva}\index{simulation!en \'evaluation} s'il a le 
m\^eme nombre de variables, et s'il  s'\'evalue de la m\^eme mani\`ere 
que la \fonb sur des entr\'ees dans \,$\{0,1\}$. 

D\'efinition analogue
pour la simulation \agq de l'application \,$g$\, par un \pol (l'anneau 
doit contenir \,$\ZZ$).
\subsubsection*{Un r\'esultat \elr} 
Le lemme suivant nous dit ce que donne la simulation naturelle d'un 
\cibo  par un \cariz: la \prof et la taille sont convenables mais les 
degr\'es peuvent r\'eserver de mauvaises surprises.
\begin{lemma} 
\label{lemAlgbri} Un \cibo 
 \,$\gamma$\,  de taille \,$\tau$\, et de 
\prof \,$\pi$\, peut \^etre simul\'e en \eva par un \cari   
\,$\psi$\, de taille \,$\leq 4\tau$\, et de \prof \,$\leq 3\pi$\, (sa 
\prof \muv reste \'egale \`a \,$\pi$\, donc le degr\'e des \pols est 
\,$\leq 2^\pi$). Cette simulation fonctionne sur tout
\acom (non trivial).
\end{lemma}
\prv Les seules valeurs des \bols sont $0$ et $1$, on a donc
$$x\land y=xy,\quad \lnot x = 1-x,\quad x\lor y =x+y-xy$$
sur n'importe quel \acom (non trivial){\footnote{~Nous 
rappelons que dans 
les chapitres \ref{chap DeterUniv}  et  \ref{chap Perma}
les seules \oparis autoris\'ees sont \,$+$\, 
et \,$\times$\, ce qui nous contraint \`a introduire des \muls par la 
constante $-1$ pour faire des soustractions. Ceci implique
que le \pol \,$x+y-xy$\, est \'evalu\'e par un \cir 
de \prof \'egale \`a $3$.}}. 
\qed

\subsubsection*{Simulation d'une \exb par une \exa} 

Le lemme suivant est une cons\'equence directe du \tho de \paran
\ref{BBrent} et du lemme \ref{lemAlgbri}.

\begin{lemma} 
\label{lemSFBFA} 
Une \exb \,$\varphi$\, de taille \,$\tau$\,  peut \^etre 
si\-mu\-l\'ee
en \eva par une \exa  de \prof major\'ee par \,$ {6\over \log \,3/2} 
\,\log (\tau +1)\leq \alb 10,26 \,\log (\tau +1)$. 
Cette simulation fonctionne sur tout \acom (non trivial).
\end{lemma}

En particulier la taille de l'\exa est \,$\leq (\tau +1)^{10,26}$
({\footnote{~Le degr\'e du \pol est major\'e par \,$(\tau 
+1)^{3,419}$.}}). 
On en d\'eduit:
\begin{proposition} 
\label{propSimul1} ~Toute \fam dans  $\BPF$ est simul\'ee
\agqt par une \fam dans $\VPF$.
Cette simulation fonctionne sur tout \acom (non trivial).
\end{proposition}

Dans \cite{Bur2} la proposition pr\'ec\'edente est \'enonc\'ee avec une 
terminologie diff\'erente: \gui{$\BPF$ est contenu dans la partie \bole 
de $\VPF$}.

\smallskip Une proposition analogue \`a la pr\'ec\'edente et qui 
voudrait relier 
de mani\`ere aussi simple les classes
$\BP$ et $\VP$ \'echouerait parce que la traduction 
naturelle d'un \cibo  
en un \cari donn\'ee au lemme \ref{lemAlgbri}
fournit en \gnl un \pol de degr\'e trop grand. 
Autrement dit, on ne conna\^{\i}t pas d'analogue satisfaisant du lemme
\ref{lemSFBFA} pour les \cibosz.

\ms Supposons maintenant que nous ayons d\'emarr\'e avec une 
\pfa double  d'\exbs 
\,$(\varphi_{n,k})$\, associ\'ee \`a une \fam de fonctions  
\,$f_n:\left\{0,1\right\}^{v_n}\rightarrow \ZZ$. La sortie est cod\'ee 
par exemple comme suit \refstepcounter{bidon}\label{codeZ} 
dans $\left\{0,1\right\}^{m}$, le 
premier bit code le signe,
et les bits sui\-vants codent l'entier sans signe en binaire (suppos\'e 
$<2^{m-1}$). 
Par exemple avec \,$m=8$\, les entiers $5$, $-11$ et $69$ sont 
respectivement cod\'es par \,$00000101$, \,$10001011$\, et 
\,$01000101.$ Il n'y a alors aucune difficult\'e \`a calculer par un
\cari ou par une \exa de \prof \,$\O(\log m)$\, la sortie dans \,$\ZZ$\, 
\`a partir de son code.

Nous pouvons alors \'enoncer la proposition suivante, qui \gni 
la proposition \ref{propSimul1}, et qui r\'esulte \egmt du lemme 
\ref{lemSFBFA}.

\begin{proposition} 
\label{propSimul2} ~
Soit \,$(\varphi_{n,k})_{1\leq k\leq a+n^h}$\, une 
\pfa double dans $\BPF$ 
qui code une \fam de fonctions  
\,$f_n:\left\{0,1\right\}^{v_n}\rightarrow \ZZ$.
Alors il existe une \pfa \,$(\gamma_n)$\, d'\exas dans  
$\VPF$ qui simule en \eva la \fam \,$(f_n)$\, sur n'importe quel anneau 
contenant  \,$\ZZ$. 
\end{proposition}

\subsubsection*{Description d'un \cibo par une 
\exa} 

Nous pouvons faire une synth\`ese des lemmes \ref{lemSCF} et 
\ref{lemAlgbri} pour obtenir une description \agq (au sens de la 
d\'efinition \ref{defForEvDes}) d'un \ciboz. 

\begin{lemma} 
\label{lemCompt} 
Soit \,$\gamma$\, un \cibo de taille \,$\tau=\ell+1$\,
qui calcule une \fonb \,$f:\{0,1\}^{n}\, \rightarrow \, \{0,1\}$.  Il 
existe une \exa \,$\psi(x_1,\ldots,x_n,y_1,\ldots ,y_\ell)$\, de taille 
$\leq 14\,\tau$\, et de \prof  
$\leq 4+\esup{\log(3\,\tau)}$  v\'erifiant:
\begin{equation} \label{eqlemCompt}
 \forall \s{p}\in\{0,1\}^{n}\qquad  
f(\s{p}) =\sum\nolimits_{\s{r}\in \{0,1\}^\ell}\; 
\psi(\s{p},\s{r})
\end{equation}
Cette \exa utilise les seules constantes $0$, $1$ et $-1$ et l'\egt est 
valable sur tout \acom (non trivial).\end{lemma}
\prv 
On applique la simulation donn\'ee dans le lemme \ref{lemAlgbri} \`a 
l'\exb \,$\varphi$\, en forme normale conjonctive construite au lemme 
\ref{lemSCF}. 
On doit simuler \agqt chacune des \exbs de base qui sont du type 
\,$(\lnot c\lor a\lor b)\land(c\lor \lnot a)\land(c\lor \lnot b)$ 
ou du type
\,$(c\lor\lnot a\lor\lnot b)\land(\lnot c\lor a)\land
(\lnot c\lor b)$.
Dans ces \exbs  \,$c$\, est un litt\'eral positif et \,$a,b$\,
des litt\'eraux positifs ou n\'egatifs.
L'examen pr\'ecis montre que la taille maximum pour une telle simulation 
est \,$11$. Il reste ensuite \`a faire le produit de
\,$3\tau$\, expressions (chacune correspond \`a l'un des composants dans
les deux types ci-dessus). On obtient alors
 une \exa \,$\psi$\, de taille $\leq 14\,\tau$\, et de \prof  
$\leq 4+\esup{\log(3\,\tau)}$ v\'erifiant:
\begin{equation} \label{eqlemAlgbri2}
 \forall \s{p}\in\{0,1\}^{n}\quad \left( f(\s{p}) =1\;  
\Longleftrightarrow  \; 
\bigvee\nolimits_{\s{r}\in \{0,1\}^\ell}\; 
\psi(\s{p},\s{r}) = 1\,\right)
\end{equation}
En outre dans le second membre il y a une seule affectation des 
variables  
\,$r_1,\ldots,r_\ell$\, dans \,$\{0,1\}^\ell$\, 
qui rend l'\exn \,$\varphi$\, vraie (lorsque $f(\s{p})=1  $), \cad
que \,$\psi(\s{p},\s{r})$\, est nulle pour tout \,$\s{r}\in 
\{0,1\}^\ell$\, \`a l'exception de cette valeur.
D'o\`u l'\egt \ref{eqlemCompt}:
$ f(\s{p}) =\sum\nolimits_{\s{r}\in \{0,1\}^\ell}\; 
\psi(\s{p},\s{r}).$
\qed

\ms Nous en d\'eduisons les corollaires suivants.

\begin{proposition} 
\label{propCompt} 
Toute \fam \,$f_n:\{0,1\}^{v_n}\, \rightarrow \, \{0,1\}$\, dans $\BP$
est simul\'ee en \eva par une \fam dans $\VNPF$,
ceci sur tout \acom (non trivial). Toute \fam \,$g_n:\{0,1\}^{v_n}\, 
\rightarrow \, \NN$\, dans  $\DieseBP$,  est simul\'ee en \eva par une 
\fam dans $\VNPF$ sur tout \acom contenant \,$\ZZ$. 
\end{proposition}

\begin{theorem} 
\label{propBPDieseBP} 
Si $\VP(\ZZ)=\VNPF(\ZZ)$ alors $\BNC^3=\BP=\BNP= \DieseBP$.
\end{theorem}
\prv
Supposons $\VP(\ZZ)= \VNPF(\ZZ)$ et soit \,$(f_n)$\, une \fam dans 
$\DieseBP$. Remarquons que \,$(f_n)$\, est \pbor en taille.
Par la proposition pr\'ec\'edente, cette \fam est simul\'ee en
\eva par une \fam dans $\VNPF(\ZZ)$ donc par une \fam dans 
$\VP(\ZZ)=\VNC^2(\ZZ)$.
Or une telle \fam s'\'evalue par une \fam dans
$\BNC^3$ d'apr\`es le lemme \ref{lemEvbZQ}.
\qed

\medskip En fait, en utilisant des techniques nettement plus subtiles, 
B\"urgisser a montr\'e les r\'esultats suivants (cf. \cite{Bur2}).

\begin{theorem} 
\label{thBur}~
\begin{enumerate}
\item Soit \,$\FF_q$\, un corps fini, si $\VP(\FF_q)=\VNPF(\FF_q)$ alors 
$\BNC^2=\BP=\BNP$.
\item Soit \,$\K$\, un corps de \cara nulle. Supposons
que l'hypoth\`ese de Riemann g\'en\'eralis\'ee est vraie.
Si $\VP(\K)=\VNPF(\K)$ alors $\BNC^3=\BP=\BNP$.
\end{enumerate}
\end{theorem}

\subsection{Formes \agqs d\'eploy\'ees}
\label{secfald}

\subsubsection*{Forme \agq d\'eploy\'ee d'une \fonb}

Pour traiter les questions de taille d'\exns ou de \cibos il est a 
priori prometteur d'interpr\'eter une \fonb par un \pol \agq usuel.
Une traduction par\-ti\-culi\`erement simple consiste \`a 
\emph{\'etaler}  certaines valeurs de la \fonbz:
 on remplace la \fonb \`a \,$m+n$\, variables \,$f(p_1,\ldots ,\alb 
p_m,\alb q_1,\ldots ,q_n)$\, par le \pol suivant, 
en \,$m$\, variables, avec pour seuls exposants $0$ ou $1$ dans les 
mon\^omes 
\begin{equation} \label{fald}
F(p_1,\ldots ,p_m)=\sum\nolimits_{\mu\in  \left\{0,1\right\}^n}\, 
f(p_1,\ldots ,p_m,\mu)\, x^\mu
\end{equation}
 o\`u \,$\mu=\mu_1,\ldots,\mu_n$\, et 
\,$x^\mu=x_1^{\mu_1}\cdots x_n^{\mu_n}$.

Nous dirons que
\emph{le \pol  \,$F$\, est la \fad\index{forme alg\'ebrique 
d\'eploy\'ee!d'une fonction bool\'eenne} (sur les
variables \,$q_1,\ldots ,q_n$) de la  \fonb \,$f$.} 
Lorsque \,$m\neq 0$\, chaque \coe de \,$F$\, est une
\fonb de \,$p_1,\ldots ,p_m$\, qui doit \^{e}tre simul\'ee \agqtz.
Lorsque \,$m=0$\, on a une \emph{\fad pure} et les \coes de \,$F$\,
sont tous \'egaux \`a \,$0$\, ou~$\,1$.

Une d\'efinition analogue est \'egalement valable si on remplace 
\,$f:\{0,1\}^{m+n}\rightarrow\{0,1\}$\, par une application
\,$g:\{0,1\}^{m+n}\rightarrow\ZZ$.

Si la  \fonb \,$f$\, est facile \`a calculer, le \pol correspondant  
\,$F$\, aura ses \coes faciles \`a \'evaluer, mais il risque d'\^etre 
difficile \`a \'evaluer, puisqu'il y aura en \gnl un nombre trop 
grand (exponentiel en \,$n$) de \coes non nuls. 

On a alors comme cons\'equence des r\'esultats pr\'ec\'edents.
\begin{lemma} 
\label{lemDeployfb} 
Soit une \fonb \,$f(p_1,\ldots ,\alb p_m,\alb q_1,\ldots ,q_n)$\, 
\'e\-valu\'ee par un \cibo \,$\gamma$\,  
de taille \,$\tau$. Sa \fad  \,$F$\, sur les variables \,$q_1,\ldots 
,q_n$\, admet une description (au 
sens de la d\'efinition \ref{defForEvDes}) par une \exa de \prof 
\,$\leq 5+\esup{\log(3\tau)}$\, et de taille
\,$\leq 14\tau +4n\leq 18\tau.$
Cette \exa utilise les seuls constantes $0$, $1$ et $-1$ et  est valable 
sur tout \acom (non trivial).
\end{lemma}
\prv Cela r\'esulte du lemme \ref{lemCompt} et de la 
constatation suivante. 
On a pour \,$\mu_1,\ldots ,\mu_n\in\{0,1\}$\,
$$
x_1^{\mu_1}\cdots x_n^{\mu_n}= \prod\nolimits_{i=1}^n{(\mu_i(x_i-1)+1)}
$$
qui s'\'ecrit comme une \exn de \prof \,$\leq 3+\esup{\log{n}} 
<4+\esup{\log(3\tau)}$\, 
et de taille \,$4n-1.$  
Donc si la \fonb \,$f$\, est d\'ecrite
par l'\exa \,$\psi(\s{p},\s{q},\s{r})$ (lemme \ref{lemCompt}), 
le \pol \,$F$\, est \'egal \`a
$$  
\sum\nolimits_{(\s{\mu},\s{r})\in\{0,1\}^{n+\ell}} 
{\;\;\psi(\s{p},\s{\mu},\s{r})\;\cdot\;
\prod\nolimits_{i=1}^n{(\mu_ix_i+1-\mu_i)}
} 
$$
 et il admet pour description l'\exa \`a \,$m+2n+\ell$\, variables
$$ 
\theta(\s{p},\s{x},\s{\mu},\s{r})\;=\;\psi(\s{p},\s{\mu},\s{r})\;
\cdot\;\prod\nolimits_{i=1}^n{(\mu_ix_i+1-\mu_i)} 
$$
de \prof \,$\leq 5+\esup{\log(3\tau)}$\, et de taille
\,$\leq 14\tau+4n.$  
\qed

\subsubsection*{Forme \agq d\'eploy\'ee d'une \fam de \fonbs} 
 
Une \fam de \fonbs 
\,$f_n:\left\{0,1\right\}^{v_n+w_n}\rightarrow \left\{0,1\right\}$\,
admet pour \fad\index{forme alg\'ebrique 
d\'eploy\'ee!d'une \fam de \fonbsz} (sur les  \,$w_n$\, derni\`eres 
variables) la \fam des \pols \,$F_n$\, 
qui sont les \fads des fonctions \,$f_n$.
M\^{e}me chose pour la \fad d'une \fam \,$g_n:\{0,1\}^{v_n+w_n}\, 
\rightarrow \, \NN$.

On a  comme corollaire du lemme \ref{lemDeployfb}.
\begin{theorem} 
\label{propDeployP} \emph{(Crit\`ere de Valiant)}
Toute \fam de \fonbs dans  $\BP$ admet pour \fad une \fam de \pols dans  
$\VNPF$, qui convient pour tout \acom (non trivial).
En cons\'equence une \fam dans  $\DieseBP$ admet pour \fad une \fam dans 
$\VNPF(\ZZ)$, et cette \fam convient pour tout anneau contenant \,$\ZZ$.
\end{theorem}

Dans le cas d'une \fonb cela peut sembler 
un peu d\'e\-ce\-vant, puisqu'a priori $\VNPF$ est une 
classe r\'eput\'ee difficile \`a calculer (elle simule 
$\DieseBPF=\DieseBP$), mais il y a une tr\`es bonne 
raison \`a cela. En effet, supposons qu'on d\'eploie toutes les 
variables, alors si on calcule \,$F_n(1,\ldots,1)\,$
on trouve le nombre total des solutions de l'\'equation 
 \,$f_n(\s{p})=1$, \cad la somme \,$\sum_{\s{p}}f_n(\s{p})$. 
Et ce n'est donc pas surprenant que \,$F_n$\, soit a 
priori plus difficile \`a calculer que ses \coesz. 
De mani\`ere \gnlez, on ne peut gu\`ere esp\'erer que 
l'int\'egrale d\'efinie d'une fonction soit en \gnl aussi simple 
\`a calculer que la fonction elle-m\^eme.

Le crit\`ere de Valiant, malgr\'e la simplicit\'e de sa preuve, est un 
moyen puissant pour fabriquer des familles dans $\VNPF$.

Comme toutes les preuves que nous avons donn\'ees dans les chapitres 
\ref{chap DeterUniv} et \ref{chap Perma}, la preuve du crit\`ere de 
Valiant est clairement uniforme. Donc si \,$(g_n)$\, est une \fam dans 
$\DieseP$
(on prend pour entr\'ee le mot form\'e par l'entier b\^{a}ton \,$n$\, 
suivi d'un $0$ puis du mot \,$\s{p}$), alors la \fad de \,$(g_n)$\, 
admet pour description une \fam \emph{uniforme} de \caris dans
\,$\NC^1$\,
qui utilise les seules constantes $0$, $1$ et $-1$ et qui
donne le r\'esultat correct sur tout anneau contenant \,$\ZZ$.

\section{Complexit\'e binaire versus \com \bole} 
\label{secCobiCobo}
\subsubsection*{Famille de \fonbs associ\'ee \`a un \pb \algq}   
Notons \,$\left\{0,1\right\}^{\star}$\, l'ensemble des mots \'ecrits sur 
l'alphabet \,$\left\{0,1\right\}$. Nous pouvons voir cet ensemble comme 
la r\'eunion disjointe des  \,$\left\{0,1\right\}^n$.

Consid\'erons un \emph{\pb \algqz} \,$P$\, qui est cod\'e sous 
forme binaire: autrement dit, toute instance de ce \pb correspond \`a 
une question cod\'ee comme un \elt de 
\,$\left\{0,1\right\}^n$\, pour un certain entier \,$n$\, et la 
r\'eponse \`a la question, du type oui ou non, est elle-m\^eme cod\'ee 
comme un \elt de \,$\left\{0,1\right\}$.

On peut interpr\'eter ce \pb \,$P$\, comme fournissant, pour chaque 
\,$n$, une \fonb   \,$f_n:\left\{0,1\right\}^n\rightarrow 
\left\{0,1\right\}$.
Nous dirons que la famille \,$(f_n)$\, est \emph{la famille de \fonbs 
associ\'ee au \pb \algq \,$P.$}\refstepcounter{bidon}\label{fafonbass} 

Supposons maintenant que le \pb  \,$P$\, porte sur les graphes 
orien\-t\'es. Un code naturel pour un graphe orient\'e \`a \,$n$\, 
sommets est sa \emph{matrice d'ajacence} qui est une matrice dans
\,$\{0,1\}^{n\times n}.$ Cette matrice contient 1 en position 
\,$(i,j)$\, \ssi
il y a une ar\`ete orient\'ee qui va de \,$i$\, \`a  \,$j$\, 
dans le graphe consid\'er\'e. Dans ce cas, on voit que la famille de 
\fonbs associ\'ee au \pb  \,$P$\, est plus naturellement d\'efinie comme 
une \fam \,$f_n:\left\{0,1\right\}^{n^2}\rightarrow \left\{0,1\right\}$. 

On dira que le \pb \algq \,$P$\, est dans une classe de \com \bole
\,$\cC$\, si la \fam de \fonbs qui lui est naturellement attach\'ee est
dans \,$\cC$.

\subsubsection*{Famille d'\apbs associ\'ee \`a une fonction \algq} 

Consid\'erons maintenant une
\emph{fonction \algqz} \,$F,$  une fonction qu'on aurait envie de faire 
calculer par un ordinateur:
l'entr\'ee et la sortie sont cod\'ees
en binaire, \cad consid\'er\'ees comme des \elts de
 \,$\left\{0,1\right\}^{\star}$.

Supposons que  \,$t(n)$\,  est une majoration de la taille de la sortie 
en fonction de la taille  \,$n$\, de l'entr\'ee et que la fonction  
\,$t$\, \gui{n'est pas plus difficile \`a calculer} que  \,$F.$

Nous pouvons alors recalibrer la fonction  \,$F$\,  de mani\`ere que la 
taille de sa sortie ne d\'epende que de la taille de son entr\'ee.
Par exemple nous prenons la fonction
 \,$G$\, qui, pour un mot \,$\mu$\, en entr\'ee de taille \,$n,$ 
calcule le mot \,$F(\mu)$\, pr\'ec\'ed\'e d'un 1, lui-m\^eme 
pr\'ec\'ed\'e du nombre de 
0 \ncr pour atteindre la longueur $1+t(n)$.
Il est clair qu'on r\'ecup\`ere facilement \,$F$\, \`a partir de  \,$G.$

Cette convention nous permet d'associer \`a toute 
fonction \algq \,$F$\, une \fam d'\apbs 
$$\,f_{n}: \left\{0,1\right\}^n\rightarrow \left\{0,1\right\}^{1+m}$$  
o\`u \,$m=t(n)$\,  est une majoration de la taille de la sortie en 
fonction de la taille de l'entr\'ee. 
La famille d'\apbs associ\'ee \`a  \,$F$\, d\'epend donc de la fonction 
de majoration \,$t$\, que l'on consid\`ere.

\begin{definition} 
\label{defFamfbass} 
Dans les conditions ci-dessus, nous dirons que \emph{la famille 
\,$(f_n)$\, est la famille d'\apbs 
associ\'ee \`a la fonction \algq \,$F$\,  avec la fonction de majoration 
\,$t.$} 
Si nous ne pr\'ecisons pas cette fonction de majoration, nous disons 
simplement que la famille \label{fafoqass} \,$(f_n)$\, est \emph{une} 
famille d'\apbs associ\'ee \`a la fonction \algq \,$F.$
On dira que la fonction \algq \,$F$\, est dans une classe de \com \bole
\,$\cC$\, si la \fam \,$(f_n)$\, est
dans~$\,\cC$.
\end{definition}

Lorsque la fonction  \,$F$\, est calculable, dans une classe de \cobi 
connue, on choisira toujours la fonction de majoration suffisamment 
simple, de fa\c{c}on que la fonction  \,$G$\,  reste dans la m\^eme 
classe de \comz. 

Tout ce que nous venons de dire s'applique par exemple \`a une fonction 
\,$F$\,  de \,$\left\{0,1\right\}^{\star}$ ou  $\NN$  vers $\NN$ ou 
$\ZZ$, modulo des codages binaires naturels convenables.

\subsubsection*{Familles uniformes de \cibos}

Consid\'erons un \pb \algq  \,$P$\, qui est cod\'e sous forme 
binaire:  pour chaque \,$n$, une \emph{\fonbz}  
\,$f_n:\left\{0,1\right\}^n\rightarrow \left\{0,1\right\}$\, donne la 
r\'eponse pour les mots de longueur \,$n$.
Cette \fam \,$(f_n)$\, peut \^etre r\'ealis\'ee 
sous forme d'une \fam d'\exbsz, ou sous forme d'une \fam de \cibosz.

En \cobi on s'int\'eresse \`a la fois \`a la taille (de ces \exns ou de 
ces \cirsz) et \`a la difficult\'e proprement \algq qu'il peut y avoir 
\`a produire l'\exn (ou le \cirz) \num$n$ en fonction de l'entr\'ee 
\,$n$\, (cod\'ee en unaire).
 Ce deuxi\`eme aspect correspond \`a la question: 
la \fam est-elle uniforme~?

Un \tho pr\'ecis
donne l'interpr\'etation de la calculabilit\'e \etp en termes 
de \fams de \cibos (cf. \cite{Bal} \tho 5.19).

\begin{theorem} 
\label{thFamCib1} Soit \,$P$\, un \pb \algq  cod\'e sous forme binaire. 
\begin{itemize}
\item [$(1)$] Si le \pb  \,$P$\,  est r\'esoluble en temps \,$T$\, par 
une Machine de Turing \`a une seule bande, on peut construire en temps 
\,$\O(T^2)$\, une \fam de \cibos qui repr\'esente la \fam  de \fonbs 
\,$(f_n)$\, associ\'ee \`a \,$P$.
\item [$(2)$] Le \pb  \,$P$\,  est r\'esoluble \etp par une
Machine de Turing \ssi il existe une \famu de \cibos
 \,$(\gamma_n)$\, qui repr\'esente la \fam de \fonbs \,$(f_n)$\, 
associ\'ee \`a \,$P$. 
\end{itemize}
\end{theorem}

Ce \tho est important, m\^eme s'il a l'air de se mordre un peu 
la queue, puisque la \fam \,$(\gamma_n)$\,  doit \^etre uniforme, \cad 
calculable \etp par une Machine de Turing. 

Il n'est pas trop compliqu\'e \`a d\'emontrer. Il est reli\'e \`a 
l'existence d'une Machine de Turing universelle qui travaille \etpz. Le 
fait que le r\'esultat du calcul (sur une entr\'ee \,$\mu$\, de taille 
\,$n$) au bout de \,$t$\, \'etapes est bien celui affich\'e peut \^etre 
v\'erifi\'e en ex\'ecutant soi-m\^eme le programme \`a la main, et on 
peut certifier la totalit\'e du calcul en certifiant le r\'esultat de 
chaque \'etape \itmdz. Quand au bout du compte, on dit \gui{la 
sortie  a \'et\'e correctement calcul\'ee}, on peut aussi l'\'ecrire en 
d\'etail sous forme d'un \cibo 
\,$\gamma_n$\, qui fonctionne pour toute entr\'ee de taille \,$n$. 
Il faut  un peu d'attention pour v\'erifier que tout ceci reste dans le 
cadre de la taille \pollez. C'est le m\^eme genre d'argument qui a 
permis \`a Cook de fournir le premier et le plus populaire des \pbs 
\NPcsz,
celui de la satisfiabilit\'e des \exbs (\'etant donn\'e une 
\exbz, existe-t-il une fa\c{c}on d'affecter les variables \boles 
en entr\'ee qui donne \`a l'\exn la valeur Vrai~?), \pb plus parlant 
que le \pb \NPc universel que nous avons expos\'e page 
\pageref{NPcomp} dans la section \ref{secPNP}.

\smallskip Le  \tho \ref{thFamCib1}  nous donne imm\'ediatement. 

\begin{proposition} 
\label{propBoBoU} Soit \,$P$\, un \pb \algqz,
cod\'e sous for\-me binaire. Soit 
\,$F:\left\{0,1\right\}^\star\rightarrow \NN$\, une fonction \algqz.
\begin{itemize}
\item Si le \pb \,$P$\, est dans la classe $\P$, alors il est dans  
$\Ppoly$.
\item Si le \pb \,$P$\, est dans la classe $\NP$, alors il est dans  
$\NPpoly$. 
\item Si la fonction \,$F$\, est  dans la classe $\DieseP$, alors elle 
est dans  $\DiesePpoly$. 
\end{itemize}
\end{proposition}

\ss La signification intuitive importante du \tho 
\ref{thFamCib1} est que \emph{la classe $\Ppoly$ est l'exact analogue 
non 
uniforme de la classe $\P$:} soit en effet  \,$P$\, un \pb \algq qui 
correspond \`a une \fam \,$(f_n)$\, de \fonbsz,
\begin{itemize}
\item le \pb \,$P$\, est dans la classe $\P$ 
 signifie que \,$(f_n)$\, est calculable par une \emph{\famuz} 
\,$\gamma_n$\, de \cibos (\`a fortiori la taille de \,$(\gamma_n)$\, est 
\polle en \,$n$), 
\item le \pb \,$P$\, est dans la classe  $\Ppoly$ signifie que 
\,$(f_n)$\, est calculable par une \fam \,$(\gamma_n)$\, de \cibos dont 
la taille  est \polle en \,$n.$ 
\end{itemize}

\'Etant donn\'e que $\NPpoly$ et $\DiesePpoly$ sont d\'efinis \`a partir 
de  $\Ppoly$  de mani\`ere similaire \`a la d\'efinition de 
 $\NP$ et $\DieseP$  \`a partir de  $\P$, une autre  signification 
intuitive importante du \tho 
\ref{thFamCib1} est que les classes $\NPpoly$  et  $\DiesePpoly$ sont
les exacts analogues non uniformes des  classes  $\NP$ et $\DieseP$.

La preuve qu'un \pb \agq donn\'e $P$ est $\NP$-complet
 donne en g\'en\'eral un proc\'ed\'e uniforme de 
r\'educ\-tion d'une \fam arbitraire dans
$\BNP$ (uniforme ou non) \`a la \fam de \fonbs attach\'ee \`a $P$. 
En cons\'equence on obtient l'implication \,$\P=\NP\,\Rightarrow\, 
\BP=\BNP$. Autrement dit la conjecture non uniforme $\BP\neq \BNP$ est 
plus forte que la conjecture classique  $\P\neq \NP$.

La m\^{e}me remarque vaut en rempla\c{c}ant $\NP$ par $\DieseP$.

\section{Le caract\`ere universel du permanent}    
\label{secPermUniv}
\subsubsection*{Le permanent}  
 
Par d\'efinition le \emph{permanent} d'une \maca \,$A=(a_{ij})_{1\leq 
i,j\leq n}$\, sur un anneau commutatif \,$\A$\, est le \pol en les 
$a_{ij},$ not\'e \,$\per(A)$, 
d\'efini par l'expression analogue \`a celle du \deter obtenue en 
rempla\c{c}ant les signes $-$ par les $+$: 
\begin{equation} \label{eqPerm}
 \per(A)\,=\,\per_n((a_{ij})_{1\leq i,j\leq 
n})\,=\,\sum_{\sigma\in\cS_n} \quad { \prod\nolimits_{i=1}^n 
{a_{i,\sigma(i)}}}
\end{equation}
o\`u $\sigma$ parcourt toutes les permutations de 
$\{1,\alb\ldots,\alb n\}$.
Nous consid\'erons \,$(\per_n)$\, comme une \fam de \pols \`a 
\,$n^2$\, variables sur l'anneau \,$\A$.

On ne conna\^{\i}t pas de mani\`ere rapide d'\'evaluer le
permanent d'une matrice \`a \coes entiers, ni sur aucun corps de
\cara distincte de \,$2$\, (en \cara $2$ 
le permanent est \'egal au \deter
et se laisse donc \'evaluer facilement).

Lorsque les \coes sont tous \'egaux \`a $0$ ou $1$ on peut
interpr\'eter la matrice \,$A$\, comme donnant le graphe d'une relation
entre deux ensembles \`a \,$n$\, \elts \,$F$\, et \,$G$.
Par exemple les \elts de \,$F$\, sont des filles et ceux
de \,$G$\, sont des gar\c{c}ons, et la relation est la relation 
d'affinit\'e (ils veulent bien danser ensemble). Alors le permanent de 
la matrice correspondante compte le nombre de mani\`eres distinctes de 
remplir la piste de danse sans laisser personne sur le bord.
Ainsi la \fam \,$f_n:\{0,1\}^{n^2}\rightarrow \NN$\, d\'efinie
par \,$f_n(A)=\per_n(A)$\, est une famille dans~$\DieseP$.

Le crit\`ere de Valiant (\tho \vref{propDeployP}) montre par ailleurs
que la \fam de \pols \,$(\per_n)$\, est dans $\VNPF$ sur n'importe quel 
anneau commutatif: en effet la \fam  \,$(\per_n)$\, n'est autre que la 
\fad de la \fam des \fonbs qui testent si une matrice dans 
\,$\{0,1\}^{n{\times}n}$\, est une matrice de permutation.

\subsubsection*{Deux th\'eor\`emes de Valiant sur le permanent} 
Valiant a \'etabli l'\'egalit\'e  $\VNP=\VNPF$ et il a montr\'e
le caract\`ere universel du permanent, \`a la fois en \cobi et
en \coagz.
\begin{theorem} 
\label{thVper1} 
Le calcul du permanent pour les \macas \`a \coes dans \,$\{0,1\}$\,
est $\DieseP$-complet. 
\end{theorem}

\begin{theorem} 
\label{thVper2} 
Sur un corps de \cara \,$\neq 2$, et plus \gnlt sur un anneau dans 
lequel \,$2$\, est inversible, la \fam \,$(\per_n)$\, est universelle
pour la classe $\VNP$: toute \fam dans $\VNP$ est une \,$p$-projection
de la \fam \,$(\per_n)$.  
\end{theorem}

Les preuves de ces deux \thos sont d\'elicates.
Pour le deuxi\`eme nous recommandons \cite{Bur2}.

\section{La conjecture de Valiant} 
\label{secVNP}

Le petit tableau  ci-apr\`es r\'ecapitule les analogies entre 
diff\'erentes classes de complexit\'e. Dans les colonnes  {\sc 
Bool\'een} et {\sc Alg\'ebrique}
interviennent des familles \emph{non uniformes} d'\exns ou circuits.
Dans la colonne {\sc Sim} nous indiquons si la simulation \agq
du cas \bol est connue comme \'etant sur la m\^{e}me ligne: deux points 
d'interrogation signifient
qu'on ne le croit gu\`ere possible.

Rappelons que dans la premi\`ere colonne (\cobiz) toutes les inclusions
en descendant sont conjectur\'ees \^{e}tre strictes, et
que les inclusions correspondantes dans le cas \bol (2\`eme colonne)
sont aussi conjectur\'ees  strictes.

\pagebreak
\begin{center}
{\textbf{Petit r\'ecapitulatif}} \\[1mm]\nopagebreak    
{(Analogies entre \cobiz, \bole et \agqz)} \\[3mm]\nopagebreak
\noi
\setlength{\extrarowheight}{1mm}
\begin{tabular}{|m{1.5cm}|m{2.6cm}|m{2.6cm}|m{.7cm}|
}
\hline {\sc Binaire} & {\sc ~~~Bool\'een} &   
{\sc ~Alg\'ebrique} & {\sc Sim} 
\\[1mm]
\hline
$\quad \NC^1$     &$\;\;\BPF=\BNC^1$ &$\;\VPF=\VNC^1$ &  oui 
\\[1mm]
\hline
$\quad \NC^2$   &$\quad \quad \BNC^2$ &$\quad \quad \VNC^2$ &  ~?? 
\\[1mm]
\hline
$\quad \NC$   &$\quad \quad \BNC$ &$\;\VNC=\VNC^2$ &  ~?? 
\\[1mm]
\hline
 $\quad \P$   & $\quad \quad \;\BP$  &$\;\;\VP=\VNC^2$ & ~??  
\\[1mm]
\hline
    & 
&$\VQP=\VQPF$ &   
\\[1mm]
\hline
 $\quad \NP$   & $\BNP=\BNPF$ & &  
\\[1mm]
\hline
 $\quad \DieseP$   & $\DieseBP=\DieseBPF$ &
$\VNP=\VNPF$ & oui 
\\[1mm]
\hline
\end{tabular}
\end{center}

\bigskip   Valiant  a propos\'e
la conjecture: 

\medskip 
\centerline{\it Pour  tout  corps  $\K,  \quad \quad \VP\not= \VNP.$}

\medskip 
 Cette conjecture est un analogue
 \agq non uniforme de la conjecture \algq
$\P\not= \NP$ ou plus pr\'ecis\'ement de $\P\not= \DieseP$.

Sur un corps de \cara \,$\neq 2$, vu le \tho \ref{thVper2},
cette conjecture s'\'ecrit purement en termes d'\exasz:

\medskip 
\centerline{\it Le permanent n'est pas une $p$-projection du 
d\'eterminant.}

\medskip 
C'est sur les corps finis que la conjecture semble le plus
significative, parce que la situation \agq y est le plus proche
du cas \bolz: elle n'est pas perturb\'ee par la pr\'esence d'\elts de 
taille arbitrairement grande dans le corps.

Si on disposait d'une \pcd uniforme qui r\'eduise la famille $(\per_n)$ 
\`a une \fam dans $\VP$, alors  le calcul du
permanent d'une matrice dans $\{0,1\}^{n^2}$
serait dans la classe $\P$  et donc on 
aurait $\P= \DieseP$ par le \thoz~\ref{thVper1}.

Plus g\'en\'eralement, 
le \tho \vref{thBur} montre que
$\Ppoly \neq  \DiesePpoly$ implique
$\VP(\FF_q)\neq \VNP(\FF_q)$ pour tout corps fini,
et sous l'hypoth\`ese de Riemann g\'en\'eralis\'ee,
$\VP(\K)\neq \VNP(\K)$ pour tout corps de \cara nulle.

Par ailleurs si on avait  $\P= \DieseP$, 
le calcul du
permanent d'une matrice dans $\{0,1\}^{n^2}$
serait dans la classe $\P$, donc a foriori dans  $\Ppoly$
et on aurait $\DiesePpoly=\Ppoly$, mais peut-\^{e}tre
pas pour autant $\VP=\VNP$.

La conjecture de Valiant est qu'il n'existe aucune \pcdz,
m\^eme sans l'hypoth\`ese restrictive d'uniformit\'e, qui r\'eduise
 la famille $(\per_n)$ \`a une \fam dans $\VP$. 

L'avantage de la conjecture de Valiant est qu'elle est un \pb
purement \agqz, qui parle uniquement de la taille
de la \rpn d'une certaine \fam de \pols
par des \fams de \carisz.

Comme un des aspects les plus myst\'erieux de la conjecture
 \,$\P\not= \NP$\, (cela n'a pas toujours
repr\'esent\'e un million
de dollars{\footnote{~Un milliardaire am\'ericain qui aimerait devenir 
c\'el\`ebre a propos\'e en l'an 2000 un prix d'un million de dollars 
pour celui ou celle qui r\'esoudrait le \pb  \,$\P= \NP~?$. 
Six autres  conjectures
math\'ematiques importantes sont dot\'ees d'un prix analogue. 
Un million de dollars n'est d'ailleurs pas grand chose
compar\'e \`a ce que gagne un bon joueur de football, et rien du tout
 par rapport \`a un avion furtif. 
Ceci tendrait \`a dire qu'un milliardaire peut esp\'erer devenir 
c\'el\`ebre avec un investissement tr\`es modeste. 
Notez que si vous d\'emontrez
que $\DieseP\not= \NP$, vous aurez droit \`a l'admiration de
tou(te)s les math\'e/infor-maticien(ne)s, mais vous n'aurez pas le
million de dollars correspondant \`a  $\P\not= \NP$. C'est
certainement injuste, mais c'est ainsi.}} mais cela a toujours
sembl\'e tr\`es excitant)  tient \`a la question de l'uniformit\'e
des \fams de \cibos en jeu, on contournerait  cet obstacle si on
d\'emontrait la conjecture analogue non uniforme et plus forte
$\Ppoly\not= \NPpoly$. 

Et  la forme purement \agq $\VP\not= \VNP$  serait plus \`a notre 
port\'ee.
Une preuve de  $\VP\not= \VNP$  serait un pas important
qui \'eclairerait le chemin pour une preuve de 
$\Ppoly\not= \DiesePpoly$,
qui implique $\P\neq \DieseP$. Cela pourrait
sugg\'erer enfin une preuve de $\P\neq \NP$.

Un petit ennui dans cette suite de consid\'erations
informelles: les deux points d'interrogation sur la ligne $\P$ du petit
tableau r\'ecapitualtif.

Comme  $\VQP=\VQPF$ et $\VNP=\VNPF$, la
\emph{conjecture de Valiant \'etendue}, \`a savoir:

\medskip 
\centerline{\it Pour  tout corps  $\K,  \quad \quad  \VNP\not\subset 
\VQP$}

\bni   
est regard\'ee par certains auteurs comme encore plus instructive
pour la compr\'ehension du \pb \algq analogue
$\P\not= \DieseP$. 

Sur un corps de \cara \,$\neq 2$, cela \'equivaut \`a:

\medskip 
\centerline{\it Le permanent n'est pas une $\,qp$-projection du 
d\'eterminant.}

\medskip 
Notons que B\"urgisser a d\'emontr\'e que  
\,$ \VQP\not\subset \VNP$\, 
sur les corps de \cara nulle (voir~\cite{Bur2}).

\newpage \thispagestyle{empty}

\def\chaptername{Annexe}

\chapter*{Annexe : codes Maple}
\addcontentsline{toc}{chapter}{Annexe : codes Maple}  
\markboth{Annexe}{Codes Maple}

\def\polgemin   {\textsf{\textbf{polgenmin}}}
\def\barmodif   {\textsf{\textbf{barmodif}}}
\def\Faddeev    {\textsf{\textbf{faddeev}}}
\def\interpoly  {\textsf{\textbf{interpoly}}}
\def\charpoly   {\textsf{\textbf{charpoly}}}
\def\berkodense {\textsf{\textbf{berkodense}}}
\def\berksparse {\textsf{\textbf{berksparse}}}
\def\chistodense{\textsf{\textbf{chistodense}}}
\def\chisparse  {\textsf{\textbf{chisparse}}}
\def\berkomod   {\textsf{\textbf{berkomod}}}
\def\barmod     {\textsf{\textbf{barmod}}}
\def\fadmod     {\textsf{\textbf{fadmod}}}
\def\chismod    {\textsf{\textbf{chismod}}}
\def\linalpoly  {\textsf{\textbf{linalpoly}}}
\def\polmod     {\textsf{\textbf{polmod}}}
\def\inversf    {\textsf{\textbf{inversf}}}
\def\hessenberg {\textsf{\textbf{hessenberg}}}
\def\kalto      {\textsf{\textbf{kalto}}}
\def\stre       {\textsf{\textbf{stre}}}
\def\Stra       {\textsf{\textbf{stra}}}
\def\somme      {\textsf{\textbf{somme}}}
\def\devlim     {\textsf{\textbf{devlim}}}
\def\sommod     {\textsf{\textbf{sommod}}}
\def\promod     {\textsf{\textbf{promod}}}

Nous donnons, dans les pages qui suivent, les codes 
{\sc Maple} des \algos qui calculent le \polcar et dont nous 
avons test\'e les performances.

Les codes sont \'ecrits ici dans la version {\sc Maple 6}, mais
les tests ont \'et\'e faits avec la version {\sc Maple 5}.
Les diff\'erences sont les suivantes. 
Premi\`erement la version {\sc Maple 6}
a grandement am\'elior\'e son calcul standard de \polcar
(en le basant sur l'\agbz~?). Deuxi\`emement, dans {\sc Maple 6},
le dernier objet calcul\'e est d\'esign\'e par \,\texttt{\%}\,
alors que dans {\sc Maple 5} il \'etait d\'esign\'e par \,\texttt{"}.
Enfin, dans {\sc Maple 6} une \pcd se termine par 
\,\texttt{end proc:}\, tandis que dans {\sc Maple 5} 
elle se termine par 
\,\texttt{end:}\,

\ms  
Les \algos que nous avons compar\'es sont ceux de  
\Ber am\'elior\'e (not\'e \,\berkodense), 
de \JB modifi\'e (\barmodif), de Faddeev-Souriau-Frame 
(\Faddeev), de Chistov (\chistodense) 
et leurs versions modulaires respectives (nous donnons ici 
\,\berkomod), 
ainsi que les \algos correspondant \`a la \miLz, 
celle de Hessenberg 
et celle de Kaltofen-Wiedemann
(not\'es respectivement \,\interpoly, \,\hessenberg\, et \,\kalto),  
en plus de la fonction \,\charpoly\, faisant partie du  
package \,\texttt{linalg}\, de {\sc Maple} que nous avons not\'ee 
\,\linalpoly\, dans nos tableaux de comparaison.  
Nous avons \egmt adapt\'e \,\berkodense\, 
et \,\chistodense\, au cas des matrices creuses (voir les 
codes \,\berksparse\, et \,\chisparse\, d\'eriv\'es)

Les mesures du temps CPU et de l'espace-m\'emoire  
pour chaque \algo test\'e sont prises \`a l'aide des 
fonctions \,\texttt{time()}\, et \,\texttt{bytesalloc}\, du noyau de 
{\sc Maple}.

\newpage  
{\tt
\noi \#\#\# Somme des \elts d'une liste \#\#\# 
\hsz \#\#\# ~ de fractions rationnelles ~~ \#\#\# 
\hsz \somme{\af}proc(suite{\ty}list(ratpoly)) 
\hsu  normal(convert(suite,`+`)) 
\hsz end proc:
\hsz \#\#\#

\bni
\#\#\#\#\# Berkowitz dans le cas d'une matrice dense \#\#\#\#\#
\hsz \berkodense{\af}proc(A{\ty}matrix,X{\ty}name) 
\hsz local n,r,i,j,k,V,C,S,Q; 
\hsu n{\af}coldim(A); 
\hsu  V{\af}table([1{\=}-1,2{\=}A[1,1]]); C[1]{\af}-1; 
\hsu  for r from 2 to n do 
\hsd   for i to r-1 do S[i]{\af}A[i,r] od; C[2]{\af}A[r,r]; 
\hsd   for i from 1 to r-2 do 
\hst    C[i+2]{\af}\somme([seq(A[r,k]*S[k],k=1..r-1)]);  
\hst    for j to r-1 do 
\hsq     Q[j]{\af}\somme([seq(A[j,k]*S[k],k=1..r-1)]) 
\hst    od; 
\hst    for j to r-1 do S[j]{\af}Q[j] od; 
\hsd   od; 
\hsd   C[r+1]{\af}\somme([seq(A[r,k]*S[k],k=1..r-1)]); 
\hsd   for i to r+1 do 
\hst    Q[i]{\af}\somme([seq(C[i+1-k]*V[k],k=1..min(r,i))]); 
\hsd   od; 
\hsd   for i to r+1 do V[i]{\af}Q[i] od;  
\hsu  od;  
\hsu  \somme([seq(V[k+1]*X\^{\,}(n-k),k=0..n)]); 
\hsu  collect({\%},X)  
\hsz end proc: 
\hsz \#\#\#\#\#

\newpage
\noi 
\#\#\#\#\# Berkowitz dans le cas d'une matrice creuse \#\#\#\#\# 
\hsz \berksparse{\af}proc(A{\ty}matrix,X{\ty}name) 
\hsz local n,r,i,j,k,V,C,S,Q,N; 
\hsu  n{\af}coldim(A); 
\hsu  V{\af}table([1{\=}-1,2{\=}A[1,1]]); 
\hsu  N{\af}vector(n);
\hsu  for i to n do N[i]{\af}\{\} od;
\hsu  C[1]{\af}-1; 
\hsu  for r from 2 to n do 
\hsd   for i to n do 
\hst    if A[i,r-1]<>0 then N[i]{\af}N[i] union \{r-1\} fi 
\hsd   od; 
\hsd   for i to r-1 do S[i]{\af}A[i,r] od; C[2]{\af}A[r,r]; 
\hsd   for i from 1 to r-2 do 
\hst    C[i+2]{\af}\somme([seq(A[r,j]*S[j],j{\=}N[r])]);  
\hst    for j to r-1 do 
\hsq     Q[j]{\af}\somme([seq(A[j,k]*S[k],k{\=}N[j])]); 
\hst    od; 
\hst    for j to r-1 do S[j]{\af}Q[j] od; 
\hsd   od; 
\hsd   C[r+1]{\af}\somme([seq(A[r,j]*S[j],j{\=}N[r])]); 
\hsd   for i to r+1 do 
\hst    Q[i]{\af}\somme([seq(C[i+1-k]*V[k],k=1..min(r,i)]); 
\hsd   od; 
\hsd   for i to r+1 do V[i]{\af}Q[i] od; 
\hsu  od; 
\hsu  \somme([seq(V[k+1]*X\^{\,}(n-k),k=0..n)]); 
\hsu  collect({\%},X) 
\hsz end proc;
\hsz \#\#\#\#\#
} 

\bs 
Nous avons \egmt adapt\'e les codes Maple  
ci-dessus, correspondant \`a l'algorithme 
am\'elior\'e de \Berz, au cas o\`u les 
\coes appartiennent \`a un anneau-quotient du 
type $\aqo{\ZZ_p[lisvar]}{Ideal}.$ On obtient une 
\pcdz, not\'ee \,\berkomod\, dans nos tableaux 
de comparaison, qui prend en entr\'ee un entier 
positif \,$p\,$, une liste d'ind\'etermi\-n\'ees 
\,$lisvar\,$, une liste \,$Ideal$\, de \pols en 
\,$lisvar$\, et la matrice 
\,$A\in (\aqo{\ZZ_p[lisvar]}{Ideal})^{n\times n}\,$
pour donner en sortie le \polcar de \,$A\,$.

La \pcd \,\berkomod\, ainsi que les versions 
modulaires des autres \algos utilisent comme  
sous-\pcd la \pcd \,\polmod\,  
qui prend en entr\'ee un nombre entier \,$p$,  
un \pol \,$P$\, de \,$\ZZ[lisvar]\,$,  
et donne en sortie un repr\'esentant simple de  
l'image canonique de \,$P$\, dans l'anneau-quotient  
\,$\aqo{\ZZ_p[lisvar]}{Ideal}$.  
  
{\tt\mni 
 \#\#\#\#\# R\'eduction d'un \pol  modulo un id\'eal \#\#\#\#\#    
\hsz \polmod{\af} 
\hsz proc(P{\ty}polynom,lisvar{\ty}list,Ideal{\ty}list,p{\ty}posint)   
\hsu  local i, Q;   
\hsu  if nops(lisvar)<>nops(Ideal) then 
\hsd    ERROR(`The number of polynomials must 
\hsd    be equal to the number of variables`)  
\hsu  fi;   
\hsu  Q{\af}P;   
\hsu  for i to nops(lisvar) do   
\hsd    Q{\af}rem(Q,Ideal[i],lisvar[i]);   
\hsd    Q{\af}Q mod p   
\hsu  od;   
\hsu  sort(Q);   
\hsz end proc:   
\hsz \#\#\#\#\#  
}

\ms On en d\'eduit les deux calculs de base modulo l'id\'eal
consid\'er\'e, la somme d'une liste et le produit de deux \eltsz.
  
{\tt\mni 
 \#\#\#\#\# Somme d'une liste modulo un id\'eal \#\#\#\#\#
\hsz \sommod{\af}proc(s{\ty}list(polynom),
\hsix lsv{\ty}list(name),lsp{\ty}list(polynom),p{\ty}posint)
\hsu \polmod(\somme(s),lsv,lsp,p)
\hsz end proc:
\hsz \#\#\#\#\#  

\mni \#\#\#\#\# Evaluation d'un produit modulo un un ideal \#\#\#\#\#
\hsz \promod{\af}proc(P,Q{\ty}polynom,
\hsix  lsv{\ty}list(name),lsp{\ty}list(polynom),p{\ty}posint)
\hsu \polmod(P*Q,lsv,lsp,p)
\hsz end proc:
\hsz \#\#\#\#\#  
}  

\ms Il ne reste plus qu'\`a r\'e\'ecrire \,\berkodense\, en y 
rempla\c{c}ant les op\'erations somme d'une liste de \pols
et produit de deux \pols par les calculs modulaires 
donn\'es par \,\sommod\, et \,\promod.

{\tt\bni
\#\#\#\#\# Berkowitz  modulaire \#\#\#\#\#
\hsz \berkomod{\af}proc(A{\ty}matrix,X{\ty}name,lsv{\ty}list(name),
\hsu \hspace{4cm} lsp{\ty}list(polynom),p{\ty}posint)
\hsz local n,r,i,j,V,C,S,Q; 
\hsu n{\af}coldim(A); 
\hsu  V{\af}table([1{\=}-1,2{\=}A[1,1]]); C[1]{\af}-1; 
\hsu  for r from 2 to n do 
\hsd   for i to r-1 do S[i]{\af}A[i,r] od; C[2]{\af}A[r,r]; 
\hsd   for i from 1 to r-2 do 
\hst    [seq(\promod(A[r,k],S[k],lsv,lsp,p),k=1..r-1)];  
\hst    C[i+2]{\af}\sommod(\%,lsv,lsp,p);  
\hst    for j to r-1 do 
\hsq     [seq(\promod(A[j,k],S[k],lsv,lsp,p),k=1..r-1)] 
\hsq     Q[j]{\af}\sommod(\%,lsv,lsp,p) 
\hst    od; 
\hst    for j to r-1 do S[j]{\af}Q[j] od; 
\hsd   od; 
\hsd   [seq(\promod((A[r,k],S[k],lsv,lsp,p),k=1..r-1)]; 
\hsd   C[r+1]{\af}\sommod(\%,lsv,lsp,p); 
\hsd   for i to r+1 do 
\hst    [seq(\promod((C[i+1-k],V[k],lsv,lsp,p),
\hst     \hspace{6cm} k=1..min(r,i))]; 
\hst    Q[i]{\af}\sommod(\%,lsv,lsp,p); 
\hsd   od; 
\hsd   for i to r+1 do V[i]{\af}Q[i] od;  
\hsu  od;  
\hsu  \somme([seq(V[k+1]*X\^{\,}(n-k),k=0..n)]); 
\hsu  collect({\%},X)  
\hsz end proc: 
\hsz \#\#\#\#\#
}

\bs Voici maintenant sans plus de commentaire les codes
{\sc Maple} des \algos
\chistodense, \chisparse, \barmodif, \Faddeev, \interpoly,
\hessenberg, \kalto.

\newpage

{\tt
  \mni 
\#\#\#\#\# (Chistov. Cas des matrices denses) \#\#\#\#\#
\hsz \chistodense{\af}proc(A{\ty}matrix,X{\ty}name) 
\hsz local n,r,i,j,k,a,b,C,V,W,Q; 
\hsu  n{\af}coldim(A); 
\hsu  a{\af}array(0..n,[1]); C{\af}array(0..n,[1]);  
\hsu  for i to n do a[i]{\af}normal(a[i-1]*A[1,1]) od; 
\hsu  for r from 2 to n do  
\hsd   for i to r do V[i]{\af}A[i,r] od; C[1]{\af}V[r]; 
\hsd   for i from 2 to n-1 do 
\hst    for j to r do 
\hsq     W[j]{\af}\somme([seq(A[j,k]*V[k],k{\=}1..r)]); 
\hst    od; 
\hst    for j to r do V[j]{\af}W[j] od; C[i]{\af}V[r]; 
\hsd   od; 
\hsd   [C[n]{\af}\somme(seq(A[r,k]*V[k],k{\=}1..r)]); 
\hsd   for j from 0 to n do 
\hst    b[j]{\af}\somme([seq(C[j-k]*a[k],k{\=}0..j)]); 
\hsd   od; 
\hsd   for j from 0 to n do a[j]{\af}b[j] od; 
\hsu  od; 
\hsu  Q{\af}\somme([seq(X\^{\,}k*a[k],k{\=}0..n)]);
\hsu  Q{\af}X\^{\,}n*subs(X{\=}1/X,\inversf(Q,X,n));
\hsu  Q{\af}collect((-1)\^{\,}n*{Q},X)
\hsz end proc: 
\hsz \#\#\#\#\#

\mni
\#\#\# Calcul de l'inverse modulo $z^{(n+1)}$ d'un \pol en $z$ \#\#\# 
\hsz \inversf{\af}proc(P,z,n) 
\hsu collect(convert(series(1/P,z,n+1),polynom),z,normal) 
\hsz end proc:
\hsz \#\#\# cette \pcd utilis\'ee dans les \algos  
\hsz \#\#\# de Chistov sera aussi utile dans l'\algo \kalto

\newpage\noi
\#\#\#\#\#  Chistov. Cas des matrices creuses  \#\#\#\#\# 
\hsz \chisparse{\af}proc(A{\ty}matrix,X{\ty}name)
\hsz local n,r,i,j,k,a,b,C,N,V,W,Q;
\hsu  n{\af}coldim(A);
\hsu  a{\af}array(0..n,[1]); C{\af}array(0..n,[1]); 
\hsu  N{\af}array(1..n); for i to n do N[i]{\af}\{\} od;
\hsu  \#\# N[i];  \elts non nuls de la i-\`eme ligne 
\hsu  for i to n do 
\hsd   for j to n do 
\hst    if A[i,j] <> 0 then N[i]{\af}N[i] union \{j\} fi 
\hsd   od
\hsu  od;\,\,\,\,\#\#\#\# Fin de la construction de N
\hsu  for i to n do a[i]{\af}normal(a[i-1]*A[1,1]) od;
\hsu  for r from 2 to n do 
\hsd   for i to r do V[i]{\af}A[i,r] od;
\hsd   C[1]{\af}V[r];
\hsd   for i from 2 to n-1 do
\hst    for j to r do
\hsq     [seq(A[j,k]*V[k],k{\=}\{\$ 1..r\} intersect N[j])]; 
\hsq     W[j]{\af}\somme({\%})
\hst    od;
\hst    for j to r do V[j]{\af}W[j] od;
\hst    C[i]{\af}V[r];
\hsd   od;
\hsd   [seq(A[r,k]*V[k],k{\=}\{\$ 1..r\} intersect N[r])];
\hsd   C[n]{\af}\somme({\%});
\hsd   for j from 0 to n do
\hst    b[j]{\af}\somme([seq(C[j-k]*a[k],k{\=}0..j)]);
\hsd   od;
\hsd   for j from 0 to n do a[j]{\af}b[j] od;
\hsu  od;
\hsu  Q{\af}\somme([seq(X\^{\,}k*a[k],k{\=}0..n)]);
\hsu  Q{\af}X\^{\,}n*subs(X{\=}1/X,\inversf(Q,X,n));
\hsu  Q{\af}collect((-1)\^{\,}n*{Q},X)
\hsz end proc: 
\hsz \#\#\#\#\#

\newpage
\noi
\#\#\#\#\# Jordan-Bareiss Modifi\'e \#\#\#\#\# 
\hsz \barmodif{\af}proc(A{\ty}matrix,X{\ty}name)
\hsz local B,n,p,i,j,piv,dencoe;
\hsu   den{\af}1; n{\af}coldim(A); B{\af}copy(A); 
\hsu   B{\af}evalm(B-X*array(identity, 1..n,1..n)); piv{\af}B[1,1];
\hsu   for p from 1 to n-1 do  
\hsd     for i from p+1 to n do        
\hst       coe{\af}B[i,p];    
\hst       for j from p+1 to n do    
\hsq         B[i,j]{\af}normal((piv*B[i,j]-coe*B[p,j])/den)    
\hst       od    
\hsd     od;  
\hsd     den{\af}piv; piv{\af}B[p+1,p+1] 
\hsu   od;
\hsu   sort(collect(piv,X),X)   
\hsz end proc: \,\,\,\,\, \#\#\#\#\#

\sni
\#\#\#\#\#  Faddeev-Souriau-Frame  \#\#\#\#\#
\hsz \Faddeev{\af}proc(A{\ty}matrix,X{\ty}name)
\hsz  local n, k, a, C, B, Id, P;
\hsu  n{\af}coldim(A); a{\af}array(1..n); 
\hsu  Id{\af}array(1..n,1..n,identity); B{\af}copy(Id);
\hsu  for k from 1 to n do  
\hsd     C{\af}map(normal,multiply(A,B));
\hsd     a[k]{\af}trace(C)/k; 
\hsd     B{\af}map(normal,evalm(C-a[k]*Id))
\hsu  od;
\hsu  P{\af}\somme([seq(a[k]*X\^{\,}(n-k),k=1..n)]);
\hsu  sort((-1)\^{\,}n*(X\^{\,}n-P,X);
\hsz end proc:
\hsz \#\#\#\#\#

\sni 
\#\#\#\#\# Interpolation de Lagrange \#\#\#\#\# 
\hsz \interpoly{\af}proc(M{\ty}matrix,X{\ty}name)
\hsz  local n,Id,i,j,N,d,L;
\hsu  n{\af}coldim(M); Id{\af}array(identity, 1..n, 1..n);
\hsu  for i to n+1 do  d[i]{\af}det(evalm(M-(i-1)*Id)) od;
\hsu  L{\af}[seq(d[j], j{\=}1..n+1)]; 
\hsu  interp([`\$`(0 .. n)], L, X);
\hsz end proc: 
\hsz \#\#\#\#\#

\newpage

\noi 
\#\#\#\#\# M\'ethode de Hessenberg  \#\#\#\#\# 
\hsz \hessenberg{\af}proc(A{\ty}matrix,X{\ty}name)
\hsz local jpiv, ipiv, iciv, i, m, n, piv, c, H, P;
\hsc      \# Initialisations
\hsu  n{\af}coldim(A); P[0]{\af}1; H{\af}copy(A);
\hsc      \# R\'eduction de H \`a la forme de Hessenberg
\hsu  for jpiv from  1 to n-2 do
\hsd   ipiv{\af}jpiv+1; iciv{\af}ipiv; piv{\af}normal(H[iciv,jpiv]) ;
\hsd   while  piv{\=}0  and iciv < n  do
\hst    iciv{\af}iciv+1 ; piv{\af}normal(H[iciv,jpiv])
\hsd   od ;      
\hsd   if  piv <> 0  then
\hst    if  iciv > ipiv   then
\hsq     H{\af}swaprow(H,ipiv,iciv); \# Echange de lignes
\hsq     H{\af}swapcol(H,ipiv,iciv) \, \# Echange de colonnes
\hst    fi;
\hst    for i from  iciv+1 to n  do
\hsq     c{\af}normal(H[i,jpiv]/piv) ;
\hsq     H{\af}addrow(H,ipiv,i,-c);\# Manipulation de lignes
\hsq     H{\af}addcol(H,i,ipiv,c)  \, \# Manipulation de colonnes
\hst    od;
\hst    H{\af}map(normal,H) 
\hsd   fi
\hsu  od ;
\hsc      \# Calcul du \polcar
\hsu  for m from  1 to n do
\hsd   P[m]{\af}normal((H[m,m]-X) * P[m-1]) ; c{\af}1 ;
\hsd   for i from  1 to m-1 do
\hst    c{\af}normal(-c * H[m-i+1,m-i]) ;
\hst    P[m]{\af}normal(P[m]+c * H[m-i,m]* P[m-i-1])
\hsd   od
\hsu  od ;
\hsu  collect(P[n],X)  \# le \polcar de A.
\hsz end proc;
\hsz \#\#\#\#\#

\newpage
\noi \#\#\#\# developpement limit\'e \`a l'ordre n  \#\#\#\#
\hsz \devlim{\af}proc(s{\ty}ratpoly,u{\ty}name,n{\ty}integer)
\hsu   convert(series(s,u,n+1),polynom); collect(",u,normal)
\hsz end proc;
\hsz \#\#\#\#\#

\mni \#\#\#\#\#  Kaltofen-Wiedemann  \#\#\#\#\#
\hsz \kalto{\af}proc(A{\ty}matrix,X{\ty}name) 
\hsz local n,i,j,k,a,b,bv,bw,c,B,C,P,u; 
\hsu  n{\af}coldim(A); 
\hsu  \#\#\# Initialisation  
\hsu  a{\af}\stre(n); C{\af}\Stra(n); 
\hsu  b{\af}vector(2*n); B{\af}evalm(C+u*(A-C)); 
\hsu  \#\#\# Calcul des b\_i  
\hsu  b[1]{\af}a[1];  
\hsu  bv{\af}copy(a); bw{\af}vector(n);   
\hsu  for i from 2 to n+1 do  
\hsd    \#\# multiplication de  B par bv 
\hsd    for j to n do 
\hst      bw[j]{\af}\somme([seq(B[j,k]*bv[k],k{\=}1..n)]); 
\hsd    od; 
\hsd    for j to n do bv[j]{\af}bw[j] od; 
\hsd    b[i]{\af}bv[1] 
\hsu  od; 
\hsu  for i from n+2 to 2*n do  
\hsd    \#\# multiplication de  B par bv  
\hsd    for j to n do 
\hst      bw[j]{\af}\somme([seq(B[j,k]*bv[k],k{\=}1..n)]);  
\hsd    od; 
\hsd    for j to n do bv[j]{\af}bw[j] od; 
\hsd    b[i]{\af}\devlim(bv[1],u,n);   
\hsu  od;  
\hsu  P{\af}\polgemin(b,X,u,n);  
\hsu  P{\af}sort(subs(u{\=}1,(-1)\^{\,}n*res),X);
\hsu  P{\af}collect(P,X,normal) 
\hsz end proc: 
\hsz \#\#\#\#\#

\newpage\noi
\#\#\#\#\# Sous-\pcds utilis\'ees dans \kalto\,\,\, \#\#\#\#\#

\bni
     \#\#\# \polgemin: Proc\'edure de Berlekamp-Massey ~~~ \#\#\# 
\hsz \#\#\# pour le calcul du \polgmin  \#\#\#
\hsz \#\#\# ~~~~~~ d'une \srl ~~~~~~  \#\#\#
\hsz \#\#\# ~~~Ici, l'anneau de base est l'anneau des  ~~~       \#\#\#
\hsz \#\#\# ~~~  \dlis A[z]/<z\^{}(n+1)>   ~~~      \#\#\#
\hsz \polgemin{\af}proc(b{\ty}vector,X{\ty}name,z{\ty}name,n{\ty}integer) 
\hsz local i,lc,ilc,ill,R1,R2,R3,V1,V2,V3,Q; 
\hsu R1{\af}\somme([seq(b[2*n-k]*X\^{\,}k,k{\=}0..2*n-1)]); 
\hsu Q{\af}quo(X\^{}(2*n),R,X,'R2'); 
\hsu V1{\af}1; V2{\af}-Q; ill{\af}1; 
\hsu for i from 2 to n do 
\hsd \#\#\# traiter R2  
\hsd R2{\af}collect(R2,X,normal);
\hsd lc{\af}lcoeff(R2,X); 
\hsd ilc{\af}\inversf(lc,z,n);  
\hsd R2{\af}\devlim(ilc*R2,z,n); 
\hsd Q{\af}quo(R1,R2,X,'R3');  
\hsd Q{\af}\devlim(Q,z,n);  
\hsd R3{\af}\devlim(R3,z,n);   
\hsd V3{\af}\devlim(ill*V1-ilc*V2*Q,z,n);  
\hsd ill{\af}ilc; 
\hsd V1{\af}V2; V2{\af}sort(V3,X);  
\hsd R1{\af}R2; R2{\af}sort(R3,X); 
\hsu od;  
\hsu V2{\af}collect(V2,X); 
\hsu lc{\af}lcoeff(V2,X); ilc{\af}\inversf(lc,z,n);
\hsu V2{\af}\devlim(ilc*V2,z,n);
\hsu V2{\af}collect(V2,z,normal) 
\hsz end proc:  
\hsz \#\#\#\#\#

\newpage

\noi \#\#\#\# Vecteur du centre  d'\elid \#\#\#\# 
\hsz \stre{\af}proc(n) 
\hsz local i,a; 
\hsu a{\af}vector(n);  
\hsu for i to n do  
\hsd  a[i]{\af}binomial(i-1,floor((i-1)/2))  
\hsu od;
\hsu eval(a) 
\hsz end proc:
\hsz \#\#\#\#\#

\bni
\#\#\#\# Matrice du centre d'\elid \#\#\#\# 
\hsz \Stra{\af}proc(n)  
\hsz local i,C;  
\hsu C{\af}array(1..n,1..n,sparse); 
\hsu for i to n-1 do  
\hsd  C[n,i]{\af}(-1)\^{\,}floor((n-i)/2) * 
\hsc      binomial(floor((n+i-1)/2),i-1); 
\hsd  C[i,i+1]{\af}1 
\hsu od; 
\hsu C[n,n]{\af}1; evalm(C) 
\hsz end proc: 
\hsz \#\#\#\#\# 

}

\newpage \thispagestyle{empty} 
 
\refstepcounter{chapter}
\addcontentsline{toc}{chapter}{Tables, bibliographie, index.}
\label{listeAg}
\refstepcounter{section}
\addcontentsline{toc}{section}{Liste des \algosz, circuits
et \prevsz} 
\listof{agh}{Liste des \algosz, circuits et \prevsz}

\markboth{Liste des \algosz}{Liste des \algosz}
\newpage \thispagestyle{empty} 


\refstepcounter{section}
\addcontentsline{toc}{section}{Liste des Figures}  
\listoffigures  
\newpage \thispagestyle{empty} 


\newpage\thispagestyle{empty}

\cleardoublepage \thispagestyle{empty} 

\refstepcounter{section}
\addcontentsline{toc}{section}{Index des notations}  
\markboth{Index des notations}{Index des notations}
\printindexnota  

\newpage \thispagestyle{empty}  
\cleardoublepage \thispagestyle{empty} 


\refstepcounter{section}
\addcontentsline{toc}{section}{Index des termes} 
\markboth{Index des termes}{Index des termes}
\printindex  
\newpage \thispagestyle{empty}

\end{document}